\documentclass{amsart}
\usepackage{amssymb,amsmath,theorem,euscript}

\usepackage{epigraph}
\usepackage{longtable}

\usepackage[X2,T2A]{fontenc}
\usepackage[cp1251]{inputenc}
\usepackage[english,russian]{babel}

\newcommand{\и}{{\fontencoding{X2}\selectfont\cyrii}} % и с точкой малая

\newcommand{\е}{{\fontencoding{X2}\selectfont\cyryat}} %ять малая

\newcommand{\Ф}{{\fontencoding{X2}\selectfont\CYROTLD}} % фита большая
\input{epsf}

\renewenvironment{quotation}%
{\begin{list}{}{}\item[]} {\end{list}}

\begin{document}

\newcounter{sec}
\renewcommand{\theequation}{\arabic{section}.\arabic{equation}}

\def\sm{\smallskip}

%\newcounter{fact} \def\fact{\addtocounter{fact}{1}{\scc \arabic{fact}}}

\newcounter{nomer}
\def\nomer{\refstepcounter{nomer}{\arabic{nomer})  }}

\newcounter{punct}[sec]
\renewcommand{\thepunct}{\arabic{sec}.\arabic{punct}}
\def\punct{\refstepcounter{punct}{\arabic{sec}.\arabic{punct}.  }}

\def\COUNTERS{\addtocounter{sec}{1}
              \setcounter{punct}{0}
          \setcounter{equation}{0}
                  }

 \def\ov{\overline}
\def\wt{\widetilde}
\def\wh{\widehat}
 \newcommand{\rk}{\mathop {\mathrm {rk}}\nolimits}
\newcommand{\Aut}{\mathop {\mathrm {Aut}}\nolimits}
\newcommand{\Out}{\mathop {\mathrm {Out}}\nolimits}
\newcommand{\Abs}{\mathop {\mathrm {Abs}}\nolimits}
\renewcommand{\Re}{\mathop {\mathrm {Re}}\nolimits}
\renewcommand{\Im}{\mathop {\mathrm {Im}}\nolimits}
 \newcommand{\tr}{\mathop {\mathrm {tr}}\nolimits}
  \newcommand{\Hom}{\mathop {\mathrm {Hom}}\nolimits}
   \newcommand{\diag}{\mathop {\mathrm {diag}}\nolimits}
   \newcommand{\supp}{\mathop {\mathrm {supp}}\nolimits}
 \newcommand{\im}{\mathop {\mathrm {im}}\nolimits}
 \newcommand{\grad}{\mathop {\mathrm {grad}}\nolimits}
  \newcommand{\sgrad}{\mathop {\mathrm {sgrad}}\nolimits}
 \newcommand{\rot}{\mathop {\mathrm {rot}}\nolimits}
  \renewcommand{\div}{\mathop {\mathrm {div}}\nolimits}

\def\Br{\mathrm {Br}}
\def\Vir{\mathrm {Vir}}

 %%%1. ClASSICAL GROUPS
 \def\Ham{\mathrm {Ham}}
\def\SL{\mathrm {SL}}
\def\Pol{\mathrm {Pol}}
\def\SU{\mathrm {SU}}
\def\GL{\mathrm {GL}}
\def\U{\mathrm U}
\def\OO{\mathrm O}
 \def\Sp{\mathrm {Sp}}
  \def\Ad{\mathrm {Ad}}
 \def\SO{\mathrm {SO}}
\def\SOS{\mathrm {SO}^*}
 \def\Diff{\mathrm{Diff}}
 \def\Vect{\mathfrak{Vect}}
\def\PGL{\mathrm {PGL}}
\def\PU{\mathrm {PU}}
\def\PSL{\mathrm {PSL}}
\def\Symp{\mathrm{Symp}}
\def\Cont{\mathrm{Cont}}
\def\End{\mathrm{End}}
\def\Mor{\mathrm{Mor}}
\def\Aut{\mathrm{Aut}}
 \def\PB{\mathrm{PB}}
\def\Fl{\mathrm {Fl}}
\def\Symm{\mathrm {Symm}} 
 \def\Herm{\mathrm {Herm}} 
  \def\SDiff{\mathrm {SDiff}} 
 
 \def\cA{\mathcal A}
\def\cB{\mathcal B}
\def\cC{\mathcal C}
\def\cD{\mathcal D}
\def\cE{\mathcal E}
\def\cF{\mathcal F}
\def\cG{\mathcal G}
\def\cH{\mathcal H}
\def\cJ{\mathcal J}
\def\cI{\mathcal I}
\def\cK{\mathcal K}
 \def\cL{\mathcal L}
\def\cM{\mathcal M}
\def\cN{\mathcal N}
 \def\cO{\mathcal O}
\def\cP{\mathcal P}
\def\cQ{\mathcal Q}
\def\cR{\mathcal R}
\def\cS{\mathcal S}
\def\cT{\mathcal T}
\def\cU{\mathcal U}
\def\cV{\mathcal V}
 \def\cW{\mathcal W}
\def\cX{\mathcal X}
 \def\cY{\mathcal Y}
 \def\cZ{\mathcal Z}
%%% END MATHCAL %%%%%%%%%%%%%%%%%%%%%%%%%%%%%%%%% %%%%%%%%%%%%%%%%%%%%%%%%%%%%%%%% %%%
\def\0{{\ov 0}}
 \def\1{{\ov 1}}
 %%%%%%%%%%%%%%%%%%%%%%%%%%%% %%%%%%%%%%%%%%%%%%%%%%%%%%%%%%%%%%% %%% BEGIN GOTIC
 \def\frA{\mathfrak A}
 \def\frB{\mathfrak B}
\def\frC{\mathfrak C}
\def\frD{\mathfrak D}
\def\frE{\mathfrak E}
\def\frF{\mathfrak F}
\def\frG{\mathfrak G}
\def\frH{\mathfrak H}
\def\frI{\mathfrak I}
 \def\frJ{\mathfrak J}
 \def\frK{\mathfrak K}
 \def\frL{\mathfrak L}
\def\frM{\mathfrak M}
 \def\frN{\mathfrak N} \def\frO{\mathfrak O} \def\frP{\mathfrak P} \def\frQ{\mathfrak Q} \def\frR{\mathfrak R}
 \def\frS{\mathfrak S} \def\frT{\mathfrak T} \def\frU{\mathfrak U} \def\frV{\mathfrak V} \def\frW{\mathfrak W}
 \def\frX{\mathfrak X} \def\frY{\mathfrak Y} \def\frZ{\mathfrak Z} \def\fra{\mathfrak a} \def\frb{\mathfrak b}
 \def\frc{\mathfrak c} \def\frd{\mathfrak d} \def\fre{\mathfrak e} \def\frf{\mathfrak f} \def\frg{\mathfrak g}
 \def\frh{\mathfrak h} \def\fri{\mathfrak i} \def\frj{\mathfrak j} \def\frk{\mathfrak k} \def\frl{\mathfrak l}
 \def\frm{\mathfrak m} \def\frn{\mathfrak n} \def\fro{\mathfrak o} \def\frp{\mathfrak p} \def\frq{\mathfrak q}
 \def\frr{\mathfrak r} \def\frs{\mathfrak s} \def\frt{\mathfrak t} \def\fru{\mathfrak u} \def\frv{\mathfrak v}
 \def\frw{\mathfrak w} \def\frx{\mathfrak x} \def\fry{\mathfrak y} \def\frz{\mathfrak z} \def\frsp{\mathfrak{sp}}
 %% This is Lie algebra %%% END GOTIC
%%%%%%%%%%%%%%%%%%%%%%%%%%%%%%%% %%%%%%%%%%%%%%%%%%%%%%%%%%%%%%%%%
%%% BEGIN MATHBF
 \def\bfa{\mathbf a} \def\bfb{\mathbf b} \def\bfc{\mathbf c} \def\bfd{\mathbf d} \def\bfe{\mathbf e} \def\bff{\mathbf f}
 \def\bfg{\mathbf g} \def\bfh{\mathbf h} \def\bfi{\mathbf i} \def\bfj{\mathbf j} \def\bfk{\mathbf k} \def\bfl{\mathbf l}
 \def\bfm{\mathbf m} \def\bfn{\mathbf n} \def\bfo{\mathbf o} \def\bfp{\mathbf p} \def\bfq{\mathbf q} \def\bfr{\mathbf r}
 \def\bfs{\mathbf s} \def\bft{\mathbf t} \def\bfu{\mathbf u} \def\bfv{\mathbf v} \def\bfw{\mathbf w} \def\bfx{\mathbf x}
 \def\bfy{\mathbf y} \def\bfz{\mathbf z} \def\bfA{\mathbf A} \def\bfB{\mathbf B} \def\bfC{\mathbf C} \def\bfD{\mathbf D}
 \def\bfE{\mathbf E} \def\bfF{\mathbf F} \def\bfG{\mathbf G} \def\bfH{\mathbf H} \def\bfI{\mathbf I} \def\bfJ{\mathbf J}
 \def\bfK{\mathbf K} \def\bfL{\mathbf L} \def\bfM{\mathbf M} \def\bfN{\mathbf N} \def\bfO{\mathbf O} \def\bfP{\mathbf P}
 \def\bfQ{\mathbf Q} \def\bfR{\mathbf R} \def\bfS{\mathbf S} \def\bfT{\mathbf T} \def\bfU{\mathbf U} \def\bfV{\mathbf V}
 \def\bfW{\mathbf W} \def\bfX{\mathbf X} \def\bfY{\mathbf Y} \def\bfZ{\mathbf Z} \def\bfw{\mathbf w}
 %%% END MATHBF
%%%%%%%%%%%%%%%%%%%%%%%%%%%%%%% %%%%%%%%%%%%%%%%%%%%%%%%%%%%%%%%%
 %%% BEGIN MATHBB
 \def\R {{\mathbb R }} \def\C {{\mathbb C }} \def\Z{{\mathbb Z}} \def\H{{\mathbb H}} \def\K{{\mathbb K}}
 \def\N{{\mathbb N}} \def\Q{{\mathbb Q}} \def\A{{\mathbb A}} \def\T{\mathbb T} \def\P{\mathbb P} \def\G{\mathbb G}
 \def\bbA{\mathbb A} \def\bbB{\mathbb B} \def\bbD{\mathbb D} \def\bbE{\mathbb E} \def\bbF{\mathbb F} \def\bbG{\mathbb G}
 \def\bbI{\mathbb I} \def\bbJ{\mathbb J} \def\bbL{\mathbb L} \def\bbM{\mathbb M} \def\bbN{\mathbb N} \def\bbO{\mathbb O}
 \def\bbP{\mathbb P} \def\bbQ{\mathbb Q} \def\bbS{\mathbb S} \def\bbT{\mathbb T} \def\bbU{\mathbb U} \def\bbV{\mathbb V}
 \def\bbW{\mathbb W} \def\bbX{\mathbb X} \def\bbY{\mathbb Y} \def\kappa{\varkappa} \def\epsilon{\varepsilon}
 \def\phi{\varphi} \def\le{\leqslant} \def\ge{\geqslant}

\def\UU{\bbU}
\def\Mat{\mathrm{Mat}}
\def\tto{\rightrightarrows}

\def\F{\mathbf{F}}

\def\Gms{\mathrm {Gms}}
\def\Ams{\mathrm {Ams}}
\def\Isom{\mathrm {Isom}}

\def\Gr{\mathrm{Gr}}

\def\graph{\mathrm{graph}}

\def\O{\mathrm{O}}

\def\la{\langle}
\def\ra{\rangle}

%\begin{document}

 \def\ov{\overline}
\def\wt{\widetilde}

\renewcommand{\Re}{\mathop {\mathrm {Re}}\nolimits}
\def\Br{\mathrm {Br}}

 %%%1. ClASSICAL GROUPS
 \def\Isom{\mathrm {Isom}}
 \def\Hier{\mathrm {Hier}}
\def\SL{\mathrm {SL}}
\def\SU{\mathrm {SU}}
\def\GL{\mathrm {GL}}
\def\U{\mathrm U}
\def\OO{\mathrm O}
 \def\Sp{\mathrm {Sp}}
  \def\GLO{\mathrm {GLO}}
 \def\SO{\mathrm {SO}}
\def\SOS{\mathrm {SO}^*}
 \def\Diff{\mathrm{Diff}}
 \def\Vect{\mathfrak{Vect}}
\def\PGL{\mathrm {PGL}}
\def\PU{\mathrm {PU}}
\def\PSL{\mathrm {PSL}}
\def\Symp{\mathrm{Symp}}
\def\ASymm{\mathrm{Asymm}}
\def\Asymm{\mathrm{Asymm}}
\def\Gal{\mathrm{Gal}}
\def\End{\mathrm{End}}
\def\Mor{\mathrm{Mor}}
\def\Aut{\mathrm{Aut}}
 \def\PB{\mathrm{PB}}
 \def\cA{\mathcal A}
\def\cB{\mathcal B}
\def\cC{\mathcal C}
\def\cD{\mathcal D}
\def\cE{\mathcal E}
\def\cF{\mathcal F}
\def\cG{\mathcal G}
\def\cH{\mathcal H}
\def\cJ{\mathcal J}
\def\cI{\mathcal I}
\def\cK{\mathcal K}
 \def\cL{\mathcal L}
\def\cM{\mathcal M}
\def\cN{\mathcal N}
 \def\cO{\mathcal O}
\def\cP{\mathcal P}
\def\cQ{\mathcal Q}
\def\cR{\mathcal R}
\def\cS{\mathcal S}
\def\cT{\mathcal T}
\def\cU{\mathcal U}
\def\cV{\mathcal V}
 \def\cW{\mathcal W}
\def\cX{\mathcal X}
 \def\cY{\mathcal Y}
 \def\cZ{\mathcal Z}
%%% END MATHCAL %%%%%%%%%%%%%%%%%%%%%%%%%%%%%%%%% %%%%%%%%%%%%%%%%%%%%%%%%%%%%%%%% %%%
\def\0{{\ov 0}}
 \def\1{{\ov 1}}
 
 %%%%%%%%%%%%%%%%%%%%%%%%%%%% %%%%%%%%%%%%%%%%%%%%%%%%%%%%%%%%%%% %%% BEGIN GOTIC
 \def\frA{\mathfrak A}
 \def\frB{\mathfrak B}
\def\frC{\mathfrak C}
\def\frD{\mathfrak D}
\def\frE{\mathfrak E}
\def\frF{\mathfrak F}
\def\frG{\mathfrak G}
\def\frH{\mathfrak H}
\def\frI{\mathfrak I}
 \def\frJ{\mathfrak J}
 \def\frK{\mathfrak K}
 \def\frL{\mathfrak L}
\def\frM{\mathfrak M}
 \def\frN{\mathfrak N} \def\frO{\mathfrak O} \def\frP{\mathfrak P} \def\frQ{\mathfrak Q} \def\frR{\mathfrak R}
 \def\frS{\mathfrak S} \def\frT{\mathfrak T} \def\frU{\mathfrak U} \def\frV{\mathfrak V} \def\frW{\mathfrak W}
 \def\frX{\mathfrak X} \def\frY{\mathfrak Y} \def\frZ{\mathfrak Z} \def\fra{\mathfrak a} \def\frb{\mathfrak b}
 \def\frc{\mathfrak c} \def\frd{\mathfrak d} \def\fre{\mathfrak e} \def\frf{\mathfrak f} \def\frg{\mathfrak g}
 \def\frh{\mathfrak h} \def\fri{\mathfrak i} \def\frj{\mathfrak j} \def\frk{\mathfrak k} \def\frl{\mathfrak l}
 \def\frm{\mathfrak m} \def\frn{\mathfrak n} \def\fro{\mathfrak o} \def\frp{\mathfrak p} \def\frq{\mathfrak q}
 \def\frr{\mathfrak r} \def\frs{\mathfrak s} \def\frt{\mathfrak t} \def\fru{\mathfrak u} \def\frv{\mathfrak v}
 \def\frw{\mathfrak w} \def\frx{\mathfrak x} \def\fry{\mathfrak y} \def\frz{\mathfrak z} \def\frsp{\mathfrak{sp}}
 %% This is Lie algebra %%% END GOTIC
%%%%%%%%%%%%%%%%%%%%%%%%%%%%%%%% %%%%%%%%%%%%%%%%%%%%%%%%%%%%%%%%%
%%% BEGIN MATHBF
 \def\bfa{\mathbf a} \def\bfb{\mathbf b} \def\bfc{\mathbf c} \def\bfd{\mathbf d} \def\bfe{\mathbf e} \def\bff{\mathbf f}
 \def\bfg{\mathbf g} \def\bfh{\mathbf h} \def\bfi{\mathbf i} \def\bfj{\mathbf j} \def\bfk{\mathbf k} \def\bfl{\mathbf l}
 \def\bfm{\mathbf m} \def\bfn{\mathbf n} \def\bfo{\mathbf o} \def\bfp{\mathbf p} \def\bfq{\mathbf q} \def\bfr{\mathbf r}
 \def\bfs{\mathbf s} \def\bft{\mathbf t} \def\bfu{\mathbf u} \def\bfv{\mathbf v} \def\bfw{\mathbf w} \def\bfx{\mathbf x}
 \def\bfy{\mathbf y} \def\bfz{\mathbf z} \def\bfA{\mathbf A} \def\bfB{\mathbf B} \def\bfC{\mathbf C} \def\bfD{\mathbf D}
 \def\bfE{\mathbf E} \def\bfF{\mathbf F} \def\bfG{\mathbf G} \def\bfH{\mathbf H} \def\bfI{\mathbf I} \def\bfJ{\mathbf J}
 \def\bfK{\mathbf K} \def\bfL{\mathbf L} \def\bfM{\mathbf M} \def\bfN{\mathbf N} \def\bfO{\mathbf O} \def\bfP{\mathbf P}
 \def\bfQ{\mathbf Q} \def\bfR{\mathbf R} \def\bfS{\mathbf S} \def\bfT{\mathbf T} \def\bfU{\mathbf U} \def\bfV{\mathbf V}
 \def\bfW{\mathbf W} \def\bfX{\mathbf X} \def\bfY{\mathbf Y} \def\bfZ{\mathbf Z} \def\bfw{\mathbf w}
 %%% END MATHBF

%%%%%%%%%%%%%%%%%%%%%%%%%%%%%%% %%%%%%%%%%%%%%%%%%%%%%%%%%%%%%%%%
 %%% BEGIN MATHBB
 \def\R {{\mathbb R }} \def\C {{\mathbb C }} \def\Z{{\mathbb Z}} \def\H{{\mathbb H}} \def\K{{\mathbb K}}
 \def\N{{\mathbb N}} \def\Q{{\mathbb Q}} \def\A{{\mathbb A}} \def\T{\mathbb T} \def\P{\mathbb P} \def\G{\mathbb G}
 \def\bbA{\mathbb A} \def\bbB{\mathbb B} \def\bbD{\mathbb D} \def\bbE{\mathbb E} \def\bbF{\mathbb F} \def\bbG{\mathbb G}
 \def\bbI{\mathbb I} \def\bbJ{\mathbb J} \def\bbL{\mathbb L} \def\bbM{\mathbb M} \def\bbN{\mathbb N} \def\bbO{\mathbb O}
 \def\bbP{\mathbb P} \def\bbQ{\mathbb Q} \def\bbS{\mathbb S} \def\bbT{\mathbb T} \def\bbU{\mathbb U} \def\bbV{\mathbb V}
 \def\bbW{\mathbb W} \def\bbX{\mathbb X} \def\bbY{\mathbb Y} \def\kappa{\varkappa} \def\epsilon{\varepsilon}
 \def\phi{\varphi} \def\le{\leqslant} \def\ge{\geqslant}

\def\UU{\bbU}
\def\Mat{\mathrm{Mat}}
\def\tto{\rightrightarrows}

\def\Gr{\mathrm{Gr}}

\def\B{\bfB} 

\def\graph{\mathrm{graph}}

\def\O{\mathrm{O}}

\def\gl{\mathfrak{gl}}

\def\la{\langle}
\def\ra{\rangle}

\begin{center}
{	\bf\huge
 Nikolay Luzin, his students, adversaries, and defenders}

\bigskip

{\bf \large (notes on  the history of Moscow mathematics, 1914-1936)}

\bigskip

{\sc \large Yury Neretin}
 \end{center}

 {\small This is  historical-mathematical and historical
 	notes on Moscow mathematics 1914-1936. Nikolay Luzin was a central
 	figure of that time. Pavel Alexandroff, Nina Bari, Alexandr Khinchin, Andrey Kolmogorov,
 	Mikhail Lavrentiev,
 Lazar Lyusternik, Dmitry Menshov, Petr Novikov, Lev Sсhnirelman,
 Mikhail Suslin, and Pavel Urysohn were his  students. We discuss the time
 	of the great intellectual influence of  Luzin (1915-1924), 
 	the time of
 	 decay of his school (1922-1930), a moment of
 	his  administrative power (1934-1936), and his fall in July 1936.}
 	
 	\bigskip
 
 \epigraph{
 	But the thing which served as a source  	
 	of Luzin's inner drama turned 	
out to be a source of his subsequent fame...}
{Lazar Lyusternik \cite{Lyu-3}}

\epigraph{
Il est temps que je m'arr\:ete: voici que je dis, ce que j'ai d\'eclar\'e,
et avec raison, \^etre inutile \`a dire.}
{Henri Lebesgue,
Preface to Luzin's book, {\it Le\v{c}ons sur les ensembles analytiques et leurs applications}, \cite{Lebeg}}

\epigraph{
	Прошло сто лет — и что ж осталось
	\newline
	От сильных, гордых сих мужей,
	\newline
	Столь полных волею страстей?
	\newline
	Их поколенье миновалось }{Alexandr Pushkin 'Poltava'}

 There is a  common idea  that a life of Nikolay Luzin can be a topic of
 a Shakespeare drama. I am agree with this sentence but I am extremely far from
 an intention to realize this idea. The present text   
 is  an impassive historical-mathematical and historical investigation
 of Moscow mathematics of that time. On the other hand,
 this is more a story of its  initiation
 and turning moments  than a history of achievements.
  The  notes contain English overview and the bibliography. The rest 
  is written in Russian.
  %
  %\footnote{There are different reasons for this. In particular, I am afraid to translate well-written texts of
  %such persons as  Luzin,  Alexandroff, Bari,  Golubev,
  %Lyusternik, etc., also Lebesgue, Sierpinski, Weil (there are their letters, which are known only in Russian translations).
  %Also, a precise translation of 
  %speachs of 1936 is a nonobvious technical problem}. 
  
  \bigskip

 {\sc 1. Overview} (English)

% {\sc 1. Overview} (Russian)
 
 {\sc 2. Was Luzin a victim of political persecutions?}
 
 {\sc 3. Moscow mathematics and Luzin's thesis}
 
  {\sc 4. Luzin's school: the creation and the chain reaction}
  
  {\sc 5. Decay}
  
  {\sc 6. The fate of set theory}
  
  {\sc 7. Mathematics and philosophy}
  
  {\sc 8. Students attack the teacher}
  
  {\sc 9. Over the top}
  
  {\sc 10.  Young workers, left mathematicians, and red professors}
  
  {\sc 11. The year, which disappeared from annals}
  
  {\sc 12. From a dialectic republic to the MechMath}
  
  {\parindent=2cm {\sc Addendum. Few remarks on historical context}}
  
  {\sc 13. In the next inning}
  
  {\sc 14. Historians and mathematicians}
  
  {\sc 15. Playing with fire}
  
  {\sc 16. After the shock-wave}
  
  {\sc 17. Autumn of Luzin}
  
  {\sc 18. Mathematics for  non-mathematicians.}

 % {\sc 16. Luzin's tree (English)}
 
 {\sc 19. Afterword}
  
  {\sc Bibliography}

 \section{Overview}
 
 \COUNTERS
 
 {\bf\punct Overview.} In 1914  Nikolay Luzin%
 \footnote{The complete name: {\it Luzin Nikolay Nikolaevich}, Russian: Лузин Николай Николаевич
 	(in the old Russian orthography: Лузинъ Николай Николаевичъ),
 	the French transliteration used in 1920-30s was {\it Lusin Nicolas}.} (9.12.1883--28.2.1950),
 a student of Dmitry Egorov,
 returned to Moscow after a 3-year trip to G\"ottingen and Paris.
In 1915 he published the thesis '{\it Integral and trigonometric series}'
and in a short time initiated a mathematical 'chain reaction' (a sentence 
of Lavrentiev \cite{Lavr1}) in Russia.
 Since the fall 1914 he gathered a group of young mathematicians
 for a work in theory of  functions  of a real variable and in descriptive set
 theory. The group  
 included (in different periods of time) Luzin's friends
 Vladimir Golubev, Ivan Privalov, Wac{\l}aw Sierpi\'nski, Vyacheslav Stepanov
 and Luzin's students Pavel Alexandroff, Nina Bari, Alexandr Khinchin, Andrey Kolmogorov,
 Mikhail Lavrentiev,
 Lazar Lyusternik, Dmitry Menshov, Petr Novikov, Lev Sсhnirelman,
 Mikhail Suslin, Pavel Urysohn
 %\footnote{I hope that these  heros of our story
 %	are rembered in moderm mathematics. On other persons, which appear
 %	in the paper, I present short biographical background in footnotes. The final part
 %	the bibliography contains a list of biographical references}
  (also, Valery Glivenko, Lyudmila Keldysh, Alexey Lyapunov, 
 and some others). During the fall 1920-1921 the group was called {\it Lusitania}, in the historical-mathematical literature
 this term usually is extended to the whole Luzin's school 1915-1935.

 It is commonly recongnized that Luzin was the founder of the famous Moscow mathematical school of
 XX century
 and that Luzin and Lusitanians reached in 1920-30s a great mathematical breakthrough.
Quite soon the domain of investigations was extended to probability (Khinchin, Kolmogorov,  Glivenko),
 topology (Alexandroff, Urysohn, Pontryagin), number theory (Khinchin, Gelfond, Schnirelman), fluid dynamics
 (Lavrentiev, Golubev, M.~Keldysh),
 variational and functional analysis (Lyusternik, Sсhnirelman, Kolmogorov, Raikov, Gelfand), 
 differential equations (Stepanov, Petrovsky, Lyusternik),
 logics (Kolmogorov, Glivenko, Maltsev, Novikov),
 theory of function  of a complex variable  (Golubev, Luzin, Privalov, Lavrentiev).
 Menshov and Bari (and ocassionaly other Lusitanians) continued investigations in  theory of functions of a real variable.  
Luzin himself with Novikov, L.~Keldysh, Lyapunov continued works in descriptive set 
theory. Below we discuss this process,  points of its branching, and the roles of Khinchin, Alexandroff, and Kolmogorov
in this development.

Starting late 1920s - early 1930s we observe  the first students of Luzin's students
as Lev Pontryagin, Alexandr Gelfond, Anatoly Maltsev, Andrey Tikhonov, Alexandr Kurosh, Mstislav Keldysh
(but a continuation of the story is beyond the scope of these notes).
 
 However, Luzin's group decayed  after a series of mathematical and personal conflicts
 in 1922-1930. Apparently, this process must be  regarded as natural and positive,
 but conflicts between 
 Luzin and his students were unusually hot. 
 
 \bigskip
 
 Let us say few words on another side of mathematical history of that time.
 Many Soviet scientists who became widely known as pure mathematicians 
 parallelly worked in  development of high technologies and in applications
 of mathematics for natural sciences. This was one of reasons 
 of favor of mathematics in the Soviet society. For mathematicians
 parallel activities gave a possibility for a wider view to science. On a social
 level this gave to individuals certain degrees of independence of
 a pure mathematical establishment. Also, this provided many positions 
 for ma\-the\-ma\-ti\-ci\-ans.
 
 Apparently, this process was initiated in
 the famous TsAGI (The central Aero-Hydrodynamical Institute),
 which  became an initial point of a development of the aero-cosmic industry of the Soviet
 Union.
 The Institute was founded in December 1918 by the mathematician,
 plane designer and professor of mechanics Nikolay Zhukovsky, 
 he was enlisted by
 the aerodynamicist and mathematician Sergey Chaplygin. In the mid-20s, the  complex analyst 
 Golubev  joined to aerodynamical and aero\-me\-cha\-ni\-cal  investigations.
 In 1933, Golubev 
 became the first dean of the famous Mechanical and Mathematical department
 of the Moscow State University (below {\it MechMath})...
 
 All biographies of Luzin published in 1950-1985 'forgot' his employment in TsAGI (1930-1936?).
 However, in this story, we will observe  several well-known  technologists and applied mathematicians
 who
 supported Luzin in heavy moments of his biography or mantained Luzin's memory
 after his death%
 \footnote{In particular, Alexey Krylov, Sergey Chaplygin,  Viktor Kulebakin, Alexandr Nekrasov,
 	Leonid Sretensky.}.
 
 \bigskip
 
 The time of our story was complicated and stormy. A picture of the Soviet politics of 1920-30s
 contains  a dense sequence of turning point,  mistakes, crimes, corrections 
 of mistakes, apologizes to suffered, new mistakes, new crimes, big catastrophes and big successes.
 If to draw this politics  as a graph on the paper, we get a constructive
 proof of existence of a nowhere differentiable function.  
 There were two post-revolutionary tsunamis. The first was called 'The Year of the great break', 
 it was happened in 1929-1930. The second
 tsunami came suddenly in May--July 1937 and finished in  September--November 1938
 (in Soviet Union, this was called The Year Thirty Seven).

 Some Lusitanians were political radicals, in particular they participated  in the 'Great break'
 and a local revolution in the Institute for Mechanics and Mathematics 
 of the Moscow University (1929-1930)...

 \bigskip
 
 At the first glance, at the end of 1920s fortune favoured Luzin. 
 %In 1927 he invented a sieve construction in descriptive set theory.
 In  1927
 he was elected as the corresponding member of the Soviet Academy of sciences,
 in 1929 as a real member.
 In 1928 Luzin was a vice-president of the International Mathematical Congress in Bolognia.
 In
 1930 he published the monograph on descriptive set theory, {\it Le\v{c}ons sur les ensembles analytiques et leurs applications,}
 which exposed researches of Luzin  and his students Suslin, Alexandroff,  Novikov,
 Lavrentiev, Keldysh, Selivanovsky.
 In the same  1930 there appeared the first version of Luzin's textbook on {\it Differential and Integral calculus}
 for  technical universities, later it 
 it was
 revised many times (13 editions 1930-1942, the total size of editions was 270 000, and
 7 editions of a completely new version 1946-1961, the total size 550 000).
 Apparently during 25 years Luzin's books  were the most popular mathematical text-book in Soviet engineering education...
 
 However, in mid 20s Luzin lost his significance as an intellectual leader of
  Moscow mathematics. It came the time of Khinchin, Kolmogorov, Alexandroff, Pontryagin 
 (Alexandroff's student), Gelfond (a student of Khinchin and Stepanov)... Luzin
 did not understand  new mathematics, which was created  in Moscow in  1920-30s.

 \bigskip
 
 After a long trip abroad, Luzin returned to Moscow in 1930 and met various problems in 1930-1933:
 a crisis of his own scientific program,
 conflicts with former students, semi-fantastical  situation in PhysMath Department of Moscow University%
 \footnote{We discuss this story below, due to a next turn of the politics of the Government and 
 efforts of Alexandroff, Khinchin, Golubev, and others, in  1932-33 a situation  was quickly normalized and 1930s 
 were a time of florescence of mathematics in Moscow.} in 1930-1931,
 heavy illnesses... However, after a reform of the Academy of Sciences 1933-34, his formal social weight
 started to grow...
 
 \bigskip
 
 In March 1935 Luzin was confirmed as the President of the Mathematical group of the Soviet Academy Sciences.
 We must explain the sense of this assignment.
 
 First. The Academy of Science of the Soviet Union was a high
 collegium, which leaded a development of sciences and high technologies. 
 The position of an academician was life-long.
 A social status (and salary) of an academician
 was very high during the whole Soviet time. At that time the number of members of the Academy of Sciences 
 was small (in 01.01.1936 the number of real members was 98 persons, this covered  all natural
 and humanitarian sciences and also high technologies, each person had a great weight).
 
 Second. In 1930-40s a person on an important position had a wide
 freedom of initiative and a big influence (but he was seriously responsible 
 for results of his activity and his solutions).
 
 In fact, Luzin became the main mathematician of the Soviet Union.
 He was deemed unfit this position and his activity implied a negative reaction
 of a majority of Moscow mathematical community.
 
 \bigskip
 
 Unfortunately now the name  'Luzin' in the  common opinion is associated mainly with the story
 of his fall in 1936 (and this story is painted in too bright colors), for this reason we 
 have to analyse it in these notes (Sect. \ref{s:attack}--\ref{s:gran}, \ref{s:fire}) and in the Overview.
 
 Alexandroff, the main opponent of Luzin%
 \footnote{They had  a heavy private conflict, details were hidden.}, attempted to expel Luzin from the Academy.
 Recall that a position of an aca\-de\-mi\-ci\-an was life-long. However, the Academy charter
 included an item claiming that an academician can be expeled if he does not carry
 his duties or if  his activity
 harms the Soviet state%
 \footnote{This item was introduced to the charter under a pressure of the Goverment.
 	Literally: Действительный член Академии Наук лишается своего звания, если он не выполняет обязанностей, налагаемых на него этим званием, или если его деятельность направлена явным образом во вред Союзу ССР.}. Alexandroff decided to incriminate to Luzin
 numerous scientific-ethical and scientific-organizational  accusations. 
 Serpi\'nski, who was a friend of Luzin, knew about Alexandroff's intention in July 1935 (sic!).
 
 In 3 July 1936 the main Soviet newspaper '{\it Pravda}' published an anonymous paper againist Luzin,
 apparently (for a detailed discussion see below Subsections \ref{ss:efremovich}--\ref{ss:luzin-krylov}),
  the attack was initiated by
 Alexandroff (but certainly he was not an author of the paper, on the other
 hand the text reflected points of view of
 some other mathematicians). This was a signal for a general offensive.
 The Academy organized a special Committee on 'Luzin affair' that worked 7-15 July 1936. In
 the stenographic record of the Committee,  former Luzin's students Alexandroff, Khinchin, Kolmogorov,
 Lyusternik, Schnirelman and  other mathematicians (including Alexandr Gelfond, Otto Shmidt, and Sergey Sobolev)
 tried to
 prove that Luzin was hardened in sin. At the same time, Luzin was attacked by mathematicians on meetings
 in the Steklov Institute and the Moscow University (in particular, by Lev Pontryagin,
  Felix Gantmakher, mechanician Nikolay Buchholtz).
  
 Accusations of  Alexandroff and Kolmogorov
  were  scientific and ethical%
  \footnote{According the Stenographic Record, a participation of Kolmogorov was minimal,
  Alexandroff several times defended Luzin from political accusations.}. 
 It is sad to see that 
  some participants of the attack pronounced  speachs in a spirit of a 'political mistrust'
 and   invented new political charges.
 
 Our knowledge of the future in past 
 shows that this was a playing with fire. The expelling of Luzin from the Academy in 1936
 could produce extremely heavy consequences in the next, Thirty Seven, year.
 However this future was unknown to mathematicians and politicians
 (as main heroes of this year Nikolay Ezhov, Joseph Stalin, Robert Eihe,  Efim Evdokimov, Nikita Khrushchev, Stanislav Kosior,  Levon Mirzoyan, Pavel Postyshev, Mikhail Frinovsky...).
 Certainly, all%
 \footnote{All but one...} participants of the attack had no intention
 to imprison Luzin  and could  not imagine such an opportunity
 (the Code Penal ot that time was rigid but the list of  'crimes' of Luzin was  outside  it).
 But a danger of a further development of 'Luzin affair' in a political direction
 was serious even in 1936 (see a discussion below). Now it is known that
 the editor-in-chief of 'Pravda' Lev Mekhlis%
 \footnote{Mekhlis (1889-1953) {\bf later} became famous as a pitiless, fanatical 
 	(and  honest a certain sense)
 	figure.
 	 In 1940-1950 he was  a minister [нарком] of the Soviet Control (this ministry was the main
 	  anti-corruption
 	structure in the Soviet Union until 1950), in 1938-9.1940 and 7.1941-5.1942 he was the main political commissar in the Red Army.
 	After mid 50s he became one of reference demons in expositions of the  Stalin
 	period of the  Soviet history (apparently, his negative fame is  a bit stretched).
 	 In 1936 Mekhlis was not a person, who solved the
 	problem under the discussion, however there was a potential danger of his influence.}
  was displeased with the final resolution 
 of the Academic Committee (supported by higher authorities), since its statements were essentially softer
 than 'Pravda's clams.

 Sergey Bernstein and the naval engineer Alexey Krylov capably defended Luzin
 in the Academic Committee.
 Not all Lusitanians joined to the attack. We know that Menshov
 and Bari resisted openly. Apparently, their reactions were predictable.
 Initiators of the attack reckoned on a joining of Lusitanians Novikov and Lavrentiev,
 who were offend by Luzin. They avoided a joining, as a result
 accusers could not prove one of their key ethical charges.

  Obviously, the mathematical squabbling was potentially dangerous for 
  the whole Academy%
  \footnote{From the letter of Vladimir Vernadsky to Alexandr Fersman
  	(the most famous Russian geologists of that time):  {\it we will go down the sloping plane}.}.
  We know that several powerful scientists tried to stop the dangerous developments.

 \bigskip
 
 Finally, Luzin lost his position of the 'main mathematician', 
 in 1939 Kolmogorov occupied it. Certainly, the both events were positive
 for mathematics. Fortunately, Luzin was not expelled from the Academy.
 The anti-Luzin company had a character of a shock-wave.
 Its aftershocks were minimal (but they were strongly overestimated in the common
 opinion after 1991).
 
 \bigskip
 
 Reading documents of 1936, we observe lot of famous respectable persons on the both sides of the front.
 Later the Soviet mathematical community tried to forget this story. Positions of the both sides
 seem clear. The party of winners (which was more powerful)
 had reasons to hide details (see below quotations of the Stenographic record of 1936). On the other hand, Luzin's
 side  preferred to forget the unhappy Luzin's mathematical kingship 
 1935-1936 and numerous preceding  heavy conflicts with Luzin  participation.
 This produced a white gap in the literature around Luzin%
 \footnote{Of course, an oral tradition existed.
 	The author of these notes enrolled at MechMath
 	in	1975. Quite soon, I observed that all speeches about Luzin are ultra-positive
 	but some worthful faces became  wry hearing this (by mathematical reasons I joined to 
 	doubters). I remember 
 sentences that Luzin was very powerful and was	dethroned by Alexandroff
 in 1936. I remember approximate quotations <<печатает в СССР ерунду, а свои лучшие работы отправляет за границу>>
 (from 'Pravda') and <<высокая оценка научной работы  Лузина
 не меняет нашего  отношения
 к его аморальным поступкам>> from Alexandroff (in Subsect. \ref{ss:detail}.2 we observe the true origin of this sentence). 
 There was a general certainty that after 1936 foreign publications of Soviet mathematicians were stopped
 (it is not a true, they were interrupted later by other reasons). 
  I had no interest to learn more.}.
  However, several papers commented a prehistory
 of 1936 were published by Lyusternik (1967),
 Lavrentiev (1974), Alexandroff (1979-1980).
 
 A hidden
 dissent
 	around Luzin's name was stopped by Kolmogorov only after 
 the death of Alexandroff (1982). During
 the next three years (1983-1985) the Soviet mathematical community
 widely celebrated  the centenary of Nikolay Luzin. This was
 a  tardy but normal and human end of the story.
 
 \bigskip

 In my opinion, an intellectual history of the  Moscow mathematical school
 1915-1940 is more interesting and more pleasant topic than a story
 of a heavy scientific and personal conflicts.
 
 Unfortunately, the end in 1983 was not an end.

 \sm

 {\bf\punct Historical sources of 'Luzin affair'.}  Our story before 1930 is a usual history of mathematics
 (with its own problems and tragedies).
 However the conflict of 1936 is a problem of both history of mathematics and  history in the usual sense.
 It requires an investigation from the both sides.
 
 Historical sources about 1936 that are available now are not perfect. But the collection of sources 
 is huge and it gives a possibility of investigation by logical tools of historicas sciences.
 
 The list of historical sources contains the followng groups of  texts.
 
 \sm
 
 a) The {\bf  Stenographic Record}  of the meetings of the special Committee  of the Academy of Sciences 7, 9, 11, 13, 15 July 1936
 on the Luzin affair, \cite{Sten}. It is a 200-pages battle scene of 1936 containing lot of records about 
 Moscow mathematical life 1911-1936.
 
 \sm
 
 b) Papers of Lazar Lyusternik \cite{Lyu-yubilej}--\cite{Lyu-4}, which are a mixture of recollections
 and a historical investigation. The author carefully avoids the 'dangerous'  1936 year. 
 However, he  tell us a lot about its prehistory. Apparently, \cite{Lyu-3}
 presented an official version of the leading group of Moscow mathematics of 1960s. 
 
 \sm 
 
 c) The biography of Luzin \cite{TyulinaA} composed by Nina Bari%
 \footnote{The typewritten manuscript was not signed. However, the autorship
 is doubtless, it is explained in the introductory paper of the editor.} in 1950.  
 It was discovered at the end of 90s and published by Anna Tyulina in 2007.
 
 \sm
 
 d) Recollections of participants of the conflict, Pavel Alexandroff \cite{Alex-auto1}, \cite{Alex-auto2} 
 and Lev Pontryagin
 \cite{Pon1}, \cite{Pon0}. Recollections of beholders, Vadim Efremovich \cite{Efrem} and Boris Gnedenko \cite{Gned}.
 
 \sm
 
 e) Reminiscences of Mikhail Lavrentiev-jr., the son of Mikhail Lavretiev, \cite{Lavr-jr},
 and of Sergey Novikov, the son of Petr Novikov and Lyudmila Keldysh, \cite{NovS}.
 Also, late reminiscences of Alexei Gladky \cite{Gla}, a student of Petr Novikov.
 
 \sm
 
 f) Letters of Wac{\l}aw  Sierpi\'nski, Henri  Lebesgue, Arnaud Denjoy
 of 1936
 published by the French historian Pierre Dugac \cite{Dug} in the Russian periodical
 'Istoriko-Matematicheskie Issledovaniya'.
 Also, a letter by Andr\'e Weil to  Dugac, 1978 (in \cite{Dug}). It seems that 
 these letters are known only in Russian translations%
 \footnote{There is the Polish translation from Russian \cite{Duda}.}.
 Also, a comment  \cite{Weil} of Andr\'e  Weil to one of his papers.
 
 \sm
 
 g) Papers in the Soviet political and scientific periodicals of 1936-1937,
 namely,  Pravda \cite{Pravda}, Russian mathematical surveys, \cite{Uspehi-37-1}, \cite{Uspehi-37-2},
 Bulletin of the Academy of Sciences \cite{Vestn2},
 Front Nauki i techniki \cite{Front1936-1}--\cite{Front1936}.
 
 \sm
 
 h) Two letters of Lev Mekhlis, the editor-in-chief of <<Pravda>>, to Stalin and Molotov of July 1936,
 published in \cite{Delo}.
 
 \sm
 
 i) Several letters of 1936: Nikolay Luzin -- Aleksey Krylov 12.02.1936 \cite{Erm1},
 Nikolay Luzin to high authorities July 1936 \cite{Delo},
 Petr Kapitsa -- Vyacheslav Molotov 06.07.1929, \cite{Kap},
 Vladimir Vernadsky -- Alexandr Fersman 03.07.1936, 07.07.1936 \cite{Vernad-fersman},
 Nikolay Nasonov -- Alexandr Fersman  07.07.1936, \cite{Delo},
 and  Sergey Chaplygin-- Vladimir Vernadsky, 11.07.1936 
(see \cite{Yush2}, \cite{Vernad-perepiska}). I do not think that all letters of this kind were discovered by historians.

j)   Vladimir Vernadsky diary \cite{Vernad-dnevniki}.
 
 \sm
 
 k) Numerous papers about Luzin, 1946-1985 (biographies, reminiscenses, surveys of works, etc.).
 This family of texts naturally splits into several subfamilies:
 
 \sm
 
 %--- Two texts of 1920s, \cite{,} \cite{Kryl-phil} related to his
 %winning of elections to the Academy.
 
 --- Two texts 1946-1948, a congratulation paper \cite{Uspehi-46} and brochure \cite{N.N.L}.
 
 \sm
 
 --- Obituaries 1950,  \cite{Obit-2}, \cite{Obit-1}, \cite{Kuleb}, \cite{Luzin-mat-v-shkole}.
 
 \sm
 
 --- Memorial papers, 1951-1953, \cite{GoBa}, \cite{BaLyu1}, \cite{BaLyu2},
 \cite{GolKuz}, \cite{KeldNov},  \cite{Fedor}, \cite{L-Sret}, \cite{Lyap-vvedenie}, also a preface 
 \cite{Nov-red} published in 1958.
 
 \sm
 
 --- Ocassional papers 1960 -1966, \cite{??},  \cite{Bari-lyudi-nauki}, \cite{Mink}, \cite{Kuz1965}, 
 \cite{Romanov-Luzin}
 
 \sm
 
 --- Papers of 1974 in occassion of Luzin's 90th birthday, \cite{Lavr1}, \cite{Kuz-Luz}, \cite{Keld-Luzin}.
 
 \sm
 
 --- Two papers outside a classification:  paper of Alexandroff \cite{Alex-kvant}-\cite{Alex-v-shkole}, 1977,
 and recollections of
 Menshov
 \cite{Menshov-L}, 1983.
 
 \sm
 
 --- Papers of 1983--1985 after Kolmogorov's turn,
 \cite{Kuz2}, \cite{Besk}, \cite{Menshov-L-2}, \cite{Shchegol},
 \cite{Kolmogor},  \cite{Vestn3}, \cite{UspKan}, \cite{Kuz2} (a collection of papers),
 \cite{Usp}, \cite{Kanov}, \cite{Ulyan}, \cite{Ulya2}, \cite{MSU}, \cite{Usp}, \cite{Kvant},
 and also  \cite{Ulyanov-Kemerovo}, \cite{Halamaizer}, \cite{Ulyanov-mshk}.
 
 \sm
 
 --- Electronic archives of the Siberian branch of the Russian
 Academie of Sciences contain three typescripts of Lyapunov 
 \cite{Lyapunov-rukopis},   \cite{Lyap-Luzin}, \cite{Lyap-Luzin-Novikov}, 
 which were not published in that time,
 and unpublished typescript of Petr Kuznetsov
 \cite{Kuz-neopublikovano}.

\sm
 
 At a first sight, this collection of apologies of Luzin contains nothing about 1936.
 This  is not the case!
 
 \sm
 
 l) In 1966 Ivan Petrovsky (who was a rector of the Moscow State University
 at that time) initiated a collecting of  audio interview
 of Soviet scientists on  history of the first third  of XX century.
 In 2004 several mathematical interview (including Alexandroff, Boris Delauney,
 Lyusternik, Menshov,
 Pavel Yushkevich)  were published 
 in \cite{Duvakin}. They do not contain revelations,
 however they are interesting from our point of view.
 
 \sm
 
 {\bf \punct The historiography.} An epic scientific history of Moscow mathematics of 1910-1940
 was firstly represented by Boris Gnedenko \cite{Gnedenko-book}, and by Alexandroff, Stepanov, 
  Gnedenko
 in \cite{AGS}, 1948. Later there  were  many other expositions  developing this paper and its details.
 The work \cite{Lyu-1}--\cite{Lyu-4} of  Lyusternik  
 was a kind of a monograph split into a series of papers 
 (also his joint works with A.~F.~Lapko 
 \cite{LL-conferences-1}--\cite{LL} are important). In general, we follow 
 the traditional scheme (but interrupt our expostion imeddiately
 after   the start of the mathematical 'chain reaction'),
 in addition we  discuss a heavy pressure of non-solvable set-theoretical problems 
 to Luzin and his group after 1925 and the opinion of Luzin on this topic (in early 40s
 Novikov
 broke through to logic  but this is far outside our notes).
 
 The second part of our bibliography \cite{Alexandroff-Hopf}--\cite{Zhukovsky-golubev-2} contains biographical  references
 about Russian mathematicians  of 1900--1936. 
 
 Numerous biographical and semi-biographical papers about Luzin were published in 1950-1985
 (the list is Subsect. \ref{ss:1950-1953}--\ref{ss:1960-1985}).
 However these works carefully avoid details of Luzin's biography after 1929.
 Apparently, this idea belonged to Golubev,  it was perfectly realized
 in the Luzin biography composed by Golubev and Bari \cite{GoBa} 1951 (a type-written outline  of 
 Bari for that paper was accidentally preserved \cite{TyulinaA}, it is the only known
 non-self-censored big text  about Luzin written by a contemporary). 
This produced a big white piece in Luzin's biography (even now the formal CV of Luzin
is not completely clarified). 

\sm

After deaths of all participants of the squabbling of 1936, 
 due the
Restructing (the Perestroika, which also was a restructing of the 'unpredictable past
of Russia')  
a tragical intellectual hero Luzin 
 was transformed to a person of  political history.

 In 1989  the  well-known historian of mathematics Adolf Yushkevich 
 ({\bf 1906}--1993) published a paper    {\it 'The academician N.~N.~Luzin affair'} \cite{Yush1}.
 It was the first paper telling about events on the Physical-Mathematical 
 department of the Moscow University in 1929-1931%
 \footnote{This story was   officially 
 	hidden  due to the fall, arrest, exiling, and death of Dmitry Egorov
 	(the Soviet State of 1950-70s criminated polytical arrests that happend starting December 1934,
 	Egorov was arrested in 1930).
 	Recall that Egorov was 
 	the adviser of Luzin and Ivan Petrovsky, the author of 'the Egorov theorem', and the main executive in Moscow mathematics
 	in 1923-1929. The story was also non-officially prohibited due to the 
 	roles of several distinguished  mathematicians in these years
 	 (and also due to the roles of founders of the Soviet
 	school of history of mathematics). The 
 	 fall of Egorov is discussed in Subsections \ref{ss:egorov-fall}, \ref{ss:egorov-arest}, \ref{ss:losev}.
 	We relatively confidently can show that the arrest of Egorov was not related 
 	to his professional activity and professional problems.}.
  Yushkevich 
 said that the story of 1936 was  inspired by the Government and this was a lesson
 of Stalin for the Soviet intelligence.  Claims of this type were obligatory
 for  texts published during the Perestroika.
 In general, Yushkevich presented  a good exposition of some events around mathematics.

 In 1991 he published under the same title  {\it 'The academician N.~N.~Luzin affair'}
 a significantly different paper
 \cite{Yush2} (see a comparison in Subsect. \ref{ss:yushkevich-89-91} and
 \ref{ss:kolman-growth}). In particular,  he claimed that  Luzin's
 story was supervised by  Arnost Kolman%
 \footnote{An adventurer with a little-investigated biography (1892--1979),
 	see below Subsect. \ref{ss:kolman-memoirs}, \ref{ss:kolman-biography}. An Austro-Hungarian citizen, got into
 	the Russian captivity,
 in 1917-18	joined to Bolsheviks. An aggressive  Marxist philosophical publicist
 in 1930s.
  In 1931 - early 1932 he occupied the position 'the President of the Association of Institutes for Natural Sciences of the Communist Academy'.
  There were few institutes in the Association, but this position in
  1931 was quite serious. The Association finished its existence in 1932 (and 1932 was a beginning
  of the end of the Communist Academy itself). 
He was the head of the Science Department of the Moscow City Committee of the Communist Party in Spring (?) 1936 -- Spring 1937.
 He leaved memoirs \cite{Kolman-vospominaniya}, I compared his stories of 1920-30s with some
 formal historical facts and immediately observed lot of obvious bragging (see Subset. \ref{ss:kolman-memoirs}).
\newline
In the modern historical literature he is regarded as '{\it one of the main ideologists of the Stalin time}'. However
in \cite{Sten} we read two ireful speaches of Gleb  Krzhizhanovsky (a companion-in-arms of Lenin, a power engineer, a member of the Central Committee
of the Communist Party, and a vice-President of the Academy) addressed to Luzin with the following sense: 'Who can believe that you, a
Soviet academician, afraid of Raikov, Kolman,  and Kagan?? This is irresponsibly!' See Subsect. \ref{ss:krzhizh}. Apparently, such claims indicate the real level of Kolman in a social hierarchy.} and draws a picture of semi-innocent
 mathematicians. For the latter purpose he assiduously selected quotations from the Stenographic report
 \cite{Sten}, which was not published in that time...

  The stenographic report was published in  in 1998, its publication required a new transformation
  of the 'non-predictable past'. This was done by the historian of mathematics Sergey Demidov and
 the archivist
  Vladimir Esakov in the paper {\it The case of Luzin in a collective memory of scientific community.} \cite{DeEs}.
   Now  Alexandroff, Gelfond, Khinchin, Kolmogorov, 
  Lyusternik, Schnirelman
  were transformed to  {\it young friends and ideological team-mates of Kolman} who participated in an  intrigue of  Kolman and Lev Mehlis and in the Stalin lessons for scientists.
  The '{\it collective memory of scientific community}' before \cite{DeEs} never heard such things.
  
  Numerous later historical texts on 'Luzin affairs' usually are 
  popularization
  of
  \cite{Yush2} and \cite{DeEs}.
   
  New works about Moscow mathematical conflict
  of 1936 were published by Semen Kutateladze \cite{Kut1}--\cite{Kut3}, he returns to the traditional
  idea of an attack to Luzin from the mathematical community.
  
  The author of the present notes have to analyze  logical constructions of
  \cite{Yush2} and  \cite{DeEs}. We also use numerous historical sources
 published during the last 20 years (the Dugac collection \cite{Dug}, Vernadsky diaries \cite{Vernad-dnevniki},
  the manuscript of Bari \cite{TyulinaA}, and
 Pontryagin \cite{Pon1}, Efremovich \cite{Efremovich}, Lavrentiev-jr. \cite{Lavr-jr} recollections),
  which confirm the traditional point of view
 existing before 1989.
    We accept an essential part
  of Kutateladze's logical argumentation. However, in my opinion 
  he  completely loses the scientific side of the conflict,
  also his historical view contains some aberrations related 
  to a linear interpolation of the nowhere differentiable function mentioned above.

  \bigskip
  
  These notes contain lot of quotations. They have different nature. On one hand
  we are trying to tell history by voices of its participants and creators,
  on the other hand  we must discuss certain non-obvious cases and this
  requires quotations  for analyzing of information. 
  
   Now many historical sources
  are easily available in the Internet. The most of Russian mathematical journals are available
  via the cite {\it math-net.ru} created by the Steklov Institute 
  \begin{center}
  	\begin{verbatim}
  http://www.mathnet.ru/index.phtml?&option_lang=eng
  \end{verbatim}
  \end{center}
Papers from other mathematical journals 1910-1940 usually can be easily find in a standard way.
  
  A collection of  sources related to the history of Russian mathematics of  1920-30 now is present on my cite in Vienna:
  \begin{center}
  	\begin{verbatim}
 http://www.mat.univie.ac.at/~neretin/misc/luzin/notes.html
  \end{verbatim}
\end{center}
 
 I  met lot of unexpected data  in new historical sources during
 all time of work on these notes, these data broke  some my conjectures,
 so I do not think that this can not happened again.

 \section{Был ли Лузин объектом политических преследований?%
 \label{s:victim}}
 
 \COUNTERS
 
 \righthyphenmin=2
 
 Этот раздел крайне формален. Если человек был жертвой политических преследований или репрессий, 
 то это должно наблюдаемо. 
 Проверим это по внешним признакам.
 В июле 1936г. он подвергся травле со стороны газеты <<Правда>>.
  Известно, что в 1930--1931 Лузин был объектом выпадов
 около-математической партийной тусовки (С.~А.~Яновская, Э.~Я.~Кольман, М.~Я.~Выгодский).
 %, которые вряд ли были очень эффективны
 %(например, они без какого-либо успеха пытались громить учебник Грэнвиль-Лузин, о %чем ниже). 
 Документы на Лузина
 были у ОГПУ по делу П.~А.~Флоренского в начале 30х (но в ход не пошли).
 А были ли еще какие-нибудь проявления преследований?  И как  эти нападки отозвалось на его биографии?
  
  \sm
  
  {\bf \punct Лузин, 1931-1936.%
  \label{ss:luzin-do-1936}} В статье С.~С.~Демидова и В.~Д.~Есакова  \cite{DeEs}
 Лузин представлен как предмет политических атак с 1930 года.
  В качестве доказательства авторы приводят письмо Э.~Я.~Кольмана 1931 со словами <<секретно>> 
  (мы обсудим деятельность Кольмана ниже пп.\ref{ss:kolman-memoirs}, \ref{ss:kolman-biography}, письмо тоже приведено ниже
  в п.\ref{ss:kolman-tsagi}). 
  Затем они интерполируют это письмо с атакой на Лузина в 1936 году... Посмотрим на биографию Лузина в промежутке 1930-1936.

  \sm
  
$\bullet$  16.12.1930  --  уход из МГУ.
  
$\bullet$  В 1929–1934 гг. работал в Математическом отделе Физико-математического института АН СССР
в Ленинграде%
 \footnote{В 1934г. Физико-математический институт был разделен на МИАН и ФИАН. В тот же год МИАН разделился на  основной
 	институт в Москве (МИАН) и  Ленинградское отделение (ЛОМИ), которое
 	в какой-то момент стало административно самостоятельным.};
 6.02.1931г  датировано удостоверение заведующего отделом   \cite{arhiv};
 в 1934–1937 гг. — заведующий Отделом 
теории функций действительного переменного МИАН (Ленинград),
\cite{Koz}. В 1934-1936 - заведующий отделом теории функций в МИАН в Москве, \cite{TyulinaA}
(я  не уверен, что в  перечисленных данных нет какой-то путаницы.).
По-видимому, Лузин продолжал жить в Москве\footnote{Из \cite{TyulinaA}:
{\it Н. Н. жил в Москве и до 1934г. только эпизодически появлялся в Ленинграде, где была Академия.}}.

%$\bullet$.
   
   \sm
  
  $\bullet$ С 20.11.1930 работал в ЦАГИ (теоретическая группа Чаплыгина), ушел оттуда не ранее 1936,
 см. \cite{Sten} и примечания там.
  
  \sm
  
  $\bullet$ В 1930 году выходит первое издание учебника Грэнвиль%
  \footnote{На обложках разных изданий встречается как написание <<Грэнвиль>>, так и <<Гренвиль>>.}-Лузин <<Курс дифференциального
  и интегрального исчисления>>, который, по-видимому, стал самым популярным учебником высшей математики
  для инженерных институтов. Список дальнейших изданий   1930, 1930, 1930, 1931, 1933, 1933, 1934,
  1934,  1935, 1935-1937; {\bf 1937}, {\bf 1938},  {\bf 1942}, общий тираж с 33 года - 270 000 экземпляров
  \cite{Kol-savvina-2}, кроме того было издание 1931 года на украинском языке. В начале 1936 года
  учебник был рекомендован конференцией Всесоюзного комитета по высшему техническому образованию при ЦИК СССР
  (всего в списке было 7 учебников и задачников%
  \footnote{Список не безынтересен: И.~И.~Привалов «Аналитическая геометрия»,
  	 В.~Э.~Гренвиль и Н.~Н.~Лузин
  	«Курс дифференциального и интегрального исчисления», В.~И.~Смирнов
  	«Курс высшей математики для физиков и техников», т. 1—2,
  	 О.~Н.~Цубербиллер «Задачи и упражнения по аналитической геометрии», Я.~С.~Дубнов
  	«Задачи и упражнения по дифференциальному исчислению», «Сборник
  	задач по высшей математике» под редакцией Н.~М.~Гюнтера и Р.~О.~Кузьмина, Г.~М.~Фихтенгольц «Математика для инженеров».}
  	\cite{Lapko}).

  \sm

  $\bullet$  В 1932 году Хинчин (только что ставший директором НИИ механики и математики  при МГУ)
  приглашал Лузина в МГУ, Лузин отказался.
  Из письма Лузина Выгодскому 14.09.1934:
  \begin{quotation}
  	Никакого обмена письмами между Дирекцией института [механики и математики] и мной не
  	было, кроме обычных, {\bf делавшихся еще и Хинчиным и другими десятки раз
  	приглашения вступить в Институт} и моего ответа на это, ответа, который я 
  	только и мог дать ввиду хотя бы только состояния моего здоровья.
  \end{quotation}
  
   \sm
   
  $\bullet$ В 1932г.    Лузин не поехал на математический конгресс в Цюрихе (5-12.09.1932).
  В  \cite{DeEs} обсуждается несколько документов, связанных с разрешением или не разрешением Лузину
  поехать на Конгресс, очевидно кто-то противодействовал его поездке. Кончается это так:
  \begin{quotation}
  	23 августа 1932 г. Секретариат ЦК разрешил командировку в Цюрих академику Лузину, 
  	и 25 августа это решение было утверждено Политбюро ЦК. Но на конгрессе H.~Н.~Лузин не был.
  	Дорога за границу ему была закрыта. 
  \end{quotation}
  
  Процитируем письмо Колмогорова Александрову%
  \footnote{Переписка Колмогорова и Александрова была опубликована А. Н. Ширяевым
  в \cite{Kolmogorov-sh2}. Ниже мы многократно эту книгу используем.} от 15.10.1932:
  \begin{quotation}
  	Еще необычайные происшествия со съездом (по рассказам В.~В.~Степанова).
  	Н.~Н.~Лузина 31 августа вызвали в МК и сказали <<Поезжайте на съезд>>.
  	Н.~Н. ответил, что уже поздно. 3 сентября его вновь вызывали в МК партии
  	и сказали <<Поезжайте на аэроплане>>. Н.~Н. ответил: Теперь и на аэроплане поздно.
  \end{quotation}
  
  Выходит, что окончательное решение было принято самим Лузиным,
  и стоит иметь в виду, что в том же сентябре 1932 он надолго лег на лечение.

   \sm
  
  $\bullet$ В 1932-33гг. тяжело болел, провел много месяцев в санатории в Крыму и в Кремлевской больнице
  (см. примечания к \cite{Sten}, из Кремлевской больницы было также отправлено
  одно из писем Лузина Крылову \cite{Erm1}).

  \sm
  
  $\bullet$ В 1933, 1935 Лузин публикует две небольшие монографии,
  <<Современное состояние теории функций действительного переменного>> \cite{Luzin-teor-funkts}, 60 стр.,
  и <<О некоторых новых результатах дескриптивной теории функций>> \cite{Luzin-descript}, 88 стр.

  \sm
  
  $\bullet$, 1935-1936, по приглашению декана Л.~А.~Тумаркина работал на Мехмате МГУ%
  \footnote{Интересно, что он был членом НИИ механики, а не НИИ математики. 
  Понятно, что у Лузина были в среднем плохие отношение с мехматскими математиками, а с механиками хорошие.}.

  \sm
  
  $\bullet$  1934-36, председатель квалификационной математической комиссии  \cite{DeEs} (квалификационная комиссия - предшественник ВАКа).
  Напомню, что    ученые степени были отменены в 1918г.% 
  \footnote{Декрет Совета Народных комиссаров <<О некоторых изменениях
  в устройстве государственных учебных и высших учебных заведений РСФСР>> от
  01.10.1918 (за подписью  ЗамНаркомПроса 
Мих. Покровского,
Управ.делами Совнаркома
В.~Л.~Бонч-Бруевича и
секретаря Совнаркома
Фотиевой, см. соответствующий отрывок в 
\cite{LL}.} и восстановлены в 1934г. На 1934-1936 годы 
  пришлось массовое присуждение степеней, в том числе и без защиты диссертаций
  (позже подобные случаи были редкостью). В ведение комиссии входили и звания доцента и профессора,
  то есть Комиссия должна была работать в согласовании со структурами Наркомпроса.

  \sm
  
  $\bullet$ С марта%
  \footnote{В марте Лузин был утвержден на общем собрании АН СССР, \cite{Shpar}} 1935
  - по июль 1936, председатель математической группы АН СССР, т.е., главный математик
  Советского Союза. В это время ему была выделена квартира из 5 комнат, у него был персональный автомобиль с шофером%
  \footnote{Стоит иметь в виду, что в Москве 30х годов было очень тесно, жилищное строительство сильно отставало
  	от промышленного. Стоит также иметь в виду, что это было эгалитаристское общество.}
  (см. ниже п.\ref{ss:serpinski}).
  
  \sm
  
  До начала <<политических преследований>> Лузин был профессором МГУ, а также вице-президентом
  Московского математического общества.... Посмотрим, что было после лета 1936 года.

  \sm
  
  {\bf\punct Лузин, 1936-1950.%
  \label{ss:luzin-after-1936}}
  Следующее утверждение, С.~С.~Кутателадзе \cite{Kut2}:%, \cite{Kut-rus}:
  \begin{quotation}
  	Лузин был подвергнут социальному остракизму и четырнадцать лет до самой смерти жил с клеймом врага в советской маске. Для сталинщины он стал показательным социальным изгоем — враг, а на свободе. 
  \end{quotation}
  
  В массовых цитатах к этому еще часто добавляется <<остался без средств к существованию>>. 
  Это, разумеется, неверно. За ним, по крайней мере, осталась  зарплата советского академика
  (зарплата весьма высокая). Некоторые строчки из предыдущего
  пункта автоматически переходят в этот. 
  
  $\bullet$ Лузин до 1937 продолжал заведовать отделом в Ленинградской Стекловке.
  
  \sm
  
  $\bullet$ Учебник Грэнвиль--Лузин продолжал издаваться, 1937, 1938, 1942. Перерыв в изданиях
   объясняется не зажимом учебника, а наоборот, большим тиражом  1938 года (60 000)
  
  \sm
  
  $\bullet$ В 1936 году Лузин <<под влиянием всех нападок \cite{TyulinaA}>> ушел из МГУ и Московской Стекловки.
  По-видимому, тогда же его увольняют из ЦАГИ.
  
  \sm
  
   $\bullet$ В 1938–1948 гг. — старший научный сотрудник, руководитель теоретического отдела Института автоматики и телемеханики АН СССР, \cite{Koz}.
   Дата 1938  вызывает определенные сомнения, иногда называется 1939 год.
   Однако Бари \cite{TyulinaA} утверждала, что это произошло в том же 1936 году.
   Возможно, что это противоречие связано с какими-то реорганизациями учреждений. 

\sm

$\bullet$ 
 Лузин остался в редакции Мат.сборника (и оставался там до своей смерти).

\sm
  
$\bullet$ В 1940 году издается новая книга Лузина  {\it Теория функций действительной переменной},  318pp, \cite{Luz-real}.
Книга была рекомендована Наркомом просвещения как учебное пособие для педвузов,  переиздание -- 1948.

\sm

$\bullet$ Всегда было общеизвестно, что Лузин вернулся в Стекловку заведовать отделом.
Приводимые разными (причем, хорошо осведомленными) авторами даты расходятся. Бари \cite{TyulinaA},
П.~Л.~Ульянов \cite{Ulyan}, А.~П.~Юшкевич \cite{Yush2} указывают как дату возврата 1941 год. С.~М.~Никольский
утверждает \cite{Niko}, что это произошло в 1943 году, а стекловский справочник
  \cite{Koz} указывает 1946 год. Возможно, что противоречия связаны с перетрясками института во время войны,
  а также с инфарктом Лузина осенью 1941 года \cite{TyulinaA}. Дата окончания работы тоже указывается разной,
  1947 у П.~Л.~Ульянова и <<до смерти>> в стекловском справочнике. По-видимому, права писавшая в 1950г.
  Бари, сообщая, что это случилось в 1947 году после <<резкого столкновения с одним из академиков>> (Колмогоровым).
 
 \sm
 
$\bullet$ 1943-1947 -- профессор Московского университета.
  
  \sm 
  
$\bullet$ С 1944-1950 -- зав. теоретического отдела Сейсмологического института.  

\sm
  
$\bullet$ 1945 Награжден орденом Трудового Красного знамени,  \cite{Obit-1}.

\sm

$\bullet$  1946. Издает новую версию курса математического анализа в виде двух томов,
<<{\it Дифференциальное исчисление}>> и <<{\it Интегральное исчисление}>>. Переиздания 1949, 1952, 1953, 1955, 1958, 1961.

\sm

$\bullet$  1948. Академия наук издает небольшую книжку <<{\it Николай Николаевич Лузин.
Материалы к библиографии ученых СССР}.>>

\sm

 $\bullet$ Умер 28 февраля 1950.

  \sm
  
  Это не вполне похоже  на биографию <<изгоя>> или на биографию <<жертвы режима>>. 

\sm

В 1945г. в статье 	<<Развитие математики в нашей стране>>
(по поводу празднования юбилея АН СССР)
Александров писал:
\begin{quotation}
\dots и если первая его [Лузина] замечательная работа по теории функций действительного
переменного и примыкала непосредственно к уже упомянутой работе Егорова  (1911), то в своих    дальнейших    исследованиях он
сразу же выходит на широкую дорогу вполне самобытного творчества и в несколько же ближайших лет  делается    математиком
крупного масштаба. Однако    не    только    первоклассное научное
дарование самого Н.~Н.~Лузина осуществило переход московской математики на новый этап в смысле ее веса и влияния на мировую науку,— не меньшее    значение    имели и заслуги Лузина по
созданию своей собственной математической школы. В  1915 году
Лузин сплотил вокруг себя первое поколение    своих учеников, к
которому принадлежали, в частности, А.~Я.~Хинчин, Д.~Е.~Меньшов,   автор этой статьи, М.~Я.~Суслин (1894—1919) и немного
позднее — П.  С. Урысон   (1898—1924).	
\end{quotation}

В 1947 году был издан юбилейный библиографический сборник <<Математика в СССР за 30 лет>>
(обзорные статьи+библиографии). Лузин там - один из самых упоминаемых авторов.
В 1946 Б.~В.~Гнеденко опубликовал <<Очерки по истории математики в России>>. Там Лузину  воздается должное. 
В 1948 году Г.~Ф.~Рыбкин и А.~П.~Юшкевич издали первый выпуск <<Историко-ма\-те\-ма\-ти\-че\-ских исследований>>.
Он открывался обширной статьей Александрова, Гнеденко и Степанова <<Математика в Московском университете
в текущем столетии>>. И здесь Лузин -- основатель Московской математической  школы.  

\sm
  
 {\bf \punct Некоторые детали.%
 \label{ss:detail}}
 Представляются существенными детали биографии Лузина и его положения  в ближайшие полтора года 
после атаки июля 1936г.

\sm

1) Известно, что тогда Лузин сразу вышел из Математического отделения академии
и перешел в техническое. Согласно \cite{Erm2},
Колмогоров%
\footnote{Колмогоров в то время не был ни академиком, ни член-кором, но он входил
в Математическую группу АН СССР.} писал ему: «{\it более правильным, если бы Вы в основном работали в 
нашем Отделении, а Вашу деятельность в области механики 
считали дополнительной нагрузкой}», на что Лузин 19 марта 1937г. ответил отказом (цитируем частично):
\begin{quotation}
... По вопросу об основном Отделении моей работы ни состояние
моего здоровья, ни мои лета не позволяют мне в равной степени
отдавать мои силы и время нескольким Отделениям Академии. В
Отделении технических наук я работаю более двух лет и, по-видимому,
вошел в жизнь и интересы этого Отделения. Таким образом,
основным для себя я считаю это Отделение, возможно, что само
Отделение так же смотрит на мою в нем работу.

Но это не означает отклонения мною возможности принести
посильную пользу Физико-математическому отделению. Так, в 
ведущемся мною на Отделении технических наук семинарии по 
системам дифференциальных уравнений с частными производными,
предпринятым мною с целью сдвинуть с места ряд проблем 
дифференциальной геометрии, принимают участие деятели 
Физико-математического отделения.

Примите уверение в моем глубоком уважении

Н.~Лузин»
\end{quotation}

2) Другое известие, из переписки Колмогорова и Александрова \cite{Kolmogorov-sh2}:
\begin{quotation}
Александров  -- Колмогорову, 28 июня 1937 г.

%..... Я усиленно рекомендовал пригласить еще 
%Соболева или Петровского по анализу, Шнирельмана или Гельфонда или
%Чудакова по теории чисел и А.~Я.~Хинчина по теории вероятностей.
%Не знаю, что получится, вероятно, срам. 
Пришлось в моей статье
(где довольно много о Лузинских работах) сделать примечание, что
все эти работы, конечно, не меняют оценку его последующих 
аморальных и антиобщественных поступков. Кстати, по поводу моей
заметки о Лузине в БСЭ я буду говорить с Рыбкиным в том смысле,
что, конечно, дело редакции помещать расширенную мою заметку
без моей подписи или с моей подписью, но что я, во всяком случае, 
совершенно не хочу впечатления, что я боюсь дать свою подпись, тем
более, что уже ходят по разным поводам разговоры, что <<П.~С. 
сейчас, конечно, стремится к хорошим отношениям с Н.~Н., так как Н.~Н
голосовал за П.~С.>>...
\end{quotation}

Кроме текста, стоит посмотреть на дату. Известное роковое решение, по-видимому, было принято в этот день.
Редакторы того знать не могли, и не могли
(как впрочем и кто-либо другой) знать, что будет дальше. Но за окнами уже два месяца творилось черти что
(разборки в верхах и истерия в газетах, в эти дни в Москве заседал пленум ЦК, санкционировавший арест четверти состава ЦК,
об этом, надо думать, вести были).
У редакторов были основания для осторожности.

Что касается подразумеваемой в письме статьи, то
это публикация \cite{Alexandrov-1937} <<Развитие теоретико-множественной математики в СССР>>
в журнале
<<Фронт науки и техники>>, 1937, номер 8-9. В статье много и уважительно говорится про Лузина
и про Егорова, правда там Александров продолжает конфликт относительно истории
аналитических множеств, знакомый современному читателю по Академической стенограмме
1936г. В статье имеется сноска 15а:
\begin{quotation}
Само собой разумеется, что наша характеристика научной работы акад. Лузина
ни в какой степени не меняет нашего резко отрицательного отношения
к тому явлению, которое вскрыто и разоблачено советской научной общественностью
как <<лузинщина>> (раболепие перед иностранной наукой и отдельными ее представителями,
присвоение исследований своих учеников, лицемерие и двурушничество и т.д.).	
\end{quotation}	
	
	В конце статьи стояла приписка 
	\begin{quotation}
		ОТ РЕДАКЦИИ.
В области теории множеств{\bf a} [в журнале так] и непосредственно связанных с ее идеями дисциплин
советской школе математики принадлежит одно из первых мест в мировой науке.

Великая Пролетарская революция, двадцатилетие которой мы празднуем теперь, приобщила к творческой научной работе
новые кадры молодых ученых, завоевавших мировую известность.
Эти ученые явились не просто продолжателями работ своих учителей или заграничных математических школ, но и творцами
ряда новых понятий и методов, давших возможность не только решить ряд старых, казавшихся почти непреодолимыми проблем, но и поставить новые, сыгравшие большую роль во всем дальнейшем развитии математики. Они воспитывались в СССР, были свободны от довлеющего над многими буржуазными математиками и паразитирующего в особенности на теории множеств идеалистического (механистского или родственного ему) мировоззрения (вроде так называемого эффективизма Бореля). К сожалению, на идеологической стороне вопроса,
имеющей, кстати сказать, обширную литературу, автор статьи,
рисующий яркую картину развития математики в СССР основных идей теории множеств{\bf a} [еще раз] совсем не остановился.
	\end{quotation}

В выходных данных журнала были указаны ответственные редакторы биологи А.~Н.~Бах, Б.~А.~Келлер.
 Очевидно, что приписка (такого рода приписки в журнале встречались)
была сочинена не ими.  То ли не вся редакция была указана в выходных данных
журнала, то ли эти слова принадлежали какому-то рецензенту. Приписка весьма осторожна,
только немногие читатели могли бы понять, что (по мнению автора приписки)
<<эффективизм>> есть в СССР в лице прежде всего  Лузина.

 Статья <<Лузин>> в БСЭ не появилась (кстати, не было и статьи <<Колмогоров>>,
но была короткая статья <<Люстерник>> с упоминанием Шнирельмана). В 39 томе энциклопедии (585-591) вышла  статья
Александрова <<Множеств теория>>, где, в частности, говорилось:
\begin{quotation}
 В то же время определились и основные течения буржуазной философии математики: {\it логистика} (см.),
 формализм, {\it интуиционизм} (см.), а также эффективизм, возглавляемый Борелем, Лузиным.
\end{quotation}

Вставка весьма неприятная, но не следует преувеличивать ее долговременного значения: журналы 30х годов были
весьма неполиткорректны и полны брани, ругани на злобу дня хватало и в Энциклопедии.
Трудно представить, что Александров сделал эту вставку по собственной инициативе (он играл в сложные игры, но в такие не играл, см. ниже п.\ref{ss:aleksandrov-polytics}). Однако, заметим, что это подозрительно близко
к редакционной приписке к статье Александрова.

Ниже мы будем разбирать нападки на Лузина в июле 1936. Интересно, что обличения философского характера там фактически отсутствовали. Очевидно, что теперь в игру включились 
партийные философы от математики, круг которых был невелик%
\footnote{Тут первым приходит в голову Кольман (цитированные сентенции его достойны,
	а ругательство <<эффективизм>> он употреблял), но стоит иметь в виду, что
	 в марте 1937г. по нему отстрелялась
газета <<Правда>> \cite{Maximov-Kolman}. В момент разговора Александрова с Рыбкиным Кольман, если верить ему самому,
был безработным (в  \cite{Kolman-vospominaniya} он утверждал, что стал таковым после Московской партийной
конференции, т.е., в мае 1937г.).
Это его не полностью исключает, в принципе он мог еще оставаться в широкой редакции журнала, если он
в нее вообще входил.}.

В 1938г. вышла целая книга В. Н. Молодшего <<Эффективизм в математике>> \cite{Molodshy-Luzin}
\begin{quotation}
	Эффективизм Бореля нашел убежденных сторонников даже среди некоторых математиков Советского Союза, в числе которых находится акад. Н.~Лузин.
\end{quotation}

Это из предисловия, вообще в книге много ругани в адрес Лузина. Кроме того, 
об эффективистах там говорится во множественном числе, возможно, что имелись в виду
Новиков и Л.~Келдыш.

Интересно, что Лузин (который вообще на рожон не лез) продолжал в это время 
(1938-1939гг.)
философствовать (см. ниже сноску \ref{fo:uchebnik}).

\sm

3) Можно посмотреть на то, как упоминался Лузин в математической печати того времени.
В единственном номере Успехов за 1937 публикуется статья <<Изжить лузинщину в научной среде>>
\cite{Uspehi-37-1}
(что было необходимо по правилам игры). В следующем выпуске \cite{Uspehi-38} (уже за 1938г.) выходит редакционная статья 
<<Советская математика за  20 лет>>. Лузин там упоминается 4 раза, П.~С.~Александров - 4, Люстерник -3, Гельфонд -5,
Колмогоров - 9. Всё нормально.

В 1938г. выходит сборник <<{\it Математика и естествознание в СССР за 20 лет}>>. Там есть статья Л.~В.~Канторовича и Г.~М.~Фихтенгольца
 {\it Успехи теории функций вещественной переменной 
  	и функционального анализа,} \cite{Fiht-Kantor}. Лузину там воздается должное (ровно там, где нужно  воздавать).
  
  \sm

4) В описи архива Лузина есть бумага с названием:
{\it 	
Докладная записка ак. Н.~Н.~Лузина, чл.-корр. Н.~И.~Мусхелишвили, профессоров М.~А.~Лаврентьева и Б.~И.~Сегала 
(сотрудники Математического института АН СССР) в ЦК ВКП(б) о реорганизации высшего образования от {\bf 20.07.1936г.}}

\sm

5) Представляет определенный интерес время
работы Лузина над учебником, изданным в 1940г. В описи \cite{arhiv} есть  заявления с просьбой предоставить отпуск с этой целью
от
2.09.1938 --	1.08.1939.

В описи архива есть
письма Н.~Н.~Лузина Вейфелю Ною Ильичу ([редактор Учпедгиза]),  20.10.1937г.
(относится  ли это к работе над учебником, по описи понять нельзя).

\sm

По-видимому, призрак газетной кампании июля 1936г. ощущался в течение года-полутора, постепенно рассеваясь
(и, видимо, будучи оживляемым околоматематическими философами).  
Стоит иметь в виду, что написанное в газетах не принято слишком долго помнить,
так обстоит сейчас, так было и тогда. К тому же за эти полтора года произошло столько событий...  

\sm

{\bf\punct Защитники Лузина.%
\label{ss:za-list}} Перед началом дальнейшего обсуждения познакомимся
 с участниками защиты и нападения.

 Кстати, защитники, по-видимому, противоречили воле Партии.
Заодно посмотрим, как Сталин покарал их. Приведем  обрывки из биографий этих людей
за 1936-1953гг (до смерти Сталина), показывающие их социальное положение.
Начнем с тех, кто сопротивлялся открыто.

  \sm
  
  БЕРНШТЕЙН Сергей Натанович (1880-1968) -- Сталинская премия (1942),
  Награждён 2 орденами Ленина (1945, 1953), орденом Трудового Красного Знамени (1944).

  \sm
  
  КРЫЛОВ Алексей Николаевич%
  \footnote{Прикладной математик и инженер-кораблестроитель. Книга «Теория качки корабля».
  В 1898 году награждён золотой медалью Британского общества корабельных инженеров.
  В 1908—1910 годах — главный инспектор кораблестроения. Вместе с И.~Г.~Бубновым проектировал
  линкоры типа <<Севастополь>>. И т.д и т.п.} (1863-1945) -- Сталинская премия, 1941.
  Награждён 3 орденами Ленина (1939, 1943, 1945). Герой Социалистического Труда (1943).
  В 1951-1956гг. было издано собрание сочинений Крылова в 12 больших томах
  (фактически в 17).

\sm
   
   МЕНЬШОВ Дмитрий Евгеньевич -- Сталинская премия, 1951. Заведующий кафедрой Теории функций и функционального анализа Мехмата
   МГУ, 1941-1977.

 \sm
 
 БАРИ Нина Карловна (1901-1961) -- профессор МГУ.
 
 \sm
 
 Отказались играть предложенную им роль:
 
 \sm
 
 НОВИКОВ Петр Сергеевич (1901-1975). С 1945г. заведовал кафедрой в Московском педе.
 
 \sm 
 
 ЛАВРЕНТЬЕВ Михаил Алексеевич (1900-1980). В 1939—1941гг. и 1945—1948гг. — директор Института математики 
 Академии наук УССР.
 Сталинские премии 1946, 1949. Орден Отечественной войны II степени — 1944.
Орден Трудового Красного Знамени — 1945, 1948, 1953; Орден Ленина — 1953.
 
 \sm
 
 Против был и
 
 ЧАПЛЫГИН Сергей Александрович (1869-1942), начальник Теоретического отдела ЦАГИ. Герой социалистического труда, 1941.
 В 1942 Академия наук СССР учредила премию им.  Чаплыгина. В 1948 году город Ранненбург был переименован в Чаплыгин.
 
  \sm
 
 По-видимому, на стороне Лузина были Голубев  (Лузин несколько раз просил пригласить его на заседания
 академической комиссии) и Л.~Келдыш, но в опубликованных документах нет ничего
 об их участии%
 \footnote{О том, что Л. Келдыш была на стороне Лузина, есть позднее известие А. В. Гладкого \cite{Gla}.}. Келдыш в то время была без чинов и наград, а у Голубева и чинов и наград было много.
 
  \sm
  
  Так или иначе, репрессии были не очень суровыми.
  
  \sm
 
{\bf \punct Судьба нападавших.%
\label{ss:protiv-list}} Посмотрим для контроля на тех, кто участвовал в облаве на Лузина.
Сначала лузитане (непосредственные ученики Лузина).

\sm

  АЛЕКСАНДРОВ Павел Сергеевич (1896-1982), Сталинская премия, 1943.
  
\sm
  
   КОЛМОГОРОВ Андрей Николаевич (1903-1987). Председатель математической группы АН СССР с   1939. Сталинская премия 1941.
   Орден Трудового Красного знамени, 1940.
  Орден Ленина 1944, 1945.

\sm

    ЛЮСТЕРНИК Лазарь Аронович (1899-1981). Сталинская премия, 1946.
   
\sm   
    
    ЛЯПУНОВ%
    \footnote{Единственное упоминание Ляпунова в числе нападавших - у Бари \cite{TyulinaA},
    сомневаться в ее показаниях нет оснований (и есть некоторые косвенные доводы
    в пользу этого, см. п.\ref{ss:bari}). Кстати, он - родственник
    знаменитого математика Ляпунова Александра Михайловича
    и (говорят) племянник в двадцатом колене 
    полевого командира Ляпунова Прокофия Петровича. Работал в дескриптивной теории множеств, с 50х годов
    в computer science.}
    Алексей Андреевич. Уволен из Стекловки в 1937 году по сокращению штатов, восстановлен в 1939г. 
    Участник Войны (командир топографического взвода в артиллерии). Орден Красной Звезды, 1944.
    В 1945-49 -- преподаватель Артиллерийской академии. В 1946-49 - докторант Стекловки. С 1951г. - сотрудник Стекловки.

  \sm
  
  ХИНЧИН Александр Яковлевич (1894—1959), Сталинская премия, 1941. Орден Трудового Красного знамени 1944, 1945. Орден Ленина 1953.
  
  \sm
  
  ШНИРЕЛЬМАН Лев Генрихович. Самый ярый из нападавших. Покончил с собой (отравившись газом) в сентябре 1938 года. Есть две версии, обычное самоубийство, или
  результат вербовочной беседы в НКВД. См. ниже п.\ref{schnirelman-gibel}
  
  \sm
  
Теперь прочие. Приводим список тех, о чьих выступлениях есть какая-либо информация. 

  \sm
  
    АРШОН Соломон Ефимович%
  \footnote{Автор статей в Мат. Сборнике,  42:1, 44:4.  Директор издательства ОНТИ.} 
   (1892-?).
  Арестован, по-видимому, летом 1938 года. Погиб.
  
  \sm
  
    БУХГОЛЬЦ Николай Николаевич%
  \footnote{Механик,
  	профессор Московского Физмата, потом Мехмата, автор 
  	 учебника <<Основной курс теоретической механики>>, он многократно издавался в 1930-1945гг.
  	и 1965-1972гг.} (1981-1943).
  Заведующий кафедрой теории упругости (1932–1938); заведующий кафедрой теоретической механики (1938–1942) Мехмата. 
  Лауреат Сталинской премии, 1943.
  Ордена Ленина (1936), Трудового Красного Знамени (1940), «Знак Почёта» (1940).
  Генерал-майор авиационно-технической службы (1943).
  
  \sm

    ГЕЛЬФОНД Александр Осипович (1906-1968). Ордена Трудового Красного знамени,  1945, 1945. Орден Ленина, 1953.
  Член-корреспондент АН СССР, 1953
  
    \sm
  
  ПОНТРЯГИН Лев Семенович (1908-1988). Сталинская премия 1942. Орден Трудового красного знамени, 1945.
  Орден Ленина,  1953.
  
  \sm
  
    СЕГАЛ Бенцион Израилевич (1901-1971).
  Секретарь парткома Стекловки, член Президиума Математической группы АН СССР,
  1936--1938(?)%
  \footnote{Сегал - теоретико-числовик. Один из руководителей кампании против Н.~М.~Гюнтера в Ленинграде,
  	см. п.\ref{ss:lenmatfront} в 1929-30гг. (в бытность его аспирантом), позднее переехал в Москву.
  	Упомянутый Президиум Математической группы состоял из Бернштейна, Виноградова и Сегала,
  очевидно, что Сегал был там политическим представителем. На Втором Всесоюзном математическом съезде
  1934г.
Был выбран в Президиум Совета Всесоюзной математической ассоциации, Президиум состоял из 5 человек,
Сегал был секретарем (фактически, эта ассоциация, видимо, не сыграла предполагавшейся в тот момент роли). Естественно предполагать в Сегале главного политического представителя
ВКП(б) в математических науках в 1934-38(?)гг . Участвовал в заседанииях академической комиссии по делу Лузина
1936г., был немногословен, но весьма опасен. К сожалению, каких-либо биографических источников об
этом деятеле мне найти не
удалось.}. В 1935--1971  зав. кафедрой математики в СТАНКИНе (Московский станкоинструментальный институт).
  
  \sm

 СОБОЛЕВ Сергей Львович (1908-1989), Сталинские премии  1941, 1951, 1953. Герой социалистического труда 1951.
 Ордена Ленина  1945, 1949, 1951, 1953.
%  Sergey SOBOLEV (1908-1989). Stalin prizes, 1941, 1951, 1953. Hero of Socialist Labour, 1951.
 % Orders of Lenin 1945, 1949, 1951, 1953.

 % Lev PONTRYAGIN (1908-1988). Stalin prize, 1942.
  %Order of the Red Banner of Labour, 1945. Order of Lenin, 1953.
  
  \sm
  
  ХВОРОСТИН Гавриил Кириллович (1900-1938)%
 % Gavriil HVOROSTIN(1900-1938)%
  \footnote{Предводитель радикальных студентов на Физмате МГУ. в 1924-29гг. Участник свержения Егорова в 1929-30гг.
  	О нем ниже в пп.  \ref{ss:hvorostin-1}, \ref{ss:hvorostin-2}.
  Единственное известие о его участии -- у С.~П.~Новикова-мл., см. ниже п.\ref{ss:novikov-hvorostin}}, (ди)ректор Саратовского университета 1935-1937.
  Арестован 02.08.1937. Расстрелян.
 % 	Hvorostin Gavriil Kirillovich. A leader of radical communistic students of PhysMath of Moscow University
 % in 1924-1929. A participant of the fall of Egorov 1929-1930. 
 % In local histories of Saratov University, he is regarded as a good rector who tried to build 'G\"ottingen on Volga'. See below ????. }. 
 % The (di)rector of the Saratov University, 1935-1937. In 1937, 2 August, he was arrested.

 % Bentsion SEGAL%
 % \footnote{
 % 	A representative of radical communists; one of leaders of a company againist
 % Nickolay G\"unter in Leningrad, 1928-1930. However, he was a respectable number-theorist in 1930s.}(1901-1971). 
%Секретарь парткома
 % The Secretary of Communist Party Group of the Steklov Institute Ученый секретарь.
  %The head of the chair of Mathematics in Moscow Machine-Instrument Institute[СТАНКИН].

  %Лауреат Государственной премии СССР (1943).
  %Nikolay BUCHHOLZ (1981-1943)%
  %\footnote{The author of a long-living text-book on theoretical mechanics. The head of Chairs of
  %Elastisity 1933-1938 and Theoretical mechanics, 1938-1943. A non-hiding  Christian. } (1881-1943). Stalin prize, 1943.
  %A major-general of technical aviation service, 1943.
  % Order of the Red Banner of Labour, 1940.
   
   \sm
   
   ШМИДТ Отто Юльевич (1891-1956). Член ВЦИК 1935--1937.
    Начальник ГлавСевМорПути, 1932-1938, Герой Советского Союза 1937,
   Ордена Ленина 1932, 1937, 1953,
   Ордена Трудового Красного Знамени 1937, 1945, в
   1939--1942  вице-президент АН СССР.
   
   \sm

   Здесь мы тоже видим блестящих людей с успешными карьерами, 
   однако фортуна отнюдь не ко всем была благосклонна,
   и вероятность сильной неблагосклонности была весьма велика.

   \sm
   
   В день, когда была опубликована статья в Правде в Стекловке состоялось собрание, заклеймившее Лузина.
   Подробности нам не известны \cite{Front1936}: выступали
    И.~М.~Виноградов, чл.-кор. АН СССР
С.~Л.~Соболев, чл.-кор. АН СССР Л.~Г.~Шнирельман, А.~О.~Гельфонд, 
Б.~И.~Сегал, Л.~А.~Люстерник, Ф.~Р.~Гантмахер, С.~Е.~Аршон, Н.~Е.~Кочин, А.~Ф.~Бермант, М.~В.~Келдыш.

Большая часть лиц уже была названа. Виноградов входил в академическую комиссию по делу Лузина, и там молчал.
Остаются Гантмахер, Кочин, Бермант, М.~В.~Келдыш. 
Гантмахер имел на Лузина вполне конкретный зуб, 
о нем упоминается ниже (но, что он говорил на собрании, неизвестно). Бермант едва ли ругался меньше остальных.
Про Кочина и М.~В.~Келдыша мы ничего не знаем.

Из нематематиков на заседании академической комиссии самым агрессивным был академик Горбунов Николай Петрович (1892-1938),
непременный секретарь АН СССР. Он был расстрелян в 1938 году.
Остальные члены комиссии нематематики (А.~Н.~Бах, А.~Е.~Ферсман, Г.~М.~Кржижановский) завершили жизненный путь при чинах и почестях.

Ладно, мы познакомились с лицами, и это знакомство важно для дальнейшего обсуждения. Что касается статистики смертности, то она ничего не доказывает, 
   но уж точно мы не  подтверждает того,
   что выступавшие за Лузина выступали против воли Вождя, а против Лузина - в соответствии с оной. 

   Перейдем к более содержательным наблюдениям.

   \sm

{\bf \punct Линия раздела в 1936 году.%
\label{ss:razdel-1936}}
Давайте посмотрим, какой математикой занимался каждый из лузитан в первой половине 30х годов. Сначала нападающая сторона

\sm

Александров --- топология.

Колмогоров --- теория вероятностей, логика.

Люстерник -- вариационное исчисление, функциональный анализ, дифференциальные уравнения.

Ляпунов --- дескриптивная теория множеств.

Хинчин -- теория вероятностей, теория чисел.
 
Шнирельман -- теория чисел, вариационное исчисление. 

\sm

Сторона Лузина:

\sm

Бари --- ТФДП%
\footnote{Теория функций действительной переменной. Сейчас чаще говорят <<вещественный анализ>> (real analysis),
мы используем принятый тогда термин, а также его сокращение, использовавшееся в устной речи (тээфдэпэ).}.

Келдыш -- дескриптивная теория множеств.

Лаврентьев -- гидродинамика.

Новиков -- дескриптивная теория множеств.

Меньшов -- ТФДП.

\sm

Выходит, что из пятерых представителей про-лузинской стороны четверо (кроме Лаврентьева) работали в лузинской тематике, а из семерых
анти-лу\-зин\-цев шестеро (кроме Ляпунова) тематику радикально сменили. Что заставляет задуматься над тем, а не 
стояла ли за конфликтом научная составляющая.

\sm

{\bf\punct От смерти Лузина до 1953 года.%
\label{ss:1950-1953}}
a) Я нашел четыре некролога Лузина. Один был опубликован во второй по значению советской газете <<{\it Известия}>>, \cite{Obit-2}
и был подписан почти всеми институтами, где Лузин работал (кроме ЦАГИ). Есть постановление Президиума РАН, которое поручает 
написание некролога Ляпунову, Л.~Келдыш и Новикову%
\footnote{
http://odasib.ru/OpenArchive/Portrait.cshtml?id=Xu$_-$pavl$_-$635188979823906250$_-$3056}.
Известен рукописный текст некролога, написанный рукой Ляпунова, см. \cite{Lyap}.

Второй некролог был опубликован в <<{\it Успехах математических наук}>> \cite{Obit-1} без подписей, что для этого журнала не характерно, но все же случалось%
\footnote{Например, без подписей был опубликован некролог Шнирельмана в 1939г. (но это случай совсем особый), некролог Стефана Банаха (Stephan Banach) в 1946г., без подписей
(и без некролога) было извещение о смерти Александрова в 1982г.}.
Кроме того был некролог в <<{\it Математике в школе}>> \cite{Luzin-mat-v-shkole} (тоже без подписей) и в журнале <<{\it Автоматика и телемеханика}>> \cite{Kuleb},
подписанный  В.~С.~Кулебакиным%
\footnote{Виктор Сергеевич Кулебакин (1891— 1970) —   академик АН СССР ( c1939),
генерал-майор инженерно-авиационной службы. Основатель и первый директор (1939-1941) Института автоматики и телемеханики
(ИПУ АН СССР).}.

\sm

b) В 5-ом номере <<{\it Успехов>>} за 1950 год был опубликован цикл статей по дескриптивной теории множеств (Е.~А.~Щегольков, Ляпунов, М.~Я.~Арсенин), 
то есть по лузинской тематике.
Скорее всего, основные статьи статьи готовились до смерти Лузина (они  следующем номере после номера с некрологом;
они
хорошо написаны, и позже даже вышли отдельной книгой по-немецки%
\footnote{Ljapunow, A. A.; Stschegolkow, E. A.; Arsenin, W. J. {\it Arbeiten zur deskriptiven Mengenlehre.}  Berlin, 1955. iii+108 pp.}; 
в номере никак не говорится, что он мемориальный).
Небольшая вводная статья
Ляпунова \cite{Lyap-vvedenie} была прямой апологией Лузина.

\sm

c) В 1951-1953гг в Успехах был опубликован цикл из 9 статей, посвященных памяти Лузина. Он включал
три посмертных публикации Лузина \cite{Luzin-add1}--\cite{Luzin-add3}, отзыв Егорова на диссертацию Лузина
\cite{Egor-lus}, а также статьи, посвященные обзору работ Лузина, а именно,  Бари, Люстерник \cite{BaLyu1},
Келдыш, Новиков \cite{KeldNov}, Гольцман%
\footnote{Гольцман В.~К., автор статьи,  Матем. сб., 38:3-4 (1931). Других сведений мне найти не удалось.},
Кузнецов%
\footnote{Кузнецов Петр Иванович (1911--?). По-видимому, ученик Привалова. Соавтор Привалова и Лузина. Работы по специальным функциям
и уравнениям с частными производными. Автор/редактор нескольких публикаций о Лузине \cite{Kuz1965}, \cite{Kuz-Luz}, \cite{Kuz2}, а также о Егорове \cite{KuzEgor}.}
\cite{GolKuz}, Сретенский%
\footnote{Сретенский Леонид Николаевич (1902-1973), ученик Егорова и Чаплыгина,
основные работы по гидродинамике, чл.-корреспондент АН СССР с 1939,
проф. Мехмата МГУ.} 
\cite{L-Sret}, Федоров%
\footnote{Фёдоров  Владимир Семёнович (1893—1983), ученик Лузина. Во время Гражданской войны переехал в Иваново-Вознесенск вслед за
Лузиным, и остался там работать. Основные работы по теории функций комплексного переменного.}
\cite{Fedor}.

\sm

d) В 1951 году была издана в виде отдельной книги диссертация Лузина <<{\it Интеграл и тригонометрический ряд}>>.
Она предварялась биографией Лузина, составленной Голубевым и Бари \cite{GoBa}, вводной статьей Бари и Люстерника \cite{BaLyu2}
и сопровождалась тщательными комментариями Бари и Меньшова.
 
 В 1953г. Новиков, Келдыш и Бари издали русский перевод \cite{Luz-anal-rus} книги Лузина  
 <<{\it Le\v{c}ons sur les ensembles analytiques et leurs applications}>>, книга содержит важное предисловие.
 
 \sm
 
 e) В 1953, 1958, 1959 годах были изданы три тома собрания сочинений Лузина \cite{Luz-collected-1}--\cite{Luz-collected-1}.
 В последний том были включены включены несколько упомянутых выше мемориальных статей и библиография работ Лузина.
 Во втором томе, посвященном дескриптивной теории множеств (744с.), также прилагалась вступительная статья и комментарии
 Новикова и Келдыш. Кстати, почти весь том был переведен с французского.
 
 Комиссия по изданию трудов Лузина состояла из Лаврентьева (председатель), Некрасова%
 \footnote{Некрасов Александр Иванович (1883--1957), специалист по гидро-аэро-динамике, занимался флаттером,
 академик с 1946 года, Сталинская премия 1952 года. В 1937 во время арестов в ЦАГИ был арестован (кажется, единственный
 из известных людей, связанных с московским Мехматом), В 1938, после прихода Берии в НКВД (так же как и другие
 осужденные деятели ЦАГИ) переведен в закрытое Туполевское КБ, освобожден в 1943.}, Меньшова, Новикова, Сретенского, Бари и Келдыш.

\sm

Так или иначе, посмертного остракизма Лузина при Сталине (как, впрочем, и прижизненного) не видно.

Если смотреть на линию раздела, то со стороны Лузина мы по-прежнему видим тех же Бари, Меньшова, Новикова, Келдыш, Голубева,
к которым из нападавших зубров присоединяется Люстерник. Кроме того, на стороне Лузина оказывается Ляпунов, но в больших издательских проектах он,
почему-то, не участвует%
\footnote{Со стороны Ляпунова была также рецензия \cite{Lyap-retsenziya} с критическими комментариями в отношении издания первого
тома Собрания сочинений Лузина.}. На Лузинской стороне находятся также несколько известных прикладников, в основном, связанных с ЦАГИ.

\sm

{\bf \punct  Статьи о Лузине, 1960-1985.%
\label{ss:1960-1985}} Их было довольно много. Читать их небезынтересно, а 
смотреть на список, наверно, скучновато

1960. \cite{??}

1961. Бари \cite{Bari-lyudi-nauki}

1963. Минковский%
\footnote{Минковский Владимир Львович  (1911-1978), известный методист.}
\cite{Mink}

1965. П.~И.~Кузнецов \cite{Kuz1965}.

1974. Лаврентьев \cite{Lavr1}, Келдыш \cite{Keld-Luzin}, П.~И.~Кузнецов \cite{Kuz-Luz}.

1977. Александров \cite{Alex-kvant}, \cite{Alex-v-shkole} (это две версии одной статьи)

1983. Меньшов \cite{Menshov-L}, \cite{Menshov-L-2}, Колмогоров \cite{Kolmogor}, Е.~А.~Щегольков \cite{Shchegol},  Бескин \cite{Besk},
П.~Л.~Ульянов, В.~М.~Беляков \cite{Ulyanov-Kemerovo},

1984. П.~И.~Кузнецов \cite{Kuz2} (сборник статей), В.~А.~Успенский, В.~Г.~Кановей \cite{UspKan}, \cite{Vestn3} (видимо, Успенский),
компиляция в Кванте \cite{Kvant} и полустраничная заметка
П.~Л.~Ульянова \cite{Ulyanov-mshk}.

1985. П.~Л.~Ульянов \cite{Ulyan}, \cite{Ulya2}, Кановей \cite{Kanov}, Успенский \cite{Usp}, Меньшов, В.~А.~Скворцов, П.Л.Ульянов \cite{MSU}.

\sm
 
 Апологии 1950-1985 весьма разнородны -- некрологи, обзоры работ, биографии, воспоминания, гибридные тексты.
 Но их много (почти четыре десятка). В связи с эти {\bf можно задать не только вопрос о том, 
 что в них есть, но и вопрос о том, чего в них нет}.
 
 \sm
 
 {\bf\punct Чего нет в статьях о Лузине 1950-85гг. 1.%
 \label{ss:net-1}}
 Во-первых, нет ни одного намека на то, что Лузин был объектом атаки сверху. Хорошо, при Сталине об этом
 едва ли можно было  говорить. Но дальше наша родная Партия осудила <<усиливавшиеся с декабря 1934 года политические 
 репрессии>>%
 \footnote{Я не делаю высказывания о репрессиях от себя, а только напоминаю, что  Партия осудила.}.
 Если бы даже открытое признание в печати было бы невозможным, то <<политическое дело Лузина>>,
 если бы его помнили  как политическое, породило бы намеки%
 \footnote{Официальное отношение к репрессиям  конца 30х годов было отрицательным.
 Для Хрущева это было средство легитимизации собственной власти и
 борьбы с политическими противниками. При Брежневе это отношение, видимо, было еше более отрицательным, но об этом меньше говорили.
 В любом случае намеки проблемы  не составляли.}, легальные, полулегальные или нелегальные воспоминания,
 о том, что кто-то Лузина защищал, кого-то что-то заставляли и т.д. и.т.п. Нет ничего менее достоверного, чем воспоминания такого типа,
 но важно то, что никаких таких воспоминаний не было и в помине. Вплоть до появления статей А.~П.~Юшкевича 1989-91 года \cite{Yush1}, \cite{Yush2}.
 
 Здесь есть возможности для <<контрольных опытов>>. 
 
 \sm
 
 A) В биологии была  известная сессия ВАСХНИЛ августа 1948 года. Ответная реакция научного мира была довольно быстрой.
 Уже через 2-3 года И.~И.~Презент (видимо, в результате интриг научного сообщества и его собственных выдающихся <<достоинств>>)
  потерял свой деканский пост в Московском университете и вынужден был покинуть Ленинградский университет. Вскоре после смерти Сталина
 последовало открытое выступление ученых против Т.~Д.~Лысенко. А сколько потом было фольклора, неопубликованных рукописей,
 да и опубликованных материалов тоже!
 
 \sm
 
 B) Как стало общеизвестно в Перестройку, Президент Московского математического общества 
 и директор НИИ математики и механики при 1ом МГУ
 Егоров был арестован в октябре  1930года и умер  от язвы желудка  в сентябре 1931 года  после освобождения. 
 Егоров до Войны рассматривался как реакционер, но хороший математик
 (см. Академическую Стенограмму, где математика смешана с политикой); весьма уважительно говорилось о Егорове, например, в
 статьях Александрова и Колмогорова  \cite{Alexandrov-1937} и \cite{Kolmogorov-matematika} 1937-38гг.),
 а после войны --
 как важная положительная фигура и в организационном плане тоже. О нем много упоминали, но текстов, посвященных лично ему, почти не было. 
 Можно посмотреть, как отразились события 1930-1931 года  в советской математической печати.
 
 Процитируем единственную статью в Успехах, посвященную Егорову \cite{KuzEgor}, П.~И.~Кузнецов, 1971:
 \begin{quotation}
 Д.~Ф.~Егоров скончался 10 сентября 1931 г. в Казани («{\it Известия}» от
 25.IX. 1931 г.) и похоронен там же на Арском кладбище.
 \end{quotation} 
 
 Странная фраза без дальнейших пояснений, призванная оставить недоумение у читателя%
 \footnote{Я, кстати, не нашел сообщения в этом номере <<{\it Известий}>>, скорее всего, недостаточно внимательно искал.}.
 
 Цитируем  воспоминания Александрова \cite{Alex-auto1},
 
 \begin{quotation}
  \dots все двадцатые
годы Д.~Ф.~Егоров руководил и Институтом математики и механики и Математическим
обществом, и поэтому стоял во главе всей московской математики.
Это время, несомненно, было не только одним из самых продуктивных,
но и одним из самых светлых периодов жизни математической школы Московского
университета.

В начале осени 1931 г. жизнь Д.~Ф.~Егорова тяжело и горько закончилась
в одной из клиник Казанского медицинского факультета, куда он был
переведен лишь за несколько дней до своей смерти. Могила Д.~Ф.~Егорова
находится на казанском кладбище рядом с могилой Лобачевского. Заботами
казанских математиков могила Егорова долгое время поддерживалась в полном
порядке (как и могила Лобачевского). Я надеюсь, что это продолжается
и сейчас.
 \end{quotation}

Что касается вообще событий 1930-31гг в Московской математике (в этой статье речь  о них пойдет ниже в \S \ref{s:annals}),
то в математической печати, начиная с 1942 года, многократно подчеркивалось, что происходило что-то, о чем не
очень удобно говорить.
Была история до 1930 года,
 а потом снова с 1932 года. А в середине была какая-то неведомая диалектика, см. Александров \cite{Alex-1942}, \cite{Alex-MMO1},
 Александров, О.~Н.~Головин \cite{Alex-MMO3}, А.~Г.~Курош \cite{Kurosh}. 
 
 Так или иначе, намеки при желании делались. 
 В рассказах о Лузине намеки тоже были. Но не о политике, а о конфликтах Лузина с учениками.... 
 
 \sm
 
 В общем, {\bf в публикациях московского математического мира
 	 идея политической атаки на Лузина сверху до 1989 года не наблюдалась.}
 
 \sm
 
 Следующее замечание менее очевидно.
 
 \sm
 
 {\bf\punct Чего нет в статьях о Лузине. 2.%
 \label{ss:net-2}}
 В тщательно изданном <<{\it Собрании сочинений}>> Лузина не было раздела <<Основные даты жизни и деятельности>>,
 что несколько необычно, но в принципе допустимо.  Нет этого раздела и в сборнике статей о Лузине,
 изданном П.~И.~Кузнецовым в 
\cite{Kuz2} 1984г..
 Открывавшая в 1985году номер <<{\it Успехов}>>, посвященный Лузину, статья П.~Л.~Ульянова  \cite{Ulyan}
 фактически состояла в Curicullum Vitae Лузина. Но этот CV не полон (и, возможно, не точен). Но, что более удивительно,
 {\bf все опубликованные до 1985 года статьи о Лузине
 вместе взятые, а также (насколько я могу судить)
 весь справочный материал в околоматематической литературе, не дают возможности составить его простейшего 
 послужного списка!}
 
 Биография Лузина идет как обычно до 1930 года, а потом расплывается.
 Тон был задан в 1951 году биографической статьей Голубева и Бари  \cite{GoBa}:
 
 \begin{quotation}
 В 30 и 40-х годах, кроме Института им. В.~А.~Стеклова, Н.~Н.~Лузин работал и в других институтах Академии наук:
 в Сейсмологическом и в Институте автоматики и телемеханики. 
 \end{quotation}
  (когда??)
 \begin{quotation}
  В эти годы работы в институтах Академии наук Н.~Н.~Лузин уже не был связан с университетом систематически.
  Однако иногда он возобновлял там работу, и это неизменно оказывало влияние на молодых математиков,... 
 \end{quotation}
(из фразы можно понять, что возобновлял не один раз, но когда?)
\begin{quotation}
 Связь с Институтом [Стеклова] стала более прочной с 1934 г., когда Академия наук и ее
Математический институт были переведены в Москву. Н.~Н.~Лузин продолжал руководство отделом теории функций до конца жизни;
все сотрудники этого отдела являются его ближайшими учениками.
\end{quotation}
(да-да, до конца с жизни, но с перерывом, причем об этом перерыве 40 лет спустя все студенты Мехмата знали;
кстати, согласно рыбе к этой статье, написанной Бари \cite{TyulinaA},  не до конца жизни).

\sm

Все бесчисленные математические  статьи о Лузине следовали этому стилю....

\sm

Но картина еще более забавна. {\bf Ни один из перечисленных выше текстов} (которые, кстати, не все написаны про-лузинцами) 
{\bf не содержит ни упоминания, ни намека на 
формально самую высокую социальную позицию Лузина} -- а именно его должность главного математика Советского Союза в 1935-1936 годах.
 Объяснение событий 1936 года должно одновременно объяснять эту странность...
 
 \sm
 
 Но я чуть-чуть погрешил против истины. Дело в том, что была еще серия статей Люстерника \cite{Lyu-1}--\cite{Lyu-4}
 <<Молодость Московской математической школы>> (1967-1970), объемом 113 стр. И в ней между делом говорится:
 \begin{quotation}
  Правда, как это иногда бывает, по мере сужения его [Лузина] роли в реальном
руководстве увеличивалась его роль в официальном руководстве математикой
и ее представительстве: в период Лузитании он не входил в узкий
президиум Математического общества, в период распада Лузитании он был
вице-президентом Общества, в первой половине 30-х годов, в период своей
возросшей изоляции, он был председателем математической группы Академии
наук. Но одно для него не могло заменить другого. В самом деле, все помнят
Лузитанию, но кто помнит, какие посты занимал Лузин? 
 \end{quotation}
 
 Разумеется, Люстерник лукавил. Во-первых, 35-36 годы, это не <<первая половина 30х>>. 
 И он-то очень хорошо помнил, <<какие посты занимал Лузин>>. И помнили это все математики, которым в 1936 году было за двадцать,
 а таких 30 лет спустя было еще много. Собственно, Люстерник, предлагает о Лузинском царстве 1935-36гг не вспоминать....
 И вся про-лузинская литература, начиная с 1951 года, следовала этой идее... 
 
 \sm
 
 {\bf\punct Беседы с Меньшовым, 1975.%
 \label{ss:menshov-1975}} По инициативе Петровского в МГУ с 1966г. записывались фонограммы воспоминаний
 известных ученых
 о первой трети XX века. Работали над этим В.~Д.~Дувакин, а потом В.~Ф.~Трейдер. Опубликованы они были много позже, в частности, записи
 Александрова, Б.~В.~Гнеденко, Б.~Н.~Делоне, Люстерника, Меньшова, О.~А.~Олейник, Л.~И.~Седова, А.~П.~Юшкевича появились в 2005г. \cite{Duvakin}. 
 
 Меньшов, конечно, рассказывал 
 о своем знакомстве с Лузиным, о сотоварищах по Лузитании...
 \begin{quotation}
  В.~Д.[Виктор Дувакин]: Да. Лузина уж давайте сегодня начинать не будем, это большая тема, а вот мы так сейчас...
  
Д.~М.: Да. О Лузине я бы так сказал — я бы стал рассказывать о Лузине с периода с 14-го года по 30-й год....

[следующая беседа]

В.~Д.: Может, и о 30-х годах расскажете?

Д.~М.: Нет, это мне не хотелось бы. Лузин тут отошел, у него были неприятности, и с математиками некоторыми. 
Это не для записи я говорю. Вы это записываете?

В.~Д.: Сейчас я ваши слова записал, да.

Д.~М.: Нет, вот последние?

В.~Д.: Да.

Д.~М.: Лучше сотрите их. (Перерыв в записи.)
 \end{quotation}

 \sm
 
{\bf\punct Последние очевидцы.%
\label{ss:ochevidtsy}}  Вот даты смерти последних участников событий:

Ляпунов  1973

Новиков     1974

Л.~Келдыш -- 1976

М.~Келдыш - 1978

Лаврентьев - 1980

Люстерник -  1981

Александров -1982

Виноградов - 1983

Колмогоров  - 1987

Меньшов -    1988

Понтрягин -- 1988

Соболев 03.01.1989

На момент смерти Соболева в математической литературе%
\footnote{В нематематической литературе было упоминание в { \it Письмах В.И. Вернадского А.Е. Ферсману},
	опубликованных в 1985г.} о конфликте 1936 года был опубликован лишь один абзац в воспоминаниях Понтрягина \cite{Pon0}
 в 1978 году. Звучал он так:
\begin{quotation}
 Но 1936 г. отмечен для меня и другим важным событием.
К этому времени трения в отношениях с моим учителем П.~С.~Александровым
привели к тому, что я открыто взбунтовался против него. В то время
Н.~Н.~Лузин подвергся резкой критике со стороны математической общественности.
И мой бунт против П.~С.~Александрова выразился в том, что,
выступая на обширном собрании математиков против Н.~Н.~Лузина,
я в довольно резкой форме указал на некоторые недостатки в действиях
П.~С.~Александрова. После моего выступления Павел Сергеевич подсел ко мне
и поблагодарил меня за правильную критику в его адрес....

Два мои выступления: первое в 1936 г. с критикой Н.~Н.~Лузина и попутно
П.~С.~Александрова и второе в 1939 г., направленное против решения
партгруппы Московского математического общества, оба тщательно подготовленные
и довольно резкие, быть может, даже агрессивные, были первыми
проявлениями моего боевого духа в борьбе за правое дело.
\end{quotation}

Кроме известных нам участников столкновения, в 1989 году оставались еще очевидцы. В принципе, их было не так уж мало,
но не все они писали воспоминания. Четверых  мы  упомянем.

В.~А.~Ефремович 1903 -1989

А.~П.~Юшкевич 1906-1993

Б.~В.~Гнеденко 1912-1995

С.~М.~Никольский 1905-2012

Что касается <<разделительной линии>>, то (если не считать ухода из жизни многих участников событий),
то в 1974 году она оставалась в принципе той же, что и в 1951 году -- три статьи в Успехах написаны все теми же лицами,
Ларентьевым, Л.~Келдыш и П.~И.~Кузнецовым. Но дальше картина начала меняться.

\sm

По-видимому, тогда же Александров задумался над тем, что после его смерти подробности 1936 года всплывут на поверхность.
В 1974 году в статье в <<{\it Науке и Жизни}>> \cite{Alex-nauka-i-zhizn'} он вспоминает Лузина%
\footnote{Статья не об этом. Она малоизвестна и, по-моему, интересна.},
а 1977 он году он публикует  \cite{Alex-kvant}, \cite{Alex-v-shkole}  восторженный рассказ о Лузитании.
В воспоминаниях \cite{Alex-auto1}, \cite{Alex-auto2} он (тщательно избегая деталей) пишет о своем конфликте с Лузиным
(см. ниже п.\ref{ss:aleksandrov-1936}).

\sm

Последние видимые отсветы конфликта относятся к 1980-1983году. Об этом и 
о том, как Колмогоров в 1983 году объявил мир, будет сказано в конце данных записок.

 \section{Московская математика и диссертация Лузина}
 
 \COUNTERS
 
 {\bf\punct CV Лузина до 1930 года.%
 \label{ss:CV}}
 
 9 декабря (27 ноября) 1883 г -- дата рождения;
 
 1893-1901 - учился в гимназии в Томске (один год  в Иркутске);
 
 1901 -- поступил на Физмат Московского университета;
 
 ноябрь или первые числа декабря 1905 -- вторая половина лета 1906 -- поездка в Париж;
 
 1907 женился на Н. М. Малыгиной;
 
 1909 сдал магистерские экзамены;
 
 1910 становится приват-доцентом;
 
 осень 1910-1912 работает в Гёттингене;
 
 1911 -- публикует первую работу;
 
 1912 -- публикует <<теорему Лузина>>;
 
 1913 -- май(?) 1914 - работает в Париже;
 
 С осени 1914 преподает в МГУ в звании приват-доцента;
 
 В 1915 - готовит диссертацию  <<{\it Интеграл и Тригонометрический ряд}>>. Отпечатана в типографии в 1915г.
 и  опубликована как статья в <<{\it Математическом Сборнике}>> в 1916 году  (242с);
 
 27.04.1916 защищает диссертацию. Совет присуждает ему степень доктора, минуя степень магистра;
 
 1917 ?  -- декабрь 1930  -- профессор Московского Университета;

 1918, октябрь -- осень 1920. Лузин  работает в Иваново-Вознесенском политехническом институте;
 
% 1919—1921 -- профессор Иваново-Вознесенского политехнического
%института.

 1919, март-октябрь, - и.о. декана физико-математического факультета МГУ%
 \footnote{Известно, что тогдашний декан Физмата М.~М.~Новиков в марте был избран ректором. Лузин упоминается в качестве и.о.декана, 15 октября
 	1919г. новым деканом был избран  А.~Н.~Реформатский, \cite{Ilchenko}.};

 1920-1925 - профессор Лесотехнического института \cite{Ryb-lesoteh};

август 1925-июль 1926 - командировка от Наркомпроса за границу;

январь 1927 - избран член-корреспондентом АН СССР;

май 1927 - командировка от Наркомпроса за границу;

 1.07.1928г. по 1.09.1930г.
%1928 (август?) - начало 1930 (?)
командировка от Наркомпроса за границу%
\footnote{Из Академической стенограммы:
\newline
ГОРБУНОВ.
  У  меня  такие  данные:  1905—1906  год — в  Париже.  1910—1912  год — три
года в  Германии,  частично в  Англии.  1912—1914 — три года в Париже.  1925—1926  год —
9  месяцев  в  Париже.  1926—1927  год — 5  месяцев в  Париже.  1928—1929  год - два года в
Париже. Вот только такие данные имеются.};

1928 -- вице-президент Международного математического конгресса в Болонье;

1929, январь. Избран в действительные члены Академии наук СССР, сначала
по кафедре философии, а затем по кафедре математики;

с 1929  -- cотрудник Физико-математического
института им. В.~А.~Стеклова АН СССР,  Ленинград;

1930 -- во Франции выходит монография Лузина <<Лекции об аналитических множествах и их приложениях>>;

В сентябре (?) 1930 года возвращается в Москву

с 20.11.1930 - устраивается на работу в ЦАГИ;

  16.12.1930 - уход из МГУ (по собственному желанию).

 \sm
 
 {\bf \punct Ранняя биография Лузина.%
 \label{ss:rannyaya}}
 Лузин о себе \cite{Ovanesov}:
 \begin{quotation}
 Томск стоит на реке Томь.
 За рекой - медвежьи берлоги. Я учился в <<классической гимназии>>. У моих сверстников 
 был культ физической силы; вполне понятно почему: близость столь сильных зверей, 
 как медведи, и рассказы о них заставляли видеть в физической силе высшее благо.
 
 Я был физически слабым (хотя и нормальным) и робким. Товарищи хотели меня приохотить 
 к их интересам и злоупотребляя своею силою, достигали обратного... Был в Томске единственный книжный магазин, в который тогдашняя <<Европейская Россия>> пересылала всякую ненужную литературу: была полная мешанина, разбираться к которой я быстро приохотился. Я был
  единственным ребенком своих родителей. Отец -- торговый служащий. Оба они были
   полуобразованы, но все читать и писать могли. Понимая мое отчуждение от товарищей, они предоставляли мне полную свободу действий внутри меня самого. За эту свободу -- вечное им спасибо! Я рос <<сам из себя>>, и моя голова была полна миром фантазии. Я читал буквально все,что видел на прилавке магазина, и, конечно, ничего не понимал. Я читал Канта <<Критику чистого разума>> ранее романов Жюль Верна. Книги философские больше всего меня привлекали, потому что я их не мого понять, и я искал тайного смысла их, хотел понять его, похитить
    даже <<силом>>
 \end{quotation}

  Из биографии Лузина, составленной  Голубевым и Бари \cite{GoBa}:
 \begin{quotation}
  В старших классах гимназии Н.~Н.~Лузин читал очень много
  и в самых разнообразных направлениях; книги по чистой философии увлекали его, давая воображению обильную пищу.
  Но математику до самых последних лет гимназии Н.~Н.~Лузин недолюбливал и боялся,
  так как царившая тогда всюду система преподавания ее была построена более на механической памяти:
  нужно было безукоризненно заучивать наизусть формулировки теорем и в точности памятью воспроизводить доказательства,
  по возможности не отступая от текста книги («Геометрия» Давидова, «Алгебра» Киселева).
  Для Н.~Н.~Лузина это было трудно переносимой мукой, так как механической памятью он совершенно не обладал; по этой же причине для него были закрыты история, география и языки, требовавшие запоминания времени, места и форм.
  Его занятия по математике шли в гимназии хуже и хуже, так что он утратил репутацию хорошего ученика
  и отец вынужден был взять для него «репетитора». К счастью, это был весьма талантливый  
  студент только что тогда открывшегося в г.~Томске Политехнического института; он произвел на Н.~Н.~Лузина сильнейшее впечатление тем, что показал ему математику не как систему механического заучивания, а как систему рассуждений, направляемую живым воображением. С тех пор он до некоторой степени утратил неприязнь к математике, перерешал самостоятельно все имевшиеся тогда задачники по элементарной математике и, естественно, в этом отношении стал в гимназии на первое место.
 \end{quotation}
 
 Теперь цитируем Бари \cite{TyulinaA}:
\begin{quotation}
 Было решено сделать из него 
инженера......

К поступлению в инженерную школу он особенно не тяготел, но и 
возражений каких-нибудь с его стороны не было.

Но было одно формальное препятствие - окончил гимназию
без особого блеска, не получил ни медали, которая могла бы
облегчить поступление в высшую школу, ни особенно блестящих
знаний, которые позволяли бы надеяться на успех в трудных 
конкурсных экзаменах, не привила гимназия Н.~Н. и особой 
решимости и настойчивости, которые позволили бы ему преодолеть 
встретившиеся трудности. Поэтому было решено сначала поступить на
Физико-математический факультет университета, окончание 
которого позволяло поступить без экзамена сразу на старшие курсы
Политехнического института....

Осенью 1901 г. Н.~Н.~Лузин был зачислен студентом первого
курса математического отделения физико-математического 
факультета Московского университета. Таким образом был сделан, как
думал его отец, первый шаг к выгодной и блестящей карьере 
инженера; но дальше все сложилось совсем не так, как предполагали
родители.

Как рассказывал сам Н.~Н., первая же лекция по высшей 
математике решила дело....

Большее влияние
оказало на Н.~Н. общение с молодым, тогда только начинавшим
свою многогранную преподавательскую и научную деятельность
профессором Д.~Ф.~Егоровым, ученым чрезвычайно разносторонним
и бывшим в курсе всех научных новинок. Наконец, несомненное
влияние на развитие Н.~Н. оказал ряд способных 
товарищей-студентов; вместе с ними Н.~Н. организовал математический
студенческий кружок, в котором с большим интересом 
разрабатывали вопросы, представлявшие в то время особую научную 
актуальность. Из участников этого кружка прежде всего необходимо
отметить оригинальную и самобытную личность Флоренского .....
В то время он [Флоренский] с необычайным интересом занимался 
вопросами философии, пытался обосновать на каких-то соображениях
Аристотеля теорию актуальных бесконечно-малых, 
многочисленные доклады на подобные темы и делались им в студенческом 
математическом кружке.
\end{quotation}

Так или иначе, наступила осень 1905 года (то есть Лузин уже четыре года проучился в Университете).
Продолжаем цитировать  Бари:
\begin{quotation}
Университет шумел, как улей; занятия осенью 1905 года
то начинались, то прекращались.
Аудитории превратились в
место сходок и массовой агитации....

Первые годы университетской учебы Н.~Н. жил в знаменитой
«Кокоревке», где он снимал номер, там же жили и его родители.
Теперь он, увлеченный бурным потоком общественного подъема,
тоже пытается принимать какое-то участие в революционном
движении. Проживание при таких условиях у всех на виду, в
большой гостинице, было явно нецелесообразным и по 
рекомендации кого-то из товарищей (В.~А.~Костицына%
\footnote{Костицын Владимир Александрович (1883—1963), естествоиспытатель,  окончил Московский Физмат,
	работы по математике, астрономии, климатологии, геофизике, экологии
	(в основном, применение математических методов в этих науках).
	 Известен также
бурной биографией.  Член РСДРП, в 1904г. получил шашкой в голову при разгоне демонстрации, участник Декабрьского восстания 1905г. в Москве и боев на Пресне.
Арестовывался в 1907-1908 году, исключен из МГУ, учился в Сорбонне. В 1917
году - заместитель комиссара Юга-Западного фронта. По некоторым утверждениям, арестовывал Деникина.
По минимуму, Деникина конвоировал. Из воспоминаний Грей, дочери Деникина,
\newline
{\it В конце митинга, около 5 часов вечера, два человека вошли в камеру Деникина:
Костицын, помощник Иорданского, и капитан Бетлинг,
командующий батальоном юнкеров, назначенный охранять заключенных...
Честный Костицын, в отличие от Иорданского, не был по природе своей кровожаден...}
\newline
(обрываю на полуслове,  Деникин тоже об этом рассказывает).
\newline
В бумаге, утвержденной Политбюро от 23 марта 1928г. (где обсуждается предстоявшая  операция
по советизации АН СССР), предлагается выяснить возможность включения Костицына 
в  список лиц, которых было бы желательно внедрить в АН \cite{Politburo}.
Однако в  1928-29гг. он эмигрировал.
Упоминается его сотрудничество с Вольтерра (Vito Volterra).
Утверждается, что в 1940 г. писал обращение о возврате в СССР.
 В 1941-1942гг находился в немецком концлагере.
Участвовал во французском подполье. После Войны получил советский загранпаспорт,
но на родину не вернулся.
В СССР  в 1984 году была издана его книга <<Эволюция атмосферы, биосферы и климата>>.
От него остались интересные воспоминания, которые после его смерти были переданы
в Советское посольство, а оттуда --  в Институт
марксизма-ленинизма. По-видимому, до сих пор они полностью не опубликованы,
отрывок есть в \cite{Kostitsyn}.}%
. - Авт. рук.) Н.~Н.
снимает комнату на Арбате, в семье вдовье врача Малыгина.
Семья состояла из старушки-вдовы Малыгиной и ее незамужней
дочери, Надежды Михайловны. Дом был тихий, внимание 
полиции не привлекал и в бурные дни октября 1905 г. перед 
появлением знаменитого виттевского манифеста «17 октября» в комнате
Н.~Н. не только ночевали нелегальные лица, но под его кроватью
был даже склад бомб...

Все это время Н.~Н. не прерывал занятий под руководством
Д.~Ф.~Егорова, и при создавшейся обстановке Д.Ф.Егоров 
посоветовал Н.~Н., чтобы не терять времени%
\footnote{Можно допустить и иную причину, учитывая увлеченность нашего героя <<бурным общественным подъёмом>>...}, уехать за границу и
учиться во время перерыва в университетских занятиях в одном
из заграничных университетов. Намечены были или Гёттинген,
или Париж, как наиболее крупные математические центры.
Д.~Ф.~Егорову удалось найти другого студента, который бывал за
границей и немного владел разговорным языком, французским и
немецким [В.~В.Голубев], и в первых числах 
декабря, перед самым вооруженным восстанием в Москве, Н.~Н. со 
своим спутником уехали в Париж.
\end{quotation}

Теперь цитируем записи Голубева \cite{Golubev-Omega}, <<Тетрадь Омега>>, 1942:
\begin{quotation}
 Я помню, когда Лузин и я собирались в конце 1905 года ехать за границу, то мы зашли
перед отъездом к Егорову. В разговоре Лузин попросил у Егорова указаний, что читать по
философии. Я тоже попросил таких указаний и по наивности предупредил Егорова, что я
материалист. Не помню, что мне указал Егоров, кажется, ничего не указал, но помню, что
он поморщился. Как я узнал много позднее, Егоров не только был свободомыслящим,
разносторонним человеком, как я себе представлял профессоров, а человеком верующим,
что было неплохо, и сверх того, до известной степени сектантом, что уже было много
хуже.

Я вспоминаю вечера, когда мы с Лузиным жили в Париже. Когда к нам в пансион
переселился мой товарищ по гимназии, Всеволод Вячеславович Елагин, то почти каждый
вечер происходило одно и то же. Часов в 10, когда мы с Лузиным кончали пить вечерний
чай, являлся Елагин и начинался бесконечный спор на религиозные темы. Лузин
доказывал существование Бога, бессмертие души и т. д., а Елагин доказывал, что нет ни
Бога, ни души. Я обычно не вмешивался в спор, только от времени до времени
иронизировал и подливал масла в огонь, когда спор начинал затихать.
\end{quotation}

Снова Бари:
\begin{quotation}
 В Париже Н.~Н. пробыл до конца летнего семестра 1906 года, и
все эти полгода пребывания за границей пройти в упорной и 
систематической работе. Лекций он слушал немного. Слушал Бореля в
Сорбонне, который читал теорию целых функций; так же слушал
знаменитого Пуанкаре, он читал разложение в ряды 
пертурбационных функций небесной механики; кроме того, в College de France
слушал Адамара, который читал теорию распространения волн.
Иногда ходил на лекции Дарбу по теории поверхностей. Но он
упорнейшим образом работал над изучением математической 
литературы в библиотеке Св. Женевъевы. Упорному изучению научных
вопросов посвящалось буквально все время. Над размышлениями над
научными вопросами просиживал целые ночи; часто поздно 
восходящее зимнее солнце заставало его еще за работой....

Н.~Н. вернулся в Россию во второй половине лета 1906 года.
\end{quotation}

Теперь Голубев и Бари:
\begin{quotation}
 ... В конце того же года он сдал государственный экзамен и был оставлен Д.~Ф.~Егоровым
 при университете «для приготовления к профессорскому званию». 
 В 1907 г. Н.~Н.~Лузин женился на Надежде Михайловне Малыгиной. За время обучения в университете
 Н.~Н.~Лузиным было прочитано и изучено много труднейших и глубоких трактатов по самым различным областям
 математики, так что он был хорошо подготовлен к магистерским экзаменам еще на студенческой скамье.
 «Время же оставления при университете он употребил на слушание лекций {\bf на медицинском факультете},
 куда намеревался поступить, чтобы впоследствии идти в народ, но потом был вынужден оставить этот план,
 так как работа в анатомическом театре оказалась ему не по силам. Тогда он перешел к слушанию лекций
 {\bf на философском отделении историко-филологического факультета}, 
 который через год оставил, потому что лекции по философии не давали указания на возможность творчества». 
\end{quotation}

Егоров Лузину 21 июня 1907, \cite{Kol-savvina-1}:
\begin{quotation}
...Я очень о Вас беспокоюсь и страшно жалею, что не удалось повидать Вас и поговорить с Вами. В письме всего не скажешь,
да к тому же я слишком мало знаю о о причинах Вашего угнетенного состояния духа. 
Так бы хотелось по мере сил помочь Вам и внушить хоть немного бодрости и более спокойного отношения ко всему окружающему.
Помните, что у Вас вся жизнь впереди, что Вам ни в коем случае не приходится отчаиваться: все еще можете вернуть, 
и изменить и наверстать. Если Вам временно показалось, что наука ничего Вам не дает, то это, во-первых, ошибка,
а во-вторых, перед Вами открыта вся жизнь, которая слагается из научного труда и многого другого.

Надеюсь, что Вы успокоитесь и что в скором времени я Вас увижу; тогда подробнее  поговорим.	
	\end{quotation}
	
%	В опубликованных письмах Лузина ???

Снова Бари:
\begin{quotation}
... В то время курс в Университете был четырехгодичным и, 
следовательно, ему нужно было окончить Университет в 1905 году, но
один год пропал из-за перерыва занятий. Между тем 
обстоятельства складывались так, что тянуть с окончанием курса было нельзя:
Н.~Н. существовал на средства родителей, а родители быстро 
дряхлели, да и торговые дела отца из-за революционных событий 
находились в самом плачевном состоянии. Поэтому первой заботой
Н.~Н. по возвращении из-за границы было сдать государственный
экзамен и получить соответствующие права. Экзамен был 
благополучно сдан в течение 1906-1907 года и Н.~Н. был оставлен, по 
обычаю того времени, на два года «для приготовления к 
профессорскому званию». Это несколько упрочило его материальное 
положение, так как оставленные при Университете получали в течение
двух лет по 600 рублей в год казенной стипендии. Для Н.~Н. это
было весьма существенно....

К этому времени относятся попытки обеспечить себя 
преподаванием математики в средней школе. Н.Н. взял уроки математики
в старших классах частной гимназии Флерова, которая 
помещалась в переулке около Никитского бульвара, но из этого ничего не
вышло. Опыта преподавания у него никакого не было, дисциплину
поддерживать в классе не умел и, несмотря на то, что отдельные
интересующиеся математикой учащиеся очень ценили уроки Н.~Н.,
подавляющее большинство были избалованные и плохо 
дисциплинированные дети зажиточных родителей, которые на уроках Н.~Н.
«ходили па головах»....Н.~Н. пришлось
отказаться от преподавания в средней школе....

К 1909 году Н.~Н. сдал так называемые магистерские 
экзамены, получил существовавшее тогда звание «магистранта» и 
права преподавания в высшей школе по прочтении двух пробных 
лекций, одной по собственному выбору и второй по назначению
факультета. Н.~Н. прочел пробные лекции и предполагал с осени
1909 г. читать в Университете курс по теории функций 
действительного переменного, но оказалось, что такой курс уже был
объявлен С.~С.~Бюшгенсом, который держал экзамены 
одновременно с Н.~Н. Тогда по совету Б.~К.~Млодзеевского Н.~Н. объявил курс
по теории интегральных уравнений. Читать этот курс Н.~Н. не
пришлось, так как Н.~Н. как раз в это время получил от 
факультета заграничную командировку для усовершенствования в 
математических науках в Гёттинген и Париж.
\end{quotation}

Так или иначе, осень 1910 года, Лузину почти 27 лет, у него ни одной опубликованной
работы, и у нас никаких сведений, о том, что у него было что-то неопубликованное.

\sm
 
{\bf\punct  Московская математика начала XX века.%
\label{ss:moskva}}
Люди, учившиеся математике в Московском университете начала XX века
обычно (хотя и не всегда%
\footnote{Некоторые из них в20-30е годы были  предметом поливания грязью со стороны передовых
кругов, в особенности Некрасов и Бугаев. И тот, и другой много и реакционно философствовали.
Не знаю как Бугаев, а Некрасов был крупный математик, и в этом качестве поношения он точно не заслужил.
Впрочем, передовым кругам всегда видней.}) тепло вспоминали своих профессоров, среди них

\sm

Андреев Константин Алексеевич (1848-1921)

Бугаев Николай Васильевич (1837—1903)

Власов Алексей Константинович (1868—1922) 

Егоров Дмитрий Федорович (1869-1931)

Жегалкин Иван Иванович (1869-1947)

Жуковский Николай Егорович (1847-1921)

Лахтин Леонид Кузьмич  (1863-1927)

Млодзеевский Болеслав Корнелиевич  (1858-1923)

Некрасов Павел Алексеевич  (1853—1924) 

Чаплыгин Сергей Алексеевич (1869-1942)

\sm

Если же мы будем вспоминать математиков и их математические результаты, то вспомним Егорова и его  теорему
из учебников анализа (действительного или функционального): для любой сходящейся п.в.
последовательности измеримых функций на отрезке найдется 
множество сколь угодно малой меры, на дополнении которого последовательность сходится равномерно%
\footnote{Об ее истории,  см. недавнюю статью В.~И.~Богачёва \cite{Bogachev-Egorov}, см. также \cite{Bogachev}.}.
Интересно, что Егоров занимался, в основном, дифференциальной геометрией, а вспомним мы именно теорему из ТФДП
(по ТФДП у Егорова была еще одна заметка).
Понятно, что люди из России должны помнить Жуковского и Чаплыгина, которые вообще-то были сильными математиками,
но прославились они не 
 как математики, а как аэродинамики, организаторы ЦАГИ и отцы-основатели аэрокосмической промышленности СССР.

А если вспоминать знаменитых русских математиков начала XX века, то окажется, что все они, кроме Егорова,
жили не в Москве. Процитируем Александрова \cite{Alex-MMO2}.
 \begin{quotation}
  При всех отдельных достижениях Московской математической школы
второй половины прошлого [XIX] века можно с достаточной уверенностью сказать,
что научные результаты московских математиков того времени были ниже
их творческих возможностей. В математической жизни Европы Москва
прошлого века занимала провинциальное положение.
 \end{quotation}

 И еще раз  Александров, \cite{Alex-MMO3}:
 \begin{quotation}
  Геометрическая школа была первой серьезной математической школой,
возникшей в Москве. Ее работы были тесно связаны с работами по классической
теории уравнений в частных производных в духе методов, развивавшихся
во Франции, начиная с Коши и кончая Дарбу. Вся эта область широко
разрабатывалась в Московском университете (в частности, Д.~Ф.~Егоровым
и С.~А.~Чаплыгиным). Но в полном смысле слова на мировое математическое
поприще Московское математическое общество и вместе с ним вообще московские
математики вышли уже в десятых годах текущего столетия, незадолго
до первой мировой войны, трудами Д.~Ф.~Егорова (1869—1931) и Н.~Н.~Лузина
(1883—1950), положивших начало московской школе теории множеств
и теории функций (сначала действительного, а потом и комплексного переменного).
 \end{quotation}

 \sm
 
 {\bf \punct Петербургская школа начала XX века.%
 \label{ss:peterburg}} Эта школа доминировала в России
 и к ней, например, относились такие знаменитости как
 А.~А.~Марков-старший, А.~М.~Ляпунов, А.~Н.~Коркин, Н.~Я.~Сонин, Г.~Ф.~Вороной.
 Картина однако была не столь уж однозначной. Приведем несколько примеров.

 В самом конце 1890х Ф.~Э.~Молин \cite{Molien-book}, \cite{Molien-Kanunov}
 искал место работы в российских университетах. У него были на руках письма
 поддержки от разных известных европейских людей, например, от Фробениуса
 (Ferdinand Georg Frobenius). Попытки, однако, оказывались неудачными.
 Опубликовано одно заключение конкурсной комиссии, куда входили А.~М.~Ляпунов, В.~А.~Стеклов,
 М.~Ф.~Ковальский и (с совещательным 
голосом) Л.~Струве.
 
  \begin{quotation}
 К 
сожалению, комиссия не могла составить себе самостоятельного
суждения о степени оригинальности и научного значения
работ г.~Молина, так как работы... относятся к области,
с которой, как стоящей в стороне от важнейших научных
дисциплин, члены комиссии знакомы лишь поверхностно..
Теория высших комплексных чисел представляет весьма;
сложное и искусственное построение, вызванное 
известным стремлением к обобщению понятия о числе, не
оправданному насущными потребностями....
  \end{quotation}
 
 Понятно, что отклонение заявок - дело обычное, как и обычна проявляемая при этом несправедливость.
 Но здесь интересна формулировка,
 и то, к чему она относится. <<Гиперкомплексные системы>> -- это на современном языке 
 ассоциативные алгебры. Результаты Молина были (в частности) такие
 
 -- любая конечномерная алгебра есть полупрямое
 произведение полупростой алгебры на радикал, а любая полупростая алгебра есть прямая сумма
 матричных алгебр. 
 
 --- соотношения ортогональности для матричных элементов 
 представлений конечных групп, 
 
 ---  групповая алгебра конечной группы разлагается в прямую сумму матричных
 алгебр.
 
 Cейчас это частично или полностью  читается  в университетских курсах алгебры,
 и входит в общие учебники алгебры (возможно, что и понятие фактор-алгебры пошло от Молина)%
 \footnote{Кстати, заявка в Харьков была на место астронома, а Молин по образованию был астроном,  его первые две работы
 были посвящены движению комет.}.
 
 Конечно, современному читателю, в связи с вышесказанным, приходят  на пямять фамилии Веддербарна и Фробениуса. 
 Но они о работах Молина (Theodor Molien) знали и на них ссылались.
 Веддербарн (Joseph  Wedderburn)  дальше обобщал теорему Молина на произвольные поля,
 а Фробениус параллельно c Молиным занимался теорией представлений. Собственно,
  Молин -- один основателей теории представлений наряду с Фробениусом и Бернсайдом.

 В итоге он в 1900г. нашел место профессора в Томском университете и в дальнейшем поднимал науку в Сибири
 (в нашей повести он еще появится).
 
 \sm
 
 Теперь цитируем статью Лаврентьева о Лузине:
 \begin{quotation}
  В книге [в диссертации Лузина, которая  была книгой, и по размерам, и по сути] ставились задачи
с наброском доказательств, в этих случаях попадались такие фразы:
«мне кажется», «я уверен».

Этот стиль не вписывался в классические традиции математических
работ; ленинградцы не признали монографию Н.~Н.~Лузина существенным
вкладом в науку. Академик В.~А.~Стеклов при чтении монографии Н.~Н.~Лузина
делал на полях много иронических замечаний: «ему кажется, а мне не
кажется», «гёттингенская болтовня» и т. п.
 \end{quotation}

А вот из письма петербургского академика В.~Я.~Успенского А.~Н.~Крылову 1926:
\begin{quotation}
 Относительно Лузина я знаю, что он хороший специалист в своей области (теория множеств 
 и связанная с нею канторовско-лебеговская дребедень), блестящий npoфeccop, создавший в Москве школу
 своих учеников и своим влиянием упразднивший настоящую математику в Москве\dots 
 Лично я Лузина почти не знаю и потому, 
 если я позволяю себе высказываться неодобрительно о его направлении, то только потому, что чувствую 
 к этому направлению глубокое
 омерзение и твёрдо уверен, что мода на это скоро пройдёт.
\end{quotation}

Успенский был прекрасный математик, в статье \cite{Uspensky1920} 1921г.,
он вывел формулу
для  асимптотики  числа разбиений%
\footnote{Напомним, что для натурального $n$ величина $p(n)$ определяется как число представлений $n$ в
виде суммы $n=k_1+\dots + k_m$, где $k_1\ge \dots \ge k_m$.},
\begin{multline*}
p(n)=\frac{\exp\Bigl\{\pi \sqrt {\frac23(n-1/24)}\Bigr\}}{4\sqrt 3 (n-1/24)}
	\Bigl[1 - \frac{\sqrt 3}{\sqrt{ 2\pi(n-1/24)}}\Bigr]
	+
	O\Bigl(\exp\Bigl\{\frac{\pi}{\sqrt 6} \sqrt n\Bigr\}\Bigr)
	\sim\\\sim
	\frac{e^{\pi \sqrt {2/3} \sqrt n}}{4\sqrt 3 n}
\end{multline*}
Эта замечательная формула была опубликована чуть раньше Харди
и Рамануджаном [Godfrey H. Hardy, Srinivasa Ramanujan], \cite{Hardy-Ram}, 1918,
но у них там гражданской войны не было). 

Работал Успенский в теории чисел и теории вероятностей%
\footnote{Из воспоминаний Б. А. Розенфельда \cite{Rozenfeld}:
	\newline
	{\it
		Много
		рассказывал мне Борис Николаевич [Делоне] о своем друге Якове 
		Викторовиче Успенском (1883—1947), ленинградском математике, академике,
		который во время одной из поездок за границу женился, привез жену
		в Ленинград и, так как ленинградская жизнь ей не понравилась,
		уехал с ней снова за границу.}
	\newline
	Уехав, Успенский сложил с себя должность академика, на освободившееся место с кафедры философии был переведен Лузин.}. Шел 1926 год,   задним числом 
видно, что Хинчин, занимавшийся вслед за Лузиным <<дребеденью>>, тогда уже заслужил титул классика теории вероятностей.
Да и Колмогоров уже подтягивался... Впрочем, на Крылова призывы Успенского не повлияли.

Это можно было бы продолжать. А.~П.~Юшкевич, рассказывая о русской математике конца  XIX -- начала XX века  
\cite{Yush1917}, сообщает (и это хорошо аргументирует) 
\begin{quotation}
 Прежде всего началось усвоение новых идей и направлений, главным образом, в периферийных университетах.
 \end{quotation}

 Вот Московский университет и был лучшим <<периферийным университетом>>.
 
 \sm
 
{\bf\punct  Москва и новые веяния.%
\label{ss:moskva-new}} 
 Снова цитируем Александрова \cite{Alex-1955}:
 \begin{quotation}
 Преподавание
математических наук в Московском университете в конце прошлого века
не было лишено черт провинциализма. С этими-то традициями и Б.~К.~Млодзеевский
и Д.~Ф.~Егоров повели решительную борьбу. Они, а также И.~И.~Жегалкин (1869—1947)
впервые стали читать курсы лекций, посвященные новым
областям математики, а старые области науки излагать в соответствии с ее
действительным состоянием в то время. Так, первое вполне современное
изложение вариационного исчисления, а в значительной степени и дифференциальных
уравнений принадлежит Д.~Ф.~Егорову (первые годы XX века).
Первый курс по теории функций действительного переменного был прочитан
Б.~К.~Млодзеевским (1901). Б.~К.~Млодзеевскому и Д.~Ф.~Егорову принадлежит
огромная заслуга введения в университетское преподавание научных
семинаров. Это нововведение оказалось в высшей степени удачным: семинары,
как известно, стали в настоящее время основной формой научного воспитания
студентов старших курсов и аспирантов. Пишущий эти строки, еще будучи
студентом второго семестра, участвовал в семинаре Д.~Ф.~Егорова (1914),
темой которого были «бесконечные последовательности». В том же семинаре
участвовали А.~Я.~Хинчин, Д.~Е.~Меньшов и др. Участники этого семинара
получили первое соприкосновение с настоящей живой математической наукой
— впечатление, которое уже не может никогда изгладиться!
\end{quotation}
 
 Стоит иметь в виду, что этот автор этого отрывка впоследствии поссорился с Лузиным. А в 1914-1915 годах
 он (так же как  Хинчин и Меньшов) находился под влиянием Лузина...
 
 Теперь цитата из Бари \cite{TyulinaA}:
 \begin{quotation}
  Среди... преподавателей [университета] появился ряд молодых энергичных и 
талантливых ученых, приват-доцентов и молодых профессоров, как
С.~А.~Чаплыгин, И.~И.~Жегалкин, А.~П.~Поляков и другие. А вместе с
тем коренным образом изменился и характер преподавания. Если
до 1905 года объем преподавания определялся в основном 
выполнением обычных министерских программ, то теперь ведущая роль
все более и более стала переходить в чтение факультативных, 
необязательных для студентов курсов, относящихся к современным,
актуальным вопросам науки, в которых лекторы ставили своей 
задачей вовлечь учащихся в решение вопросов, стоящих в то время
перед наукой... Механический кабинет ...
Н.~Е.~Жуковского, все более и более превращался в эти годы в исследовательскую
лабораторию, где велась напряженная разработка самых 
актуальных задач теоретической и прикладной механики, гидромеханики
и авиации...
 \end{quotation}

 \sm
 
 {\bf\punct Жегалкин.%
 \label{ss:zhegalkin}} Процитируем воспоминания  Н.~М.~Бескина%
\footnote{Бескин  Николай Михайлович, деятель школьного математического образования,
	автор задачников и методических пособий, см. Математика в школе, 1965, 1.}
 (который поступил
 на Физмат МГУ в 1921 году) о Лузине, опубликованные в Математике в школе \cite{Besk}, 1983:
 \begin{quotation}
\dots  (напомню, что Лузин считал своими учителями Д.~Ф.~Егорова, Ивана Ивановича Жегалкина и Жюля Таннери)\dots
 \end{quotation}
 
Забавно, что сентенция, брошенная между делом
 как общеизвестная,    в прочей  литературе о Лузине отсутствует.
 Но прислушаться к ней стоит. 
 
 Известно (и об этом много написано ниже) Лузин придавал огромное значение технологиям образования. Его иногда
 поругивали  за неправильность поведения как препода (но чаще этим поведением восхищались),
 и много честили как автора учебников, написанных  <<не по правилам>>.
 Однако его эффективность в обеих этих ипостасях очевидна. Что касается упомянутых героев, то:
 
 \sm
 
 Jules Tannery (1848-1910) -- французский математик, больше известный
 как деятель образования и реформатор образования. Лузин писал о нем в замечательном
 предисловии к учебнику мат.анализа Жегалкина и Слудской \cite{Luzin-Zhegalkin} (оно цитируется нами ниже 
 в п. \ref{ss:zhegalkin-preface}),
 там же Лузин высказывает солидарность с педагогическими методами Жегалкина.

 \sm
 
 Жегалкин (1869-1947) в 1902 году становится приват-доцентом Московского Университета. В 1907 публикует
 книгу <<Трансфинитные числа>> (340с) с весьма продвинутым изложением этого предмета.
 Очевидно,
 вместе с этой книгой канторовская теория множеств (теория ординалов и кардиналов, которые,
 сейчас среди среди работающих математиков не очень популярны)
 становится известной в Москве. Ясно и то, что дальше 
 многочисленные московские математики относились к трансфинитным числам как к чему-то привычному и обыденному.
 В любом случае, Лузин имел возможность еще в Москве познакомиться с продвинутой теорией множеств в качественном изложении. В 20х годах 
 Жегалкин публикует в <<{\it Математическом сборнике}>> несколько статей по математической логике, 
 а в 1930м или в самом начале тридцатых начинает вести в Московском университете семинар по логике.
 Он также был автором нескольких учебников и задачников по анализу. Когда в конце 1946
 года был высочайший указ о введении предмета <<Логика>> в школах, Жегалкин начал писать школьный учебник 
 (см. \cite{za30}, статья С.~А.~Яновской), который
 остался незавершенным.
 
 \sm

 {\bf \punct Московское математическое общество и Егоров.%
 \label{ss:egorov}}
 Советская власть не любила организаций граждан, с нелегальными организациями она обходилась решительно,
а легальных  было немного.
 Среди них, однако были различные профессиональные общества. Московское математическое
 общество в течение досоветского и всего советского периода играло большую положительную
 роль, каковая почему-то начала постепенно рассеиваться с концом Советской власти.
 
 В 1905-1921 году президентом Общества был Жуковский, он умер, не пережив смерти от туберкулёза своей единственной дочери.
 
 Новым президентом стал Млодзеевский. С этим избранием был связан первый 
 известный нам конфликт с участием Лузина (по-видимому, он был заметным событием, поскольку о нем говорят
 несколько независимых источников, \cite{TyulinaA}, \cite{LL}, \cite{Alex-MMO3}). Лузин хотел провести в 
 президенты Егорова, тот, по-видимому, сам отказался.  Лузин некоторое время был в оппозиции, 
 потом стороны помирились. Как будто ничего особенного, все-таки не было.
 
 Млодзеевский также стал первым директором основанного (по инициативе О.~Ю.~Шмидта) в 1921 году НИИ математики и механики при Физмате МГУ.
 Умер он в самом начале 1923 года.
 
 До середины тридцатых это имя для левых было  ругательным, но 
 в статье <<Математика>> в БСЭ \cite{Kolmogorov-matematika} 1938, Колмогоров,  говоря о Московской геометрической
 школе, произносит слово <<первоклассная>> и упоминает Млодзеевского,
 кроме того в Энциклопедии выходит и вполне положительная статья о самом Млодзеевском.
В положительном ключе о нем пишет
в статье 1946 года \cite{Alex-MMO1} Александров,
 01.12.1948 в Московском матобществе состоялось <<{\it заседание памяти Болеслава Корнелиевича Млодзеевского в связи
 с 25-летием со дня его кончины}>>. \cite{Uspehi-49}. Это к вопросу о неоднородности тогдашних времен.
 
 В 1923 году было избрано новое руководство Общества, 
 Егоров - президент, Лузин -- вице-президент, Привалов - секретарь. Таким
 оно оставалось до осени 1930 года. Тогда же Егоров становится директором НИИ математики и механики.
 В 1924 году его выбирают членом-корреспондентом Российской Академии наук.

В 1922 году Егоров возобновил издание Мат.Сборника. Как вспоминал Люстерник
\begin{quotation}
 Выход
в свет выпуска журнала в 1922 г. после четырехлетнего перерыва был отмечен 
на заседании Математического общества как праздничное событие.
Мне пришлось тогда поздравлять Математическое общество от имени студенческого математического кружка при нем%
\footnote{Похоже, что с начала Гражданской войны и до этого момента основным советским математическим журналом были
<<{\it Ученые записки Иваново-Вознесенского политехнического института}>>.}.
\end{quotation}

Журнал стал принимать статьи на иностранных языках (французский, немецкий, английский), что делало
журнал отчасти понятным зарубежным читателям (и тем самым увеличивало сферу распространения журнала).
Кроме того, стали приниматься статьи иностранных авторов (подробнее о журнале того времени, см. \cite{Demidov-sbornik}).

В 20е годы <<{\it Сборник}>> был единственным регулярным центральным математическим
журналом в СССР, хотя были менее заметные повторяющиеся издания в Ленинграде, Казани, Харькове, Ташкенте.

\sm

 Время Егорова потом долго (открыто - после Войны) вспоминали как золотой для математики век.
В 20е годы в МГУ было неспокойно, Егоров как <<реакционер>> подвергался атакам революционной молодежи (о чем ниже п.\ref{ss:egorov-war}),
но, так или иначе, до конца 1929 года он удерживался на посту директора Института механики и математики.

\sm

Среди учеников Егорова, кроме Лузина, известны:

\sm

Голубев Владимир Васильевич (1884-1954)

Привалов Иван Иванович (1891-1941)

Степанов Вячеслав Васильевич (1889-1950)

Петровский Иван Георгиевич (1901-1973)

Сретенский Леонид Николаевич (1902-1973) (в качестве научного руководителя указывается также Чаплыгин))

Фиников Сергей Павлович (1983-1964)

Бюшгенс Сергей Сергеевич (1882-1963) (в качестве научного руководителя указывается также Млодзеевский)

Вениаминов Владимир Николаевич (1895-1932)

Костицын Владимир Александрович 1883—1963 (?%
\footnote{Упоминания о том, что Костицын -- ученик Егорова,  встречаются в нематематических источниках.
	 Судя по первым публикациям Костицына,
	это может соответствовать истине, правда доучивался Костицын уже в Париже.}) 

\sm

Привалов, Степанов и Вениаминов%
\footnote{О Вениаминове известно мало. С начала 20х был начальником кафедры в Академии воздушного флота им. проф. Н.~Е.~Жуковского.
Работал также в химико-технологическом институте. Автор нескольких работ по теории функций 
комплексного и вещественного переменного и по гидродинамике, а также нескольких вузовских учебников. Покончил с собой.}
числятся в списках Лузитании. 
Судя по всему, Лузин также оказал существенное  влияние
на Голубева, его магистерская (1916) и несостоявшаяся (по причине отмены ученых степеней) докторская
диссертация (обе были переизданы в \cite{Golubev-trudy}), по-моему, отразили влияние  ТФДП на комплексный анализ%
%В любом случае,  работы Голубева 1914-1924гг.,
%посвященные комплексному анализу, по духу  родственны работам московской школы теории %функций%
\footnote{Александров \cite{Alex-1945}, 1945:
	\newline
	{\it
	 Под значительным влиянием методов Лузина развивались и первые, собственно математические исследования В. В. Голубева (род. 1884), относящиеся к весьма тонким и трудным вопросам теории функций комплексного переменного.}}.
 В воспоминаниях о Лузитании Голубев не значится, потому что в 1918-30гг.
 он работал в Саратове.

Петровский, судя по всему, не имел прямого отношения ни к Лузину, ни к Лузитании.

Фиников и Бюшгенс -- профессора Мехмата, представители классической дифференциальной геометрии в духе Егорова.
Сретенский -- известный механик.

Интересно, что двое из учеников Егорова, как и он сам, приобрели известность также как положительные
административные деятели. Петровский -- декан Мехмата 1940-44 и знаменитый ректор МГУ 1951-1973.
Голубев был деканом физмата Саратовского университета (1918-19), потом ректором (1921-23),
первым деканом Мехмата МГУ (1933-1934) и снова деканом (1944-1952), директором НИИ механики (1936–1941, 1942–1950),
заведующим многочисленными кафедрами в разных институтах, гражданских и военных, имел чин генерала-майора
инженерно-авиационной службы (1944). Внес большой вклад в постановку обучения на Мехмате.

\sm

 {\small
{\bf \punct Геометрические работы Егорова.}
По работам Егорова был обзор П.~И.~Кузнецова \cite{KuzEgor} 1975г.
О Егорове очень много писали последние 25 лет, однако 
содержательного обсуждения его работ (кроме истории <<теоремы Егорова>> \cite{Bogachev-Egorov}) почти не было,
см. интересную статью \cite{Fomenko-Egorov}, к сожалению, слишком краткую. 
Автор ни в коей степени не берет на себя непростой труд
закрыть эту брешь.

По-видимому, основной его работой Егорова по геометрии был 240-страничный мемуар (фактически книга)
<<Об одном классе ортогональных систем>> \cite{Egorov-main} 1901г. Сначала он рассматривает римановы метрики
на двумерных поверхностях вида
$$
ds^2=\frac{\partial\omega}{\partial x} dx^2+ \frac{\partial\omega}{\partial y} dy^2
$$
где $\omega$ -- фиксированная функция. Такие римановы метрики он называет {\it потенциальными}.
Понятно, что семейство кривых $x=\alpha$ ортогонально семейству кривых $y=\beta$ во всех точках.
Это выполнено для всех метрик вида $f(x)dx^2+g(y)dy^2$, а класс потенциальных метрик Егоров выделяет
из следующего соображения: при параллельном переносе вдоль кривых $x-y=const$ (в смысле римановой связности)
касательная плоскость движется параллельно относительно системы координат.
Смысл мотивировки не очень ясен (по крайней мере автору этих записок), но дальше оказывается, что
такие системы обладают разнообразными интересными свойствами, что оправдывает их определение задним числом.
Егоров
исследует условия существования таких метрик на данной поверхности (при фиксированной римановой метрике
они есть и зависят от трех функций одной переменной).

Далее Егоров вводит {\it потенциальные поверхности} в $\R^3$. Для этого рассматривается сетка из линий кривизны,
риманова метрика записывается в этой сетке, 
и требуется, чтобы эта метрика была потенциальной относительно 
этой сетки (при подходящем выборе координат на двух ортогональных линиях кривизны).
Оказывается, например, что линии кривизны потенциальной поверхности при гауссовом отображении
переходят в потенциальную систему на сфере.

Далее рассматриваются потенциальные системы координат в пространстве. Оказывается, что координатные
поверхности сами будут потенциальными, а любая потенциальная поверхность в $\R^3$
может быть реализована как координатная поверхность трехмерной координатной системы.
А если, например, взять координатные линии трехмерной системы, 
то в любой точке их соприкасающиеся плоскости имеют общую прямую. Это свойство характеризует трехмерные потенциальные системы...

Мы примерно обрисовали объект исследования (и очень далеки от идеи излагать его содержание).
Это сложная работа с нелинейными УрЧП и со странными и необычными (во всяком случае, для автора этих записок)
красотами.

Результаты Егорова частично излагалась  во втором издании книги 
{\it Le\c{c}ons sur les syst\`emes orthogonaux et les coordonn\'ees curvilignes} Дарбу (Gaston Darboux) \cite{Darboux},
потенциальные системы он назвал (E)-системами (фактически Дарбу ссылается на заметки Егорова в Compt. Rendus).

По-видимому, это деятельность Егорова не имела прямых продолжений. Впрочем, в 90х годах было обнаружено,
что условия плоскости потенциальной метрики $\sum_j \frac{\partial \omega}{\partial x_j} dx_j^2$
связаны с WDVV-уравнением, см. \cite{Dubrovin}.}

\sm

{\bf\punct Старт Лузина.%
\label{ss:start}} Так или иначе, в конце 1910 года Егоров отправляет Лузина в 
командировку в Гёттинген. Цитируем Бари: 

\begin{quotation}
 ... насколько можно
судить по рассказам Н.~Н., наиболее тесная связь у него была ....
с молодым и весьма энергичным ученым, который
только что начинал свою работу в Гёттингене, профессором 
Эдмундом Ландау [Edmund Landau]. С Ландау хорошо был знаком Д.~Ф.~Егоров, учитель
Н.~Н. Несомненно, Н.~Н. имел от Д.~Ф.~Егорова рекомендации к
Ландау; вероятно, в семинарах Ландау и началась работа Н.~Н. и
его общение с молодыми заграничными учеными, которые всегда
собирались в большом количестве в Гёттингене, особенно в летние
месяцы. Во всяком случае, первая печатная работа Н.~Н., 
несомненно, носит на себе некоторое влияние Ландау....

Но уже работы, напечатанные Н.Н. в следующем 1912 году,
ярко свидетельствуют об исключительной самостоятельности
научного творчества Н.~Н....
\end{quotation}

Лузин -- Флоренскому, 06.04.1912 \cite{Luz-Flor}:
\begin{quotation}
Во время командировки мною предпринята работа о тригонометрических рядах Fourier, которая почти доведена до конца. 
Вероятно, это и будет первой диссертацией [в письме вскользь упоминается результат \cite{Luzin-start} 
и высказывается  гипотеза о том, что ряд Фурье сходится <<кроме отдельных точек>> (это не так)]...
Работа была предпринята главным образом для меня самого, так как после неудачи с
Continuum problem мне просто хотелось проверить свои силы и выяснить причины неудачи (и их происхождение)%
\footnote{Еще по крайней мере двое выдающихся математиков -- героев настоящих записок -- вскоре столкнулись с неудачами
и внутренним кризисом при атаках на континуум-гипотезу...}.
Так как работа, делаемая не случайно, поверхностной аналогией и <<вдохновением>> (как я хотел это ранее), 
но упорным напряжением, удалась почти вся (сейчас заканчиваю) -- то я получил некоторую бодрость
и лучше чувствую себя (духовно). Если все доведу до конца, то дальше не буду нервничать и буду спокойно работать.
Абсолютная неуверенность в себе давно уже гнела меня. Огромный недостаток было отсутствие школы и простого умения работать.
В этом смысле я многому обязан немцам (Hilbert'у, Landau).
\end{quotation}

В Гёттингене случилось еще одно событие. Цитируем книгу Александрова и В.~В.~Немыцкого о Степанове
\cite{Stepanov-AN}:
\begin{quotation}
По окончании университета в 1912 году
В.~В.~Степанова и И.~И.~Привалова оставляют при
университете, т. е., употребляя современную 
терминологию, они поступают в аспирантуру и вскоре 
направляются в научную командировку в Гёттинген.
Там они слушают лекции знаменитого математика
Д.~Гильберта, а также и другого немецкого 
профессора Э.~Ландау, о которых. Вячеслав Васильевич
впоследствии много вспоминал...

В Гёттингене происходит знакомство В.~В.~Степанова и И.~И.~Привалова
 с Н.~Н.~Лузиным,
представлявшим собой в то время самую яркую
восходящую звезду на небосклоне московской, да
и не только московской математики....

Первая самостоятельная математическая 
работа В.~В.~Степанова «К принципу Дюбуа-Реймона в
теории роста функций»  была сделана им в
1915 году и опубликована в Математическом 
сборнике. В этой работе была решена задача, которую
поставил перед В.~В.~Степановым Н.~Н.~Лузин....
\end{quotation}

Любопытно, что Привалов и Степанов, надо думать
знакомые с Лузиным ранее, попадают в сферу его влияния в Гёттингене. 

Снова Бари:
\begin{quotation}
За два года заграничной командировки им была проделана, 
огромная работа и между прочим было напечатано восемь научных
работ в лучших русских и заграничных научных журналах....

Отчет Н.~Н. о заграничной командировке, пересланный в 
Министерство народного просвещения профессором Д.~Ф.~Егоровым,
вызвал высокую оценку; ученый комитет Министерства постановил
признать занятия Н.~Н. во время заграничной командировки 
успешными, отчет о них вполне удовлетворительным и признал 
желательным продолжить заграничную командировку еще на один
год с той же стипендией (2000 руб. в год), что 22 июня и было 
утверждено министром....

Этот дополнительный год заграничной командировки Н.~Н.
провел в Париже. Здесь он систематически работал в семинаре
Адамара, завел личное знакомство с крупнейшими 
математиками, в особенности с молодежью (Пикар, Адамар, Борель, Лебег,
Данжуа [Charles  Picard, Jacques  Hadamard, \'Emile Borel, Henri  Lebesgue] и ряд других) и полностью проникся идеями 
французской математической школы, в особенности духом неподражаемого
изящества, легкости формы изложения, блеска лекций, которые
отличали французскую математическую науку.
\end{quotation}

Бари  -- источник пролузинский, ее можно было бы заподозрить в излишней патетичности,
но все это согласуется с  другими данными (можно лишь усомниться в <<молодости>> части
упомянутых французов). Во время Первой Мировой Войны Лузин будет переправлять в Comptes Rendus
академику Адамару (статьи в Comptes Rendus представлялись академиками)
заметки никому не известных молодых людей из Москвы. В заметках обычно не помещаются
полные доказательства. Понятно, что их публикация зависела во многом от честного слова Лузина. В 1936 году
Лебег и Данжуа решительно окажутся на лузинской стороне, а советские математики  поимеют
проблемы с публикациями во Франции.

\sm

{\bf \punct Первые работы.%
\label{ss:first}}
О работах по ТФДП, см., например, \cite{BaLyu1}-\cite{BaLyu2}, \cite{Ulya2}).
Мы лишь упомянем три работы Лузина этого периода.

\sm

Первой публикацией \cite{Luzin1911}, 1911, был  пример аналитической функции $F(z)=\sum с_k z^k$, у которой коэффициенты Тейлора
стремятся к нулю, а ряд Тейлора расходится во всех точках окружности $|z|=1$. Это было в теме тогдашних
дискуссий теории функций, и предполагалось, что наоборот, должна быть сходимость почти всюду. 

Контрпример строился примерно так. Рассмотрим сумму 
$$T_n(z)=1+z+z^2+\dots+z^n, \qquad z=e^{i\phi}$$
Вблизи точки $z=1$ выражение $|T_n(z)|$ близко к $n+1$, возьмем маленький
интервал $\phi\in (-\epsilon,\epsilon)$, где $|T_n(z)|>n$.
Далее возьмем сумму
$$
S_n:=\sum_{l\ge 0:\,2\pi l\epsilon \le 1} \Bigl[z^{l(n+1)} T(e^{2\pi l\epsilon} z)\Bigr]
$$
В любой точке окружности хотя бы одно суммы $S_n$ будет по модулю больше $n$.
Теперь выпишем сам ряд:
$$
F(z)=\sum_{n=0}^\infty a_n z^{r_n} S_n(z), 
$$
где $a_n$ --  последовательность, стремящаяся к нулю, причем
$a_n>1/n$,
 а 
последовательность натуральных чисел $r_n$ растет достаточно быстро,
так, что при $p>q$ любой показатель степени при слагаемом $z^j$ в $S_p$ будет больше 
любого такого показателя для $S_q$.
 В частности, при раскрытии всех скобок и получении ряда 
$\sum c_k z^k$ ни разу не случится приведения подобных членов.
Выражения, стоящие выше в квадратных скобках, имеют вид 
$
\sum_{k=N}^M c_k z^k
$, 
и по построению эти выражения не стремятся к нулю... Контрпример вышел весьма изысканный.

\sm

Вторая работа, 1912, содержала  <<теорему Лузина>> \cite{Luzin-start}, которая вскоре вошла в учебники:
для любой измеримой функции $f$ на отрезке и любого сколь угодно малого $\epsilon>0$
существует  непрерывная функция $g(x)$, такая, что  $f(x)\ne g(x)$ 
на множестве $C$ меры $\ge\epsilon$. Множество $C$ можно выбрать совершенным.
Лузин сам же немедленно начал применять эту теорему. 

\sm

В том же 1912 году Лузин показывает, что у любой конечной п.в. измеримой функции  $f(z)$ есть первообразная,
то есть {\it непрерывная} функция $F(x)$ такая, что $F'(x)=f(x)$ п.в.%
\footnote{Сравнительно недавно был обнаружен
	многомерный аналог этой теоремы (Moonens, Pfejfer, 2008):
\newline
{\it Пусть $U$ – область в $\R^n$ и $F: U \to R^n$ – измеримое по Лебегу
	отображение. Тогда найдется такое непрерывное отображение 
	$G:U\to \R^n$, что $G$ почти всюду в $U$ дифференцируемо и 
	его градиент совпадает с $F$ почти всюду
	в $U$.}} 

Теорема странная, и, скорее всего, тогда она была неожиданной.
Возьмем, например, функцию $x^{-2}$ на отрезке $[-1,1]$. Обычная первообразная $-x^{-1}$ тут не подходит,
она разрывна. С другой стороны, возьмем какую-нибудь <<канторовскую лестницу>> $k(x)$, то есть непрерывную функцию, производная которой почти всюду равна нулю. Тогда
$(F(x)+k(x))'=F'(x)$ почти всюду, т.е., первообразная в этом смысле определена не однозначно.
Теперь мы можем подобрать канторову лестницу  так, что особенность $-x^{-1}+k(x)$ в нуле станет устранимой.
Если же взять производную от $-x^{-1}+k(x)$ в смысле обобщенных функций, то получится $x^{-2}$ 
плюс (знаконеопределенная) сингулярная мера бесконечной вариации.

 Лузин показывает,
что такая первообразная есть не только в этом тривиальном примере, но и всегда.
Первое используемое соображение -- что для любой непрерывной функции
$g(x)$ и любого $\epsilon>0$ мы можем выбрать первообразную $G(x)$ так, что $|G(x)|<\epsilon$. Дальше применим теорему Лузина
и возьмем совершенное множество $C$, на котором функция
данная измеримая функция $f$
непрерывна. Возьмем непрерывную функцию $g$, равную $f$ на $C$ и линейную
на дополнительных интервалах. Возьмем у нее маленькую первообразную $G$. 
Дальше надо <<разыграть>> увеличение множества $C$ (см. \cite{Luzin-1916}, \S\S14-15).

\sm

{\bf\punct Возвращение в Москву.%
\label{ss:return}}
В 1911г. Лузин публикует упомянутый контрпример, в 1912ом - аж целых 7 статей...
Вскоре следует какой-то очередной кризис. 

Егоров -- Лузину, 7 апреля 1914 \cite{Kol-savvina-1}:
\begin{quotation}
Вы все молчите. Очевидно, с Вами неладно. Но как бы дело ни стояло, Вы должны помнить, 
что Вам дарован математический талант, и Вы должны его сохранить, а равным образом и постараться,
чтобы пришедшие Вам счастливые идеи не остались бесплодными. А поэтому возьмите себя в руки и займитесь прежде всего
своим моральным лечением. Без нравственной уравновешенности научная работа не может успешно идти.

Думаю, что Вам лучше бы вернуться в Москву, раз в Париже Вы оказались банкротом. Здесь Вы подтянетесь,
да и поддержку со стороны, может быть, найдете; а там далее и за работу приметесь.... 
\end{quotation}

Из интервью Меньшова \cite{Menshov-Duvakin}:
\begin{quotation}
 Егоров знал, что Лузин приедет осенью, и решил, насколько я теперь понимаю, подготовить для него хорошую аудиторию.
 И поскольку Лузин был специалист по теории функций, то в весеннем семестре 14-го года Егоров объявил специальный семинар по теории функций. Но это был семинар все-таки довольно официальный. 
 Было шесть групп. В каждой из групп было по нескольку человек, не больше шести, по разным отделам теории функций
\end{quotation}

Скорее всего, замысел Егорова - это домысел Меньшова.
Но так или иначе, <<хорошая аудитория>> для Лузина (как мы видим чуть ниже, лучше не придумаешь),
была подготовлена.

\sm

{\bf\punct  Диссертация Лузина.%
\label{ss:disser}} 
Весной 1914 года Лузин возвращается в Москву и в 1915 публикует диссертацию 
<<{\it Интеграл и Тригонометрический ряд}>>, Москва, Типография Лисснера и Собко.
Он же была опубликована в 1916 году как 242-страничная статья в Мат. Сборнике \cite{Luzin-1916}.

Работа содержит предшествующие результаты Лузина. 
Дальше извлекаются следствия из теоремы о существовании первообразной.
 Показывается,
что любая измеримая функция (конечная п.в) на окружности
является граничным значением (в смысле некасательного предела) некоторой гармонической функции.
Кроме того, доказывается, что любая 
измеримая конечная почти всюду функция есть (с  точностью до почти всюду)
сумма тригонометрического ряда (надо
иметь в виду, что коэффициентов Фурье у этой функции, вообще говоря,  нет, потому что
задающие их интегралы расходятся).
Сумма ряда вычисляется по Абелю--Пуассону, т.е.,
$$
\sum:=\lim_{t\to 1-0} \sum c_k t^{|k|} e^{ik\phi}. 
$$
Как показал  позже Меньшов, без оговорки про Абеля--Пуассона можно обойтись.

\sm

Значительная часть диссертации посвящена обобщениям интеграла Лебега
(Данжуа, Бореля, Перрона[Oskar Perron]), взаимоотношениям между разными известными 
обобщениями, свойствам получаемых таким образом первообразным. Сейчас,
сто лет спустя эта тематика (тогда очень модная) выглядит несколько экзотично
(впрочем за прошедшее время это отчасти упростилось).
 Лузин и здесь находит изящные утверждения,
например дает простое описание функций, которые могут первообразными
по Данжуа.
 
 \sm
 
Кроме этого, диссертация содержит много различных утверждений, логических связок, контрпримеров. Приведем два образца. 

\sm

Как известно, любая гармоническая функция
$u(x,y)$ двух переменных является вещественной частью
некоторой голоморфной функции $u(x,y)+iv(x,y)$,
гармонические функции $u$ и $v$ называются {\it сопряженными}.
Рассмотрим теперь гармонические функции в единичном круге с $u(0,0)=0$.
Они (с некоторыми оговорками) определяются своими граничными
значениями. Поэтому переход от функции к сопряженной определяет линейное
преобразование в пространстве функций на окружности $e^{i\phi}=x+iz$,
$$
J (\sum_n (a_n \cos n\phi+ b_n\sin n\phi)=\sum_n (- b_n \cos n\phi+ a_n\sin n\phi).
$$
Лузин замечает, что этот оператор задается формулой
\begin{multline*}
J f(\phi)=\frac 1{\pi} \mathrm {v.p} \int_0^{2\pi} \cot \Bigl(\frac{\phi-\psi}2  \Bigr) \,
f(\psi)\,d\psi:=\\:= \lim_{\epsilon\to +0} \frac 1{\pi} \int_{\epsilon}^{\pi}
\bigl(f(\phi+\alpha)-f(\phi-\alpha)\bigr) \cot \frac \alpha 2 \, d\alpha.
\end{multline*}
и показывает, что последний предел определен для почти всех
$\phi$, если $f$ принадлежит $L^2$.

\sm

Это, конечно, преобразование Гильберта, и приоритет 
тоже у Гильберта
 (знаменитую книгу  \cite{Hilbert-book}  Лузин то ли не знал,
 то ли не заметил там соответствующего отрывка). Однако, возможность
 определить это преобразование с помощью сингулярного интеграла {\it для всех}
 функций из $L^2$ --  утверждение тонкое, и оно было новым.
 
\sm

А вот другой пример (кстати, выводится с помощью преобразования Гильберта):
\begin{quotation}
... всякое измеримое множество всеми своими частями,
включая части меры бесконечно малой порядка выше 
первого, расположено симметрично относительно почти каждой точки области, 
если пренебрегать бесконечно малыми
расстояниями порядка выше первого%
\footnote{Здесь стоит вспомнить замечательную теорему плотности Лебега:
для любого измеримого множества $A$ отрезка 
$\lim_{\epsilon\to 0} \frac 1{2\epsilon}
 \mu(A\cap [x-\epsilon,x+\epsilon])$
равен 0 или 1 для почти всех точек отрезка. Утверждение Лузина весьма удивительно,
я, впрочем, не знаю, удалось ли в итоге из него извлечь дальнейшие следствия.}.
\end{quotation}

В диссертации обсуждается состояние теории функций действительного переменного на тот день, рассказывается о том, что  понятно и что непонятно,
какие задачи можно решать. Доказательства там иногда есть, а иногда лишь намечаются,
это позволяет при небольшом объеме сделать работу разнообразной и насыщенной. 
Если ее читать, то даже при отсутствии склонности к ТФДП (у автора настоящих заметок ее
уже давно
нет),  книга начинает завораживать.

\medskip

Диссертация на степень магистра защищалась в начале 1916 года, оппонентами были Егоров и Лахтин. 
Забавно, что стиль написания статей
за 100 лет изменился заметно, а стиль отзывов на диссертацию 
 в России за 100 лет больших изменений не претерпел. 
 \begin{quotation}
  Обращаясь к недостаткам работы Н.~Н.~Лузина, замечу, что, конечно,
отдельные случайные погрешности, не имеющие влияния на результаты,
можно указать, как во всяком другом труде, и я, например, со своей стороны,
могу отметить определение функции с ограниченным изменением для
совершенного множества....
 \end{quotation}

 Концовка отзыва Егорова \cite{Egor-lus} была необычна:
\begin{quotation}
 Высокие достоинства труда Н.~Н.~Лузина не оставляют во мне никакого
сомнения в том, что было бы только справедливо оценить его дарованием
высшей учёной степени. Прибавлю к этому, что у Н.~Н.~Лузина есть ряд
и других ценных работ и что его имя пользуется почётной известностью
в математическом мире. Замечу ещё, что разбираемый труд в его настоящем
виде содержит столько ценного материала, что его свободно хватило бы на
два отдельных сочинения, особенно при дальнейшем развитии некоторых
указаний автора, которые часто имеют характер простого намёка.

Ввиду всего вышеизложенного я полагал бы безусловно отвечающим
обстоятельствам дела, допустив Н.~Н.~Лузина до защиты настоящего сочинения,
ходатайствовать, в случае удовлетворительной защиты, перед Советом
Императорского Московского Университета об утверждении Н.~Н.~Лузина
в степени доктора чистой математики.

13 марта 1916 г. \hfill Орд. проф. Д.~Егоров
\end{quotation}

Лузин в самом деле получил высшую ученую степень. 
Но главное было не в этом. Его книга оказалась программным сочинением и заворожила московский математический мир.

 \section{Школа Лузина: рождение и цепная реакция%
 \label{s:school}}
 
 \COUNTERS
 
 {\bf\punct Школа Лузина.%
 \label{ss:school}}
 В числе участников группы Лузина 1915-1935 прежде всего называют
 его друзей:
 
 \sm
 
 Голубев Владимир Васильевич (1884-1954)
 
 Привалов Иван Иванович (1891-1941)
 
 Серпинский Вацлав [Wac\l aw Franciszek Sierpi\'nski] (1882-1969)
 
 Степанов Вячеслав Васильевич (1889-1950)
 
 \sm
 
\noindent
 и учеников:
 
 \sm
 
   Александров Павел Сергеевич (1896-1982)
   
   Бари Нина Карловна (1901-1961)
    
  Гливенко Валерий Иванович (1897-1940)
  
  Келдыш Людмила Всеволодовна (1904-1976)
  
   Колмогоров Андрей Николаевич (1903-1987).

   Лаврентьев Михаил Алексеевич (1900-1980)
   
    Люстерник Лазарь Аронович (1899-1981). 
    
    Ляпунов  Алексей Андреевич (1911-1973)
    
    Меньшов Дмитрий Евгеньевич (1892-1988)
    
    Новиков Петр Сергеевич (1901-1975)
    
    Суслин Михаил Иванович (1894-1919)
    
    Урысон Павел Самуилович (1898-1924)
  
  Хинчин Александр Яковлевич (1894—1959)
  
  Шнирельман Лев Генрихович (1905-1938)
  
  \sm
  
  Эти люди появлялись в числе учеников Лузина не одновременно, начиная с первой порции 1914-1916 годов
  --
  Александров, Меньшов, Суслин, Хинчин, -- кончая поздними учениками, Новиковым и Келдыш (1923-1924)
  и Ляпуновым (1933).
  
  Что касается Серпинского, то он был германский подданный, работавший в Львовском университете
  (тогда Львов входил в Австро-Венгрию),
  оказавшийся в августе 1914 в Белоруссии. Он был выслан (интернирован) в Вятку. Благодаря усилиям Егорова и
  Млодзеевского 
  ему разрешили жить в Москве%
  \footnote{Серпинский: {\it <<Но первые две мои работы по теории множеств, имеющие 
некоторое значение, были напечатаны только в 1916 г.>>}, \cite{Sierpinski2}.}. В 1917 году он уехал в Польшу, в дальнейшем  
  продолжал сотрудничать с Лузиным. В 1920  году вместе с С.~Мазуркевичем
  (Mazurkiewicz, Stefan) и 3.~Янишевским (Janiszewski, Zygmunt) основал журнал Fundamenta Mathematica,
  где  Лузин и его ученики опубликовали много работ.

  Понятно, что группа Лузина не могла состоять только из будущих знаменитостей. В рассказах участников
  упоминаются разные другие лица. Упомянем наиболее известных из них, а также просто  тех,  у кого удается найти математические работы%
\footnote{\label{fo:maximov}%
	В описаниях Лузитании никогда не упоминается Максимов Исай Максимович (1889--1976). 
Он, однако, присутствует в библиографическом справочнике \cite{za30},
его упоминал также сам Лузин, по-видимому, он же представлял статьи Максимова в Доклады (я не проверял).
Максимов закончил Казанскую духовную семинарию, работал священником, учителем, преподавателем
педучилища. В 1928 поступил в аспирантуру МГУ, судя по всему, к Лузину. С 1930
работал доцентом в Чувашском пединституте. Первую журнальную статью опубликовал 
в 1935 году (в 46 лет!). Публиковался в Сборнике, Compt. Rendus, Compositio, Bull. AMS, Annals,
Annali di Pisa, Tohoku. Тематика статей соответствует поздней Лузитании. Загадочная картина рассыпается, когда выясняется, что его статьи 
полны ошибок и несообразностей, см. рефераты MR0186765, MR0176935,  MR0186560, MR0000656, MR0001814, MR0001813
в MathReview
(чтобы не было сомнений в грамотности рецензентов, замечу, что  среди них J.~W.~Tukey и E.~Mendelson).}%
$^,$%
\footnote{Из Александрова, \cite{Alex-auto1}:
<<{\it В Лузитании мужской и женский пол были представлены в количественном
отношении примерно поровну. В этом коллективе выдвинулось одно
действительно яркое женское дарование: Нина Карловна Бари. Несколько
лет позже, при новом значительном пополнении состава учеников Н.~Н.~Лузина
прибавилась еще Людмила Всеволодовна Келдыш.}>> 
\newline
Еще отмечу, что из упомянутых 
в следующем списке лузитанок две в дальнейшем
стали женами известных ученых.}
  (в списке могут быть ошибки в датах)
  
 Ф\"едоров   Владимир Сем\"енович (1893—1983) 
 
Селивановский [Линде-Селивановский] Евгений Августович%
\footnote{Автор работ  C. R. 184, 1311-1313 (1927); Матем. сб., 35:3-4 (1928),  379–413; Fundam. Math. 21, 20-28 (1933).
Новиков, \cite{LL} сообщает, что он заболел тифом и в состоянии тифозной горячки выбросился из окна.}

Немыцкий Виктор Владимирович (1900-1967) 

Вениаминов Владимир Николаевич (1895-1932)

Ковнер Сем\"ен Самсонович (1896—1962)

Селиверстов Глеб Александрович  (1905–1944)

 Хлодовский Игорь Николаевич (1903—1951),

Селиванов Н.А.

Богомолова Вера С.

Рожанская Юлия Антоновна ( 1901 - 1967) 
 	
Лихтенбаум Люциан Михайлович (1900?-1968?)

Леонтович Евгения Александровна (1905-1997).

Лисенков Н.~М.

До нас дошло довольно много воспоминаний и лузинской школе и о Лузине того времени.
Ниже мы их обсудим, но начнем все же с математической части истории.

\sm

{\bf\punct Годы 1916-1917.%
\label{ss:1916}}
В 1916--1917 годах в парижских Comptes Rendus появляются заметки четырех молодых людей из Москвы,
Александрова, Хинчина, Меньшова и Суслина.

\sm

a) {\sc Александров.} Как мы все помним, в теории меры Лебега возникает вопрос об описании измеримых множеств.
Напомним, что класс измеримых множеств должен образовывать сигма-алгебру, т.е. быть замкнутым относительно
счетных объединений, счетных пересечений, и относительно дополнений.
Есть два стандартных подхода,  мы можем взять лебеговскую или борелевскую
 сигма-алгебру. По формальному определению, борелевская сигма-алгебра на прямой -- 
наименьшая сигма-алгебра, 
содержащая все интервалы. Взяв пересечение всех сигма-алгебр, содержащих все интервалы,
мы получим борелевскую алгебру. Но это уж слишком неконструктивно%
\footnote{Это доказательство в рамках системы аксиом Цермело-Френкеля, оно использует множество подмножеств множества
подмножеств континуума, что вообще пугает и удивляет, потому что сама борелевская сигма-алгебра континуальна.
И вообще в обычном анализе множества больше континуума не появляются, а здесь огромное множество появилось и исчезло.}. 

Можно рассуждать так. 
Рассмотрим все  множества, получаемые счетными объединениями интервалов.
Потом берем счетные пересечения множеств, полученных таким способом.
Полученный класс называется $G_\delta$. Мы берем счетные пересечения полученных
множеств. Так делаем  до бесконечности и получаем новый класс множеств, который все равно не подходит.
Мы снова берем объединения, потом пересечения и т.д. Как показал Лебег \cite{Lebeg-nomme}, процесс этот трансфинитен,
и успеха мы достигнем, лишь перечислив все счетные трансфинитные числа, а этажи конструкции будут нумероваться первым несчетным трансфинитом,
$\aleph_1$.
Впрочем, это умозрительно,
потому что возможность строить явно представителей получаемых классов быстро становится проблематичной.

В общем, выглядит муторно, и не во всех отношениях эстетично.
Определение измеримости по Лебегу выглядит проще, часто без борелевской сигма-алгебры можно обойтись, но 
- увы! - это не всегда возможно. К примеру, она участвует в теореме Рисса--Маркова.
Да и ситуации, когда нужно рассматривать много мер на одном пространстве вполне обычна
(а тогда от борелевской сигма-алгебры не уйти, потому что лебеговское пополнение по каждой
мере будет свое).

\sm

Вернемся в 1916 год. Александров решает стоявшую тогда (и поставленную перед ним Лузиным) проблему о возможной мощности борелевских множеств.
А именно показывает, что несчетное борелевское множество континуально. Это получается
как следствие другого утверждения: несчетное борелевское множество содержит совершенное подмножество.
Но наиболее замечателен был подход Александрова к доказательству. Он вводит так называемую
$A$-операцию. А именно, рассматривается набор замкнутых множеств $X_{k_1,\dots, k_n}$, 
индексируемых конечными наборами натуральных чисел. А далее берется
$$
\bigcup_{k_1, k_2, k_3, \dots} \bigcap_{n} X_{k_1,\dots,k_n}
$$
Оказывается, что любое борелевское множество может быть получено как результат применения ОДНОЙ  $A$-операции
(и, тем самым, борелевские множества могут быть получены без трансфинитной индукции). И
оказывается, что множество, полученное применением $A$-операции содержит совершенное подмножество
(что решало проблему о мощности).

По-видимому, это соответствовало чему-то, витавшему тогда в воздухе, потому что в том же 1916 году,
близкие утверждения были получены Хаусдорфом (Felix Hausdorff).

\sm

2) {\sc Меньшов.} Он в 1916 году опубликовал пример ненулевого тригонометрического ряда,
который почти всюду сходится к нулю. Результат был сенсационным, и Меньшов на том 
стал знаменит. Но более интересно то, что конструкция Меньшова по сути контрпримером 
не была. А именно, он рассматривал сингулярную меру $\mu$, сосредоточенную на некотором канторовском множестве.
У этой меры, как мы сейчас знаем, есть коэффициенты Фурье
$$
c_n=\int_0^{2\pi} e^{in\phi} \,d\mu(\phi).
$$
Соответствующий ряд Фурье $\sum c_n e^{in\phi}$, как мы сейчас знаем, сходится в смысле обобщенных функций
к самой мере $\mu$. Если  взять стандартное канторовское множество, то ряд окажется всюду расходящимся
в обычном смысле слова. Но изменив процедуру выбрасывания серединок, можно получить желаемый результат
(это связано с диссертацией Лузина и с замечанием выше о том, что производная у Лузина - не
совсем настоящая).

Зигмунд (Antoni Zygmund) 60 лет спустя писал:
\begin{quotation}
 ... [эта работа] может рассматриваться как начальная точка современной теории тригонометрических рядов...
\end{quotation}

По-видимому, разложения сингулярных мер в ряды тогда тоже витали в воздухе.
Скажем, произведения Рисса (Frederic Riesz) были введены  в 1918г., \cite{Riesz}.

\sm

3. {\sc Хинчин.} Заметка Хинчина 1916 года  в наибольшей степени соответствовала тогдашнему mainstream,
но выглядит несколько экзотично с точки зрения сегодняшнего дня. 
После того, как Лебег 
изобрел свой интеграл, у многих авторов последовало желание придумать интеграл еще более общий
(желание обусловленное не только тем, чтобы <<превзойти Лебега>>, 
но и тем, что несобственный интеграл Римана интегралом Лебега не обобщается и тем, что хотелось бы иметь
операцию взятия первообразной, обратную к дифференцированию).
Было изобретено несколько интегралов, Данжуа, Перрона (Oskar Perron) и два интеграла Бореля.
Уже в двадцатые годы было обнаружено, что интегралы Перрона и Данжуа эквивалентны,
а в диссертации Лузина было показано, что два интеграла Бореля не дают ничего нового.
Определение интеграла Данжуа утомительно. Хинчин его обобщает \cite{Hinchin-start}, но определение все равно утомительно.
В том же 1916 году Данжуа публикует эквивалентное определение того же интеграла.
Поэтому то, что получилось, называется <<интегралом Данжуа--Хинчина>> или 
<<широким интегралом Данжуа>>.

Если читать эту работу Хинчина
 \cite{Hinchin-start}--\cite{Hinchin-start2}, то там обнаруживается много изящных деталей
 и возникает желание эту работу продолжать.
Люди и сейчас  занимаются подобными вещами (и есть серьезные упрощения, Jaroslav  Kurzweil),
но все это как-то плохо инкорпорировано в нынешнюю математику. 

Впрочем, Хинчин тогда же вводит понятие {\it асимптотической производной}, которое иногда
используется в анализе, а также понятие {\it абсолютной непрерывности}, которое  вошло в учебники анализа.

\sm

4. {\sc Суслин.} В начале XX века было впечатление, что хоть множеств в теории множеств много,
но все множества, которые могут появиться в анализе -- борелевские. Потому что открытые множества -
борелевские, а что ни конструируй из 
 борелевских множеств, всегда   будет получаться борелевское множество.  В частности, была теорема Лебега
о том, что проекция борелевского множества -- борелевское множество. Суслин в 1916 году
нашел ошибку в доказательстве Лебега%
\footnote{Лебег утверждал, что для убывающей последовательности множеств проекция пересечения есть пересечение проекций. 
	Очевидно, Лебег написал это сгоряча и не подумавши. Со всеми бывает.}...
 А потом показал, что неверно и само утверждение.

В заметке, опубликованной в январе 1917 года,
Суслин вводит $A$-множества, как множества, которые могут быть получены $A$-операцией.
Он показывает, что это то же самое, что проекции борелевских множеств, а также то же самое,
что проекции $G_\delta$-множеств. А также то, что этот класс больше, чем класс борелевских множеств%
\footnote{Вместе с заметкой Суслина
Лузин посылает в Comptes Rendus свою собственную заметку \cite{Luzin-1917} с доказательством измеримости 
аналитических множеств.}.

Лузин потом называл такие множества <<аналитическими>>, видимо, имея в виду, что 
это множества, которые могут по делу появиться в анализе. Но слово <<аналитический>>
имеет в математике другой смысл, термин не очень прижился, а такие множества этого класса чаще называют
<<суслинскими>>. 

\sm

Так или иначе, в 1915--16гг. Москва внезапно стала центром изучения теории функций. Пока только теории функций.

\sm

{\bf\punct Замечания и свидетельства.%
\label{ss:1916-svid}}
 Хинчин поступил в Университет в сентябре 1911 года, Меньшов - 1912 года, а Александров и Суслин
в сентябре 1913 года. Весной 1914 года Лузин появляется в Москве.
К январю 1917ого года все упомянутые заметки были 
уже опубликованы в Comptes Rendus в Париже (который, кстати, был отделен от Москвы фронтами Мировой войны).

Здесь интересна не только быстрота с которой делались отнюдь не простые работы, но и скорость, 
с которой молодые люди входили в тему. В этом была безусловная заслуга Лузина, но, отчасти,
это было и свойство темы. ТФДП того времени требовало острого ума, но все же не требовало долгого обучения.

Очень быстро продвижения, полученные одними участниками группы вводятся в оборот другими.
В диссертации Лузина уже обсуждается введенная Хинчиным асимптотическая производная
и упоминается результат Меньшова об интеграле Данжуа%
\footnote{Из диссертации Лузина:
	\newline
{\it 	<<Перейдемъ теперь къ р\ешен\ию поставленной выше (§ 58) задачи Б. К. Млодз\еевскаго
	относительно общего опред\елен\ия производной. Мы даемъ
	это р\ешен\ие, сл\едуя мысли А.~Я.~Хинчина>> Далее сноска
	<<Предложенной имъ въ сообщен\ии студ. Математическому Кружку.
	6 ноября 1914.>> 
	\newline
	 <<Во время просмотра корректуръ этой работы, Д. Е. Меньшовымъ
	былъ указанъ на зас\едан\ии Московскаго Математическаго Общества
	(1914 г. декабрь), прим\еръ функц\ии, интегрируемой методомъ Dirichlet
	и, однако, неинтегрируемой въ смысл\е ($B''$).>>}
	\newline
	 Осень 1914,
	студенты Меньшов и Хинчин уже серьезно работают в лузинской тематике }%
. Работа Суслина  существенно опиралась
на изобретение Александрова.

Приведем две цитаты. Первая из Александрова \cite{Alex-nauka-i-zhizn'}:
\begin{quotation}
 Я впервые встретился с ним [Лузиным] будучи студентом 2-го
курса. Впечатление от этой встречи было, можно прямо
это сказать, потрясающим, и я запомнил его на всю жизнь.
Обратившись к нему после окончания лекции за советом,
как мне заниматься математикой дальше, я был прежде
всего поражен внимательностью и — не могу найти 
другого слова — уважением к собеседнику — как ни странно
звучит это слово, когда речь идет о беседе уже 
знаменитого, хотя и молодого еще, ученого с 18-летним 
студентом. Выслушав меня, Лузин посредством умело 
поставленных вопросов очень скоро разобрался в характере
моих математических склонностей и сразу же в доступной
мне форме обрисовал основные направления, которые он
мог мне предложить для дальнейших занятий; очень 
осторожно он сам меня склонил к выбору одного из этих
направлений, причем все это было сделано тонко, без
всякого нажима и — как я теперь могу сказать — 
правильно. Я стал тогда же учеником Лузина, и это было
в эпоху его наивысшего творческого подъема.

Лузин жил тогда совершенно один в меблированных
комнатах, жил только наукой. Мне запомнилась одна его
фраза, сказанная в одну из многочисленных наших
встреч: «Я дни и ночи думаю над аксиомой Цермело%
\footnote{Аксиома выбора. Важная героиня математики того времени и нашего повествования тоже.
Позже Лузин стал  относился к ней настороженно.}
(такая
есть в математике знаменитая аксиома, которая была
тогда — и еще много лет спустя — в центре исследований
по логическим основаниям математики). Если бы только
кто-нибудь знал, что это за вещь!».

Мое знакомство с Лузиным пришлось довольно точно
в середину того десятилетия, в котором он получил все
самые значительные свои результаты.
Видя Лузина в эти годы, я видел действительно то, что
может называться вдохновенным отношением к науке,
и я не только учился у него математике, я получил и урок
того, что такое настоящий ученый, а также и того, чем
может и должен быть профессор университета.
\end{quotation}

%\begin{quotation}
% Естественно, возникла обратная задача: всякое ли 
%множество, получаемое применением А -операции к 
%замкнутым множествам, является борелевским? Лузин правильно
%увидел в этой задаче центральную проблему 
%дескриптивной теории множеств того времени и с большой 
%настойчивостью ставил ее П.~С.~Александрову и М.~Я.~Суслину
%(1894—1919), тогда студенту 3-го курса.
%\end{quotation}

%\begin{quotation}
% «Но здесь тоже надо мужество — сказать
%что я занимался не только как переписчик Вашей ноты и не только
%тем, чтобы сокращать и вдавливать коленом этот материал, но мною
%было сделано очень много, и эта таблица была сделана мною 
%одним».).
%\end{quotation}

Вторая - из  популярной книги Серпинского \cite{Serp} <<О теории множеств>>, 1964:

\begin{quotation}
 В 1916 г. Николай Лузин, в то время молодой профессор 
Московского университета, предложил своему ученику Михаилу Суслину изучить
работу А.~Лебега об аналитически представимых функциях, 
опубликованную в 1905 г. в известном французском журнал е Journal des Mathematigues.
В этой работе Лебег приводит, между прочим, теорему о том, что если
аналитически представимая функция обратима, то обратная ей функция
тоже является аналитически представимой. Доказательство этой теоремы
Лебег основывает на нескольких леммах и в их числе на лемме, согласно
которой проекция на прямую общей части бесконечной стягивающейся
последовательности множеств, расположенных на плоскости, является общей
частью проекций этих множеств. ........
Лебег не привел доказательства этой леммы, считая
ее очевидной, и того же мнения придерживались, видимо, многочисленные
читатели работы Лебега в течение 10 лет. Так вот, Суслин, продумывая эту лемму, 
установил, что она неверна....

Мне довелось быть свидетелем того, как Суслин сообщил Лузину
свое замечание и вручил ему рукопись своей первой работы. Лузин
очень серьезно отнесся к сообщению молодого студента, и подтвердил,
что тот действительно нашел ошибку в труде известного ученого. Я также
читал рукопись Суслина непосредственно после Лузина и знаю, как Лузин
помогал своему ученику и как направлял его работу. Некоторые авторы 
называют аналитические множества множествами суслинскими; правильнее
было бы называть их множествами Суслина—Лузина.
\end{quotation}

\sm

{\bf \punct Поворот Лузина.%
\label{ss:povorot}} Работы Лузина 1911-1914гг. и его диссертация 1915г.
 были посвящены теории функций действительной переменной.
С 1917 года статей по ТФДП у него больше не было (или почти не было%
\footnote{Была статья Sur un mode de convergence de l'int\'egrale de Dirichlet в Трудах Физ-матем. ин-та Казанского университета, 1934,т.6, сер. 3, стр. 1-4 и была неоконченная рукопись
	<<Об одном особом интеграле>> середины 20х, опубликованная
в \cite{Luz-trig}.}), правда он написал несколько работ по теории функций
комплексного переменного (они собраны в конце первого тома <<Собрания сочинений>>), среди них очень 
известная работа с Приваловым о граничных значениях аналитических функций \cite{Luzin-Privalov} (Александров \cite{Alexandrov-1937} писал, что работы была сделана в 1919г).

Теперь 
его основной темой становится дескриптивная теория множеств, т.е. исследование множеств, которые могут
естественным образом возникнуть путем применения естественных операций к простым множествам
(по смыслу, это множества, которые могут возникнуть сами собой в анализе).
Начала подобной теории существовали к 1917 году
в виде исследований французской школы по борелевским множествам и классам Бэра,
а по большей части эту науку еще предстояло создать.
Занимались этим в 20х прежде всего Лузин и члены его группы%
\footnote{Конечно, в этом направлении работали и другие математики, упомянем уехавшего из этой группы в Варшаву
Серпинского, а также Кондо (Motokiti Kond\^o).}, фактически эта область тогда отделилась от ТФДП%

По-видимому, решающей причиной  поворота Лузина был успех, достигнутый Александровым и Суслиным.
В силу результатов Суслина (см. выше),   определение аналитических множеств в итоге оказывалось более простым
и конструктивным, чем изначальное определение борелевских множеств.

Лузин находит несколько описаний других описаний множеств этого  класса,
например, как  непрерывные образы пространства $\R\setminus \Q$.
Интересно, что непрерывные взаимно однозначные 
(с точностью до счетных множеств) образы этого множества являются борелевскими множествами.
А другое определение - множество является борелевским тогда и только тогда, когда
оно аналитично вместе со своим дополнением.

Лузин увлекается этой темой, придумывает разные интересные вопросы, решает их с учениками,
и это составит основную часть написанной Лузиным монографии \cite{Lus-anal}, 1930 года.
В наши цели не входит обсуждение этих результатов, Лузин 
излагал себя сам, были   обзоры этой науки в Успехах 1950года,
 \cite{Schegol1950}, \cite{Arsenin}, а также статьи \cite{KeldNov}, \cite{Keld-Luzin},
 \cite{Usp}, \cite{Kanov}. Мы же пока оставим чисто научную часть биографии Лузина до следующей ее критической точки
 (1924-25гг).

\sm

Важно, что  все, сделанное для пространств $\R$ или $\R^n$,
фактически имеет значительно б\'ольшую общность - а именно переносится на польские пространства%
\footnote{Польское пространство -  полное сепарабельное метрическое
пространство. Метрика обычно является в той или иной степени искусственной, поэтому чаще  говорят 
о топологических пространствах, гомеоморфных полному сепарабельному топологическому пространству.
Это большой естественный класс пространств, например, $G_\delta$-множество
(счетное пересечение открытых множеств) в польском пространстве само является польским пространством
(Александров, \cite{Alexandrov-G-delta}).
\newline
В частности, польским пространством является
и упомянутое чуть выше $\R\setminus \Q$. Оно гомеоморфно <<пространству Бэра>>
$\N^\N$, т.е. пространству всех последовательностей $(l_1,l_2,  \dots)$, состоящих из натуральных чисел.
Возьмем две  последовательности $(l_1,l_2,  \dots)$, $(m_1,m_2,  \dots)$. Пусть $k$ -- первый номер,
в котором эти последовательности различаются. Положим расстояние между ними равным $1/k$.
Мы получаем полное вполне несвязное метрическое пространство. Отображение из множества иррациональных
точек отрезка $[0,1]$ в $\N^\N$
-- это  разложение числа в цепную дробь (понятно, что $\R\setminus \Q$ гомеоморфно $(\R\setminus \Q)\cap [0,1]$).
Рациональные числа выпадут, потому что у них это разложение 
конечно.
\newline
Лузин часто работал именно с пространством Бэра, то есть с довольно общим
вполне несвязным польским пространством. Но сам он, скорее всего, в 20е годы этого не понимал,
и в начале книги \cite{Lus-anal} объясняет резкое изменение дескриптивных возможностей при переходе от $\R$ к $\R\setminus\Q$
(обнаруженное им, видимо, в результате немалых усилий)
разной природой рациональных и иррациональных чисел, и тем, что лучше рассматривать 
отдельно иррациональные числа как <<однородную>> часть континуума. Странность этого взгляда была отмечена
в комментариях Новикова и Келдыш к русскому изданию \cite{Lus-anal}
(если мы выкинем из $\R$ произвольное счетное плотное множество, то выйдет то же самое).
По-видимому, здесь на  собственных математических изысканиях Лузина отразился его
личный конфликт с Александровым
и нежелание понимать топологию по Александрову.
\newline
Упоминаемый ниже в п.\ref{ss:name} пример аналитического неборелевского множества
из статьи \cite{Lus-baire} <<Об арифметическом примере...>>, основанный на цепных
дробях, фактически никакого отношения к теории чисел (<<арифметике>>) в обычном ее понимании не имеет.
\label{fo:baire}%
}.
В частности, все (с тривиальными оговорками) польские пространства
имеют изоморфные 
борелевские структуры. См. <<Топологию>> Куратовского \cite{Kuratwski}.

Дескриптивная теория множеств в дальнейшем  в неформальном смысле раздвоилась. С одной стороны 
есть чисто логическая часть, где доказывается недоказуемость (см. \cite{Kanov0}, \cite{Kanov}). Но 
осталась на удивление большая общематематическая часть с позитивными результатами, 
см. \cite{Kechris}, см. также \cite{Gao}, причем в последние два десятилетия здесь
произошло определенное оживление.

Стоит отметить, что аналитические множества появляются во вполне реальных задачах.
Рассмотрим группу $\mathrm{Ams}(M)$ преобразований пространства $M$ с вероятностной
 мерой $\mu$, оставляющих меру инвариантной.
Группа $\mathrm{Ams}(M)$ имеет естественную сепарабельную метризуемую топологию,
$g_j\in \mathrm{Ams}(M)$ сходится к  $g$, если для любых измеримых
$A$, $B\subset M$ имеет место сходимость $\mu\bigl(g_j(A)\cap B\bigr)\to \mu\bigl(g(A)\cap B\bigr)$
(так, что получается польское метрическое пространство).
Огромная литература (начиная с 30х годов) посвящена изучению классов
сопряженности в $\mathrm{Ams}(M)$. Оказывается, что отношение сопряженности
в $\mathrm{Ams}(M)$ не является борелевским.
 А именно, рассмотрим
в $\mathrm{Ams}(M)\times\mathrm{Ams}(M)$ подмножество $\Xi$, состоящее из
пар $(g,h)$, таких, что $g$ и $h$ сопряжены. Оказывается,
что это подмножество является аналитическим, но не борелевским;
более того,  утверждение остается в силе, если ограничиться
лишь эргодическими  $g$, $h$, см. \cite{Hjo}, \cite{FRW}. Это,
в частности показывает безнадежность вопроса о классификации
эргодических преобразований с точностью до сопряжения
(что, впрочем и до того едва ли могло вызывать сомнения).

\sm

{\bf \punct Математика в Иванове-Вознесенке.%
\label{ss:ivanovo}} Так или иначе, пришла Гражданская война,
жизнь в Москве стала очень тяжелой,  и математики
стали разъезжаться в разные провинциальные центры, где жизнь, видимо, была полегче, а также
(удивительным образом)
появились рабочие места. 

Голубев и Привалов уезжают в Саратовский университет, где в 1917-1918 открылся физико-математический факультет,
Голубев там становится заведующим кафедрой, потом деканом, а потом  -- ректором (1921-1923).

Александров в своих воспоминаниях рассказывал, что Лузин ему предложил в качестве задачи континуум-гипотезу,
он в связи с этим столкнулся с тяжелым кризисом и далее вел жизнь театрального деятеля
в Новгород-Северске и Чернигове. Когда в Чернигов вошли деникинцы он, как советский работник, был арестован,
но нашлись люди, которые за него хлопотали, и он был освобожден. Потом он  читал лекции в Черниговском
пединституте по литературе и математике.

В 1915 году в Иваново-Вознесенск был эвакуирован Рижский политехнический институт. 
В 1918 году на его базе создается Ивановский политехнический институт и Ивановский педагогический
институт.
Осенью 1918 года туда едет Лузин, а за ним несколько его учеников, Меньшов, Суслин, Хинчин, Фёдоров%
\footnote{Меньшов \cite{LL} называет много других имен московских математиков, работавших
тогда в Иванове.}.
Лузин одновременно сохраняет должность профессора в Москве и читает там лекции
(Иваново недалеко от Москвы, но передвижение по железным тогда было связано с определенными сложностями,
рассказывается, что Ивановский институт имел свой вагон, и связь с Москвой была приемлимой).

Хинчин на некоторое время уехал из Иванова и был деканом в Нижнем Новгороде, потом вернулся и стал
деканом в Ивановском педагогическом институте. 

Суслин после конфликта с Лузиным
уехал к себе в деревню под  Саратовым и умер  от тифа в конце 1919 года.

Егоров и Степанов оставались в Москве.

В конце 20 года по окончании Гражданской войны   члены лузинской группы стали возвращаться,
впрочем, вернулись не все.
 Голубев еще долго работал в Саратове, Фёдоров так и остался  в Иванове, а Хинчин  задержался там до  1922 года.

\sm

{\bf\punct Поворот Хинчина.%
\label{ss:povorot-H}} Мы помним, что 1920-40годы были героической эпохой теории вероятностей,
связанной с именами Колмогорова, Хинчина и Леви (Paul L\'evy). А как это начиналось в Москве?
Наша более широкая цель в пп.\ref{ss:povorot-H}--\ref{ss:tfdp} -- посмотреть как начиналась математическая <<цепная реакция>>,
о которой писал Лаврентьев.

Как мы видели, Хинчин задержался в Иванове до 1922 года.
Список его трудов  содержит такие строчки:

\sm

{\small
Об одном свойстве непрерывных дробей и его арифметических приложениях, Известия
Иваново-Вознесенского политехнического ин-та, 1922, 5, 27--41.

\sm

К вопросу о представлении числа в виде суммы двух простых чисел, Известия Иваново-Вознесенского политехнического  ин-та, 1922, 5, 42--48.

\sm

Об одном вопросе теории диофантовых приближений. Изв. Иваново-Вознесенского
политехнического ин-та, 1925, 8, 32--37.
}

\sm

Не знаю, многие ли библиотеки мира могут похвастаться наличием этого журнала, но, так или иначе, мы видим, что к 1922 году
Хинчин занялся теорией чисел (и, кстати, темы всех трех статей - разные). Тем же годом подписана следующая его статья по 
по скорости сходимости цепных дробей \cite{Khinchin-chain}. Приведем два более поздних результата 
Хинчина на эту тему
\cite{Khinchin-copositio1}--\cite{Khinchin-copositio2}. Разложим 
 число $\alpha$ в непрерывную дробь, пусть $a_1$, $a_2$, \dots -- неполные частные, пусть
$q_1$, $q_2$, \dots -- знаменатели подходящих дробей. Тогда следующие пределы существуют почти всюду
и равны (п.в.) постоянным величинам
$$
\lim_{n\to\infty}\sqrt[n]{a_1\dots a_n}=C, \qquad \lim_{n\to\infty}\sqrt[n]{q_1\dots q_n}=D,
$$
При этом
$$
C=\prod_{j=1}^\infty \Bigl[\frac{(j+1)^2}{j(j+2)} \Bigr]^{\log_2 j},
$$
а
вторая величина была впоследствии вычислена Леви,
$$
 \qquad D=e^{\frac{\pi^2}{12\ln 2}}
.
$$

А вот пример результата Хинчина \cite{Khinchin-Diofant}, 1926, по диофантовым
приближениям.
Пусть $a$ пробегает иррациональные числа.
Пусть $\phi(t)$  -- положительная функция на полуоси $t>0$, 
 такая, что $t^2\phi(t)$ монотонно убывает.
 Пусть интеграл $\int_1^\infty t\phi(t)\,dt$ расходится. Тогда для почти всех 
  $a$ существует бесконечно много рациональных чисел $p/q$, 
  таких, что 
  $$
  \Bigl|a-\frac pq\Bigr|<\phi(q).
  $$
 Если же интеграл сходится, то для почти всех $a$
 есть лишь конечное количество таких рациональных чисел.
 
 Так или иначе уже в 1922 году Хинчин начал внедрять методы теории функций действительной переменной в теорию чисел.
 Приведем две цитаты. Одна -- из Люстерника \cite{Lyu-3}, 
 
 \begin{quotation}
  В 1922 г. А.~Я.~Хинчин переехал в Москву (но еще 4 года наезжал
в Иваново). Он объявил спецкурс по теории чисел. Постоянным слушателем
этих лекций был Л.~Г.~Шнирельман.
 \end{quotation}
 
 Вторая цитата -- из Александрова \cite{Alex-1955}:
 \begin{quotation}
 В этом смысле огромное
значение имел семинар по аналитической теории чисел, который А.~Я.~Хинчин
вел в Московском университете в 1925—1926 гг. В этом семинаре работали
такие, тогда еще совсем молодые математики, как А.~О.~Гельфонд.
(бывший по комплексному переменному учеником И.~И.~Привалова),
Л.~Г.~Шнирельман, Н.~Г.~Чудаков, Н.~П.~Романов и др. Все эти математики
стали, как известно, крупными специалистами по теории чисел, и в этом
семинар А.~Я.~Хинчина, несомненно, сыграл значительную роль.
\end{quotation}

О продолжении теоретико-числовой истории  в начале 30х, см.п. \ref{ss:ShGe}.

Что касается истории теории чисел в СССР, то стоит помнить, что в Ленинграде
в 20х годах была сильная теоретико-числовая школа, включавшая в себя 
Я.~В.~Успенского,  Р.~О.~Кузьмина, И.~М.~Виноградова,
 Б.~Н.~Делоне. В 1934 году Виноградов и Делоне вместе
 Институтом Стеклова  переехали в Москву. Между ленинградской
 и московской школой была явная неприязнь, но я думаю, что в действительности
 стороны внимательно следили друг за другом.

Тогда же, в первой половине двадцатых Хинчин появляется в другой ипостаси.
В 1924 году он публикует \cite{Hinchin-logarifm} знаменитый закон двойного логарифма.
Пусть мы случайным образом бросаем монетку. Пусть $s(n)$ -- число <<гербов>> минус число
<<решек>>, выпавших
за первые $n$ бросаний. Насколько быстро может расти это число?
Это классическая и много исследовавшаяся задача.
Напомним, что ответ Хинчина такой
$$
\limsup_{n\to\infty} \frac{s(n)}{\sqrt{2n\ln\ln n}}=1.
$$
Слабую версию того же утверждения Хинчин опубликовал годом раньше \cite{Khinchin-logarifm1}.

В 1925 году была опубликована знаменитая теорема Хинчина и Колмогорова о суммах независимых случайных величин.
Напомним ее формулировку. Пусть $\sum \xi_j$ -- ряд, составленный из независимых случайных величин. Пусть, для простоты
они имеют нулевые средние. Тогда сходимость ряда из дисперсий влечет сходимость ряда  $\sum \xi_j$ 
почти всюду%
\footnote{Здесь очевидна $L^2$-сходимость, а сходимость почти всюду факт очень тонкий.}. 

%Цитирую Колмогорова \cite{Kolmogor}.
%
%\begin{quotation}
% В 1925 г. я окончил Московский университет как студент и поступил
%в университетскую аспирантуру. Моим руководителем в аспирантуре попрежнему был Н. Н. Лузин.
%(Напомню, что пребывание в аспирантуре не
%завершалось тогда диссертацией, как сейчас: ведь ученые степени были
%введены лишь в 1934 г.) Еще в 1924 г. я начал интересоваться теорией вероятностей.
%Моя первая работа в этой области относится к тому же 1924 г.
%Она была выполнена совместно с А. Я. Хинчиным (также учеником Н. Н. Лузина).
%Все мои занятия по теории вероятностей совместно с А. Я. Хинчиным,
%весь вообще первый период занятий этой теорией отмечен тем, что мы применяли
%методы, разработанные в метрической теории функций. Такие темы,
%как условия для применимости закона больших чисел, условие сходимости
%ряда независимых случайных величин, велись по существу методами, выкованными
%в общей теории тригонометрических рядов, т. е. методами, разрабатывающимися
%Н. Н. Лузиным и его учениками.
%\end{quotation}

Перед тем как оставить рассказ о Хинчине,  процитируем Люстерника
\cite{Lyu-3}:
\begin{quotation}
 Добавлю, Лузитания и пост-Лузитания состояли из большого количества
пересекающихся дружеских дуэтов, трио, квартетов и т. д. А.~Я.~Хинчин
не входил ни в один из них и представлял собой отдельную компоненту.
В процессе распада Лузитании и ее превращения в комплекс направлений
большую роль сыграл отход А.~Я.~Хинчина от традиционной тематики
по теории функций и переход к новым для московской школы темам из теории
чисел и теории вероятностей. Однако вследствие, очевидно, большой изолированности
А.~Я.~Хинчина в Лузитании этот уход не был в эмоциональном
отношении таким драматическим, как откол топологической школы.
\end{quotation}

Стоит отметить, что в научном отношении
поворот Хинчина имел не меньшее значение, чем поворот Александрова.

\sm

{\bf \punct Студент Колмогоров.%
\label{ss:kolmogorov}} В 1920 году Колмогоров поступил на Физмат Московского Университета
(кстати, вступительных экзаменов в тот год не было). Цитируем Ширяева \cite{Shiryaev}:

\begin{quotation}
 Параллельно с университетом Андрей Николаевич поступает на 
математическое отделение Химико-технологического института им. Д.~И.~Менделеева 
(где, как раз, требовался вступительный экзамен по математике).
Сам Андрей Николаевич так объясняет это свое решение: <<{\it Не бросал
мысль о технической карьере, почему-то меня увлекала 
металлургия}>>, <<{\it Техника тогда воспринималась как что-то более
серьезное и необходимое, чем чистая наука}>>.

Продолжает он и свои серьезные занятия историей, участвуя в семинаре
по древнерусской истории известного историка профессора С.В. Бахрушина
на историческом факультете Московского университета в качестве 
вольнослушателя. В этом семинаре в 1920 году Андрей Николаевич делает свой
первый научный доклад%
\footnote{Эта работа была издана уже после смерти Колмогорова, \cite{Kolmogorov-Novgorod}.} — о земельных отношениях в Новгороде на основе
анализа писцовых книг XV—XVI веков, в котором <<{\it использовались... некоторые приемы математической теории}>>
\end{quotation}

С сентября 1921 года Колмогоров занимался в семинаре Степанова по тригонометрическим 
рядам. Цитируем интервью с Колмогоровым \cite{Kolmogor}, взятое В.~А.~Успенским в 1983г.:

\begin{quotation}
Когда и каким образом состоялось Ваше личное знакомство с Н. Н. Лузиным?

— Знакомство состоялось, когда я был студентом второго курса. На
этом курсе я начал заниматься в семинаре В.~В.~Степанова. Работая в этом
семинаре, я решил задачу, которой интересовался Н.~Н.~Лузин. Возможно,
что и сама эта задача была им поставлена. Во всяком случае, именно со ссылкой
на Н.~Н.~Лузина формулировка задачи обсуждалась на семинаре
В.~В.~Степанова. Речь шла о построении ряда Фурье со сколь угодно медленно
стремящимися к нулю коэффициентами. Мне удалось решить эту
задачу (эта была моя первая самостоятельная работа%
\footnote{См. \cite{Kolmogorov-slowly}.}). Когда об этом рассказали
Н.~Н.~Лузину, он обратился ко мне (помню, это было на университетской
лестнице) и предложил регулярно приходить к нему на занятия.
\end{quotation}

 В том же  1922 году Колмогоров построил пример функции из $L^1$, ряд Фурье которой расходится почти всюду,
 \cite{Kolm-divergent}. В самом начале того же года под влиянием Александрова,
 Колмогоров вводит новую операцию над множествами \cite{Kolmogorov-operatsii}. 
 Цитируем воспоминания Колмогорова об Александрове \cite{Kolmogorov-Alexandrov}: 
 \begin{quotation}
 Мне запомнились вытащенные откуда-то Павлом Сергеевичем огромные
листы бумаги со схемами образования множеств все более высоких классов,
созерцание которых в конце концов привело Павла Сергеевича к тому результату,
что все B-множества  любого класса являются A-множествами%
\footnote{Любое борелевское множество может быть получено операцией Александрова.}. Эти
листы раскладывались по полу и Павел Сергеевич вместе со мной ползал по
ним, желая сделать наглядным получение B-множеств высоких (хотя бы
и трансфинитных) порядков в результате однократного применения А-операции...

...написанные в 1921—1922 гг. мои дескриптивные работы пролежали в письменном
столе Н.~Н.~Лузина, находившего их методологически неправильными,
без всякого движения до 1926 г.

Именно Павел Сергеевич добился того, что мои работы по дискриптивной
теории были все же опубликованы. Он же был инициатором моего оставления
на работе после аспирантуры в Институте математики и механики при Московском
университете.
\end{quotation}

 Снова цитируем интервью \cite{Kolmogor}:
 \begin{quotation}
 В 1925 г. я окончил Московский университет как студент и поступил
в университетскую аспирантуру. Моим руководителем в аспирантуре 
по-прежнему был Н.~Н.~Лузин. (Напомню, что пребывание в аспирантуре не
завершалось тогда диссертацией, как сейчас: ведь ученые степени были
введены лишь в 1934 г.) Еще в 1924 г. я начал интересоваться теорией вероятностей.
Моя первая работа в этой области относится к тому же 1924 г.
Она была выполнена совместно с А.~Я.~Хинчиным (также учеником Н.~Н.~Лузина).
\end{quotation}

В 1925 году Колмогоров опубликовал свою первую работу по логике <<О принципе tertium non datur>>,
\cite{Kolmogorov-tertium}. Напомню, что в начале XX века шли бурные споры о так называемых 
<<основаниях математики>>. Среди предлагавшихся проектов реформ был интуиционизм, 
выдвинутый Брауэром (Luitzen Egbertus Jan Brouwer, 1881-1966). В частности
говорилось, что <<закон исключенного третьего>> верен лишь для конечных совокупностей,
и не может применяться в большей общности. Впрочем, слабая версия <<$A$ влечет не не $A$>>
остается. Колмогоров показал, что если высказывание $Z$ выводимо из $A$, $B$, \dots
в рамках обычной логики, то <<не не $Z$>> выводимо из <<не не $A$>>, <<не не  $B$>>, \dots
в рамках правил, разрешенных интуиционистами. 

Интуиционизм
многократно обсуждался у Лузина на семинаре, и интерес Колмогорова, надо думать, произошел оттуда.
Но это уже математическая логика, которой сам Лузин не занимался.

Радикальный отрыв Колмогорова от Лузинской тематики начался с уже упомянутой совместной с Хинчиным 
 работой по теории вероятности 1925.
Следующая статья Колмогорова на эту тему была опубликована  в 1927 году, \cite{Kolmogorov-grand},
а в 1932 году после публикации 60-страничной книжечки 
Grundbegriffe der Wahr\-schein\-lich\-keits\-rechnung, 
\cite{Kolmogorov-osnovy}
он становится живым классиком.

На этом оставим Колмогорова.

\sm

{\bf\punct Александров и Урысон.%
\label{ss:uryson}} Цитируем книгу Гнеденко <<Очерки по истории математики в России>> \cite{Gnedenko-book}, 1948:
\begin{quotation}
В 1921г. Павел Самуилович Урысон начал
заниматься топологией — современной весьма важной 
ветвью геометрии. В следующее году к нему присоединился
П.~С.~Александров. Этим было положено начало советской
топологической школе. Сейчас весьма уместно сказать 
несколько слов об её основателе и одном из самых крупных
топологов нашего века.

Научная жизнь П.~С.~Урысона началась с работ по 
экспериментальной физике, выполненных им в 15-летнем 
возрасте под руководством академика. П.~П.~Лазарева.
В 1918—1919 гг. он выполнил несколько математических
работ и с этих пор не прерывал своих занятий в этой 
отрасли науки. 
\end{quotation}

Отрывок из Люстерника:
\begin{quotation}
Первые топологические
работы П.~С.~Урысона, о которых он рассказывал в конце 1921 г. и начале
1922 г.%
\footnote{
Доклады в Московском математическом обществе:
\newline
П. С. У р ы с о н «О размерности множеств».
\newline
6. П. С. У р ы с о н «Общие теоремы о размерности»
\newline
П. С. У р ы с о н «Индексы ветвления».
}, как и работы А. Я. Хинчина по метрической теории
чисел%
\footnote{А. Я. X и н ч и н «О двоичных дробях».}, воспринимались как «ответвления» внутри школы теории функций.
Но уже весной 1922 г. произошло отпадение от Лузитании «двух
П.~С.» — П.~С.~Александрова и П.~С.~Урысона. Помню их совместное выступление
в маленькой аудитории на 3-м этаже возле лестницы с докладами
о некоторых свойствах плоских континуумов. Присутствовал В.~В.~Степанов
и студенты В.~А.~Ефремович, В.~В.~Немыцкий, Ю.~А.~Рожанская,
Л.~Г.~Шнирельман и др. Мы почувствовали, что присутствовали при рождении новой школы.
\end{quotation}

Последняя фраза может быть эвфемизмом: ничего, противоречившего по духу ТФДП,
в тогдашних работах Александрова и Урысона не было%
\footnote{В связи с этим стоит отметить слова Меньшова из интервью \cite{Menshov-Duvakin}:
\newline
{\it  Вот он [Лузин] в 1916-м – 17-м году читал второй раз спецкурс, по вопросу, далекому от него, — 
о теоретико-множественной топологии. Причем этот курс он читал очень своеобразно. Он к нему не готовился.
У него был один мемуар известного польского математика Янышевского, если я не путаю фамилию, и вот он,
этот мемуар, так сказать, был основой, началом теоретико-множественной топологии. И вот лекция состояла в том, 
что мы вместе разбирали этот мемуар. Лузин смотрел в этот мемуар и начинал доказывать теорему — ну, какую-то – 
ту или другую теорему. Но там было написано очень сжато, Лузин это читал в первый раз или, может,
так, накануне чуть-чуть перелистал, и, значит, поэтому не сразу восстанавливал доказательство. 
И мы все помогали ему восстановить доказательство, и доказательство восстанавливалось общими силами.
Так что мы были как бы соавторами этого спецкурса, как будто мы с ним совместно читали этот спецкурс.}
\newline
По-видимому, речь идет о Janiszewski, S.
{\it Sur les continus irr\'eductibles entre deux points}. 
J. de l’\'Ec. Pol. (2) 16, 79-170 (1912).}. 

Первая топологическая работа Александрова и Урысона была опубликована
в 1923 году в \cite{Alex-Urysohn1}. Она была первой из заметок,
анонсировавших <<Мемуар о компактных топологических пространствах>> \cite{Alex-Urysohn2},
он был написан в 1922, но опубликован лишь в 1929 году.

\sm

Урысон погиб во Франции, в Бретани в 1924 году, купаясь в штормовом море. Необычна судьба 
его последней
работы об универсальном метрическом пространстве, опубликованной Александровым в 1927 году
\cite{Urysohn}.

 Урысон строит полное сепарабельное метрическое
пространство $\mathbb U$, удовлетворяющее следующему свойству:
 Пусть $A$ – конечное метрическое пространство, а $B$ – некоторое метрическое пространство, 
 полученное из $A$ добавлением одной точки. Тогда
любое изометрическое вложение $A\to \mathbb U$ продолжается до изометрического
вложения $B\to \mathbb U$ (тот же объект получится, если вместо всевозможных
конечных пространств брать произвольные компактные).

Пространство, обладающее эти свойством, единственно и обладает разными другими 
удивительными свойствами. Например,
любое сепарабельное метрическое пространство $M$ вкладывается в $\mathbb U$
изометрично. Таких вложений много. Зато любые два изометричных вложения {\it компактного}
метрического
пространства $K$ в $\mathbb U$ переводятся одно в другое изометрией пространства $\mathbb U$
(С.~А.~Богатый, 2000).
В частности, $\mathbb U$ имеет огромную группу изометрий; оказывается, что любая группа
с польской топологией может быть вложена в группу изометрий пространства $\mathbb U$
(В.~В.~Успенский, 1990)
.

% Более того, вложить его можно так, что любая изометрия $M$
%будет продолжаться до изометрии $\mathbb U$.

Урысон строит этот объект явно, в том смысле, что он получается в результате 
выполнения точной процедуры построения (но реально конструкция очень громоздка и использовать ее проблематично,
фактически используется факт существования и свойство добавления точки).
Работа, видимо, не привлекла в свое время внимания. Она вошла в собрание
сочинений Урысона, опубликованное в 1951 году, и снова осталась незамеченной.
По-видимому, она была извлечена на свет божий лишь в 1986 году Катетовым
(Miroslav Kat\v etov) и сейчас является популярным объектом исследований.
Кстати оказалось, что случайное выписывание метрики на счетном пространстве,
приводит к пространству, пополнение которого -- урысоновское (А.М.Вершик, 2004)%
\footnote{Одна из возможных точных формулировок такая: возьмем пространство $\cM$ всех
	метрик на счетном пространстве. Тогда для плотного $G_\delta$-подмножества
	в $\cM$ пополнение  является урысоновским.}.

\sm

Так или иначе, в 1922-1924 Александров и Урысон занимались созданием общей топологии%
\footnote{В \cite{Lyusternik-Shnirelman-2} говорилось
	\newline
\it	Покойный П. С. Урысон читал в 1924 г. в Московском 
	математическом обществе доклад, в котором приводил доказательство 
	существования двух замкнутых геодезических на выпуклой поверхности
	жанра [рода] 0. К сожалению его работы не могли быть опубликованы за
	его смертью.} и в  объеме учебной дисциплины математики
с ней знакомятся. У людей моего поколения и людей несколько более старших
Александров запомнился именно как деятель общей и теоретико-множественной топологии,
когда рассматриваются весьма общие точечные множества,
общие топологические пространства или пространства
мощности больше континуума. Но надо иметь в виду, что   интересы 
Александрова были шире и
имели  серьезные последствия  за пределами этих  тем.

В 1927г. Александров публикует заметку в Сompt. Rendus
\cite{Alexandrov-Cheh0} (см.также \cite{Alexandrov-Cheh}, \cite{Alexandrov-Cheh1}),
в которой вводит понятие {\it нерва покрытия}.
А именно, возьмем конечное открытое покрытие $\cO$ компактного
пространства $X$,
$$
X=\cup \Omega_j.
$$
По этим данным строится симплициальный клеточный комплекс $C(X,\cO)$.
Его вершины $\omega_j$ нумеруются элементами покрытия.
Если пересечение множеств $\Omega_{i_1}$, \dots, $\Omega_{i_k}$
не пусто, то мы рисуем симплекс 
с вершинами $\omega_{i_1}$, \dots, $\omega_{i_k}$.

Если всевозможные  пересечения множеств $\Omega_j$
связны и односвязны, то комплекс $C(X,\cO)$
гомотопически эквивалентен начальному пространству
$X$, и мы получаем возможность считать (и определять)
гомологии в через покрытия пространства.
Для сложных пространств мы можем заниматься измельчением
покрытий и рассмотрением пределов соответствующих групп.

Это то, что называется {\it когомологиями Чеха}.
  Чех (Eduard \v Cheh) в 1932г. использовал работу Александрова,
 несколько ее обобщив. Термин <<когомологии Александрова--Чеха>>
 был бы более справедливым. Но, кажется, он используется лишь в русской литературе%
 \footnote{Отметим, что это
 результат эволюции математической терминологии, а  не взаимоотношений
 между двумя авторами.}.

\sm

В 20е годы Александров  начинает создание свой
научной школы. Среди его учеников этого времени известны 
А.~Н.~Тихонов, В.~В.~Немыцкий, А.~Г.~Курош, В.~А.~Ефремович, Н.~Б.~Веденисов,
 Л.~А.~Тумаркин,  А.~Н.~Черкасов, 
и самый знаменитый  -- Л.~С.~Понтрягин.

Научная экспансия Московской математической школы -- отдельная большая история,
см., например, \cite{Gnedenko-book}, \cite{AGS}, \cite{Lyu-3}, \cite{Philips}. Мы оставляем этот предмет в момент
начавшейся цепной реакции. Но все же необходимо сказать несколько слов о 
математических работах двух важных героев
последующего изложения \S\S\ref{s:attack}--\ref{s:after}.

\sm

{\small

{\bf\punct Московские теоретико-числовики: Шнирельман и Гельфонд.%
\label{ss:ShGe}}
Напомним, что трансцендентные числа впервые построил Лиувилль (Joseph Liouville) 1844, а именно
он заметил трансцендентны разреженные ряды вроде
$\sum_{n=0}^\infty 10^{-n!}.$
 Кстати, в этом случае доказать трансцендентность легко (это приятное упражнение для студента-математика любого курса), но надо было быть Лиувиллем, чтобы догадаться искать
трансцендентные числа в таком виде. Понятно, что это число было специально придумано, чтобы быть
трансцендентным, доказывать же трансцендентность конкретных встречающихся в математике чисел тяжело. В 1873г. Эрмит (Charles Hermite) наконец доказал трансцендентность числа $e$.
В 1882 Линдеман (Ferdinand von Lindemann) доказал, что трансцендентны числа вида
$e^\beta$, где $\beta$ алгебраично. В частности, число $e^{i\pi}=-1$
не трансцендентно, а потому трансцендентно $\pi$. Это влечет неразрешимость задачи
о квадратуре круга. Линдеман (по другим утверждениям Вейерштрасс,
но по-моему, у Линдемана есть равносильное этому высказывание) также доказал, что числа вида $e^{\beta_j}$
алгебраически независимы над алгебраическими числами (если $\beta_j$ -- алгебраические числа, линейно независимые
над $\Q$). После этого,
по-видимому, прогресс в этой области приостановился. 

Проблема Гильберта  состояла в вопросе: доказать (или опровергнуть), что
все числа вида $\alpha^\beta$, где $\alpha$, $\beta$ алгебраичны, причем $\alpha\ne 1$,
а $\beta\notin \Q$. 

Во второй половине 1920х началось оживление в этой области, подходы к изучению
трансцендентности были связаны
с исследование поведения подобранных голоморфных функций. Автор не в состоянии
проследить это движение. Были, например, работы Д.~Д.~Мордухай-Болтовского (см. Мат. сборник,
1927, 34:1), быть может оказавшие влияние на Гельфонда. 

 Согласно \cite{Gelfond-60},
\begin{quotation}
 Александр Осипович Гельфонд, один из выдающихся советских математиков, родился в Петербурге 24 октября 1906 г. в семье врача. В 1924 г.
он поступил в Московский университет и в 1927 г. его закончил. В 1930 г.
закончил аспирантуру. Его работой в аспирантуре руководили А. Я. Хинчин и В. В. Степанов. В 1929—1930 гг. Александр Осипович преподавал
в МВТУ. С 1931 г. ... работает [написано в 1966г.] на механико-математиче­ском факультете%
\footnote{С весны 1931г. - профессором.} МГУ. 
\end{quotation}

Напомним вкратце формальную историю решения проблемы Гильберта,
(доказательства, см., например, в \cite{Baker}).
В 1929г. Гельфонд показывает \cite{Gelfond-1929} (см. также \cite{Gelfond-marksizm}), что трансцендентно любое число вида $\alpha^{i\sqrt n}$, где
$\alpha$ алгебраично, а
$n$ - целое число, не являющееся полным квадратом (например, $(-1)^i=e^\pi$),
в этот момент он становится знаменитостью.
На следующий год  Р.~О.~Кузьмин  \cite{Kuzmin} доказывает такую же теорему
для чисел вида $\alpha^{\sqrt n}$. В 1933г. К.~Боле (Karl Boehle, Math.Ann., 1933)
получает такую теорему. Пусть $\alpha$ алгебраично. Пусть $L$ -- алгебраическое расширение $\Q$
степени $k$. Пусть $\beta_j$ -- базис пространства $L$ над  $\Q$.
Тогда хотя бы одно из чисел $\alpha^{\beta_j}$ трансцендентно.

Наконец, в 1934 Гельфонд решает в положительном смысле 
проблему Гильберта \cite{Gelfond-1934}, \cite{Gelfond-izv}. В тот же год (несколько позже) теорема   была независимо доказана
Т.~Шнайдером (Theodor Schneider, J. Reine Angew. Math. 172 (1934).).
По-видимому, решающий вклад в достигнутый в 1929--1933гг. успех был внесен Гельфондом.

Много позже, в 1967г. Бейкер (A. Baker) доказал трансцендентность всевозможных произведений
$$
\alpha_1^{\beta_1}\dots \alpha_k^{\beta_k},
$$
где $\alpha_j$ -- алгебраические числа, $\ne 0$, $1$, а $1$, $\beta_1$, \dots, $\beta_k$ линейно
независимы над $\Q$.

\sm

Шнирельман - автор нескольких приобретших сразу широкую известность результатов из разных разделов математики.
Кажется, все они относятся к 1928-1930г.,  см. их обсуждение
в \cite{Gelfond-Shnirelman}.

\sm

a) Люстерник и Шнирельман
в 1929г. аннонсировали 
доказательство следующей гипотезы Пуанкаре: на сфере с любой
гладкой римановой метрикой существует по крайней мере три замкнутых несамопересекающихся
геодезических \cite{Lyusternik-Shnirelman}.
Эта междисциплинарная работа (вариационное исчисление, топология, дифференциальная 
геометрия, функциональный анализ) имела большое значение. Быть  может, с нее начался
выход московской математики в функциональный анализ.

Что касается доказательства \cite{Lyusternik-Shnirelman-2}, то в нем были обнаружены дырки. Люстерник
несколько раз публиковал версии этой работы, самая подробная 
-- большой мемуар \cite{Lyusternik49} 1949г. По-видимому, окончательно
дырки были заделаны к началу 1990х%
\footnote{W. Klingenberg,  Trans. Amer. Math. Soc. 356 (2004), no. 6, 2545–2556:
	\newline
	\it In a short note, Lusternik and Schnirelmann
	[1929] sketched a proof of the fact that actually three simple closed geodesics exist
	on any orientable closed surface of genus 0, i.e., on a surface given by $S^2$, endowed
	with an arbitrary Riemannian metric. A complete proof of the theorem of
	Lusternik-Schnirelmann is notoriously complicated. Over the years, many papers
	have appeared on this subject and only a few are satisfactory; cf. Taimanov [1992].},
примерно тогда же стало известно,
существование счетного числа замкнутых геодезических, правда, вообще говоря, самопересекающихся (подробнее об этой задаче,
 см. \cite{Taimanov}).

\sm

b) Шнирельман в 1929г. доказал такую теорему \cite{Kvadrat}: дана жорданова кривая на плоскости с кусочно-непрерывной кривизной. Тогда в нее можно вписать квадрат (то есть существуют четыре точки кривой,
лежащих в вершинах квадрата). Из его публикации непонятно, знал ли Шнирельман о том,
что проблема существования вписанного квадрата для {\it непрерывной} жордановой
кривой была сформулирована О.~Тёплицем (Otto Toeplitz) в 1911г.
(см. \cite{Toeplitz}). Вскоре после этого вышли
две статьи А.~Эмха (Arnold Emch. Amer. J. Math. 35 (1913), Amer. J. Math. 37 (1915)),
где теорема была доказана для кусочно-аналитических кривых. В частности, это верно для ломаных.
Дальше можно  снижать требование на  гладкость:  аппроксимировать кривую ломаными
и брать предельное положение квадрата. Но при этом нужно обезопаситься,
	чтобы квадрат в пределе не перешел в точку. В общем, для непрерывных кривых
проблема остается нерешенной. 

В любом случае статья Шнирельмана была одной из
ранних  работ по дифференциальной топологии и способствовала интересу к этому предмету в СССР%
\footnote{Рассмотрим четверки точек $(A,B,C,D)$ плоскости. Пространство всех четверок есть $\R^8$. Множество четверок, лежащих на кривой образует четырехмерное подмногообразие, гомеоморфное тору. Множество четверок, удовлетворяющим следующим условиям на расстояния,
$AB=BC=CD=DA$ образует четырехмерную поверхность. Нужно показать, что две четырехмерных
поверхности имеют общую точку.}.

\sm

c) Третий известный результат Шнирельмана был связан с проблемой Гольдбаха:
верно ли, что любое четное число есть сумма двух простых чисел, а любое нечетное - сумма трех простых. Шнирельман в 1929г. доказал, что существует такое конечное
$N$, что любое натуральное число представимо в виде суммы $\le N$
простых \cite{Shnirelman29}. Это было, видимо, первым положительным утверждением в этом 
направлении. Но число $N$, извлекавшееся из доказательства, было немаленьким (800 000).

В 1937г. Виноградов решил проблему Гольдбаха для достаточно больших нечетных чисел,
а именно для чисел, больших, чем
$3^{3^{15}}$, по порядку величины это $10^{6 846 16}$.
Эту границу потом долго понижали (не приблизив ее сколько-нибудь к астрономическим
числам%
\footnote{Для сравнения объем Вселенной, измеренный в кубических
нанометрах, $10^{74}$ или что-то вроде этого. C другой стороны,
математики обычно не замечают, с какими конечными величинами они порой работают.
Автор данных записок как-то поймал себя на том, что он оперирует с конкретным конечным
вероятностным пространством, содержащим $10^{3 000 000}$ точек,
а потом сообразил, что для вероятностного пространства это совсем немного.
Я не знаю, пытался ли кто-либо оценивать число коэффициентов,
в системах многочленов, задающих алгебраические многообразия
(полученные в результате долгих алгебро-геометрических манипуляций).
В общем, число Виноградова великовато для возможностей прямого перебора
в нашей Вселенной,
но называть его очень большим не стоит.}.
Сейчас гипотеза Гольдбаха для трех чисел, быть может, доказана без оговорки <<для достаточно больших>> 
(Harald  Helfgott).
}

\sm

%10^{16} -cветовой год
%5\cdot 10^{10} световых лет
%10^9 милимикрон

{\bf\punct ТФДП и московская математика.%
\label{ss:tfdp}}
Цитируем Люстерника \cite{Lyu-3}:
\begin{quotation}
 Как же все-таки школа со сравнительно узкой тематикой стала
базой для развития универсальной математической школы? Оказалось, что
«прорыв» научного фронта, хотя бы на узком участке, имеет гораздо более
широкое значение, чем может показаться с первого взгляда...

Высокий уровень на «узком основании» в начальный период развития
научного центра представляет часто более благоприятную ситуацию для
дальнейшего его развития, чем более низкий и равный уровень на более широком
основании; 
\end{quotation}

Как будто, так и вышло. Мастерское владение ТФДП   оказалось тараном 
в руках молодых лузитан, оно пригодились и в теории чисел, и в
теории вероятностей, и в функциональном анализе, и в комплексном анализе.
А дескриптивная теория множеств постепенно (усилиями, главным образом, Новикова),
переходила в логику.

Опять цитирую Колмогорова:
\begin{quotation}
Все мои занятия по теории вероятностей совместно с А.~Я.~Хинчиным,
весь вообще первый период занятий этой теорией отмечен тем, что мы применяли
методы, разработанные в метрической теории функций. Такие темы,
как условия для применимости закона больших чисел, условие сходимости
ряда независимых случайных величин, велись по существу методами, выкованными
в общей теории тригонометрических рядов, т. е. методами, разрабатывающимися
Н.~Н.~Лузиным и его учениками.
\end{quotation}

Ну, и повторю еще раз, ТФДП в то время была новой наукой, и начинать ей активно заниматься
можно было без долгого обучения. Мы уже заметили, с какой скоростью входили в тему 
Александров, Хинчин, Меньшов, Суслин, Колмогоров. Так что ТФДП еще сыграло роль учебной 
площадки, на которой учились самостоятельному мышлению люди, впоследствии ставшие знаменитыми в других науках.

\sm

{\bf\punct Лузитания.%
\label{ss:luzitaniya}}
Про Лузитанию впоследствии было много ностальгических воспоминаний, я, по-видимому, должен что-то процитировать.
Александров \cite{Alex-auto1}:
\begin{quotation}
 Лузитания считалась «орденом», «командором» которого был Н.~Н.~Лузин, а «гроссмейстером» — Д.~Ф.~Егоров. 
 Лузитания была действительно уникальным и неповторимым коллективом молодежи, жившей не только напряженной,
 насыщенной математической жизнью, но и жизнью, которая была непосредственно радостной и веселой. 
 Такой коллектив мог возникнуть лишь в самые первые годы революции, когда вся страна переживала единственный в истории,
 неповторимый подъем во всех областях своей жизни. Веселые и необыкновенно оживленные лузитанские собрания,
 на которых, кстати сказать, не было ни капли вина, происходили при непременном участии учителей 
 всей этой молодежи — Д.~Ф.~Егорова и Н.~Н.~Лузина, и это говорит о том, 
 насколько простыми и непринужденными могут быть отношения между учителем и его учениками.
\end{quotation}

Лаврентьев \cite{Lavr1}:
\begin{quotation}
   Лузитанцы  признавали  двух  начальников  «бог-отец-Егоров  и  бог-сын-Лузин»; 
   Н.~Н.~Лузин  новичкам-лузитанцам  говорил:  «главный  в  нашем  коллективе  Егоров,  окончательная  оценка  
работы,  открытия  принадлежит  Егорову».  Новички  быстро  ориентировались:  
Д.~Ф.~Егоров  — форма,  а  содержание  — Н.~Н.~Лузин,  но  все  основные  и новички  соблюдали  форму  и  три  раза  в  году  приходили  домой  к Д.~Ф.~Егорову  
(пасха,   рождество,  именины).  Как  правило,  говорили  старшие,   остальные   
молчали  и  ждали  конца   визита.   

С  Н.   Н.  Лузиным  отношения  были  много  проще,   и  Н.~Н.~Лузин   сам   
веселился,  когда  лузитанцы  вытворяли  студенческие  озорства.  В  управлении  
«Лузитанией»  главными  помощниками  Н.~Н.~Лузина  были  три  Павла  (из  которых  один  Славочка)  со  своими  функциями:  П.~С.~Александров  — создатель 
тайн  «Лузитании»,  П.~С.~Урысон  — хранитель  тайн  Лузитании  и  В.~В.~Степанов
— глашатай   тайн   «Лузитании»  — ПСЫ). 
Был  такой случай.  Пришли на  лекцию Н.~Н.~Лузина,  ждали час;  по инициативе 
 (кажется,  П.~С.~Александрова)  все  20 пошли  к  Н.~Н.~Лузину  домой.  
Дверь  открыла   жена  Н.~Н.~Лузина   и  сказала,  что   Командора   похитила   
девушка  — увела  в  Малый  театр.  Всеобщее   возмущение   (особенно  лузитанок  —
 все    были    влюблены    в    учителя).    Составили    под    руководством    
П.~С.~Александрова  «грозное»  письмо  с  порицанием.  
Вышли   на   улицу   (Арбат)  — что   делать?   Решили   все   идти   в  театр.

В складчину  купили
 2 билета  и приемом «прошли двое,  один остался,  а другой 
вышел  с  2  билетами»  прошли  в  театр  все  20.  Когда  наступил  антракт,  трое  
прошли  за  командором  и,  под  угрозой  большого  шума,  привели  Н.~Н.~Лузина
в  фойе,  где  все  его  обступили,   стали  качать  и  петь  лузитанскую  песню  
«Наш  бог Лебег,  кумир  интеграл  рамки  жизни  сузили,  так  приказал  нам наш 
командор  Лузин...»,  вечер  закончился  на  Арбатской  площади,  где танцевали 
фокстрот  под  «гребенки».  
\end{quotation}

Люстерник:
\begin{quotation}
 Этот    веселый    математический    маскарад    продолжался    весь   1921   г.   
А  уже   в  следующем  году,  когда  жизнь  начала   стабилизироваться,   все  это  
постепенно   отпало.   Осталось   просто   лузинская   школа   теории   функций...

Лузитания,   к  счастью,  была  добровольным  обществом,  независимым 
  в   «административном    отношении»   от   своего   метра:   в  нее   входили   
и   из  нее  уходили   по   собственному   желанию.
\end{quotation}

И еще Люстерник:
\begin{quotation}
 В Лузитании были все признаки секты — свой социальный микроклимат
со своими ответами на вопрос «что такое хорошо и что такое плохо», с сектантской
непримиримостью и узостью. Некоторые молодые члены Лузитании
не хотели ничего признавать, кроме тех вопросов теории функций, которые
там исследовались. Для других разделов математики придумывались
шуточные названия: «уравнения в несчастных производных», «теория
неприятностей», «разные конечности» и т. д. (Это относилось к самым молодым
лузитанцам.) Но, очевидно, в тех условиях создание своего микроклимата
было необходимо для того, чтобы увлечься математической наукой.
\end{quotation}

Вообще рассказы о Лузине при их художественности во многом штампованы.
Эта цитата относится к 1967 году, но буквально то же самое теми же словами писал Гнеденко \cite{Gnedenko-book}
(который не был очевидцем Лузитании) в 1946 году, и в сокращенном варианте это писалось
в статье
Александрова, Гнеденко и Степанова в \cite{AGS}, 1948. Сомневаться в картинке с натуры нет оснований,
но в отношении напрашивающегося вывода о сектантстве лично Лузина стоит проявить осторожность. 
 Лузин,  конечно, был фанатиком своего научного направления,
и считал, что его ученики должны заниматься тем же. Но, скажем, списки
докладов на лузинском семинаре, приведенные в \cite{Lyu-2}-\cite{Lyu-3},
довольно разнообразны. Идеи, что ТФДП как-то отрицает обычный анализ, 
в сочинениях Лузина не видно. Что касается дескриптивной теории
множеств, то Лузин наблюдал несоответствие ее классической математике,
и считал, это существенным недостатком теории множеств
и основанием для ее реформирования. Сам Люстерник приводит данную ему Лузиным программу
аспирантского (магистерского) экзамена.
\begin{quotation}
 \tiny
1. Классический анализ (в основе лежал 3-й том известного курса анализа Пикара).

2. Теория функции действительного переменного (книги Бореля и Лебега). 

3. Теория функций комплексного переменного (книга Озгуда).

4. Обыкновенные дифференциальные уравнения (кажется, Гурса, 2-й том; позже
стал использоваться курс Бибербаха).

5. Уравнения в частных производных (кажется, тоже Гурса).

6. Интегральные уравнения (Гурса, 3-й том, иногда рекомендовали и классическую
книгу Гильберта, а позже курс Ловита).

7. Вариационное исчисление (курс Больца, Тонелли).

8. Уравнения с функциональными производными (это были элементы функционального анализа).

9. Геометрия 1 («Основания геометрии» Гильберта).

10. Геометрия 2 (большой курс дифференциальной геометрии Бьянки, конечно,
не полностью).

11. Теория чисел (Минковский, «Диофантовы приближения»).

12. Алгебра (курс Вебера).

13. Теория инвариантов (для этого предмета использовались записи лекций по
теории инвариантов, кажется, читанные Егоровым).

14. Теория вероятностей (не помню, что было рекомендовано).

15. Механика 1 (общая, курс Аппеля, 1-й том).

16. Механика 2 (сплошная среда, другой том Аппеля, а также статья
Н.~Е.~Жуков­ского «Видоизмененные методы Кирхгофа» о струйном обтекании).
\end{quotation}
(это много отдельных предметов). И где тут сектантство?
\begin{quotation}
 Сдачу экзамена по какому-нибудь предмету можно было заменить
докладом по этому предмету, содержащим новый результат, новое доказательство и т. п....
Кажется, наибольшее
число зачтенных докладов было у А.~Н.~Колмогорова — семь. 
\end{quotation}

\sm

{\bf \punct Лузин-лектор.%
\label{ss:lector}} Описывать манеру преподавания - дело сложное. Так или иначе,
лузинская школа самостоятельного мышления 
это важная часть нашей истории, поэтому приведем отрывки из воспоминаний разных авторов.

Известный инженер-нефтяник В.~Н.~Щелкачёв \cite{ShCh}:
\begin{quotation}
... Как он читал лекции нам? Он говорил, 
 что его долг – никогда не готовится к лекции, чтоб не «давать студентам жвачку». Он начинал
 лекции примерно так: «Мы вчера с вами доказали одну теорему, а на этой лекции попробуем доказать другую.
 Надо задуматься, как от той теоремы перейти к этой». Это была лекция-рассуждение. 
 Я, конечно, никогда не мог позволить себе так читать лекции, для этого нужно
 было быть действительно гениальным ученым. А он не стеснялся в конце лекции говорить:
 «Знаете, а ведь мы пошли по неправильному пути. Для доказательства попробуйте прочитать такую-то книжку».
....
  \end{quotation}
  
Рассказ Н. М. Бескина
\cite{Besk}, Лузин в воспоминаниях профессионального педагога:
\begin{quotation}
 Я не ученик Лузина и именно поэтому взялся писать о нем. Его учениками были десятки, а слушателями
 -- тысячи человек... А роль Лузина как университетского преподавателя огромна. Одни ученики
 не могут воссоздать его полный портрет... 
 
 Прежде, чем рассказывать о лекциях Лузина, хочу напомнить приметы того далекого времени.
 Университет, как и вся Москва до 1921 года не отапливался. Зимой 1920/21
 учебного года Лузин читал лекции в валенках, в шубе с поднятым воротником и меховой шапке
 (в следующем году шуба была уже не нужна). Стипендии были лишь единичные, поэтому все студенты
 где-нибудь работали и посещали только те лекции, которые их интересовали. Они не имели
 иных стимулов к учению, кроме интереса к математике. Математика как профессия казалась
 бесперспективной (в отличие от инженерных специальностей). Диплом был не нужен. Презрение к документам
 было столь велико, что студенты, оканчивавшие университет, не находили нужным
 зайти в канцелярию и получить полагавшийся им диплом...
 
 Чем измеряется влияние лекций Лузина? Тем, что они запоминались на всю жизнь. Лучше сказать не запоминались, а входили
 в сознание слушателя как элементы его математической культуры и математической методологии.
 А в этом и заключается высшая цель лектора...
 
 Будучи студентом второго курса я посещал лекции (необязательные)
 по теории функций действительного переменного. Я не записывал лекции, а только внимательно слушал.
 Экзамена по этому курсу не было, и специально я теорией функций действительного переменного
 не занимался. Она впервые понадобилась мне через несколько лет после окончания
 университета. Я попытался, не обращаясь к литературе восстановить по памяти курс Лузина.
 Изо дня в день я сидел над листом бумаги и вспоминал. Мне удалось все восстановить.
 В доказательствах теорем основная идея всегда легко вспоминалась, а подробности
 приходилось достраивать самостоятельно...
 
 Когда Лузин читал обязательные курсы, большинство его лекций не имело отношения к программе.
 Это был свободный разговор, большей частью экспромтом, на математические темы.
 В 1921/22 учебном году он читал на первом курсе высшую алгебру.
 Я хорошо помню о математическом творчестве, об аксиоме Цермело,
 об основаниях математики. Он разбирал недавно вышедшее английское
 издание Б.~Рассела <<Введение в философию математики>>. 
 В результате программа выполнялась лишь в малой степени. Первый семестр посвещался теории определителей.
 Лузин изложил не более трети программного материала.. Остальные две трети предстояло самостоятельно
 изучить (для сдачи экзамена) по книге Нетто <<Теория определителей>>...
 
 Теперь я понимаю педагогическое кредо Лузина так: он преподавал в первую очередь математику,
 и лишь во вторую -- какой-нибудь ее раздел. Он был математиком, а мы - будущим математики.
 Он вводил нас в науку, преподавал свой подход к проблемам математики.
 Разумеется, нельзя забывать, что Лузин преподавал профилирующий предмет взрослым людям...
 
 Вначале нас удивляло, что Лузин не готовился к лекциям. Доказывая теорему, он часто натыкался на препятствие
 и пробовал другой путь. Он видел свою задачу не в предъявлении готового доказательства,
 а в том, чтобы показать, как мыслит математик, как он ищет истину и преодолевает препятствия.
 В большинстве случаев он действительно решал вопрос у доски. Однако бывали случаи, когда он ошибался нарочно.
 Зачем?.......

 Доказывалась теорема: множество точек открытого квадрата и множество
 точек его стороны равномощны [понятно как]
 
 -- Это очень просто, -- говорил Лузин, мы разбираем десятичные знаки
 как будто сдаем колоду карту: одну карту одному партнеру другую другому. Вот и все доказательство.
 
 Действительно, как просто...
 
 Изложив это доказательство Лузин вдруг принял очень озабоченный вид и сказал:
 
 Прошу прощения, я ошибся. Это доказательство неверно. [понятно, где]....
 
 Он задумался, а мы все огорчились. Почему? Потому что было досадно, чтоб такое хорошее
 доказательство приходится забраковать из-за такого ничтожного, непринципиального,
 но все же несомненного дефекта. Возникли разные предположения, как его устранить.
 Все предложения Лузин внимательно слушал и обстоятельно критиковал...
 
 Во многих доказательствах идея проста, но она непосредственно не проходит.
 Преодоление же препятствий бывает гораздо сложней, чем основная идея. При изложении доказательства
  в готовом виде эти подробности заслоняют идею... Лузин расчленял эти моменты.
  Простая идея запоминается навсегда, а подробности пусть забудутся...
\end{quotation}

Наконец, рассказ ученика, 
Люстерник \cite{Lyu-3}:
\begin{quotation}
На  лекциях  часто  ставились  
задачи.   По  этому  поводу   один  из   бывших   студентов   университета   сказал:   
«Другие   профессора   показывают   математику   как   завершенное   прекрасное   
здание  — можно  лишь  восхищаться   им.  Лузин  же  показывает   науку   в   ее   
незавершенном   виде,   пробуждает   желание   самому   принять   участие   в   ее   
строительстве».  Не  знаю,  была  ли  это  случайная  обмолвка,   но  после  одного  
замечания  Н.~Н.~Лузина  мы  с  Ниной  Бари  «состязались»  в  решении  якобы  
нерешенной   задачи,   решение   которой   легко   было   найти   во   многих    более    
подробных    курсах....

Читает  Николай  Николаевич  лекцию,  вдруг  задумывается:   «Я  не  могу  
восстановить  доказательства,   может   быть,  кто-нибудь  из  коллег  мне  поможет?»
Царит   напряженная   тишина,   все  думают,   кто-то   подходит   к   доске,   
пробует   доказать,   не   выходит,   и  покрасневший   возвращается    на    место.    
Наконец,  счастливец  под  завистливые  взоры  остальных  встает  и  рассказывает  
у  доски  найденное  доказательство.   Метр  говорит:   «Спасибо,  коллега»,   благосклонно  улыбается
.  Победитель  радостный  садится  на  место.  В  самом  ли  
деле  профессор   «утерял»  доказательство  или  это  была  хорошо   проведенная   
игра,  прием  для  пробуждения   активности   и  самостоятельности?   Я   слышал   
от   бывших   студентов,   слушавших   общий   курс   алгебры   лет   пять   спустя,   
что  как-то  позже  Лузин  предложил  аудитории  найти  «новый»  метод 
доказательства  какой-то  теоремы  из  теории  детерминантов  и,  искусно  «дирижируя»  
аудиторией,   «получил»   от  нее  это  доказательство,   и  это прошло  с  большим  
подъемом   и   интересом....
\end{quotation}

Напоследок, отрывок из А. П. Юшкевича:
\begin{quotation}
Каждая его лекция представлялась нам вдохновенным творческим процессом поиска и открытия истины.  То задумчивой паузой, то размышлением вслух, то прямым предложением слушателям вместе с ним разрешить встретившуюся трудность, Лузин приковывал внимание аудитории, делая ее активной участницей решения рассматриваемого вопроса. Этот последний прием Лузин с особым успехом применял в своих специальных курсах...

Быть может, при изложении высшей алгебры подобные эффекты представляли собой лишь педагогический прием, или же наш профессор не успевал  продумать  все перед лекцией, так или иначе, мы испытывали необыкновенное увлечение. Говорили также, что некоторые приемы повторялись из года в год, как у знаменитого историка Ключевского, но для тех, кто был их свидетелем, и в какой-то мере соучастником впервые, это не имело значения. 
\end{quotation}

{\bf\punct Лузин - научный руководитель.%
\label{ss:ruk}}
Цитируем интервью с Колмогоровым, \cite{Kolmogor}:

\begin{quotation}
В чем состояли эти занятия с Н.~Н.~Лузиным?

— Каждый ученик приходил к Николаю Николаевичу Лузину в его
арбатскую квартиру раз в неделю вечером — в постоянно выделенный для
него день недели. Мой день был общий с Петром Сергеевичем Новиковым....

Для Москвы, для Московской математической школы важное
значение имел новый подход к работе с молодежью. Существенным в этом
подходе было вполне индивидуальное личное руководство, а также умение
придавать избранной тематике особенную значимость. Н.~Н.~Лузин настойчиво
внедрял следующий метод работы (он и сам работал таким образом,
и приучал к этому своих учеников): берясь за какую-либо проблему, надлежит
смотреть на нее с различных точек зрения. Надо пытаться доказывать
гипотезу и одновременно опровергать ее. Если доказательство не выходит,
надо переходить к опровержению гипотезы, к построению противоречащего
примера. Если не получается построение — надо снова вернуться к доказательству.
И пока не получится результат, нельзя покидать данную область.
В теории функций действительного переменного такая установка двойного
видения (поиск доказательства — поиск опровержения), такой подход к делу
естественно привел к культивированию чрезвычайно высокой техники построения
примеров (или, как теперь принято говорить, контрпримеров).
\end{quotation}

Люстерник \cite{Lyu-3}:
\begin{quotation}
Лузин   смело   провел   коренную   реформу   в   деле   подготовки   молодых   
научных  работников  математиков.  Считалось  необходимым  до  начала   самостоятельной 
научной  работы  изучить  предварительно  много  толстых  фолиантов..... 
Лузин   сумел  избежать   в  подготовке   математической   молодежи   опасностей   
«переучивания»    и   дилетантизма,   призывая    к   ранней    самостоятельности,    
сочетающейся  с  продолжающимся  математическим  образованием.  Так  именно 
действовали  позже  и  другие  научные  руководители   в  МГУ.  Но,   такова   
диалектика  жизни,  это  привело  к  быстрому  распаду  Лузитании.   Лузитанцы   
стали   проявлять   самостоятельность   и   в  других   областях   математики...   
\end{quotation}

Лаврентьев \cite{Lavr1}:
\begin{quotation}
 В  этот период  расцвета  ярко  выявилась  основная  черта  школы Н.~Н.~Лузина
 — это  была  школа  развития  самостоятельного  мышления,  способности  
расчленять  проблемы,  искать  обходные  пути,  ставить  новые  проблемы.  
Развитие  этих  способностей  интеллекта,  характерная  для  школы  Н.~Н.~Лузина,  
было  очень  важно  тогда  (полстолетия  тому  назад)  и  приобретает  особое  значение
сегодня  в  эпоху  возросшей  роли науки  и  ее  значением  в  научно-техническом
прогрессе.  Большую   роль  в  развитии  «Лузитании»  играли
   лекции   Н.~Н.~Лузина: к  своим  лекциям  Н.~Н.~Лузин  готовился  только  вчерне  и  его  
лекции  были далеки  от стандарта  «хорошей  лекции». Н.~Н.~Лузин  часто опаздывал,  
но  слушатели   все  или  почти   все  лузитанцы  приходили   во  время  
и  в  ожидании  в  коридоре  велись  разговоры  по  проблемам  и  около  проблем....

  Лузин  заботился,  чтобы  лузитанцы  особенно  те,  которые  проявляли  
самостоятельность  в  мышлении,  не  теряли  времени  на  подготовку  к  экзаменам 
по  областям,  далеким  от  теории  функций  (астрономия,  физика,  химия,  
механика). Н.  Н.  Лузин давал  советы  — надо  хорошо знать  оглавления
 выучить  выборочно,  20—30\%.  Если  будет  задан  вопрос  из  незнакомой  части,  
то  (не  стесняясь)  надо  быстро  начать  рассказывать  из  знакомого   раздела.   
Преподаватель,  как правило,  останавливать  не будет;  если  остановит,  то надо 
быстро  начать  говорить  из  чего-нибудь  другого.  Астрономы  имели  привычку
оставлять  студентов  (для  подготовки)  одних  — в  этом  случае  надо  было  
осторожно  пронести  с  собой  учебник  и  списать  то,  что  задано.  
Для  аспирантов  было  «негласно»  установлено  такое  правило:  если  у 
аспиранта  по  теме  экзамена  есть  самостоятельный  результат,  то  спрашивают  
только  по  этому  результату.  Мы  все  стремились,  вместо  изучения  толстой  
монографии  200—300  стр.  (как  правило,  на  иностранном  языке),  придумать  
новую  постановку   (обобщение)  задачи...  

 Основная черта лузинской школы - развитие самостоятельного мышления - стала для меня главенствующей, где бы я ни работал.
\end{quotation}

\sm

{\bf\punct Международные контакты.} Входившим в науку молодым людям еще предстояло приобрести международную известность.
Мы уже отмечали, что первые публикации Александрова, Меньшова, Хинчина и Суслина шли в Compt. Rendus 
через Адамара (сопровождаемые письмами Лузина). Еще одним путем (начало 20х) был журнал Fundamenta Mathematicae,
издававшийся Серпинским с 1920г.,
этот путь тоже, надо думать, поддерживался Лузиным, впрочем, старшие из лузитан были лично знакомы с Серпинским.
В дополнение к этому приведем известие А.~П.~Юшкевича, 1982, \cite{Yushkevich-Duvakin}:
\begin{quotation}
 Вот Лузин рекомендовал Клейну этих своих молодых учеников, просил им оказать всякое содействие. 
 Тут начались контакты московской школы, москвичей вообще с Гёттингеном, которые продолжались вплоть до мрачных лет,
 когда к власти пришли нацисты и уничтожили, по существу, гёттингенскую школу Феликса Клейна и Гильберта. Вот.
\end{quotation}

Во второй половине двадцатых Александров, Хинчин, Колмогоров  приобрели международную известность и сами стали рекламой
советской математики. Но сначала
к этому должен был приложить
руку  сам Лузин. 

\sm

{\bf\punct  Из интервью Люстерника.} Вот еще несколько эмоциональный отрывок, из интервью 1970г.,
записанного на диктофон \cite{Lyusternik-Duvakin}:
\begin{quotation}
 Л.~Л.: Да, я вам скажу. Это был человек очень такой своеобразный, такой повышенной нервности. Но дело в том, 
 что более спокойные люди и трезвые, например, Чаплыгин, не могли в трудных условиях создавать вокруг себя…
 
В.~Д.[Виктор Дувакин]: Коллектив.

Л.~Л.: …молодежь. Нужна была вот эта нервозность какая-то. Словом, может быть... как-то... он именно был героем таких трудных времен. 
Когда жизнь наладилась, его влияние как-то стало меньше. Вот это очень любопытно, что максимум его влияния... вот 
его влияние — это было военные годы и первые послереволюционные годы 
— 14-й — 24-й год. Дальше начало падать. Чем труднее времена, тем относительно его влияние возрастало. 
Такая нервная личность — она больше как-то привлекала, чем такие... такой крупный ученый, как Чаплыгин... 
Понимаете, увлечь нужно было... надо было как-то эмоционально заразить, нужно было иметь эмоциональный потенциал. 
Ну, мысль понятна.
Он был очень нервный человек, человек был очень противоречивый, очень сложный и противоречивый, 
и поэтому он возбуждал очень противоречивые эмоции, иногда переходившие от восхищения к слишком одностороннему отрицанию... 
Целая, можно сказать, это была интеллектуально-эмоциональная драма его отношений с учениками. 
Это было очень сложно.
\end{quotation}

\section{Распад%
\label{s:raspad}}
 
 \COUNTERS
 
 \epigraph{Фигура в полном смысле слова трагическая, в высшей степени противоречивая, которой были доступны,
 совершенно независимо от ее творческих взлетов... Лузин, конечно, был гораздо более крупным ученым,
 гораздо более крупным талантом, чем то, что ему удалось осуществить.
 Его возможности как математика были, несомненно, гораздо больше, но... 
 Теперь, его возможности как человеческой личности были тоже очень велики,
 но тут были чрезвычайно... Как человек, как человеческая личность он был способен к чрезвычайно большим взлетам и вдохновениям, 
 но, как часто бывает у крупных людей, он был способен и к глубоким падениям.}
 {Из диктофонной записи Александрова 1971г., \cite{Alexandrov-Duvakin}.}
 
 Так или иначе, в нескольких цитатах предыдущего параграфа видны отсветы событий 1936 года.
 Несколько участников событий, избегая подробностей, впоследствии рассказывали о конфликтах
 в Лузитании.
 
 \sm
 
 {\bf\punct Показания Люстерника.%
 \label{ss:raspad-lyusternik}}
 Опасную и сложную  тему распада школы Лузина он предпочитает излагать в стихах (\cite{Lyu-3}. Почему-то во всех
 изданиях присутствуют многоточия (относительно полный  текст есть в 
 \cite{Lust-memoirs}), см. также \cite{Lyusternik-Duvakin}.

\begin{quotation}
 А дальше все как будто просто —
 \newline
Процесс естественного роста,
\newline
Тематика все расширялась,
\newline
Своей дорогой каждый шел —
\newline
И школа Лузина распалась
\newline
На ряд блестящих новых школ.
\newline
Но был мучительно тяжелым
\newline
Процесс распада этой школы.
\newline
Держалась крепко Лузитания
\newline
Огромным шефа обаянием.
\newline
Один не может интеллект
\newline
Эмоциональный взять барьер,
\newline
Большую роль играл аффект —
\newline
(Вот чувство дружбы, например,
\newline
Что двух П. С.%
\footnote{П.~С.~Александров и П.~С.~Урысон}
соединяло,
\newline
От Лузина их отдаляло,
\newline
И вот они пошли вдвоем
\newline
Топологическим путем.)
\newline
Наплыв эмоций (в плане личном),
\newline
Пересыщение привычным,
\newline
Желание самому стать первым
\newline
Иль расшалившиеся нервы,
\newline
Да мало ли что, но кто куда,—
\newline
Птенцы уходят из гнезда,
\newline
А это Лузин, хоть скрывал,
\newline
Болезненно переживал...
\end{quotation}

Вот рассказ Люстерника о себе:

\begin{quotation}
 Впервые задал я вопрос:
\newline 
К чему и театральность поз,
\newline
И тон его такой слащавый?
\newline
Не мы ль его раздули славу?
\newline
От Лузитании тогда
\newline
Стал отходить, хоть и не сразу,
\newline
А в аспирантские года
\newline
Я с Лузиным был мало связан.
\newline
Сначала просто так болтался:
\newline
За слишком многое уж брался.
\newline
Потом работал в направленьях,
\newline
В Москве тогда отнюдь не модных —
\newline
В вариационном исчислении,
\newline
В задачах в частных производных.
\newline
Я метод сеток развивал.
\newline
(К вопросу о приоритете —
\newline
Году то было в двадцать третьем.)%
\footnote{Рассматривается уравнение Лапласа в двумерной области с гладкой границей
и его приближение разностной системой на прямоугольной сетке. Было доказано,
	что при измельчении сетки решения разностных систем сходятся к уравнению Лапласа. 
	Статья Люстерника на эту тему вышла в 1926 году в <<Мат.сборнике>> \cite{Lyu-Dirichlet}.
	В списке выступлений на Мат.обществе за 1924 год упоминается
доклад	Люстерника <<Задача Дирихле>>.  Вопрос о приоритете действительно стоял,
в 1925 Курант (R.~Courant) опубликовал статью на ту же тему. Впрочем на эту тему были и более ранние работы
(J.~le~Roux, 1914; R.~G.~D.~Richardson, 1917), по-видимому, ни Люстернику, ни Куранту в тот момент не известные,
см. \cite{CFL}. }
\end{quotation}

И немножко прозой:
\begin{quotation}
 Но уже в конце 20-х годов начала чувствоваться некоторая изолированность
Н.~Н.~Лузина в коллективе московских математиков. Это выявилось
в эпизоде сравнительно незначительном: при выдвижении кандидатур
в Академию в московских математических организациях кандидатура Лузина
не была поддержана большинством математиков%
\footnote{Из Академической стенограммы:
\newline
ХИНЧИН. {\it Я считаю, что основных причин ухода H.~H. из Московского университета
было две, и первая из них та, о которой сегодня не говорилось, а именно, то обстоятельство, что Институт математики,
в котором H.~H.~Лузин работал, высказался в свое время
против избрания его академиком, когда выставлялась его кандидатура в Академию наук.
Это, как я мог наблюдать совершенно точно, создало первую трещину, которая, по-моему, 
и до настоящего времени не заросла вполне между H.~H.~Лузиным и Московским
университетом. Здесь, стало быть, причина личного характера.}}
\newline
\phantom{.}\qquad
Подразумеваемая Хинчиным вторая причина обсуждается ниже в п.\ref{ss:prompartiya}., входивших в Лузитанию,
и не получила официальной поддержки; практического значения это не
имело, но Н. Н. Лузин почувствовал себя «развенчанным королем» (результат
выборов был неожиданным: Н. Н. Лузин был избран академиком... по философии)
\end{quotation}

{\bf\punct Рассказ Бари.%
\label{ss:raspad-bari}}
 \begin{quotation}
 Талантливые ученики редко похожи на своих учителей. 
Поддаваясь авторитету, научному обаянию, творческой мощи своего
учителя, они вместе с тем вносят в науку и свои личные идеи, свое
личное понимание задач, стоящих перед наукой, путей ее 
развития. И это новое, вносимое талантливыми учениками, часто весьма
сильно расходится с научным мировоззрением учителя. В этом
расхождении проявляется диалектика развития науки, а 
возникающие противоречия представляются движущей силой научного 
прогресса.

За восемь лет, протекших с начала работы Н.Н. в Московском
университете, многое изменилось. Научная молодежь, которая с
таким благоговением слушала вдохновенные лекции Н.Н. в 1914
году, превратилась в 1922 г. в зрелых ученых, которые шли уже
своими собственными научными путями, вырабатывали свое 
самостоятельное научное мировоззрение. Наконец, и самая школа 
приняла совершенно другой вид. Если в 1914 году ближайших 
учеников Н.Н. было 4-5 человек, то в 1922 году в направлении, 
намеченном Н.Н., работали уже десятки молодых, талантливых и 
энергичных ученых. Чтобы сплотить эту многочисленную группу 
инициативной и талантливой молодежи, уже недостаточно было 
одного научного таланта руководителя, умевшего зажигать творческую
инициативу молодежи. Нужен был и другой, более прозаический,
но также необходимый, чисто административный талант главы
школы, умение держать в равновесии центробежные силы, 
вызываемые различием научного мировоззрения, умение превратить эти
центробежные силы в мощный фактор дальнейшего развития 
школы, где возникающие противоречия способствуют только 
успешному развитию пауки.

Н.Н. обладал исключительным талантом учителя; с 
поразительным успехом умел он своим преподаванием зажигать любовь к
научному творчеству, вовлекать молодежь в трудную, 
напряженную, а часто и неблагодарную работу по исканию научной истины.
{\bf Но Н.Н. в весьма малой степени обладал административными 
талантами, необходимыми для руководителя большой научной 
школы.}

..............Еще в 1919 г. отошел от 
работы в этом коллективе талантливый и многообещающий М.~Я.~Суслин; он вскоре преждевременно погиб, заразившись сыпным 
тифом. Ряд самых старых участников коллектива отошел от работы в
нем и пошел своим собственным путем....
\end{quotation}

Надо сказать, что распад группы Лузина и возникновение нескольких
новых школ были явлениями положительными (здесь стоит согласиться с Люстерником), да и вряд ли это можно было остановить.
Но важно, что вполне пролузинский источник подтверждает отсутствие у Лузина 
<<административных талантов>>. 
Для сравнения, из Академической стенограммы, 1936:
\begin{quotation}
	ХИНЧИН.
... И я считаю, что Н. Н. Лузин, при всех его
	прежних заслугах, которых я не отрицаю, {\bf к организаторской деятельности — при всех его
	личных свойствах — не способен.}
	\end{quotation}

Что ж поделать, при многих
выдающихся качествах этих способностей у Лузина не оказалось.
Ничего страшного в этом не было бы, но этому человеку вскоре предстояло занять должность
главного советского математика.

Продолжим цитировать Бари:
\begin{quotation}
 Конкурентом Н.~Н. в академии был его учитель Д.~Ф.~Егоров.
Его избрали «почетным академиком». Но, как сказал острый на
язык А.~Н.~Крылов, разница между «почетным академиком» и 
просто «академиком» в те годы была такая же, как между 
«милостивым государем» и просто «государем». Отсюда начался холодок в
отношениях с Д.~Ф.~Егоровым. С Б.~К.~Млодзеевским отношения 
испортились еще раньше, когда после смерти Н.~Е.~Жуковского в
1921 г. проводились выборы президента Математического 
общества и «Лузитания» весьма напористо и достаточно бесцеремонно
пыталась провести своего Д.~Ф.~Егорова против Б.~К.~Млодзеевского, 
который по мнению «лузитанцев» постарел и не соответствовал
духу науки. Президентом был выбран тогда Б.~К.~Млодзеевский%
\footnote{Об этом инциденте пишут также Александров \cite{Alex-MMO3} 
и Т.~А.~Млодзеевская, \cite{LL}, согласно ей в последний год жизни Млодзеевского
стороны помирились.}\dots
\end{quotation}

 \begin{quotation}
  С 1925 года начинается новая полоса жизни и деятельности
Н.~Н.~Лузина. Как в период с 1918 по 1924 г. его работа в 
Университете прерывалась его длительными поездками в Иваново, так и в
период с 1925 по 1930 г. он проводит большую часть времени в 
поездках за границу, отрываясь все более и более от работы в 
Московском университете, от руководства созданной им школой, от 
работы с советской молодежью. За границей в командировке от Наркомпроса он пробыл с августа 1925 г. по июль 1926 г. на 
Всероссийском математическом съезде весной 1927 г. (27.04-4.05) Н.~Н.
сделал один из основных докладов «Современное состояние теории
функций действительного переменного». В мае 1927 г. он опять
уезжает в командировку за границу, участвует в математической
конференции во Львове. В 1928 г. он снова уезжает за границу,
участвует в августе 1928 г. в Международном [математическом]
конгрессе в Болонье, читает лекции в Брюсселе, а затем в течение
двух лет почти непрерывно проживает в Париже, где работает
над своей книгой <<Lecons sur les ensembles analytiques>>\dots
 \end{quotation}

 \begin{quotation}
  Есть ученые-одиночки, ученые-отшельники. В своих 
уединенных кабинетах и лабораториях делают они великое дело науки;
всякое общение с людьми для них помеха в их непрерывной и 
напряженной работе мысли, им не нужно общение с молодежью, с
учениками, они стремятся избавиться от чтений лекций, от 
выступлений. Таков был А.~М.~Ляпунов — один из крупнейших русских
ученых математиков и механиков.

И есть ученые-коллективисты. Для их творчества необходим
коллектив, товарищи, молодежь, студенты. Только в живом 
общении с ними, в непрерывном обмене мыслей, планов, в обсуждении
достигнутых итогов продуктивно и напряженно работает мысль,
зреют научные идеи. Для таких ученых лекция, семинар - это 
радость творчества. Таков был Н.~Н.~Лузин.

Как мы видели, обстоятельства сложились так, что в течение
шести лет с 1924 по 1930 год в силу разных причин Н.~Н. был
оторван от той привычной обстановки Московского университета,
русской молодежи, своих учеников, которые теперь выросли, 
научно возмужали и делали первые решительные и самостоятельные
шаги. Шесть лет он был оторван от этой, такой необходимой для
его научного творчества среды, и каких шесть лет.
 \end{quotation}
 
 {\bf\punct Рассказ Лаврентьева.%
 \label{ss:raspad-lavrentiev}}
 
 \begin{quotation}
   Уже в первые годы внутри «Лузитании» возникали
  конфликты, но эти конфликты задевали лишь отдельных лузитанцев и не оказывали сколько-нибудь
  существенного влияния на общую деловую обстановку.
  
  Внутренние противоречия, которые привели к развалу, начались
  в 1925—1928 гг. Здесь надо отметить два обстоятельства: 1) сравнительно
  доступные задачи были к этому времени решены, и в главной тематике
  Н.~Н.~Лузина остались задачи, над которыми бились безуспешно много лет,
  и не только у нас но и за рубежом. Многие лузитанцы сами стали искать
новые направления. Столпы «Лузитании»-— П.~С.~Александров и П.~С.~Урысон
  начали успешно развивать топологию и ряд других новых направлений.
  Ряд сильных математиков, при участии В.~В.~Степанова, устремились
  в область дифференциальных уравнений (И.~Г.~Петровский, А.~Н.~Тихонов).
  Дольше других в области теории функций работал П.~С.~Новиков, 
  который открыл в ней новые пути и связал свои исследования с логикой.
  Сильное развитие получила теория функций комплексного переменного
  с выходами в геометрию и гидроаэродинамику.
  
  С переездом Академии наук в Москву, начала интенсивно развиваться
  теория чисел во главе с крупнейшими в мире специалистами в этой области.
  Распад школы был обусловлен также тем, что сам Н.~Н.~Лузин ряд лет
  посвятил второй своей большой монографии по дескриптивной теории функций
  и оторвался от молодых. 
 \end{quotation}

 {\bf \punct Каноническая точка зрения.%
 \label{ss:raspad-canonical}}
 Цитируем Академическую стенограмму
 \begin{quotation}
БУТЯГИН.%
	\footnote{Бутягин Алексей Сергеевич (1881-1958), ректор МГУ 1934–1941, 1943. Известно о нем мало.
	Окончил московский Физмат в 1906, математик (публикаций найти мне не удалось). И. о. ректора МВТУ им. Н.~Э.~Баумана в 1924–1929.
	Управление МГУ при нем было, видимо, весьма жестким.}
  Для  полноты  освещения  вопроса  об  уходе  из  Университета  следует  еще
указать и на следующее. С этим моментом совпало как раз обострение отношений Лузина
с целым рядом учеников, в частности, с П.~С.~Александровым. Я думаю, что эта 
последняя причина была чрезвычайно веской. 
 
............................................ 
 
БУТЯГИН.
  Не  было ли в этом  причиной то  обстоятельство,  что  новая  молодая  советская  школа 
  уже  перерастала  Лузина  и  что  он  чувствовал,  —  что  он  отходит  на  вторые
позиции?

ШМИДТ.
  Да,  вот о причинах ухода нужно учесть в  формулировке то,  что говорил тов.
Бутягин. Может быть, есть смысл добавить к начальной характеристике,  когда мы  говорим,
что Лузин возглавлял в течение такого-то периода школу, не было ли смысла в соответствии 
с этими фактами сказать:  «Однако эта школа распалась как только ученики Лузина 
стали сколько-нибудь самостоятельными учеными, ибо он не давал им возможности
самостоятельно  расти?»

АЛЕКСАНДРОВ.
  В  1922  г.  H.  H.  мне  прямо сказал:  «Пока вы занимаетесь топологией,
между  нами  не  может быть  никакого  научного контакта».

ШМИДТ.
 И не только по отношению к вам это было:  Вы были наиболее яркой  фигурой среди его учеников. Но это было по отношению ко всем.
 
 ...........................
 
 СОБОЛЕВ. Нельзя здесь говорить об отношении учителя и ученика, поскольку Новиков
 давно перерос Лузина. Никто в этом не сомневается.
   \end{quotation}

 Из воспоминаний Понтрягина:
 \begin{quotation}
 Талант Лузина как Учителя «с большой буквы», по-видимому, далеко превосходил его талант как творческого математика.
 Отсюда его трагедия. Его ученики часто начинали быстро превосходить его в своих достижениях
 и уходили из его области в более значительные разделы математики. 
 Он ревниво относился ко всему этому, и возникали враждебные отношения с учениками. 
\end{quotation}

И снова Люстерник и снова в стихах:
\begin{quotation}
 И, может быть, он не был избавлен от того, от чего предостерегал Тютчев
 
От желчи горького сознанья,

Что нас поток уж не несет

И что другие есть призванья,

Другие вызваны вперед.
\end{quotation}

\sm

{\bf\punct  Из интервью Юшкевича.%
\label{ss:raspad-yushkevich}} Вот отрывок из интервью Юшкевича \cite{Yushkevich-Duvakin}, взятого в 1982г.
и опубликованного лишь в 2005г.

\begin{quotation}
(Перерыв в записи)

А. Ю.
...крупные ученые, чем я, чрезвычайно высоко ценили деятельность Николая Николаевича Лузина. Вы начинаете записывать?

В. Т.[Валентина Трейдер]: Ну, да.

А.Ю.: Да, ну, пишите. Дело в том, что я сказал, что Лузин был очень сложной фигурой и в некотором смысле — трагической фигурой.
Дело было в следующем. Лузин, сам математик талантливый... ну, я не знаю, какого класса, ну, он был не идеальный [?] человек,
но это был человек большого таланта... Лузин, человек весьма талантливый, был совершенно исключительным учителем,
то есть воспитателем молодых ученых, я имею в виду. Он их заражал своим энтузиазмом, 
он совершенно переменил стиль обращения профессора со студентами, он был доступен, он держал себя, ну,
как равный или как первый, но среди равных; к нему можно было всегда прийти, его провожали гурьбой после лекции до дому, 
он рассказывал, он ставил задачи, он ставил вопросы, он с удовольствием читал работы их, он представлял их для публикации,
в «Доклады Парижской Академии наук» сперва, потом — в другие места. Он рекомендовал их, значит,
геттингенским математикам (он сам перед этим лет за десять был в Геттингене и провел там около года или полутора лет)
и так далее, и так далее. В общем, это был человек, который умел очаровывать. Он умел очаровывать,
но он мог и разочаровывать. Дело в том, что тематика, которой он занимался,
сама по себе очень важная и интересная, была все-таки ограниченной — 
теория функций действительного переменного в том направлении, которое разрабатывал он, — 
и складывалось дело так, что целый ряд молодых ученых, не чувствуя уже интереса к проблематике,
казалось бы, в основном изученной и малоперспективной, во всяком случае,
с трудом поддающейся дальнейшей глубокой разработке, переходили на другие темы.
Школа, которую они получали у Лузина, в высшей степени содействовала их занятиям
в другой области. Они получали великолепные навыки для последующего творчества:
и технические средства, и разработка интуиции математической. 
Они были готовы заняться многим другим. Вот они выходили из-под его, так сказать, влияния,
из сферы его влияния, и постепенно его ученики создавали собственные школы.
Вот Павел Сергеевич Александров и Павел Самуилович Урысон, два друга, — они ушли в топологию, 
Александр Яковлевич Хинчин, его ученик был тоже, — он ушел в основном в теорию вероятностей,
понимаете ли, другие уходили в теорию дифференциальных уравнений, и так далее, и так далее, и так далее.

{\bf И вот трагедия его заключалась в том, что постепенно близ него оставались лишь очень немногие, верные,
так сказать, ему и математике люди.}

Ну, они шли дальше, но... Такой был Дмитрий Евгеньевич Меньшов: он шел все в том же направлении, правда,
далеко и глубоко, в одной специальной области — теории тригонометрических рядов и ортогональных рядов.
Или была такая Нина Карловна Бари, одна из самых талантливых женщин-математиков, — 
она тоже занималась близкой к Меньшов тематикой и оставалась в ее пределах.
А иные уходили. Вот Лаврентьев, Михаил Алексеевич, будущий академик, — он перешел 
на прикладную тематику в значительной мере, сначала занялся некоторыми проблемами теории аналитических функций,
потом он перешел просто к вопросам прикладной математики, математической физики,
и, конечно, Лузину обязан многим, но он уже от него отошел.

И вот в этом была большая трагедия личная Лузина.
\end{quotation}

\sm

{\bf\punct В заключение.}
Люстерник был одним из гонителей Лузина в 1936 году, позже он, надо думать, изменил свою точку зрения.
Как мы видели, он -- единственный из нападавших -- в 50х годах участвовал в издании работ Лузина и мемориальных статей.
Бари -- источник вполне пролузинский, Лаврентьев относился к Лузину <<сложно>> (о чем ниже п. \ref{ss:plagiat}), но
с большим уважением.

Концы с концами сходятся. Молодые люди в самом деле перерастали Лузина (или достигали не меньшего уровня),
уходили из лузинской тематики (ТФДП плюс дескриптивная теория множеств), вооруженные владением ТФДП 
и навыками мышления, полученными в Лузитании. Этот уход не нравился Лузину,  на этой почве происходили столкновения. Люди становились самостоятельными учеными через эти конфликты...

Лаврентьев пишет о противоречиях 1925-1928, которые привели к развалу. Но мы знаем и о более ранних
конфликтах Лузина с Суслиным и Урысоном...

Бари и Люстерник в качестве дополнительной причины развала школы называют долгие поездки Лузина за границу
и его отрыв от молодежи. Лаврентьев вскользь упоминает еще одну причину, и речь о ней пойдет в следующем параграфе.

 \section{Судьбы теории множеств%
 \label{s:fate}}
 
 \COUNTERS

 \epigraph{
 	Впрочем, именно усталость от <<канторовского рая>> сыграла известную роль: огромное количество новых понятий без приложений и всегда без связи
 	с другими ветвями классической математики, многочисленные работы на явно искусственные темы и по крайней мере несоразмерные с их полезностью -- все это требовало крайней осторожности.} 
 	{	Лузин \cite{Lus-bologna}, 1927.}

 {\bf \punct  Проективные множества.%
 \label{ss:projective}} Мы прервали научную биографию лично Лузина на 1924 годе. 
 Открытие Суслина повлекло за собой исследование <<аналитических множеств>> и дополнений до аналитических множеств.
 В 1924 году Лузин вводит более широкий класс множеств -- <<проективные>>. Возьмем польское
 пространство $X$, фактически Лузин брал $\R$, а в своей монографии $\R\setminus\Q$.
 Возьмем всевозможные его конечные степени $X^n$ и всевозможные борелевские множества
 в них. Далее мы разрешаем две операции: дополнение до множества  и проектирование множества из
 $X^{k+l}$ на $X^{k}$. Хотелось бы добавить операцию объединения, но она не нужна, полученный класс замкнут относительно
 объединения автоматически. Кроме того, можно с самого начала брать не все борелевские множества, а лишь множества типа
 $F_\sigma$.
 
 Этот класс ({\it проективные множества}) является объединением возрастающей цепочки подклассов, нумеруемых числом операций 
 дополнения и проектирования, эта цепочка является строго возрастающей. 
 
 Лузин пытается понять, можно ли перенести на эти множества результаты об аналитических множествах и дополнениях
 до них. Он же придумывает изуверский метод решета \cite{Lusin-sieve}.
 
 Программа началась с заметки \cite{Lus-nonsolvability}, 1925, и мемуара \cite{Luzin-sbornik-1926}, 1926.

 Если я не ошибаюсь, новые ученики Лузина (Новиков, Селивановский, Л.~Келдыш, Ляпунов) 
 направляются на дескриптивную теорию множеств.
 
 \sm
 
 {\bf\punct Взгляд издалека.%
 \label{ss:sets}} В 1940 году Гёдель доказал, что отрицание кон\-ти\-ну\-ум-ги\-по\-те\-зы не доказуемо
 в системе аксиом Цермело--Френкеля \cite{God}. В 1963 году П.~Коэн установил, что континуум-гипотеза вообще независима от этих аксиом.
 Это было первое утверждение такого рода, а дальше результаты посыпались как из рога изобилия,
 оказалось, что почти все давно стоявшие вопросы теории множеств являются недоказуемыми и неопровержимыми.

 С другой стороны, когда-то давно была знаменитая книга Хаусдорфа <<Теория множеств>> \cite{Haus}, 1914, 1927, состоявшая из многочисленных
 положительных утверждений. Что происходило между 1927 и 1963 годом?
 
 \sm
 
 {\bf \punct Взгляды Лузина на неразрешимость, 1925--1930.%
 \label{ss:nerazreshimost}}
 Книга Лузина <<Лекции об аналитических множествах>>, 1930, была книгой с многочисленными  результатами
 по теории множеств в рамках
 аксиоматики Цермело--Френкеля,
 возможно, это была последняя  такая оригинальная книга.
 Работа, однако, еще содержит   умопостроения
 о связях между разными нерешенными проблемами.
И вот, что Лузин писал в Заключении:
 
\begin{quotation}
Автор этой книги склонен рассматривать проективные
множества как объекты, определение которых не может быть
вполне закончено: это чисто отрицательные понятия, которые
ускользают от всякого вида положительного определения....
{\bf Автор рассматривает как неразрешимый
вопрос о том, все ли проективные множества, измеримы или нет},
так как, на его взгляд, возможности самих способов определения проективных
множеств и меры в смысле Лебега несравнимы и, следовательно, лишены
логических взаимоотношений»
\end{quotation}
 
На самом деле это он говорил уже в своей первой первой опубликованной заметке
\cite{Lus-nonsolvability}, 1925,  по
проективным множествам:
\begin{quotation}
 «Усилия, которые я сделал, чтобы
решить этот вопрос, привели меня к следующему неожиданному заключению:
существует семейство эффективных множеств, имеющее отображение
на континуум [то есть континуальное семейство], таких, что неизвестно и никогда не будет известно, имеет
ли произвольное множество этого семейства (предполагаемое несчетным)
мощность континуума..., ни также —
измеримо ли оно. ... это семейство проективных
множеств»
\end{quotation}

Он также говорил, что для проективных множеств невозможно установить,
является ли любое такое множество объединением открытого множества и множества первой категории.

Цитируем Новикова и Келдыш:
\begin{quotation}
 На страницах своей книги Лузин высказывает мысль,
что трудности, связанные с проблемой о мощности 
CA-множеств [дополнений до аналитических множеств], носят принципиальный характер и что она не может
быть решена, исходя из принципов классической теории
множеств [вопрос в том, верно ли, что несчетное множество этого класса континуально].
\end{quotation}

Во всех этих случаях Лузин оказался прав, эти вопросы  независимы от аксиоматики 
теории множеств (ссылки и обсуждение см. у В.~А.~Успенского \cite{Usp}). 

По всей видимости, Лузин был первым человеком, который понял, что нерешаемые задачи теории
множеств в действительности являются неразрешимыми.
Продвижение в этих вопросах началось вскоре после смерти Лузина:
в 1951 году Новиков показал, что измеримость по Лебегу проекций $CA$-множеств не доказуема.
Впрочем доказательства независимости появились уже после работы Коэна.

Что касается положительной части, то Лузин с учениками продвинулись до пределов возможного.
Цитируем слова В.~А.~Успенского \cite{Usp} о работах В.~Г.~Кановея:
\begin{quotation}
\dots
Результаты в отношении двух последних проблем близки к тому, чтобы их
можно было назвать уникальными в дескриптивной теории: кажется, еще
никому не удавалось обнаружить разрешимую проблему среди серьезных
проблем дескриптивной теории, обсуждавшихся, но оставленных открытыми
в изысканиях классиков — Н.~Н.~Лузина и П.~С.~Новикова.
\end{quotation}

{\bf \punct Континуум-гипотеза.%
\label{ss:continuum}} Сначала об истории со Шнирельманом,
который является важным героем нашего дальнейшего изложения. Вот как рассказывает Люстерник:
\begin{quotation}
 Несколько  необычным  был  прием  у  Лузина  впервые  пришедшего  к  нему  
совсем  юного  Льва  Генриховича  Шнирельмана.  Николай  Николаевич  посадил  его  и  сказал:
«Сидите,  думайте  —  я  буду  смотреть  на  вас».  Оказалось,  
ему  снился  сон,  что  придет  мальчик   с  анкетными  данными  Л.~Г.   и  решит  
континуум-проблему.  Это уже из  области  иррационального  и  «Достоевского»...  
\end{quotation}

Теперь Понтрягин:
\begin{quotation}
Шнирельман рассказывал мне, что Лузин едва не загубил его как математика в самом начале 
его пребывания в университете. Лузин читал на первом курсе «Высшую алгебру». 
Хотя это не была его специальность, но он делал это для привлечения к себе студентов.
Лузин обратил внимание на Шнирельмана и предложил ему заняться решением континуум-проблемы. 
При этом он сказал: «Бросьте все лекции, ничему не учитесь и только думайте об этой проблеме».
Шнирельман, конечно, ничего не мог придумать по континуум-проблеме, а занятия он прекратил на целый год.
При встречах Лузин говорил ему: «Ну, что? Вы думаете? Думайте! Думайте!» Шнирельман не смел сказать, что он не знает, что думать. 
Занятия в университете он прекратил на целый год и с большим трудом вошёл потом в курс нормального обучения. 
 \end{quotation}
 
Речь, видимо, идет о 1921 годе, когда Лузин еще не осознал неразрешимость задач теории множеств.
Сам он над континуум-гипотезой много размышлял
и активно искал разные к ней разные подходы, и, не надо себе представлять дело так,
что Шнирельман пытался решить ее в одиночку.

А вот, что Лузин говорил на Первом Всероссийском съезде математиков (опубликовано в 1928 году)

\begin{quotation}
 Так прошло около двух десятков лет. За это время произошло следующее:
не только не выяснилась математическая природа рассуждения Zегmelо, но,
напротив того, к проблеме установления места continiium'a на алефической шкале
были прибавлены еще и другие проблемы, частью связанные с трансфинитными
числами, частью свободные от них — проблемы, вся обстановка которых
представляет ту же самую степень безнадежности в смысле разрешения конечным
рассуждением, каковая имеется в вышеуказанной проблеме об алефе continuum'a.
При этом нисколько не смягчившемся положении дела вполне следовало ожидать
возобновления спора. Не нужно думать, что здесь доминирующую роль играет
психологический момент, связанный с неудавшимися попытками разрешения и
 немногими имеющимися в науке прецедентами <<доказательств невозможности>>.
Всякая наука имеет свои собственные временно недоступные проблемы, решения
которых приходится иногда ждать веками — теория чисел, например. Но 
обстановка таких проблем всегда в высшей степени сложная, и внимательный глаз, 
исследующий её, вскоре подмечает многочисленные нити, тянущиеся от
исследуемой проблемы к многим вопросам и проблемам — решенным или нет — самых
разнообразных математических областей, на первый взгляд совершенно чуждых
области рассматриваемой проблемы. {\bf Ничего похожего нет в обсуждаемой 
обстановке: рассматриваемые проблемы теории функций%
\footnote{Лузин имеет тут в виду дескриптивную теорию множеств} характеризуются крайней
простотой и полной изолированностью обстановки. Нет никаких нитей,
связывающих их с проблемами вне теории функций: каждая из этих проблем
формулирована в терминах теории функций и должна решиться в этих же терминах и
методами и рассуждениями только теории функций}. Решение должно произойти
внутри области почти простыми средствами элементарной логики, причем число
испытуемых комбинаций обычно весьма ограничено. {\bf И если такого решения уж
вскоре же не находится — то в этом и повод к тревоге за прочность основ области.}
\end{quotation}

А вот, что Лузин писал в 1935 году \cite{Luzin-descript}:
\begin{quotation}
	В математической литературе не раз проскальзывали 
	намеки на то, что возможны и другие гипотезы континуума.
	Среди них наиболее интересной мне представляется та, 
	которая выражается алефическим равенством
	$$
	2^{\aleph_0}=2^{\aleph_1}
	$$
	
	Мы не станем доискиваться, с каким именем связано
	самое первое начертание этого равенства автором, который,
	действительно, серьезно мыслил его возможность. Мы 
	назовем его пока <<второй гипотезой континуума>>.
	
	Вторая гипотеза континуума мне представляется в той
	же самой степени изъятой от противоречий, как это имеет
	место и для первой гипотезы континуума. Мне думается,
	что со временем прогресс теории Hilbert'a будет настолько
	значительным, что позволит с успехом атаковать 
	доказательство непротиворечивости и этой второй гипотезы 
	континуума.
	
	А тогда пред нами предстанет необходимость выбирать
	между различными гипотезами континуума, в равной мере
	изъятыми от противоречия. И этот выбор, вне всякого 
	сомнения, будет продиктован одним только наблюдением фактов.
	Мы видим, насколько были проницательны слова E. Borel'я, писавшего:
	<<следует различать подлинную математику от чисто словесных логических спекуляций, в которых
	озабочены лишь одним совершенно отрицательным 
	качеством: свободой от словесного противоречия>> (Теория меры
	и теория интегрирования, \S 3, 1914).
\end{quotation}

Понятно, что две континуум-гипотезы находятся в противоречии друг с другом,
и тем самым Лузин считал, что они обе независимы от стандартной системы аксиом.

\sm

{\bf\punct И снова о кризисе школы.} Так или иначе, с середины 20х годов 
Лузин
понимал и открыто говорил, что задачи теории множеств, которыми он занимался,
по большей части неразрешимы. Если в 1916-1922гг он мог увлекать молодых
людей необозримыми далями и выводить их в люди на задачах с простыми формулировками,
причем задачах относительно простых, то теперь ситуация менялась на в точности противоположную.
Удивительно скорее, что молодые люди при Лузине еще оставались, хоть в меньшем
числе. Трудности нарастали на глазах,  сам Лузин и его программа к 30ому году
находились в кризисе. Понятно, что и интерес к его деятельности
со стороны угасал. Окружающие видели много более интересных
и содержательных тем для исследования. Новиков продолжил движение вперед и в итоге прорвался в 
математическую логику. Но это уже было в начале 40х годов.

\sm

{\bf\punct Чего хотел Лузин?%
	\label{ss:intention}}
Еще раз цитируем Заключение к <<Лекциям по аналитическим множествам>>:
\begin{quotation}
	Но если допускать все множества, измеримые B, то 
	необходимо допустить проективные множества. Следовательно,
	если желать ограничивать математический анализ лишь 
	изучением вполне законченных объектов и вполне определенных
	взаимоотношений, то нужно пожертвовать некоторыми 
	множествами, измеримыми В [борелевскими множествами], и даже некоторыми 
	иррациональными числами. В конечном итоге вполне определимых 
	иррациональных чисел имеется лишь счетное множество, хотя их
	перенумерование и не может быть осуществлено при помощи
	математического закона. Таким образом, арифметический
	континуум заведомо содержит неопределимые точки.... 
	
	В конечном счете вопрос будет разрешен определенно лишь
	усилиями научной мысли, т. е. наблюдением математических
	фактов. Философские рассмотрения служат лишь для того,
	чтобы отличить истинно плодотворное направление от
	бесконечного множества других. Только два случая возможны:
	Или дальнейшие исследования приведут когда-нибудь к 
	точным соотношениям между проективными множествами, а также
	к полному решению вопросов относительно меры, категории
	и мощности этих множеств. С этого момента проективные
	множества приобретут в математике право гражданства наравне
	е наиболее классическими из множеств, измеримых В.....
	
{\bf	Или указанные проблемы из теории проективных множеств
	останутся навсегда нерешенными и к ним добавится 
	множество новых проблем, столь же естественных и столь же 
	недоступных. В этом случае ясно, что пришло время произвести
	реформу в наших идеях об арифметическом континууме.}
\end{quotation}

Вопрос о <<реформе наших идей об арифметическом континууме>>
встает совершенно естественно при чтении Лузина. Пытался ли он сам
ее провести?

Цитируем Келдыш, 1974:
\begin{quotation}
	Остается еще добавить, что в ту эпоху, когда Н.~Н.~Лузин занимался
	дескриптивной теорией множеств, вокруг него объединялась группа его
	учеников, активно работавших в этой области. Н.~Н.~Лузин умел исключительно
	увлекательно говорить и при обсуждении иногда даже, казалось бы,
	самых простых вопросов указывать очень интересные далекие перспективы
	и связи с важными проблемами. Новизна его идей, важность их, особенно
	вопросов о существовании неразрешимых проблем, не только для теории
	множеств, но для всей математики, создавала среди окружавшей его молодежи
	обстановку большого энтузиазма. У всех, работавших в его группе учеников,
	было ощущение участия в очень важном и исключительно интересном
	деле. Задача, поставленная Н.~Н.~Лузиным, как уже было сказано выше,
	сводилась к следующему: решить те проблемы дескриптивной теории множеств,
	которые можно решить средствами теории множеств, и сформулировать
	те, для которых это невозможно................... 
	%И сейчас уже можно сказать,
	%что в течение 20-х—30-х годов эта программа действительно была выполнена.
\end{quotation}

\begin{quotation}
	В течение
	последующих пятнадцати лет Н. Н. Лузиным с группой его учеников эта
	программа была выполнена, хотя следует отметить, что вторая ее часть,
	касающаяся выделения неразрешимых проблем, тогда носила чисто эвристический
	характер. В то время еще не существовало методов, позволяющих
	устанавливать недоказуемость математических утверждений.
\end{quotation}

Но, что касается вопроса реформы, то  Новиков и Келдыш, \cite{KeldNov}
давали отрицательный ответ на заданный вопрос:
\begin{quotation}
	Подойдя в своих работах к исследованию границ применимости принципов теории множеств,
	Н.~Н. живо интересовался вопросами оснований
	математики. Сам он никогда не работал в этой области, но внимательно 
	следил за ней и критически оценивал то, что в ней происходило.	
\end{quotation}

Сомнения на этот счет могут оставаться. Цитата из <<Лекций об аналитических множествах>>:
\begin{quotation}
	По-видимому, при современном
	состоянии науки является преждевременным нападать на дедекиндовскую
	теорию иррациональных чисел, если только,
	желать, чтобы это нападение оказалось плодотворным, а 
	именно дало новые положительные результаты, ускользающие от
	нас в области этой теории.
	Таким образом, мы ограничимся тем, что примем ее и
	будем рассматривать как промежуточный инструмент, считая
	возможным в дальнейшем указывать на некоторые трудности
	этой теории....
	
	Целью теории множеств является вопрос чрезвычайной
	важности: можно или нет рассматривать линейную 
	протяженность атомистическим образом как множество точек; вопрос
	этот, кстати, уже не нов и восходит к эллинам.
\end{quotation}

Но это опять вопрос о поиске реформы...
 Читатель может
воспринимать вопрос Лузина как странный, но это во многом вопрос  восприятия математики
и математического опыта. Понимание, что в теории множеств задачи обычно неразрешимы,
естественно должно было повлечь вопрос о том, <<правильна>> ли сама теория множеств.
 Вроде бы мы можем назначать 
 произвольные ответы на некоторые вопросы (например ставить континуум в разные места
 <<алефической шкалы>>), что необычно для прочей математики, а с другой стороны в прочей
 математике ничего не будет меняться от того, какой мы ответ назначим.
 С другой стороны,
 почти во всех реальных математических
построениях с континуумом на нем присутствует борелевская структура,
причем она одна и та же, откуда бы мы континуум не брали; она уважается (почти всегда) 
при всех реальных математических манипуляциях. То есть выходит, что континуум
-- больше, чем множество.
С другой стороны, борелевская структура - объект малопонятный (в этом смысле
она  хуже континуума,
который хоть кажется понятным, если не задумываться над опасными вопросами).
И - увы! - Лузин с Суслиным обнаружили, что борелевская структура сохраняется не всегда...

Опять цитируем Келдыш:
 \begin{quotation}
 Он относился отрицательно к теории Брауэра -- интуиционизму как к теории
  разрушительной. Однако и теорию Гильберта, и в частности её подход к решению
 труднейших задач теории множеств он не считал полноценной
 [на теорию Гильберта Лузин высказывал и оптимистические взгляды].
 
 Н.~Н. считал, что доказательство непротиворечивости континуум-гипотезы не является удовлетворительным решением континуум-проблемы. Он
 считал, что целью теории множеств является накопление фактов, которые
 позволили бы выбрать определённую гипотезу. Арифметический континуум, —
 говорил он, — должен быть единственным. 
 \end{quotation}

 Замечу, что в отношении теоретико-множественного скептицизма Лузин был отчасти
 близок  к другим сомневающимся начала XX века, Лебегу, Бэру, Борелю,
 Г.~Вейлю [Hermann Weyl], а также раннему интуиционизму. Основными предметами сомнения были аксиома выбора,
 правомерность операции множества подмножеств, и то, что континуум можно рассматривать
 как множество.

 Вот еще цитата из Лузина:
\begin{quotation}
Мы не будем стараться дать определение слову {\bf назвать}, нам кажется,
 что это достаточно примитивное
понятие для того, чтобы
его определение было по меньшей мере бесполезно; можно только постараться объяснить при помощи
примеров смысл этого слова. Назвать функцию, множество это значит однозначно указать
эту функцию индивидуальным образом без каких бы то ни было возможных двусмысленностей.

Одна из основных идей, которыми мы обязаны Лебегу, - это точное различение двух понятий: 
называемое и неназываемое множество....

{\bf Совокупность всех точечных множеств не имеет в объективной науке, каковой являются классические части
математики, другого существования как чисто словесное; следовательно невозможно ее рассматривать
как истинный математический объект, который может быть включен в рассуждения.

... Нужно включить в область незаконного самую совокупность действительных чисел}... Одна из наиболее
важных проблем теории функций действительной переменной - в преобразовании совокупности чисел 
в законную совокупность... или в том, чтобы назвать при помощи конечного определения класс чисел,
инвариантных относительно всех преобразований классического анализа и алгебры, 
этот класс чисел будет составлять практический континуум, который используют математики.
\end{quotation}

Похоже, что у Лузина была идея ужесточения научного законодательства (вслед за Брауэром%
\footnote{Последняя из цитированных сентенций близка к Борелю. Марксистские философы того
времени не без некоторых оснований объединяли Бореля и Лузина под вывеской
<<эффективизм>>.}).
При этом классическая математика не должна была пострадать
(это вполне ясно из лузинских речей), а с другой стороны,
у Лузина были весьма высокие эстетические запросы (что явно не оговоривалось, но
тоже ясно). Краевые условия
на реформу крайне жесткие,
и попытки реформы не последовало%
\footnote{Во всяком случае ее не видно из опубликованных
работ Лузина,  я не видел и не слышал никаких известий о такой попытке.}.
 А, возможно, -- эта точка зрения скользит в публикациях Лузина -- реформа считалась
 им преждевременной.
 
 \sm

{\bf \punct Скептицизм в отношении натурального ряда.%
\label{ss:N}}
У Лузина были и более радикальные сомнения.
Например, Лузин(см. письма Выгодскому \cite{Demidov-pismo})  считал возможным введение
актуальных бесконечно малых, сейчас мы знаем, что это в самом деле возможно.

Высказывал он и скептицизм в отношении натурального ряда.
Приведем цитату из \cite{Luz-fund}, эта статья -- часть письма Куратовскому:
\begin{quotation}
	\dots Если я трачу время на рассмотрение этих вещей, то не потому, что считаю их действительно серьезными, а потому, что через множество чисто словесных существований, слишком легких, чтобы принимать их всерьез, я вижу слабый свет настоящей интуиции, могущим привести нас к совершенно неожиданным фактам, которые мы обнаружим, если следовать другому пути.\dots
	
	Несмотря ни на что, я не могу рассматривать как {\it данное} множество целых положительных чисел, потому что сама идея актуальной бесконечности мне кажется мало естественной, чтобы рассматривать ее в себе. Однако критика бесконечности ужасно трудна, и мы начинаем понимать, какие трудности видел Борель, когда он писал:
	
	<<Кажется неоспоримым, что математики делают или думают, что делают, из этой неограниченной последовательности
	$$
	\text{1,2, 3, 4, 5, 6,\dots}
	$$
	совершенно ясную идею>>.
	
	Слова <<или думают, что делают>> очень характерны для точки зрения натуралиста.
	
	Представляется трудным пройти мимо новых идей, к которым нас приводят кванты и исследования Дэ Ситтера [Willem De Sitter] и Лемэтра [Georges Lema\^itre]. Фундаментальная проблема состоит в том, чтобы выяснить, является последовательность положительных чисел {\it вполне} объективной? Кажется, что она почти объективна, и что имеются следы несомненной субъективности, такой, что нельзя говорить о последовательности целых положительных чисел всегда, во всех случаях, в одном и том же смысле. Однако в данный момент слишком преждевременно ставить жгучую проблему о {\it единственности последовательности целых положительных чисел} и говорить о конечных недостижимых числах, отправляясь от 1. Современная теория чисел весьма униформизирована методом полной индукции. И тем не менее наука о числах благополучно существовала в XVIIв. В то время существовали другие умственные привычки, которые не дошли до нас.
	
	С точки зрения натуралиста, трансфинитные числа являются только очень большими конечными числами. Исследования о множествах -- классификация Бэра и другие -- приобретают совершенно новый аспект. Это исследования о конечных выборах, но надо принимать указанную предосторожность относительно слова <<конечное>>\dots
	
	{\bf Я не считаю бесполезными усилия в изучении классов Бэра, аналитических и проективных множеств: они не являются в точности объектами, имеющими самостоятельную ценность, но движениями мысли, которые,будучи перенесенными в область целых положительных, позволят нам {\it по аналогии} открывать новые феномены, недостижимые другими методами.} 
	
Но я Вам уже писал об этих вещах.

Москва, 14 мая 1933г.
\end{quotation}

Цитата относится к 1933г. 
Ниже в п. \ref{ss:kolman-tsagi}, из несколько необычного исторического источника мы узнаем о том, что  Лузин был уверен в неполноте  арифметики на несколько лет раньше.

Автор не в состоянии понять, что именно Лузин имеет здесь в виду, да и сам он говорит
о <<слабом свете>>. Может, он в самом деле хотел переделать рассуждения о проективных
множествах в рассуждения о натуральных числах.

 Эти сомнения так или иначе запомнились в московском математическом
сообществе, см., например, статью П.~К.~Рашевского \cite{Rashevski} <<О догмате натурального числа>>. К сожалению, через 10 лет после смерти Лузина на эту тему началась известная клоунада А.~С.~Есенина-Вольпина, см. \cite{Neretin-1968}, которая в итоге (уже позже и по иным причинам) имела 
для московской математики печальные последствия. Интересно, что ее  поддержала
 гонительница Лузина С.~А.~Яновская (которая едва ли была в состоянии понять суть дела)
и, до какой-то степени, главный противник Лузина Александров (который, к сожалению, в дело вникать не стал).

\sm

\section{Философия и математика%
\label{s:filosofiya}}

\COUNTERS

Отчасти, цель этого раздела - комментарии к книге Graham, Kantor
<<Naming infinity>> \cite{GrK}, 2009 (там 
голословно   излагается историческая схема развития дескриптивной 
теории множеств,
 но так как это целая книга, причем по-анг\-лий\-ски,
приходится ее комментировать). С другой стороны, это
повод упомянуть дискуссии 20--30х годов об основаниях математики, не последним из участников которых был Лузин.

\sm

{\bf \punct Лузин, философия, философствования и философы.%
\label{ss:philosophy}}
Об интересе Лузина к философии писалось выше, его философствующий взгляд на 
математику виден из многих цитат предыдущего параграфа.

Приведем длинную цитату из предисловия Лебега \cite{Lebeg} к книге Лузина 1930г.
Тут интересно и самом Лузине, и об отношениях математики и философии:
\begin{quotation}
После первых крупных успехов теории множеств философы и математики
и сочли, что настал момент протянуть друг другу руку через разделяющую
их глубокую пропасть. Завязавшийся разговор с самого начала походил на
игру в бессвязные слова; считалось, что тут всего лишь минутное недоразумение, еще одно усилие и взаимопонимание будет достигнуто. Но вспомнили
о Зеноне Элейском
и о парадоксе лжеца.

В мои студенческие годы мы охотно атаковали за чашкой кофе глобальные
 идеи: беседа разгоралась и казалась нескончаемой, когда один из нас
восклицал: «Вот ты, для начала, существуешь ли ты? Я говорю «ты» для
удобства, ведь существую один я...». Но тут мы осознавали, что надо идти
работать, и расставались до завтра. С появлением Зенона и лжеца математики
поняли, что им следует удалиться; но на следующий день они не вернулись,
и не было кофе. Некоторые даже поклялись, что больше они на этом не
попадутся; они приложили старания к тому, чтобы не сказать ничего такого,
что позволило бы втянуть их в философскую дискуссию,— сожалея единственно
 о том, что им не удается, подобно Эрмиту, спрятать свои идеи в уравнения. Эрмит, мэтр, оставивший столь глубокий след в науке, которого
сумели процитировать, лишь извлекая фразы из частной переписки.

Глубоко было разочарование среди математиков. Обратились к
мужам, которых страстное желание все знать сделало философами и
которые были одарены столь блистательно, что казались на самом деле 
способными все объять и все постичь; но заставить их выйти за рамки традиционной философии не удалось, и потому при создании научной философии,
пригодной для науки, ученые положительно должны были рассчитывать
только на себя самих. Для тех, кто долго и упорно трудился над техникой,
это явилось бы делом, завершающим карьеру.

Однако постоянно стараться говорить только на техническом языке —
разве значит хорошо подготовиться к этому делу? Не рискует ли проявить
убогость мысли тот, кто, опасаясь дискуссий, никогда не говорит ничего
хоть сколько-нибудь общего? Человек науки, занявший слишком сдержанную,
недоверчивую позицию, никак не может надеяться, что в подходящий
момент он будет содействовать созданию столь желанной для него новой
философии.

К счастью, здесь, у г-на Лузина, другая позиция, обнадеживающая.
У него тоже, может быть, есть некоторое недоверие к философии, но у него
нет перед ней страха; у него есть философские интересы, и он в них сознает
ся. Математические требования и философские требования у него постоянно
соединены, даже, можно сказать, сплавлены. Хотя его книга — сочинение
по математике, написанное математиком для математиков, почти на каждой
ее странице отчетливо проступает эта тесная связь философских и математических мыслей, что придает монографии исключительную значительность
и совершенно необычайную привлекательность. 
\end{quotation}

Из письма Лузина к Выгодскому (речь идет где-то о 1906-1907 годе):
\begin{quotation}
	Лишь после этих 1/2 года на историко-философском факультете,
	где я слушал всех известных тогда философов и не вынес,
	почему-то страсти к их рассуждениям, я возвратился к математике,
	и приняв предложение Егорова, принялся за пополнение математического образования,
	не спрашивая советов и руководствуясь, как и раньше, случаем.
\end{quotation}	

Но интерес Лузина к философии сохранился, на что есть много прямых и косвенных
свидетельств%
\footnote{\label{fo:uchebnik}Лузин -- Д. М. Петрушевскому, 26.8.1939 \cite{Erm2}:
\newline
{\it 
Десять дней посвятил на учебник, который Учпедгиз все равно хочет выпустить и к которому я не только равнодушен, но и который вызывает у меня неприязненное чувство. Многое, очень многое, касающееся философского обоснования пространства, чисел, величины, времени, количества, бесконечности, категории <<конечного>>, движения, изменения, счета и т.д. просто выкинули или сделали такие примечания, которые мне и в голову никогда бы не пришли... Поэтому я глубоко равнодушен к дальнейшим судьбам книжки, хотя ей и предсказывают некоторое педагогическое будущее.}
\newline
Речь идет об учебном пособии <<Теория функций действительного переменного>>
для педвузов. Вообще-то в учебных пособиях по естественно-научным предметам не принято философствовать и точку зрения Учпедгиза можно понять. Кстати, неприязнь к книге после долгой работы над ней -- дело тоже обычное.}.
 
 \sm
 
 {\bf \punct Академик-философ.}
Забавно, что в  1929 году Лузин был выбран в Академию по кафедре философии.
Это было результатом интриги Крылова (который возглавлял выборы по
только что объявленной кафедре философии). Предполагалось, что будут два философа, один
- естествоиспытатель, другой - гуманитарий.
В тот момент в АН СССР было лишь два академика-математика,
Крылов и Я.~В.~Успенский. Цитируем Стенограмму:
 \begin{quotation}
 КРЫЛОВ. Я просил бы приобщить к делу все материалы по представлению Лузина в
 	академики.
 	Я помню, что представленные работы Лузина не казались мне чисто математическими
 	 работами, а скорее работами философскими. Но мы провели его по кафедре математики. Успенский дал тот же отзыв.
 	 
 	Затаенная мысль у нас была такая в Комиссии, что оттого, что мы изберем его по
 	кафедре философии, он не перестанет быть математиком, и, следовательно, мы получим
 	лишнего математика, а философов и без того там довольно.
 	
 БЕРНШТЕЙН. А как произошел переход его с кафедры философии на кафедру математики?
 	
 КРЫЛОВ. А когда Яков Викторович [Успенский] уехал в Германию
 	и отказался от звания академика, освободилась вакансия, и постановлением Общего собрания он был переведен.
 \end{quotation}

Напомним, что это были экстремальные академические выборы 1929г. На второе место философа
(гуманитария) шел тогда А.~М.~Деборин,
продавливаемый верховной властью.

По математике академиками тогда были выбраны академиками Бернштейн, Виноградов, Н.~М.~Крылов,
член-корами Александров,
Б.~Н.~Делоне,
Н.~Г.~Чеботарёв.

% Щелкачев
 
% Тахо-Годи
 
 {\small 
{\bf\punct Книга Грэхема и Кантора.%
\label{ss:graham}}
Теперь о работах Грэхема [Loren Graham] и Кантора [Michel Kantor],
статье <<Два подхода к оценке математики как феномена культуры: 
Франция и Россия, 1890-1930г>>
и книге <<Имена бесконечности>>. Авторы утверждают, что дескриптивная
теория множеств в России была связана с имяславием.

Имяславие -- православное религиозное течение, говоря о нем,
обычно приводят цитаты типа:
\begin{quotation}
«Бог — везде; и, конечно, Он находится и в Своём имени».

\sm

 «Имя Божие есть Бог; но Бог не есть имя. Существо Божие выше энергии Его, хотя эта энергия выражает существо Имени Бога.»
 \end{quotation}
 
 В том, что  такое имяславие на самом деле, надо разбираться по специальной богословской литературе
 (что само по себе представляет непростую задачу).  В число сторонников этого  течения 
 входили, постоянно или временно, ряд известных философов,  
 Владимир Эрн, Павел Флоренский, Сергий Булгаков, Алексей Лосев. Имяславцами
 были также математик Егоров и (временно) механик Бухгольц. Течение 
 зародилось в 1908-1909 году, его поддержали русские монахи с горы Афон.
 В 1913 году Синод осудил имяславцев. Русские военные корабли 
 вывезли монахов-имяславцев в Россию.
Однако у имяславцев нашлись влиятельные сторонники,
в итоге имяславие не было осуждено Церковью, а вывезенные
монахи начали пропаганду (или проповедь, кому что больше нравится) имяславия.

Теперь цитируем книгу Грэхэма и Кантора, а затем их же статью \cite{GrK0}.
По ходу мы делаем несколько замечаний, но обсуждение по сути -- ниже. 
\begin{quotation}
Теорию множеств, которую блистательно разрабатывали во Франции, постиг сильнейший кризис, после чего на сцену вышли русские. Мы расскажем, как две разные духовности, связанные с разными культурами, привели к противоположным результатам: с одной стороны французскому скептицизму и неуверенности, с другой - к русскому творчеству и прогрессу. Главная идея книги заключается в том, что своеобразная религиозная ересь%
\footnote{На всякий случай отметим, что православная церковь  вердикта <<ересь>> в отношении
	имяславие не выносила.} оказалась своеобразными акушерскими щипцами, вызволившими на свет новую область математики.
\end{quotation}

Продолжаем цитировать:
\begin{quotation}
	В свои зрелые и профессионально активные годы все трое -- Егоров, Флоренский и Лузин -- были глубоко религиозными людьми \dots В конце концов Флоренский и Егоров были арестованы коммунистическими властями по обвинению в смешении математики и религии и закончили жизнь в неволе%
	\footnote{Дела эти исследовались, Егоров (о его деле ниже п. \ref{ss:losev}) и Флоренский обвинялись в участии в контрреволюционных организациях, а вовсе не в религиозно-математическом смешении. Подавляющее большинство  участников дела имяславцев \ref{ss:losev} никакого отношения к математике не имело, да едва ли математика интересовала тогдашнее ОГПУ. Егоров в итоге  был подвергнут административной высылке в Казань. <<Суд>> над Лузиным подробно обсуждается в следующих двух параграфах,
		там тоже ни о каком религиозно-математическом смешении и  помина нет.}. (Одна из жестоких шуток истории заключается в том, что обвинение Флоренского и Егорова в смешении математики и религии было справедливым, и несмотря на противоположное предположение коммунистов, смешение это оказалось удивительно плодотворным для математики.) Третий из них едва избежал заключения, хотя и подвергся <<суду>> за идеологические отклонения.
		Из всех троих [т.е., Егоров, Лузин и Флоренский] именно священнику Павлу Флоренского
		предстояло предать новому синтезу математики и религии наиболее яркую форму,
		сыгравшую огромную роль в исследовательских математических работах Егорова,
		Лузина и их учеников...

		\sm
		
Так как непрерывные функции описывают <<детерминированные>> процессы%
\footnote{Отметим, что данная посылка ошибочна, даже если не иметь в виду нигде недифференцируемую функцию Вейерштрасса.}, то экспансию детерминизма в психологию, социологию и религию Флоренский рассматривал как разрушительный итог временного увлечения такими функциями в математике... Флоренский воззвал к <<заре нового, прерывного миросозерцания>> и призвал своих коллег, таких как Лузин и Егоров, развить этот новый подход, который как он думал, совместить математику, религию и философию%
\footnote{Непонятно, как можно называть Егорова коллегой Флоренского. Продвинутый курс ТФДП в Московском университете с 1900-01 года читал Млодзеевский. Разрывными функциями до того занимался Бугаев, кстати научный руководитель Флоренского. Егоров в 1902-1903гг. год пробыл в Берлине, Париже (где, в частности, слушал Лебега), Геттингене. Он еще несколько раз выезжал за границу.}.		
\end{quotation}

\begin{quotation}
Во время большевистской революции Павел Флоренский жил в Сергиеве
Посаде, городке неподалеку от Москвы, и религиозно и интеллектуально был
близок к диссидентам из числа имяславцев. Он формулировал свои идеи на
математическом языке и делился ими с Лузиным и Егоровым. В 1920-е гг. в
Москве действовал «имяславческий кружок», в котором соединялись идеи религиозных диссидентов и математические концепции. Кружок состоял 
из пятнадцати или шестнадцати философов, математиков и религиозных мыслителей.
Иногда они собирались на квартире Егорова, и на нескольких из этих
встреч Флоренский представлял свои доклады. Здесь он излагал идею о том,
что «точка, в которой встречаются божественная и человеческая энергия,
есть «символ», который больше чем он сам»...

Флоренский понимал, что благодаря имяславцам вопрос «именования» 
получил совершенно новое звучание. Дать имя чему бы то ни было значило
 породить новую сущность. Флоренский был убежден, что математика – продукт
свободного творчества человеческих существ и что она имеет религиозный
смысл. Люди могут использовать Свободную Волю и обратить ее в перспектив
ную математику и философию....

Развитие теории множеств стало для Флоренского драгоценным примером
того, как наименование и классификация могут приводить к крупным 
научным открытиям. «Множество» являлось просто присваиванием названий в соответствии с произвольной ментальной системой, а не распознаванием типов
онтологически сложившихся объектов. Когда математик создает «множество»,
 давая ему название, он дает жизнь новой математической сущности%
 \footnote{Сравните со скептицизмом Лузина в отношении теории
 	множеств, <<множество чисто словесных существований, слишком легких, чтобы принимать их всерьез>>.}.
Наименование множеств является математическим актом, точно так же как,
с точки зрения имяславцев, именование Бога – актом религиозным.
Флоренский утверждал, что приходит время нового вида математики...

{\bf Таким образом коренное слово «имя» присутствовало в русском
языке одновременно в качестве математического термина, обозначавшего
новые виды множеств, и в качестве религиозного термина «имяславия» («имяславцы»). Действительно, большая часть работы Лузина над теорией множеств
включала изучение эффективных множеств («называемых множеств»). Для Флоренского это означало, что и религия, и математика двигались в одном и том же направлении...}

Перед Флоренским стоял вопрос аналогичный: «Возможно ли убедить себя в существовании Бога, не давая Ему определения?» Ответ для Флоренского и, позднее, для Егорова и Лузина был
следующим: акт именования сам по себе дает объекту жизнь. Таким образом
«именование» стало ключом и для религии и для математики. Имяславцы давали определение Богу, поклоняясь Его имени, математики давали определения
множествам, присваивая им имена...

Изучая теорию множеств,
французы хотели ясно идентифицировать ее философские, математические и
психологические составляющие и сохранять их по отдельности. В противоположность
 этому, Лузин и некоторые его друзья верили в то, что математика
связана с религией, но они не могли искренне говорить об этих связях в послереволюционный период из-за враждебного отношения к таким вопросам
советского окружения. Они знали о проблемах с властями, которые возникнут
немедленно, если взгляды, обсуждавшиеся на встречах «имяславского
кружка», станут известны.	
	\end{quotation}

На эти работы было несколько русских рецензий, А.~Геронимус, С.~С.~Демидов, А.~Н.~Паршин \cite{GDP},
 В.~А.~Бажанов \cite{Bazhanov}, В.~Н.~Катасонов \cite{Katasonov} (написана с историко-религиозной точки зрения%
\footnote{Избегая рассуждений на эту тему и предпочитая рассматривать
	вопрос со стороны математики, приведу две цитаты из этой статьи
	\newline
 <<\it Прежде всего, нужно сказать, что идея творения не имеет никакого отношения к имяславию. Имяславцы, призывая Бога, не творят Его, и нам не
 известен ни один из имяславцев, утверждавший нечто подобное... Поклоняясь
 имени Божьему, повторяя Иисусову молитву, имяславцы приближаются к Богу,
 как бы актуализируют Его присутствие для себя. Хотя подобное высказывание
 следует понимать скорее не как страдательность Бога, но как человеческое обращение к Богу... Но никогда речь не идет о творении.>>
 \newline
 <<Богословски достаточно грубое сближение идей имяславия с идеей творения
 проистекает у авторов рецензируемой книги из их желания связать имяславие и
 математические конструкции теории множеств.>>}),
Г.~Синкевич \cite{Sink}. В рецензиях указаны многочисленные несообразности, ошибки,
голословные утверждения и натяжки. Основная идея книги не выдерживает критики.
Следующий пункт содержат некоторые дополнительные доводы.

\sm

{\bf\punct Исторические несообразности.%
\label{ss:nesoobraznosti}} Нам говорят, что Флоренский, 
Егоров и Лузин были имяславцами, математиками и занимались дескриптивной теорией множеств. Это умопостроение ложно.

\sm

a) {\sc Флоренский}  был сторонником имяславия (или, по крайней мере,
высказывал к нему симпатии). Он учился на Московском
Физмате, и это повлияло на его мышление, 
часто его называют математиком,  
но никаких его научных математических работ (опубликованных или
рукописных) не известно. Есть философствования на основе математического  
образования, но это не профессиональная математика%
\footnote{В связи с теорией множеств у него была статья 1904г. <<О символах бесконечности>> \cite{Florenski}. Это что-то
вроде религиозно-философской популяризации  канторовской теории.  Там говорится о мощностях, их сравнении, о том, как можно строить
большие мощности (путем последовательного применения операции степени).
Сказано вскользь определение вполне упорядоченного множества и сказано несколько слов об ординалах (трансфинитных числах). Сейчас часто произносятся сентенции типа того,
что это было <<первое развернутое изложение теории трансфинитных
 множеств Кантора>> в России. Это
сильное преувеличение как в отношении развернутого
изложения, так же, как и сильное преувеличение в отношении слова 
<<первый>>, см. \cite{Bazhanov}. Стоит отметить, что сохранились 
записи лекций Млодзеевского 1902г., сделанные Флоренским, 
в которых последовательно
излагалась теория множеств  до теории ординалов включительно  \cite{Medvedev-Mlo}.}. Следов  научного влияния 
Флоренского на математиков 
тоже не видно (и вообще это было бы странным при отсутствии научных 
работ).

Отношение Лузина и Егорова к математическим сочинениям Флоренского
было  отрицательным.
Из письма Лузина жене от 3 августа 1908 \cite{Luz-Flor}:
\begin{quotation}
Как только он показал свои работы по математике –
опять старое зашевелилось во мне мнение: все его работы не имеют
цены в области математики. Намёки, красивые сравнения – что-то
упивающее и обещающее, дразнящее, манящее и
безрезультативное%
\footnote{Весьма понятная характеристика, если посмотреть, например
книжку Флоренского <<Мнимости в геометрии>>, 1922.}. И под конец я перестал понимать, что же такое
Флоренский? 	
	\end{quotation}
	
 А вот рассказ А.~Ф.~Лосева:
	 \begin{quotation}
	 	Далее [речь идет о работе <<О мнимостях в геометрии>>] отец Павел говорит... а если тело движется со скоростью, большей, чем скорость света? Тогда-то что получается?! Согласно этой формуле [для энергии] получится не единица, не нуль вместо подкоренного выражения, а мнимая величина...
	 	То есть [тело]	становится идеей, не имеющей объема....
	 	
И все же я знал хороших математиков, которые дрожали от ненависти при чтении этих страниц.

Хотя Флоренский там и говорит, что он советовался со своим университетским товарищем Николаем Николаевичем Лузиным, -- они вместе в 1904 году окончили математический факультет, товарищи по учебе. Так вот тот, пишет Флоренский, прочитал и одобрил.

Но я близко знал учителя академика Лузина -- это Дмитрий Федорович Егоров, крупнейший математик с мировым именем. Так вот он нервно-отрицательно относился к этим умозаключениям.

{\bf Когда я с ним беседовал -- а он всегда был корректный человек, джентльмен, -- при упоминании об этом рассуждении он напрягался и терял самообладание}%
\footnote{Вот что писал \cite{Kolman-Losev} о Лосеве Кольман, который упоминался в \S1
и о котором еще много будет говориться ниже.
\newline
{\it	То, что математика служила мракобесию не только в прошлом, но и оказывает ему услуги в настояшем, всем известно. Такой черносотенный идеалист как Лосев, в книге, вышедшей в Москве в 1927году (<<Античный космос и современная наука>>) ссылается на другого попа Флоренского, на его книгу <<Мнимость в геометрии>>. Сам от себя Лосев добавляет рассуждение о том, как теория относительности
	приводит неизбежно к платонизму. Вот основная формула теории относительности: известное смещение
	\newline
	$\beta=\sqrt {1-v^2/c^2}$
	\newline
	Здесь, говорит он, весь платонизм в этой формуле.}
	\newline
	Забавно, как подобные джентльменские дискуссии преломляются в современной литературе. Вот что пишет Тахо-Годи в биографии Лосева (книге, вообще, весьма качественной):
	\newline
{\it Книги А. Ф. Лосева были теснейшим образом связаны с современностью. Он писал не просто об античном космосе, но о достижениях современной науки, самых последних, наиболее интересных, но и опасных в 20-е годы, да и не только тогда (например, теория относительности Эйнштейна, знаменитая формула Лоренца, математические теории П.~А.~Флоренского).}	
	\newline
О, Русская Земля!
}.... 
\end{quotation}

Между двумя отрывками есть известное противоречие, но мы не можем знать,
что на самом деле говорил Лузин Флоренскому, а  из многих других источников знаем
о чрезмерной вежливости Лузина. А отношение Егорова к математике Флоренского описано 
ярко.
	
Теперь мнения Егорова и Лузина   об отношениях математики и метафизики.
	Из воспоминаний Н.~М.~Бескина (описываются 20е годы):
\begin{quotation}
Н.~Н.~Лузин рассказывал нам об интересном человеке П.~А.~Флоренском. Они окончили университет одновременно... Н.~Н.~Лузин считал Флоренского талантливым математиком. Он сожалел, что П.~А.~Флоренский,  окончив университет, отверг предложение остаться в аспирантуре по математике и посвятил себя богословию. Он считал книгу Флоренского <<Столп и утверждение истины>> метафизической... Следующее высказывания Лузина я помню почти дословно

--- Метафизика отделена от науки непроницаемой стеной. Между ними нет никаких связей. В метафизике может быть застой, а может быть расцвет и бурная кипучая деятельность. Ни в том, ни в другом случае, это не окажет ни малейшего влияния на процессы, происходящие в науке.
\end{quotation}

Письмо Егорова Лузину от 27.2(12.3)1906:
\begin{quotation}
	... Вы говорите о перемещении центров в миросозерцании. Мне думается, что научная работа, ее направление и даже вкусы в этой области не должны бы зависеть от миросозерцания. Миросозерцание - само по себе, а наука - сама по себе!
\end{quotation}
	
	б) {\sc Егоров}. Он был имяславцем, был математиком, но работ по дескриптивной
	теории множеств у него не было. По ТФДП он опубликовал лишь две 
 заметки в Comptes Rendus
	\cite{Egorov1911}, \cite{Egorov1912} 1911 и 1912 года
	(вторая заметка - реакция на теорему Лузина о существовании первообразной),
	обе вышли до вывоза имяславцев с Афона (и до начала распространения имяславческого учения в России).
	Список работ Егорова есть в 
	\cite{KuzEgor}. 

Егоров отделял научную деятельность от общественных и философских взглядов.
В 20-30е годы ультралевые  склоняли его налево и направо в связи с его реакционными
общественными взглядами и (открытым)
 отрицательным  отношением к марксизму (много таких
сентенций есть в Стенограмме). Но отношение к нему как к математику оставалось уважительным
и никто не обвинял его в  смешении общественных взглядов и математики (разве что Хворостин и Райков
обвиняли его в занятиях чистой математикой, что с их зрения было предосудительно; Яновская громила его за отсутствие марксисткой диалектики
на заседаниях мат.общества; Люстерник с сотоварищами обвиняли Егорова  в отсутствии 
политических докладов на Матобществе и недостаточный вклад математики в строительство социализма).

\sm
 
 в) {\sc Лузин}. Он был основной фигурой дескриптивной теории множеств в 1920х
 годах. Он был православным (о чем есть  многочисленные известия из разных периодов его жизни).
 
 Но утверждения о его  имяславии (или склонности к имяславию) и влиянии религиозного мистицизма на его
 математические работы    ни на чем не основано.

Надо иметь в виду, что объем <<Собрания сочинений>> Лузина 2000 страниц,
есть сотни страниц работ, в него не вошедших
(и, кстати, Лузин много в своих работах философствует). Опубликована переписка 
Лузина со многими людьми. Он сам -- герой многочисленных рассказов современников,
писем, дневников. Многие из этих источников носят неподцензурный характер.
...
Никаких следов  иррационального мышления в математике у Лузина никто не предъявлял.
Нет и никаких сведений о его имяславстве.

 Есть тонкое (или тончайшее)
понимание предметов, которыми он занимался, но нет ничего, выходящего за эти рамки.
Да, он был православным, интересовался философией, в том числе религиозной.
Он дружил с Флоренским (по крайней мере, в молодости), интересовался
(и увлекался) его философией.
Но не видно никаких оснований утверждать что-либо большее.

\sm

Никаких упоминаний <<имяславства>> Лузина в источниках нет. Конечно, всегда можно
настаивать на неполноте (но здесь источников слишком много, будь у Лузина профессиональный мистицизм,
он где-нибудь бы прорвался).
Но некоторые источники дают прямой отрицательный ответ на этот вопрос.

\sm

1) Тахо-Годи \cite{Taho-Godi} (вторая жена А.~Ф.~Лосева) перечисляет участников московского имяславского
кружка, их немного, человек двадцать пять, Лузина среди них нет%
\footnote{Там нет и никого из лузитан.
	Есть один математик, кроме Егорова, Николай Михайлович Соловьёв, по-видимому он совпадает
с автором двух статей по геометрии 1902-1903гг., видимых Zentralblatt'ом. Кроме того
в кружке было несколько людей, связанных с Московским Физматом,
механик Бухгольц,  В.~М.~Лосева
(жена Лосева),
защитившая диссертацию по небесной механике (про возмущение двойной звезды пролетающей мимо звездой),
и молодые люди В.~Н.~Щелкачев (?),
 В.~Л.~Олсуфьев (ум. 1929), П.~А.~Черемухин (1905-1991), человек со сложной биографией,
 по словам Щелкачева перед Войной он стал кфмн в Сталинграде. Как будто не видно, чтобы они имели
 какое-либо отношение к дескриптивной теории множеств.
К 1929г. московский имяславческий кружок, полураспался.}. При этом Лузин (с характеристикой <<великий математик>>)
ей, разумеется, был известен, и про общение Лосева с Лузиным она упоминает.

\sm

2) Имяславцев привлекали по делу 1930 года (под него попал  Егоров;
 в деле упоминался  Бухгольц, который тогда уже отошел от имяславия). Документы эти видел историки С.~С.~Демидов,
 \cite{Demidov}),
 Лузин там не фигурирует...

\sm

г) {\sc Дескриптивная теория множеств.}  Лузин занимался  ей не в одиночку:
а со своими учениками
Александровым, Суслиным, Колмогоровым, Лаврентьевым,
Гливенко, Новиковым, Келдыш, Селивановским, Ляпуновым%
\footnote{К 1940-41гг. в дескриптивной теории множеств работали
	также  ученики Новикова В.~Я.~Арсенин, З.~И.~Козлова,  В.~Янков.}.
Выход в эту область был совершен совместными усилиями Александрова,
Суслина и самого Лузина.
В Ленинграде этой наукой занимались ученики Фихтенгольца
Л.~В.~Канторович и Е.~М.~Ливенсон. Если я не ошибаюсь, никаких
следов ни имяславия, ни христианского мистицизма  у перечисленных
людей не известно. 
Более того, про многих из них мы можем твердо заключить, что имяславцами они не были.

В Польше дескриптивной теорией множеств занимались Серпинский и Куратовский,
в Японии -- Кондо (Kond\^o), были и статьи фон Неймана. Кажется, никого из них мы не тоже не
подозреваем в причастности к имяславию.

\sm

{\sc  д) Понятие <<называния>> у Лузина.} Об этом пойдет речь в последнем пункте настоящего параграфа.
}

 \sm

{\bf \punct Дискуссия об аксиоме выбора.%
\label{ss:choice}} Эта аксиома была основным около-философским  пунктом, около которого ломали копья 
 математики первых трех десятилетий XX века.
Сейчас работающие математики обычно мало интересуются этой дискуссией
(о ней есть большие  интересные исследования
\cite{Moore}, \cite{Medvedev-choice}),
а историко-математические комментарии  на этот счет часто содержат  неточности
или аберрации зрения.

Напомню, что  аксиома выбора (она же аксиома Цермело) состоит в следующем. Пусть есть семейство
(множество) непересекающихся  множеств $A_x$, где $x$ пробегает некоторое множество $X$. 
Тогда существует множество $S$, содержащее по одному элементу $a_x$ из каждого множества
$A_x$. Если множества $A_x$ пересекаются, нужно сказать чуть аккуратнее, например, что
есть отображение $\psi:X\to \cup A_x$, такое, что
$\psi(x)\in A(x)$ для всех $x\in X$. 

Если просто взять и выбросить Аксиому выбора из аксиоматики теории
множеств, то последствия для математики окажутся довольно тяжелыми.
Но острые возражения касались, в основном, случая, когда $X$ несчетно.
Аксиому счетного выбора критики, скорее, были склонны принимать.
Есть  более сильная версия аксиомы счетного выбора, так называемая
{\it Аксиома зависимого выбора}%
\footnote{Пусть $F$ -- многозначное всюду определенное отображение
	множества $X$ в себя (формально: рассматривается подмножество $\Phi$ в
	$X\times X$, такое, что проекция $\Phi$ на первый сомножитель 
	сюръективна). Тогда для любого $x_0\in X$ существует бесконечная последовательность $x_j$,
	такая, что
	$x_{j+1}\in F(x_j)$. (Paul Bernays, 1942).} Она сильнее аксиомы счетного выбора и для обычной математики,
видимо, достаточна (мы привыкли произносить использующие ее рассуждения не задумываясь).
	%\footnote{%Пусть дано отношение $R$ в множестве $X$, то есть подмножество $R$ в $X\times X$,
	%	пусть область определения $R$, то есть проекция $R$ на первый сомножитель $X$,
	%	совпадает с $X$. Тогда для любого $x_1\in X$ существует
	%	последовательность $x_1$, $x_2$, \dots, такая, что для любого $j$
	%	выполнено $(x_j,x_{j+1})\in R$ (Paul Bernays, 1942)},

\sm

{\bf Ниже под словами <<аксиома выбора>> мы имеем в виду ее сильную форму,
обеспечивающую несчетный выбор.}

\sm

На первый взгляд, аксиома выглядит вполне безобидной и самоочевидной и для несчетных
множеств $X$. Но сформулировавший
ее в 1904г.  Цермело
сходу доказал <<теорему Цермело>>: любое множество допускает вполне упорядочение%
\footnote{Упорядоченное множество называется вполне упорядоченным,
если любое его подмножество имеет минимальный элемент.}.

Эта теорема сразу вызвала резкое
неприятие со стороны ряда математиков, в особенности Бореля, Бэра и Лебега
(дискуссия началась в 1905г. в <<Пяти письмах о теории множеств>>. \cite{5}).
Вот первое возражение Бореля \cite{Borel}, цитируем перевод по 
\cite{Medvedev-school}:
\begin{quotation}
Приведя затем такое рассуждение: чтобы вполне упорядочить 
произвольное множество $M$, достаточно выбрать в нем некоторый
 элемент, которому мы припишем ранг 1, потом отличный от 
него другой элемент, которому припишем ранг 2, и т. д. 
трансфинитно, т. е. до тех пор, пока не исчерпаем все элементы
 множества $M$ последовательностью трансфинитных чисел,— Борель 
продолжал: «Но ни один математик не станет рассматривать 
это рассуждение как законное. Мне кажется, что возражения, 
которые можно выставить здесь, действительны и для всякого 
рассуждения, в котором предполагается произвольный выбор, 
совершенный {\bf несчетное} множество раз; такие рассуждения 
находятся вне пределов математики» 
\end{quotation}

Кстати, сразу подчеркнем, аргументы Бореля относятся не к безобидно звучащей аксиоме,
а к следствию из нее (теорема Цермело: любое множество может быть вполне упорядочено).

Обсудим некоторые другие следствия.
Рассмотрим теорему Хана--Банаха (которая, впрочем, тогда
еще не была известна). Пусть дано пространство $Y$, подпространство
$V$
и линейный ограниченный линейный функционал $\ell$ на $V$.
Мы хотим показать, что функционал можно продолжить на $Y$
без увеличения нормы. Добавим к $V$ один вектор $y_1$.
Доказывается лемма, что функционал можно продолжить
на подпространство $V+\C y_1$. Добавим второй вектор $y_2$.
Дальше можно сказать <<и так далее>>, и сказать, что поэтому функционал продолжается 
на все пространство $Y$. Если объемлющее пространство $Y$ сепарабельно, 
то мы можем можем выбрать последовательность $y_1$, $y_2$, \dots,
чья линейная оболочка
плотна в $Y$, в итоге мы получим ограниченный линейный функционал на плотном
подпространстве, дальше мы его продолжим по непрерывности.
Но говорить <<и так далее>> можно и для несепарабельных пространств,
имея в виду трансфинитную индукцию (предварительно вполне упорядочив
пространство $Y$), и рассуждение все равно будет по сути верно (в рамках аксиоматики,
включающей аксиому выбора). 

Результат, однако, для несепарабельных пространств будет странным.
Если в качестве $Y$ взять пространство $L^\infty[0,1]$, то
ограниченные линейные функционалы имеют вид
$$
\ell(f)=\int_{[0,1]} f(x)\,d\nu(x),
$$ 
где $\nu$ -- {\it конечно-аддитивный} заряд
(знаконеопределенная мера)
конечной вариации. Однако ни одной конечно-аддитивной меры,
которая не была бы счетно-аддитивной, предъявить невозможно.
То есть пространство, двойственное к $L^\infty[0,1]$, состоит
из <<осязаемой>> части, которая состоит из счетно-аддитивных
зарядов, или проще, линейных функционалов вида
$$
\ell(f)=\int_{[0,1]} f(x)\,g(x)\,dx,
\quad \text{где $g$ интегрируема},
$$
и большой <<неосязаемой>> части, <<существование>> которой
обеспечивается трансфинитной индукцией (сейчас чаще ссылаются на эквивалентную
ей лемму Цорна).

Ситуация становится еще более странной, если рассмотреть банахову алгебру 
$L^\infty[0,1]$. Согласно теореме Гельфанда, множество $\Omega$ всех характеров
(непрерывных гомоморфизмов этой алгебры в $\C$) является компактным пространством.
Каждой функции $f\in L^\infty[0,1]$ ставится в соответствие функция $F$ на
$\Omega$, определенная формулой
$$
F(\chi):=\chi(f),\qquad \text{где $\chi$ пробегает $\Omega$}.
$$
Это соответствие устанавливает изоморфизм между банаховой алгеброй
$L^\infty[0,1]$ и алгеброй всех непрерывных функций на $\Omega$.

Все прекрасно, но ни одной точки из $\Omega$ предъявить нельзя
(кто никогда не задумывался над такими вопросами, может попробовать). 

\sm

Приведем описание дискуссии в изложении Лузина \cite{Luzin-teor-funkts}:
\begin{quotation}
	\dots впервые в математике встал вопрос во 
	всей его остроте о существенной необходимости делать различие 
	между <<тем, что есть на самом деле>> (= реальность) и „тем, 
	что может быть доказуемо для человеческого интеллекта (=  
	познание)— вопрос мало новый для философии, прибавим. Едва 
	только появился на горизонте этот вопрос, как не стало  
	обычного единомыслия математиков, тотчас же разбившихся на толки. 
	Для одних, получивших название идеалистов, <<все это — слишком 
	антропоморфично: надо ни на одну минуту не забывать о самом 
	предмете и поменьше говорить о человеческой возможности его 
	охватить, — говорит же нам палеонтология о фауне третичной 
	эпохи или астрофизика — о строении всей поверхности луны, 
	хотя в первом случае тогда не было человека, и во втором 
	случае имеется невидимая половина у нашего спутника. Это 
	мнимое злоупотребление категориями времени и пространства не 
	мешает быть упомянутым дисциплинам <<истинными науками>>. 
	<<Мало ли вещей на свете, о которых мы не знаем и о которых 
	мы никогда не узнаем. Вот, например, капля воды. Число  
	молекул конечно, — но какое оно? Я не знаю и никогда не узнаю — 
	и это еще не резон уничтожать кинетическую теорию материи. 
	Все это — слишком антропоморфично>> (Адамар). 
	
	Для других, по иронии судьбы получивших прозвище  
	реалистов, необходимо избегать жонглирования с символами,
	 которые не соответствуют ничему. Истинная наука никогда не будет 
	 сборищем пустых слов и чисто логических понятий без  
	 концепций. {\it Наука не есть логомахия}.. За словами всегда должна 
	 скрываться сама реальность. Рассуждение Цермело — только 
	 греза, так как каждый из идеалистов, говоря о выборах,  
	 выбирает и грезит по-своему, и нет не только возможности  
	 сообщить своему собеседнику о проделанных в бесконечном  
	 количестве выборах, но и быть согласным даже с самим собою. 
	 Построение, которое нельзя описать, рассуждение, немогущее
	  быть выполненным во всех его шагах до конца, — все 
	 это находится вне науки. Победы и горести в мире призраков 
	 не окажут никогда никакого эффекта — ни задерживающего, 
	 ни ускоряющего — на ход науки, занятой открытием конечных 
	 соотношений между вещами мира чувств или мира концепций 
	 (Борель). 
	 
	 Эта, к изумлению самих математиков, обнажившаяся разница 
	 взглядов с течением времени не смягчилась, но была просто 
	 отмечена и констатирована. По молчаливому соглашению  
	 решено было не настаивать на ней и не углублять ее, но  
	 подождать новых фактов \dots. 
	 
	 Возобновление спора произошло в два-три последние года, 
	 и это произошло в столь подчеркнутой и ни перед чем не  
	 останавливающейся форме, что дало повод Адамару говорить о  
	 подлинном <<землетрясении в математике>>. И при том это  
	 потрясение отнюдь не ограничивается <<окраинными владениями  
	 математического царства>> — что было бы простою болезнью роста  
	 науки, связанною с неладами и неустройствами при недавней  
	 ассимиляции какого-нибудь нового понятия,—-но идет в самую 
	 глубь его, опрокидывая и разрушая все, что считалось прочным 
	 со времени Ренессанса. 
	 
	 Атака повелась <<реалистами>> на этот раз на принцип исключенного третьего, и
	 авторами ее явились Брауэр (Brouwer) 
	 и Вейль (Weyl). Этот принцип поставлен ими под категорическое 
	 сомнение. Так как он все-таки действенен для текущей жизни, 
	 то Брауэр считал его законным лишь для конечных множеств, 
	 отрицая его применимость к каким бы то ни было бесконечным 
	 множествам — включая сюда и натуральный ряд. И дабы  
	 объяснить его действенность в области обыкновенной логики, Брауэр 
	 изменяет порядок вещей, объявляя логику лишь частью теории 
	 конечных множеств....
	 
	 Интуиционизм 
	 в тех его формах, которые он принял у последователей Брауэра 
	 и Вейля и в которых, прибавим, неповинны его инициаторы, 
	 представляет из себя скорее игру ума, чем науку. Он есть, 
	 в сущности, смесь обостренной критики и догматизма  
	 теологического характера. Он возводит конечное в догму; но та  
	 <<революция>>, которую он обещал внести в математику, является 
	 не освобождением мысли, но загромождением путей  
	 исследования, т. е. в сущности, одним из реакционных запретов мысли. 
	 
	 Лично мне кажется, что мнение Адамара о необходимости 
	 прекращения бесполезных и не могущих ничего выяснить 
	 споров — вполне справедливо. Необходимо в самом деле  
	 подождать новых фактов, касающихся возможности полной 
	 арифметизации континуума. {\bf Из всех тезисов Брауэра и Вейля самым 
	 интересным является утверждение, что континуум не есть  
	 множество точек.} 
\end{quotation}

В этом отрывке 
и в большинстве более поздних описаний дискуссия выглядит как философско-логическая.
Это было и так, и не так.

 Надо думать,   участники дискуссии имели разные философские взгляды. Но у них были и различия в математическом опыте. 
Приведем слова Лебега в дискуссии на международном съезде математиков
в Цюрихе, 1932, см. \cite{Gonseth} (цитируется по \cite{Medvedev-choice}):
\begin{quotation}
	Когда 35 лет тому назад пороховая бочка взорвалась от подожженного г.Цермело фитиля, 
	то {\bf обнаружился один замечательный факт: все те, кто до того времени пользовался множествами,
	были против аксиомы выбора; все, кто не пользовался, были за аксиому выбора.
	Это -- любопытный факт, и можно было бы поставить вопрос:
	то ли нас сделала реакционерами закоренелая привычка обходиться без явно сформулированной аксиомы выбора,  
	то ли здесь была другая причина.} 
\end{quotation}

Ровно так оно и было. {\bf Именно Борель, Бэр и Лебег применяли
нетривиальную теорию множеств к функциям действительного переменного%
\footnote{То же замечание относится  к Лузину, который занимался нетривиальной теорией множеств.
	Кажется, Хаусдорф все же был исключением.}
(в прочих областях математики множества, если и использовались, то как
удобный язык)}. Эти французы были против
следствий, которые извлекались из аксиомы выбора. И именно они понимали
этот предмет лучше кого-бы то ни было, дело было даже не в странности
следствий -- вроде парадокса Банаха--Тарского (которого, впрочем, в момент начала диспута еще не было), -- а том, что
эти следствия относились к какому-то потустороннему меру, не имеющему отношения 
к обычной математике). Но коль скоро следствия
согласно всем правилам логики извлекались
из новой аксиомы, то они пытались найти логические и философские
 доводы против нее самой. И здесь уж конкретный человек
 применял ту философию, к которой был горазд. 

Аксиома выбора утверждает существование множества из представителей,
В связи с этим развернулись споры о теоремах существования
и о том имеют ли смысл теоремы существования, которые не дают возможность 
найти сам существующий объект. Я не буду их цитировать, лишь отмечу,
что слова <<теорема существования>> могут иметь в математике разный смысл. 

\sm
 
 Первый пример -- теорема существования и единственности для обыкновенных
 дифференциальных уравнений. Обычное доказательство (со сведением к интегральному уравнению и решением его методом последовательных приближений Лиувилля--Пикара), по существу,
 конструктивно и дает способ приближенного нахождения решения с любой точностью%
 \footnote{Вполне конструктивно и доказательство со сходимостью ломаных Эйлера. 
 В качестве средства для реального вычисления они сомнительны, но теоретически  ответ дают.}.
 Точно так же конструктивно обычное доказательство теоремы существования для
 решения задачи Дирихле для уравнения Лапласа путем сведения задачи к интегральному уравнению.
 
 \sm
 
 Второй пример - доказательство Кантора существования трансцендентных чисел.
 А именно, множество алгебраических чисел счетно, а множество всех
 вещественных чисел несчетно. Поэтому неалгебраические числа существуют.
 Это доказательство дает способ построения трансцендентного числа. А именно,
 мы точно формулируем правило, нумерующее алгебраические числа отрезка $[0,1]$. 
 Составляем из них список и пишем действительное число, которое
 отличается от $n$-ного числа списка в $n$-ном знаке.
 Этот пример, по-видимому, отличен от предыдущего. Там мы имели возможность
 получать решение с любой степенью точности, а здесь вопрос о точности бессмысленен.
 
 \sm
 
 Третий пример. Доказательства от противного. Например, докажем,
 что гладкое отображение $f$ замкнутого
  единичного шара $B$ в себя имеет неподвижную точку
 (ослабленный вариант теоремы Брауэра).
 Допустим, что такой точки нет. Для $x\in B$
 рассмотрим луч, идущий из точки $f(x)$ в точку $x$. Пусть $\psi(x)$
 -- пересечение луча со сферой $S$. Получается непрерывное отображение
 $\psi:B\to S$, тождественное на сфере (ретракция). Мы можем стянуть
 сферу в точку по шару. Применяя к этой гомотопии отображение $\psi$,
 мы получаем, что сфера стягиваема по себе, т.е. гомотопна точке. Поэтому все замкнутые формы
 на сфере должны быть точны. Это касается, в частности, стандартной формы
 объема. Но тогда интеграл от этой формы по сфере равен 0, что, очевидно,
 не верно. Поэтому неподвижная точка существует. Однако никакого способа поиска
 этой точки наше доказательство не дает.
 
 \sm
 
 Четвертый пример -- это упомянутые выше характеры алгебры $L^\infty[0,1]$,
 существование которых доказывается с помощью аксиомы выбора.
 
 \sm
 
 На первый взгляд, третий и четвертый пример однотипны, в том смысле,
 что мы устанавливаем существование объекта, не предъявляя способа его нахождения.
 Однако различия начинаются на следующем шагу.
 Зная наличие неподвижной точки, мы можем с ней работать и при желании мы
 можем ее какими-то другими способами искать. Эта точка так или иначе 
 включается в общую картину анализа. Но мы не найдем характера
 алгебры $L^\infty[0,1]$ ни точно, ни приближенно, мы знаем только то, что они есть,
 и что их много. Следствия для обычного анализа из их существования
 извлекать проблематично, а если это и получится, то, скорее всего,
 окажется, что это можно делать напрямую, без привлечения аксиомы выбора и банаховых
 алгебр.

Приведем на этот счет две цитаты из Лузина. 
Сначала из \cite{Luzin-teor-funkts}:
\begin{quotation}
Если теперь сойти с принципиальной почвы и перейти к практике 
текущей математической работы, то, с одной стороны, в пользу 
аксиомы Цермело говорит ее «очевидность», величайшая простота 
ее применения в математических рассуждениях, чрезвычайная 
легкость получения при ее помощи самых разнообразных 
примеров, носящих хотя иногда исключительно парадоксальный 
характер, но все же не повлекших до сих пор ни к какому 
противоречию (например, Хаусдорфу удалось разделить шар радиуса $R$ 
на такие четыре непересекающиеся попарно части $A$, $B$, $C$, и $D$, 
что, двигая, как твердые тела, части $A$ и $B$ и прикладывая 
надлежащим образом их друг к другу, он получил опять шар радиуса $R$, 
и, двигая так же точно и прикладывая части $C$ и $D$ друг к 
другу, он получил второй такой же шар радиуса $R$. Множества 
$A$, $B$, $C$, и $D$ суть неизмеримые и образованы Хаусдорфом при 
помощи аксиомы Цермело). С другой же стороны, против нее 
говорит именно эта самая чрезвычайная легкость ее применения 
и немедленность даваемых ею ответов, так как {\bf  математические 
сущности, сформированные при помощи ее, не крепки, не 
обладают устойчивостью, имея слишком расплывчатые, 
неопределенные свойства, чтобы практически служить затем точкой опоры 
для математических рассуждений, направленных уже на 
классические математические предметы}. Напротив, образование 
математического предмета без аксиомы Цермело часто представляет 
чрезвычайные трудности, зато такой математический предмет
будучи построен, почти всегда имеет большую ценность для 
дальнейших изысканий». 
\end{quotation}

Мы видим, что Лузин  рассматривает принятие или непринятие аксиомы выбора не 
только и не столько как
 априорный логическо-философский  вопрос, сколько вопрос о следствиях,
 которые могут быть из этой аксиомы извлечены. В целом его отношение
 было скептическим. 
Например,
\cite{Luzin-sbornik-1926}:
\begin{quotation}
	применять свободный выбор — это значит, по моему 
	мнению, жонглировать соединениями пустых слов, смыслу 
	которых не соответствует никакой интуитивно доступный факт.	
\end{quotation}	

Однако Лузин сознательно использовал ее в части своих работ.

\sm

{\bf \punct Что случилось с аксиомой выбора в дальнейшем, и почему сейчас 
против нее мало кто возражает?}
В 1940г. Гёдель доказал непротиворечивость этой аксиомы (то есть, если
теория множеств без аксиомы выбора непротиворечива, то она остается непротиворечивой
и при ее добавлении). Но, скорее всего, решающую роль сыграло не это.

Если выкинуть сильную форму аксиомы из аксиоматики (оставив, вместо нее аксиому зависимого выбора),
то ничего страшного с математикой не случится.
Возникнет, однако ряд неудобств технического характера.
В условиях некоторых теорем придется делать оговорки, например,
в теореме Хана--Банаха или теореме Крейна--Мильмана придется накладывать
условие сепарабельности пространств. Соответственно, при применении 
теорем придется оговорки проверять (что иногда бывает не совсем приятно).
Иногда при доказательствах дешевле сослаться на аксиому выбора, чем ее обойти
(например, существование меры Хаара).
 Можно заметить, что кратчайшее определение аменабельности группы
 имеет в виду аксиому выбора (хотя на самом деле, у этого понятия есть эквивалентные
 формы, не имеющие этой аксиомы в виду). Тоже самое касается теоремы о короне...
 В общем есть разные возможности для неудобств...
 
 Есть неудобства и другого плана. При преподавании математики
 для математиков
 с отменой аксиомы выбора останется слишком много неразрешимых  вопросов.
 Например,
 существует ли неизмеримое множество? Существует ли разрывный линейный функционал
 на гильбертовом пространстве? (и т.д и т.п.). Без конца запинаться на тему о том, что
 данный вопрос неразрешим как-то неуютно.
 
 Неудобства аксиома выбора тоже создает (однако кому она не нравится могут ее и не применять), но неудобства при ее отмене их перевешивают.
 Несмотря на  осмысленность критики начала XX века
 вопрос на сегодняшний день решается таким прагматически-эстетическим образом.

 \sm

{\bf\punct Называемые множества.%
\label{ss:name}} В связи с бурными дискуссиями
начала XX века естественно возник вопрос о том что такое <<явно описываемые
множества>> (есть много других слов, употребляемых в этом контексте,
<<конструктивные>>, <<эффективные>>, <<хорошо определенные>> и т.п.). Лузин часто употреблял слова
<<называемое множество>>, <<называть>>. В основном этот термин встречался
в его работах на французском, <<ensemble nomm\'e>>(см. выше  в п. \ref{ss:intention} длинную
сентенцию Лузина на эту тему). Называемое множество -- это множество, которое можно
однозначно описать конечной фразой русского языка в соответствии с правилами
описания, принятыми в теории множеств. Понятно, что таких множеств
может быть лишь счетное число. В частности, можно  <<назвать>>
лишь счетное число действительных чисел.

Слово 
<<ensemble nomm\'e>> появились в 1905 году в мемуаре Лебега
\cite{Lebeg-nomme}. Это та самая статья, где был ляп,
найденный Суслиным, и которая в остальных отношениях замечательна
(и оказалась отправной точкой для разных размышлений Лузина).
 Там Лебег доказывает непустоту всех классов Бэра%
 \footnote{Нулевой класс Бэра состоит из непрерывных функций.
 	Первый класс - из поточечных пределов последовательностей непрерывных функций.
 	Для счетного трансфинитного числа $\alpha$ класс номер $\alpha$
 	состоит из поточечных пределов последовательностей функций,
 	лежащих в  классах $<\alpha$. Классы Бэра -- двойник классов борелевских
 	множеств.}. Для этого он по каждому счетному трансфиниту
 $\alpha$ строит  универсальную (Бэровскую) функцию двух переменных $F_\alpha(x,t)$ такую,
 что для любого фиксированного $t$ функция $f_t(x):=F_\alpha(x,t)$
 содержится в классе с номером $<\alpha$ и любая функция класса номера меньше
 $\alpha$ имеет такой вид для некоторого $t$. Интересно, что универсальная функция строится явно.
 С помощью этой универсальной функции Лебег доказал непустоту всех классов Бэра
 (равносильным образом все трансфинитные классы борелевских множеств не пусты).
 
Это не отменяет вопроса о явном построении примеров функций из классов Бэра.
 Сам Бэр построил  функции первых трех классов. Позже Келдыш построила примеры множеств 4-ого класса%
 \footnote{Из книги Лузина \cite{Lus-anal}:
 \newline
{\it ... элемент класса 4 есть множество тех 
иррациональных точек $(\alpha_1, \alpha_2,\dots,\alpha_\nu)$}[это обозначение для неполных частных
при разложении числа в цепную дробь], {\it содержащихся между 0 и 1,
у которых среди неполных частных $\alpha_\nu$, есть бесконечное
множество различных, каждое из которых повторяется
бесконечное множество раз.}
 \newline
 Отметим, что цепные дроби здесь ни причем, фактически работа идет в <<пространстве Бэра>>
 $\N^\N$, см. сноску \ref{fo:baire}.%
 \label{fo:keldysh}}
 а потом около 1940г. \cite{Keldysh-B0}, \cite{Keldysh-B}
 и любых трансфинитных классов (в не большей степени, в которой  сами эти трансфинитные числа являются явными),
 по-видимому, это был один из последних позитивных результатов теории множеств, наступало
 уже 
 время доказательств недоказуемости.
 
Вернемся к мемуару Лебега. Он строит (<<называет>>)  неборелевское множество (формально, небэровскую
функцию). Понятно, что неборелевские множества существуют (потому что
борелевских множеств всего континуум, а множество подмножеств континуума
имеет б\'ольшую мощность), однако изобретение конкретного неборелевского множества было замечательным результатом.
В принципе слово <<явно>> в данном случае может быть предметом сомнений.
	В конструкции Лебега участвует множество всех счетных трансфинитных чисел,
	то есть непонятный $\aleph_1$. В тот момент это еще могло выглядеть относительно безобидным,
	но в последовавшие три десятилетия степень недоумения математиков при виде этого
	алефа возрастала. 
	
Завершает Лебег свой мемуар вопросом:
\begin{quotation}
Можно ли {\bf назвать} неизмеримое множество?%
\footnote{Как мы сейчас знаем, нельзя.}
\end{quotation}

Для Лузина вопрос о назывании различных неборелевских множеств
(без аксиомы выбора и без использования множества всех счетных трансфинитных чисел)
приобретал особое значение: это был вопрос о реальности предмета
развиваемой им науки.

Он к нему неоднократно возвращался,
и вот один из его примеров, замечательный по изяществу. 

Разложим число из отрезка $[0,1]$ в цепную дробь. Пусть $n_1$, $n_2$, \dots --
неполные частные. Мы берем множество $E$ всех чисел, таких, что
существует бесконечная подпоследовательность $n_j$, такая, что
$n_{k+1}$ делится на $n_k$ для всех $k$. Тогда $E$
-- аналитическое неборелевское множество.

Лузин много обсуждал <<называние>> в разных задачах теории множеств.
Приведем еще пример; цитата из  статьи Успенского \cite{Usp}: 
\begin{quotation}
...	На действительной прямой
	$\R$ существует множество из $\aleph_1$ точек. Это наблюдение оставляет открытым
	следующий вопрос: {\it можно ли назвать совокупность точек континуума $\R$,
	имеющую мощность ровно $\aleph_1$}? Мы будем говорить об этом вопросе как о проблеме
	Лебега, так как Н.~Н.~Лузин
	 [\cite{Luz-anal-rus}, стр.202] приписывает его
	Лебегу и даже приводит из Лебега соответствующую цитату; на самом же деле
	первая четкая формулировка рассматриваемого вопроса принадлежит самому
	Н. Н. Лузину, и цитата из Лебега это только подтверждает.
\end{quotation}

Под словом <<назвать>> подразумевается однозначно описать
подмножество континуума из $\aleph_1$ элементов.

\sm

{\small
Если же отвлечься от дела и вернуться к <<именам Бога>>, то во втором случае
  еще можно пытаться что-то заподозрить (хотя это вполне естественный и правильный математический вопрос 
  и попытка атаковать континуум-гипотезу,
  со стороны Лузина было много таких атак),
то первый пример все же вне подозрений. Между <<именем Божьим>> у имяславцев
и <<ensemble nomm\'e>> у Лебега и Лузина единственным сходством является слово
<<имя>>, точнее, его корень. Еще раз напомним, у Лузина присутствовал сильнейший скептицизм 
в отношении теории множеств, а вовсе не вера в ее божественность
и божественность операции наименования.}

 \section{Ученики атакуют учителя%
 \label{s:attack}}
 
 \COUNTERS
 
 \epigraph
 {Эти собрания носили самый безобразный характер.}
{ Н.~К.~Бари \cite{TyulinaA}.}
 
 \epigraph
{Я присутствовал на этих собраниях. Общее впечатление было неприятное, даже омерзительное.}
{Л.~С.~Понтрягин \cite{Pon1}}

\epigraph 
 {Для меня была невыносима эта травля...}
 {Б.~В.~Гнеденко \cite{Gned}}
 
 \epigraph{
 В стенах Института было организовано побоище  на Н.~Н. — 
 собрали сотрудников Института и других 
 математиков, каждый из них пинал Н.~Н. как мог.}
 {С.~М.~Никольский \cite{Niko}}

 Смысл приведенных цитат  различается. К сожалению, автор данных записок
 здесь должен согласиться со  сказанным и с той, и с другой, и с третьей стороны. 
 
 Так или иначе, июльская история 1936 года много обсуждалась в научной литературе,
 обсуждалась не очень точно,
   в литературу популярную она вышла уже в совсем превратном виде. Поэтому нам придется
 взять 
 на себя (не очень приятный) труд по ее разбору.
 
 \sm
 
 Если брать отсчет с 3 июля 1936 года, то в следующие
 дни сторону разума представляла сторона,  противостоявшая
 изгнанию Лузина из Академии: из математиков это Бернштейн, Крылов,
 Чаплыгин, Меньшов, Бари, Новиков, Лаврентьев. Из нематематиков мы знаем об усилиях геолога В.~И.~Вернадского, физика С.~П.~Капицы, зоолога Н.~В.~Насонова, химика Н.~С.~Курнакова.
 По-видимому, противодействовали опасному развитию событий также биохимик А.~Н.~Бах и геолог А.~Е.~Ферсман.
 Но, не забывая этого, мы все же  должны попытаться понять позиции разных участников
 атаки.
 
 \sm
 
 {\bf\punct Начало облавы.%
 \label{ss:nachalo-oblavy}}
 Юшкевич \cite{Yush1} описывает это так:
 \begin{quotation}
 	Следует признать, что инсценировка «дела Н.~Н.~Лузина» была искусно продумана с самого начала. Поводом явился, казалось бы, незначительный эпизод.
 	В конце июня 1936 г. Лузин был приглашен на выпускной экзамен в школу № 16
 	Дзержинского района Москвы. Затем его попросили публично высказать свои
 	впечатления. 27 июня в «Известиях» появилась его статья «Приятное разочарование», где говорилось, что приехал он в школу с несколько предвзятым
 	мнением о плохой математической подготовке выпускников, но в ходе экзамена
 	был приятно разочарован: все ответы, даже на вопросы, выходившие за рамки
 	школьной программы, были правильны, свидетельствовали о глубоком понимании предмета, а слабых учеников не оказалось ни одного.
 	
 	Конечно, отзыв был чрезмерно хвалебным и, быть может, Лузин хотел
 	снискать благоволение высших сфер. Как вскоре разъяснил директор Г.~И.~Шуляпин, часть ответов была очень хорошей или хорошей на самом деле, но имелись и «троечники», примерно одна треть. Тем не менее приглашение на экзамен
 	являлось несомненной ловушкой%
 	\footnote{Все же стоит уточнить. Из статьи Лузина в Известиях:
 		{\it
 			В числе многих московских ученых я был приглашен присутствовать на выпускных
 			испытаниях в одной из школ... Я выбрал 16-ю школу Дзержинского района. Пришел туда во время испытаний выпускников по тригонометрии.}
 		Ловушкой (и, скорее всего, пробным шаром тоже), было предложение написать статью в <<Известиях>>.}: {\bf  в математических кругах было известно, что
 	Лузин всегда дает только хорошие отзывы на любые представленные ему работы} и на экзаменах чрезвычайно снисходителен. {\bf Предложение опубликовать
 	в печати свои впечатления явилось следствием элементарного расчета. Было
 	заранее ясно, что выскажется Лузин весьма похвально, за что его можно будет
 	подвергнуть критике, после чего приступить к более серьезным обвинениям}.
 	К тому же в Наркомпросе он высказался о преподавании математики в средних
 	школах в целом отрицательно, правда особенно подчеркивая низкое качество
 	общепринятого тогда учебника геометрии Гангнуса%
 	\footnote{Об этом учебнике, см. ниже п.\ref{ss:kiselev}}.
 	
 	Через пять дней началась газетная кампания против Лузина. 2 июля «Правда» поместила «Ответ академику Лузину», подписанный директором 16-й школы
 	Г.~И.~Шуляпиным. В статье говорилось, что в школе есть разные ученики,
 	в том числе слабые, что школьники не владеют навыками математического
 	мышления, не умеют работать с книгой. Школа нуждается не в лицемерных
 	похвалах, а в товарищеской критике. Далее следовал вопрос:
 	 «Не было ли Вашей целью замазать наши недостатки и этим самым нанести нашей школе
 	вред?».
\end{quotation}

 На другой день, 3 июля 1936 года в Правде появилась анонимная статья 
 <<О врагах в советской маске>>, где Лузину были предъявлены следующие обвинения:
 
 1) Многочисленные хвалебные отзывы, не соответствующие истине (приводятся конкретные случаи):
 \begin{quotation}
 	Фабрикация заведомо ложных похвальных отзывов — это частица линии академика Лузина, линии на засорение советской математической науки людьми неподготовленными. Такую
 	линию академик Лузин последовательно ведет вот уже несколько лет и до самых последних месяцев и дней, когда он рекомендует на научную работу в Академию наук людей,
 	которым в лучшем случае надо начать прохождение университетского курса.
 \end{quotation} 
 
 2) Лузин публикует лучшие свои работы за границей, а в СССР печатает 
 полную  халтуру:
 \begin{quotation}
 	Это своеобразное вредительство в науке академика Лузина видно и на примере его
 	собственных работ. Чтобы казаться активным членом Академии наук, он печатает в различных изданиях в СССР многие свои якобы научные статьи. Но научная ценность большинства этих статеек ничтожна. Сам Н.~Лузин в дружеских беседах не стесняется заявлять, что он ежегодно выбрасывает на рынок науки несколько мемуаров, внутренне
 	смеясь над их содержанием. Он называет эти работы белибердой, которую надо печатать исключительно для «устрашения количеством». Эти, с позволения сказать, «труды»
 	нарочито разводнены, и аспирант Академии наук Ф.~Р.~Гантмахер доказал на втором всесоюзном математическом съезде, что результаты серии работ Н.~Лузина о методе академика А.~Н.~Крылова (составление векового уравнения), занимающие 160 страниц, легко
 	умещаются на ... трех страничках. Под видом научных работ печатает Н.~Лузин некрологи, мелкие и бессодержательные заметки-комментарии и тому подобное. А более или
 	менее самостоятельные свои работы Н.~Лузин отсылает для печатания за границу — во
 	Францию, Польшу и даже... Румынию.
 \end{quotation}
 
 3) Обвинение в преследовании Суслина и его смерти:
 \begin{quotation}
 Н. Лузин сделал все возможное, чтобы
 выставить своего ученика Суслина из Москвы, не дать ему возможности работать. И как
 только умер (в 1919 г.) М. Суслин, его открытием Н.~Лузин не замедлил воспользоваться
 в напечатанных за границей работах, выдавая открытие погубленного им ученика за
 свое собственное.	
 	\end{quotation}
 
 .
 
 4) Присвоение результатов учеников. Названы Суслин (см. выше) и Новиков:
 \begin{quotation}
 	Не стесняется и теперь Н.~Лузин печатать работы своих учеников под
 	своей фамилией, как это было в прошлом году с книгой «О некоторых новых результатах
 	дескриптивной теории функций». На обложке этой книги красуется фамилия Лузина, а
 	внутри — работа его ученика П.~Новикова.
 \end{quotation}
 
 Статья завершается потоком политических ругательств:
 \begin{quotation}
 	Но полупочтенный академик забывает, что большевики хорошо умеют распознавать змей, в какие бы шкуры они ни рядились. Мы хорошо знаем, что {\bf Н.~Лузин — антисоветский человек.}
 	Академик Лузин пускает в ход лесть (иногда тонкую лесть) советской научной молодежи —
 	 будто она все уже знает, умеет, научно подкована и прочее и прочее. А исподтишка
 	ухмыляется, зубоскалит, рассказывая по секрету дружкам своим, что все это не всерьез и
 	есть надежда на то, что время молодежи подходит к концу, и подсиживает, изгоняет из
 	Академии действительно талантливых молодых ученых%
 	\footnote{Эта фраза несколько загадочна, по-моему в Стенограмме она не отражена. Это может
 		свидетельствовать о чем-то, нам не известном, а может  быть и оговоркой газетчиков.
 		Речь также могла идти не о самой Академии, а о Математической группе Академии,
 		которая включала дополнительных лиц.}\dots
 	
 	Мы знаем, откуда вырос академик
 	Лузин: мы знаем, что он один из стаи бесславной царской «Московской математической школы», философией которой было черносотенство и движущей идеей — киты российской реакции: православие и самодержавие. Мы знаем, что и сейчас он недалек от
 	подобных взглядов, может быть, чуть-чуть фашистски модернизированных. Но социальная почва, взрастившая Лузиных, исчезла, ушла и уходит из-под ног. Академик Лузин
 	мог бы стать честным советским ученым, каких из старого поколения много. Он не захотел этого; он, Лузин, остался врагом, рассчитывая на силу социальной мимикрии, на
 	непроницаемость маски, им на себя напяленной.
 	Не выйдет, господин Лузин! Советская научная общественность срывает с вас маску
 	добросовестного ученого и голеньким, ничтожным предстаете перед миром вы, ратующий якобы за «чистую науку» и продающий интересы науки, торгующий ею в угоду
 	прежним хозяевам вашим, нынешним хозяевам фашизированной науки\dots
 \end{quotation}
 
 Две последовательных статьи в Правде были сигналом. Трудящиеся 
 должны были ответить на них митингами с осуждением Лузина и 
 лузинщины. Первый митинг состоялся уже 3 июля в Институте Стеклова.

 7 июля Президиум АН СССР создал специальную комиссию: 
 \begin{quotation}
	Для разбора дела ак. Н.~Н.~Лузина в связи со статьями «Правды» — «Ответ
	академику Н.~Лузину» и «О врагах в советской маске» и в связи с фактами, оглашенными на собрании научных работников Математического Института АН, образовать Комиссию под председательством Вице-Президента Академии Наук ак. Г.~М.~Кржижановского в составе академиков: А.~Е.~Ферсмана, С.~Н.~Бернштейна, О.~Ю.~Шмидта,
 	И.~М.~Виноградова, А.~Н.~Баха, Н.~П.~Горбунова, чле\-нов-кор\-рес\-пон\-ден\-тов АН: 
 	Шнирельмана,  Соболева,  Александрова и профессора  Хинчина.
 \end{quotation}
 
 Председатель комиссии -- Кржижановский (1872—1959) -- соратник Ленина,  член ЦК ВКП(б)
 1924-1939, один из немногих членов ЦК, переживших 37ой год. В 1929г. под давлением власти 
 был избран в АН СССР и стал ее вице-президентом. По совместительству,
 инженер-энергетик, занимался разработкой многих хозяйственных проектов. 
 
  Из математиков -- членов комиссии, на стороне Лузина
 был лишь Бернштейн. Шнирельман,  Соболев,  Александров и  Хинчин
 участвовали в атаке. Отдельно от них в нападениях участвовал Шмидт.
 Позиция Виноградова не ясна, он молчал. Атака поддерживалась присутствовавшими (или появлявшимися)
 на заседаниях комиссии
 Люстерником, Гельфондом, Сегалом и Аршоном.
 
 В комиссию также входили влиятельные нематематики.

 Горбунов (1892—1938) -- партийно-государственный деятель,
 один из участников операции по советизации АН СССР в 1928-1930г.
 В 1929-30г. вице-президент ВАСХНИЛ.
 С 1935г. -- академик АН СССР, также занимал
  должность <<непременного секретаря>> АН СССР. 
 Принимал активное участие в  политических обвинениях
  Лузина,
 научно-этические вопросы его интересовали мало.
 
 Ферсман Александр Евгеньевич (1883—1945) -- знаменитый геолог,
  Бах Алексей Николаевич (1857—1946) -- биохимик-революционер.
 Оба участвовали в обвинениях Лузина, но оба, по-видимому, были противниками
 радикального решения.

 Комиссия заседала 7, 9, 11, 13, 15 июля.
 До нас дошла стенограмма заседаний этой комиссии,
 эта стенограмма была  чудесным образом <<обретена>> в конце 80х сразу после смерти последних математиков -- участников атаки на Лузина,
 \begin{quotation}
{\it А.~П.~Юшкевич} \cite{Yush1}, 1989:	
 Стенограммы, о которых идет речь, обнаружила сотрудница Архива АН СССР И. П. Староверова. 
 
 \sm	
 	
 {\it С.~С.~Демидов, В.~Д.~Есаков} \cite{DeEs}:
 Т. А. Токаревой совместно с работниками Архива АН СССР был обнаружен случайно сохранившийся и, весьма вероятно, неполный экземпляр машинописного текста стенограммы заседаний Комиссии
 Академии наук СССР по делу академика H.~H.~Лузина. Обязательные экземпляры, которые должны были храниться в Архиве в положенном для них месте, оказались уничтоженными заинтересованными лицами.
\end{quotation}

 Сторону Лузина представлял Бернштейн (присутствовал 7 и 15 июля)
 и Крылов (присутствовал 15 июля). В остальные дни, 9, 11 и 13 июля Лузин
 должен быть обороняться в одиночку (судя по всему, ситуация разрешилась в целом
 благополучно для Лузина к 13 числу).
 
 Но вернемся к началу.
  9 июля последовала новая анонимная статья <<Традиции раболепия>> в <<Правде>>, где критике были подвергнуты уже участники атаки на Лузина (Александров, Колмогоров, Хинчин),
  а также Бернштейн. Об этом чуть ниже.
 
 \sm
 
 {\bf \punct Рассказ Ефремовича.%
 \label{ss:efremovich}} В.~А.~Ефремович (тополог, ученик Александрова) в 1986-87гг давал интервью неизвестному корреспонденту.
 Магнитофонные пленки позже были расшифрованы  его дочерью \cite{Efremovich} и  опубликованы. Мы частично 
 приводим большой отрывок, посвященный истории с Лузиным:
 \begin{quotation}
 {\it Сколько лет она [Лузитания] длилась? Года до 24-25-го?}
 
  \sm
 
 Потом она распалась, но какая-то часть ее оставалась при Лузине, и в этой оставшейся части были Людмила Всеволодовна Келдыш и ее муж П.~С.~Новиков. Оба они сохранили о нем очень хорошие воспоминания и любовь к нему. А старшие ученики постепенно с ним расходились. Например, Урысону он запретил заниматься топологией. Урысон не послушался и, стало быть, перестал быть его учеником. Он был одновременно учеником и Лузина, и Егорова.
 .......................................
 
 Что можно еще сказать о Лузине: он был психически не вполне нормальным. Могу рассказать историю, которой я лично не был свидетелем, но знаю, что это так было. Одному из математиков, Александру Осиповичу Гельфонду, он сказал, что с ним хотел поговорить Бухарин. Бухарин в то время был в полной славе. Гельфонд пошел к Бухарину, записался на прием. Бухарин спрашивает: «Зачем Вы пожаловали?» Гельфонд: «Лузин сказал, что Вы хотели со мной говорить». Оказывается, что ничего подобного не было%
 \footnote{Вот как этот эпизод описывал А.~И.~Яглом \cite{Yaglom}, знакомый с Гельфондом и его окружением:
 	\newline{\it
 	Лузин, как известно, любил издеваться над людьми. А Гельфонд был, во-первых, партийный, а, с другой стороны, без чувства юмора. И когда Гельфонд решил седьмую проблему Гильберта, то Лузин объявил ему, что его лично вызывает Николай Иванович Бухарин (тогда главный редактор <<Правды>>), чтобы напечатать в <<Правде>> большую статью о решении седьмой проблемы Гильберта. И Гельфонд пошел к Бухарину на прием. Бухарин, естественно, его не принял. Гельфонд сказал что Бухарин лично интересовался его достижением, на что секретарша ответила, что Бухарин ничего про это не знает, но если это действительно такое великое достижение, то он может послать корреспондента, который с Гельфондом поговорит.}
 	\newline
 	Для справки: Бухарин был ответственным редактором газеты <<Правда>> 7.1918 - 4.1929 и газеты <<Известия>> 2.1934 - 27.2.1937.
 	Более правдоподобно, что речь идет о <<Известиях>>. Седьмую проблему Гильберта Гельфонд решил в 1934,
 а первый подход сделал в 1929х.}%
 . Странное такое поведение. Так что какие-то отклонения от нормы у него, несомненно, были.
 
  \sm
 
 {\it По Вашим рассказам получается, что Лузин не был образцом ученого старой школы, у него были какие-то отклонения от этических норм. История с Суслиным, история с Новиковым.}
 
  \sm

 Да, два случая были — он пытался присваивать результаты своих учеников....
 
  \sm
 
 {\it Есть широко распространенное мнение, что отношения власти с учеными в двадцатые — тридцатые годы шли по такой схеме, что под удар попадали те, кто воплощал в себе тип независимого, наиболее честного и благородного ученого. Тем не менее произошла история с Лузиным. Ведь Лузин, действительно, из математиков подвергся наиболее громким гонениям. Он не сел, остался жив, но громкая кампания шла почему-то вокруг его имени. Остальных арестовывали тихо.}
 
  \sm
 
 {\bf Это поношение, насколько я понимаю, шло не сверху, а это было возмущение части математической общественности.}
 С одной стороны были его ученики, которые его почитали и защищали, но с другой стороны были те, кто помнил Суслина, 
 помнил, как он с Урысоном поступил, запретил ему заниматься топологией.
 
 Привнесение идеологических моментов было духом времени, так тогда было принято формулировать обвинения, 
 это все присоединилось, когда вопрос о Лузине поднимался, но уже с другой стороны, не математической.

 Я думаю, что у него были довольно антисоветские взгляды, но он их, естественно, скрывал. Старался даже слишком «подкоммунивать».
 Я помню его выступление в Богословской аудитории (потом она называлась Коммунистической), это самая большая аудитория в здании МГУ 
 на Моховой. Он там рассказывал о своей поездке во Францию, как там плохо, что там нищие, бледные, еле живые, стоят с протянутыми руками 
 и т. д. Такая вот агитация, это он старался делать, его никто за язык не тянул.
 
 В 29-м году его избрали в академики, но не по математике, а по философии, для того чтобы избежать избрания некоторых других людей,
 которых не хотели избирать. Их и не следовало избирать, но они по указаниям сверху должны были быть избраны. И вот, чтобы уменьшить 
 число таких «назначенцев», по философии был избран Лузин, как меньшее зло. Это не явилось причиной затаенного недовольства.
 «Лузинщина» появилась позже, а в 29-м году он был еще в славе.
 
  \sm
 
 {\it Если считать, что хотели нападать на математиков старого закала, то кто тогда был математиком старого закала?}
 
  \sm
 
 Егорова к тому времени уже освободили, и он вскоре умер – в 31-м году. Сидел Каган Вениамин Федорович, кажется, он год только,
 и его отпустили. Но это было раньше, в 30-31-м году, а «лузинщина» — это 35-36 год.
 
 Из математиков в этих гонениях принимал участие П.~С.~Александров. Понтрягин занимал двойственную позицию.
 Лузин каялся, конечно. {\bf Математики, которые произвели атаку на него в то время, 
 недооценили общего политического положения.} Думаю, что злонамеренной их акция не была.
 Некоторые очень неблаговидные поступки Лузина действительно имели место, и их следовало разоблачить.
 {\bf Я думаю, что если бы они понимали, что человек может здесь сложить голову, то они, конечно, этого не сделали бы.}
 
  \sm
 
 {\it Такой идеологический накал кампании был, может быть, неожиданностью для многих?}
 
  \sm
 
 Думаю, что да. Я про себя могу сказать. Я не принимал участия в этой атаке на Лузина, я был только свидетелем,
 но я вполне сочувствовал этому нападению. Но, если бы я понимал в то время, чем это чревато, то, ясное дело, я не сочувствовал бы.....
 
  \sm
 
 {\it Разгар шельмования Лузина — это какой примерно год?}

  \sm
 
 Это раньше — 35-36-й годы. Мы тогда еще как-то не поняли. {\bf Я лично очень торжествовал, когда Лузина клевали},
 потому что его поведение было очень противным. Говорили, будто бы это дошло до Сталина и он сказал: ну, это ничего.
 Раз Сталин махнул рукой, то это сошло. Его выгнали из МИАН-а, он устроился в другом месте.
 
 \end{quotation}
 
 Напомню, что в 1937 году Ефремович был
 арестован и получил 10 лет. Колмогоров, Александров и Понтрягин, а также Соболев, А.~А.~Андронов, Петровский
 многократно предпринимали усилия по его освобождению \cite{GoncharovNehotin},
 и, возможно, отчасти преуспели - освободили Ефремовича чуть раньше положенного срока. Это
 к тому, что не стоит подозревать его в наивности, он просто был честный рассказчик. 
 
 \sm

{\bf  \punct Из публикации Дюгака.%
\label{ss:serpinski}}
Пьер Дюгак \cite{Dug} опубликовал несколько писем Серпинского,  Данжуа и Лебега,
имевших отношение к лузинской истории. Самое важное для нас мы приводим полностью:
\begin{quotation}
	 Письмо В.~Серпинского к А.~Данжуа (Arnaud Denjoy).
	 Варшава, 30 июля 1936 г
	
\bigskip

Дорогой коллега и друг,

Вернувшись с Конгресса в Осло, в котором не принял участия
ни один русский, я узнал ужасную новость, касающуюся нашего
друга Лузина. В московской газете «Правда» от 14 июля помещена
статья, озаглавленная «Враг, с которого сорвана маска». Посылаю Вам французский перевод — слово в слово — этой статьи.
Оказывается, что {\bf  по инициативе Павла Александрова и его сотрудников} 
была образована комиссия из членов Академии  СССР
{\bf для изучения очень тяжких обвинений, предъявленных Лузину его
бывшими «благодарными» учениками, которые требовали изгнания Лузина из Академии}.
Эти возражения были в основном политического толка. Лузина упрекают в том, что он публиковал самые
важные результаты во французских изданиях, в том, что он отказался подписать провокационный адрес к французским ученым — 
слово «Франция» при этом явно не называлось по очевидным мотивам. Все это верно, но это было именно в то время, когда
г-н Александров находился в Германии и публиковал свои результаты по-немецки и в немецких сборниках. 
Более того, только несколько месяцев тому назад, в Германии вышла большая книга г-на
Александрова и г-на Хопфа — который сам является немцем, живущим в Швейцарии.

Затем Лузина упрекают за его поведение до 1918 г., говоря,
что он был на стороне правых и что он так с правыми и остался.

Чтобы уничтожить Лузина морально, ему вменяют то, что он
систематически крал идеи у своих учеников, особенно у Суслина.
Каждый, кто хоть немного знает Лузина, понимает, что подобный
упрек не имеет никаких оснований.

{\bf Действия г-на Александрова и его компаньонов, конечно, есть
такая подлость, на которую, что даже невозможно представить,
был бы способен ученый.}

В своем письме от {\bf 27 июля 1935 г. — то есть вот уже год назад — г-н Лузин писал}:

«Возвращаясь теперь к очень трудной для меня самозащите по
поводу приписывания Суслину тех результатов, на которые он не
имел никакого права и которых у него даже в мыслях не было, я
должен сказать, что эта самозащита спровоцирована очень большой и совершенно реальной опасностью. Г-н {\bf Александров имеет
виды войти в Академию наук в качестве действительного члена,
сместив меня. С этой целью он требует пересмотра моих работ, заявляя, что я не имею права быть членом Академии, поскольку мои
идеи все украдены у Суслина. Такой пересмотр вполне возможен и
реален}».

Когда я был в Москве, в сентябре 1935 г., г-н Александров
заверил меня, что опасения Лузина — чисто мнимые, и что он очень
уважает Лузина, своего бывшего учителя. В моем присутствии
Александров протянул руку Лузину и объявил, что всегда будет
его другом.

Теперь имеются две проблемы.

Первая такова: как можно помочь Лузину? К сожалению, ввиду явных политических причин, всякое действие, происходящее от
польских ученых, могло бы только повредить Лузину. Я думаю,
что, может быть, французские ученые могли бы помочь Лузину,
заявив в совершенно дружеской манере послу СССР в Париже,
что упреки Лузину, кроме тех, что касаются его дружеского расположения к Франции, не имеют никаких оснований, и что, если 
Лузина исключат из Академии и заменят Александровым, 
это повредит престижу Академии [наук] СССР. В действительности,
 невозможно понять, что же выиграют Советы от всей этой омерзительной истории.

Что касается нас, поляков, мы, конечно, должны выступить с
резкой декларацией, но только в случае, если Лузин будет исключен из Академии [наук] СССР. Это надо потому, что Лузин —
иностранный почетный член Польской академии и что он является единственным советским доктором honoris causa польского
университета. Эта декларация, естественно, разорвет все связи
между польскими учеными и учеными СССР.

Вторая проблема — это проблема с Александровым.

Г-н Александров был приглашен Варшавским университетом
читать лекции осенью 1936 г. Естественно, теперь это приглашение
будет отозвано.

Однако я придерживаюсь того мнения и того же мнения мои
польские коллеги, что {\bf присутствие господ Александрова, Хинчина,
Колмогорова, Шнирельмана, которые самым нечестным образом
выступили против своего бывшего учителя и ложно обвинили
его, — нельзя терпеть ни в каком собрании честных людей. Они не
смогут, таким образом, участвовать ни в одном международном
конгрессе и ни в какой комиссии}. Я думаю, что подобная декларация, подписанная математиками всех цивилизованных стран, будет
необходима. Это совсем не является политической или социальной
акцией. Г-н Александров и компания совершили подлость, которая
должна рассматриваться как таковая, независимо от каких-либо
политических или социальных режимов. Терпимость ввиду такой
подлости была бы непростительной.
Примите, дорогой коллега и друг, выражение моих самых преданных чувств.

\sm

В. Серпинский

\sm

P.S. Мой адрес до 25 августа будет такой: Prof.W.Sierpi\'nski,
Czerchawa kolo Sambora, Польша.
\end{quotation}

Понятно, что Серпинский высказывает точку зрения самого Лузина
на подоплеку событий. Конечно, Лузин мог ошибаться, 
но, заметьте, он писал Серпинскому о предстоящей попытке исключить его
из Академии за год до того, как она началась! Единственным способом
изгнать Лузина из Академии было объявить его приносящим вред советскому государству.

Последовала переписка между Серпинским, Данжуа, Лебегом и Монтелем
(Paul  Montel).
Стороны были едины в отношении Александрова и сотоварищей
и обсуждали, чем они могли бы помочь Лузину. 12 августа 
Борель и Ланжевен (Paul Langevin) передали письмо в советское посольство
в Париже. В Москве, впрочем, к тому времени все уже кончилось.
Приведем лишь две цитаты на другую тему:
\begin{quotation}
Из письма А.~Данжуа к В.~Серпинскому
Мимизан-Пляж, 5 августа 1936 г.
Вилла Moнa Лиза (Ланд)

\sm

{\bf Замечательный взлет, которым Лузин ознаменовал предпоследний год (так как он стал теперь персоной, руководящей математиками в СССР)},
совпал со взаимным, очень явным
сближением между Францией и Советами. Разве 
не Лузин реализовал в Москве французское направление в науке, и разве франко-советское согласие не было бы нарушено в России, если бы он
пострадал?

\sm
	
	Из письма А.~Данжуа к А.~Лебегу
	Мимизан-Пляж, 7 августа 1936 г.
	Вилла Moнa Лиза (Ланд)

\sm
	
Примерно год назад Лузин получил большой пост. Он фактически руководил советским математическим отделением [Академии
наук]. {\bf У него была квартира из пяти комнат, автомобиль с шофером.}
	\end{quotation}

{\bf\punct Андре Вейль о деле Лузина.%
\label{ss:weil}}
И еще одна цитата из Дюгака:
\begin{quotation}
А.~Вейль (Andre Weil), которому я послал машинописное досье в сотню страниц о «деле» Лузина,
ответил 2 декабря 1978 г.:

«{\bf По разным поводам я не уставал говорить русским математикам, что к тому времени, о котором идет речь в Вашей работе, влияние Лузина в русской математической
жизни вот-вот станет гибельным (с чисто научной точки зрения), и что самые серьезные русские математики, и притом наименее подозреваемые в том, что они имеют неподобающие политические взгляды, оказались единодушны в этом отношении; в
сущности на это указывает то досье, которое Вы собрали}; но это также то, что такие
люди, как Данжуа, Монтель, Борель были совершенно не в состоянии оценить:
Монтель и Борель были в ту пору абсолютно чужды направлению современной математики, а Данжуа по другой причине был также этому чужд, так как он не знал
ничего из того, что делали другие. {\bf  Печально, конечно, что для того, чтобы положить
конец диктатуре, которую Лузин пытался ввести для математиков России, обратились к политическим и малопривлекательным методам. Однако надо заметить, что
это человек, чьи реакционные в политическом смысле взгляды были общеизвестны,
продолжал занимать очень влиятельную должность до самого 1935 г., и что, как это
вытекает из Вашего досье, с ним не случилось ничего худшего, чем серьезное предупреждение, которое он получил.} В частности, он остался членом Академии наук (это
положение не только почетное, как это происходит у нас, но весьма комфортное и
очень хорошо оплачиваемое)
\end{quotation}

{\bf\punct О том же. Другие источники.%
\label{ss:pontryagin}}
Вот Понтрягин \cite{Pon1}:
\begin{quotation}
Одним из мест кулуарных обсуждений был «математический салон», который в то время держала одна дама, весьма почтенного возраста.
Фамилию я её не помню, помню только имя и отчество: Софья Моисеевна. В этом «салоне» принимались различные математики, 
но в основном это были Л.~А.~Тумаркин, А.~О.~Гельфонд, Л.~Г.~Шнирельман, а затем присоединился и я....

Именно здесь мы обсуждали событие с Лузиным, и у некоторых возникло подозрение, что {\bf среди организаторов} письма мог быть П.~С.~Александров. 
Мы хорошо знали, что у Александрова с Лузиным были отвратительные отношения. Хотя Александров был ученик Лузина,
они грубо враждовали между собой. Так как источником информации об этой вражде был Александров, естественно, я был его сторонником....

Александров был учеником Лузина. Он, несомненно, очень многому научился от Лузина и в течение длительного
времени был под его сильнейшим влиянием. Но к тому времени, когда я стал разговаривать с Александровым о Лузине,
эти два человека находились уже в состоянии непримиримой ненависти%
\footnote{Ко всему этому была сноска Понтрягина:
	\newline{\it 
	 Тогда мы ещё не знали, что грядёт 1937 год. Я описываю события так, как я и мои товарищи воспринимали происходящее тогда — в 1936 году. Позже я понял, что Советскому правительству нужно было разогнать школу русского математика Н.Н. Лузина. Уничтожить его самого они не решились.}
	\newline
	Стоит отметить, что разгонять школу Лузина в 1936г. не было никакой необходимости: Лузина атаковали именно представители его школы.}....
\end{quotation}

Процитируем также дневники Вернадского \cite{Vernad-dnevniki}:
\begin{quotation}
6.VII.936 
	(Узкое) 
	
	Вчера вечером приезжал Зел[инский]. Лузин в отчаян[ном] положении]. Он был в ред[акции] «Правды» ред[ актора] не было.
	Заместитель ред[актора] заявил, что они в курсе всего, за все «факты» ручаются, давно за ним следят. 
	Он подал в отставку от (должности) председателя] группы%
	\footnote{В то время академиками были Бернштейн, Лузин, Виноградов, А.~Н.~Крылов, Н.~М.~Крылов, Шмидт (он был избран по отделению геологии). Члены-корреспонденты Александров, Голубев, А.~И.~Некрасов, В.~И.~Смирнов, Соболев, Шнирельман. Быть может, я кого-то пропустил.
	Кроме того, в состав группы входили лица, не имевшие академических чинов.}. Сейчас же было собрано заседание группы.
	Отставка принята подавляющим большинством%
	\footnote{То есть кто-то был против.}. Поднимается вопрос об удалении его из Акад[емии]. {\bf Л[узин] думает, что кампанию против него ведет чл[ен]кор[респондент] Александров, желающий занять место председателя группы.} 
	
	Сервилизм и страх царят. Н.~Д. (Зелинский) и Нас[онов] думали об обращении] в През[идиум] (АН СССР).
	Мне кажется, надо в президиумы] отделений и копию непременному] с[екретарю] по каждому отделению иначе «скоп» и сговор. 
	Кржиж[ановский] Лузина не принял, был «занят заседанием»! 
	Партийные здесь вполне в таких обстоятельствах бессильны и боятся всего больше нас. Л[узин] должен обратиться к Сталину.
	В сущности, вопрос стоит о газете «Правда», т[о] е[сть] введенной (в заблуждение) или пошедшей на обман. 
	Здесь вчера разговор с Архангельским. Он считает, что это уже решение фактически правительства, т[ак] к[ак] неподписанная]
	статья «Правды» органа партии имеет такое значение. 
	
6.VII.936 Узкое	
	
Все боятся и боязнь имеет основание, т[ак] к[ак] защититься без большой смелости нельзя. Газетная статья «Правды» полна инсинуаций и передержек. Он  должен был бы подать прокурору и написать Сталину. Но едва ли у него хватит решимости: он неврастеник. Боюсь трагедии. Людей не берегут и на каждом шагу это видно. 	
\end{quotation}
 
 Независимо от Серпинского, приводится точка зрения Лузина 1936года. Про значение неподписанной статьи не очень ясно,
 она была на третьей газетной странице.
 
 Теперь цитируем статью Юшкевича <<Дело академика Лузина>>, но не ту которую все цитируют,
 \cite{Yush2},
 а другую  \cite{Yush1}:
\begin{quotation}
 Как видим, молодое поколение математиков энергично выдвигалось на руководящие позиции. Процесс этот, сам по себе естественный,
 в тогдашней обстановке приобрел весьма зловещий характер. Часть выдвинувшихся на рубеже 20-х и 30-х годов молодых ученых считала, 
 что старшее поколение — 
 Д. Ф. Егоров и Н. Н. Лузин в особенности — отстают от новых направлений науки, не должны далее возглавлять математическую жизнь ни в Москве, 
 ни в Советском Союзе вообще. {\bf Тут играли роль и чисто личные моменты, например, всем известные напряженные отношения между Н. Н. Лузиным и
 П. С. Александровым.} 
\end{quotation}

И еще раз напомним, ни один пролузинский текст 1950-1985гг не содержит  намеков на то, 
что атака на Лузина велась сверху.

\sm

{\bf \punct Показания С. П. Новикова-младшего.%
\label{ss:novikov-hvorostin}}

\begin{quotation}
Расследование было проведено тогда отцом (кажется, вместе с Люстерником и
 Лаврентьевым, знающим партийные круги). Они установили, что было письмо П. Александрова
к влиятельному человеку по имени Хворостин, с возмущением излагающее мерзости Лузина.
 Хворостин находился в Саратове и имел большие связи в ЦК. Лузина он ненавидел, это
было известно. Хворостин-то, как они решили, и передал материалы в ЦК и инициировал
статью. Павел Сергеевич был великий мастер биллиардного удара%
\footnote{
	Для сравнения, взгляд Лузина на Александрова.
	Лузин -- М.Я.Выгодскому 14.09.1934 \cite{ErmTok}.
	\newline
{\it ... Правда, нравы теперь иные и
	вместо скрытых протоколов, где фиксируется полный отвод меня
	от науки, от моей дорогой Теории Функций, я имею лестное 
	заявление в «Правде» и «Известиях»
	\newline
	Но разве это меняет дело и разве я не вижу тени великого Талейрана.}
\newline
В первой фразе речь идет о дележе квот в Математической энциклопедии, которая должна была издаваться 
под редакцией Выгодского. В «Известиях» была короткая заметка Колмогорова (как директора НИИММ)
о встрече советских и французских математиков в Москве, она полностью приведена в 
 \cite{ErmTok}. Что за заявление в <<Правде>> -- не ведомо. K кому относится слово <<Талейран>>,
 как будто, сомнений не вызывает. В том же письме 
 <<{\it ... разоблачение лжи, коварства, притворства шагов Талейрана -- Вы
 знаете кого я разумею под этим именем}>>.}! 
\end{quotation} 
По-моему, нет оснований не доверять С.~П.~Новикову \cite{NovS} в данном вопросе. Относительно
малая политизированность первой статьи в Правде указывает, скорее, на то, что она исходила из
математических кругов. Некоторая дополитизация статьи должна была произойти сама собой
(уже при проходе через Хворостина, о котором речь еще пойдет, 
а потом и неизвестного нам автора
статьи).

 Конкретные обвинения из статьи почти соответствует
направлению, которое имела лузитанская атака, застенографированная в \cite{Sten} (мы увидим, что по одному пункту позиции лузитан различались, а именно Александров и Хинчин возражали против обвинения в печатании
лучших статей за границей).  

\sm

{\bf\punct Что говорил сам Александров?%
\label{ss:aleksandrov-1936}}  В 1977 году он написал статью о Лузитании в журнале <<Квант>> (математический 
журнал для школьников):
\begin{quotation}
В августе 1924 года этот коллектив пережил тяжелый удар:  погиб один из самых
талантливых членов Лузитании   и вообще один из самых талантливых
советских математиков - Павел Самуилович Урысон.

Члены Лузитании восприняли этот удар судьбы как личное горе каждого из них, а Лузитания
как целое ненадолго его пережила...
 \end{quotation}
 
 Идет ли речь о причине и следствии или просто о последовательных событиях - понять
 невозможно. В принципе, конспирология по поводу гибели Урысона была -- что
он не просто так пошел купаться в шторм, но все равно не ясно при чем тут может быть
Лузин. 

Потом Павел Сергеевич опубликовал свои воспоминания в двух номерах <<Успехов>>
\cite{Alex-auto1}--\cite{Alex-auto2}:
\begin{quotation}
Однако, по моему убеждению, годами высшего расцвета Н.~Н.~Лузина
как математика да и как человеческой личности были не годы Лузитании,
а непосредственно предшествовавшие годы: 1914—1915.
Узнав Лузина в эти самые ранние творческие его годы, я узнал действительно
 вдохновенного ученого и учителя, жившего только наукой и только
для нее. Я узнал человека, жившего в сфере высших человеческих духовных
ценностей, в сфере, куда не проникает никакой тлетворный дух.
Выйдя из этой сферы (а Лузин потом вышел из нее), человек неизбежно
подпадает под власть тех сил, о которых Гёте сказал:

Ihr f\"uhrt in's Leben uns hinein,

Ihr lasst den Armen Schuldig werden

Dann \"uberlasst Ihr ihn der Pein,

Denn jede Schuld racht sich auf Erden.

Вы вводите нас в жизнь,

Вы делаете беднягу виновным.

Потом вы предаете его на муку,

Потому что на земле отмщается всякая вина.
\end{quotation}

Во второй части говорится следующее:
\begin{quotation}
Позже, когда
я вернулся к математике и стал снова продуктивно ею заниматься, нарушенные отношения  между мною и Н.~Н.~Лузиным уже не смогли восстановиться,
но эти нарушения отношений имели характер совершенно не типичный
для отношений между учеником и учителем и для нас здесь интереса не
представляют.
\end{quotation}

В общем, Александров подтверждает факт <<особых>> отношений,
уклоняясь от обсуждения каких-либо подробностей.

И, наконец, отрывок из интервью Александрова 1971г., \cite{Alexandrov-Duvakin}.
\begin{quotation}
 Высший расцвет Лузина — это, я бы сказал, десятилетие: с 1915-го по 1925-й год. Потом пошло уже...
 
В. Д.: Под уклон, да?

П. А.: ...под уклон. А уж когда его избрали в академики, то это было резким уклоном, 
как и некоторой деградацией всей личности в целом.

В. Д.: Ах, так? С чем же это было связано?

П. А.: Ну, не стоит это все подымать, это, так сказать, не предназначено для истории.
\end{quotation}

\sm

{\bf\punct За пять месяцев до июля.%
\label{ss:luzin-krylov}}

\begin{quotation}
	Н.~Н.~Лузин — А.~Н.~Крылову, Москва, Арбат, д. 25. кв. 8 
	12 февраля l936 г.

	Лично я с огромной охотой ушел бы из председателей 
	Математической группы, но упустил момент. Сейчас ввиду 
	отстранения Ивана Матвеевича [Виноградова] и  
	неожиданной для всех нас замены его ленинградским коллегою, 
	не проживающим в Москве, мне бросить все и уйти  
	представляется затруднительным. Но помимо моей печальной 
	деятельности «философа от математики», у меня ведь есть 
	и технические работы (движение поезда, периодограммы, 
	etc), и я с огромным бы удовольствием ушел просто в 
	члены Технической группы, чем {\bf быть в центре  
	столкновения математиков, имеющего случайный и чисто личный  
	характер}. По проблема ухода очень трудна, и я не знаю, как 
	ее решить, кроме разве уйдя вообще из числа членов  
	академии. Я уже начал присматриваться
	 к нынешнему Московскому университету и дал согласие на чтение в нем 
	этою весною эпизодического курса. 
\end{quotation}

Не известно, о каком <<столкновении математиков>>, в которое втянут Лузин, идет речь. Но, во всяком случае, речь идет о напряженной обстановке, 
и, по-видимому, Лузин считает, что должен кого-то защищать... Что за отстранение Виноградова (с должности директора Стекловки??) 
и кто такой <<его ленинградский коллега>>\footnote{Естественно предположить, что 
	это начальник ленинградской части
Стекловки, наличие каковой должности до 1940г. анналами отрицается.}, история тоже умалчивает.  Директором Стекловки Виноградов в итоге остался.

\sm

{\bf \punct Важное уточнение.} Вроде бы всё сходится на том, что инициатором атаки был Александров.
Но запущенный им процесс уже шел через партийные круги, имевшие свои взгляды и свои интересы.
Надо иметь в виду, что сентенция <<публикует  свои лучшие работы за границей>> к Александрову 
могла быть прицеплена к Александрову, еще в большей степени, чем к Лузину%
\footnote{У Александрова в 1930-35гг. вообще не было опубликованных в СССР научных статей, 
	кроме выступлений на математических съездах
(но были две   книги).}.
Ниже, но не в этом параграфе, мы увидим, что такая точка зрения на то, где нужно публиковаться, присутствовала
в определенных московских математических кругах. В связи с этим не исключено, что 
в подготовке выступления <<Правды>> участвовал еще какой-то <<теневой>> математик (который не был вычислен Лузиным).

\sm

{\bf\punct Позиция нападения.%
\label{ss:napadenie}} Перейдем к обсуждению Академической стенограммы.

\begin{quotation}
	ХИНЧИН. А как статья Устава Академии наук говорит об исключении?
	
	ГОРБУНОВ. Статья Устава об исключении гласит следующее: «{\bf Действительные члены},
	почетные члены, члены-корреспонденты Академии наук... {\bf лишаются своих званий 
	постановлению Общего собрания, если деятельность их направлена во вред Советскому Союзу}».
	
	БАХ. Ну, да, это конечно будут эксплуатировать за границей, будут сопоставлять с
	Конституцией.
	
	ГОРБУНОВ. Да, но мы все так по Конституции и делаем.
	
	БАХ. Но все-таки тут у нас нет квалификации со стороны такого органа как Наркомвнудел [НКВД].
	{\bf  Все это вращается в той области,
	что мы имеем здесь перед собою врага, а на
	меня скорее это производит впечатление, что это недоброжелательно, что это не друг.
	Ведь враг — это нечто активное.}
	
	КРЖИЖАНОВСКИЙ. Мы сейчас живем в такую полосу, когда те, кто не с нами, те против
	нас.
	
	СОБОЛЕВ. В сущности, сейчас, когда мы видим широчайшую демократию нашего государства в целом,
	мы не можем подходить с той точки зрения, что здесь нужны доказательства Наркомвнудела и т.д. 
	Сейчас важнейшие мероприятия обсуждает весь народ в целом
	и сейчас пора подходить к такому вопросу, как исключение, необязательно по постановлению Наркомвнудела,
	а по нашей собственной инициативе. Разница во внутреннем положении сейчас такова, 
	что мы имеем право решать в этих случаях сами, без всяких инстанций. По-моему, это очень важное обстоятельство, что мы имеем собственную инициативу.
	
	БАХ. Существует определенное государственное учреждение для этого. Ведь у нас
	нет еще всеобщего суда, т.е. ставится на общественное обсуждение тот или другой вопрос, и что общественное обсуждение поставит, то и будет.
	Это неправильно. Пока мы
	существуем в таком государстве, где государственная власть вполне организована, и она
	сможет даже не послушаться общественного мнения, как можно ставить таким образом
	вопрос? Мы никаких оснований не имеем так его ставить. У нас есть крепкое государство. Поэтому нельзя так ставить вопрос.
	
	СОБОЛЕВ. {\bf Но ведь мы же не арестовывать его собираемся.}
	
	БАХ. Тов. Соболев говорит, что при нашей демократии этого достаточно, а мне кажется, что этого не достаточно.
	
	ХИНЧИН. По Уставу для исключения достаточно установить, что деятельность направлена во вред Союзу ССР. 
	И неужели Академия наук в иных случаях не может самостоятельно, без помощи других органов, установить,
	что деятельность такого-то лица направлена во вред Союзу? {\bf Деятельность того или другого лица может быть направлена во вред
	Союзу и в тех случаях, когда эта деятельность такого рода, что НКВД не сочтет своей
	обязанностью вмешиваться в нее.}
\end{quotation}	
	
Итак, Лузина предполагалось обвинить в нанесении вреда Советскому Союзу,
и на этом основании лишить звания академика. 

Быть может, взгляды части нападавших проясняются следующими словами Хинчина:
\begin{quotation}
	ХИНЧИН.
	Я не считал и не считаю до сих пор, что {\bf звание академика — это кресло в
	Совете, и это не должно быть даже синекурой или титулом, или оплатой за прежние заслуги, а это должно быть постом организатора. И я считаю, что Н. Н.Лузин, при всех его
	прежних заслугах, которых я не отрицаю, к организаторской деятельности — при всех его
	личных свойствах — не способен.}
\end{quotation}

Стоит иметь в виду, что академиков тогда было мало (98 душ на весь СССР, и это число
покрывало все естественные и гуманитарные науки, а также инженерных академиков).

\sm

{\bf\punct Позиция защиты.%
\label{ss:zashchita}} Она была сразу озвучена Бернштейном: 

\begin{quotation}
БЕРНШТЕЙН. Если угодно, я выскажу свое отношение к Николаю Николаевичу, прежде всего как к ученому. Я считаю, что мы, прежде всего, должны учесть именно то, что
мы имеем в данном случае человека исключительно крупного научного значения, имя
которого пользуется большим авторитетом на Западе — вполне заслуженным.
 Доказывать сейчас его научные заслуги на этом заседании я не считаю необходимым. Важно
только, чтобы это вообще было зафиксировано, чтобы мы собрали соответствующий
материал, который в полной мере это подтвердил бы. \dots

Второе положение, которое является чрезвычайно важным, чрезвычайно существенным, которое мы должны выяснить, это вопрос о том, 
{\bf можно ли считать, что в действиях
H. H. имеется элемент сознательного вредительства? Имеется ли в тех материалах,
 которыми мы располагаем, достаточно для этого указаний? Я ставлю вопрос именно о сознательном вредительстве в отличие от того, что, может быть, фактически его выступления были вредными. И в этом отношении я, прежде всего, должен подчеркнуть — я считаю, что его выступления, его легкомысленные отзывы, являющиеся результатом недостаточно продуманного отношения к той ответственности, которую несет академик,
высказывая определенные суждения, — такого рода выступления являются, безусловно,
объективно вредными\dots
 у него ни в
какой мере не было стремления при этом поддерживать людей мало талантливых, ни в
какой мере не было желания повредить этим развитию математических наук}. К развитию
математики в нашем Союзе H. H. вообще относился и относится всегда с чрезвычайно
большим интересом и делал все возможное для того, чтобы математика у нас развивалась, возможно, лучше так, как он это понимает.
Во всяком случае,  прошлая его деятельность в Университете [и] та роль, которую он
сыграл в развитии математической школы в Москве, является исключительно крупной, и
этого, мне кажется, не следует забывать, и это, вероятно, в большей степени будет подчеркнуто и другими, когда мы будем выносить какое бы то ни было окончательное суждение о Николае Николаевиче.

Поэтому {\bf я считаю, что те обвинения, которые возводятся против H. H., имеют основание в некоторых его поступках}.
Но те выводы, которые из этих обвинений делаются,
являются необоснованными, поскольку они перерастают, как мы, математики, говорим,
экстраполируют, выводят нечто гораздо большее, чем то, что можно извлечь из этих
фактов. Те утверждения, которые я сейчас делаю, нуждаются, конечно, в более полном
обосновании. Необходимо различные факты, которые здесь отмечены, соответствующим образом проанализировать,
и каждый в отдельности соответствующим образом
оценить. {\bf Я уверен в том}, что Комиссия, отнесясь к этому делу с полной беспристрастностью, сумеет в этих фактах выделить то, 
что является действительно следствием тех личных недостатков характера H. H., о которых я говорил, и сумеет доказать, {\bf что те, более
или менее широкие выводы, которые здесь пытаются сделать, указывая, что H. H. является врагом советской власти, сделать ни в коем случае невозможно.}
\end{quotation}

То есть защита согласна признать часть обвинений, отрицая то, что
для исключения Лузина из Академии есть какие-либо основания. 
Бернштейн попытался высказать и такой тезис:
\begin{quotation}
	БЕРНШТЕЙН. \dots
	Но я подчеркиваю, что
	независимо от этого, независимо от того, что подобного рода вещи могли бы отрицательно характеризовать H.~H.,
	тем не менее, у него есть ряд неоспоримых научных заслуг, где
	никто не претендует на то, что он заимствовал, начиная с его известной докторской диссертации, и которые, во всяком случае, 
	являются {\bf достаточным основанием для того,
	чтобы признать его крупнейшим советским математиком.}
\end{quotation}

Но удержаться на этой позиции (<<крупнейший советский математик>>) было невозможно, и под напором атакующих
она была оставлена.

\sm

{\bf\punct Предварительные замечания по поводу отдельных участников.}
На заседаниях комиссии Лузина так или иначе атаковали почти все
члены комиссии (кроме Бернштейна и молчавшего Виноградова) и почти все участники прений.

\sm

1) Лузин уже был сведен с должности председателя Математической группы, и понятно, 
что Комиссия должна была его в той или иной форме осудить.
Основной вопрос был в том, будет ли Лузин изгнан из Академии. 

Бах, судя по разным его репликам, не был доброжелателен по отношению
к Лузину. Но мы видели, что именно по ключевому вопросу он сдерживал  нападение. 
Не выглядит доброжелательным и Ферсман.
Но в какой-то момент он произносит длинную тираду с таким выводом
\begin{quotation}
	ФЕРСМАН. ...  Но это вопросы, которые требуют не
	внешнего решения, не внешней защиты, а глубокой какой-то внутренней переработки.
	Вы сами, вероятно, найдете пути, которые смогут обеспечить это. Здесь дело не во внешних словах.
\end{quotation}
Это не самое страшное из того, что мог бы услышать Лузин.

\sm

2) На Лузина обрушился поток научно-этических и политических
обвинений со стороны математиков. Оказывается, что разные участники нападения предъявляли
политические обвинения с разной степенью решительности. Некоторые
же их и вовсе не предъявляли, а один из участников нападения защищал
Лузина в этом отношении.

\sm

Позиции нападавших математиков довольно сильно различались, были две группы,
Шнирельман--Люстерник--Гельфонд (а также Б.~И.~Сегал и С.~Е.~Аршон) и Александров--Колмогоров, 
кроме того,
по отдельности выступали
Хинчин, Соболев, Шмидт.

Перейдем к обсуждению научно-этических обвинений, политические обсуждаются в следующем параграфе.

\sm

{\bf \punct История с Суслиным.%
\label{ss:suslin}}
В своих воспоминаниях Александров излагает эту историю так:
\begin{quotation}
В 1918 г.
Н.~Н.~Лузин на время переехал в Иваново (тогда еще Иваново-Вознесенск).
По его совету туда переехали А.~Я.~Хинчин, Д.~Е.~Меньшов и М.~Я.~Суслин;
они, как и Н.~Н.~Лузин, преподавали в Ивановском Политехническом институте. Однако Суслин в Иванове не ужился и скоро потерял там свою работу.
Ввиду этого, по инициативе В.~В.~Голубева и И.~И.~Привалова, возник план
о предоставлении Суслину профессорского места в Саратовском университете.
Ожидалась рекомендация Н.~Н.~Лузина. Лузин ее не дал и не поддержал предоставления 
Суслину преподавательского места в Саратовском университете.
Не получив этого места, Суслин уехал к себе в деревню (в Саратовской
губернии). Вскоре он заболел сыпным тифом и умер. В историю советской мате­матики
была вписана одна из самых трагических ее страниц. На письменном
столе Лузина до конца его жизни стоял портрет Суслина, единственный портрет Суслина, который я видел. 
\end{quotation}

Об этом пишет  и Понтрягин в своих воспоминаниях,
ссылаясь на независимый от лузитан источник:
\begin{quotation}
	При поисках работы он ездил по провинциальным университетам, но оказалось,
	что во всех этих университетах уже имеется письмо Лузина, в котором он резко
	отвергает кандидатуру Суслина как преподавателя. О наличии такого письма рассказывал мне А.~А.~Андронов%
	\footnote{Андронов Александр Александрович (1901--1952) --  математик
		(известный работами по качественной теории дифференциальных уравнений) и механик. Работал в Нижнем Новгороде/Горьком.}. 
	Оно было в университете в Горьком.
\end{quotation}

Вернемся к обсуждению 1936 года:
\begin{quotation}
		ХИНЧИН.
{\bf  Но мне кажется, что в «Правде»
имеется фактически не совсем верное утверждение.} Во-первых, неверно, что Лузин 
сделал все возможное, чтобы выставить своего ученика Суслина из Москвы. Лузин перевел
Суслина из Москвы в Иваново-Вознесенск, но перевел его в качестве своего сотрудника.
Причем надо сказать, что условия работы в Иваново-Вознесенске были тогда во многих
отношениях лучше, чем в Москве. Я сам там был. После Суслина я через год приехал.
Нельзя это рассматривать как выживание. Наоборот, это можно рассматривать как услугу. 
Пробывши год в Иванове, Суслин отсюда уехал. Я доподлинно не знаю, в чем тут
дело. Но когда я приехал после него в Иваново, то единодушное мнение в Иванове было
такое, что из Иванова Лузин действительно его выжил. Суслин выехал в Саратов на свою
родину. И известно, что Н.~Н.~Лузин послал телеграмму ректору Саратовского университета,
настаивая на том, чтобы Суслину не была предоставлена профессура, несмотря
на то, что Суслин имел уже мировую известность. Суслин вскоре после этого умер там
от сыпного тифа...........

ХИНЧИН. \dots Я знаю, что там были очень резко враждебные отношения.
 H.~H. очень скрытен в этом отношении, хотя он и бывал со мною откровенным.
Но каждый раз, когда речь заходила о Суслине, у него против воли прорывалась такая
инстинктуозная, резкая враждебность, так что приходилось верить на слово. Но тут ходили
всякие слухи. Говорили, что тут супруга Н.~Н. сыграла большую роль, выжила Суслина.

АЛЕКСАНДРОВ. Я считаю, что в борьбе за организацию славы H.~H., которую взяла на
себя супруга H.~H., это был источник всяких зол%
\footnote{Вообще и эту идею, и просто отрицательное отношение к
	<<супруге Н.~Н.>> высказывают многие источники, в том числе доброжелательные 
	к Лузину.}.

ХИНЧИН.
{\bf Суслина называют учеником Н.~Н.~Лузина, загубленным H.~H. Ну, когда человек
умирает от сыпного тифа, то это слишком резкое выражение. Ведь он мог заболеть
сыпным тифом и в Иванове. Но общее мнение такое, что из Иванова H.~H. выжил Суслина.}
\end{quotation}

Можно сомневаться в речах недоброжелателей, но вот что говорит сам Лузин:
\begin{quotation}
	 ЛУЗИН.
 Будучи в годы разрухи приглашен в Иваново-Вознесенский Политехнический институт, я взял с собою Суслина и устроил ему должность экстраординарного профессора,
причем Институт потребовал от него сдать в течение двух лет магистрантские экзамены.
К сожалению, я не смог заставить Суслина всецело посвятить себя научной работе, и по
прошествии двух лет, ввиду не сдачи им магистрантских экзаменов у него возникли 
трения с дирекцией Института, и он вынужден был, покинув Иваново-Вознесенск, уехать к
своим родителям. Не скрываю, что на почве моих требований о подготовке к экзаменам
и отсутствия у него достаточного импульса это сделать, у нас произошли трения. Я и
сейчас, 17 лет спустя, считаю недопустимым, для столь талантливого человека, как Суслин,
загружение педагогической работой в ущерб научной работе. {\bf Когда Суслин захотел
стать профессором в Саратове, я выразил свое мнение о том, что ему не следует давать
там кафедры, так как я находил, что это был единственный способ для того, чтобы он
был вынужден, хотя бы под давлением внешних обстоятельств, приступить, наконец, к
серьезной научной работе.} При этом я хорошо знал, что родители Суслина вполне
обеспеченные люди, и он в заработке не нуждается».

ФЕРСМАН. Ну, а если бы он не остался в Саратове, если бы ему в Саратове отказали в работе в университете, он мог бы вернуться в Иваново, за ним там оставалось
место?

ЛУЗИН. Тогда была предложена такая формулировка, что до тех пор, пока он не сдаст
экзаменов, он, стало быть, увольняется.
\end{quotation}

Так что был и конфликт (мы не знаем его подробностей) с выживанием Суслина из Иванова,
  и антирекомендательные письма.
А оправдания Лузина не вызывают сочувствия к нему.

\sm

Что касается обвинения Лузина в смерти Суслина, то как
верно сказал Хинчин, оно несправедливо.
На время Гражданской войны приходятся эпидемии
испанки, трех болезней, называвшихся словом <<тиф>>,
сыпной тиф, брюшной тиф, возвратный тиф,
была также высокая заболеваемость 
малярией и дизентерией. Принято считать, что эпидемия испанки 1918 года унесла
в Европе больше жизней, чем Первая мировая война%
\footnote{Кстати,  упоминавшийся выше Janiszewski умер от испанки
в начале 1920г.}.
Огромная смертность в России времен Гражданской войны объясняется
прежде всего этими эпидемиями%
\footnote{Общей статистики на этот счет по понятным причинам нет,
хотя есть данные по отдельным местностям и отдельным промежуткам 
времени.
 Cтатистика по Красной армии проводилась
 Генштабом, она известна
\cite{Krivosh}:
Общее количество умерших от инфекционных болезней
283079. Убито и умерло от ран 146488 (в два раза меньше). Если у ним добавить
пропавших без вести, попавших в плен, погибших
в результате происшествий, покончивших самоубийством, и т.п.,
то цифра возрастает до 376269 (из коих не все погибли). То есть даже по воюющей армии
смертность от эпидемий составляла почти половину общих потерь.}

\sm

Биография Суслина исследовалась В.~И.~Игошиным,
нашедшим в архивах ряд документов, связанных с Суслиным,
\cite{Igo1}, \cite{Igo2}. Наличие конфликта с администрацией Иваново-Вознесенского политеха 
понятно из опубликованных документов, но детали не ясны.
Суслин в этот период имел проблемы со здоровьем, и много его
сил отнимали хозяйственные заботы.
  О его поездках по провинциальным
университетам после увольнения из Иваново-Вознесенска.
(о которых говорит Понтрягин) данных нет, перед самой смертью
Суслина в его деревню пришло письмо, подтверждающее его оставление  в Московском
университете, оно было в связи со смертью получателя отправлено обратно
(и было найдено в архиве Московского университета)%
\footnote{Из Воспоминаний Александрова:
\newline
		«М.~Я.~Суслин уже в самые ранние студенческие годы проявил 
		себя как интересный и своеобразный человек.
		Уже в 18—19 лет он составил 
		себе особую программу своего дальнейшего интеллектуального 
		развития. Математика была только началом этой программы. 
		Вторым этапом должны были быть физика и химия, за которыми 
		должна была последовать биология. В качестве завершения 
		программы мыслилась медицина, которой М.~Я.~Суслин и  
		предполагал посвятить всю свою дальнейшую жизнь.}.

\sm

{\bf\punct Обвинения в присвоении чужих результатов.%
\label{ss:plagiat}}
На комиссии Лузин обвинялся в присвоении результатов
своих учеников Новикова, Лаврентьева, Суслина. В июле
1936 года Новиков и Лаврентьев вяло подтвердили факт плагиата,
сделав это так, что это
  не усилило позиции нападающих. 
Позднее они оба отстаивали память Лузина, а Лаврентьев
был председателем комиссии по изданию трудов Лузина.

Однако на этот счет есть показания их сыновей.
Из воспоминаний М.~М.~Лав\-рен\-тье\-ва-млад\-ше\-го \cite{Lavr-jr}:
\begin{quotation}
	Отец с восхищением вспоминал, каким блестящим лектором был Лузин и какие оригинальные задачи придумывал для учеников.
	Но у него был и большой недостаток. Часто происходил такой казус: когда ученик успешно решал поставленную учителем задачу,
	учителю почему-то становилось обидно. И он говорил: на самом деле задачу решил не ученик, а он сам, Лузин.
	Из-за этого у Лузина случалось много конфликтов. Обидел он и моего отца. Однажды отец нашел красивое решение для задачи, поставленной Лузиным.
	Лузин сказал: «Отправляйте свой результат во Францию» (французская математика в те годы была авторитетнее нашей).
	А сам отослал французам и свою статью, которую напечатали быстрее. И написал французским коллегам, что молодой Лаврентьев его обокрал. 
	
	Для отца это было моральной травмой. Он сразу сменил научную тему, занялся аэродинамикой в группе Чаплыгина и отошел от Лузина,
	перестал с ним общаться. \dots
	
	У отца осталось неоднозначное отношение к Лузину, хотя в своих публичных выступлениях он всегда вспоминал об учителе с восхищением.
	Он говорил, что перенял от Лузина способность ставить интересные задачи и привлекать талантливую молодежь. 
	Но отец никогда не обижал учеников и скорее готов был отдать им свои научные результаты, чем присвоить чужой труд.
\end{quotation}

 В своих воспоминаниях Лаврентьев-старший \cite{Lavr2} много  положительно
 говорил о Лузине. Но там есть и такая фраза:
\begin{quotation}
Мой дипломный результат Лузину понравился, он даже включил его в свою книгу, но про меня забыл.
\end{quotation}

Из записок С. П. Новикова-младшего:
\begin{quotation}
	\dots
Он [Лузин] крал у П.~С.~Новикова, моего отца\dots
\end{quotation}

Понятно, что этих воспоминаний в руках членов комиссии не было.
Новиков и Ляпунов написали в Комиссию большую бумагу, где плагиат
вроде бы и подтверждался, но так, что это  не усиливало нападавших.
Далее Новиков исчез.
\begin{quotation}
КРЖИЖАНОВСКИЙ. Новиков в отпуске... по телефону не отвечает.
\end{quotation}

Нападающие пытались 
 обосновывать плагиат в отношении Новикова, и, отчасти,
 в этом преуспели.
\begin{quotation}
АЛЕКСАНДРОВ. \dots Я начну с работы Новикова. Действительно, была очень неудачная
	публикация, заключающаяся в том, что некоторый результат, опубликованный в работе,
	совместной работе Вашей и П.~С.~Новикова, был затем опубликован в Вашей собственной ноте. 
	Я читал эту работу Новикова и вашу собственную действительно с большим
	интересом. И все-таки их содержание более или менее идентично. А теперь вы приводите объяснение проф. Бореля, 
	где он пишет, что те 3 строки, в которых дается такое важное
	признание с Вашей стороны, им вычеркиваются. (Зачитывается французский текст.)
	Я должен сказать, H. H., что здесь какое-то недоразумение. Я отказываюсь, насколько
	я знаю обычаи математических публикаций, верить, чтобы Борель или кто-нибудь другой, 
	мог вычеркнуть три строчки такого содержания «...» указав при этом «...» (французский текст)
	\end{quotation}
Лузин пытается защищаться, по-моему, совсем не убедительно
(желающие могут прочитать Стенограмму).

О другом случае говорит Колмогоров:
\begin{quotation}
	КОЛМОГОРОВ. \dots Но
	я хочу поставить более существенный общий вопрос и получить на него общий же ответ: считает ли H.~H., что такое положение, когда в его работах печатались работы Новикова, когда докторская диссертация Новикова до сих пор не напечатана, — является нормальным положением?
	
	ЛУЗИН. Отвечу следующим образом. Андрей Николаевич, Вы ставите мне в вину, что
	диссертация Новикова не была опубликована. Содержание этой диссертации исчерпывающе полно было опубликовано мною с авторским именем в моей книжке [имеется в виду \cite{Luzin-descript}].
	
	КОЛМОГОРОВ. Я и спрашиваю: разве это нормально?
\end{quotation}

Да ненормально, но, строго говоря, это не плагиат. Лузин излагает  результаты Новикова, воздавая тому должное  и признавая за ним авторство результатов. Непонятно только, почему фамилии Новикова нет на обложке книжки, в которой излагаются новые результаты Новикова.

Говорилось еще многое, но  нападающие смогли обосновать лишь часть своих обвинений.
\begin{quotation}
КРЖИЖАНОВСКИЙ.
  Возражений  нет?  Принято. По-моему,  не  следовало  бы  говорить  относительно  явного  плагиата,  а сказать —
«особенно недобросовестным является отношение в случае Новикова».
\end{quotation}

Заметную часть Стенограммы составляет обсуждение открытия <<аналитических множеств>>, где нападающие (в основном Александров)  стараются обвинить Лузина в частичном присвоении результатов Суслина.
Я воздержусь от обсуждения этой темы, где наличной сейчас информации,
возможно, недостаточно, если кому эта история интересна, то вот дополнительные ссылки: письмо Лузина Канторовичу   
\cite{ReshKut}, книжка Серпинского по теории множеств
\cite{Serp}, статьи В.~И.~Богачёва и А.~В.~Колесникова \cite{Bogachev-Kolesnikov}
и В.~М.~Тихомирова \cite{Tihomirov-A}.

Впрочем, дискуссия о Суслине зацепила также вопрос о лузинской манере цитирования. Приведем отрывок из предисловия Лебега
\cite{Lebeg} к книге Лузина:
\begin{quotation}
Всякий, вероятно, удивится, когда узнает, читая Лузина, что я, между
прочим, изобрел метод решета  и первым построил аналитическое множество.
 Никто, однако, не удивится так, как я. Г-н Лузин лишь тогда бывает
совершенно счастлив, когда ему удается приписать собственные открытия
кому-либо другому. Странная причуда; она мне кажется простительной,
поскольку нет опасности, что г-н Лузин создаст школу в этой области.
\end{quotation}

Замеченный Лебегом стиль
 касался не только цитирования Лебега:
\begin{quotation}
	АЛЕКСАНДРОВ.
	В течение всех выступлений H.~H., я думаю, это можно будет по стенограмме установить, H.~H. в течение этого заседания стоял на той точке зрения, что А-множества открыты Суслиным.\dots
	
	Там я прочел введение от начала до конца... В этом введении совершенно ясно, что вся теория этих множеств принадлежит Лебегу, т.е. содержится в своем зародыше у Лебега, имени Суслина
	здесь вообще совершенно нет.
	
	Теперь дальше идет первая глава, которая имеет заглавие: «Конструкция Лебега и ее
	обобщение».
	
	Второй параграф — [«Каноническое решето Лебега»].
	
	Третий параграф — «Аналитические множества Лебега».
	
	Имя Лебега встречается по несколько раз на каждой странице этой первой главы.
	
	Вторая глава «Измеримость в смысле B [Бореля]».
	
	Первый параграф — «Идеи Бореля».
	
	\dots
	
	Потом множества, измеряемые в смысле В—Бореля.
	Дальше идет целый ряд параграфов. Среди них вообще встречается в заглавиях параграфов только имена Лебега — 5—6 раз, имя Бореля — 2 раза. Значит, H.~H., если цитировать, если смотреть на это количество цитат Лебега и Бореля, можно составить впечатление, насколько Николай Николаевич тщателен в указаниях авторства даже в тех случаях,
	когда это не считается особенно обязательным.
\end{quotation}

Лузин в самом деле приписывает свои идеи и идеи своих учеников Лебегу, Борелю,
Адамару (и это потом вводило в заблуждение некоторых исследователей,
оригинальность книги Лузина  маскируется этим стилем):
\begin{quotation}
АЛЕКСАНДРОВ.	В
	тексте, не в примечаниях, а в тексте, имя Суслина, если я  не ошибаюсь, совершенно не
	фигурирует.
\end{quotation}

Это в самом деле так, в основном тексте книги <<Лекции по аналитическим множествам>> имени Суслина (который 
открыл эти самые <<аналитические множества>>) нет. Ссылки на него есть в примечаниях, которые были собраны в конце
книги.
Не упоминается в основном тексте и Александров (на него есть одна ссылка в примечаниях).
А Новиков, Лаврентьев, Келдыш и Селивановский в основном тексте  упоминаются.

Опять некрасиво, но формально авторы и работы цитировались, и дотошный читатель
мог до ссылок добраться.

В целом,  благодаря Новикову и Лаврентьеву (они должны были делать выбор, они его
сделали), успех нападения по этому пункту был весьма относительным.
Да, скользкое поведение, но без полного <<криминала>>.

%Но в оправдание Лузина могу сказать: у него было феодальное мышление: «я феодал, это мои %вассалы, и все, что они делают, — моё, но и я вассал вон того феодала, и все моё — его».

\sm

{\bf\punct Лузин, 1930 год.}  Напомним, что в 1927 году Лузин был избран членом-корреспондентом АН СССР,
в 1928 был вице-президентом Международного математического конгресса, в  1929 избран действительным членом АН СССР.
В 1928-29 году во Франции работает над монографией, которая подытоживала 
работы его самого и учеников по дескриптивной теории множеств,
она выходит в 1930. Обычно
люди после написания оригинальных книг некоторое время приходят в себя.
Мы видим, что в том же 1930 году выходит учебник Грэнвиль--Лузин.

Осенью 1930 года Лузин возвращается в Москву и сталкивается с новыми проблемами.
Конфликт с учениками. Кризис его собственной научной программы и невозможность
быстро предложить что-то новое. Если в 1922 году он был математическим богом,
то теперь он мало кому в Москве интересен. На московском Физмате творится черт
знает что. Бешеная политизация, нашествие неведомых людей, введение бредовых форм обучения. 
Атаки на чистую математику. Арест Егорова.
В конце 1930 года Лузин уходит с Физмата и переходит на работу в ЦАГИ (в теоретический отдел Чаплыгина).

В это  время он начинает публиковать прикладные работы.

\sm

{\bf\punct Обвинения в публикации халтуры.%
\label{ss:haltura}}
Выше цитировалось высказывание <<Правды>> о <<методе академика Крылова>>
 и статье Ф.~Р.~Гантмахера. Приведем несколько высказываний на этот счет.
 Бари, 1950:
\begin{quotation}
	Но сам Н.~Н. искал области 
	приложений не в развитии тех методов и идей, которые выросли из его 
	исследований и из его преподавания, а во включении в работу по 
	задачам ему совершенно чуждым. Так он занимался в эти годы 
	теорией так называемого векового уравнения, темой, на которую ему,
	несомненно, указал академик А.Н.Крылов.
\end{quotation}

Из Стенограммы:
\begin{quotation}
	СОБОЛЕВ. Вот относительно Гантмахера. Он писал в редакцию журнала, что работа
	Гантмахера ученическая, и ее не стоит печатать. Дело было так. Гантмахер посылает
	свою статью в «Математический сборник». Лузин пишет по поводу этой статьи два письма. Одно письмо Гантмахеру: «Работа Ваша замечательная, нельзя в таком дрянном
	журнале как «Математический сборник» ее печатать». Он не употребил этого слова как
	«дрянной», но что-то в этом роде. Другое письмо в редакцию «Сборника»: «Работа Гантмахера носит совершенно ученический характер. В «Математическом сборнике» помещать ее нельзя, я согласен ее поместить в каком-нибудь другом месте». Два письма по
	поводу одной и той же работы.
\end{quotation}

Представляется невероятным, что Соболев мог бы тут блефовать. История была недавней,
и документы были, скорее всего, целы. Продолжение дискуссии:
\begin{quotation}
	КРЫЛОВ. Нужно сказать, что о работе Гантмахера один отзыв дан самому Гантмахеру, а другой — в редакцию «Сборника».
	
	КРЖИЖАНОВСКИЙ. Здесь мы упоминаем вот что. Здесь говорится о своеобразном
	«великодушии».
	
	БЕРНШТЕЙН. В основном это верно.
\end{quotation}

Понтрягин высказывался на эту тему на собрании в МГУ.
Дословный текст  нам не известен, но
оно было отреферировано неизвестным автором в журнале <<Фронт науки и техники>>:
\begin{quotation}
	\dots
	Вот, например, молодой математик Гантмахер сделал на трех страницах то, что Лузин
	написал на 160. Но достаточно Лузину воспротивиться, и редакция университетского
	«Математического сборника» не помещает эту замечательную работу. Не удается также
	провести ее и в «Известиях» Академии наук. Правда, Гантмахеру дается там «удовлетворение» — отмечается,
	 что он проделал такую же работу, что акад. Лузин, но при этом
	упускается отметить, что это сделано Гантмахером на трех страницах, вместо 160.
	Почему? \dots
\end{quotation}

О чем шла речь. Крылов был судостроителем, военным инженером, 
ма\-те\-ма\-ти\-ком-при\-клад\-ни\-ком и интересовался методами вычислений. Он предложил  способ
 вычисления характеристического (точнее, минимального)
многочлена для матрицы $A$  размера $n$, который, как объяснил
Гантмахер, сводится к следующему.
Возьмем вектор $v$, пусть, для простоты, он циклический%
\footnote{Напомним, что у оператора есть циклический вектор, если
каждое его корневое подпространство является жордановой клеткой,
множество операторов  без циклического вектора имеет коразмерность 3 в пространстве
всех операторов. Если же циклический вектор есть, то множество всех
нециклических векторов - объединение конечного набора  подпространств коразмерности 1.}. Берем цепочку векторов
$v$, $Av$, \dots, $A^{n-1}v$. Они образуют базис. Раскладываем
вектор $A^{n} v$ по этому базису:
$$
A^n v=\sum_{j=0}^{n-1} \alpha_j A^j v.
$$
(для этого надо решить систему из $n$ линейных уравнений).
Записывая оператор $A$ в новом базисе, мы получаем
$$
\begin{pmatrix}
0&1&0&\dots&0\\
0&0&1&\dots&0\\
\vdots&\vdots&\ddots&\ddots&\vdots\\
0&0&0&\dots&1\\
\alpha_0&\alpha_1&\alpha_2&\dots&\alpha_{n-1}
\end{pmatrix}
$$
Понятно, что числа $\alpha_j$ -- это и есть коэффициенты характеристического многочлена%
\footnote{Крылов это говорил иными словами. 
	По матрице он писал систему линейных  дифференциальных уравнений
	первого порядка,
	а дальше приводил ее к одному дифференциальному уравнению.}
\footnote{Выписанная матрица - каноническая форма для оператора, имеющего циклический вектор;
	если циклического вектора нет то можно взять прямую сумму матриц такого вида.
В связи с этим я слышал слова <<вторая жорданова нормальная форма>>,
но никогда не встречал таких слов в печати. Кто и когда ее заметил, автор не знает.
Возможно, что это в самом деле изобретение Крылова и Гантмахера.}.
Подробнее см. в <<Теории матриц>> Гантмахера \cite{Gant}. В итоге для получения
характеристического многочлена мы должны решить одну систему
линейных уравнений порядка $n$. Если мы вычисляем напрямую,
то мы должны вычислить определитель порядка $n$, у которого диагональные 
элементы зависят от дополнительной переменной (то есть, нам придется проводить преобразовывать
матрицами, чьи элементы являются многочленами). Прием Крылова полезен
и в практических вычислениях используется.

В 1931г. Лузин написал на эту тему 4 статьи общим объемом 160 страниц,
\cite{Luzin-Krylov-1},  \cite{Luzin-Krylov-2} (кстати, объем тогдашних
<<Известий АН СССР>> был небольшим и издавался он по всем
естественным наукам вместе; для сравнения, общий объем <<Мат.Сборника>>
в 1930 и 1931гг. был 212 и 240стр.). Надо думать,
что цель Лузина состояла в поддержании
отношений с Крыловым (кстати в названии было
<<метод академика Крылова>>, в заголовках математических статей
не принято писать чины при фамилиях)%
\footnote{ЛЮСТЕРНИК. На нашем заседании [непонятно, на каком] профессор Гол[убев] ...
	 цикл его работ о работах А.~Н.~Крылова ... назвал халтурой.}.

\begin{quotation}
	СЕГАЛ. Все, близко знающие H.~H. и главные его работы, подчеркивают, что он большой мастер придавать изложению сжатость и изящество.
	В отзыве, подписанном Новиковым и Ляпуновым, подчеркнуто, что такое-то доказательство значительно упрощено
	Лузиным и что то, что другой автор написал бы на 20-30 страницах, Лузин пишет на трех
	страницах. Если мы будем рассматривать советские работы Лузина, то нигде вы не найдете такого многословия, как у него. Здесь наоборот, то, что можно было бы сказать на
	трех страницах, он говорит на 20-ти страницах; он излагает исследовательские работы
	так, как будто бы они предназначаются для людей начинающих. Это обстоятельство также
	имеет место.
	
	ХИНЧИН. К тому, что сейчас было сказано, я хочу добавить маленький факт. Не далее
	как вчера по просьбе редакции журнала «Геофизика», я взял на себя труд просмотреть
	статью Лузина, очень большую, в несколько десятков страниц — 50-60 страниц — это для
	математики очень большая статья, — предназначенную для помещения в этом журнале
	по вопросу об анализе периодограмм. Она в высшей степени подтверждает то, о чем
	сейчас было сказано. Статья написана прекрасно, очень красивым языком, с большой
	ясностью, но, что меня глубоко поразило, что человек так размазывает. Столько воды
	напустил, что я не мог себе отдать отчета, для чего это сделано.
	%
	%Хинчин. Мне пришлось давать официальный отзыв на последнюю работу Лузина.
	%Там есть один достойный внимания результат, но этот результат мог быть доказан на 4-
	%5 страницах, а статья занимает 77 страниц.
\end{quotation}

Вообще-то почти всем математикам  <<удается>> время от времени писать плохие статьи, это как-то само собой получается, без всякого на то желания авторов. Обвинение
в единичной халтурной работе вряд ли может считаться серьезным, оценка все же должна вестись по хорошим работам и по совокупности работ. Лузин не был автором, склонным к халтуре, и, кстати,  многие его работы <<второго уровня>>  достойны руки мастера, отметим биографию Эйлера
\cite{Luzin-euler},
  добавление \cite{Luzin-add-lebeg}
  к книге Лебега и предисловие \cite{Luzin-Zhegalkin}
  (впрочем, с предисловиями у Лузина не все было благополучно).
  
  С <<методом академика Крылова>> 
   Лузин подставился, а потом усугубил это блокированием
  публикации Гантмахера, выйдя за рамки приличий.
  
  Что касается <<периодограмм>>, то, по-видимому, речь идет об опубликованной
  посмертно работе Лузина
  \cite{Luzin-periodogram}. Пусть дана функция, которая представляет из себя
  сумму гармоник,
  $$
 f(t)= \sum_{j=1}^n a_j \sin\bigl(\omega_j t+ \beta_j\bigr),
  $$
частоты, амплитуды и  число $n$ (вообще говоря) неизвестны. Требуется восстановить
неизвестные параметры%
\footnote{Пример, где возникает подобный вопрос -- описание движения Луны, задача, на протяжении тысячелетий  доставлявшая
наблюдателям неба много мучений (и они продолжались, по крайней мере, весь XIX век). Другой вопрос (связанный с движениями Луны, но все же
отдельный), предсказание высот приливов в данной точке морского побережья.
Лузин, скорее всего, имел в виду попытки предсказания погоды,
занимался он этим по предложению одного из лузитан С.~С.~Ковнера, работавшего
зам. директора Института геофизики.}. Кстати, задача эта не столь уж далека от основной лузинской тематики
(возможно, на это намекает Бари в отрывке цитированном  чуть выше).
 Статья (70 страниц в Собрании сочинений Лузина) представляет из себя обзор, рассчитанный на  нематематиков.
Хинчин отчасти прав, но чересчур суров%
\footnote{Для сравнения:
\newline
Лузин-- Крылову, {\bf 22 декабря 1932г.} (в оригинале указан 1923г., очевидно неверный,
но письмо довольно точно датируется, исходя из его содержания).
\newline
{\it
... прежде всего о кандидатурах в члены-корреспонденты. Я продолжаю утверждать,
что кандидатура А.~Я.~Хинчина пока преждевременна. {\bf А.~Я.~Хинчин не работает уже два года,
не печатает ничего, и по сведениям, даже перестал делать доклады}. Неизвестно,
когда он возвратится к Науке и вернется ли он к ней вообще.
\newline
Я считаю, что звание члена-корреспондента есть важное научное отличие. А.~Я.~Хинчин
слишком молод, чтобы давать ему за его прошлое. Значит это будет аванс.
\newline
Но печальный опыт  Academie de Paris с  Painlev\'e
(я отнюдь не сравниваю научных сил Painlev\'e и А.~Я.~Хинчина,
совершенно несравнимых), избранного в Institut de France
и переставшего работать в самый момент избрания
и потом всю жизнь занимавшегося чуждыми науке делами
-- этот опыт, по моему, не стоит повторять с А.~Я.~Хинчиным...}
\newline
Однако, и занимаясь академическими интригами,
 Лузин сохраняет высокохудожественный
  стиль. Отвлекаясь от вопроса о значении тогдашних работ Хинчина
  (находившегося в расцвете своего таланта и прокладывавшего
  новые пути в теории вероятностей), отметим, что аргумент
  с Пенлеве относится к области поэзии, а не логики,
  упомянем также число научных публикаций Хинчина по годам:
  1929 --- 8, 1930 --- 1, 1931 --- 0, {\bf 1932} --- 7 (из них 5
  в легко доступном Мат. Сборнике
  - в Мат. Сборнике), 1933 --- 6 статей и монография Asymptotische Gesetze der Wahrscheinlichkeitsrechnung. 
  <<Сведения>> об отсутствии докладов тоже, по меньшей мере, преувеличены,
  во всяком случае доклад Хинчина 21.05.1932 
  состоялся на первом заседании восстановленного после смуты
  Московского математического общества (вряд ли Лузин мог этого не знать).}. 

\sm

{\bf\punct Печатанье лучших работ за границей.%
\label{za-granitsej}}
Обвинение из анонимной статьи в <<Правде>> 
в следующей публикации  9 июля  было усилено и распространены на других
математиков:
\begin{quotation}
	До сих пор считается почему-то нормальным, естественным в научной среде печатать труды советских ученых, прежде всего,
	за границей или даже только за границей.
	Взять хотя бы математику. Большинство советских ученых-математиков (Александров, Колмогоров, Хинчин, Бернштейн и др.) публикуют свои работы за границей, не
	печатая их 
	у нас, в СССР, на русском языке.
	\dots
	
	Дошло до того, что даже популярные работы (по топологии, теории вероятностей)
	профессоров Александрова, Хинчина, Колмогорова впервые напечатаны были за границей 
	на немецком языке, а затем только был «поднят вопрос» о переводе этих работ 
	советских ученых на русский язык и переиздании в СССР. Чувство глубокого возмущения
	охватывает, когда знакомишься с такого рода фактами.\dots	
\end{quotation}

Вопрос этот не имел того смысла, который сейчас принято ему придавать.
Был 1936 год, до Холодной войны и Железного занавеса было еще 10 лет.
Разъяснения члена ЦК ВКП(б) Кржижановского (уже после второй статьи в Правде, о которой речь
чуть ниже) были такими:
	\begin{quotation}
КРЖИЖАНОВСКИЙ. Ведь здесь два момента: первый о том, что по преимуществу печатают за границей, второй момент — в СССР печатается то, над чем человек внутренне
		смеется. Это очень серьезная вещь. Представьте психологически это из Вашей собственной картины. Человек представляется в лучшем случае безмерно испугавшимся обывателем, который чувствует власть и при каждом проявлении власти только трепещет. Из
		информации видно, что чтобы гарантировать себя в этом отношении, он что делает? Он
		думает: надо сунуть в рот пищу этому чудовищу, на — пишу, работаю. {\bf И главное обвинение не в том, что он печатает работы за границей, а в том, что он печатает в СССР такие работы, какие он сам называет в разговорах с друзьями явной «белибердой>>}\dots
	\end{quotation}
	
	\begin{quotation}
	КРЖИЖАНОВСКИЙ.   {\bf  Здесь никто не возражает против того, чтобы значительные работы
		печатались за границей; здесь разбирается конкретный случай совершенно определенного отбора. Здесь товарищи подчеркивают совершенно ясно: мастер стиля, мастер этого отбора во всех русских его работах, напечатанных в СССР, показывает как раз обратные качества.} Что это случайно или это делалось сознательно? Из всей характеристики,
		которая здесь перед нами прошла, для нас совершенно ясно, что это делалось сознательно.
	\end{quotation}  

Очевидно, что обсуждаемое обвинение  не могло исходить от Александрова
(он сам  публиковал почти все статьи  за границей, за границей
же была опубликована  знаменитая
книга Alexandroff, Hopf <<Топология>>. 
Но, возможно,
это обвинение было сформулировано в математическом  мире
(мы обсудим этот вопрос ниже).
Понятно, что многим работающим математикам оно должно было не понравиться.
Процитируем заседание Комиссии 7 июля до выхода второй анонимной статьи в <<Правде>>

\begin{quotation}
	СОБОЛЕВ. Теперь вопрос об изданиях. Я согласен с С.~Н. [Бернштейном], что
	обвинение в том, что он печатал за границей лучшие статьи — это обвинение, может
	быть, необоснованное. Да если бы это и было так, то это может быть объяснено тем, что
	все-таки заграничные журналы печатают несравненно скорее, несравненно лучше. Так
	что в этом смысле могло иметь место совершенно естественное стремление напечатать
	там, где выйдет скорее, поскольку вопрос о выходе в свет той или иной работы, имеет
	смысл. Это обвинение я считаю недостаточно серьезным.
	
	\dots\dots
	
	ХИНЧИН. Что же касается печатания его работ, то я должен признать, что факты, изложенные в этой статье, совершенно правильные, уместные, но они не дают материала для
	обвинения. Совершенно верно, что Н.~Н.~Лузин свои лучшие работы печатает за границей, но я думаю, что 90\% советских математиков поступают таким же образом.
	
	ФЕРСМАН. Но одновременно печатают и здесь.
	
	ХИНЧИН. Да, одновременно печатают и здесь, но именно лучшие работы печатаются
	за границей. Ну, может 90\% — это и преувеличение, однако значительная доля советских 
	математиков так поступают, и не из политических соображений, а из желания иметь
	оттиск в хорошей обложке, на хорошей бумаге и иметь поскорее, иметь без опечаток. У
	нас, к сожалению, этого нет.
	
	АЛЕКСАНДРОВ. В частности, я должен сказать, что процент иностранных работ по отношению к советским у меня гораздо больше. Я только в нынешнем году стал печатать
	мои работы здесь, потому что до сих пор у нас в совершеннейшем развале было это
	дело.
	
	СОБОЛЕВ. И притом, в то время как наше издательство дает 25 оттисков, любое иностранное издательство дает 100 оттисков.
	
	АЛЕКСАНДРОВ. {\bf  Но я должен указать, что одну из самых своих крупных работ H.~H.
	выпустил в «Математическом сборнике»} [имеется в виду первая большая работа по проективным множествам \cite{Luzin-sbornik-1926} ]. Я думаю, что этот пункт вообще не принадлежит к числу сильных пунктов обвинения.
	
	ХИНЧИН. Тем более что {\bf здесь еще один момент является для нас весьма серьезным.
	Что для нас политически более правильно: печатать наши работы здесь или за границей? Где должен быть центр тяжести, это еще спорный вопрос.}
\end{quotation}

Заметим, Александров -- главный обвинитель -- возражает (не единственный раз)
против справедливости данного высказывания Правды в отношении Лузина, а Хинчин открыто сомневается в самой идее, высказанной <<Правдой>>. 

Вторая статья в <<Правде>> ставила Александрова, Хинчина и Колмогорова
 в положение обороняющихся (а Бернштейн исчез и появился уже на последнем заседании комиссии).
  На собрании 9 июля в МГУ Александров и Колмогоров
 признавали свои ошибки и обещали в дальнейшем вести себя лучше \cite{Front1936-1}.
 
 В Стенограмме Александров говорит следующее:
 \begin{quotation}
 АЛЕКСАНДРОВ. \dots
 Я резюмирую: едва ли особенно серьезные упреки можно сделать персонально Лузину. Здесь мы все приблизительно одинаково виноваты. 
 \end{quotation}

 Однако  Хинчин и Соболев продолжали возражать против этого обвинения
 (Хинчин цитируется ниже, в п. \ref{ss:hinchin1936}).
 
 У части нападавших была другая точка зрения на этот вопрос:
\begin{quotation}
СЕГАЛ. Я подробно ознакомился с работами H.~H., начиная с 1930 года, как выпущенными на русском языке, напечатанными в советских журналах, так и за границей. И у
меня создалось совершенно определенное впечатление, что тут был сознательный плановый отбор. Ни одной серьезной работы на русском языке и в советских журналах не
было напечатано, кроме одной работы, помещенной в номере пятом «Трудов физико-математического института».	
	\end{quotation}
	
	\begin{quotation}
ШМИДТ. Самый факт, что подавляющую часть своих работ Лузин печатал за границей, а в Союзе печатал только явную чепуху, это бесспорный факт и, если математическая общественность затрудняется квалифицировать этот факт [нужным] образом, то
		объясняется это только тем, что в этом повинны и другие. Сегодняшняя статья в «Правде» сигнализирует о том, что это действительно явление общее. Это верно, но можно,
		однако сказать, что это не умаляет вины Лузина, потому что в этом стиле нашем советском, фактически антисоветском — не печатать у нас — повинен тот же Лузин, как руководитель одной из школ, как наибольший западник, как популярная фигура и т.д. Все следовали его образцу. Так ведь дело обстоит. Нужно сказать, что раньше до деятельности
		Лузина, такового взгляда у математиков российских, пожалуй, не было. Крупнейшие
		российские математики печатали преимущественно в России.
	\end{quotation}
	
	\begin{quotation}
	ЛЮСТЕРНИК. Я должен сказать, что до революции H.~H. не ориентировался вовсе на
		заграницу. Его книга «Интеграл и тригонометрический ряд», которая составила основную славу Лузину, была напечатана в России в 1915 году. Значит у H.~H. до революции,
		наоборот, лучшие работы печатались здесь.
	\end{quotation}
	
Теперь обсудим, до какой степени сентенции из <<Правды>>
соответствуют истине (не в смысле <<что такое хорошо, а что такое плохо>>, а в смысле истинности).
 Расцвет Лузина начался в 1911--1915 гг, а к 1930
году он оказался в научном кризисе (мы это подробно обсуждали в \S \ref{s:fate}). С 1933гг. года положение с изданием математических
журналов в СССР начало быстро улучшаться (см. ниже п. \ref{ss:journals}), а до того  было весьма неважным. Так что большая часть лучших работ Лузина после 1917 года в самом деле была напечатана за границей. 

Признавать научный кризис публично Лузину по понятным причинам не хотелось, 
но даже его сторонники это понимали:
\begin{quotation}
БЕРНШТЕЙН.	\dots
для
оправдания вообще несколько пониженного в научном отношении характера его нового
творчества, \dots
\end{quotation}	

На пятом по счету заседании Комиссии 13 июля уже вполне затравленный Лузин говорит:
\begin{quotation}
ЛУЗИН. \dots Тов. Александров против этого сказал, что мои лучшие работы восходят до 20
	года. Совершенно правильно. Дальше в моих теоретических работах — я не говорю о
	проективных множествах, потому что это довольно сильная вещь — но дальше в моих
	чисто теоретических работах начинаются сильные колебания, и я не считал эти чисто
	теоретические работы настолько сильными, что они только и достойны заграницы.
\end{quotation}

Но до того он выдвинул объяснение, что под давлением Райкова, Кольмана и Кагана
он стал печатать в СССР прикладные работы, а за границей -- теоретические.
Приведем три длинных цитаты: 

\begin{quotation}
	ЛУЗИН.
	Когда тов.~{\bf Райков}, который указывал на многие вещи, с которыми в
	научном отношении я не был согласен, на то, например, что теоретические работы не
	нужны и т.д., — все это на меня действовало. Я ведь плохо разбирался в обстановке.
	\end{quotation}
	
	\begin{quotation}
		ЛУЗИН.
		Всем известны, ни от кого
		этого не скроешь, статьи тов.~{\bf Кольмана}, который говорил, что мои теоретические вещи
		пропитаны насквозь идеализмом, что все это не что иное, как вредная ерунда, что, собственно говоря, в наше время должна быть совершенно другая установка. И вот все это
		мне было хорошо известно, я это читал своими глазами. А отрицательное отношение
		тов.~Кольмана, который занимает высокое положение, и с мнением которого я обязан
		был считаться, вовсе не как с человеком, который отрицательно относится ко мне, но с
		человеком, которому партия доверяет, его мнение было для меня чрезвычайно веским,
		— это отрицательное отношение заставило меня пересмотреть всю мою деятельность,
		после того как на меня эти нападки были напечатаны
		
		\dots
		
Моя цель была в том, чтобы не повредить моей
родине. И если те лица, которым видно, как все это движется, находят мою теоретическую деятельность вредной, пусть она не вредит. Но из меня перли эти теоремы. Из меня
шли эти теоремы. Что же мне было делать? Что же — в себя засовывать, что ли? И я 
решил, что все, чем я могу нашей стране служить в прикладном направлении, печатать у
нас. А в том направлении, которое признается вредным в моей практической деятельности (ведь в нескольких газетах было напечатано мнение тов. Кольмана, — и с этой стороны был употреблен прямо совершенно определенный отрицательный резкий отзыв о
моих теоретических работах в связи с математическими ошибками, которые были там
констатированы — не в моих сочинениях, а в отзывах тов. Кольмана, который в математическом отношении ошибается, потому что он в некоторых тонких понятиях теории
функций путает, и это все мои ученики знают, и все смеются над той грубой путаницей, которая может проскользнуть по нечаянности) — все это печатать за границей. Ведь
я с этим должен был считаться.

\dots
Но я решил перенести принципиально все то, что имеет приложение, или может иметь приложение, все это печатать в нашей стране, а все, что относится
к чисто теоретическим вопросам, к очень отдаленным трансфинитным вопросам, все
это перенести для заграницы, и то просто потому только, что я не мог удержать в себе
творчество. Вот почему у меня целый ряд статей есть о методе Крылова.
\end{quotation}

\begin{quotation}
	ЛУЗИН...
%Когда я послал статью из Парижа и сам приехал, то узнал, что эта статья
%[<<Дифференциальное исчисление>> для Большой советской энциклопедии] 
%вызвала жестокие нарекания, и тогда Вениамин Федорович [{\bf Каган}] заявил, что эту статью печатать нельзя и т.д. Так
%продолжалось около двух-трех лет. После этого ко мне обратились с просьбой эту статью развить. Когда я спросил, зачем это,
%раз она вызвала такое противодействие, то
%мне сказали, что я должен развить ее в том же направлении. Это мною было сделано,
%и я ее развил%
%\footnote{Это отголосок отдельной небезынтересной истории, которая обсуждается ниже. 
%Главред Энциклопедии Шмидт тогда был атакован Яновской  за неправильное издание энциклопедии
%и  несоответствие ее математических статей диалектическому материализму, а также за
% отсутствие
%в оных классового анализа,  Каган, который до  1931г. курировал естественно-научный раздел, был
%обвинен Кольманом (в том же самом). Был также предъявлен  античный образец
%того, как надо делать правильно \cite{Kolman-Yanovskaya}. Кагану было некуда деться.
%Для дальнейшего важно, что по прошествии некоторого времени
%статья Лузина  \cite{Luzin-encyclopedia-diff} была расширена <<в том же направлении>>. Кроме того, Лузин в конце 1930
%(чуть раньше истoрии с Энциклопедией) оказался на кафедре, которой заведовал Каган.}.
Я должен сказать следующее. В связи с этой статьей [имеется в виду статья Дифференциальное исчисление
для Большой Советской Энциклопедии, см. ниже п.\ref{ss:kolman-kult1}] я много думал по поводу печатания своих научных работ, и мне пришла идея:
так как установка делается на строительство, на индустриализацию, то, быть может, было полезнее для нашей страны печатать
всякие работы, которые имеют близкое или отдаленное отношение к прикладной науке.
Имея в виду историю с этой статьей для энциклопедии, я боялся повредить нашей стране,
так как эта статья была признана вредной, я боялся, так как не разбираюсь во всех тонкостях,
будучи человеком книжным и плохо понимающим политические моменты. 
Поэтому я просто решил печатать всё то, что имеет чересчур отдаленное значение, за границей,
а что имеет сколько-нибудь прикладное значение — у нас.	
\end{quotation}

Что в сказанном правда -- что в 1930 году Лузин должен был сильно испугаться,
и можно уверенно сказать, что испугался (речь об этом годе еще пойдет ниже \S \ref{s:annals}).
Однако разные другие выводы, которые можно с несомненностью сделать
из сказанного, будут далеки от истины. Атака на чистую математику в самом деле
была в  1929-1930г., но она была отбита (отбита с участием лузинских обвинителей -- чистых математиков, а также Шмидта) не позднее
начала 32 года (а, скорее всего, раньше). После этого никто не против теоретической математики не возражал. 

Другая неправда --  Лузин по-русски опубликовал ряд работ по своей основной тематике.
Это  две  тонких книжки \cite{Luzin-teor-funkts}, 1933, \cite{Luzin-descript},
1935,
кроме этого было три добавления к русскому переводу книги Лебега  \cite{Lebeg-1934},
1934.
Книги рассчитаны на более широкого читателя, чем статьи (а Сегал в цитате выше лукавит,
правильным произнесением фразы исключая книги из рассмотрения). Было еще
3 статьи по чистой математике в русских журналах.% и одна в трудах Всесоюзного съезда????.

Что правда -- что все прикладные работы писались по-русски, их на тот момент было всего 6, включая 4
статьи по <<методу академика Крылова>>. Так что если не считать этого злополучного <<метода>>,  
инкриминировать  Лузину в смысле отбора худших работ для советских
изданий было  нечего.

Но присутствовавшие на Комиссии прекрасно знали и про кратковременность атаки на чистую математику
(благо, что многие участвовали в ее отражении), и про лузинские книжки.
Своими сомнительными оправданиями Лузин  усугублял свое положение и навешивал на себя новые обвинения.   

Кржижановский произнес на это жесткую отповедь, но нам будет интересно привести ее
ниже в п.\ref{ss:krzhizh} по другому поводу. А здесь процитируем Александрова:

\begin{quotation}
	АЛЕКСАНДРОВ. Действительно, из того, что говорил Райков, было много глупостей.
	Много упреков, совершенно необоснованных, Райков делал не только вам, но и многим
	другим математикам. В той части, в которой мы считали эти упреки со стороны Райкова
	несправедливыми, мы против них спорили, иногда даже довольно резко спорили. И,
	между прочим, вы отлично знаете, что мы все продолжали заниматься вещами нисколько не менее абстрактными, чем ваши. Область, которой я занимался, приложима к практике не более чем ваша. И на заседании Института я отстаивал эту точку зрения, что моя
	область не имеет непосредственно приложения на практике, точно так же как и теория
	чисел их не имеет...
	\end{quotation}

\begin{quotation}
АЛЕКСАНДРОВ. \dots Но то, чем мы занимались особенно тогда, было чисто теоретической областью, и нам приходилось спорить, потому что тогда некоторые работники нашего Института действительно высказывали ту точку зрения: то, что немедленно неприложимо на практике — это все белиберда. Но с этой точкой зрения никто не согласился. Были бурные заседания, и мы отстаивали противоположную точку зрения, которая была оправдана. {\bf Но вы, H. H., проявили такое свойство вашей природы — я буду говорить совершенно откровенно — которое во
многих случаях мешало в жизни, проявили малодушие, граничащее с трусостью. Какой
бы человек, носящий партийный билет} [видимо, имеется в виду Райков], {\bf чтобы вам ни сказал, вашей первой реакцией
сейчас же был страх, как бы чего не вышло. Сопротивления у вас не было.}
\end{quotation}

В этом месте Стенограммы (и не только в этом, см. слова Колмогорова, цитированные в следующем пункте)
 звучит обвинение совсем другого
свойства. В сложной обстановке 1930-1932гг. Лузин  не оказался соратником тех,
кто оказывал сопротивление накату на математику, а потом восстанавливал
научную жизнь в Москве...

{\bf\punct Отзывы.%
\label{ss:otzyvy}} Спустя 30 лет Люстерник рассказывал об этом в таких
словах \cite{Lyu-3}:
\begin{quotation}
 Некто X написал работу с опровержением геометрии
Лобачевского. Его работа была послана на отзыв Н.~Н.~Лузину с просьбой
ответить на вопрос, стоит ли Х-у предоставить академический паек.
Н.~Н.~Лузин ответил: «Из работы Х-а видно, что он незнаком с основными
трудами в этой области. Поэтому я считаю, что ему следует предоставить
академический паек, чтобы он получил возможность ознакомиться с ними».
Паек Х-у был предоставлен. Позже подобные иронические отзывы были
одним из оснований для выдвинутых в прессе обвинений против
Н.~Н.~Лузина.
\end{quotation}

В 1936 году он говорил иначе:

\begin{quotation}
	ЛЮСТЕРНИК. Я редко с вами встречался, почти не встречался, и мог не говорить об
	отзывах. Я говорил много с Михаилом Алексеевичем [Лаврентьевым] и с Ниной Карловной [Бари], и они мне говорили, что они неоднократно с Вами беседовали. Эти отзывы
	и вас компрометируют, и вообще недопустимы. Вы сами говорили, что вы не знали, что
	коллектив отрицательно относится к вашим отзывам.
\end{quotation}

В многочисленных положительных отзывах на низкокачественные научные работы, книги и сомнительные заявки на должности
(applications)
Лузина хором обвиняли все присутствующие. Говорили, что таких отзывов было
вообще много (обратите внимание на цитату из Юшкевича,
с которой начинается этот параграф), и называли конкретные имена.
Стоит иметь в виду высокое научное положение Лузина в 1929-1936гг.
В 1929 он оказывается одним из пяти академиков по математике.
С 1934г. он -- председатель комиссии по присуждению
ученых степеней и званий по математике.
Вскоре он становится и  Председателем математической группы Академии.

Во многом из того, что обсуждалось, сейчас уже трудно чего-то понять
(участники знали много больше, чем мы сейчас).
Но мы предпримем сейчас такой разбор (он набран мелким шрифтом) и покажем,
что при определенных преувеличениях (в которых возможно отразились также контры нападавших
между собой и сведение счетов с отсутствовавшими) эти обвинения, в 
значительной части были справедливы.  

\sm

{\small 
Приведем сначала пару отрывков из дискуссии:
\begin{quotation}
АЛЕКСАНДРОВ. Когда я приехал в Ленинград, то буквально каждый математик рассказывал о том, что человек напечатал полнейший вздор, но явился с письмом от вас в Дом
ученых и требовал постановки доклада.

ЛУЗИН. Он приходил ко мне три раза, буквально осаждал.

АЛЕКСАНДРОВ. В Ленинграде все говорили, что человек совершенно безграмотный
приходит с письмом от вас и требует, чтобы в Доме ученых был поставлен доклад на
тему действительного бреда, и математики с большим трудом добились того, чтобы этого не было.

ЛЮСТЕРНИК. Они, к сожалению, не смогли добиться этого, и только вследствие путаницы в рассылке повесток Дома ученых доклад сорвался.

ЛУЗИН. Мое письмо можно найти. Я просил относительно отдыха, но относительно
того, чтобы ставить доклад, я не просил и не мог просить, потому что я не знал его как
лектора, а даже верные идеи — плохой лектор, излагая их, может сорваться.
\end{quotation}

Стоит попытаться представить себе, какое впечатление должны были производить подобные 
оправдания на присутствующих.%
\footnote{
Вот еще образец того, как Лузин оправдывался
(в Стенограмме чуть ниже):
\newline
{\it 
С МЕСТА. О первой работе Льва Генриховича [Шнирельмана] вы выражались, что это топорная работа. Эта работа с трудом может быть названа топорной.
\newline
ЛУЗИН. Это было бы безумием, а что элементарная работа и не доделана до конца,
Лев Генрихович и сам это признает. Он дал 1500 или 2000 этих самых простых... Это колоссальное открытие. Я готов подписаться и выступить в печати, что это гордость и украшение нашей страны.
\newline
СЕГАЛ. Я всего раз виделся с вами до этого, и вы именно так отзывались о работе
Льва Генриховича.
\newline
ЛУЗИН. Я признаю топорными работы Чебышева. Они слишком элементарны. Я
Льва Генриховича очень уважаю и ценю, и ожидаю от него очень многого.
\newline
ХИНЧИН. Вы не помните точного выражения?
\newline
ЛУЗИН. Да и зачем это? Я не первый день относительно этой работы говорю. Я думаю,
что Лев Генрихович вполне заслуживает того, что он член-корреспондент Академии наук.
А другие его работы, дальнейшие, носят гораздо большую глубину. Я здесь совершенно
не виноват. Топорной работой я признаю в этом смысле и работу Чебышева.}
\newline
В принципе Лузин имел полное личное право считать работу Шнирельмана топорной... Но беседа дикая.}

В обвинениях, по-видимому, присутствовали преувеличения, и об этом говорил
Александров, мы процитируем ответ Лузина на высказывания Александрова:
\begin{quotation}
	ЛУЗИН (13 июля).
	Относительно Романова. Я очень рад могущественной поддержке, которую мне
	оказал Павел Сергеевич, что это лицо — не идиот, что он талантливый человек, но своей
	рекомендацией я мог бы ему повредить. В этом отношении я приношу Павлу Сергеевичу
	глубочайшую благодарность, потому что я думал, что я мог в отношении его ошибиться.
	
	Относительно Безсонова Вы также меня поддержали.
	
	Относительно Эйгеса [Вы] меня также поддержали тем, что специалист по геометрии
	[В. Ф. Каган] в нем когда-то принимал участие.
\end{quotation}

Фамилия <<Романов>> упоминалась в <<Правде>> 16 июля, как будто высказываний
 Александрова о нем в
Стенограмме, предшествовавших словам Лузина, не видно
\footnote{Романов Николай Павлович (1907--1972).  О нем есть статья в <<Успехах>>, 12:3(75) (1957). Был в Москве в аспирантуре в 1929-32г. у Хинчина и Шмидта. В 30х продолжал теоретико-числовую
	работу Шнирельмана. В 1935г, работая в Томске, защитил докторскую диссертацию.
	Позже был заведующим кафедрой в Ташкенте. Скорее всего его (или его докторскую) не одобрял кто-то из (трех)  участвовавших в нападении теоретико-числовиков.}.

%\begin{quotation}
	%{\sc Александров.} Вы никогда не считали его математиком, и тут были даже шутливые
	%высказывания. Но ни отзыв о Безсонове, но ни отзыв об Эйгесе я не считаю большими
	%вашими прегрешениями, потому что Безсонова я считаю довольно хорошим математиком,
	%и я не взялся бы говорить ни за, ни против, так что я считаю, 
	%что я никогда не позволю себе в данном случае упрекать вас в сознательном вредоносном %действии, а в случае
%	с Харламовой это было действительно возмутительно.
%\end{quotation}

Разобраться  в деталях прочих обвинений, не зная ничего (или очень
мало) о действующих лицах проблематично. Попробуем это сделать.
\begin{quotation}
АЛЕКСАНДРОВ. \dots
	Н.~Н. писал предисловия ко многим книгам в этом же духе —
	смешные, смехотворные.\dots
\end{quotation}

Что это за <<предисловия ко многим книгам>>, автору данных записок не известно.
Я видел предисловие к учебнику Мат. анализа Жегалкина и Слудской
(весьма интересный текст, который, кстати, в 1985г. был перепечатан в <<Успехах>>),
к переводу математического справочника Дубельта для инженеров, и к изданиям Грэнвиля.
Очевидно, речь шла о чем-то другом.
Впрочем, дважды повторенная сентенция Соболева бросает на это некоторый свет: 
\begin{quotation}
СОБОЛЕВ. Я просил бы Вас, H.~H., со всей откровенностью ответить на такой вопрос:
	согласились ли бы вы написать предисловие к книжке ферматиста Тер-Микаэляна, если
	бы вы эту книжку должны были издавать в Париже? Согласились ли бы вы написать то
	же самое предисловие с тем, чтобы на французском языке вышла эта теорема Ферма?
	Хватило бы у вас храбрости для этого?
\end{quotation}

Книга эта  в итоге не вышла.
Л.~А.~Тер-Микаэлян (1867-1943) - крупный инженер, строитель железных дорог.
Он был не первый и не последний почтенный человек, ударившийся в ферматизм.
Почему Лузин не предоставил почетное право разбираться с ферматизмом теоретико-числовикам,
непонятно (вообще-то по каждой долго нерешаемой математической проблеме есть стандартный набор ошибок,
которые известны специалистам, совершить новую ошибку, особенно в таком деле, как доказательство теоремы Ферма, дело непростое).
Предисловие Лузина (оно было издано в 1993 году \cite{Luz-mikaelyan}) интересно как философско-математический текст,
но ведь это предисловие крупнейшего математика к 
книге, по меньшей мере, подозрительной...
\begin{quotation}
	Но главное, что следует читателю принять во внимание, - это то, что дело
	 вовсе или не в том,  чтобы даваемое доказательство Великой теоремы было верным или неверным.
	 Все дело  в {\it методе}, которым оно проводится. Самое доказательство Великой теоремы может оказаться вполне верным, и,
	 тем не менее, может не представлять ни малейшей ценности.
	 Ему тогда место в специальном журнале, где оно, в случае если совершенно не соответствует характеру XVII века, 
	 ни содержит какого-нибудь арифметического принципа, будет сохраняться как некая музейная ценность.
	 Наоборот, та или иная попытка рассмотрения Великой теоремы может быть не полной, но если она хорошо
	 дает почувствовать самый характер методов, которыми оперировал divus arithmeticus - то это уже явится большой ценностью.
	 
	 В этом - значение и смысл предлагаемой работы автора.	
\end{quotation}	

Вот так. Ни много, ни мало.

\sm

Теперь просто об отзывах и рекомендациях.
\begin{quotation}
	ЛУЗИН. \dots
	Я думаю, что среди моих отзывов было более 50\% таких, которые действительно стояли на абсолютной высоте.
	
	Затем надо считаться с индивидуальностью. Я не машина, я живой человек и у меня
	характер чересчур мягкий. Я не хочу говорить о слабости воли, — что дано, то дано, —
	но я человек очень мягкий. И все движения, которые вызваны были по тому или другому
	поводу, когда мне приходилось давать отзывы, это были движения сердца. Возьмите несчастного Шадхана.
	
	Шадхан написал мне письмо: он, мальчик, участвовал в ужасающих
	вещах, был живым свидетелем погрома, под который попала его семья. Он говорил о
	нервном состоянии его семьи, говорил о том, как пробивал себе путь, и передо мной
	совершенно живо встала эта картина несчастного юноши, который был исключительно
	безумно и горячо предан науке, в этом отношении он хотел пробиться. У меня рука не
	могла отказать ему.
	
	\dots
	
ЛЮСТЕРНИК. Я хочу сделать несколько добавлений к вопросу о Шадхане. Шадхан с
	отзывом H.~H. пошел к Тумаркину. Тумаркин попросил меня переговорить с ним, и я
	убедился, что он не знает даже, что такое производная, не имеет представления ни о чем
	и даже задач по элементарной математике, задач на логарифмические уравнения он не
	может решать, а вы его рекомендуете в Академию.
	На счет возраста — так ему сейчас 22 года, это не ребенок, ну 21 год, во всяком случае, он не ребенок. Он дифференцировать не умеет.
	
ЛУЗИН. Но он читает в Педагогическом ВУЗ'е лекции по математике и за чтение этих
	лекций получает премии и благодарности.
	
ЛЮСТЕРНИК. Возмутительно, что он читает. Надо было протестовать. Это наша задача.
	
ЛУЗИН. Кадров нет.

ЛЮСТЕРНИК. Кадров нет, это другой вопрос. Но на этом основании рекомендовать —
	это возмутительно.
	
ЛУЗИН. Ведь я его не рекомендовал ни на степень доктора, ни на степень кандидата...

СМЕХ И ГОЛОСА. Еще бы этого не доставало!

ЛУЗИН. ...ни на степень аспиранта при Академии наук. Я написал письмо, в котором
	просил посмотреть и указал относительно того, что ему переходить, конечно, в Академию наук в качестве аспиранта невозможно, но я просил, чтобы Отдел кадров посмотрел, нельзя ли что-нибудь сделать для такой бесконечной жажды знаний, которая у него
	была.
	
АЛЕКСАНДРОВ. Скажите, каким образом Шадхан мог быть свидетелем погромов, когда
	ему всего 20 лет?
\end{quotation}

По-моему, обсуждавшиеся в стенограмме случаи <<Ласаев>>, <<Харламова>> тоже являются доказанными, хотя в случае
<<Харламова>> Лузин мог опасаться лиц, стоявших за аспиранткой.

Хинчин, как и Люстерник говорит о многочисленных случаях.
\begin{quotation}
	ХИНЧИН.
	Дальше, здесь говорится об отзывах. Должен сказать, что к списку тех, которые здесь
	даются, я мог бы прибавить очень много известных мне. Мне приходилось по должности моей сталкиваться с этим,
	так как ко мне как к директору Института математики приходили люди, требовавшие внимания к своим исследованиям и приносившие с собою
	отзывы и характеристики Лузина. Я знаю отзывы значительно более резкие. Ко мне приходили люди явно психически больные, 
	работы которых представляют сплошной бред
	сумасшедших, по поводу которых, однако, H.~H. в своих отзывах говорил, что в них имеются мысли, которые,
	вне всякого сомнения, будут положены в основу будущей математики, в основу будущей физики и т.д.
	\end{quotation}}
	
Одно из немногих выступлений Колмогорова:	
\begin{quotation}
	КОЛМОГОРОВ. Нам говорили, что отзывы ваши противоречат нашим, и, пожалуй, это
	был единственный случай, когда авторитет H.~H. мог принести чрезвычайно существенный вред. 
	Было такое время, когда мы проводили аспирантуру с большим трудом, и его
	авторитет мог принести существенный вред.
	\end{quotation}

 <<Провести аспирантуру>> относится к ликвидации последствий событий 30-31гг., 
 о которых пойдет речь ниже, и получить здесь подножку от человека, в котором естественно
 было предполагать соратника, надо думать, было сильно обидно.

\sm

{\bf\punct Замечания.%
\label{ss:zamechaniya}} Участники дискуссии предъявили Лузину много иных научно-этических
обвинений, которые мы не пытаемся реферировать%
\footnote{СОБОЛЕВ.  {\it Относительно политики, которую проводил H.~H., в частности, когда дело
	касалось выборов членов-корреспондентов АН СССР и т.д., я помню, например, выборы
	1934-го года, когда Н. Н. вел странную совершенно политику. Именно, С.~Н.
	 [Бернштейн] представил довольно большой список серьезных ученых с тем, чтобы
	отнестись со всей серьезностью к такому ответственному делу и чтобы поговорить о
	том, кто действительно достоин, кто может быть избран. В этом списке была выставлена
	талантливейшая молодежь, как А.~Колмогоров, Гельфонд и т.д. Список этот был составлен, с моей точки зрения, вполне объективно. И требовалось только в группе подойти
	серьезно к оценке каждого кандидата. И вот Н.~Н. по формальным соображениям —
	когда чего-то не хватало, не был представлен своевременно какой-то отзыв, который
	можно было представить на следующий день, — словом из-за каких-то мелочей, к которым он придрался, в совершенно истерическом тоне потребовал, чтобы ни одна кандидатура, кроме той, которую он выставил, не обсуждалась. Фактически он сорвал всякое
	обсуждение выборов и поставил группу в такое положение, когда был выставлен только
	один кандидат [Голубев]. Я не буду говорить о том, что этот кандидат достоин, может быть, это и
	так, но самый факт этого отказа в обсуждении, этот срыв выборов достойных кандидатов
	говорит сам за себя. Я лично видел в этом деле то, что H.~H. было просто, может быть,
	глубоко противно, что случайно может быть окажется выбранным какой-нибудь представитель из той молодежи вроде А.~Н.~Колмогорова, которого он в Академию допустить не
	желал. Я после группы говорил Н.~Н., что это безобразие, H.~H., то, что вы делаете. А он
	мне ответил, что это священные традиции Академии наук и т.д.
\newline
Я считаю, что в Академии он проводил политику, которая, во всяком случае, шла во
вред Академии наук.}}. Очень много говорил Александров. Мы не имеем доступных источников для проверки многих свидетельских показаний 1936.
%Например, мы не можем знать, кто был прав или не прав в перепалке Александрова и Лузина,
%об их разговоре за полгода до июля.
 Часть обвинений произносившихся на комиссии, как мы видели, были натянуты или преувеличены.
Очень вероятно, что нападавшие имели против Лузина еще что-то,
что они говорить не хотели (Александров точно имел, надо думать,
имел и Гельфонд).
 Но стоит иметь в виду, что часть участников нападения (прежде всего
Александров и Хинчин) активно возражали против обвинений, которые считали недостоверными. 

Многие обвинения не вызывали сомнений для участников дискуссии,
которые имели заведомо больше информации, чем мы. 
При тех данных, которую мы сейчас имеем, можно утверждать, что часть обвинений
предъявленных на комиссии была доказана (в значительно меньшей степени это замечание
может быть отнесено к статьям в <<Правде>>).

\sm

Автор данных записок много лет наблюдал научный мир изнутри и далек от мысли,
что Лузин превосходил своими отрицательными качествами всех академиков
(хотя в 1935-36гг. он имел б\'ольшую власть, чем поздние академики). Своими 
выдающимися достоинствами он превосходил подавляющее большинство из них. Но то, что
атака на Лузина была вызвана свойствами его характера и его собственными действиями,
по-моему, не вызывает сомнений. Выпишем список учеников,
с которыми он явно или неявно поссорился: Суслин, Александров, Урысон, Хинчин, Колмогоров,
Лаврентьев, Люстерник, Шнирельман, Новиков.
Часть из них здравствовала на памяти ныне живущих. Мы можем вспомнить о них не только
хорошее, но все же помним их как людей достойных. Можно вспомнить и конфликты Лузина с его учителями, 
Млодзеевским и Егоровым, людьми тоже уважаемыми, и, судя по всему, порядочными;
был у него и конфликт с А.~В.~Костицыным. Что-то многовато.

\sm

{\bf\punct Рассказ Гнеденко.%
\label{ss:gnedenko}} Свидетельств непосредственных очевидцев об облаве на Лузина сохранилось немного.
Рассказы Ефремовича и Понтрягина уже были задействованы выше (а слова С.~М.~Никольского были буквально процитированы
в начале параграфа, они ничего не добавляют к исследуемой картине).

Гнеденко в воспоминаниях о Колмогорове писал:
\begin{quotation}
	Тяжелый осадок оставила в моей душе антилузинская кампания в газетах, журналах и даже на митингах.
	Некоторые из его учеников клеймили его, отказывались от него как от учителя, забывая, 
	что именно он был создателем знаменитой Московской математической школы. 
{\bf 	Для меня была невыносима эта травля, мне казалось, что те, кто позволял себе позорить учителя, 
	совершали бесчестный поступок.} Для меня учитель был и остается вторым родителем, 
	чье имя должно быть поставлено рядом с именами матери и отца.
	Исключительно достойно себя вели тогда Н.К.Бари и Д.Е.Меньшов.
	В ту пору я еще не знал Лузина, не слышал его лекций, не видел его ни разу.
	Я знал только, что он один из создателей замечательной научной школы московских математиков.
	\end{quotation}
	
Стоит все же отметить, что именно два учителя Гнеденко (Хинчин и Колмогоров) нападали на своего
общего учителя.	
В статье 1990г. о Степанове Гнеденко говорил чуть определеннее:
\begin{quotation}
  Мне было стыдно,
и я не был на срочно созванных общественных судилищах. Очень немногие нашли в себе
силы выступить в защиту когда-то горячо любимого учителя; среди них нужно отметить
Н.~К~Бари и Д.~Е.~Меньшова.
\end{quotation}
Кстати, слова <<среди них>> предполагает, что Гнеденко имел в виду кого-то еще%
\footnote{Кстати, из известных лузитан в этой истории (документах и поздних рассказах) не упоминается Гливенко.}...

О поведении Меньшова и Бари известно мало. Юшкевич пишет так
\begin{quotation}
	Против постановления, по словам присутствовавшего на нем профессора
	А.~Т.~Григорьяна, в то время еще студента, голосовали только два самых преданных Н. Н. Лузину 
	его ученика — профессора Н.~К.~Бари и Д.~Е.~Меньшов.
\end{quotation}
Впрочем, Юшкевич мог быть очевидцем этого собрания и сам.

\sm

{\bf\punct Из статьи Бари.%
\label{ss:bari}}
\begin{quotation}
	Тем не менее, та брешь, которая образовалась в отношениях 
	между Н.~Н. и группой московских математиков, из которой часть 
	были его учениками, не только не исчезла, по еще более  
	углубилась. Сейчас очень трудно установить причины этого расхождения. 
	Ясно только одно: многое в этих отношениях основывалось на  
	недопонимании. 
	
	Несомненно, что часть его прежних учеников и сотрудников 
	раздражала известная самоизоляция, в которую поставил себя Н.~Н. 
	после возвращения из-за границы. С другой стороны, возможно, 
	что самому Н.~Н. казалось, что его ученики не всегда оценивают ту 
	большую роль, которую в их творчестве играло его идейное  
	руководство. Это взаимное недовольство и расхождение прорвалось  
	летом 1936 года. 
	
	Началось с того, что в «Правде» появилось письмо какого-то 
	директора средней школы, в котором подвергалось чрезвычайно 
	резкой критике отношение Н.~Н. к работе школы. Суть письма  
	состояла в том, что школа пригласила Н.~Н. присутствовать на  
	экзаменах, надеясь получить от него критику работы и указание на пути 
	дальнейшего улучшения работы. Вместо этого Н.~Н. все огулом 
	хвалил, неумеренно и неискренне всем восхищался, проявил,  
	одним словом, полное безразличие к итогам работы под видом  
	добродушия и доброжелательства. А так как были иногда весьма  
	критические замечания со стороны Н.~Н. по работе средней школы,  
	делавшиеся в интимных научных кругах, то все поведение Н.~Н.  
	рассматривалось как двурушничество.
	
	Как потом оказалось, это было началом широкой кампании 
	против Н.~Н. Вслед за этим письмом появилась редакционная 
	статья в «Правде», где вся деятельность Н.~Н. рассматривалась, 
	как проявление антисоветских настроений, преклонения перед  
	заграницей, двурушничества. Н.~Н. обвинялся в присвоении чужих 
	научных достижений, например, его учеников П.~С.~Новикова и 
	ряда других, в затирании и травле талантливых учеников, 
	например, М.~Я.~Суслина, что косвенно привело даже к его  
	смерти, в подхалимстве, в протаскивании весьма слабых работ 
	раболепствовавших перед ним бездарных математиков и в других 
	смертных грехах. Статьи в «Правде» нашли свое продолжение в 
	собраниях профессоров и преподавателей в Университете, в АН 
	СССР. 
	
	Эти собрания носили самый безобразный характер. Выступали 
	ученики Н.~Н., выступала математическая молодежь, учившаяся у 
	Н.~Н., продолжавшая его работы, Л.~С.~Понтрягин, Л.~Г.~Шнирельман, А.~А.~Ляпунов и другие, причем их выступления часто  
	состояли из площадной брани, позорной для стен Университета и  
	Академии наук. Выработался даже термин «лузинщина». Это была  
	хорошо подготовленная, организованная травля. В известном  
	смысле резюме всего содержала редакционная статья в «Успехах  
	математических наук», выпуск III, 1937 г.; редакция состояла в то  
	время из А.~Ф.~Берманта, Ф.~Р.~Гантмахера, А.~Н.~Колмогорова, Н.~Е.~Кочина, В.~Д.~Купрадзе, Л.~А.~Люстерника, И.~Г.~Петровского, 
	А.~И.~Плеснера, Б.~И.~Сегала, В.~И.~Смирнова, В.~А.~Тартаковского и 
	Д.~А.~Райкова (тех. ред.). 
	
	Таким образом, две стороны деятельности Н.~Н. подверглись 
	особенно резким нападкам. Во-первых, отмечалось, что Н.~Н. не 
	всегда аккуратно и точно отмечал достижения своих учеников, так 
	что иногда получалось такое впечатление, что полученные ими  
	результаты выдавались Н.~Н. за результаты, полученные им самим. 
	Во-вторых, не менее резкому осуждению подвергалось  
	преклонение Н.~Н. перед иностранной наукой, его научный  
	«космополитизм». Известная доля правды была и в том, и в другом.  
	Несомненно, что при тех формах коллективной работы, которая велась в 
	научных семинарах Н.~Н., часто весьма трудно было до конца  
	разобрать, кому впервые принадлежала та или иная математическая 
	идея. Так, например, было с происхождением созданной М.~Я.~Суслиным,
	П.~С.~Александровым и Н.~Н.~Лузиным теории «A-множеств», которые впоследствии, во время пребывания за границей, 
	были Н.~Н. систематически разработаны в его большой работе 
	<<Lecons sur les ensembles analytiques>>. В обработке этой  
	теории Н.~Н. во время своего пребывания за границей, естественно, 
	находился в условиях исключительно благоприятных. Однако во 
	всех этих обвинениях необходимо помнить, какую ведущую роль 
	во всех исследованиях, проведенных учениками Н.~Н., играли его 
	собственные идеи и его прямые указания. Среди выступавших  
	против Н.~Н. многие просто не знали всех обстоятельств дела и их 
	суждения носили иногда особенно авторитетный вид просто по 
	полному непониманию дела. 
	
	Точно так же и «космополитизм» Н.~Н. надо понимать очень 
	условно. Несомненно, что в обращении Н.~Н. было кое-что от 
	манер «преувеличенной галантности», не всегда искренней, и это 
	естественно, раздражало и коробило иногда, в особенности, когда 
	приходилось наблюдать его обращение с иностранными учеными 
	Однако не надо забывать, что Н.~Н. и его школа заняли видное 
	место в иностранной математической науке, конечно, не формой 
	обращения, а своими крупнейшими научными результатами. В  
	течение двух или двух с половиной десятилетий работы Н.~Н. и его 
	учеников непрерывно служили украшением мировой научной  
	литературы по теории функций, а относительно польской математики 
	едва ли будет преувеличением сказать, что она в 20-х и 30-х годах 
	выросла и развивалась на идеях московской школы Н.~Н.~Лузина. 
	
	Едва ли можно сейчас во всех подробностях выяснить все  
	движущие пружины той жестокой критики, {\bf в размерах своих явно 
	несправедливой}, которой подвергался Н.~Н. систематически в  
	течение летних месяцев 1936 года. \dots	
\end{quotation}

	Довольно странен список <<Понтрягин, Шнирельман, Ляпунов>>.
	Присутствие Шнирельмана вполне понятно (см. следующий параграф).
	Неупоминание Александрова  может допускать объяснения не научно-политического характера
	(Бари, скорее всего, знала обстоятельства связанных с Лузиным конфликтов; видно,
	что она не считала Лузина во всем правым; к тому же, как мы
	отчасти видели в этом параграфе и увидим в следующем,
	поведение Александрова на Комиссии было не столь уж однозначно).
	
	\sm
	
	{\it Ляпунов}.
	В опубликованных документах 1936 года Ляпунов фигурирует лишь
	в совместной бумаге о плагиате, поданной в Комиссию вместе с Новиковым,
	а также в нескольких репликах в Стенограмме, которые ничего не значат.
	Однако напомним, что Ляпунов не участвовал в издании сочинений Лузина в  1951-1959
	(диссертация \cite{Luz-trig}, парижская монография \cite{Luz-anal-rus}, и собственно <<Собрание сочинений>>
	\cite{Luz-collected-1}--\cite{Luz-collected-3})
	и в серии мемориальных статей в <<Успехах>> 1952-53гг
	\cite{BaLyu1}, \cite{KeldNov}, \cite{GolKuz}, \cite{L-Sret}, \cite{Fedor}%
	\footnote{Напомним, что участники книжного издательского проекта - Бари, Голубев, Л.~Келдыш,
	Лаврентьев, Люстерник, Меньшов,  Некрасов, Новиков, Сретенский, В.~С.~Федоров, кроме них в <<Успехах>>
	участвовали  также В.~К.~Гольцман и П.~И.~Кузнецов. Из гонителей Лузина присутствовал лишь Люстерник,
	который, надо думать, раскаялся, и раскаяние его, надо думать, было принято. Впрочем, непосредственно
в комиссию по изданию сочинений Лузина и он не входил.}.
	Это косвенно подтверждает слова Бари.
	Была также  рецензия Ляпунова \cite{Lyap-retsenziya} на первый том собрания сочинений 
	Лузина, содержавший работы по теории функций. Из 
	рецензии видно  недовольство Ляпунова; основная суть претензий - что сопроводительные статьи, биографию
	и библиографию было решено сосредоточить в третьем томе%
	\footnote{Стоит еще отметить, что диссертация Лузина была
		к этому моменту
	образцово издана в \cite{Luz-trig} вместе с сопроводительной статьей,
	большими  содержательными комментариями, 
	биографией Лузина, библиографией, а также статьями Лузина,
	примыкающими к его диссертации (частью переведенными с немецкого и французского). Очевидно, что это издание потребовало огромного
	труда
со стороны Бари, Меньшова, Голубева и Люстерника. При издании 
тома <<Собрания сочинений>>, посвященного ТФДП, перед издателями должен был встать вопрос о том, чтобы избежать
повторения той же книжки. Поэтому рецензия Ляпунова не выглядит  джентельменской.}. Кроме того, известны 3 
	рукописи Ляпунова  о Лузине разных лет, которые не были  в своё время опубликованы. Исходя из упоминаемых
	и неупоминаемых
	в рукописях крупных научных событий две из них датируются 1950-1960гг и серединой
	60х годов, третья рукопись содержит дату, 1971. Естественно думать, что сотоварищи по каким-то
	причинам не хотели иметь с ним дела по этому поводу.
	
	\sm
	
	{\it Понтрягин.}
	То, что мы знаем о выступлении Понтрягина от него самого
	и из \cite{Front1936-1}, не дает ответа на вопрос, почему он возглавил список.
	Но мы знаем также, что склонностью к резким высказываниям он всегда обладал.

\sm

{\bf\punct Снова о <<Собрании сочинений>> Лузина.%
\label{ss:SS}}
Выше было замечено, что там нет CV.  Кроме того:

1) Предисловия Лебега нет ни в  отдельном издании <<Лекций по аналитическим множествам>>, 
ни во втором томе <<Собрания сочинения>> (куда эта книга вошла).

2) Примечания в конце <<Лекций по аналитическим множествам>> в обоих изданиях книги превращены в подстрочные
(тем самым фамилии Суслина и Александрова становятся видимыми непосредственно
при чтении книги).

3) В списке научных трудов Лузина (по-видимому, полном) не указаны страницы.

4) Список публикаций  Лузина не совсем полон.  Разумеется, там нет статьи <<Приятное разочарование>>.
Но там нет и предисловий (что с одной стороны естественно,
а другой стороны и не очень -- некоторые из них являются весьма интересными сочинениями). Скорее всего,
нам и сейчас не все публикации Лузина видимы.

5) Во втором томе, посвященном дескриптивной теории множеств,
 отсутствуют публикации 
Sur un raisonnement nouveau dans la th\'eorie des fonction descriptive
и Sur un choix d'ensemble parfait distingu\'e dans un compl\'ementaire analytique 
arbitraire ay ant des constituantes non d\'enombrables
в С. R. Acad. Sci, Paris, 201.

6) В третьем томе (куда попали прикладные
 работы) нет статей по методу академика Крылова.

Понятно, что все это, так или иначе, отражает свару 1936 года.

 \section{За гранью%
 \label{s:gran}}
 
 \COUNTERS
 
 \epigraph
 	{В <<деле>> Лузина бросается в глаза несоответствие между {\it мелочностью} предъявленных ему конкретных обвинений и {\it резкостью}
 		 их политических квалификаций.}{Бирюков~Б.~В. \cite{Biryukov-yanovsk}}

 Напомним, что статья в Правде от 3 июля  завершалась потоком ругательств
 в адрес старой московской математической школы и ее питомца Лузина.
 Эти ругательства были существенно дополнены и обогащены 
  в ходе обсуждения на Академической комиссии.
 
 \sm
 
  {\bf\punct Подача Шнирельмана.%
  \label{ss:shnirelman-podacha}}
  \begin{quotation}
  	{\sc ГОРБУНОВ.} [9 июля 1936]
  	Во-вторых, у меня есть заявление члена-корреспондента тов. Шнирельмана относительно одного случая, который дает известный политический штришок и относится к 1930
  	году. Когда после процесса Промпартии группа московских ученых по инициативе тов.
  	Люстерника решила обратиться к французским ученым с письмом протеста против
  	угрозы интервенции, то Лузин, отговариваясь болезнью, не дал своей подписи под документом, хотя его подпись была очень важна ввиду его связей с французскими учеными.
  \end{quotation}
  
  Шесть лет спустя в 1936 году Шнирельман и Яновская вспомнили об этом случае.
%  Из выступления Яновской на собрании в МГУ в тот же день 9 июля.
%  
%  \begin{quotation}
%  	В 1930 г. H. H. Лузин председательствовал на том собрании ученых, которое приняло
%  	обращение к французским ученым — протест против интервентов — в связи с делом
%  	Промпартии. Но Лузин уклонился от того, чтобы собственноручно подписываться под
%  	этим воззванием. Напрасно тогда молодой аспирант Рабинович стучался в двери Лузина. %Узнав, что он пришел из института математики за подписью, Лузин заявил, что болен, что ни %принять, ни подписать обращение не может.
%  \end{quotation} 
%  
  В следующие дни этот эпизод склонялся на комиссии так и сяк, основным свидетелем был Люстерник.   Мы вернемся к этому  ниже при обсуждении событий вокруг Лузина в 1930 году, см. п.\ref{ss:prompartiya}.
  Отголосок дискуссии был запечатлен в статье в 
   <<Правде>> от 14 июля: 
  \begin{quotation}
  	В 1930 г., во время процесса «Промпартии», физики и математики Москвы подписали обращение к ученым за рубежом с призывом поднять протест против интервенционистских намерений по отношению к СССР некоторых иностранных держав. На этом
  	обращении не было подписи академика Лузина. Его позорное недвусмысленное поведение в то время, его жалкие и нечестные увертки сейчас ссылками на болезнь, политическую неграмотность, незнание политической обстановки и прочее, разоблачают профессора Люстерник, Хинчин, Соболев и другие.
  \end{quotation}

  \sm
  
  {\bf \punct Подача Гельфонда.%
  \label{ss:gelfond-podacha}}
  \begin{quotation}
  	{\sc ГЕЛЬФОНД.} В 1928 году случился такой эпизод, который очень всех удивил. Был международный съезд в Болонье, и среди математиков узнали, что H. H. принял предложение Серпинского войти в Президиум или Оргкомитет Всеславянского математического
  	съезда. Это вызвало чрезвычайное удивление, потому что всякому ребенку ясно было,
  	что такое панславянство, что оно было орудием самодержавия, ясно было, что если
  	Польша устраивает что-то общеславянское, то не надо было быть очень политически
  	квалифицированным человеком, чтобы понять, что это такое. Это вызвало возмущение
  	многих, например, С.~Н.~Бернштейна и т.д. В конце концов, H.~H. послал телеграмму в
  	Президиум Академии наук — можно ли ему принимать участие во Всеславянском математическом съезде?

  	 ГЕЛЬФОНД. \dots Но совсем недавно мне попались
  	Труды этого съезда, и там указано, что члены Оргкомитета — Н.~Н.~Лузин и H.~M.~Крылов%
  	\footnote{Крылов, Николай Митрофанович (1879—1955), академик АН СССР.}
  	 (эти Труды имеются в библиотеке Московского университета). 
  	Там указывается,
  	что, к сожалению, в силу некоторых обстоятельств Н.~Н.~Лузин и H.~M.~Крылов не смогли принять участия в работах этого съезда. 
  	Зато там выступали представители эмигрантов-математиков. Эти труды можно достать и посмотреть. Отказались ли они официально или нет, 
  	неизвестно, во всяком случае, там сказано, что они являются членами
  	Оргкомитета.
  \end{quotation}
  
  Если отвлечься от предельной политизированности самого Гельфонда, то
  <<Всеславянский съезд>>, разумеется, был политизированным мероприятием, и принимать
  участие в нем не следовало. Можно еще добавить, что Польша была резко враждебна
    очень слабому в 1928 году СССР
  и стремившилась к гегемонии в Восточной Европе.  Но Лузин все же в съезде  не участвовал. 
  
  \begin{quotation}
  	 ХИНЧИН. По поводу этого Всеславянского съезда. Кстати, Вы его неправильно именуете. Как будто бы этот съезд назывался
  	Съездом математиков славянских стран%
  	\footnote{Эта реплика -- смягчающая.}.
  	
  ХИНЧИН. Я просто для протокола говорю, потому что нужно правильно называть документы.
  	Вы считаете себя виноватым в том, что после того как вас включили в Президиум
  	съезда, вы не протестовали и не потребовали вашего снятия. Но все-таки включили Вас в
  	этот съезд с Вашего согласия или без Вашего согласия? Было бы очень странным, и в
  	практике съездов это не принято. Или организационный комитет без вашего согласия
  	включил? Если бы с Вашего согласия, то было бы странно, конечно, после этого просить
  	о снятии. А если согласия не было, то странно, что они вас включили.
  	\end{quotation}

  	\begin{quotation}
  	ЛЮСТЕРНИК. Меня удивляет все-таки такая мотивировка. Вы должны были написать
  	письмо Борелю [это в связи с делом Промпартии], а потом решили, что «так как он так поступил, я не написал письмо». То
  	же самое и здесь, такая аргументация, что Вы решили не иметь ничего общего и поэтому
  	допустили, чтобы в издании трудов этого Съезда была поставлена Ваша фамилия и что
  	там было выражено соболезнование, что Вы не смогли туда приехать.

  	\dots 
  	
  	 ЛЮСТЕРНИК. {\bf  Кого Вы боялись, например? Советской власти? Я не понимаю этой трусости, кого Вы боялись?}
  	
  	 ЛУЗИН. Я не хотел.
  	
  	ЛЮСТЕРНИК. Неужели было кем-нибудь плохо понято, если бы вы заявили отказ поехать туда? Больше того, {\bf вы могли бы в «Правде» напечатать о Вашем отказе}.
  	
   	ЛУЗИН. Вы правы, Лазарь Аронович. Я только не догадался.
  \end{quotation}
  
  Еще один эпизод обсуждения:
  \begin{quotation}
  АРШОН. У Н.~П.~Горбунова имеются два факта, которые более или менее характерны.
  Один тут был оглашен — это приглашение Лузина в оргкомитет Всеславянского математического съезда. Это относится к 1927 г.
  
  Второй, не менее важный факт: когда здесь был Серпинский [топологическая конференция, сентябрь 1935], махровый черносотенец, то этот Серпинский предлагал Новикову напечатать свои работы, причем, по словам
  того же самого Новикова, он сказал ему, что было бы желательно, чтобы эти работы шли
  за двойной подписью — Новикова и Лузина. Это как раз те работы, о которых мы говорим, что это плагиат. Тут вот что — совершенно невероятно предполагать, что махровый
  черносотенец желает создать авторитет советскому ученому. Поэтому этот факт, с моей
  точки зрения, является очень характерным.
    \end{quotation}
  
  На непросвещенный взгляд автора данных записок
  Серпинский так же не мог быть черносотенцем, как и матерью-одиночкой.
  Но передовым людям всегда видней.
  
  \sm
 
 Закончим обсуждение цитатой  
  из газеты <<Правда>> 14 июля 1936 года.
  \begin{quotation}
  	Реакционный съезд за границей избрал
  	его членом президиума, и Лузин дал на это свое согласие.
  	Теперь в комиссии Лузин рассказывает, что впоследствии его охватил «ужас». Он
  	понял, что участвует в антисоветской манифестации. Его спросили: «Почему же вы не
  	вычеркнули своего имени?» Лузин ответил, что он «боялся», но не пояснил, кого.
  	%Антисоветский характер в 1930 г. имел уход Лузина из университета, когда были
  \end{quotation}
  
  \sm
  
  {\bf\punct Подача неизвестного.%
  \label{ss:podacha-neizvestnogo}}
  \begin{quotation}
  	С МЕСТА. Характерно, что он во время Кассо, когда ушли даже кадеты из Университета%
  	\footnote{Речь идет о массовом уходе преподавателей МГУ в знак 
  	протеста против циркуляров министра просвещения Кассо  о полицейском надзоре
  над университетами (1911).}, не ушел — был приват-доцентом и остался приват-доцентом
  	\end{quotation}
  	
  	О Кассо вспомнил Шнирельман незадолго до этого (обвиняя Егорова).
  	<<Доброжелатель>> с места проэкстраполировал обличение на Лузина, который в 1911
  	был приват-доцентом и находился в заграничной командировке
  	(1910-1914) от Московского университета. Дальше началось уже отдающее полным бредом обсуждение этого предмета (с участием, в основном, Горбунова и Кржижановского). Такая мелочь, что Лузин, не уйдя из МГУ, в итоге создал
  	Московскую математическую  школу, в расчет браться не могла. Лишь в последний 
  	день, когда снова появился Бернштейн и пришел Крылов, прозвучали возражения:
  	\begin{quotation}
  КРЫЛОВ. Но ведь он уезжал за границу при Кассо. Это 25 лет тому назад было. Нельзя
  	ему ставить в вину, что он не примкнул тогда, а остался в той группе, в которой раньше
  	был.
  	
  	 БЕРНШТЕЙН. Почему не упоминают о том, что констатировано на том заседании, на
  	котором я был, — что H. H. в 1920 году во время разрухи проявил колоссальный энтузиазм? Почему не упоминается, что с его стороны не было никаких попыток, и никто не
  	может сказать, что он выступал когда-либо против советской власти. Вы должны понимать, что если даже Александров говорит, что в него верили как в бога и были с ним
  	близки, то никто, однако, не констатировал, что у него были в то время черносотенные
  	настроения. В период, когда была уже советская власть, никто не может сказать ничего,
  	в чем проявилось бы это его отношение к советской власти.
  	\end{quotation}

 {\bf \punct Самоподача Лузина.}
 \begin{quotation}
 ЛУЗИН.	Величайшей моей ошибкой, не упоминаемой в статье, но глубоко мною пережитой,
 	был и мой уход из Университета в 1930 г. Сейчас мне трудно даже отдать себе отчет, как
 	могло случиться, что я, крупный ученый, видевший свое призвание в преподавании
 	математики, смог отказаться от преподавания в Университете и не понять, что развитие
 	нашей жизни и быстрый рост нашей страны открывает передо мною перспективы, о
 	которых я никогда не мог мечтать.
 	Сознавая вполне всю глубину политической моей ошибки в неумении оценивать
 	обстановку и вовремя переключаться, я согласен нести за нее любую ответственность.
 \end{quotation}

И это долго обсасывалось на комиссии, и лишь в последний день Лузин получил поддержку.
 \begin{quotation}
 	КРЫЛОВ. Это не о том. В 1930 г. он просто ушел из Университета, а ошибка его в том,
 	что он начал изъясняться, почему ушел. Он должен был бы сказать — захотел и ушел.
 \end{quotation}
 
 А Бернштейн находит точные допустимые формулировки и в связи с уходом
 Лузина из Университета, и его занятиями прикладной математикой. 
\begin{quotation}
	БЕРНШТЕЙН.
	Но вот когда он ушел из Университета, то, во-первых, он тогда в Университете не был
	ни в какой мере заинтересован, во-вторых, с учениками он был в явно плохих отношениях и работать в Университете для него не представлялось уже интересным. В-третьих, для
	него было важно встать на путь прикладной математики, поэтому он ушел в ЦАГИ.
	
	Я точно не помню эти мотивы, но в разговоре со мной он как будто бы именно так
	мотивировал.
\end{quotation}

{\bf \punct Один из стаи.} Понятно, что слова  из <<Правды>>
о реакционной Московской математической школе и связей Лузина с ней
тоже должны были быть обсуждены. Это было сделано с превеликой
тщательностью, итог был скорее приемлимым для Лузина (хотя ему дополнительно
приплели и Кассо). Из заключения Комиссии:
\begin{quotation}
	H. H. Лузин был питомцем старой московской математической школы, принадлежавшей к наиболее реакционному крылу профессуры. Во время университетских событий 1911 г. (при министерстве Кассо) H. H. Лузин — тогда приват-доцент — остался в
	университете вместе с наиболее правым крылом, хотя в области своей науки он стремился к перестройке преподавания по образцу заграничных университетов%
	\footnote{Эта фраза положительна. Что касается 1911г., то - напомним - он в это время был в Гёттингене.}. Сам H. H. Лузин
	воздерживался от каких-либо явных политических выступлений в том или другом направлении.
\end{quotation}
 
 {\bf\punct Александров.%
 \label{ss:aleksandrov-polytics}} Было два противника Лузина, Александров и Колмогоров,
 которые удерживались в научно-этических рамках. Более того, Александров
 возражал против политических обвинений, или пытался преобразовать их
 в этические. Приведем несколько цитат.

Раз:
\begin{quotation}
	АЛЕКСАНДРОВ.
	И я должен сказать, что в это советское время я сделался свидетелем необычайной картины, картины,
	которой я никогда и нигде на свете не видел, картины такого исключительного научного
	энтузиазма и такого исключительного научного подъема, который охватил этого человека, большой человеческий коллектив, состоявший свыше чем из 36 человек, который
	весь концентрировался вокруг личности Н.~Н. Это были годы 1920—1921—1922 гг. Тут я
	думаю, что Н.~Н. выступал не только в роли создателя школы. Н. Н. выступал действительно в роли организатора молодежи, в роли лица, умевшего в эту молодежь, уже советскую молодежь, вселить какой-то исключительный научный энтузиазм. Тут, конечно, говорить о чем-либо антисоветском, мне кажется, совершенно невозможно, потому что
	H.~H. отдавал себе отчет. Он, имея абсолютное влияние на всех нас, настолько, что каждое его слово воспринималось действительно как абсолютная истина, если бы он желал
	нас настроить антисоветским образом, имел бы для создания таких настроений все возможности, которые никогда ни одному профессору не были даны, потому что вокруг
	него было общество людей, которые буквально на него молились. И, тем не менее, я, который был очень близок с ним тогда, и не только близок, но я должен сказать, что это был
	период очень большой дружбы, и дружбы не только с моей стороны к H.~H., но и его
	дружбы ко мне, настолько большой, что H.~H. находил возможным делиться со мною
	очень интимными сторонами своей жизни, быть очень откровенным со мною; {\bf и я
	жизнь Н.~Н. этого периода действительно знаю очень хорошо, — я категорически отрицаю, что какие бы то ни было антисоветские настроения H.~H. могли бы проявляться при
	полной его откровенности со мною.} Я ни одного проявления этих настроений не могу
	констатировать. Я опять повторяю. {\bf Я не знаю тех материалов, которые находятся в распоряжении «Правды», но мне казалось, что это противоречит тому восприятию личности
	этого человека, которое у меня есть, — чтобы в последующие годы он мог бы сделаться
	врагом советской власти.}
\end{quotation} 
 
Два: 
\begin{quotation}
АЛЕКСАНДРОВ. Я
 должен  сказать в  защиту Лузина  относительно  1920,  21  и  22  гг.  Лузин
тогда  пользовался  громадным  влиянием  у  своих  учеников.  Это  действительно  верно.
{\bf Однако должен сказать, я не помню ни одного случая, когда он, пользуясь этим влиянием, 
старался бы находить такие темы для антисоветских разговоров,  или пропускания
этих ниточек}.  Этого я констатировать не могу.  Ни разу такого  случая  не было.

ШМИДТ.
  Но разговоров о  положительной  стороне  большевизма,  например, тоже не
было.  Вообще в период гражданской войны на такие темы разговоры не велись.

АЛЕКСАНДРОВ.
  У  него  была такая  аполитичность  в  самом  цветущем  ее  состоянии,  но
активно антисоветским человеком, я могу с определенностью сказать, в те времена,  когда
я был его учеником, он себя не проявлял.
\end{quotation}

Три:
\begin{quotation}
	АЛЕКСАНДРОВ.
	 считаю Н. Н., — это может быть очень резкое выражение — интриганом
	\dots И если у H. H. была политика, то эта политика, мне кажется,
	имела чисто личный характер: H. H. желал иметь «популярность» среди таких рядовых
	научных работников. Он каждому скажет по комплименту. Ему хотелось иметь исключительно за себя голоса всех, не считаясь совершенно с объективным весом.
	{\bf Мне кажется,
	что факты, о которых говорит С. Л.[Соболев], они бесспорны, но они относятся к
	личному интриганству, но не к сознательному политическому вредительству.}
\end{quotation}

В принципе вот в следующем отрывке можно усмотреть политическую составляющую:
\begin{quotation}
	АЛЕКСАНДРОВ. В одном из пунктов резолюции указывается на уничижительное отношение H.~H. к советской науке. Я считаю, что тот характер дискуссии, который ведется у
	нас, является блестящим подтверждением этого пункта, потому что этот характер дискуссии ничем иным как уничижением всех здесь собравшихся, не может быть назван.
\end{quotation}

Но это все же мелочь на фоне личного ручательства Александрова%
\footnote{Скорее всего эти слова сам Александров считал истинными,
и, скорее всего, они, по-видимому, и были истинными. Лузин, по-видимому,
не понимал достижений московской математики за 1924-1936гг., в любом случае,
читая разные его тексты, опубликованные прижизненно и посмертно,
я не видел никаких признаков обратного.}.

\sm

{\bf\punct Хинчин.%
\label{ss:hinchin1936}}  Хинчин тоже был активным нападавшим,  был не чужд 
политических обвинений и был сторонником изгнания Лузина из Академии.
 Однако  он обсуждает степень справедливости самых разных высказываний
 и возражает, когда он считает их  несправедливыми. Выше уже было несколько примеров,
 вот еще несколько высказываний с политическим оттенком:
\begin{quotation}
	ХИНЧИН.
	Но все-таки факт, что в ущерб
	мировой науке и в ущерб престижу русской науки было то, что большинство из наших
	дореволюционных математиков печатались исключительно на русском языке, печатались исключительно в России. К сожалению, с работами Ляпунова дело вышло так, что
	для нас это просто досадно теперь, когда его центральная теорема в теории вероятности,
	представляющая собою, может быть, вообще самое крупное, что в этой науке за все
	время ее существования было сделано, прежде чем она стала известна научному миру,
	была два или три раза «переоткрыта» европейскими математиками. {\bf В конце концов,
	нужно позаботиться о том, чтобы наши статьи появлялись на иностранных языках, иначе
	это грозит престижу советской математики.}
\end{quotation}

Это сказано 9 июля, после второй статьи в <<Правде>> того же дня. Хинчин возражает не 
против отдельного высказывания, и не уточняет нюансы, он возражает против основного
смысла статьи! Напомню, что он сам в этой статье был объектом нападения, но
примерно то же самое чуть, но короче он говорил и 7 июля (см. выше  п.\ref{za-granitsej})

Следующий эпизод 1928 года, обсуждаются выборы в Академию, когда НИИ механики
и математики выставил кандидатуру Егорова и не выставил Лузина:
\begin{quotation}
		ХИНЧИН.
	Я тогда, наоборот, выставил кандидатуру Д.~Ф.~Егорова
	покойного, который не был избран. И я должен сказать, что я до сих пор жалею, несмотря
	ни на что, что он не был избран.
\end{quotation}

Напомним, Егоров был арестован, пробыл несколько месяцев в тюрьме,
 был выслан в Казань, где вскоре умер. Он не был заклеймен 
как враг народа, но он числился реакционером, и выступавшие на комиссии один за другим клеймили его в этом качестве.    

И снова Хинчин возражает газете <<Правда>>:
\begin{quotation}
		ХИНЧИН.
	Мне хотелось возразить против одного пункта этой статьи: «Мы знаем, откуда вырос акад. Лузин: мы знаем, что он один из стаи бесславной царской «Московской математической школы», философией которой было черносотенство и движущей идеей —
	киты российской реакции: православие и самодержавие. Мы знаем, что и сейчас он недалек от подобных взглядов, может быть чуть-чуть фашистски-модернизированных» и т.д.
	Насчет теперешних взглядов я ничего не могу сказать, но эти сведения мне кажутся фактически неверными.
	
	Дело в том, что когда H.~H.~Лузин был молод, когда складывались его академические
	и, вероятно, его политические убеждения, насколько мне известно, (я, правда, должен
	сказать, что я не очень следил за этим) — в это время в Московском университете происходила борьба между реакционной, действительно черносотенной группой и другой
	группой, которую возглавлял Д.~Ф.~Егоров и которая стремилась к европеизации, в буржуазном смысле этого слова, Московского университета. И Лузин целиком принадлежал к этой второй группе. Он отнюдь не был связан с самодержавным управлением.
	Наоборот, его отношение к этой группе было всегда очень резкое. Он принадлежал к
	западникам%
	\footnote{Слово <<западник>> в данном контексте положительно.}.
\end{quotation}

Тема для обсуждения своеобразна, и тем не менее.
 Далее:
\begin{quotation}
	ХИНЧИН. Я считаю, что основных причин ухода H.~H. из Московского университета
	было две, и первая из них та, о которой сегодня не говорилось, а именно, то обстоятельство, что Институт математики, в котором H.~H.~Лузин работал, высказался в свое время
	против избрания его академиком, когда выставлялась его кандидатура в Академию наук.
	Это, как я мог наблюдать совершенно точно, создало первую трещину, которая, по-моему, и до настоящего времени не заросла вполне между H.~H.~Лузиным и Московским
	университетом. Здесь, стало быть, причина личного характера.
\end{quotation}

Политическое обвинение Лузина частично переводится в личное \dots

Мельком упомянем о Соболеве. Он тоже участвовал в политических
обвинениях, но и он  был неоднозначен. О старой Московской школе,
которую обличала <<Правда>>:
\begin{quotation}
	СОБОЛЕВ.
	Теперь о Московской школе. Связь с Московской школой. Я считаю, что связь с Московской школой мы не можем ставить ни в какую вину. Всякий человек в какой-нибудь
	школе воспитывается. И направление этой школы не имеет в себе никаких вредных тенденций, против которых можно было [бы] бороться. Я считаю, в частности, что дескриптивная теория функций связана с большими философскими задачами, и здесь H.~H. имеет серьезные большие заслуги. {\bf Я согласен с С.~Н. [Бернштейном] также и в том, что действительно в части создания современной Московской школы математиков H.~H. принадлежат громадные заслуги — хотел он этого или не хотел, но он был
	тем центром, который создал всю Московскую математическую современную школу.
	Это абсолютно ясно.}
\end{quotation}

Кстати, сам Соболев к Московской школе не относился.

\sm

{\bf \punct Несколько замечаний о политических нападающих.%
\label{remarks-polytics}}
Из математиков основную атаку с политической стороны вела группа
Шнирельман--Люстерник--Гельфонд. На собрании в МГУ
главным выступающим была Яновская, и по смыслу ее выступление 
было близко к речам этой участникам группы (см. ниже п.\ref{ss:yanov2}.г.
 По политической линии наседали также эпизодически появлявшиеся на комиссии Сегал и Аршон.
 
   Сам по себе нападал Шмидт. Он разбил предложенную Бенштейном
   идею, что Лузин -- крупнейший советский математик (что, впрочем, к
   политике не относилось), активно участвовал в исследовании вопроса
   о Лузине и Московской математической школы (примерно в духе процитированного выше
   заключения комиссии). Кроме того, он предлагал исследовать вопрос о поведении
   Лузина за границей.

	\sm
	
{\bf \punct Несколько замечаний о нападении.%
\label{ss:modern}}	
В связи с этой историей у современного пост-советского читателя возникает несколько вопросов.

\sm

1) {\it Понимали ли они, что Лузин мог сесть?}
Мы видели свидетельство Ефремовича о том, что не понимали,
а также видели реплики Соболева и Хинчина (<<не сажать же мы его собираемся>>).
Если бы Лузин
был исключен из Академии, это могло бы очень плохо кончиться, причем не только
для него, но и для нападавших. Однако участники истории не могли подозревать
о предстоявшем 37 годе (да и вообще о нем вряд ли кто-либо подозревал). Они привыкли произносить речи на собраниях, обличать
врагов и друг друга. За подобными обличениями в 37 году с легкостью следовали аресты обличенных
и обличателей, но этого опыта у присутствовавших тогда не было. 
Речь шла о <<политическом недоверии>>, а не об аресте.

У нас есть свидетельство о восприятии математиками тогдашней жизни, оставленное
Андре Вейлем [Andre Weil] \cite{Weil}, который приехал на топологическую конференцию  осенью 1935 года
и несколько недель пробыл в Москве:
\begin{quotation}
	«Это международное математическое собрание, впервые проведенное в СССР, было
	также, пока был жив Сталин, последним%
	\footnote{Данные о международных математических собраниях
		1930-1953гг.
		в СССР такие.
				 Первому Всесоюзному математическому съезду 1930г.
		пытались придать отчасти международный характер, см. \cite{LL-conferences-1},
		в частности иностранцы читали там половину пленарных докладов. В 1934 году в Москве была большая
		международная конференция
		по тензорной дифференциальной геометрии, организованная Каганом, с участием Э.~Картана (Eli Cartan),
		Бляшке (Wilhelm Blaschke), Кэлера (Erich K\"ahler), Схоутена
		(Jan Schouten) и др. 
		Упоминаемая Вейлем топологическая конференция организовывалась Александровым. 
		В переписке Александрова и Колмогорова за 1936 упоминается подготовка к (несостоявшейся) международной конференции
		по функциональному анализу (обсуждался недостаток средств и их
		перетягивание разными организациями).
		В 1937-38гг. у советского руководства были иные занятия, далее началась Мировая война.
		  В 1945 г., с 15 июня по 3 июля, в Москве и Ленинграде   проходила Юбилейная сессия,
	посвященная 220-летию Академии наук. В ней, в частности, участвовали,
122 гостя из 18 стран. По математике упоминаются \cite{LL-conferences-2} доклады
Ж.~Адамара, С.~Стоилова (Simion Stoilow), М.~Фреше (были ли там еще какие-либо
иностранные математики,
автор выяснить не смог). Видимо, это был широкий жест, демонстрировавший восстановление
международных научных связей после Войны. А с 1946-1947г. началась Холодная война,
и в самом деле началась жесткая политика изоляционизма. Я не понял, когда
была следующая международная математическая конференция СССР, но совсем не сразу после
смерти Сталина. В 1961г. в Киеве был Симпозиум по нелинейным колебаниям, а в 1966г. -- Конгресс в Москве.
Я что-то ничего более раннего в анналах не раскопал.\label{fo:220}}.
		То, что оно состоялось в 1935 г. при участии целой плеяды математиков всех стран, среди которых были самые выдающиеся
	из тех, которые кто близко, а кто издали имели дело с топологией, в ретроспективе
	казалось аномалией и почти чудом. В эти самые времена хватало математиков, чтобы видеть
	в этом знак начала либерализации советского режима, да я и сам был недалек от того, 
	чтобы согласиться с этой иллюзией, которую громкие московские
	процессы не замедлили жестоко развеять».
\end{quotation}  

Что касается смягчения системы, то оно в середине 30х в самом деле было,
а смягчение на математическом направлении происходило с 32 года. Кстати,
 <<Правда>> оказалась  орудием борьбы в математическом мире  не в первый и не в последний раз, о чем ниже тоже пойдет речь.
 
 Может, после следующих трех параграфов точка зрения собеседников Вейля станет более понятной читателю.

 \sm
 
 {\it Но ведь у них же на памяти был случай с Егоровым....} 
 Ниже в п.\ref{ss:losev} мы обсуждаем, что сейчас известно об обстоятельствах ареста Егорова.
 По-видимому, арест не был связан с его профессиональной деятельностью. Надо думать, что так
 думали и современники.
 
 Приведем две цитаты:
 \begin{quotation}
 	ШМИДТ...
 	 И этот взрыв [выступление математической молодежи] имел драматические последствия для Егорова. Здесь не было связи, но по времени все это совпало.
 	 Егоров был арестован, затем он был выпущен и умер.
 \end{quotation}
 
Теперь Ефремович:
\begin{quotation}
	 Егорова тоже взяли в 30-м году. Формально его обвиняли в том, что он член какой-то религиозной секты. Вероятно, так и было. Но было в этой секте чуть ли не человек 10-15. Причем она была где-то на Кавказе, а в Москве он был чуть ли не один. 
\end{quotation}
 
Интересно, что сообщение Ефремовича 60 лет спустя не вполне точно, но
 отражает реальность, см. п.\ref{ss:losev}.

 \sm
 
 {\it Как же люди могли такое говорить?} В самом деле, впечатляет.
 Но стоит иметь в виду, что русская/советская/российская интеллигенция
 всегда была очень политизирована и определенные общественно-политические вопросы
 имеют для нее фактически религиозное значение (знаменитый сборник <<Вехи>> \cite{Vehi} определял ее как секту).
 Взгляды могут меняться на формально противоположные при
 остающейся неизменной структуре мышления.
 Особое возмущение сегодня вызывает не сам факт политизации красной интеллигенции 30х годов (современная
 либеральная интеллигенция политизирована не меньше), а тем, что зеркало в данном случае,
 не переставляя голову с ногами, меняет местами лево и право.
 
 \sm
 
 {\it Неужели непонятно, что эти люди действовали под действием страха?}
 Смотря кто. Речи под принуждением или под диктовку имеют определенную специфику. Стиль <<как мать говорю и как женщина>>
 в Академической стенограмме не наблюдается, 
 люди пытаются отстаивать то, что
 считают нужным, оспаривают обвинения, которые считают несправедливыми,
 имеют разные мнения, того не скрывая, и возражают газете <<Правда>>. Больше всех 
 возражал Хинчин, но  против отдельных  правдинских сентенций 
 возражали многие. 
 
 Эти слова относятся ко всем, кто нападал в Академической комиссии.
 А так страх конечно был. Те, кто был против этого, на собраниях должны были боялись высказываться
 (а представьте себя окруженным  охваченными единым духовным порывом партийцами и математиками...).
  Наверно кто-то выступал, потому что велели... Но это все вне дошедших до нас протоколов.

 \sm

{\bf \punct Шнирельман.%
\label{ss:schnirelman-NKVD}} И был один человек, сильно выделявшийся из прочих. 

\begin{quotation}
ШНИРЕЛЬМАН.
\dots	там, где H.~H. был Председателем группы как крупный ученый — это плюс, но там, где он выступал как Председатель группы, где он вел себя необъективно, этот плюс обращается в минус. Такие факты,
	как дача отзывов, во многих случаях заведомо несправедливых, доказывают именно полную безответственность, которую обнаруживал H.~H. в таком ответственном деле, которое ему было поручено. Это является, безусловно, {\bf преступлением}. Ряд таких фактов,
	идущих вместе, дают возможность утверждать, что, во всяком случае, несмотря ни на
	какие научные заслуги, H.~H. как научно-общественная фигура за последние годы является отрицательной фигурой.\dots
	
	Последний пункт в статье, которая была опубликована в «Правде», содержит обобщение для этих фактов и квалификацию причин, обобщающих, дающих объяснение деятельности H.~H. последних лет. Не имея всех данных по этому поводу, я не считаю себя в
	праве высказать какое-нибудь суждение здесь. {\bf Единственно, что я могу сказать, это то,
	что те инстанции, которые будут иметь в своем распоряжении этот материал, проверят
	эти предположения и высказывания в «Правде», и если они оправдаются в какой бы то
	ни было мере, то я, со своей стороны, считаю, что должны быть приняты самые решительные меры, и что в этом случае нельзя считаться ни с высокой ценностью его научной роли, ни с той объективно полезной ролью, которую он сыграл в развитии советской математики.}	
\end{quotation}

Слово <<преступление>> присутствует в Стенограмме лишь один раз. Другие ораторы
говорили о <<принесении вреда>> или <<вредительстве>> (что тоже не радует глаза, но все же это предшествующая ступень).

Значительно хуже выделенный жирным текст в конце. Слово <<инстанции>>
наводит на мысль об НКВД, но в принципе может допускать какое-либо иное истолкование.
К сожалению, <<самые решительные меры>> -- это  не исключение из Академии.
Последнее -- прерогатива Комиссии (в которую входит Шнирельман) и общего собрания Академии, а не <<инстанций>>. 
То есть речь идет о чем-то более решительном, это 
может быть лишь арест или ссылка. Круг с <<инстанциями>> замыкается. Это
те самые <<инстанции>>.

Приведенное высказывание не единственно.

\begin{quotation}
	ШНИРЕЛЬМАН.
	Я думаю, что этот вывод, сделанный со всей отчетливостью, уже сам по себе и с полной очевидностью доказывает, что H.~H. нельзя доверять никакого научно-общественного дела. Это надо сформулировать так, как здесь в действительности.
	{\bf Второй вопрос: является ли он активным контрреволюционером или сознательным,
	хотя бы и своеобразным вредителем? Я думаю, что ответа на этот вопрос мы пока не
	можем дать, так как у нас нет материалов. Думаю, что для того, чтобы выяснить этот
	вопрос (а выяснить это необходимо, так как он очень важный), мы должны сделать следующее: доверить это тому авторитетному органу, у которого есть в распоряжении весь
	материал.}
\end{quotation}

Здесь никаких сомнений быть не может: Шнирельман говорит о <<передаче дела в НКВД>>.
Причем
<<мы должны>>.

\sm

{\bf \punct Гибель Шнирельмана.%
\label{schnirelman-gibel}}
Александров, Гнеденко и Степанов в первом выпуске <<Историко-математических
исследований>> \cite{AGS} писали о Шнирельмане следующее:
\begin{quotation}
	Л. Г. Шнирельман родился в 1905 г. в семье учителя 
	в г. Гомеле. Уже мальчиком 12—15 лет он обнаружил  
	необыкновенные математические способности. В 1921 г. он  
	поступил в Московский государственный университет и  
	окончил его через два с половиной года. В 1927—1928 гг., 
	вместе с Люстерником, он развивал топологические методы 
	в анализе. В 1929 г. Шнирельман уехал в должности 
	профессора математики в Новочеркасск. Там он  
	натолкнулся на изумительно простую и смелую идею  
	рассматривать общие свойства произвольных числовых  
	последовательностей. Этим, как мы говорили, был сделан шаг 
	в решении проблемы Гольдбаха, была обобщена проблема 
	Варинга и пр. Осенью 1930 г. он возвратился в Москву. 
	В 1933 г. был избран членом-корреспондентом Академии 
	наук СССР. Позднее он занимался разнообразными  
	проблемами анализа и топологии. 
	
	\dots\dots\dots
	
	Нелепая смерть в 1938 г. вырвала Л. Г. Шнирельмана 
	из рядов московских математиков\dots 	
\end{quotation}

Позже где-то было написано, что он отравился газом на кухне... Позже...
Рассказов было несколько и содержание их различно.

Приведем отрывок из интервью, взятого Е. Б. Дынкиным у А. М. Яглома \cite{Yaglom}.
Последний год перед Войной Яглом занимался у Гельфонда в математическом кружке.
 \begin{quotation}
{\sc А.М.} А потом пришел Гельфонд и сказал, что Шнирельман умер... Я не помню, плакал ли он, во всяком случае,
 он был страшно расстроен. Сказать, что Шнирельман покончил жизнь самоубийством, Гельфонд не мог.
 Он был в чем-то очень советский человек. Мог и бояться просто. 
 Во всяком случае, было сделано все, чтобы мы не понимали, что происходит.
 
{\sc Е.Б.}: Ну, разумеется. А ты знаешь, как это все было? Мне Софья Александровна Яновская рассказывала, 
что он оставил записку: <<Я умираю честным перед товарищами и советской властью>>. Его заставляли доносить.

{\sc А.М.}: По-видимому, это была громадная трагедия для советской математики,
что умерли два бесконечно талантливых и, как говорят, бесконечно порядочных человека%
\footnote{Очень типичная для подобных случаев характеристика. В подобном контексте
употреблялось также слово <<честнейший>>. Урысон,  впрочем, тут не причем.}: Урысон (совсем рано),
и Шнирельман.

{\sc Е.Б.}: Ну, да. Ну, Урысон — это был, так сказать, несчастный случай,
а Шнирельман — это преступление. Софья Александровна мне говорила, 
что работника НКВД, который его вербовал, расстреляли. Но от этого, конечно, 
Шнирельман не вернулся к жизни.
\end{quotation}

Из статьи Н. Я. Виленкина
<<Формулы на фанере>>:
\begin{quotation}
	 Утверждают, что он принял решение открыть газ на 
	 кухне после того, как был вызван на Лубянку 
	 и получил предложение «освещать» одного 
	 из видных партийных деятелей, в дом  
	 которого был вхож. Друг Шнирельмана Л.~А.~Люстерник сказал мне однажды, что видел  
	 предсмертную записку Льва Генриховича (и,  
	 насколько я помню, хранил ее). Было бы весьма 
	 интересно найти ее [Записка эта не известна].
\end{quotation}

Из Воспоминаний Юшкевича \cite{Yush-school}:
\begin{quotation}
 В 1930-е гг. Гельфонд особенно подружился с Л.~Г.~Шнирельманом.
 Это был чрезвычайно умный математик и интересный человек.
Нередко мы собирались тесной компанией в Доме ученых — оба
друга с женами и Шнирельман (редко кто-либо еще) — и он 
рассказывал нам фантастические, придуманные им истории с элементом
мистики и ужасов. Самоубийство Шнирельмана 24.09.1938 г. не
только сильно огорчило, но и обеспокоило Гельфонда. О причине
этого поступка, по-видимому, Гельфонду известной, он говорить
не хотел, я знаю только, что ему стоило больших усилий получить
разрешение напечатать в «Успехах математических наук» 
некролог.
\end{quotation}

Из записок Голубева 1942 <<О судьбе Эвариста Галуа>> \cite{Golubev-Omega}
(записки писались <<в стол>>).

\begin{quotation}
 А вот другой пример, трагический и грустный.
Среди молодых математиков, выросших около Егорова, был некий Лев Генрихович
Шнирельман. Человек не без таланта, но со всеми недостатками человека, мнящего себя
гением. В теории чисел он нашел результат очень интересный, обративший на себя
внимание [Эдмунд] Ландау  и других европейских математиков. Итак, гений был
общепризнан. А дальше начинаются вещи, о которых лучше бы не писать в его
биографии: травля его учителя Егорова, присвоение им себе каких-то «сыщицких»
функций в среде математиков, травля Лузина.

Таланту было мало, а фимиам курился вокруг трусами, прохвостами и добрыми
приятелями непрерывно. Ждали гениальных открытий, а их не было. Ждали всемирной
славы, почестей, а их тоже не было. Попытка устроиться на работу в Университет ни к
чему не привела: лектор был никудышный, никто его не слушал, администратор был
вообще ни к чему не пригодный. В результате ни славы, ни признания, ни денег. Бедняга
не выдержал и отравился газом. Здесь погубили друзья и собственное самомнение,
выросшее на почве неумеренных надежд и похвал.
\end{quotation}

Отрывок перекликается со Стенограммой.

И последнее -- воспоминания Понтрягина:

\begin{quotation}
{\bf Шнирельман был незаурядный, талантливый человек с большими странностями.
Было в нём что-то неполноценное, какой-то психический сдвиг.} Я помню, как трудно было ему уйти
от меня из гостей: он останавливался в прихожей и не мог двинуться дальше. Тогда говорили,
он не имел никаких успехов у женщин и это сильно угнетало его. Кроме того, с ним произошло 
большое несчастье в смысле научного творчества. 
Он сделал выдающееся научное открытие, дав первое приближение к решению теоретико-числовой проблемы
Гольдбаха. Этот успех грубо исказил его отношение к математической проблематике.

Ему принадлежала следующая формулировка: «Я не хочу заниматься промыванием золота, 
я хочу находить только самородки». Ясно, однако, что найти самородок можно, только промывая золото 
и подбираясь к самородку постепенно.

Он отказался от этого пути и утратил творческую инициативу. Когда это произошло, 
он впал в полное уныние и говорил часто мне: «Имеет ли право жить человек, который уже ничего не делает,
а в прошлом сделал что-то замечательное?» Я утешал его как мог. Кончилось это трагически: 
Шнирельман преднамеренно отравился. Я помню, как Люстерник встретил меня на вокзале,
когда мы с матерью возвращались с юга, и сообщил о происшедшем несчастье.

В то время Шнирельман жил уже в хорошей квартире вместе с матерью.
Она видела, что с ним происходит что-то неблагополучное, и следила за ним. 
Однажды ночью она была чем-то очень встревожена и хотела даже посмотреть, что с сыном.
Но, подумав, что он спит, не решилась пойти к нему. Утром обнаружила, что он закрылся в кухне, 
заложил все щели и пустил газ. Когда она обнаружила его, он уже был безнадежно мёртв, хотя ещё и не остыл... 
Так трагически кончилась жизнь Льва Генриховича Шнирельмана.
\end{quotation}

Насколько я понимаю, что случилось на самом деле остается не известным. У нас нет
оснований сомневаться в словах Дынкина, но он является не первичным источником,
а передает слова Яновской. Нет оснований сомневаться в свидетельстве Виленкина,
но буквально он свидетельствует лишь о записке. С другой стороны, Голубев  и Понтрягин знали то, что им
сказали, а им могли сказать не всё.

\sm

{\bf\punct Грань, которую никто не перешел.%
\label{ss:gran}}
Мы снова цитируем Стенограмму.:
\begin{quotation}
	ШНИРЕЛЬМАН.
	Мы понимаем под Московской школой то течение, положительное математическое течение, которое сыграло существеннейшую роль в развитии советской науки, развитие
	которой в значительной степени связано именно с именем Николая Николаевича.
\end{quotation}

\begin{quotation}
		ШНИРЕЛЬМАН.
	Я знаю его с 1921-го года, и я могу судить об этом только
	косвенно, поскольку я был тогда студентом. Но у меня сложилось такое впечатление, что
	с 1916—1917 гг. по 1922—23 гг. — этот период - золотой период деятельности периода деятельности H. H. в Московском университете. Это тот период, который независимо от
	того, какие отдельные теневые стороны имеются у H. H., следует оценить как период,
	принесший громадную пользу советской науке.
\end{quotation}

Заметим, Шнирельман готов отправить Лузина за Можай, но научные заслуги Лузина
он признает.

\begin{quotation}
	ШНИРЕЛЬМАН.
	Егоров — фигура весьма интересная, и мне кажется, что для выяснения
	этой фигуры ее можно разделить на две части: с одной стороны, — Егоров как политическая фигура, и с другой стороны, — Егоров как научная фигура. Политически — это
	настолько реакционная фигура, что он был реакционнее некоторых из тех зубров, которые были здесь упомянуты. Это факт, которого я не был, конечно, свидетелем, но это
	все знают. В 1911 году, в самые мрачные времена Кассо, Егоров остается в Университете в качестве одного из немногих профессоров, но гораздо реакционнее его в науке был
	Млодзеевский. И между тем он ушел оттуда и остался без государственной службы. Так
	что эту свою физиономию последовательно реакционную Егоров сохранил, кажется, до
	последних дней. Достойно удивления то обстоятельство, что очень долго, до 1930 года,
	эта реакционная политика господствовала и могла в какой-то степени проявляться.
	
	Что касается научной стороны, то здесь нужно Егорова характеризовать, как это
	было отмечено, как западника.
	
	Что касается взаимоотношений Лузина, то мне кажется, позволительно полагать, что
	Лузин являлся непрерывным продолжателем Егорова.
\end{quotation}

Нельзя сказать, чтобы этот монолог мог порадовать слух. Но Шнирельман снова отделяет научную сторону
от политической (слово <<западник>> здесь положительно). 

В целом по Стенограмме  разные ораторы шельмуют Лузина по политическим
мотивам, но никто не ставит под сомнение его заслуги как математика
(отрицалось то, что он крупнейший из советских математиков, но, с одной стороны,
отрицалось справедливо, а с другой -- на этот статус, надо думать, в глубине души претендовали многие
присутствовавшие). Егоров,
соответственно политический реакционер, но положительная научная фигура.
И никто, даже Шнирельман, не переходит эту черту.

\sm

На первый взгляд, может показаться, что признание линии между  наукой с одной стороны,
и идеологией и политикой с другой, было само собой разумеющимся. Это не так.
В советских биологических дискуссиях настоящие биологи  уже к 1929-1930 году были за
этой чертой (ниже есть немного об этом, п.\ref{ss:biologi}). За нее вышли и советские математики к 1970ому году. Боюсь, что и из-за нее так и не вернулись их наследники.

\vspace{22pt}

{\bf \punct Из дискуссия 11 июля.}
\begin{quotation}
СЕГАЛ. Разрешите, может быть и задним числом, но сделать некоторое дополнение. В
заключении сказано, что «исходя из всего этого, необходимо согласиться с квалификацией «Правды» Лузина 
как врага советской власти». Я считаю, что эта характеристика
слишком сжата, и надо добавить, по меньшей мере «... в своей деятельности приносивший вред развитию науки в СССР».
 {\bf Я предлагаю в конце добавить: «Показывает, что
H. H. Лузин своей деятельностью за последние годы приносил вред советской науке и
Советскому Союзу».}

КРЖИЖАНОВСКИЙ. Нет возражений? (Принимается.)
\end{quotation}

После перерыва на заседание был приглашен Лузин. В статье Демидова и Есакова
\cite{DeEs} утверждается, что в предложенной к обсуждению резолюции
эта формулировка исчезла.

Так или иначе, борьба выходила на более высокий иерархический уровень. О том, что
происходило в верхах, нам на сегодняшний день сообщили мало.
Имеющуюся в наличии обрывочную информацию мы обсудим ниже.
Но перед этим мы должны вернуться к бурным событиям на московском Физмате 1929-1932  годов.

\section{Рабочая молодежь, левые математики и красные профессора%
\label{s:workers}}

\COUNTERS

 	\epigraph{Математика в период пролетарской диктатуры должна
 		стать во втузе не только вспомогательным предметом преподавания, но и
 		могучим воспитательным средством для выработки
 		диалектико-материалистического мировоззрения. 
 		\newline
 		Основные положения, данные в учении Маркса, Энгельса
 		Ленина и их последователей, достаточны для того, чтобы помноженные на энтузиазм новых пролетарских научных
 		работников, овладевших и овладевающих математическими знаниями,
 		дать построение новой математики, сознательно руководимой методом
 		диалектического материализма. }{Л. А. Лейферт \cite{Leifert}, 1930}
 	
 	\epigraph{... основной задачей наших 
 		марксистских авторов и изданий в отношении математики должна была быть 
 		беспощадная борьба с математическим идеализмом, где бы и в какой бы
 		форме он ни проявлялся, неуклонное разоблачение его подчас очень
 		тонких хитросплетений, неусыпная работа над установлением и в 
 		области математики подлинного единства теории и практики на основе
 		борьбы за осуществление генеральной линии партии на основе 
 		практики социалистического строительства.}{С. А. Яновская \cite{Yanov-BSE}, 1931}

Здесь мы обсуждаем  биографии нескольких деятелей 1929-1931 годов до начала бурных событий:
три  относительно молодых человека Г.~К.~Хворостин, Д.~А.~Райков  и М.~Я.~Выгодский, два левых профессора, Шмидт и В.~Ф.~Каган,
красный профессор С.~А.~Яновская, авантюрист Э.~Кольман, а также вскользь нескольких своеобразных статистиков,
среди которых главным в тот момент был В.~И.~Хотимский, и (тоже вскользь) героев Ленинградского
математического фронта.
 
\sm

{\bf\punct Eropoву была объявлена война...%
\label{ss:egorov-war}} Цитируем статью И. Зайденвара \cite{Zaidenvar}%
\footnote{Мне не удалось найти об  авторе  статьи каких-либо данных.} 1930 года: 
\begin{quotation}
 В 70-х годах прошлого столетия, по инициативе Бугаева,
 Некрасова и др., организовалось Математическое Общество. Это
была кастовая, бюрократическая  организация, не сумевшая 
создать серьезной научной школы в математике.

Таким оно вошло в революционный период и, сохраняя
традиции академизма и философского идеализма, совершенно
не пыталось связаться с советской общественностью....

Система Егорова--Костицина%
\footnote{Как отмечалось выше, Костицын в 1928 упоминался в бумаге от Политбюро
	как один из возможных кандидатов на внедрение в АН СССР. На момент публикации этой статьи он уже был невозвращенцем,
	и ее автор идеологически объединяет Костицына с Егоровым.}
 держала в экономической 
зависимости наших молодых советских ученых: отказывалось достойным, талантливым, продвигались идеологически себе близкие.
 Общество не ставило, конечно, вопросов марксистской методологии
науки, зато в  его рядах было не мало активных
церковников. Связи с производством также не было, и Общество
проповедывало самодовлеющую «чистую» науку.
Вливавшиеся в университет широким потоком пролетарские
кадры студенчества своим классовым чутьем сразу почуяли недоброе.
Еще несколько лет назад проф. Eгopoву была объявлена
война. Вначале Егоров вынужден был уйти с должности
Председателя Предметной Комиссии, а в 1929 г. он был отстранен
от руководства Институтом....
\end{quotation}

 Стоит заметить, что в 30е годы на Егорова лево-радикальными математиками
было вылито изрядное количество помоев, но все они имели политическую природу, никто не обвинял Егорова в задвигании <<достойных и талантливых>>????.
Что касается упомянутой <<войны>>, то согласно Форду  \cite{Ford3}, вождями нападавших были Хворостин и Райков (хотя автор этих записок допускает, что
список вождей может быть не полон). Что касается Предметной комиссии, то Форд сообщает следующее

\begin{quotation}
 В протоколах собрания 6 февраля 1924 г.  упоминается
об избрании трех руководителей Предметной комиссии по 
математике. Председателем был избран Егоров, заместителем 
председателя — В.~А.~Костицын, секретарем — В.~В.~Степанов. На собрании 
присутствовало 25 человек, треть из которых составляли студенты,
включая и Хворостина. Голосование отразило эту пропорцию: 
около двух третей присутствовавших было за избрание руководства в
этом составе; треть была против. Все указывало на раскол между
преподавателями и студентами-активистами, действовавшими, 
скорее всего, по указаниям «идеологов». Уже на следующем собрании
5 марта 1924 г.  Егоров, Костицын и Степанов были 
смещены со своих постов....

На собрании 19 марта 1924 г.  на рассмотрение был
предложен новый список: председателем Предметной комиссии
выдвигался И.~И.~Жегалкин, заместителем 
председателя В.~В.~Степанов, секретарем — Г.~К.~Хворостин. Каждый был избран 
единогласно. Очевидно, было достигнуто определенное соглашение, 
благодаря которому пролетарское студенчество смогло получить в 
комиссии большее представительство: в ее руководство вошел 
Хворостин....
\end{quotation}

В конце 1925 года И.~Г.~Петровский (видимо, из-за болезни Хворостина) становится секретарем Предметной комиссии.
Фактически, сторона Егорова выходит из этого столкновения  победившей. Но мы должны присмотреться к тогдашним борцам.

\sm

{\bf\punct Хворостин.%
\label{ss:hvorostin-1}} Данные по ранней биографии Хворостина Гавриила Кирилловича (1900-1938) приводятся Фордом:

\begin{quotation}
 Хворостин [в своей автобиографии] указывает, что
он стал активно сотрудничать с большевиками во время 
революции, тогда ему было 17 лет; в Сибири он был арестован 
колчаковцами и даже посажен в тюрьму. В 1920 г. он поступил учиться на
рабфак им.~М.~Н.~Покровского, который закончил в 1923 г... Осенью
1923 г. Хворостин... станов[ится] студент[ом]
физико-математического факультета Московского университета.
Мы знаем, что уже к началу 1924 г. Хворостин активно 
действовал в Предметной комиссии по математике как представитель 
студенчества. Из-за болезни ему пришлось прервать обучение в
1925—1926 учебном году. Это обстоятельство объясняет почему в
конце 1925 г. Петровский заменил Хворостина в качестве 
секретаря Предметной комиссии по математике.
Хворостин пишет, что из-за большой загруженности 
общественной работой, он не смог завершить учебу на 
физико-математическом факультете. Хотя подробности в его автобиографии 
отсутствуют, она по своей тональности мало отличается от написанной
Райковым. Хворостин отмечает, что он потратил много времени на
создание специализации «прикладная математика», по которой он
собирался продолжить работу. С 1929 г. он был одним из 
студентов рабоче-крестьянского происхождения, обучавшихся 
специальной программе. Эти студенты часто отбирались по их симпатии к
советской системе, хотя при этом демонстрировать особый уровень
политической активности не было необходимости.
\end{quotation}

{\bf\punct  Райков.%
\label{ss:raikov-1}}
Вот автобиография Райкова Дмитрия Абрамовича (1905-1980), составленная при поступлении в аспирантуру 14.12.1929, \cite{Ford3}:
\begin{quotation}
«Отец -- кустарь До этого был рабочим пекарем. Сейчас 
— чернорабочий на мельнице. Родился я в 1905 г. в г. Одессе. Учился в
городском училище (до 3 класса). С 1918 по 1920 г. работал. В
1920 г. поступил на рабфак. В августе 1922 перевелся на рабфак
МВТУ в Москву. В конце января 1923 поступил в типографию и
перевелся на вечерний рабфак им. М.~Н.~Покровского, который и
окончил в июне 1923 г. В конце 1923 г. поступил на Физмат 1-го
МГУ на математическое отделение, но до конца 1926 г. почти не
учился, т.к. работал в типографии и вел большую общественную
работу в КСМской ячейке Физмата. Числился на 2 курсе. С 
начала 1927 г. ушел из типографии и нажал на занятия. Летом 1927 г.
перескочил со 2 на 4 курс. Весною 1927 г. учеба сорвалась 
вследствие борьбы с оппозицией (был секретарем КСМской ячейки).
Второе полугодие 1927/28 учебного года работал в семинаре
проф. А.~Я.~Хинчина по теории чисел и прочел там доклад: «Новое
доказательство трансцендентности числа е, данное Spaeht'oм».
Осенью 1928 г. также не занимался, так как был занят 
общественной работой. Во втором полугодии 1928/29 уч.г. работал по 
теории вероятностей у проф. А.~Я.~Хинчина и прочел на его семинаре
(часть курса лекций проф. А.~Я.~Хинчина была превращена в 
семинар) доклад об усиленном законе больших чисел. Кроме того дал
два новых доказательства теоремы Бернулли. В настоящее время
принимаю участие в работе семинара по алгебре и теории чисел.

С начала 1928 являюсь выдвиженцем. Состою сейчас членом
аспирантского бюро Института математики и механики 1-го МГУ в
качестве представителя от выдвиженцев.

{\bf Интересуюсь, собственно, методологией математики и ставлю
себе целью способствовать внедрению материалистической 
диалектики в математику. Но для этого нужно хорошо проработать самое
математику, а с другой стороны, развивать марксистскую 
методологию математики нужно для преобразования самой математики, а
не только для иллюстрации математическим материалом
диалектики.} Поэтому считаю нужным поступить аспирантом в 
Институт математики и механики 1-го МГУ и выбираю своей 
специальностью (на первое время) теорию чисел.
С мая 1920 г. состоял членом КСМ, с. мая 1928 г. являюсь
членом ВКП/б/».
\end{quotation}

Это начало жизненного пути известного впоследствии чистого математика. Пока же в  1929-30гг. он прославится как борец
с чистой математикой (в качестве какого он упоминался выше в п. \ref{za-granitsej}).

\sm

{\bf\punct Выгодский.%
\label{ss:vygodski-1}} Выгодский Марк Яковлевич (1989-1965) родился в Минске, учился в гимназии в Баку, в 1916
году поступил в Варшавский университет, находившийся в эвакуации в Ростове-на-Дону. Во время Гражданской войны был подпольщиком,
по одним сведениям в Ростове%
\footnote{В частности, утверждается, что он был стенографистом в штабе Деникина. Информация, по-видимому, исходит от О.~Г.~Шатуновской (надежность
этого источника оставляет желать лучшего, цитируем по книге Померанц <<Следствие ведет каторжанка>>)
\newline{\it
В Ростове жил тогда Марк Выгодский. Он работал стенографом в деникинской канцелярии и всё передавал нашим.
Марк стал плакать, перед этим одного нашего разведчика поймали, привязали к березам и разорвали пополам.
Он стал умолять меня: „Оля, не ходи. Я чувствую, что будет плохое, я тебя умоляю!“.
Я говорю: „Марк, что ты меня умоляешь! Меня послали, я должна идти“. Купили мне билет и почти до Курска я доехала.}
}%
, по другим в Азербайджане (а, возможно, и там, и там). Согласно
\cite{Kol-savvina-2}, 
<< {\it в 1923 г. окончил Московский университет и начал преподавать в Коммунистическом университете им. Я. М. Свердлова,
затем – в Институте Красной профессуры}>>.
С 1923 года учился в аспирантуре при московском Физмате, научный руководитель О.~Ю.~Шмидт
Хотя о Выгодском много писали (например, \cite{Vyg-nekrolog}, \cite{Rozenfeld}, \cite{Demidov-Vygodsky}),
его биография 1920-30х годов известна плохо. Например, лишь в последнее десятилетие
стало известно, что он в 1931-32 году был директором НИИ механики и математики в МГУ.

В Коммунистическом университете он использовал опыт, полученный у Деникина. А именно, в каталоге Российской
национальной библиотеки есть такая карточка (самой книги я не видел).

\begin{quotation}
Выгодский Марк Яковлевич (1895-1965)

Курс парламентской стенографии (по системе Гальсбергера)

М.П. Госиздтельство ``Мосполиграф``

1-ая Образцовая тип., 1923

64, 113стр., 27 см.

Перед загл. - Марк Выгодский - препод. стенографии в Коммун. ун-та им. Я.М.Сведлова
\end{quotation}

Приведем три биографических отсвета.
Из воспоминаний  Бескина \cite{Besk29} (отрывок сам по себе небезынтересен):
\begin{quotation}
Чистка в МГУ (1924 г.) касалась только студентов. Она не 
была публичной, студенты держали ответ только перед комиссией.
Она называлась <<академической чисткой>>, и ее официальной
целью было освободить университет от неуспевающих студентов%
\footnote{По рассказу Кузина \cite{Kuzin}, студентов, имевшие не более трех хвостов, чистке
не подлежали. На биологическом отделении московского Физмата преподаватели срочно (авансом)
позакрывали студентам хвосты, и на том дело кончилось. Указание на чистку, по-видимому, было общесоюзным,
или, по крайней мере, обще-РСФСР-овским.}.
Такая формулировка открывала возможность произвола, потому
что критерия успеваемости не существовало. Как я уже 
рассказывал, студент мог дойти до старших курсов; имея <<хвосты>> за 
младшие. Несмотря на это, он мог активно работать в какой-либо
области математики....

Чистка была двухступенчатой. Сначала студент являлся в
комиссию математического отделения. Эта комиссия могла 
принять только два решения: 1) оставить, 2) передать в 
вышестоящую факультетскую комиссию. В первом случае для данного
студента чистка на этом оканчивалась, во втором же случае он
Должен был явиться в комиссию Физмата, решение которой почти
всегда было отрицательным. Председатель комиссии 
математического отделения был студент (тогда это было возможно),
впоследствии профессор М.~Я.~Выгодский (1898—1965), у факультетской же комиссии было два председателя — О.~Ю.~Шмидт и
профессор физики А.~К.~Тимирязев.

Когда я предстал перед М.~Я.~Выгодским (он был один), он
спросил:

—- Читали ли вы марксистскую литературу?

— Читал, но мало.

— Напрасно.

Мое дело было передано в факультетскую комиссию. Имели
ли эти две комиссии какие-нибудь указания сверху относительно
отдельных студентов — я ис знаю. Я надеялся, что на 
заседании факультетской комиссии будет председательствовать
О.~Ю.~Шмидт, который меня хорошо знал, как своего слушателя,
но мне не повезло — на этот раз председательствовал
А.~К.~Тимирязев.

Комиссия была многочисленной. Ее члены сидели за длинным
столом. В середине длинной стороны сидел председатель, рядом с
ним усаживался чистящийся.

Как ни странно, я не помню, о чем меня спрашивали, но
помню, о чем спрашивали моего товарища по несчастью
И.~Н.~Бронштейна. Ему задали три вопроса:

1. Чем занимались ваши родители до 1917 года?

2- Какой сейчас политический строй в Афганистане?

3. Найдите производную от $x^x$.

................................

Факультетская комиссии отчислила почти всех, направленных
на ее усмотрение. Я был отчислен, как и мои сокурсники:
И.~Н.~Бронштейн — будущий профессор МГУ, Л.~А.~Тумаркин
(1904—1974), физик В.~Л.~Грановский н многие другие.
Пока я пытался апеллировать, дело неожиданно повернулось в
нашу пользу. Но это счастье было куплено дорогой ценой. Чистка
привела к трагическому результату: произошло несколько 
самоубийств.

Оказалось, что наверху к этому очень чувствительны. Шел
еще только 1923—1924 год. В университет немедленно прибыла
небольшая комиссия из ЦК ВКП(б) с неограниченными 
полномочиями. Без всякой волокиты она стала пересматривать дела
всех исключенных и тут же принимать решение о восстановлении
(часто условном).
\end{quotation}

Стоит заметить, что Выгодский в это время уже не был студентом.

А вот образец стиля самого Выгодского, 
(предисловие редактора перевода М.~Я.~Выгодского к книге Цейтен <<История математики в XVI и XVII веках>>, 1933),
где он излагает не вполне обычную точку зрения на перевод книг.

\begin{quotation}
Если в указанном выше отношении переработка Цейтена была невозможна, то в ряде других сторон {\bf русское издание значительно отличается от оригинала}....
 Принимая во внимание требования и особенности русского языка, мы должны были в ряде мест прибегнуть к очень вольной передаче оригинала,
 и нам кажется, что мы скорее виновны в недостаточно тщательном проведении этого приема, чем в злоупотреблении им.

Но этого мало. {\bf В качестве редактора я считал необходимым подвергнуть книгу обработке не только стилистической, но и по существу.} 
В целом ряде мест Цейтен трудно понимаем не только потому, что труден его стиль, но и потому, что он набрасывает картину слишком 
быстрыми движениями кисти. Для читателя, впервые знакомящегося с историей математики, такая беглая характеристика часто совершенно недостаточна.
{\bf В таких случаях я распространял изложение, а иногда и вовсе отказывался от цейтеновского текста, заменяя его другим.} Разумеется,
при этом всегда представлялось опасным разбить единство произведения. Я стремился к тому, чтобы избежать этой опасности и добиться того,
чтобы читатель не мог отличить интерполированный текст от текста оригинала.

{\bf Я счел себя также в праве сгладить те места, где, излагая по существу правильные взгляды, автор в несущественных деталях проявляет
некоторую наивность}, являющуюся следствием стихийности его материалистических воззрений. Разумеется, при этом момент субъективности всегда остается;
однако, я стремился к тому, чтобы всемерно избежать таких изменений в тексте, которые вносили бы в точку зрения автора существенные поправки.
Это, конечно, отнюдь не означает, что я согласен со всеми высказываниями Цейтена. Напротив, в ряде пунктов, и притом весьма существенных, 
взгляды автора являются, с моей точки зрения, в корне неправильными. В таких случаях, однако, я не считал возможным делать Цейтена ответственным 
за чуждые ему взгляды. 
Я предпочел сохранить в таких местах текст автора, сопровождая его редакционными примечаниями, иногда довольно пространными.
\end{quotation}

На Первом Всесоюзном математическом съезде Выгодский выступил с программным докладом по истории математики:
\begin{quotation}
 Все это однако указывает лишь на то, что классовые элементы
в математических теориях не столь очевидны, как в других 
теоретических дисциплинах. Как раз в последнее время мы имеем в математике
никем не отрицаемый кризис, и притом кризис ее основ, и сейчас мы
имеем все основания сомневаться в возможности однозначного ответа
на каждый поставленный вопрос. Не различие в устройстве мозга
у спорящих сторон обусловливает этот кризис. Его причины следует
искать в моментах социального порядка. Этим я отнюдь не хочу 
сказать, что одна группа математиков защищает, скажем, интересы 
промышленников, а другая — интересы аграриев. Нет, такое 
утверждение было бы злейшей пародией на марксизм. Но, прослеживая развитие
идеологии общества, мы должны будем в его классовой структуре
найти корни различных течений философской мысли, различных 
мировоззрений и их оттенков....

Именно признание решающей роли материальных факторов для
направления научной деятельности людей дает нам в руки критерий
для установления и проверки исторической гипотезы, а учёт влияния
классовой идеологии позволяет нам отделить в свидетельских 
показаниях современников, а также позднейших историографов и 
комментаторов, отделить в них элементы объективно исторические от 
идеологических искажений, обусловливаемых иногда непониманием эпохи,
иногда и просто апологетическими мотивами.
\end{quotation}

{\bf \punct Шмидт.%
\label{ss:schmidt}} Шмидт Отто Юльевич (1891-1956) закончил в 1913 году Киевский университет (научный руководитель Д.~А.~Граве) и 
был там оставлен для подготовки к профессорскому званию по математике. В 1916 году  опубликовал монографию
<<Абстрактная теория групп>> \cite{Shmidt-1916}, которой на следующие 20-30 лет предстояло стать основным учебником
по теории групп в СССР. Здесь можно спорить, связано ли это с тем, что Шмидт стал политической знаменитостью, 
или с достоинствами книги (автор настоящих записок склоняется ко второй точке зрения). В любом случае это была хорошая оригинальная книга
по новейшей в тот момент области математики.

Так или иначе, с конца 1916 года Шмидт втягивается в общественную деятельность. Вспоминает Делоне \cite{Shmidt-zhizn}:

\begin{quotation}
 «Как-то мы с ним пошли гулять в Голосеевский лес, и
я его спросил:

— Почему Вас так мало видно?

— А я занят организацией карточной системы в Киеве.

— Неужели это Вас интересует? — спросил я.

— Да, видите ли, Борис Николаевич, во мне два 
человека — человек науки, ума и человек действия, воли, и
эта деятельность удовлетворяет второго из них».
\end{quotation}

Летом 1917 он оказывается продовольственным деятелем в Петрограде, после Октябрьской революции
борется с саботажем в министерстве продовольствия, в марте 18ого  вступает в партию меньшевиков-интернационалистов,
далее входит  в коллегию Наркомпрода, и далее ведет бурную жизнь политического деятеля, которая постепенно смещается в сторону организации
науки, образования и книгоиздания.

Иногда его заносило и в МГУ...
Цитируем Александрова (об основании НИИ механики и математики):
\begin{quotation}
		В холодное ноябрьское утро 1920
		в одной из плохо отапливавшихся аудиторий Московского университета
		на Моховой было собрание, посвященное организации
		Московского института математических наук. Инициатором этого
		собрания был О.~Ю.~Шмидт, председательствовал Д.~Ф.~Егоров,
		а я был в числе участников.
		
		Если не ошибаюсь, это была первая попытка создания
		научно-ис\-сле\-до\-ва\-тель\-ско\-го института, по крайней мере в области
		близких мне научных дисциплин...
		Отто Юльевич был совсем молод.
		ему было 29 лет. Не удивительно, что его доклад, во многом
		расходившийся со старыми университетскими традициями, был
		встречен с осторожностью и даже недоверием. Помню те горячие
		прения, которые он вызвал, и помню, как Отто Юльевич
		непередаваемым обаянием своей личности, силой и ясностью
		своего ума, готовностью выслушивать и понимать самые различные
		точки зрения, умением убеждать своих противников
		сильной  и тонкой логической аргументацией сумел победить
		недоверие многих даже консервативно настроенных московских
		математиков.
	\end{quotation}
	
	А вот
	независимый источник, Костицын \cite{Kostitsyn}
	\begin{quotation}
	 Шмидт очень понравился Дмитрию Федоровичу [Егорову]
		тем, что с достоинством отстаивал в разговоре коммунистическую политику, отнюдь не стараясь понравиться собеседнику. Шмидт понравился
		и Б.~К.~Млодзеевскому. Через некоторое время удалось ввести его и в Московское математическое общество. Вел он себя с очень большим тактом
		и ни разу не дал никакой фальшивой ноты.
		\end{quotation}

Послужной список Шмидта поражает своей длиной.
Вот например список его должностей 
1929 года \cite{Shmidt-trudy} (причем я вижу, что по крайней мере одна должность пропущена):

\sm

{\footnotesize
1929 г. Главный редактор БСЭ [Большая советская энциклопедия], член коллегии Наркомпроса, член президиума
ГУСа[Государственный ученый совет], заместитель председателя комиссии по подготовке научных кадров, член Ученого
комитета при СНК[Совет народных комиссаров], заместитель начальника ЦСУ[Центральное статистическое управление] (до ноября 1929 r.), заведующий
Секцией естественных и точных наук Комакадемии, заместитель председателя экспертной
комиссии по Ленинским премиям (до 1930 r.), главный редактор журнала
<<Научное слово>>, член редакции журналов <<Вестник Комакадемии» и <<Естествознание
и марксизм», член главной редакции МСЭ[Малая советская энциклопедия], профессор 2-ro МГУ (до июля 1929 г.).
Апрель. На 2-й конференции марксистско-ленинских учреждений делает доклад
<<О работе Секции естественных и точных наук и задачи марксистов в области естествознания».

Присутствует на XVI конференции ВКП(б).

Июль. Назначен Правительственным комиссаром Земли Франца Иосифа%
\footnote{Шмидт, прибыв на Землю Франца-Иосифа, объявил ее советским владением.}
и начальником
экспедиции на ледоколе <<Г. Седов». Экспедиция водрузила флаг СССР на о-вах
Земли Франца-Иосифа, построила самую северную ($80^\circ20'$) научную станцию (в бухте
Тихой); совершила плавание в северную часть архипелага; прошла весь Британский
канал, провела ряд глубоководных гидрологических наблюдений  геологических
исследований.

Утвержден профессором и заведующим кафедрой алгебры 1-го МГУ, в
связи с чем оставил преподавание во 2-м МГУ.

Ноябрь. Назначен членом президиума Госплана.

Декабрь. Избран членом президиума Комакадемии.

Утвержден членом художественного совета Камерного театра.

Печатает статьи <<Задачи марксистов в области естествознания>>, 
<<Задачи СССР на Крайнем севере>>, <<Мы выполнили поручение правительства>> и др.
}

\sm

Понятно, что на части этих должностей он присутствовал по касательной. Его деятельность как математика вроде бы закончилась 
в 1916 году.... Но это было не совсем так. Цитируем Александрова, \cite{Alex-Shmidt}, \cite{Alex-Shmidt2}

\begin{quotation}
 В 1927 г. я встретился с Отто Юльевичем в Геттингене....
 В этот мировой научный центр во время летнего семестра
стекались ученые со всех концов земного шара...
Отто Юльевич приехал туда в научную командировку, оторвавшись
на короткое время от своих уже тогда очень разнообразных
занятий...

[Ему] было 36 лет... по его собственным словам, как  бы окунулся в математическую
 работу. Результат был выдающимся. Достаточно
было этих, по-существу нескольких недель досуга, чтобы Отто
Юльевич, овладев всем тем, что было сделано в области его
математической специальности за целое десятилетие, не только
оказался полностью на уровне последних достижений этой
науки, но и сразу же пополнил ее собственными первоклассными исследованиями.

Теорема теории групп, известная под именем теоремы
Шмидта, представляет собой одну из основных теорем современной алгебры...

Я помню заседание Геттингенского математического общества
под председательством Гильберта, на котором
О.~Ю.~Шмидт излагал свою теорему...

Я помню впечатление,  которое произвел этот доклад блестящий
не только  ло содержанию, но и по языку, по всей своей
внешней форме. Впечатление было огромным, несмотря на то,
что делался  он в таком месте. где люди были избалованы и
знали цену хорошим докладам....

Выступление О.~Ю.~Шмидта в Геттингене имело большой
и широкий успех. {\bf Посудите сами, приехал из Советского Союза
крупный общественно-политический деятель, делает блестящее
математическое открытие и столь же блестяще излагает его.}
Естественно, успех О.~Ю.~Шмидта стал своего рода сенсацией.
Отто Юльевича стали нарасхват приглашать не только в профессорские
семейства, но и на всевозможные официальные
приемы. Посыпались просьбы выступить на различных собраниях,
до собрания женщин-ученых включительно.... 
\end{quotation}

Несколько патетично, но по сути интересно.
Теорема Ремака (R. Remak, 1911) состоит в том, что разложение конечной группы в произведение единственно (в том смысле,
что как ни раскладывай, сомножители будут одинаковыми). У Александрова речь идет 
 о теореме Ремака--Шмидта или Крулля--Шмидта (варианты: Ре\-ма\-ка--Крул\-ля--Шмид\-та и просто Шмидта) \cite{Schmidt-1928}.
Шмидт рассматривал (вообще говоря) бесконечные группы с операторами%
\footnote{То есть берется абстрактная группа и в ней фиксируется семейство эндоморфизмов. Под это
	понятие подпадают, например,
модули (абелева группа и действие кольца на абелевой группе).}, которые удовлетворяют дополнительному условию: убывающие, а также возрастающие, цепочки инвариантных подгрупп
не могут быть бесконечными. Тогда разложение группы в прямое произведение инвариантных подгрупп единственно.

\sm

Так или иначе, возврат к математике оказался в бурной биографии Шмидта лишь кратким эпизодом, с 1929 года он начинает
свою карьеру полярного исследователя и организатора Северного морского пути 
(в часы досуга на ледоколе он писал заметки в математические журналы,
но, кажется, они не столь интересны). В качестве заведующего кафедры алгебры
(1929-1949) он  запомнился в памяти дальнейших поколений мехматян как ненаблюдаемая
на Мехмате персона (было ли это так с самого начала, неясно), однако
не вызывает сомнения, что Шмидт принес на московский Физмат продвинутую алгебру
(до Шмидта сколько-либо заметных алгебраистов в Москве не было),
и наложил свой отпечаток на деятельность
Куроша (неофициального зава) с коллегами. 

\sm

{\small
Уже в конце жизни, занимаясь проблемой происхождения Солнечной системы, Шмидт как математик тряхнул стариной. 
Есть известная проблема о поведении системы нескольких тел, притягивающихся по гравитационному закону Ньютона,
на длительном промежутке времени (очень сложная, надо сказать, проблематика).
Если тела два, то задача точно решается, и возможны лишь три вида движения тел.
Будем для простоты говорить о звезде и планете: планета вращается вокруг звезды по эллипсу, 
планета прилетает из бесконечности по гиперболе и по 
продолжению той же
гиперболы отправляется на бесконечность, и третий - промежуточный случай - когда планета (комета) прилетает 
и улетает по параболе.

Оказывается, что вероятность захвата в идеальной    задаче $n$ тел равна нулю%
\footnote{Быть может, здесь нужен неформальный комментарий. Пусть есть система из двух
тел, назовем их Солнцем и Юпитером. В систему влетает небольшое третье тело, 
назовем его астероидом. Как правило, он улетит по той же гиперболе, по которой прилетел.
Но в случае пролета вблизи Юпитера, астероид может ускориться, и тогда
он улетит по другой гиперболе. Но он может  замедлиться и перейти на эллиптическую
орбиту. На первый взгляд, астероид остался в системе. Но его орбита будет пересекаться
с орбитой Юпитера, и рано или поздно (но, быть может, очень нескоро) 
они окажутся одновременно вблизи точки пересечения орбит. Тогда астероид может
или замедлиться (и оказаться на новой эллиптической орбите, пересекающейся с орбитой
Юпитера) или ускориться с шансом уйти из системы.  Оказывается, что за бесконечное время
с вероятностью единица астероид  из системы уйдет.
\newline
Отметим, что в реальных планетных системах присутствуют иные факторы, кроме чистого  ньютоновского тяготения:
приливное трение, солнечный ветер, испарение газовых оболочек, столкновения и пр.
}. Есть идеальная теорема, см. \cite{Little}, о том, что система
из $n$ 
тел, в прошлом пребывавшая в ограниченной области пространства в течение бесконечного времени,
не может захватить прилетевшее из 
бесконечности $(n+1)$-ое тело (точнее, вероятность захвата на бесконечное время равна нулю).

Рассмотрим случай трех тел. В 1927-1932 г
Ж. Шази (Jean Chazy) \cite{Chazy} доказал теорему, о том, что движение трех тел при бесконечном отрицательном времени и при бесконечном положительном времени
<<одинаковы>>, в том смысле, что если все три тела были бесконечно далеки друг от друга, то они в итоге и останутся бесконечно далеки;
если два тела вращались вокруг друг друга, а третье было бесконечно далеко, то так же будет и в будущем и т.д.
(причем все будет в точности так, без поправки <<с точностью до событий вероятности ноль>>).
Работа была известной и считалась правильной.

В 1947 Шмидт заявил, что это неверно, и что с ненулевой вероятностью возможен частичный захват:
из трех звезд, находившихся на бесконечности может возникнуть двойная звезда, а третья звезда уйдет снова на бесконечность...
Он предложил численный пример, который по его поручению обсчитывали сотрудники. Дальше была дискуссия, о том,
до какой степени это можно считать доказательством, вопрос привлек внимание ряда математиков
(Г.~Ф.~Хильми, К.~А.~Ситников, А.~Н.~Колмогоров,  	Г.~А.~Мерман, В.~М.~Алексеев). В итоге было показано, что 
в принципе любое финальное движение может перейти в любое, но в некоторых случаях с вероятностью 0.
Как и утверждал Шмидт, с ненулевой вероятностью возможен частичный захват, а также обмен (одна планета прилетела к звезде,
а другая улетела),  подробней об этой задаче и ее истории см. об этом работы В.~М.~Алексеева \cite{Alekseev0}, \cite{Alekseev}%
\footnote{На самом деле, возможность частичного захвата или обмена с ненулевой
вероятностью показывается очень легко, для этого достаточно рассмотреть ситуации, близкие к прямому соударению. Это было сделано (занимавшимся общей топологией) Ситниковым, учеником Александрова. Чувствуется, что Шмидт
не проходил лузинскую выучку.}.}

Закончим цитатой из Александрова л 30х годах:
\begin{quotation}
Он продолжал руководить семинаром по теории групп, поддерживал тесное общение со своими учениками-алгебраистами,
оставался редактором основного математического журнала <<Математический сборник>>; как государственный деятель нередко призван 
был вмешиваться в решение многих, относящихся к математике вопросов. Отто Юльевича интересовало все, 
что происходит в математике и у математиков: он часто говорил о возросшей  интенсивности нашей математической жизни, 
об успехах математической молодежи. Живя в основном жизнью механико-математического факультета Московского университета,
я рассказывал Отто Юльевичу о том, что там происходит, часто жаловался на всевозможные неурядицы и трудности.
Помимо разумных советов и прямой помощи, Отто Юльевич никогда не упускал случая сказать: <<Не забывайте,
какие исключительные возможности вы имеете, каких ранее никогда не имели.>> И этот принципиальный оптимизм
не воспринимался как казенное утешение, а обладал внушающей силой.
\end{quotation}

{\bf\punct Каган.%
\label{ss:kagan}} Еще одним левым профессором московского Физмата был  Вениамин Федорович Каган (1869-1953).
С 1897 года он работал в Новороссийском университете (Одесса), с 1922 - в Москве. Крупных оригинальных результатов
в его биографиях, как будто, не приводится. Он больше был известен образовательной, и, особенно, издательской деятельностью.

В 1904 году  организовал в Одессе научное издательство <<Матезис>>. 
 В \cite{Kagan-kniga} приводится письмо Шмидта, руководившего в 1922  Государственным издательством,
\begin{quotation}
 ... Госиздат чувствует себя обязанным дать научную и научно-по\-пу\-ляр\-ную литературу. Волей-не\-во\-лей мы очутимся 
 перед необходимостью выпустить аналогичные книги, а отчасти даже переиздавать те же. Боюсь, не вышло бы плохо: Госиздат, как экономически подавляюще
 сильный, погубит возможность возрождения  Matheus'a, не воспользовавшись его навыками и традициями. 
 
 Поэтому предлагаю следующее. Нам с Вами объединиться. Госиздат дал бы капитал для возрождения Matheus'a
 под Вашей дирекцией. Это было бы автономное предприятие Госиздата СССР. \dots

 Предложение О.~Ю.~Шмидта было принято В.~Ф.~Каганом. Он приехал в Москву, и возглавил Научный отдел Государственного издательства. 
 Эта работа настолько его увлекла,
 что в течение почти десяти лет она занимала в его жизни не меньшее место, чем научная и университетская деятельность.
\end{quotation}

Еще одним успешным проектом Кагана был сериал <<Труды семинара по векторному и тензорному анализу>>, выходивший с 1933 года
(позднее это издание продолжал П.~К.~Рашевский).
Кстати, с Каганом связано внедрение тензорной геометрии в московской математике. 
Из его учеников стал знаменит П.~К.~Рашевский (Rashevsky--Chow theorem,
с которой началась его оригинальная деятельность, была опубликована в 1938г.
в трудах Пед. института им. Либкхнехта), еще один ученик -- В.~В.~Вагнер.

Яновская \cite{Yanov-PZM-1} в 1930 году упоминала о Кагане в таком стиле:
\begin{quotation}
Ибо даже в среде таких, подлинно близких нам и преданных советской
власти математиков, как В.~Ф.~Каган, не вполне изжиты еще конвенциалисткие
взгляды Пуанкаре.
\end{quotation}

Шмидт и Каган в июне 1930г. выступали на Первом Всесоюзном съезде математиков
с организационно-политическими речами, один открывал съезд, другой его закрывал.
В конце 1930г. Каган вслед за Шмидтом был атакован в Комакадемии и обвинен 
в идеализме (кстати, той же С.~А.~Яновской, позже идея была подхвачена также Э.~Кольманом).
В 1931 годах он на короткое время был арестован  (по-видимому, по издательским делам).
Подробностей в опубликованной литературе мне найти не удалось. По выходу он, по-видимому,
 остался влиятельным человеком (например, в 1934 году он организовывал международную конференцию
по геометрии в СССР).

\sm

{\bf \punct Яновская.%
\label{ss:yanovskaya-1}}
Как известно, год 1929 был «{\it годом великого перелома на всех фронтах социалистического строительства}».
Этот год и последовавшие за ним годы, 1930-33 были годами великих потрясений. Годы 1929-31 были временем плохо управляемого
социального переворота, когда сторонники системы получали карт-бланш на самые разные действия в соответствии
с собственным разумением. В 1929-1931гг большие потрясения пережила и система образования.

Цитирую Г.~А.~Гамова \cite{Gamov}, который после трехлетней поездки за границу вернулся  в Советской союз в разгар этой катавасии:

\begin{quotation}
 В России государственных
философов готовят в Коммунистической академии в
Москве и направляют во все учебные и исследовательские
институты, чтобы уберечь профессоров и исследователей от
идеалистических, капиталистических ересей. Государственные
философы обычно до некоторой степени знакомы с предметом
исследовательского учреждения, за которым им поручено
надзирать, будучи ранее или школьными учителями,
или прослушав в академии односеместровый курс по данной
дисциплине. Но они получают власть над научными директорами
учреждений и могут наложить вето на любой исследовательский
проект или публикацию, которые отклоняются от
«правильной идеологии».
\end{quotation}

Гамов, по-видимому, был человеком веселым и склонным к преувеличениям, но в общем что-то в этом роде и происходило.
На интересующий нас московский Физмат комиссаром  была отправлена выпускница Института красной профессуры
Софья Александровна Яновская (1896-1966). Цитируем ее автобиографию, приведенную в  \cite{Mints}:
\begin{quotation}
Родилась в 1896 г., 31 января, в местечке Пружаны,
быв[шей] Гродненской губернии в семье служащего.
Когда мне нe было еще двух лет, семья переехала в Одессу.
В 1906 г. поступила и в 1914 г. окончила с золотой
медалью Вторую городскую женскую гимназию. В том же
году поступила на естественное отделение Высших женских
курсов.
 По настоянию профессора С. О. Шатуновского
 была переведена на математическое отделение.
На курсах стала принимать участие в работе подпольного
Красного Креста [для помощи политзаключенным].
В ноябре 1918 г. вступила в подпольную
организацию большевиков. Перевозила инструкции
обкома и литературу через фронт. Была секретарем редакции
<<Коммуниста>>[печаталась в Одесских катакомбах] во время англо-французской интервенции
... По установлении Советской власти была секретарем
 редакции <<Известий>>, затем была послана с
 группой товарищей в Елисаветград [впоследствии Кировоград]
на ликвидацию последствий <<григорьевщины>> 
При отступлении из Елисаветграда вступила в ряды
Красной Армии: была политработником на фронте, завинформ.
отделом в газете <<Красная Армия>>, органе
Политуправления XII армии.

С 1920 по 1923 г. работала в Одесском губкоме партии:
завинформотделом, отделом учета, статистики и распределения.
Осенью 1923 г. была командирована губкомом на естественное
отделение Института красной профессуры,
на котором училась с 1924 по 1929 г. С января 1925 г.,
еще будучи слушательницей ИКП, начала преподавать в
ИКП математику, а в Московском государственном университете
- историю и философию математики.

С 1928 г. начала печатать работы по истории и методологии
математики.

С 1931 г. была утверждена пpoфeccopoм в Московском
государственном университете, ИКП и Академии наук...
\end{quotation}

Вот эпизод описанный арабисткой Х.~И.~Кильберг, приятельницей Яновской в 20-30е годы:
\begin{quotation}
 Во время отступления из Одессы белые захватили в
плен нескольких красноармейцев. Пленных они расстреливали
на мосту, и те падали в реку. Среди них, оказывается,
была и Софья Александровна. Пуля прострелила
высокую тулью шляпы. Софья Александровна упала в
реку, сумела выплыть и потом целую ночь отсиживалась
в воде в камышах.
\end{quotation}

У
Яновской было много учеников и ей посвящены многочисленные апологии, вот краткий список: \cite{Yanov-60}, \cite{Yanov-UMN}, \cite{Yanov-obituary} \cite{Mints},
\cite{Biryukov-yanovsk}, \cite{Birukov-Borisova}, \cite{Yanov-Bashmakova},
\cite{Yanov-BDU}, \cite{Kushner}, \cite{Levin-2011}, \cite{Yanov-reminiscences}, \cite{Yanov-anelis}, \cite{Anelis}, наиболее интересна работа Б.~В.~Бирюкова \cite{Biryukov-yanovsk}. Процитируем \cite{Yanov-Bashmakova}:
\begin{quotation}
Первые годы её работы приходятся на очень сложный и мрачный период в истории университета. Новая власть ставила под свой контроль науку
и высшую школу, одной из главных задач которой была борьба с реакционной профессурой, утверждение пролетарского студенчества, внедрение в
профессорскую и студенческую среду единственно верной марксистско-ленинской идеологии. Все это сопровождалось шумными пропагандистскими
кампаниями, чистками и поисками врагов. Борьба с «егоровщиной» 
и «дело академика Н.~Н.~Лузина» -- это лишь отдельные, хотя и наиболее известные, мрачные события в жизни московского математического
сообщества конца 20-30х гг. Особая роль в этих мероприятиях отводилась
	так называемым «красным профессорам» -- членам большевистской партии,
	носителям новой марксистской идеологии. Они должны были выступать не
	только рупором новой идеологии и политики, но и быть зоркими стражами,
	не допускающими идеологической крамолы, выявляющими и обличающими
	её, в каком бы виде -- даже самом невинном -- она ни появлялась. Эта «борьба» 
	нередко заканчивалась в застенках ОГПУ. К числу «красных профессоров» в Московском университете относилась и С.~А.~Яновская. В эти годы
	многие её поступки кажутся нам непонятными. Ведь она вместе с Э.~Кольманом -- 
	одной из наиболее одиозных фигур в советской науке тех лет --
	громила «реакционную профессуру» и в той или иной мере способствовала,
	мягко говоря, созданию тяжёлой атмосферы вокруг ряда известных математиков 
	(например, Д.~Ф.~Егорова, арест которого последовал в 1930 г.)...	
\end{quotation}

Вот что пишет биограф Яновской Б. В. Бирюков \cite{Biryukov-yanovsk}:
\begin{quotation}
Не существует никаких документов, говорящих о каком-либо участии С.~А.~Яновской в «деле Егорова». 
Трудно, однако, предположить, что она в нем не участвовала. Косвенно это подтверждается ее словами,
сказанными в Харькове на всесоюзном съезде математиков, которые будут приведены ниже, в связи с «делом» Н. Н. Лузина.
 Объяснение же тому, почему никаких документальных свидетельств ее участия в гонениях на Д.~Ф. не обнаружено,
я вижу в том, что свидетельства эти уничтожены. Как пишет Форд, С.~С.~Демидов и С.~С.~Петрова,
просматривая архивы ММО, установили, что документы, касающиеся реорганизации Общества,
которая последовала после смещения Егорова, в них отсутствуют. «Все страницы, относящиеся к заседаниям,
на которых могли обсуждаться арест Егорова и пути реорганизации Общества, вырваны». 
Кто-то заметал следы своих деяний ...

 Впрочем, не исключено, что соответствующие документы могут быть найдены в других архивах.
 Нужные протокольные записи либо извлечения из них могут сохраниться и в частных архивах. 
 Что касается тех, кто <<поработал>> над протоколами ММО, то вполне возможно составить представление о том,
 кто это был. Можно полагать, что в советское время заботиться об уничтожении этой части «исторической памяти» 
 кто-либо вряд бы стал. Значит, произойти это могло, когда началась Перестройка и особенно в постсоветское время.
 Вероятно, можно установить, кто обращался к протоколам ММО -в архиве могут сохраниться соответствующие записи.
 Поиск этот может быть успешен, так как число тех, кто заинтересован в утаивании этих данных, невелико%
 \footnote{Если все понимать буквально (и считать, что Бирюков продумывал каждое слово), то сказанное недвусмысленно. Но я, в свою очередь, в связи  испытываю некоторое (хотя и не абсолютное) недоумение.}. 
\end{quotation}

Завеса над тем, что ее интересовало  на Физмате, несколько приоткрывается следующей ее цитатой  \cite{Yanov-PZM-1} 1930г.,
\begin{quotation}
Уже не говоря о том, что ни в Московском, ни в Ленинградском математическом
обществе ни разу не было даже упомянуто слово «диалектика», и что,
если бы не год, проставленный на протоколе, заседание Математического
общества в 1929 году ничем нельзя было бы отличить от заседания в 1909 r., 
\end{quotation}

\sm

{\bf\punct Кольман Эрнест Яромирович.%
\label{ss:kolman-memoirs}} 
(Arnos\v{t} Kolman%
\footnote{В чешских источниках он также присутствует как Ernest Natanovi\v{c}.}),  06.12.1892, Прага — 22.01.1979, Стокгольм)
Кольман не относился ни к <<красным профессорам>>, ни к <<рабочей молодежи>>, но его, как
одного из ярких героев 1930-1931гг. удобно обсудить здесь.
Выше упоминались разные едва ли симпатичные детали биографий разных людей,
но это были детали, а вовсе не  полные биографии. Упоминавшиеся герои  в иных проявлениях и в иные моменты
жизни оказывались людьми
яркими и способными к положительной деятельности,  речь о чем еще пойдет ниже.
Про Кольмана последнего не скажешь, он выглядит человеком, в своем роде  <<идеальным>>.

Хотя о нем  много и с патетикой писали,  биография его не вполне ясна, и исследований по нему 
мне лично не известно, нет и его библиографии (нельзя не признать, что изучение его
творений -- удовольствие достаточно сомнительное). В 1930-1931гг он извергал поток статей
в марксистских журналах того времени, <<Под знаменем марксизма>>, <<Естествознание и марксизм>>,
<<Фронт науки и техники>>, <<Большевик>>. Предметом особого его пристрастия была математика,
но он был  широкой души человеком, и доставалось от него  физикам, химикам, биологам,
философам (разумеется), статистикам, инженерам. 
%Если судить по тому, на что есть ссылки, и тому, что я 
%разыскал, этот поток красноречия несколько  слабеет с 1933г.
%В 1938-1940гг. предметом его <<любви>> оказываются биологи. 
%Но эти сентенции могут оказаться неточными, библиографии его творений нет.

В качества образца приведем отрывок из его статьи марта 1931г., \cite{Kolm-boevye-voprosy}
<<Боевые вопросы естествознания>> (зачем мною ниже выделены жирным некоторые слова, будет объяснено потом):
\begin{quotation}
 Второй пример из более абстрактной науки, т.е математики. Существует так называемая
 московская математическая школа. Это — почтенная школа давнишних традиций,
 она всегда поддерживала три основных лозунга царизма: {\bf самодержавие, православие и народность},
 хотя бы тем, что профессора вроде Некрасова или Бугаева читали в Московском университете доклады,
 где доказывали, что анализ поддерживает православие, теория вероятности — народность,
 арифметика — самодержавие, и т. д. До сих пор правое крыло этой {\bf московской математической школы}, 
 во славе с контрреволюционером Егоровым, представляло еще достаточно живую струю 
 не только в смысле политических настроений, но и в методике, в тематике и в направлении исследования.
 Это та школа, корифеи которой с гордостью заявляют, что они занимаются только теми частями математики,
 которые не имеют никакого практического приложения. {\bf Эта школа очень тесно связана} с международной контрреволюцией,
 с французским комитетом интервенции, с такими людьми, как известный французский 
математик Борель, участвующий в работах Комитета французской тяжелой промышленности,
или итальянский математик Энриквес, — все известные {\bf фашисты}, которые не только проповедуют 
ту же философию солипсизма, что и члены московской математической школы, но, 
кроме того, разъезжают по Польше и не столько делают там математические доклады, 
сколько создают культурный интернационал борцов против Советского союза. 
Странно, конечно, но факт, что, скажем, академик {\bf Лузин}, который выдвинут при советской власти 
в Академию наук по кафедре философии, напечатал в Париже, в издании Бореля, 
в прошлом году книгу по теории аналитических множеств, книгу, 
как раз касающуюся такой области, которая не имеет никакого приложения к действительности,
которая толкует об абсолютно непрерывном. И вот эта книга снабжена предисловием, в котором говорится буквально,
что книга ценна не столько математически, сколько философически, поскольку она утверждает солипсизм%
\footnote{Для читателей-нематематиков, если кому попадется этот текст,
и для которого упомянутая книга Лузина - абсолютная абракадабра: никакого <<солипсизьма>> Лузин там не утверждал,
а о том, что это за <<абсолютное непрерывное>>, едва ли можно узнать у кого-либо, кроме Кольмана лично.},
как единственно научную философию. 
\end{quotation}

В статье этой достается также разнообразным технарям (не перечисляю), биологам,
медикам, физикам%
\footnote{Достается всем оптом и в розницу: {\it Два слова об одной общей проблеме, — о проблеме классификации наук...
Hа днях мне один товарищ рассказал,
что существует такая наука — гидробиология. Это — наука, трактующая о всем том живом царстве, 
которое населяет воду. Значит — о китах, о лягушках, о рыбах, о раках, амебах, о водорослях, о медузах. 
Все это изучает гидробиология. Это плохая штука. Существует институт этой гидробиологии,
там люди серьезно занимаются, на это уходят большие деньги. Какая это наука?
}};
 поименно достается левым  советским  ученым
О.~Ю.~Шмидту, М.~Л.~Левину, С.~Г.~Левиту,  И.~И.~Аголу, Б.~М.~Гессену, получают должное и М.~Планк,
Л.~де~Бройль, В.~Гейзенберг, Э.~Шредингер,
Р.~Милликен, А.~С.~Эддингтон, Э.~Персико, и
разумеется, Н.~И.~Бухарин, А.~М.~Деборин, Н.~А.~Карев, Я.~Э.~Стэн и И.~К.~Луппол. 

В 1976 году Кольман эмигрировал в Швецию, и там написал романтическо-обличительные воспоминания, 
читать которые особенно тошнотворно
на фоне его сочинений 30х годов. {\bf  Ниже
	приведена биография
 Кольмана} со слов его самого, в подстрочных ссылках приведены мои комментарии%
\footnote{Кольману во время бегства из СССР было 84 года. Конечно, у него могла испортиться память.
Но и в этом предположении   надо подходить к мемуарам критично. С другой стороны, стоит помнить 
и то, что мемуары чаще пишутся вовсе не для того, чтобы донести истину.}:

--- 
Австро-Венгерский подданный, родился в Праге, учился в 
двух тамошних университетах то ли философии, то ли математике. Участвовал в Мировой войне,
в 1915 году попал в плен.

--- В конце 17 года в качестве члена Комитета 
<<военнопленных революционных социал-демократов-ин\-тер\-на\-ци\-о\-на\-ли\-стов>> по приглашению Муралова оказывается в Москве.

---  Оказывается делегатом  Третьего съезда советов  (в то время 
орган, официально формировавший высшую власть) 10—18 (23—31 января) 1918  от организации военнопленных%
\footnote{В опубликованном в 1918г. Стенографическом отчете списка делегатов, к сожалению, нет.}.

---  Бухарин и Землячка лично дали дали ему рекомендацию в РКП(б) в апреле
1918.

--- Был избран членом Московского комитета РКП(б) (видимо, весной 1919г., причем избран на районной конференции Басманного района(?)). Встречался с Лениным на первомайской демонстрации и потом ходил к нему домой.

---- Участвовал в Гражданской войне, начальник советско-партийной
школы 5-ой армии Восточного фронта, далее с конца 1919г - начальник интернационального подотделения политотдела%
\footnote{В своих воспоминаниях Кольман много пишет о своем знакомстве с Гашеком.
В частности, описывая конец декабря 1919г. (никаких сомнений в датировке нет, это происходит в только что занятом Омске),
Кольман сообщает: {\it Я зачислил Гашека в интернациональное подотделение политотдела, руководство которым мне только что поручили.} Однако согласно разным биографиям Гашека (например, Радко Пытлик. Гашек. — Молодая Гвардия, 1977.), Гашек сам руководил этим подразделением с 5 сентября 1919г.}
по осень 1920г.

--- С 1920 года в Коминтерне,
в 1921 выехал в Германию, сидел в тюрьме в Дюссельдорфе,
в 1923 году возвращается в СССР, снова ездит в Германию.

--- 1924 заведующий Агитпропом Замоскворецкого
райкома%
\footnote{
{\it Среди агитаторов, которых приходилось направлять с выступлениями
на фабрики, заводы и учреждения района, был и Лазарь Моисеевич Каганович. 
Бывший рабочий-кожевник (закройщик), он работал тогда в
профсоюзе кожевников, помещавшемся в нашем районе. Но хотя он
слыл хорошим агитатором, посылать его с докладами приходилось с
разбором - говорил он тогда по-русски еще с заметной примесью еврейско-украинской 
местечковой <<говирки>>. Когда мне позднее пришлось
встречаться с ним и работать вместе, Каганович сам вспоминал об этом.}
\newline
А вот какой на самом деле был послужной список Кагановича:
\newline
12.1921 - 5.1922
член Туркестанского бюро ЦК РКП(б)
\newline
6.1922 - 10.3.1924
заведующий Организационно-инструкторским отделом ЦК РКП(б)
\newline
25.4.1923 - 23.5.1924
кандидат в члены ЦК РКП(б)
\newline
31.5.1924 - 29.6.1957
член ЦК РКП(б) - ВКП(б) - КПСС 
\newline
2.6.1924 - 30.4.1925
секретарь ЦК РКП(б)
\newline
2.6.1924 - 18.12.1925
член Организационного бюро ЦК РКП(б)
\newline
Очевидно, продвижению Лазаря Моисеевича по партийной лестнице 
помогло  мудрое руководство товарища Кольмана.}.

--- с конца 1924-1925 - заведующий Московским Губполитпросветом

--- ?-- 1929 член Московской комиссии партийного контроля.

--- 1925-1929 заведующий издательством «Московский рабочий»%
\footnote{Подтверждается справочником <<Вся Москва>> за 1927 -- 1928гг.
В справочнике 1929г. заведующий другой, что ничему не противоречит.}.

---  август 1929 года -- март 31 года  помощник заведующего агитпропом ЦК ВКП(б)%
\footnote{
Кольман говорит, что за это время сменилось три руководителя и называет Кнорина и Криницкого.
\newline
{\it Работать с ним [с Кнориным] было легко, мне было поручено следить
за идеологическими проблемами естественных наук и техники, за
работой соответствующих институтов и обществ. Моим напарником
– по части общественных наук – был Борис Маркович Таль.}).
\newline
Сведения об исчислении начальников
не верны: В.~Г.~Кнорин ушел с этой должности в 1927г. (т.е. до указанного периода и когда Кольман еще заведовал редакцией),
Заведующие отделом  были А.~И.~Криницкий (3.5.1928 - 19.11.1929) и  А.~И.~Стецкий (19.11.1929 - 5.1.1930)
Раздел 5.1.1930 разделился на два, Культуры и пропаганды (в который отошла наука), перешедший к Стецкому,
 и Агитации и массовых компаний, которым заведовали 
Г.~И.~Каминский (5.1 - 8.1930) и К.~И.~Николаева (8.1930 - 1933).
\newline
Кольман рассказывает, что он вместо заболевшего Криницкого ходил на прием к Молотову
и с оным повздорил. Но вместо заболевшего начальника бывают заместители,  их список вроде бы 
существует
(см. базу knowbysight) и Кольмана не включает. 
\newline
Упомянутый Б.~М.~Таль (Криштал) (1898-1938) как раз заведовал 
сектором науки (7.1929 - 1.1930) и, по идее, должен был быть начальником Кольмана, а не его напарником.
Кто далее заведовал этим подразделением, автору выяснить не удалось.
\newline
Забавно, что Кольман
о своей работе в Отделе не рассказывает ничего конкретного, кроме того, что ему был неприятен
Криницкий (надо думать, что воспоминания Кольмана о неупомянутом им А.~И.~Стецком  были более неприятны, о чем ниже п.\ref{ss:stetski}), 
и о каком-то конфликте в ИМЭЛ.
Кроме того, в это же время Кольман ездит на Урал в Среднюю Азию после конфликта с Криницким.}

--- в 1930г. входил в редакционную комиссию XVI съезда ВКП(б)%
\footnote{В Стенографическом отчете съезда \cite{16} список редакционной комиссии есть, а Кольмана там нет.
	\newline
{\it	Как ответственный работник аппарата ЦК, я поспешил на съезд,...}
	\newline
Дальше он рассказывает про общение со Сталиным...
	}

--- в  марте 1931%
\footnote{Скорее всего, это было в 1929г., см. предыдущую сноску. Собственно, 1930-31гг  - это
время, когда Кольман наблюдался на научном горизонте.}
отправился в  поездку на Урал и Среднюю Азию (судя по всему, длительную). 

--- до апреля 1931 обрел в Институте Маркса--Энгельса <<Математические рукописи>> Маркса
(о чем ниже).

--- Летом 1931 года -- член советской делегации на 2-м международном конгрессе по истории естествознания и техники%
\footnote{Да, так. Он числится также в утвержденном Политбюро списке советской делегации на Гегелевский конгрессов Берлине, ноябрь 1931}. 

---  С 1932 года  работал в институте Маркса-Энгельса-Ленина (ИМЭЛ). 

--- К 1932 году работал в Комакадемии, входил в ее
президиум и руководил Ассоциацией естествознания. 

--- В 1932 году  был назначен директором Института Красной Профессуры Естествознания. 

--- В 1932г. - руководитель советской делегации на Международном математическом конгрессе в Цюрихе.

--- В 1932 году входил в комиссию по проверке партийных документов на Украине%
\footnote{Автору данных записок неизвестно, была ли проверка партийных документов  в 1932.
Кажется, было два таких мероприятия, в 1935 и 1936гг (но я недостаточно силен в истпарте и в
этом не уверен).}.

--- В 1932 году был членом Комакадемии
президиума и руководил Ассоциацией естествознания.

---  1932  <<{\it я был назначен директором Института Красной Профессуры Естествознания}>>, 

---1934: <<{\it в то время я работал в комиссии по плану реконструкции Москвы}>>.

--- С начала 1936г завотделом науки [Отдел науки, научно-технических изобретений и открытий] Московского городского комитета партии%
\footnote{<<{\it С ликвидацией Комакадемии с самого начала 1936 года
я стал работать завотделом науки Московского городского комитета партии, 
первым секретарем которого был Каганович, но вскоре его сменил
Хрущев, бывший до того вторым секретарем.}>> [Кольман опять не помнит, кто у него был тогда начальством.
Каганович оставил эту должность в марте 1935 г., а Комакадемия была ликвидирована согласно
решению Политбюро от 7 февраля 1936г.]
<<{\it Из секретарей нашим отделом руководил Каганович, а потом Хрущев,
и поэтому я имел возможность, еженедельно докладывая им,
ближе узнать их... Я помню их обоих очень хорошо.>>} Дальше 
 тридцать пять тыщ одних курьеров:
{\it не говоря уже о том, что я наблюдал их поведение
{\bf на заседаниях  секретариата и бюро ЦК}, как и на многочисленных
совещаниях.} К сожалению, Кольман не видел там Сталина (об трех встречах в книге рассказывается отдельно),
видимо, замещая оного по его болезни.}.

--- В 1937 году (в мае или июне) стал безработным и был им год.

--- В 1938 году - инспектор 
Всесоюзного Комитета по Высшей Школе.

--- В марте 1939 года  перешел в Институт философии
АН в качестве старшего научного сотрудника, а затем зав.
отделом диалектического материализма. 

--- Одновременно преподавал логику в Московском юридическом институте,
а затем в Педагогическом институте имени Ленина.

%%%%%%%%%%%%%%%%%%%%%%%%%%%%%%%%%%%%%%%%%%%%%%%%%%%%%%%%%%%%%%%%%%%

\sm

Дальнейшая биография непотопляемого авантюриста нас не очень интересует. В Войну
занимался политработой среди пленных;
в 1945 г выехал в Чехословакию,
в 1948 году был арестован, сидел в в Москве до начала 1952 года, 
преподавал математику в Московском автомеханическом институте,
работал в Институте истории естествознания и техники АН СССР.
Снова был в Чехословакии  (1959—1963), вернулся в СССР на пенсию.
В СССР  публиковал историко-математические и пропагандистские книги, а также логические сборники.
В 1976 году эмигрировал и направил  
открытое письмо Брежневу. 

В описании биографии, как мы заметили, присутствуют явные нестыковки. 
Можно было бы предположить, что  нестыковки  в списках начальников возникла из-за неточности
указываемых дат. Это не так. Во время ухода к В.Г.Кнорина с зав.отдела ЦК
Кольман заведовал издательством, а Л.И.Каганович оставил должность секретаря МК
до того как в аппарате ЦК возник Отдел науки.

Одна дата легко устанавливается
\begin{quotation}
...  я вошел в состав Президиума [Комакадемии] и стал вместо «меньшевиствующего идеалиста» Отто Юльевича
Шмидта руководить Ассоциацией естествознания, объединявшей ее
естественнонаучные, технические и математические секции.
\end{quotation}
О том, как Шмидт стал меньшевистующим идеалистом, мы поговорим
ниже. Случилось это в январе 1931 года  (а Кольман в самом деле появляется на его
месте). Теперь об Институте Красной профессуры естествознания:
\begin{quotation}
Характерно, что теперь о них стараются молчать; во втором издании БСЭ 
даже нет упоминания об этом учреждении, игравшем тогда,
в тридцатых годах, в период ожесточенной политической борьбы,
значительную роль [следует рассказ об этом институте]. 
\end{quotation}

Надо сказать, что и в первом издании БСЭ (том 34, сдан в производство 1 апреля 1937г.)
на тему {\it <<Красной профессуры институты (ИКП)>>} говорится немного (меньше колонки
об учреждении, еще недавно столь грозном). Про Институт красной профессуры естествознания (ИКПЕ) удается
узнать лишь то, что он был организован в 1932г. и был закрыт до 1936г.
Сведения об этом институте удалось найти в современном университетском 
учебнике философии \cite{Alexeev-phil}:
\begin{quotation}
 В числе негативных последствий деятельности этих «руководителей» [Митин, Кольман]
в течение двух-трех лет были закрытие философского и естественно-научного отделения Комакадемии,
закрытие журнала «Естествознание и марксизм» и т. п.

Было распущено Общество воинствующих материалистов-фи\-лант\-ро\-пов.
Ликвидирован Тимирязевский научно-исследовательский институт, 
расформирована (в середине 1932 г.)
Ассоциация институтов естествознания при Коммунистической академии, 
а через некоторое время сама Коммунистическая академия влилась в состав Академии наук СССР  [1936]%
\footnote{Отметим, что влилась все же небольшая часть ее комсостава.}.
Был ликвидирован Институт красной профессуры естествознания.

На заседании дирекции Института философии Комакадемии 23 июня 1933 г. 
отмечалось: «После ликвидации ИКПЕ [Института красной профессуры естествознания]
и Ассоциации естествознания при Комакадемии, после закрытия журнала
«За марксистско-ленинское естествознание» (пришедшего в 1931 г. на смену журналу «Естествознание и марксизм»)
— сейчас по существу нигде нет непосредственного руководящего центра или такой организации,
которая занималась бы вплотную этими вопросами,
и в первую очередь вопросами теоретического естествознания».
\end{quotation}

Ассоциация институтов естествознания была закрыта в 1932 году (у меня, правда, источник википедического типа),
а Институт красной профессуры естествознания, как мы видим, к июню 1933 года был в прошлом.
В 1933-1935 годах он продолжает упоминаться  в журнале <<Под знаменем марксизма>> как работающий в Комакадемии. 

В 1936 году он опубликовал книгу <<Предмет и метод современной математики>>
(сдана в набор 2/IV 1936 г., подписана к печати 14/VII 1936), в которой он претендовал
на обзор и философское осмысление современной математики (ниже есть пара цитат).
Эта книга 17 марта 1937г. была атакована в газете <<Правда>>  философом А.~А.~Максимовым%
\footnote{Тоже небезынтересный персонаж, до Кольмана, впрочем, не дотягивающий.} \cite{Maximov-Kolman},
бывшим коллегой Кольмана по Комакадемии,
 в 1938 году  разгромная рецензия на нее была опубликована Гельфондом и Шнирельманом
\cite{GeSh} в <<Успехах математических наук>>.

Он претендовал на то, что он математик (а позже, что он -- историк математики).
Что касается научных статей, то 
поиск по современным базам дает у него два доклада на Математическом конгрессе в Цюрихе, 1932 году,
один по <<Математическим рукописям>> Маркса, о которых еще пойдет речь в п. \ref{ss:yanov2}а, другой - по кватернионному анализу,
публикация
выглядит как расширенный abstract, в котором, по-моему, ничего нет. Доклады на конгрессе не рецензировались.

Кроме того, есть статья  <<О разбиении круга>> в <<Мат.Сборнике>>
за 1937 год. Я пытался понять, что там написано, но безуспешно:
в статье  считается число разбиений круга наборами из $n$ кривых (видимо, с точностью до изотопии),
но понять, какие именно разбиения перечисляются, я не смог, непонятно ничего, начиная с первых строчек
(отсюда не следует, что это полная чушь, и, возможно, что это не так).  
Она поступила в редакцию 27 июля 1936 года, с двумя временн\'ыми совпадениями: аккурат после дела Лузина
(когда многие математики испытывали страх Божий),
а с другой -- вскоре окончания им работы над книгой, когда у него освободились руки для новых свершений.

\sm

{\bf\punct Статистики.%
\label{ss:hotimski}} Главным из них был Хотимский Валентин Иванович (1892-1938), человек с яркой
и едва ли исследованной биографией (и автор не берет на себя эту задачу). 
Эсер, участник Февральской революции, большевик с 1918г.
Участник Гражданской войны, партизанского движения, подпольщик. В 1918г. член Екатеринбургского
совета и вроде бы требовал расстрела всей царской семьи.

В 1927-32гг. заведовал математической секцией в Комакадемии. В 1932г. бригада под руководством Хотимского
издала учебник <<Статистика>> \cite{Hotimsky-book} (производящий впечатление и на морально 
подготовленного человека).
Сочинения Хотимского высоко оценивали истинные знатоки. Кольман \cite{Kolman-vospominaniya} 
(это воспоминания):
\begin{quotation}
Хотимский был подлинным ученым, творческой личностью, борцом, беззаветно любившим науку и добивавшимся истины,
критикуя ошибочные, по его убеждению, идеи. Он не боялся восстановить
против себя лиц, занимавших тогда командные места.

 В своих работах Хотимский выступал против ошибочной концепции некоторых советских статистиков, 
 будто в планируемом хозяйстве неприменим закон больших чисел%
\footnote{Я не смотрел, что Хотимский говорил о законе больших чисел, но его книга
по статистике полна борьбы с вредительскими буржуазными теориями. Его
соавторы в дальнейшем стали важными советскими статистиками.} . Из-за того, что приверженцы этой концепции занимали ведущие места, советская статистика
и экономические науки вообще на долгое время сильно отстали в
применении современных прогрессивных математических методов.
\end{quotation}

Яновская \cite{Za-povorot}, конец 1930г.:
\begin{quotation}
Скажем, задача создания нового учения о средних была разрешена именно благодаря тому, что с одной
стороны, применялась марксистская методология, а с другой, задача решалась на практике социалистического
строительства... Эта книга, к сожалению, у нас не известна, и к сожалению, выходит
с большим запозданием... [речь идет о предварительном варианте книги \cite{Hotimsky-book}]

Между тем здесь мы имеем впервые подлинное осуществление теории и практики
с марксистской методологией.
\end{quotation}

Приведем несколько цитат из статей ученых статистиков
(часть из них входила в число авторов учебника) из сборника \cite{Dialektika}
1931 на тему <<осуществления теории и практики>>:
\begin{quotation}
\footnotesize
Ястремский  Борис Сергеевич (1977-1962)

ГПУ верной рукой вылавливает вредителей. Что же нам, 
научным работникам и работникам практическим, следует делать? Нам
нужно итти на помощь этому верному стражу революции. Нужно
выявлять вредителей, которые ради якобы свободы мысли против
«стандартизации» мышления ведут свою якобы идейную борьбу. Этих
вредителей надо вылавливать, и мы в этом отношении должны
помочь ОГПУ.

\sm

Старовский, Владимир Никонович (1905-1975)

Приведенный мною пример показывает, как самые элементарные
статистические показатели в руках вредителей превращаются в 
чрезвычайно острое оружие в борьбе против социалистического плана.
Если в таких простых вещах вредители оставались неразоблаченными,
то еще легче было им протаскивать те свои статистические расчеты,
которые прикрывались вуалью сложнейших статистических формул.
Не случайно т. Хотимский посвятил очень много времени освещению
математико-статистического вредительства. В связи с необходимостью
борьбы с этими вредителями мы должны хорошо владеть 
математическим методом.

\sm

Боярский, Арон Яковлевич (1906--1985)

Задача не только в том, чтобы не 
ослабить своего напора в борьбе с ее остатками, в борьбе со всякими
антибольшевистскими теориями, но в том, чтобы удесятерить свою 
активность и непримиримость, вести и на теоретическом фронте войну
с классовым врагом вплоть до его уничтожения.

\sm

Цыпкин Л. [идентифицировать его не удается]

...перейти к некоторой общей оценке статистической методологии 
современных народников в отношении ее к социализму. Мне кажется, 
такой подход дает нам возможность увидеть всю необходимость пути
от научной книги, от может быть иногда блестящих рефератов весьма
парадоксального свойства до Лубянки.

\sm

Смит-Фалькнер, Мария Натановна (1878 - 1968)

Но еще лучше, если в области науки коммунистическая мысль,
марксистское общество, которое находится при Комакадемии, 
превратится тоже в бдительный орган революционной диктатуры 
пролетариата, в бдительный орган борьбы с немарксизмом. Ибо нужно 
сказать прямо, что никогда еще мы не видели такой 100-процентной 
корреляции между вредительством и антимарксистской теорией, какую
установило ГПУ. Мы можем гордиться тем, что по линии статистики
первые выступили в борьбе с этими антимарксистами, но эта гордость
пустяки по сравнению с той гордостью, которую мы когда-нибудь
будем иметь право себе приписать, если мы сделаемся ГПУ научной
мысли в области статистики и в ее применении к планированию.
\end{quotation}

Интересно (и даже удивительно), что эти люди благополучно пережили 37 год.

Сам Хотимский в 1932 или 33 году оказывается на должности начальника политотдела
моторно-тракторной станции на Дальнем востоке. Потом возвращается в Москву,
с 1935 
--
начальник Отдела статистики населения и здравоохранения Центрального управления народно-хозяйственного учёта
   Государственной плановой комиссии при СНК СССР, участник переписи 1937;
арестован в феврале 
 в 1938, расстрелян в1939ом.

\sm

{\bf \punct На Ленинградском математическом  фронте.%
\label{ss:lenmatfront}} Это название книги, изданной в 1931 году <<редакционной комиссией>> 
в составе 
Л.~А.~Лейферта, Б.~И.~Сегала, Л.~И.~Федорова:
\begin{quotation}
Какие же основные группировки сложились в самом Ленинградском университете,
 в стенах которого, в основном, и шла борьба на
ленинградском математическом фронте до 1928 года.

Почти с первых шагов работы и вплоть до 1928 г. образовались
в предметной математической комиссии три довольно прочные и устойчивые группировки: правая группа (Н.~М.~Гюнтер, В.~И.~Смирнов,
Г.~М.~Фихтенгольц и др.), левая группа (Л.~А.~Лейферт, А.~Д.~Дрозд,
А.~Р.~Кулишер и др.) и наконец менее устойчивая, чем крайние
группы, промежуточная группировка (И.~М.~Виноградов%
\footnote{Виноградов, по-видимому, 
	 никогда не был сторонником советской власти, и всегда считался человеком кадетских взглядов.}, А.~М.~Журавский и др.).

Эта промежуточная группировка несомненно идейно в области педагогических и частью организационных вопросов тяготела более
к левой группе, но довольно часто давала себя запугать правым и поддавалась на их маневры и ссылки на «академические традиции», «интересы науки» и~т.~д.
\end{quotation}

В <<и др.>> левой группы входил также составитель книги Б.~И.~Сегал. Кто такой Л.~И.~Федоров,
и имел ли он какое-либо отношение к математике, автору установить не удалось. Вот люди, которые упоминаются в
качестве <<левых>> в сборнике

 Гусева М.

Дрозд Антон Донатович 

Кулишер Александр Рувимович (1875 - после 1943)

Лейферт Леонид Абрамович (1992 -1938)

Люш В. В.

Милинский Владимир Иванович (1898-1942)

Мрочек, Вацлав Ромуальдович (1879-1937)

Рабинович Е.С

Сегал Бенцион Израилевич (1901-1971)

Шеин Ю.

Упомянутые <<правые>> памятны ныне здравствующим математикам.
Конечно, памятен Виноградов. 
Имена из левой группы  в науке  неизвестны (и, кстати, здесь виден резкий контраст с Москвой, 
где лево-радикалы, в большинстве своем, были людьми
яркими и дельными (в научном или организационном плане%
\footnote{При чтении материалов тех лет волей-неволей приходится помнить, кто был
	в Москве за большевиков,
	а кто за коммунистов. Для полноты процитируем Виленкина
	\newline
{\it
 	Проявлявшие особую  
 	активность руководители комсомольской  
 	организации математического факультета [МГУ] Л.~А.~Тумаркин, Д.~А.~Райков, А.~А.~Космодемьянский 
 	и Н.~А.~Слёзкин в дальнейшем стали  
 	серьезными учеными.}
 \newline
 (Космодемьянский и Слёзкин -- механики)
 }). Математические работы удается найти лишь у Сегала (который был не только председателем парткома Стекловки,
но и теоретико-числовиком). Дрозд, правда, был астрономом (и одно время был директором Пулковской обсерватории, уже после обсуждаемых событий,
и, кажется, по партийной линии).  Некоторые из этих деятелей, по-видимому,
стали выгодополучателями математического переворота 1931 года, но на короткое время.

Главой левой группы, по-видимому, был Лейферт%
\footnote{Лейферт на занятиях ярко описывается в воспоминаниях
	 Е.~С.~Венцель \cite{Grekova}. Но ее рассказ, скорее всего, сдержит преувеличения.}.

 %\begin{quotation}
	
 %\end{quotation}
 
 \section{Год, которого не было в летописях%
 \label{s:annals}}
 
 \COUNTERS
 
 \epigraph{
 	Как говорил автору когда-то 
 	один из его учителей член-корреспондент 
 	АН СССР В.~В.~Степанов, «на Мехмат МГУ 
 	военный коммунизм пришел с опозданием на 
 	12 лет».}
 	{Н.~Я.~Виленкин \cite{Vilenkin} }

{\bf\punct Лакуна.%
\label{ss:lakuna}} 
Александров несколько раз писал об истории матобщества и несколько раз (1942, 1946, 1957)
повторял с разными вариациями следующее:
 \begin{quotation}
  Однако при всем блеске этих научных достижений, ярко отражавшемся
на научных заседаниях Московского математического общества, в деятельности
Общества не могла не чувствоваться все возраставшая оторванность
от жизни, вызванная тем, что руководители Общества, и прежде всего
Д.~Ф.~Егоров, смотрели на математику абстрактно, никак не связывая ее
с происходившими в стране великими историческими сдвигами. {\bf  Разрыв
между традициями Общества, культивировавшимися во всей их неизменности,
приводившей иногда даже к комическим положениям, и требованиями
широкой научной общественности делался все глубже, и он разрешился
в конце концов тяжелым кризисом, выразившимся в том, что Общество целый год (1931)
фактически бездействовало. Естественно, что кризис
этот вскоре был преодолен путем широкой демократизации Общества.}
В него влилось много десятков новых членов, значительная часть
которых принадлежала к общественноактивной научной молодежи. Обновленное
Общество с большой широтой подошло к основным своим задачам,
видя их не только в научных сообщениях о новых математических достижениях
(которые, конечно, продолжали составлять основной стержень
работы Общества), но и в стремлении воздействовать на тематику научных
исследований посредством постановки докладов обзорного, проблемного,
а иногда и научно-популярного характера, посредством дискуссий, возникавших
как по поводу тех или иных докладов, так и самостоятельно. В Обществе,
совместно с университетским семинаром по истории и философии математики,
стали ставиться доклады, посвященные общим философско-методологическим и историческим проблемам математики и т. п.
 \end{quotation}

 Из сказанного ясно, что в этот год происходило что-то, о чем говорить было нельзя.
 Уже давно умерли все очевидцы, случилась Перестройка, потом 20 лет либеральной демократии.
 Многое было сказано (и автор не уверен, что все было сказано точно),
 но тема эта и доныне остается отчасти табуированной,
 о 1930 годе все же стало известно довольно много, год 1931, по-прежнему, скрыт завесой.
 
 Автор попытается собрать разные отблески и  составить из них мозаичную картинку.
 
 \sm
 
 {\bf \punct Фон.%
 \label{ss:break}} Вообще это был Год Великого Перелома и <<Наступления социализма по всему фронту>>. 
 Партия взяла курс на социальный переворот и сторонники революционных преобразований получали карт-бланш 
 на свои действия,
 их поддерживала вся мощь тогдашней государственной машины, мощь вообще-то небольшая,
 но и силы противников были невелики, и они не были организованы.
 
 В университетах операцию возглавляла Коммунистическая академия
 
 %%%%%%%%%%%%%%%%%%%%%%%%%%%%%%%%%%%%%%%%%%%%%%%%%%%%%%%%%%%%%%%%%%%%%%%%%%%%%%%
 
 \sm
 
 {\bf\punct Директор%
 	\footnote{Ильченко \cite{Ilchenko}: <<Должность директора была введена вместо должности
 	ректора университета в декабре 1929 г., когда управление вузами было перестроено на принципах единоначалия.>>}
 	Удальцов.%
 \label{ss:udaltsov}}  В 1925-1928гг ректором МГУ был историк В.~П.~Волгин, а 1925-28гг. --
 А.~Я.~Вышинский. Политические процессы, на которых <<Ягуарьевич>> выступал в качестве обвинителя, 
 тогда еще были впереди. В воспоминаниях выпускников МГУ я не встречал сильно отрицательных речей о Вышинском. 
 Например, так вспоминал о Вышинском
 М.~С.~Смиртюков, управляющий делами СовМина СССР \cite{Smirtyukov}: 
 \begin{quotation}
  А лекции по общим юридическим дисциплинам на младших курсах читал Андрей Януарьевич Вышинский,
  который был ректором университета. Естественно, тогда и подумать никто не мог, что этот умнейший преподаватель
  и блестящий лектор превратится в грозного прокурора Союза ССР. 	
 	\end{quotation}
 	
 Удальцов Иван Дмитриевич (1885-1958), член РСДРП  с 1905 года, с 1906 г. —  ячейки РСДРП(б) в МГУ, участник Первой мировой войны и Гражданской войны (зам.начальника политотдела армии),
 профессор экономического отделения МГУ с 1921г. И.~о. ректора МГУ 1 сент. 1928 – 22 марта 1929, далее ректор, освобожден от должности 30.06.1930.
 
 Из воспоминаний биолога Кузина
  Борис Сергеевича (1903-1973 или 1975):
 \begin{quotation}
 Но полное водворение новых порядков в университете произошло только в году 1930 и 1931. Не помню, когда достаточно либеральный и вполне приличный ректор Волгин был заменен Удальцовым.
 В этом не было ничего академического. Просто какой-то тип с вострой бородкой и бегающими глазками. Читал он, кажется историю социализма. Впрочем, он не столько читал, сколько управлял. 
 Сталкиваться с ним мне не приходилось. Но скоро поползли слухи о всяких его художествах с секретаршами и машинистками.
 Пользуясь ректорской властью, он оттяпал себе какие-то аппартаменты в одной из университетских квартир, оборудовал их казенными коврами и прочей мебелью для отдыха от трудов праведных.
 Университет же находился в явном запустении. Фасад старого здания на Моховой, как я его помнил, стоял в вечных лесах. Суммы на все потребности, как хозяйственные,
 так и другие, отпускались нищенские. Опять-таки не помню, как в верхах было решено поднять университет и поправить его дела. Удальцов так же незаметно вылетел из ректоров, как и стал им.
  Его место заступил столь прославившийся 
впоследствии А.~Я.~Вышинский, не оставивший при этом, кажется, своей должности 
начальника Главпрофобра, а возможно, что и какого-то высокого юридического поста.
Фасад университета был быстро отремонтирован и покрашен. Двинулись в ход и
кое-какие стройки. И во многом другом почувствовалось улучшение и порядок. Но
порядок этот стал уже совершенно новым, чуждым духу науки, бюрократическим и казарменным.

 ... Ценность записок, подобных тем, какие я сейчас пишу, была бы полнее, если бы даты описываемых событий приводились с полной точностью. Но здесь в Борке%
 \footnote{В 1953-1973гг Кузин работал зам. директора Института внутренних вод, Борки (Ярославская область, Рыбинское водохранилище, директором был И.~Д.~Папанин).
 Ошибка с порядком ректоров и Вышинским небезынтересна, по-видимому, <<столь прославившегося 
 впоследствии>> Ягуарьевича нельзя было поставить в эпоху до больших безобразий.}, я не имею возможности проверять их, а память сильно обманывает... 	
 	\end{quotation}
 	
 Из <<Летописи Московского университета>> \cite{Ilchenko}:
 \begin{quotation}
 	Декабрь, 27
 	
 	Издан приказ №1 ректора И.~Д.~Удальцова о перестройке управления университета на основе введения единоначалия:
 	
 	 «1. В связи с введением “единоначалия” как принципа управления университетом правление университета реорганизуется в совещательный орган при ректоре. Все распоряжения по университету будут впредь отдаваться от имени ректора или его помощников (проректоров) в пределах полномочий последних.
 	
 	2. В осуществление принципа единоначалия деканаты реорганизуются впредь в совещательные органы при деканах. Вся ответственность за состояние и работу факультета целиком и полностью возлагается единолично на декана».	
 \end{quotation}
 
 Из книги \cite{Petrovsky3}:
 \begin{quotation}
  13 февраля 1930 г. совещание проректоров и деканов 1-ro МГУ по докладу
ректора И.~Д.~Удальцова принимает постановление, в котором говорится: «В связи
с выделением из 1 МГУ факультетов физико-математического, медицинского
и химического, факультетов истории и философии, факультета литературы и искусства, 
факультета совправа, а также преобразованием их в самостоятельные
вузы, считать необходимым разработать (в основном) в 5-дневный срок проект
финансового, административного и хозяйственного расчленения 1 МГУ, поручив это тройке в составе И.~Д.~Удальцова, И.~Г.~Руфанова и В.~Ф.~Андреева»
 \end{quotation}

 В журнале <<Красное студенчество>> в январе 1930г. была опубликована статья Удальцова (цитируется по \cite{Petrovsky3}):
 \begin{quotation}
 175-ЛЕТНИЙ СТАРЕЦ
 	
 	В этом году исполнится 175 лет существования 1 Московского государственного университета. Не случайно, что на пороге этого юбилея встал вопрос о дальнейшем
 	существовании университета. На протяжении почти двух столетий рос
 	и расширялся университет, «отражая в себе» рост и продвижение науки. Из единого
 	небольшого стержня, постепенно развиваясь, он разбился на ряд факультетов, объединяющих множество научных дисциплин. Университет и сейчас сохраняет первоначальную организационную форму, 
 	подсказанную еще в момент организации университета эмбриональным состоянием науки в России,
 	университет сохраняет структуру, позволявшую в едином учреждении охватить почти все дисциплины в области
 	учебной и на\-уч\-но-ис\-сле\-до\-ва\-тель\-ской работы и тогда оправдавшую название университета.
 	
 	Сейчас эта форма, связанная с русским средневековьем, уже отжила свой век.
 	Она вступает в конфликт с жизнью. Уже с момента Октября, когда от науки требуется
 	быстрое реагирование на требования жизни, когда перед каждой отраслью науки
 	ставятся совершенно конкретные задачи - структура университета показала себя
 	неуклюжей и недостаточно эластичной. Невозможно ввести точную плановость,
 	полноценно обслужить и оборудовать хозяйственно и экономически все отделения,
 	нет возможности дать все условия для расцвета каждой специальности.
 	Университет, разросшийся в чудовищно громоздкое учреждение, объединяющее в себе шесть разнородных факультетов, не может с нужной быстротой разрешать выдвигаемые жизнью вопросы. Университет не в состоянии приспособиться
 	к условиям и темпам социалистического строительства. Университет отстает от жизни
 	и в сущности он оказывается не двигателем, а тормозом развития науки. Ноябрьский
 	пленум ЦК, поставивший вопрос о темпах, о новых методах работы, о теснейшей
 	связи с производственными заданиями - произнес свой приговор над прежней структурой университета.
 	
 	Мы считаем необходимым расчленение университета на его вполне самостоятельные
 	составные части: Медицинский институт (он должен возникнуть из слияния
 	медфаков 1 и 2 МГУ); Институт общественных наук (в него должны войти
 	факультеты: совправо, ис\-то\-ри\-ко-фи\-ло\-соф\-ский и факультет литературы и искусства);
 	Институт физико-химических наук с факультетами- химфак и физмат (он разобьется
 	на ряд самостоятельных факультетов).
 	
 	Пора старику-университету на 175-летнем юбилее своей жизни - на покой.
 	Около девяти тысяч студентов, две тысячи рабочих и служащих, полторы тысячи
 	научных работников, 14 гектаров крыш на университетских зданиях и одиннадцать
 	километров тротуаров вокруг них - сами настаивают на своем разделении ...
 \end{quotation}

 Эта программа начала быстро претворяться в жизнь \cite{Ilchenko}:
 \begin{quotation}
 	Май, 20
 	Приказ ректора И.~Д.~Удальцова: «В целях скорейшего завершения работ по выделению 
 	передаваемых в другие наркоматы факультетов и отделений университета объявить работу по выделению ударной».
 \end{quotation}
 
 Программа по разделению Университета была во многом выполнена. Гуманитарные факультеты были отделены в МИФЛИ
 (Московский институт философии, литературы и истории). Были отделены химики.
  Впрочем к концу его ректорства было полуразгромлено и  то, что еще оставалось...
 
 Химики вернулись  в МГУ  в 1932-33гг.
 В 1934 году был восстановлен истфак, в ноябре 1941 ИФЛИ был присоединен к МГУ,
 и тогда в МГУ появились филологический, философский и экономический факультеты
 (и вместе с последним МГУ получил того же Удальцова в качестве декана экономфака,
 а также экономфак, возглавляемый Удальцовым.).

 \sm
 
 %%%%%%%%%%%%%%%%%%%%%%%%%%%%%%%%%%%%%%%%%%%%%%%%%%%%%%%%%%%%%%%%%%%%%%%%%%%%%%%
 
 {\bf\punct Падение Егорова.%
 \label{ss:egorov-fall}} Сведения о Егорове  искал в 
 архиве МГУ Ч.~Форд, он нашел много интересных документов \cite{Ford-arhiv}, \cite{Ford3}.
 Однако и дата  ухода с директорского поста, и обстоятельства этого ухода, по-видимому, остаются неизвестными.
 Последняя известная подпись Егорова в качестве директора НИИММ датирована 24 октября 1929г.
 Через два месяца, 21 декабря 1929г. состоялось общее собрание аспирантов МГУ (о котором сохранились документы).
 По-видимому, в тот момент
 Егоров уже ушел (или был отправлен) в отставку. В любом случае 24.12.1929 
 в качестве директора уже присутствует Шмидт.
 % В январе директором уже был утвержден Шмидт.
 Форд предполагает, что Егоров ушел под давлением революционной молодежи (данные о такой активности,
 прежде всего Хворостина и Райкова имеются). Шмидт вскользь говорил об этом
 моменте на Комиссии  1936г.
 \begin{quotation}
 	ШМИДТ.
 Егоров, конечно, по
 	сравнению с теми [Лахтин, Бугаев, \dots], был фигурой прогрессивной. С другой стороны, по своему мировоззрению, по всему своему отношению и к философским, и к религиозным и прочим вопросам, по своему крайнему консерватизму даже в мелочах университетской жизни, —
 	он был самым настоящим душителем, именно душителем. Он умел воспитывать учеников. И в советское время это долго терпелось. Нас было мало. Талантливая молодежь
 	только подрастала, и это некоторое время продолжалось, а затем это закончилось взрывом, буквально взрывом, ибо нельзя было больше терпеть. И этот взрыв имел драматические последствия для Егорова...
 \end{quotation}
 
 Однако, стоит иметь в виду, что на Егорова было кому давить и сверху (с ди-ректората),
 и изнутри Физмата. Отсутствие данных о таком давлении заставляет вспомнить гипотезу об уничтожении документов.
 
 Еще пара штрихов о тогдашней обстановке, Статья Ч. Форда, \cite{Ford3}:
 \begin{quotation}
  Райков и Хворостин были среди восьми претендентов, поступавших в 
аспирантуру Института математики и механики осенью 1929г. Шесть
из них были рекомендованы партийной организацией 
физико-ма\-те\-ма\-ти\-че\-ско\-го факультета. Окончательное решение 
Отборочной комиссии было получено только в середине декабря. Из
шести абитуриентов, рекомендованных партийным бюро, трое
были членами партии, остальные трое — комсомольцами. Среди
членов партии были Райков и Хворостин; третьим был Кива Израилевич Защинский, состоявший в партии с 1920г. На 
общем собрании аспирантов 21 декабря 1929 г. он был одним из 
четырех, избранных в президиум собрания.
 \end{quotation}
 
 Стоит иметь в виду, что комсомол тогда был еще не очень массовой организацией, а партийных на Физмате было
 совсем мало. 
А вот из журнальной статьи начала 1930г. \cite{Dolin}:
 \begin{quotation}
 	21 декабря состоялось собрание молодых научных кадров - аспирантов I Московского университета...
 	В нем насчитывается сейчас свыше 250 аспирантов. 
 	В повестке дня стояли отчет временного бюро, утверждение плана работ общевузовского бюро на 1929/30 г.
 	и выборы бюро\dots Свое вступительное слово председатель временного аспирантского Бюро тов. Мицкевич посвятил общей политической установке в деле подготовки новых кадров из рабочих и крестьян и задачам,.
 	стоящим перед ними в университете. Затем  он остановился на характеристике профессорского, 
 	научного и аспирантского состава университета, в связи с прошедшей отчетной перевыборной кампанией профессорского персонала. 
 	Так выяснилось, что некоторые группы аспирантов и часть товарищей из руководства старого,
 	бюро в этом вопросе занимали неверную позицию. Кампания сразу разделила весь аспирантский коллектив резко на группы.
 	Одна из них открыто выступала в защиту реакционной профессуры и тем самым проявляла себя как антисоветскую
 	и реакционную группу. Вторая исподтишка сидела в «нейтральном болоте», была аполитична, и только рабоче-крестьянская,
 	советская и партийная группа аспирантов занимала правильную политическую линию в этом вопросе. 
 	
 	Такое положение лишний раз подчеркивает неблагополучие в аспирантском составе; несомненно,
 	здесь имеет место засорение его чуждыми элементами, выдвинутыми в аспиранты реакционной профессурой:
 	Особенно резко дифференциация среди аспирантов выявилась в математическом институте, 
 	когда аспиранты с выдвиженцами из студентов дали хорошую отповедь реакционному профессору Егорову,
 	критикуя его косность, оторванность, инертность и аполитичность в реформе и перестройке всей педагогической,
 	научно-исследовательской работы и методологии. 
 	\end{quotation}
 	
Так или иначе, Егоров покинул пост директора института. Президентом Московского математического
общества он остался.

\sm

 %%%%%%%%%%%%%%%%%%%%%%%%%%%%%%%%%%%%%%%%%%%%%%%%%%%%%%%%%%%%%%%%%%%%%%%%%%%%%%%
 
 {\bf\punct Математика в МГУ весной 1930г.%
 \label{ss:1930-spring}}
Теперь эпизод, который по-видимому, не документирован, но,
очевидно произвел впечатление на присутствовавших.
О нем известно из воспоминаний Бескина 1993г.
\cite{Besk29} и погромных статей статистика
Старовского \cite{Starovski} и  Кольмана \cite{Kolm1} 1931г.
\begin{quotation}
(Бескин)

После вступительной
программной речи О.~Ю.~Шмидта  выступил
Д.~Ф.~Егоров и сказал:'<<
Я всегда считал, что люди, не обладающие стандартным 
мировоззрением, могут тем не менее добросовестно трудиться и
приносить пользу науке. Я очень огорчен, что вы думаете иначе.>>

(Старовский)

Мне недавно пришлось слушать на заседании совета Института
математики и механики речь проф.~Егорова, тогда еще не 
разоблаченного вредителя. Он выступил со своего рода программной речью и
так горячо, со слезой даже в голосе, сказал: «Что вы там толкуете
о вредительстве... худших вредителей, чем вы, товарищи, нет, ибо
вы своей пропагандой марксизма стандартизируете мышление»...

 (Кольман)

 .... но характерно — если взять лишь события последнего месяца -- что признанного вождя
реакционной московской математической школы, еще в прошлом
году директора математического института, состоявшего церковным старостой,
но не желавшего быть членом профсоюза, проф. Егорова московское
математическое общество упорно не желало исключить из своего состава.
Когда же Егоров заявил, что «не что-либо другое, а навязывание стандартного
мировоззрения ученым, является подлинным вредительством», докладчик-коммунист
не только сам не дал ему отпора, но в заключительном слове
отвел предложение сделать из выступления Егорова организационные выводы,
объяснив все «недоразумением». Такова политика некоторых коммунистов,
проводимая ими в реакционнейшей профессорской среде, в среде хранителей
традиций Цингера, Бугаева, Некрасова, разрабатывавших теорию вероятностей,
науку о числе и анализ для доказательства незыблемости «православия,
самодержавия и народности», для подкрепления философии Лопатина в среде
тех людей, которые вполне последовательно на недавнем своем съезде отказывались 
послать приветствие XVI партийному съезду.
[Далее идет осуждение подобных партийных товарищей]
\end{quotation} 

 Отрывок из Кольмана посвящен не столько Егорову (который, судя по тексту, уже был
 арестован, а также
  исключен из Мат.Общества),
а об осуждении <<докладчика-коммуниста>>, т.е.  Шмидта
(Шмидт был обвинен во время <<событий последнего месяца>>,
о которых еще пойдет речь, о <<приветствии XVI партийному съезду>> - тоже ниже).
Впрочем, кольмановский текст содержит вполне интересные отсветы об обоих этих деятелях.

\sm

А.~Ф.~Лапко и Люстерник публиковали список докладов в НИИММ 1929-1930гг.

{\small
1929\dots	

14.12. Защита С. В. Бахвалова «О совместном изгибании двух связанных поверхностей». 

24.12. О назначении директором НИИММ О. Ю. Шмидта; А. И. Некрасов «О вихре в вязкой
жидкости».

1930г. 4.1. Сообщение П. С. Александрова о последней работе Л. С. Понтрягина; В. И. Гливенко «О неявных функциях».

 14.1. Н. В. Мелентьев «О методе приближенных вычислений»; 
 
 26.1. Е. Л. Николаи «О деятельности Ленинградского механического общества».
 
  3.2. В. В. Степанов «О составлении пятилетнего
и генерального плана потребности СССР в математиках».

 4.2. С. Н. Бернштейн
«Приглашение участвовать в математическом съезде в Харькове». 

14.2. Л. С. Понтрягин «О непрерывных алгебраических телах». 

24.2. П. С. Александров
«Топология замкнутых множеств и теория размерности»; Г. Н. Николадзе «О непрерывных системах геометрических фигур». 

14.3. И. В. Хлодовский «Разложение
дифференцируемых функций в ряды полиномов». 

24.3. Защита А. Н. Колмогорова
«Исследование о понятии интеграла». 

24.4. А. Н. Колмогоров «Аналитические
методы в теории вероятностей». 

4.5. А. Я. Xинчин «О рабочей книге по математике»;
О курсах ускоренной подготовки преподавателей вузов и втузов. 

20.5. А. О. Гельфонд «Исследования по теории Адамара».

 6.6. Защита Г. Б. Гуревича «О некоторых интегральных задачах тензорного анализа». 
 
 18.6. Выдвижение на премию
работ А. О. Гельфонда, А. Н. Колмогорова, Л.~С.~Понтрягина, Л.~Г.~Шнирельмана.

3.7. В. Бляшке (Германия) «Лиувиллевские элементы дуги». 
}

\sm

Мы видим при Шмидте нормальные научные доклады и нормальные защиты,
исключение составляет доклад Степанова о планировании математики,
тема эта, очевидно, появилась под давлением сверху.
Надо думать, под давлением сверху или под давлением обстоятельств
был сделан и доклад Хинчина, но это все же в рамках естественной тематики. 

\sm 

Александров и О.~Н.~Головин \cite{Alex-MMO2} в 1957 году публиковали списки заседаний Мат.Общества.
Вот список за 1930г:

\sm
 
% \begin{quotation}
{\small
\hangindent=1cm\noindent
11 января 1930 г. 1. И. И. Привалов «{\it Об областях точек, изображающих
общие значения пары аналитических функций}».

2. А. Я. Xинчин «{\it Об одной вариационной задаче волновой механики}».

\hangindent=1cm\noindent
1 февраля 1930 г. 1. Е. М. Ливенсон «{\it О os-функциях Hausdorff'а}».

2. Л. В. Канторович «{\it О проективных множествах}».

3. Ц. Г. Бурстин «{\it Об изгибании гиперповерхностей}».

\hangindent=1cm\noindent
22 февраля 1930 г. 1. К. Н. Шапошников «{\it Механика квант}».

\hangindent=1cm
2. Ю. А. Рожанская «{\it О непрерывном разбиении поверхности на одномерные
континуумы}».

\hangindent=1cm\noindent
1 апреля 1930 г. 1. П. С. Александров «{\it О разбиении пространства и общей
теории размерности}».

\hangindent=1cm
2. Е. Е. Слуцкий «{\it Основания теории периодограмм для случая зависимых
испытаний}».

\hangindent=1cm\noindent
11 мая 1930 г. 1. К. Ф. Огородников «{\it О первой теореме теории ошибок
наблюдений}».

2. Д. Ф. Егоров «{\it О разложении резольвенты на главные части}».

3. Л. Г. Шнирельман «{\it Об аддитивных задачах в теории чисел}».

\hangindent=1cm\noindent
10 октября 1930 г. 1. Л. С. Лейбензон «{\it Приложение теории гидравлического
удара к исследованию работы глубоких насосов}».

2. М. А. Лаврентьев «{\it Об одной минимальной задаче в конформном отображении»}. 
}
% \end{quotation}

\sm

Мы видим, что весной 1930 года Общество заседало в обычном режиме. Арест Егорова, по-видимому, произошел
в ночь с 9 на 10 октября 1930г.

\sm
 
 %%%%%%%%%%%%%%%%%%%%%%%%%%%%%%%%%%%%%%%%%%%%%%%%%%%%%%%%%%%%%%%%%%%%%%%%%%%%%%%
 
  {\bf \punct Математика должна быть прикладной.%
  \label{ss:applied}}
 Атака на чистую математику без подробностей мелькает в многочисленных 
 источниках. Вот связный абзац из воспоминаний Н.~М.~Бескина:
 \begin{quotation}
 На математическом отделении  [МГУ]  был введен крайний прикладной
 уклон. На младших курсах преподавание стало очень элементарным. 
 Главными предметами считались приближенные вычисления,
и  Учение логарифмической линейки, номография, начертательная
геометрия, черчение. На старших курсах культивировались 
прикладные методы теории функции комплексного переменного, 
приближенное решение дифференциальных уравнения. Теория 
функций действительного переменного, теория чисел были в загоне.
Во главе математического отделения стоял Б.~С.~Дворкин,
посредственный ученый, позднее работавший в педагогическом
институте... Была
подведена теоретическая база под это утилитарное направление.
\end{quotation}

Не совсем ясно, о каком в точности периоде времени рассказывал Бескин.
Очевидно, что это нашествие было одной из проблем, с которыми столкнулся Лузин,
вернувшийся из поездки во Францию осенью 1930г. По разным отсветам ясно, что
под вопрос ставилось вообще существование чистой математики (и, кстати, одним из борцов с ней был
Райков), на какое-то время была отменена аспирантура по этой специальности.
Отчасти защите чистой математики была посвящена речь Шмидта на Первом всесоюзном
математическом съезде (см. следующий пункт).

Искались и другие решения. В статье И.~Зайденвара
<<Октябрь в математическом обществе и в институте
механики и математики>> 1930г. (весьма соответствующей духу момента
в худшем смысле этого слова), среди прочего,
сообщается следующее:
\begin{quotation}
 В жизни Института [математики и механики] также наметился резкий перелом в
сторону, активного участия в производственных задачах дня. Институт
связался уже с Государственным Институтом Сооружений, ведутся переговоры о связи с другими институтами 
(Тепло­технический Институт, ЦАГИ, ВЭИ%
\footnote{Кроме ЦАГИ, который тогда сыграл важнейшую роль,
названы
\newline
Центральный научно-исследовательский институт промышленных сооружений  (строительная механика)
\newline
Всесоюзный научно-исследовательский теплотехнический институт (энергетика и теплотехника)
\newline
Всесоюзный электротехнический институт (электротехника)}
),
созданы плановые группы
 по проблеме кадров, по связи с промышленностью, учреждаются 
 подготовительные курсы для подготовки в аспирантуру
института людей из рабочего класса;  в общем <<чистая наука>>
пошла на службу сегодняшнему творческому дню.
\end{quotation}

Аспирантура для людей, не прослушавших нормальных курсов, по тогдашним (и по более поздним)
нравам была делом необычным.
Но поиск прямого сотрудничества с ведущими технологическими исследовательскими институтами
был правильным  ходом со стороны НИИММ.

Еще одно известие \cite{Glagolev}:
\begin{quotation}
 Осенью 1930 г. Нил Александрович [Глаголев%
 \footnote{Глаголев (1888—1945) занимался в основном (теоретической) номографией. Это
 было такое обобщение логарифмической линейки. На
 бумаге изображается кривая или набор кривых. Дальше посредством взгляда на картинку
 или простой операцией, вроде приложения линейки, мы приближенно решаем конкретное уравнение или вычисляем
 значение какой-нибудь функции. В 20-30е годы <<номограммы>> были полезным средством для массовых типовых 
 приближенных вычислений,
 конкурировавшим с механическими и электрическими арифмометрами (а задача изобретения номограмм была
 весьма забавна и содержательна), с развитием электронной вычислительной техники номографические
 вычисления, так же как и табличные вычисления (и заодно с ними старинные <<счёты>>), стали архаичными. 
 Глаголев работал в НИИММ со времен Егорова (и, по крайней мере лично, не относился к <<сонму>>, о котором чуть выше
 была цитата из Бескина), он также известен
 составлением и редактированием учебников геометрии для средней школы. Видимо, был учеником А.~К.~Власова.}] был 
командирован Институтом математики Московского 
университета в Германию [Берлин и Гёттинтен] для ознакомления с 
состоянием там прикладной математики.
\end{quotation}

Так или иначе, это было время жестокого наката против чистой математики, и самые разные 
люди - левые и правые -
занимались как прямой защитой математики, а также поисков выходов на приложения - на приложения
для реального дела, а не правильного писания заявок на гранты (математики тогда в свою науку верили).
И в том, и другом они в итоге преуспеют... 

\sm

{\bf\punct Математический съезд. Речь Шмидта.%
\label{ss:sezd-schmidt}}
В Харькове
 24—29
июня 1930 г.
проходил Первый Всесоюзный съезд математиков (председатель 
оргкомитета - Бернштейн).
Открывался он  программным выступлением 
Шмидта «Роль математики в строительстве социализма»,
а закрывался выступлением Кагана «О теоретической и прикладной математике в СССР».
Текст выступления Кагана мне не удалось найти (по-видимому его как-то отражает резолюция съезда),
выступление Шмидта  было в разных отношениях интересно, мы приведем из него
большие выдержки.

\sm

а) {\sc Рабочие места для математиков.}
\begin{quotation}
	\small
  Я скажу хотя бы то, что сводный бюджет
научных учреждений всего Союза на будущий год составляет 450 млн. рублей. Я смею утверждать, что ни в одной стране, за исключением США,
такого бюджета научных учреждений нет. 450 млн. рублей-это бюджет
абсолютно больший, чем в любой стране, кроме США, а относительно
к народному доходу это больше, чем в какой бы то ни было стране было и
есть. Я обращаю ваше внимание еще на то, что мы имеем сейчас в настоящий
момент шесть тысяч аспирантов - молодых ученых, получающих
от государства стипендию не менее 100 р. в месяц, а в большинстве случаев
и больше...

Мало того, при развитии высшего образования в связи с подготовкой
необходимых для страны кадров специалистов приходится проектировать
еще большие темпы подготовки научных работников... 

... Как
вы знаете, пятилетка осуществляется успешно. Наиболее узким местом
ее осуществления является недостаток специалистов. Это узкое место
будет преодолено в течение ближайших четырех, может быть, пяти лет.
За это время вольются в народное хозяйство больше миллиона новых
специалистов. Соответственно с этим происходит расширение высшего
образования. Как раз в день моего отъезда сюда Совнарком утвердил
цифру приема в вузы на будущий год в 120 тыс. человек, против пятидесяти
с небольшим тысяч приема прошлого года и менее сорока тысяч позапрошлого.

 Из этих 120 тыс.
не менее чем 100 тысячам мы, математики, будем преподавать высшую математику
в технических, сельскохозяйственных, экономических
и педагогических вузах....
 Для вас ясно по опыту каждого из нас, что преподавание и научная работа тесно переплетаются и не только в том смысле,
что большое количество часов преподавания мешает научной работе. Мы
знаем, что в процессе преподавания мы подымаем новые пласты научной
молодежи и часто сталкиваемся с задачами, которые нас углубленно заинтересовывают. 

....  у нас именно инженеров более всего не 
хватает и, как прямое следствие, на рынке преподавателей высшей школы
более всего похватают математиков. Молодой человек, который занимается
нашей наукой, имеет все шансы стать профессором в 25 лет. Такая большая нужда. 
\end{quotation}

По отношению к приводимым цифрам стоит быть осторожным.
Слово <<наука>> в СССР включал <<high technologies>>.
В других странах, где технологии разрабатывались,
 эти расходы в большей степени, чем в СССР,
проходили по линии разрабатывавших их промышленных фирм.
 Прием в вузы то увеличивался, то уменьшался, отсев, кажется, был очень большим \cite{Lapko}.
Но то, что количество рабочих мест для математиков в среднем быстро  увеличивалось,
не вызывает сомнений%
\footnote{Стоит иметь в виду, что в СССР математики, работавшие в технологических
организациях, имели возможность заниматься чистой наукой и публиковать свои
результаты. На памяти автора для этого нужно было при публикации подписывать
 определенную форму о неразглашении тайн (<<Акты экспертизы>>) и собирать
под ней положенный набор подписей. Бумага раздражала, но альтернатива
этому (как мы сейчас знаем) -- фактический запрет на публикации тем, кто работает
в неакадемических организациях.}. При просмотре разных биографий в самом деле оказывается,
что разные люди по окончанию физматской или мехматской аспирантуры сразу становились, например,
заведующими кафедрами.

Следует также иметь в виду, что стране в самом деле были нужны специалисты,
 их в самом нужно было учить (что в корне отличается от современной ситуации
раздутого имитационного образования).

\sm

{\sc b) О чистой и прикладной математике.}
\begin{quotation}
	\small
 Математика из всех наук имеет
наибольший соблазн считать себя наукой надклассовой и вообще наукой
<<безгрешной>>. Тут есть и тот соблазн, что приложение математики далеко
не адекватно всему математическому творчеству, что значительная часть
математиков далека от приложений. Тут есть тот соблазн, что математические истины вытекают, 
казалось бы, из совершенно особых свойств
нашего ума или создаются нашим умом независимым образом, и поэтому
математика может и должна держаться в стороне от классовой борьбы
и от социалистического строительства. Что касается приложения математики, то все вероятно знают заявление одного из крупнейших
ученых, что ему нравится теория чисел потому, что это та часть математики, которая еще не запятнана приложением. Несомненно, эта точка
зрения довольно популярна. Есть тенденция превратить математику
в особый мир, который не обязан быть ни в каком соответствии с миром
действительным....

\dots мы
хотим математику в целом видеть среди коренных научных сил, участвующих в нашем строительстве. На деле, разумеется, математика никогда
не была так далека от приложения, это больше фикция тех математиков,
которые заняты наиболее абстрактными ее ветвями и, занимаясь этой
абстракцией, потеряли сами чутье, откуда эта абстракция берется. Постепенно в традиции теряются первоначальные корни абстракции, первоначальным толчком которой было действительное явление. На протяжении истории математики это вскрывается довольно ярко. Ведь математика XVIII в. не говорила о кризисе чистой математики, и бессмысленно было бы тогда об этом говорить при очень тесной связи с реальным миром, когда математики радовались тому, что находимые ими абстрактные методы дают блестящий результат и позволяют охватить, как
тогда думали, все мироздание формулой чисел. И тогда, и сейчас - если
говорить откровенно не о философской тенденции, заложенной в реакционном крыле математиков, а о действительном положении - нужно
сказать, что современная математика вовсе не так далека от приложений.
Позвольте привести в доказательство два примера: во-первых, те явления
которые сейчас происходят в физике - углубление и развитие 
теории относительности, с одной стороны, комплекс вопросов, связанных
с квантовой механикой и волновой механикой, с другой стороны - показывают, что самые абстрактные в философском смысле математики,
занимающие на философском фронте давление от нас позиции, в то же
время участвуют в разрешении этих задач, выдвигаемых физикой, как
бы они ни желали замкнуться в своем блестящем одиночестве. Независимость 
математики является фикцией. На деле математика умерла бы
давно, если: бы она пребывала в этом блестящем одиночестве! На деле
она и теперь питается другими науками...

 Все мы должны вооружить технику новым оружием, должны
мобилизоваться для разрешения новых задач, которые к этому моменту
встанут. Это в отношении техники, но не техникой исчерпываются наши
задачи и не тем еще, что мы входим нашей наукой интегральной частью
в общее образование, и не тем только, что будущие педагоги и статистики
будут впитывать ее как часть образования. Вы помните лозунг: <<социализм -- это учет>>. Нечего говорить, что для учета нужно уметь считать
и вовсе не элементарно. Учет в условиях страны социализма это есть
решение большой сложной задачи, одно только правильное размещение 
промышленности приводит к очень сложным задачам подсчета%
\footnote{По нынешнему - применение математики в экономике.}. Нам
нужно однако не только это. В стране, где строится социализм, где нужно
уметь считать, нужно, чтобы это умение математически формулировать
стоящие перед каждым задачи, уменье подойти во всеоружии науки к каждой конкретной задаче, умение руководить наиболее экономно и точно, -
чтобы это уменье
было всеобщим достоянием. Нам необходимо трудиться
над тем, чтобы общая математическая культура у нас была выше чем
у других.

Совершенно неверно думать, что советская власть и руководящая
наша партия ждут от науки только непосредственного практического
приложения. Это не верно. Это опровергается хотя бы нашим съездом,
который создан при поддержке правительства...

{\bf Математика представляет собой одну из наиболее разветвленных
наук, настолько разветвленных, что когда кто-нибудь делает математический доклад, то очень небольшое число присутствующих может следить за ним компетентно. Это таит в себе опасность общего разброда, опасность отрыва от питающих корней, опасность заблуждений в дебрях леса,
откуда не выбраться, и где напрасно расточаются силы и гибнут таланты.}  
 Успех в том, чтобы рядом со специальными
работами была и работа синтетическая, чтобы основа математики — методология -- была в центре нашего внимания, чтобы каждый отдельный математик,
удаляясь в лес специальных изысканий, видел спасательные огни. 
\end{quotation}

Так или иначе, говоря о необходимости приложений, Шмидт проводит защиту чистой математики...

\sm

{\sc c) Классовая борьба и диалектика.}
\begin{quotation}
	\small
Неверно также утверждение о том, что современная математика
за рубежом или здесь есть математика внеклассовая. Этого нет. Разумеется, 
конечно, есть отдельные математики, которые занимаются вопросами о таких-то пространствах или о таких-то уравнениях и при
этом субъективно, в отношении науки по крайней мере, далеки от того,
чтобы с этим связывать какой-либо классовый интерес, но мы говорим
об огромных явлениях, в которых участвуют сотни людей, и вправе
брать общие средние и общие выводы, независимо от субъективных позиций
 тех или иных лиц. Как же используется математика? Товарищи статистики 
 например знают очень хорошо, что делают с нашей наукой. Всем
известны формулы Пирсона,- математическое содержание они имеют 
небольшое, это также всем известно, нужно это только откровенно сказать.
Однако использование этих формул принимается как нечто глубоко научное
	и обоснованное, и если какое-либо социальное явление на каком-то
	отрезке времени располагается по кривой, которая похожа на кривую
		Пирсона, то постулируется, что это явление происходит согласно такой-то
		кривой Пирсона, и делается предсказание о том, что и в дальнейшем будет
			так. Американская литература заполнена такими предсказаниями,
			и в частности американская экономическая литература, применяя такого
			рода кривые, предсказывала, как известно, дальнейшее развитие и
			процветание в Америке, а вот там произошел кризис, которого формула
			Пирсона как раз не предусматривала...
			
Разве случайно, что 20 лет тому назад руководящие
научные журналы не поместили бы статьи, касающейся бога, религии
и веры, тогда как сейчас вы можете ее найти в научных журналах?
Нет, это не случайность, это отражает изменившееся положение буржуазии,
которая сейчас в 30-х годах гораздо ближе к своей конечной катастрофе,
чем в начале века, которая видела своими очами революцию,
которая видит победу революции в одной огромной стране, дрожит за
свои классовые позиции и соответствующим образом свою идеологию заостряет в реакционном духе....

Никогда не удастся
нам, математикам, добиться общих больших результатов, если мы не
будем пользоваться этим всеобщим методом, всеобщим мировоззрением.
Наша наука должна догнать те темпы, которыми идет строительство социализма. 
Она развивается у нас не такими темпами, а для того, чтобы
догнать, нужно прежде всего иметь общетеоретическое вооружение.
Каков этот метод? Это метод диалектического материализма. Мы можем
смело сказать - метод современной науки, смело проводя знак равенства
между методом современной науки и диалектическим материализмом,
ибо главнейшие успехи, которые действительно мир переделывают, мы
можем отнести за счет диалектического материализма, метод которого
пока применяется в сфере наук общественных и осуществляется там
под руководством материалистической диалектики Маркса, Энгельса и
Ленина. Теперь очередь наступает за другими науками. Мы убеждены,
что и в других науках настоящий, подлинный расцвет наступит именно
тогда и именно в связи с тем, что и их руководящим общим методом станет
диалектический метод Маркса, Ленина и Энгельса. Только при наличии
этого метода проблема связи абстракции с действительностью для математики
перестанет быть чем-то зыбким, а получит реальное подробное
разъяснение. Только в свете этой теории станет ясно, какое же значение
имеет математика в среде других наук. Я не имею основания излагать
здесь основы марксисткой диалектики, считая, что они в значительной
мере известны. Марксистская диалектика есть прежде всего более высокая
ступень мышления, нежели формальная логика. Формальная логика
есть орудие хорошее, но орудие более грубое, годное лишь для более
упрощенных случаев. Для более же тонких случаев, где выступают внутренние
противоречия, метод формальной логики пасует, а диалектическая
логика дает орудие для разъяснения этих явлений. Это есть более
высокая ступень мышления, и нужно прямо признать, что успехи во всей
нашей научной деятельности будут прямо перевернуты, когда к ним будет
приложен метод диалектического материализма. Это не значит, что мы,
зная переход количества в качество, должны всюду-кстати и некстати применять
эту формулу. Дело не в этом. Дело в глубокой проработке
диалектики, в глубокой проработке основ нашей науки под этим углом
зрения. Это очень большая и сложная работа, это работа, для которой
нужно много участников, в которой каждый из нас должен принять участие,
ибо мимо этой задачи нельзя пройти, жизнь выдвигает ее...

 Мы рассматриваем математику
как весьма сильное орудие. Мы считаем, что это орудие пока
служило всяким господам, мы хотим, чтобы это орудие служило в нашей
стране одному господину и служило не как раб, а чтобы носители этой
науки сами творчески, кан хозяева, участвовали в общем деле строительства
социализма...	
\end{quotation}

Шмидт -- политик, представляющий власть на съезде и не последний деятель в Комакадемии
в тот момент.
Он должен балансировать на разных лезвиях, которые мы частично видим,
и которые мы, скорее всего,  видим не все. 
До какой степени он верил  в свои рассуждения о марксистской диалектике
и о предвещаемом ей расцвете,
сразу не ясно (он вообще-то был марксист и революционер),
но стоит отметить, что Шмидт представлял более мягкую точку зрения,
чем Кольман, Яновская, Хотимский, Райков и пр. Предлагается диалектику
<<прорабытывать>> (и говорится, что это будет полезно), делать это долго, но не непосредственно внедрять
в математику.

\sm
 
 {\bf\punct Математический съезд. Приветствие съезду ВКП(б).%
 \label{ss:sezd-VKP}}
 Цитируем книгу про <<Ленинградский
математический фронт>> \cite{Lenmatfront}:
 \begin{quotation}
Так как съезд математиков совпал по времени с началом XVI партийного съезда, 
то огромное большинство участников математического
съезда естественно пожелало послать приветствие коллективному
вождю пролетариата и руководителю социалистического строительства —
Съезду партии. При постановке вопроса о принятии приветствия Съезду
партии выявилось, как далеко отошли от настроения масс ее прежние
вожди, такие люди как академик Бернштейн, проф. Егоров, проф.
Гюнтер и другие. В своих попытках тормозить принятие приветствия
Съезду партии они оказались совершенно изолированными, и от них
отошли даже их бывшие сторонники. Когда же московские и ленинградские
математики-материалисты объявили открытое заседание
фракции материалистов, то больше половины съезда пришли на эту
фракцию и потребовали немедленного созыва пленума съезда
для принятия приветствия. Фракция материалистов настояла и на том, чтобы
доклады философского характера представителей Коммунистической
академии в Москве тт. Выгодского и Яновской были перенесены с
последнего дня съезда на один из ближайших дней.
\end{quotation}

Не обязательно все было ровно так (источник не абсолютно надежный),
но что-то в этом роде случилось. Нежелание Бернштейна, Егорова и Гюнтера
отправлять приветствие%
\footnote{В Стенограмме Съезда ВКП(б) сообщается о 1180 приветствиях, полученных Съездом по телеграфу,
	из них 43 от научных организаций. Приветствие от математиков было опубликовано в каком-то журнале.} на фоне происходившего погрома университетов и тяжелых картин коллективизации более, чем понятно...

На съезде состоялись 11 пленарных докладов
(среди них иностранные участники
 Монтель, Бляшке, Картан, Адамар, Данжуа, Блюменталь, Лихтенштейн).
 Кроме того был доклад Шмидта, открывающий съезд, доклад Кагана 
 «О теоретической и прикладной математике в СССР» и <<два соединенных доклада 
 всех секций>> назначенных под давлением <<фракции материалистов>>%
\footnote{В недавней историко-математической статье Т.~А.~Токаревой \cite{Tokareva-2001} это интерпретируется так:
	\newline
	{\it 
		... на Харьковском съезде 1930 г. 
		наблюдается заметное оживление интереса к истории и философии науки.
		На пленарное заседание Съезда выносятся два, связанных с этой
		проблематикой, доклада....}%
	}
(версии этих докладов есть в \cite{Dialektika}):
\begin{quotation}
\footnotesize
М. Я. Выгодский «Проблемы истории математики с точки зрения
методологии марксизма».

С. А. Яновская «Критика основных современных течений в области
оснований математики с точки зрения философии диалектического материализма»%
\footnote{Доклад посвящен тогдашним направлениям в математической логике.
Предваряется он эпиграфом <<Кто хочет знать врага, должен побывать во вражеской стране>>
Ленин (это немецкая пословица, ссылка на Ленина подчеркивает важность немецкой сентенции).
Например, о Лузине говорится следующее:
\newline{\it В так называемой Московской математической школе интуиционистские
идеи были встречены доброжелательно даже группой молодых советских математиков
во главе с А.~Я.~Хинчиным (см. статью последнего в Вестнике Комакадемии, N16) 
 Один из признанных вождей этой школы проф. Лузин является вполне законченным интуиционистом.
В докладе о современной борьбе течений в обоснованиях математики на первом конгрессе
по теории познания в точных науках (Прага, 1922) проф. А.~Френкель среди представителей этого направления упоминает и проф. Лузина.}
\newline
У Лузина в самом деле были определенные пересечения с интуиционистами
(а Гейтинг относил его к <<полуинтуиционистам>>). Но интуиционистом, и тем более законченным
интуиционистом,  не был.
}.
\end{quotation}
 
 %%%%%%%%%%%%%%%%%%%%%%%%%%%%%%%%%%%%%%%%%%%%%%%%%%%%%%%%%%%%%%%%%%%%%%%%%%%%%%%
 
 {\bf\punct Арест Егорова.%
 \label{ss:egorov-arest}}
 %В 1929-1930г. ОГПУ разрабаьывало большое дело <<Истинно православной церкви>>.
Егоров был арестован, согласно рассказу Щелкачёва \cite{ShCh}, в ночь с 9 на 10 октября.
Арестован он был по разрабатывавшемуся ОГПУ общеполитическому делу
контрреволюционной организации церковников <<Истинно-Православная Церковь>>.
У нас нет данных о том, что арест Егорова  был как-либо связан с его
профессиональной  деятельностью. У нас также нет  сведений 
о каких-либо доносах людей из математической или около-математической тусовки
на Егорова (хотя отсутствие таких доносов представляется весьма маловероятным).  Поэтому вопрос о том, что известно об обстоятельствах его ареста,
 вынесен в добавление к следующему параграфу.
 
  Он просидел в тюрьме несколько месяцев%
 \footnote{Возможно, почти 11, а, возможно, меньше.
 Из интервью Александрова 1971г. \cite{Alexandrov-Duvakin}: {\it Правда, он был освобожден, но уже за несколько месяцев до своей смерти.} },
был приговорен к 5 годам лагерей, который были заменены на 5 лет ссылки.
 Был выслан в Казань, где по-видимому, должен был работать профессором Казанского университета. Умер
  10 сентября 1931г. в больнице Казанского института усовершенствования врачей \cite{Demidov}.
 По-видимому, этому предшествовала голодовка Егорова. Виленкин в 1991г. излагал это так
 \begin{quotation}
Егорова арестовали и сослали в Казань (по одной
из версий, за сбор денег в пользу 
репрессированных единоверцев — он принадлежал
к религиозному течению «христославцев»).
Некоторые из его учеников, например,
B. В. Степанов, ездили к нему. По 
дошедшему до автора рассказу видного 
казанского алгебраиста Н. Г. Чеботарева, Дмитрий
Федорович объявил голодовку, сильно 
подорвавшую его здоровье. Он был помещен в
больницу, и, хотя его изолировали, 
Чеботареву удалось проникнуть к нему и 
побеседовать. В 1931 г. Егоров умер в возрасте
62 лет.
 \end{quotation}
 
Приведем также казанский источник -- В.В.Морозова\footnote{Морозов  Владимир Владимирович  (1910-1975) --известный алгебраист,
ученик Н.~Г.Чеботарёва. Его воспоминания о Чеботарёве \cite{Morozov}, которые мы цитируем, были написаны в 1963г. и изданы,
видимо, лишь в 2002г.}:
 \begin{quotation}
  Несколько лет спустя, в 1931 г., Н.~Г.[Чеботарёв] и Д.~Ф.~Егоров – тогда уже почетный член АН СССР –
встретились в Казани при обстоятельствах, не менее печальных, Д.~Ф. находился перед этим в
заключении, здоровье его сильно пошатнулось, возникла необходимость госпитализации, и
он был отправлен в Казань и помещен в клинику ГИДУВа; родоначальника Московской
математической школы навещали казанские математики, в первую очередь Н.~Г. Впрочем,
лечение было уже бесполезным, Д.~Ф. вскоре скончался и был похоронен на церковной аллее
Арского кладбища, вблизи памятника Лобачевскому.
 \end{quotation}

В 1932г. в Большой Советской Энциклопедии, том 24,
появилась такая статья
\begin{quotation}
Егоров Дмитрий Федорович (1969-1931), московский профессор математики,  имеет работы по анализу, теории чисел и др.,
не содержащие, однако, выдающихся научных открытий. Видный представитель реакционной московской математической школы
(идеалистической). Активно боролся против мероприятий Советской власти по реорганизации высшей школы
и научных институтов. После разоблачения <<егоровщины>>
был снят в 1929 с поста директора Института математики и механики
и в 1930 исключен из членов Московского математического общества.
\end{quotation}
Том был сдан в печать 16.5.1931г. в разгар свистопляски, о которой речь идет  в настоящем параграфе. Редакторская работа была закончена
15.6.1932г. Отделом математики в БСЭ в то время руководил Кольман, о чем чуть  ниже.

Отношение к Егорову было разным у разных математиков и со временем оно быстро менялось.
Мы видели, что  в 1936г. даже  самые радикальные из работающих математиков воспринимали Егорова сложно.
 В 1938г. в статье <<Математика>> для БСЭ Колмогоров писал:
\begin{quotation}
 В это же время в Москве (где ранее господствовали глубоко провинциальные
 и реакционные научные течения, возглавлявшиеся Н.~В.~Бугаевым)  создается научная школа 
 в области дифференциальной геометрии [Петерсон (1828-81), Д.~Ф.~Егоров(1869-1931), Б.~К.~{\it Млодзеевский}, (см.) 
 (1858-1923)]%
 \footnote{Квадратные скобки и курсив из оригинала. Характерно, что на энциклопедическую статью 
 о Егорове ссылки нет.}%
 и особенно в теории функций действительного переменного [Егоров, Н.~Н.~Лузин (род. 1883) и др.]
\end{quotation}
Далее упоминается слово <<первоклассные научные школы>>,  относящееся именно сюда.
Весьма положительно писал о Егорове и Александров в статье 
\cite{Alexandrov-1937} 1937г.

Превращение Егорова в однозначно положительную фигуру  и как администратора тоже было
завершено сразу после Войны.

 Вот еще интересный документ, приведенный в книге  \cite{Ilchenko}, дата не сообщается. Автор письма Петровский был деканом в 1940-44г..
 Кажется, упоминаемое <<литерное снабжение>> -- реальность, относившееся уже к военному времени.
 \begin{quotation}
ИЗ ПИСЬМА ДЕКАНА МЕХАНИКО-МАТЕМАТИЧЕСКОГО ФАКУЛЬТЕТА
И.~Г.~ПЕТРОВСКОГО ПРОРЕКТОРУ МГУ В.~И.~СПИЦЫНУ

Прошу Вас возбудить ходатайство о предоставлении литерного снабжения ...
Анне Ивановне Егоровой ... Анна Ивановна Егорова -вдова почетного члена
Академии Наук СССР, одного из старейших профессоров нашего факультета Дмитрия
Федоровича Егорова. Дмитрий Федорович был основателем современной Московской
математической школы, получившей мировую известность. Долгое время он
состоял Президентом Московского Математического Общества. Анна Ивановна
Егорова достигла уже преклонного возраста. {\bf Она получает персональную пенсию},
не имеет возможности по состоянию здоровья служить и находится сейчас
в тяжелом положении».

Архив МГУ. ф. 1, оп. 28, россыпь 
\end{quotation}

 %%%%%%%%%%%%%%%%%%%%%%%%%%%%%%%%%%%%%%%%%%%%%%%%%%%%%%%%%%%%%%%%%%%%%%%%%%%%%%%
 
 {\bf\punct Переворот в Мат.Обществе.%
 \label{ss:perevorot}} В списке Александрова--Головина (см. п \ref{ss:1930-spring})
 последний доклад на старом мат.обществе помечен 10 октября.
 Так или иначе \cite{Deklaratsiya},
 \begin{quotation}
 ... На арест своего председателя Об-во%
 \footnote{Лузин оставался вице-президентом, Привалов -- секретарем.} никак не реагировало
и назначило обычное деловое заседание  [на 21 ноября 1930] с докладами ближайшего 
соратника Егорова по Ин-ту и Об-ву - Финикова
и исключенного только что из комсомола Куроша%
\footnote{Возможно, что Курош был исключен из Комсомола потому, что усомнился 
в процессе Промпартии. См. Рассказ И.~Р.~Шафаревича \cite{Sha} о катании на лыжах на Кавказе
(но у меня нет уверенности в именно этой интерпретации рассказа).}.
 \end{quotation}
 Но реакционеры не прошли,
  \cite{Sbornik-revolution}:
  \begin{quotation}
{\bf В Математическом об-ве была создана инициативная 
группа по его реорганизации, куда вошли ряд видных московских 
математиков}; декларация инициативной группы
[см. следующий пункт], направленная как 
против активно реакционных элементов в среде математиков, так и 
против мнимо советской позиции буржуазно-де\-мо\-кра\-ти\-че\-ских 
попутчиков первого периода НЭПа, была принята Математическим
об-вом. Об-во исключило из своих рядов Егорова и других реакционеров, пополнило свой состав, прежде всего, за
счет аспирантуры. Новый устав, выработанный об-вом, ставит
его целью, прежде всего, поворот советской математики на 
обслуживание социалистического строительства. Был реорганизован также
Математический кружок, влившийся в Об-во в качестве его 
педагогической секции.

Последнее распорядительное собрание Об-ва избрало новые
выборные его органы: президиум, комиссии по исследовательской, 
педагогической и массовой работе, редакции журналов «Математический
сборник» и массового журнала.
\end{quotation}

Перейдем к исследованию упомянутой декларации, а потом вернемся к картинам переворота,
описанного его участниками.

\sm

 %%%%%%%%%%%%%%%%%%%%%%%%%%%%%%%%%%%%%%%%%%%%%%%%%%%%%%%%%%%%%%%%%%%%%%%%%%%%%%%
 
 {\bf \punct Декларация инициативной группы.%
 \label{ss:deklaratsiya}}
 В номере 11-12 журнала <<Работник просвещения>> была опубликована  <<Декларация инициативной группы
 по реорганизации математического общества>> (\cite{Deklaratsiya}, она перепечатана в \cite{Tok-white} и легко доступна).
 
 \begin{quotation}
  Обострившаяся классовая борьба в СССР толкнула правую часть профессуры в лагерь контрреволюции. 
  Реакционная профессура возглавляла все раскрытые в последнее время вредительские организации
  и контрреволюционные партии. Благодаря блестящей деятельности ОГПУ разоблачены преступления целого ряда
  научных бонз, умевших искусно скрываться за разными масками - от холодной лойяльности до шумно рекламируемой
  горячей приверженности советской власти. И в среде математиков выявились активные контрреволюционеры.
  Арестован за участие в контрреволюционной организации проф. Егоров, признанный вождь московской математической
  школы,...
 \end{quotation}

 Этот хорошо известный, но малоисследованный документ, который 
 со своей трескучей политической риторикой
 производит дикое впечатление. 
 Однако автор
 данных записок, имевший сомнительное удовольствие прочитать много разных сочинений 1929-1931гг.,
 должен отметить наличие в декларации иных черт,  выделяющих  ее на комакадемическом фоне
 кольмановского типа
 (и выделяющих ее в лучшую сторону).

 \sm
 
 {\sc Статус документа.} Автору пришлось просматривать разные околонаучные публицистические  издания того времени,
 и он не мог бы сказать, чтобы математике, а особенно событиям  на московском Физмате, уделялось особое внимание. Страна большая,
 городов много, институтов много, наук много. Эта декларация полностью занимает 5 страниц журнального текста,
 в конце стоят подписи Л.~А.~Люстерник, Л.~Г.~Шнирельман, А.~Гельфонд, Л.~Понтряши (в статье так), Некрасов.
 
 Подписи стоят последней строчкой, уже за пределами обычной типографской страницы, инициалы Некрасова
 обрезаны, а инициалы Гельфонда и Понтрягина представлены частично, в строчке не остается места ни для одной буквы.
 Очевидно, текст был важнее точных подписей
 (можно допустить также какую-то вставку в текст на этапе корректуры и обрезание части подписей).

 По тогдашним временам определенного рода газетные и журнальные статьи имели инструктивный статус.
 Очевидно, этот статус имела и обсуждаемая статья. В то время под руководством Комакадемии проходила советизация
 университетов, частью этой операции была советизация и идеологизация профессиональных сообществ 
 (итоги преобразования <<обществ>> к весне 1932 года окажутся плачевными,
 о чем говорится в статье Стецкого, приведенной ниже в п.\ref{ss:stetski}).  Данная публикация имела значение <<Делай как я!>>.
 Чтобы это поняли недостаточно понятливые, к заголовку была поставлена сноска:
 
 \begin{quotation}
  Печатая настоящую декларацию инициативной группы Московского математического общества,
  редакция считает необходимым привлечение
  внимания научной общественности к работе всех обществ с целью расширения
  борьбы с реакционными элементами в них 
  и перевода их на рельсы активного содействия социалистическому строительству.
 \end{quotation}

 В заголовке отсутствует слово <<московского>>, что экономит одну строчку.
 Возможно, что это тоже было ужиманием текста,
 хотя сам заголовок можно было бы подвинуть вверх. Возможность сдвинуть заголовок, возможно,
 свидетельствует о вставке в конце статьи.
 
 Статья была перепечатана в журнале <<Научное слово>>, в номере 1 за 1931г. 
 К подписавшимся приставили должности: проф. Люстерник~Л.~А., проф. Шнирельман~Л.~Г.,
 проф. Гельфонд~А., доцент Понтрягин~Л., проф. Некрасов, не добавив опущенных инициалов%
\footnote{Т.~А.~Токарева \cite{Tok-white} (на основании титула <<проф.>>)
 твердо утверждает, что инициалы  
последнего из подписавших декларацию были <<А.И>>. Но
Гельфонд окончил аспирантуру в конце 1930г., автор не смог найти подтверждения того, что 
тот  к моменту подписания Декларации был профессором (в начале 1931г. он им стал). В связи с этим,
автор не уверен в истинности предлагаемой атрибуции.
 Форд \cite{Tok-white} предлагал инициалы К.П.
(это был аспирант),
но его аргументация, видимо, не публиковалась и автору неизвестна. В любом случае, ясно, что при перепечатке статьи
точностью данных об авторах не слишком озабачивались. Также непонятен принцип, по которому
перечислены авторы. Первое место у Люстерника может свидетельствовать как о лидерстве, так и о том, 
что он был
<<старшим по званию>>.
\newline
В декларации есть фраза:
{\it Большинство из нас, подписавших эту декларацию, 
беспартийные научные работники, до сих пор не принимавшие участия в этой
борьбе.} Первые четверо были беспартийными, то есть пятый вроде бы был  партиен
(известно, что А.~И.~Некрасов в 1923г. вступил в РКП(б), а в начале
1934г. был беспартийным, \cite{NekrasovAI}). Однако в \cite{Reorganizatsiya} говорится:
{\it 
...инициативная группа по реорганизации Общества, куда вошел
ряд математиков как партийцев, так и беспартийных, большая часть
молодых советских ученых.}
Скорее это подтверждает гипотезу об обрезании   списка фамилий.
Вот еще довод против подписи А.~И.~Некрасова: на выборах в АН СССР 1932г. Лузин решительно
поддерживал кандидатуру А.~И.~Некрасова, см. его письмо Крылову от 22.12.1932 \cite{Erm1}:
{\it Ведь труды Некрасова очень сильные, стиля Poincar\'e,
широты и размаха крупного мастера, и редкого изящества.}
\newline
Напомним, что <<Инициативная группа>> свергала Лузина.
\newline
Есть статья  Nekrassoff, V. A.
Nomography in applications of statistics.  
Metron 8, No. 3, 95-99 (1930). В статье нет никаких данных 
о ее авторе
(быть может, он жил в Праге).}.

\sm
 
 {\sc Об авторах.} Люстерник  работал в то время в Комакадемии, по этому своему положению
 он должен был (наряду с Яновской) входить в <<группу захвата>> московского Физмата.
 Особый статус первых трех авторов (а четвертый был их другом)
 в московской математике 1929-1937гг. достаточно очевиден. 
 Например, из статьи Яновской \cite{Yanov-PZM-1}, 1930:
 \begin{quotation}
 В СССР в текущем году сделаны работы,
которыми гордилась бы любая крупнейшая западно-европейская «школа», и
сделаны как раз молодыми советскими математиками, группирующимися
вокруг Коммунистической Академии (тт. Гельфонд, Шнирельман, Понтрягин).
\end{quotation}

Тут важно не то, что из себя представляли эти математики (ими в самом деле гордилась бы любая
западно-европейская школа, как и еще более неназванными Колмогоровым и Хинчиным), а то, что это говорит Яновская,
взгляды которой на правильность той или иной математики формировались тогда на диалектической основе.

А вот воспоминания Понтрягина: 
  \begin{quotation}
Одной из основных задач, которую ставил перед собой И.~М.~Виноградов, руководя институтом
 [когда институт Стеклова переехал в Москву], 
 было привлечение в него молодых талантов, математиков с хорошей, разумной, по его мнению, математикой...

Среди вновь привлекаемых в институт москвичей назывались шесть, которые рассматривались тогда, как молодые и талантливые.
В том числе был и я. 
  Любопытно отметить, что эти шесть человек классифицировались на три пары по их «качеству».
  На первом месте стояли А.~О.~Гельфонд и Л.~Г.~Шнирельман, на втором месте — М.~А.~Лаврентьев и Л.~А.~Люстерник, 
  а на третьем месте — Л.~С.~Понтрягин и А.~И.~Плеснер. То, что на первом месте стояли {\bf Гельфонд и Шнирельман},
  было в то время очень естественным. {\bf Эта пара составляла обойму, которую называли всегда и везде,
  когда хотели указать наиболее «талантливых» 
  молодых советских математиков.}%
  \footnote{Понятно, что такого рода списки никогда не бывают вполне справедливым.
  Однако отметим отсутствие Колмогорова и Новикова. Александров и Хинчин тоже были довольно молоды.
   Отметим, что среди шести -
  четыре подписавших обсуждаемую декларацию и политиммигрант 1932г. Плеснер Абрам Иезекиилович  (Plessner Abraham,  1900—1961),
  впрочем, родившийся в Российской империи. Он сыграл определенную роль в продвижении функционального анализа в СССР,
  но стоит иметь в виду, что Люстерник, Шнирельман и Колмогоров сами около 1930г. вышли в функциональный анализ из ТФДП. Еще стоит отметить, что список из 6 самых перспективных математиков включает 4 авторов <<Декларации>>.}%
  $^,$%
  \footnote{Например Кольман, 1934 \cite{Kolman-mass},
  \newline
  {\it ...советских ученых новой формации, давших уже ряд выдающихся исследований и открытий, как Гельфонд, Шнирельман...} }
 \end{quotation}
 
 В конце концов, о  политическом радикализме Люстерника, Шнирельмана и Гельфонда
 нам известно из Стенограммы 1936г., которая цитировалось выше и
 будет цитироваться чуть ниже.
 Определенной оценкой их социального-политического веса является наличие статей в <<Правде>>.
 У Шнирельмана была статья от 4.3.1937 года (о ней ниже \ref{ss:kiselev}), а у Гельфонда от 15.1.1941
 (<<Математическая наука в СССР>> \cite{Gelfond-pravda}.
 
 Есть такой отсвет об их деятельности в Москве 1931г.(кстати, Шнирельман в 1931г. наблюдается
 в Московском индустриально-педагогическом институте им.  Либкнехта, хотя по его биографиям он был Новочеркасске):
 \begin{quotation}
Колмогоров -- Александрову, 23 марта 1931 г.
 
 Интересно, что вся
компания (Лазарь Аронович, Лев Генрихович и Александр Осипович),
видимо, и сейчас много думает и работает в математике,
несмотря на большую активность в <<организационных>> вопросах и весьма
солидную нагрузку по специальным ВУЗам.
 \end{quotation}
 
 \sm
 
 Ну и наконец, снова слова Люстерника, но уже из Стенограммы-1936
 (непонятно, относится ли это к декабрю 1929 или ноябрю 1930г.):
 \begin{quotation}
    ЛЮСТЕРНИК.
 Я несколько раз делился со своими товарищами одним  фактом,  которому
я до сих пор объяснения не нахожу. Это то, что польская пресса%
\footnote{Кстати, выходит, что польская пресса -- потенциальный исторический источник.
Другое дело, что найти данные по такому поводу в польской прессе не так уж просто.}
очень подробно освещает
математическую жизнь Советского Союза. Дело в том, что я сам родом из Польши [Белосток], и мои
родители живут в Польше, и письма, которые я иногда получаю, говорят такие вещи:  «Мы
читали в газетах о таком-то твоем выступлении».  Здесь об этом в  нашей прессе ничего
нет. Например, во время дела Егорова — что они читали о моей роли в этом деле,  между
тем как в нашей прессе это не освещалось.... Причем мое выступление там конечно освещается
с той точки зрения,
что мол такой-то карьерист,  который выступил для того, чтобы сделать себе научную
карьеру.
 \end{quotation}

 \sm
 
 {\sc Взгляд на дореволюционную русскую математику в декларации.}
Кого-то%
\footnote{Например, Т.~А.~Токарева \cite{DeEs}:
\newline
{\it Здесь можно поспорить с авторами «Декларации», и еще раз,
усомниться в причастности подписавшихся к ее составлению. Могли,
ли члены ММО не знать года его основания? А утверждение о том,
что до революции не удалось создать «серьезной научной школы в
математике»? Ведь уже в последней трети XIX в. существовала одна
из лучших математических школ того времени — Петербургская школа
 П.~Л.~Чебышева (1821-1894). В это же время в недрах 
Московского математического общества возникли и стали играть заметную роль
в жизни математической Европы научные школы: по прикладной 
математике (Н.Е.Жуковского (1847-1921) и С.~А.~Чаплыгина
(1869-1942)), к которой принадлежал и их ученик А.~И.~Некрасов; по
дифференциальной геометрии (К.~М.~Петерсона (1828-1881),
Б.К.Млодзеевского (1858-1923), Д.~Ф.~Егорова). Наконец, в начале
1910-х гг. зародилась одна из наиболее значимых в XX в. 
математических школ — школа теории функций действительного переменного
Егорова-Лузина. Могли ли «лузитане» второго поколения - Люстерник и Шнирельман
- и ученики учеников Лузина - Гельфонд и Понтрягин — отрицать это?}
\newline
Автор позволит себе это прокомментировать: упомянутые московские школы
прикладной механики и дифференциальной геометрии не оказали сколько-либо
заметного
влияния на математическую  жизнь Западной Европы того времени
(хотя имена Петерсона, Егорова, Лахтина упоминались иногда западноевропейскими математиками
начала XX века).
 \newline
Процитируем Выгодского \cite{Vygodsky-1948}: {\it Работы Петерсона не нашли себе должной оценки даже
в среде [Московского] Математического общества, членом которого он был.
При жизни Петерсона ни один автор не поместил в 
Сборнике ни одной статьи, развивающей идеи Петерсона.} Несколько
позже 
эти работы приобрели большую популярность у московских геометров.
В Западной Европе работы К.~М.~Петерсона 1853, 1866-83гг.
тоже стали известны  с большим запаздыванием,
когда в значительной своей части были уже повторены. Они были замечены в 1884-86гг. Штаккелем (Paul St\"ackel),
а в 1905г. переведены
на французский. Изданная  самим Петерсоном по-немецки книга \"Uber Kurven und Fl\"achen, Moskau und Leipzig (1868),
по-видимому, прошла незамеченной (С.~П.~Фиников \cite{Finikov}: {\it книжка ... затерялась среди математической литературы}).
 \newline
Жуковский с Чаплыгиным по сути не были провинциалами и, в самом деле, наследие
их научной и организационной работы (в виде аэрокосмической промышленности СССР)
оказало заметное влияние на мысль (и не только на мысль) Европы. Но всё это будет
после 1929 года.}
 могут удивить такие слова из Декларации:
\begin{quotation}
 «Математическое общество возникло с другими научными 
обществами в 70 -х годах прошлого столетия, в период роста 
буржуазной культуры в старой России. Однако серьезной научной школы в
математике... создать в царской России не удалось.»
\end{quotation}
Однако о тогдашних взглядах математической тусовки
рассказывал  сам Люстерник 35 лет спустя:
\begin{quotation}
Репутации «основателей ММО [Московского мат. общества]» 
среди тогдашней математической молодежи
сильно повредили наукообразные статьи профессора МУ [Московского университета] П.~А.~Некрасова,
которые он печатал в «Математическом сборнике» в начале XX века
и в которых он пытался обосновать необходимость царивших в то время
порядков. В этих бредовых статьях П. А. Некрасов выступал от имени
«основателей Московского математического общества», к числу которых
он сам не принадлежал и никого из которых не было тогда в живых. Так
создалось долго бытовавшее и отразившееся в некоторых публикациях
мнение о реакционности «основателей ММО». И лишь когда в 1940 г. 
отмечался 185 (!)
  (тогда праздновали
и такие «некруглые» юбилеи)%
\footnote{Да уж, непонятно, если не вспомнить статью Удальцова по поводу 175-летия МГУ...}-летний 
юбилей Московского университета, Академия праздновала 220-летие%
\footnote{См. сноску \ref{fo:220}.}), появился
ряд исследований по истории математической жизни университета и Общества,
произошла «реабилитация» основателей Общества (этому была посвящена
статья М. Я. Выгодского 1940г. [Ссылка у Люстерника: Математика и ее деятели в Московском университете во второй
половине XIX века, Ист.-матем. исследования 1 (1948), 141 —183. Дата другая, но о той же статье Выгодского говорил Юшкевич
\cite{DemidovTokareva}, она не вышла в 41г. из-за Войны])%
\footnote{Читатель, которого удивит, что состояние математики в дореволюционной Москве
оценивалось по философским статьям одного автора в <<Мат.Сборнике>>, может вспомнить, что обстановка
в московской же математике 30х гг. сейчас оценивается по статьям Кольмана 1930--31гг.
(еще вспомнить отношение либеральной интеллигенции к <<совку>> в 90е годы)
\newline
Кстати, Некрасов  Павел Алексеевич (1852-1924), который побывал и ректором МГУ и попечителем Московского учебного округа,
был также  крупным и оригинальным математиком. Но передовые молодые люди едва ли читали
его  научные работы. Даже когда заслуги предреволюционной Московской математической
школы уже были официально признаны, Некрасова это признание не коснулось.}.
\end{quotation}
Остается вопрос об отношении к Егорову и к Лузину. Ну тут с одной стороны возникает парадокс кучи:
а два человека - это уже школа? Или еще нет? А с другой стороны стоит вспомнить  цитированные выше стихи Люстерника
про Лузина:
<<{\it Не мы ль его раздули славу}?>>.  

 Может остаться удивление в отношении Петербурга: все ж таки
Чебышёв, Марков-ст., Ляпунов... Нехорошо как-то.  Между тем Люстерник в той же статье
рассказывает
про оборвавшуюся связь между научными центрами... А вот комментарий с ленинградской стороны,
Б.~Н.~Делоне \cite{Delone}, 1973:
\begin{quotation}
 Но ленинградцы эту лузинскую французскую математику даже математикой не считали — считали философией какой-то, 
 и его избрали в пику москвичам по философии, вместе с Дебориным....
 
И вот между школой Эйлера-Чебышева петербургской и школой Лузина московской — собственно, французской,
парижской — все время был такой антагонизм, 
что те этих не понимали, эти — этих, пока Академию не перевели в Москву.
Когда в 34-м или 35-м, в начале, перевели Академию в Москву, мы начали сближаться, сближаться,
и вот из этого сближения обеих школ и получилось, то, что мы сейчас называем «советская математика».
\end{quotation}

Ну и добавим слова Колмогорова из статьи <<Математика>> в БСЭ, сданной в 1938г.:
\begin{quotation}
Таким образом русская дореволюционная наука, оставаясь в целом отсталой и слабой по сравнению с наукой 
передовых буржуазных стран, все же выдвинула в 19в. двух математиков самого высокого ранга, 
а к 20в. пришла также к созданию
в некоторых областях М., пусть немногочисленных и иногда чрезмерно замкнутых, но все же 
первоклассных научных школ.
\end{quotation}
С 30г. утекло уже много воды, Московская школа (Петерсон, Егоров, Млодзеевский, Лузин) уже <<реабилитирована>>, 
в частности, в самой статье Колмогорова. Но в отношении <<отсталости>>, <<слабости>>
и недоброго старого времени инерция все еще сохранялась%
\footnote{Для сравнения, из редакционной статьи в Успехах <<Советская математика за 20 лет>>, тот же 1938г.
\newline
{\it 
Русский
народ в тяжелой обстановке царской России выдвинул ряд первоклассных
математиков — в первую очередь Лобачевского и Чебышёва, далее
Остроградского, Ляпунова, Маркова, Вороного, Золотарева и др. Эти ученые
были фактически изолированы...}
\newline
О Лузине там уважительно, а о Егорове - ни слова. Впрочем, основные 
достижения Егорова были совершены до 1917г.} (или, по крайней мере, Колмогоров сохранял некоторую осторожность
в своих высказываниях).
В статье Александрова \cite{Alex-1945} 1945г. уже  излагалась нормальная история,
да и Люстерник \cite{Lyusternik-1946} в 1946г. писал уже иначе.

\sm

Цитированная сентенция из Декларации, конечно же, несправедлива по отношению к дореволюционной математике
России.
Но не видно, чтобы это противоречило тому, что мы знаем о взглядах тогдашней лево-радикальной тусовки,
скорее мы видим соответствие этим взглядам.
Нет причин на основании этих слов ставить под сомнение авторство подписавшихся.

\sm

 Возможно, что глаз современного читателя могут резануть слова:
 \begin{quotation}
  Общество... окружало невероятной помпой его [Егорова] выступления, по существу, скромного научного значения.
 \end{quotation}
 Но тут нужно посмотреть на список подписей. Эти молодые люди к тому времени 
  сделали блестящие работы,
 и, скорее всего, если убрать эмоционально-усилительные интонации,
 скорее всего,  доклады Егорова  уступали их последним достижениям%
\footnote{В связи с особым интересом к Егорову как к жертве политических репрессий и как к человеку,
пострадавшему
за христианские убеждения, современные читатели и авторы часто неправильно понимают значение Егорова.
Роль его была выдающейся,  без него Московская математическая школа едва ли могла бы сложиться в том виде, 
в котором она сложилась (это относится к его деятельности в в 1900-1917гг). Он был хорошим 
(и, видимо, жестким)
администратором в сложные 20е годы,   расцвет московской математики  начинался в  пору его администрирования.
В плане истории
идей (истории математики в строгом смысле) он был прежде всего предтечей Лузина,
а ученики Лузина (по крайней мере, Хинчин, Колмогоров, Александров) превзошли учителя 
уже в 20е годы. Подписавшие письмо молодые люди имели определенные основания смотреть
на Егорова сверху вниз.
\newline
Добавим, что  упоминающаяся в историко-математической литературе <<Школа
Егорова--Лузина теории функций действительного переменного>>
отчасти соответствует истине, но отчасти способствует неправильному пониманию 
картины. У Егорова были две статьи по ТФДП, 1911, 1912 (причем вторая, по-существу, была добавлением
к работе Лузина). Первая публикация
Лузина была в 1911г., а в 1912г. Лузин опубликовал 7 статей (но дело-то не в количестве). То, что было позже,
(Хинчин, Александров, Колмогоров, Меньшов, Бари),
было школой Лузина (да и кому было бы  под силу состязаться
с такими молодыми людьми). Егоров же занимался дифференциальной геометрией.}.
С другой стороны
 self-esteem этих молодых людей был чрезвычайно высок (а о том, что где-то у них были ляпы,
 а где-то приоритет был не за ними, они просто не знали).

\sm

{\sc Политическая часть декларации.} 
 Б\'ольшая часть декларации посвящена проблеме политической лояльности научных кадров, <<классовая борьба
 в среде научных работников>> обсуждается  в весьма жестких выражениях.

 Интересно, что самому Егорову предъявлялись обвинения в политизированности и идеологизированности%
 \footnote{{\it Освобождаясь от идеологического плена егоровщины, большинство математиков...}}, в неприеме
 в общество ма\-те\-ма\-ти\-ков-ком\-му\-нис\-тов%
 \footnote{Шмидт был членом Общества с 1922г. Известные автору тогдашние коммунисты:
 Райков (тогда еще не математик), Хворостин, 
 Яновская, Хотимский (возможно, с другими статистиками), Выгодский
 и Кольман.},
 а также каких-то <<крупнейших молодых ученых>>. Интересно, что позднее
 (скажем в Стенограмме 1936г.) таких обвинений от математиков слышно не было. Скорее всего, они не были справедливы%
 \footnote{Впрочем есть  такое свидетельство (С.~С.~Демидов \cite{Demidov})
 \newline
 {\it
 В московских математических кругах до сих пор вспоминают
об его отказе читать лекцию по математике в бывшей церкви,
приспособленной под аудиторию. Этот поступок был расценен как
«пощечина пролетарскому студенчеству». Владимир Николаевич
Молодший} [1906-1986, философ от математики, математических работ не известно, однако доктор физ.мат.наук], 
{\it бывший в 20-е годы
 аспирантом, рассказывал мне, как
однажды, встретив Д.~Ф.~Егорова в коридоре университета, 
попросил разъяснить одно непонятное место в курсе Э.~Гурса. Егоров
сразу согласился, но по ходу своего разъяснения увидел у Молодшего
значок Коммунистического Союза Молодежи. Он сразу, 
рассказывал Молодший, изменился в лице и, сославшись на то, что
ему некогда, прекратил разговор. Разумеется, столь 
демонстративное поведение не могло быть терпимым слишком долго, хотя 
поначалу-власти смотрели на все сквозь пальцы.}
 \newline
 Если верить известию В.~Н.~Молодшего, 
 то профессиональное поведение Егорова не было политически индефферентным.
 Автор был бы осторожен в этом отношении (для примера, Егоров мог идентифицировать
 конкретного аспиранта, увидев этот значок). Есть объективные данные: судя по списку
  докладов на Мат.обществе времен Егорова \cite{Alex-MMO2},  ни о  какой политизированности 
  подбора докладчиков говорить не приходится (кстати, и один из руководителей Общества Костицын был человеком
  весьма левым).
  \newline
 Усомнившись в оценке сообщения В.~Н.~Молодшего, нельзя не позабавиться 
 преемственностью эпох: судя по интонации, современный комментатор  считает, что
 демонстративная дискриминация студентов 
 по политическим взглядам - дело нормальное, лишь бы вектор дискриминации был правильным.
 И власти должны смотреть на это сквозь пальцы.}.

 О словах про <<блестящую работу ОГПУ>>.
Современный российский читатель, который возмутится такими речами своих
 праотцов,   испытает чувство своего духовного
 превосходства над ними, или, наоборот, заподозрит, что к праотцам могло быть
 применено насилие, может освежить в своей памяти статью \cite{Gadina}%
  \footnote{По-видимому, пройдет еще 50 лет, и все будут знать,
 что подписи о солидарности с убийствами 1993г. и с  призывами к 
 политическим арестам  собирались под дулами танковых
 орудий. Этот процесс постижения истории уже начался.
 \newline
 Кстати, это совсем не проблема <<поведения людей при тоталитаризме>>,
 как часто говорится. Декларация \cite{Gadina} никакого отношения к оному поведению не имеет.
 Хотя проблема эта интересна, а само явление  неистребимо.}.
 
 Ну и общеисторическое замечания. В представлениях современной передовой 
  пост-советской интеллигенции, все политические дела 20-40х годов были основаны
 на фальсифицированных
 данных (и о том же заявляют результаты политических реабилитаций 1950-1960х и  1980-2000х годов).
 Принимая данный факт за истинный, мы обнаруживаем, что у Советской власти не было противников, 
 и уж, конечно, не было  противников
 организованных и вооруженных. Это можно объяснить лишь тем, что Советская власть
 устраивала всех, в том числе и тех, чьи дела миллионами (десятками миллионов и т.д.)
 фальсифицировались. Передовая интеллигенция 20-30 годов (по крайней мере, до Тридцать Седьмого
 года) таких взглядов не разделяла, а, напротив, считала, что 
 чекист -- благородный герой,  защищающий честных граждан Страны Советов
 и показывающий врагам  кузькину мать.
 
 И вот представьте себя на месте тогдашних честных
 лево-радикальных математиков. Глава Московского математического общества, человек, с которым они сами
 столько лет общались, оказался вовсе не <<попутчиком>>%
 \footnote{Слово было такое.}, а активным
 контррреволюционером! С одной стороны, это должно было стать духовным шоком,
 с другой -- люди могли испытать облегчение и радость от того, что враг, наконец, разоблачен.
 Потом, по-видимому, станет известно, что активным контрреволюционером он не был,
 а чекисты просто не сразу во всем <<разобрались%
 \footnote{Слово было такое}>>. Но случится  это через несколько месяцев.
 
Так или иначе, не видно, чтобы эта бумага не соответствовала политическим взглядам
лиц, ее подписавших. Скорей хорошо подходит (вспомните, что первые трое авторов говорили
в Академической комиссии 1936г., и вспомните двукратное предложение
одного из подписавших декларацию передать дело Лузина в НКВД).
 
\medskip

{\sc Диалектика и математика.}
А теперь главное, {\bf что в бумаге необычно на фоне комакадемической продукции 1930-1931гг.}

Первое. При безумной политизированности,
в документе почти отсутствует исчисление диалектическо-математических ангелов на острие иглы.
Точнее, в списке из восьми планируемых
мероприятий есть  такое
\begin{quotation}
 4. Борьба за марксистское революционное созерцание в вопросах математики. Организация секции 
 марксистской методологии и истории%
 \footnote{История математики в то время была идеологическим предметом, см. ниже   \S\ref{s:history}.} математики.
\end{quotation}

Кроме того, есть фраза:
\begin{quotation}
 Общество не ставило, конечно, вопросов марксистской методологии науки. 
 Зато в его рядах процветали поповщина и мракобесие.
\end{quotation}

Для пламенной 5-страничной бумаги это очень немного.
Упоминается <<партийность науки>>, но упоминается в качестве вопроса кадрового.
В программу входит борьба с идеалистической философией,
но не борьба с идеализмом в математике, а  политическое перевоспитание математиков 
в интересах социалистического
строительства. Даже в последней цитированной фразе говорится, что поповщина была, а марксизма не было.

В связи с этим стоит отметить слова Яновской 1930г. \cite{Yanov-PZM-1}:
\begin{quotation}
Но с особенным презрением [на Матобществе] отзывались о диалектике в математике:
в лучшем случае, в среде наиболее нам близких математиков, работавших
в Коммунистической Академии%
\footnote{Известны в таком качестве Люстерник, лузитанин Лихтенбаум,
	а также Гельфонд (он работал в биологическом Тимирязевском институте,
	входившем в систему Комакадемии).
	 про <<близких>> 
к Комакадемии было чуть выше.},
говорили о собрании математических иллюстраций
 для законов диалектики, более или менее удачных, и зло издевались
над философами за малейшую математическую  ошибку%
\footnote{Отвлекаясь от темы этого пункта, процитируем сентенцию Яновской из следующего абзаца:
\newline
{\it
Даже такой «левый» математик, как Мюнц, специально
выписанный Коммунистической Академией из-за границы, где выступал
действительно как подлинный сторонник советской власти, в докладе в
Ленинградском математическом обществе «показывает», что никакого кризиса
современной математике нет, что кризис специально изобретен философами,
помощь которых науке может состоять лишь в изобретении для нее таких
мнимых кризисов, на разрешении которых зато и наживается философия.
}
\newline
О том же эпизоде из книги <<На Ленинградском математическом фронте>>:
\newline
{\it 
... лекция была весьма и весьма бледной. Сколько ни
старался лектор, игнорируя все то, что сказано Энгельсом о математике,
что сказано В.~И.~Лениным о кризисе буржуазного естествознания,
доказать, что «на Шипке все спокойно», — ничего убедительного не
получилось, и к концу лекции чувствовалось сильное разочарование.}}. 
\end{quotation}
Похоже, что в обсуждаемой бумаге <<марксистская  методология>>
проговорена сквозь зубы (а про <<поповщину и мракобесие>>, конечно, искренне)%
\footnote{Кажется маловероятным, чтобы  эту бумагу не просматривали представители политсостава.
Также  не исключено, что опубликованный текст бумаги мог быть <<улучшен>> 
по сравнению с начальным содержанием.}.

Второе. Бумага содержит элементы позитивной программы, которая, впрочем, перемешана с политической риторикой.
\begin{quotation}
1... Необходима пропаганда математики в широких массах пролетариата, выдвижение на научную работу
 талантливых рабочих. Необходима организация популярного журнала и сети корреспондентов при нем.
 Все это должно
 быть предметом деятельности секции массовой работы при обществе...
 
 2... Выдвижение новых тем математической работы. Организация секции математических методов
 в технике и экономике.
 
6 ... Организация  секции преподавания  математики, куда должен влиться реорганизованный столь же архаичный 
Московский математический кружок.

7. Организация математической работы в провинции%
\end{quotation}

Из 8 пунктов программы четыре предполагают (если очистить их от трескотни) положительные действия%
\footnote{Еще один посвящен политическому перевоспитанию кадрового состава, один - борьбе с идеологическим
влиянием буржуазной науки (и ответному влиянию). Последний пункт: {\it Инициатива в деле планирования
научной работы в области математики.} Планирование научной работы тогда жестко навязывалось сверху,
оно уже упоминалось выше в связи со списком докладов в НИИММ.
О том, кто взял на себя упомянутую инициативу, мы вскоре увидим. Не те, кто подписывал Декларацию.}

Кстати, Люстерник, Гельфонд и Шнирельман  займутся распространением  математических знаний
в ближайшие годы...

\sm

{\sc Сравнение с ленинградским аналогом.}
 В марте 1931г. долго артачившимся 
ленинградским математикам приехавшие из Москвы товарищи таки показали кузькину мать
(о чем чуть ниже, п. \ref{ss:moskva-leningrad}).
Тогда появилась <<Декларация
инициативной группы по реорганизации Ленинградского физико-математического
общества>> \cite{Lenmatfront}, под которой хошь-не-хошь людям неподобающих взглядов приходилось подписываться.
Сочинение это содержит прямые  заимствования из московской бумаги, но
в целом оно вполне комакадемично, т.е. оно посвящено
 <<перестройке математики на марксистской основе>>.

Что касается списка задач обновленного Общества, то главная комакадемическая задача (которая чуть выше была 
была под номером 4) идет первой,
а всякая ерунда с развитием образования и поиску приложений для математики 
из списка выкинута%
\footnote{
	Где-то в бумаге есть  одинокая фраза <<содействие росту
	новых математических кадров из людей рабочего класса,
	которые могли бы с успехом работать в области естествознания и техники>>.}.
Обе декларации легко доступны в интернете (поэтому я цитирую выборочно).

\sm

{\sc Подытоживая.} При всех элементах дикости, 
{\bf московская} декларация (во всяком случае ее окончательный вариант)
 не могла быть продукцией авторов типа Кольмана, Хотимского или Яновской.
Она ничем не противоречит тому, что мы знаем о взглядах и о предшествующем и дальнейшем положении первых трех
подписавших ее людей, скорее она неплохо с ними согласуется. И, кроме того,
бумага содержит   элементы,  положительные как по отношению к комакадемической продукции
того времени, так и вообще положительные.

 \sm
 
 {\bf\punct Лузин в дни переворота 1930г.%
 \label{ss:luzin-perevorot}} Так или иначе, Лузин, закончив в Париже работу над 
 вторым своим главным трудом - <<Лекции по аналитическим множествам и их приложениям>>
 и учебником по анализу для инженеров, окунается в новую московскую жизнь. Егоров только что был арестован,
 в апреле был арестован А.Ф.Лосев, который приходился Лузину по крайней мере хорошим знакомым.
 Из Декларации мы знаем о заседании Совета Института математики единогласно исключившим из своего состава
 профессора Финикова (но, быть может, лишь из состава Совета).

 Автору почти неизвестно исторических источников 1930-31гг. о Лузине осенью 1930г., но многое говорилось на заседаниях Академической комиссии
 1936г. Дальше мы цитируем  Стенограмму:
 \begin{quotation}
 ХИНЧИН. 
Это причина [ухода Лузина из Университета], бывшая одной из главных в этой цепи, но не единственной. Я бы
сказал так, что все эти бурные в политическом отношении события, годы,
которые знаменуются процессом Промпартии и всем, что с ней было связано, в Институте математики
Московского университета переживались особенно бурно, так, как ни в одном из других
университетских Институтов не переживались.
 Обуславливалось это тем, что у нас
была исключительно темпераментная математическая молодежь, темпераментная в
политическом отношении. Я говорю, имея в виду, отчасти присутствующих здесь Лазаря
Ароновича [Люстерника] и Л[ьва] Г[енриховича] Шнирельмана, отчасти партийную
молодежь, как Хворостин
 и др. {\bf Была необычайно живая, напряженная и бурно политическая атмосфера} после долгих лет затишья в этом отношении.
 
ШМИДТ С МЕСТА. Впервые прозвучали политические речи.

...............

 ХИНЧИН. И на меня произвело совершенно определенное впечатление, которое остается незыблемым
 и до сегодняшнего дня, что на H.~H.~Лузина все это произвело впечатление абсолютно запугивающее: он испугался на всю жизнь, он дрожал от страха. В это
время произошел и арест Д.~Ф.~Егорова, и когда H.~H. выступал, то на его лице было написано, 
что началось светопреставление, что теперь всех нас заберут и т.д., и он дрожал.
Этот страх и дрожание остались у него и до настоящего времени...

 Достаточно вспомнить выступление Николая Николаевича на том заседании, которое
 	было посвящено процессу Промпартии — заседании Математического общества, которое велось под председательством Н.~Н.~Лузина.
 	Многие из здесь присутствующих —
 	присутствовали там и помнят как жалко, как несчастно он себя в этой роли чувствовал и
 	как он сам в качестве председателя, в конце концов, мог только промямлить, что «да,
 	конечно, вредить — это очень нехорошее дело, быть вредителем это очень нехорошо».
 	Ничего другого он не нашел сказать. Впечатление такой страшной обывательской запуганности%
 	\footnote{Очевидно, Хинчин считал страхи Лузина ни на чем не основанными.}.
 	Это и тогда уже у него чувствовалось. Тогда уже чувствовалось (это было в
 	30-м году все-таки) отношение к большевикам как к пришедшей какой-то страшной
 	силе, которая все разрушает, всех схватывает, арестовывает и т.д. Это дрожание за себя
 	чувствовалось во всем и, в частности, сказалось в этой истории с документом, которую
 	я тоже очень хорошо помню и где дело обстояло так, как рассказал Л[азарь] А[ронович].
 	Между прочим, насколько я помню, этот документ [о нем чуть ниже] стал составляться, и первые подписи
 	стали собираться%
 	\footnote{Об этом чуть ниже.} на том же заседании Математического общества, посвященного процессу Промпартии.

.................................. 	
 	
	ГЕЛЬФОНД. Так как у меня сохранилось довольно яркое впечатление об этом собрании в Университете, 
	которое было посвящено процессу Промпартии, то мне хотелось
	бы добавить кое-что к этому впечатлению. Я помню очень хорошо, когда была прочитана резолюция молодых работников 
	Института с рядом общеполитических установок и с
	конкретными выводами относительно Егорова, возглавлявшего тогда Институт математики, Лузин был председателем%
	\footnote{Очевидно, что это тоже самое собрание, вряд ли Лузин мог оказаться председателем
	где-либо, кроме Мат. Общества. Но Егоров уже Института не возглавлял.}.
	И действительно, в этом собрании он принимал, будучи
	председателем, крайне малое участие. Тогда особенно хорошо я помню впечатление,
	произведенное на него выступлением Аппельрота%
	\footnote{Аппельрот Герман Германович (1866–1943). Во время  свистопляски  ему было 74 года. 
	Все же он был не только <<очень старый человек>> в смысле Гельфонда (что его отчасти извиняло
	в глазах высокого собрания),
	но также и человек, чьи работы 125-летней давности
	и сейчас иногда  вспоминают (и, что, по-видимому, главное в современной околонаучной
	российской парадигме, цитируют;
	впрочем, его индекс Хирша по базам WoS \& Scopus маленький). В связи с ним стоит вспомнить
	еще один сюжет из московской математической истории.
	\newline
	Напомним, что есть классическая задача о движении твердого
	тела с закрепленной точкой в поле тяжести. Как известно, есть три случая, когда задача интегрируема
	(скажем, не уточняя смысла: явно решаема). Первый - случай Эйлера--Пуансо (Louis Poinsot, 1834). В этом случае закрепленная
	точка совпадает с центром тяжести. Это та же задача, что и движение в невесомости типа космического
	аппарата (что Пуансо вряд ли имел в виду) или вариаций  вращения планеты Земля. Второй случай -- волчок  Лагранжа, когда тело
	осесимметрично, а закрепленная точка находится на оси (это всем знакомая юла). Третий - странный случай
	Ковалевской, опубликованный лишь в 1888г. (полуоси эллипсоида инерции тела удовлетворяют условиям $a=b=2c$,
	а закрепленная точка находится в экваториальной плоскости эллипсоида). Это вызвало
	 большой интерес, вскоре было выяснено, что во всех других случаях задача не интегрируема. 
	Работа Ковалевской сильно заинтересовала  русских математиков и механиков, в частности, потому, что 
	в задаче вводилось комплексное время $t\in\C$, и в исследовании движения важнейшую роль играло
	поведение <<движения>> при невещественных значениях времени. В принципе, дифференциальные уравнения
	на функции комплексной переменной
	уже давно были в ходу, а работа Ковалевской делалась под влиянием Вейерштрасса и его работ.  
	То, что такие вещи  касаются классической 
	механики, вызвало в Москве удивление, а также интерес к комплексному анализу 
	(который, в частности, проявился у Жуковского) и к голоморфной теории дифференциальных уравнений.
	В конце XIX-начале XX века разные утверждения 
	о движении твердого тела получали, в частности,
	Аппельрот, Н.~Б.~Делоне-ст., Д.~Н.~Горячев, Жуковский, А.~М.~Ляпунов, Млодзеевский, П.~К.~Некрасов, Стеклов, Чаплыгин.  
	В конце  XIX-начале XX века были найдены несколько дополнений к классическим случаям, когда задача
	частично интегрируема,
	в этих новых <<случаях>> нужно  налагать условия не только на эллипсоид инерции и на положение закрепленной точки,
	но и на начальную скорость. Чаще всего упоминаются  в этой науке случай
	Гесса(W. Hess)--Аппельрота и случай Горячева--Чаплыгина. Про Аппельрота известно, что он с 1895г. работал
	в Московском сельскохозяйственном институте, в 30е годы он еще продолжал публиковаться. Видимо, последняя  его
	 работа -- 100-страничный мемуар <<Не вполне симметричные тяжелые гироскопы>> по движению твердого тела
	в сборнике памяти Ковалевской, изданном Чаплыгиным и Н.И.Мерцаловым, 1940. Там вместо поля тяжести говорилось
	об ускоренном движении гироскопа, что, видимо, было стилистикой предвоенного 1940г.}.
	
	С МЕСТА. Кто это такой?
	
	ГЕЛЬФОНД. Да это очень старый человек. Выступление было следующее: расценивая
	определенным образом то, что делается в это время, Аппельрот выступил и заявил, что
	{\bf «Николай Николаевич, в наше смутное время наша задача пронести светоч науки сквозь
	мрак того, что делается. И так как вы возглавляете это дело, на вас падают последствия»%
	\footnote{Отдавая должное словам Аппельрота, стоит вспомнить на то, что среди прочих
		обвинений Александрова и Колмогорова в 1936г. было то, что Лузин не оказался среди
		тех, кто в 1932-33гг. восстанавливал нормальную математическую жизнь, и, скорее, им мешал.}.}
	
	Я помню, как Лузин опустил голову и молчал.
	
	ЛЮСТЕРНИК. Да, тяжелое положение было председателя.
\end{quotation}

\sm

{\bf \punct <<Это не обращение Наркоминдела, а обращение общественности>>.%
\label{ss:prompartiya}}
Кроме отсветов об этом заседании, мы кое-что знаем о событиях ближайших дней.
Дело Промпартии вызвало возмущение во французских научных и инженерных кругах, а передовая советская 
математическая общественность
возмутилась против этого возмущения и составила возмущенное письмо. Письмо, как положено, начали подписывать,
и его авторы пытались получить подпись Лузина:
 \begin{quotation}
ГЕЛЬФОНД.
Теперь относительно ухода [Лузина] из Университета. Дело было таким образом. Был составлен
и согласован с ВОКС'ом [Всесоюзное общество культурной связи с заграницей] и Наркоминделом текст этого обращения, все его подписывали.
	Ближайшие ученики Лузина сначала сомневались, потом подписали. Когда дело
	дошло до Лузина, то к нему была отправлена аспирантка — теперь преподавательница
	Казанского университета — Рабинович%
	\footnote{Из воспоминаний Понтрягина:
\newline	
{\it	В начале нашего знакомства он [Левшец] пригласил нас с мамой в США на один год. 
	Кажется, в 32-м году я получил уже это приглашение, но из него ничего не вышло. Меня не пустили.
	Очень лёгкие до этого поездки за границу советских математиков стали к этому времени уже труднее.
\newline	
К отказу в поездке мне, по-видимому, приложили руку моя приятельница по университету студентка Виктория Рабинович 
и наша преподавательница философии Софья Александровна Яновская. Во всяком случае, однажды Яновская сказала мне:
\newline
— Лев Семёнович, не согласились бы Вы поехать в Америку с Витей Рабинович, а не с матерью?
\newline
Я ответил Яновской резким отказом, заявив: «В какое положение Вы хотите поставить меня?
Кто мне Витя Рабинович? Она же мне не жена».
\newline
Такая совместная поездка в Америку на год с Витей Рабинович могла бы кончиться браком с ней,
к чему я вовсе не стремился. Яновская в то время была влиятельным партийным деятелем, и я могу себе представить,
что от неё многое зависело, в частности, если она предлагала мне поехать с Витей Рабинович, 
то она, вероятно, имела основания думать, что может организовать эту поездку. Но я на это не согласился.
\newline
Так намечавшаяся на 33-й год поездка в Соединённые Штаты на год не состоялась.}
	}
	. Сколько раз стучала она Лузину, но ей не открывали и каждый раз ее просто не пускали.
	Так он и не подписал. После этого он ушел из
	Университета. У всех было впечатление, что поводом к уходу из Университета было то,
	что французским ученым физики и математики написали письмо.
	
	ШМИДТ. Это было воспринято, вероятно, как последствие процесса Промпартии.
	
	ГЕЛЬФОНД. Это было воспринято так, что он, несомненно, французских математиков
	ставит выше. Что касается процесса Промпартии, то там мало говорилось о Промпартии, но говорилось главным образом об интервенции.
\end{quotation}

\begin{quotation}
ЛУЗИН.
	Это неверно, это просто категорически неверно. Дело в том, что ко мне звонили по
	телефону, если не ошибаюсь, звонил Л.~А.~Люстерник. На это я ответил следующее: документ я подписывать не отказываю[сь], 
	но ввиду крайне сложной дипломатической
	обстановки могу это сделать только тогда, когда будет разрешение Наркоминдела. Я просто сказал: 
	пусть ко мне будет звонок из Наркоминдела (конечно, не от самого Литвинова). Пусть будет это сделано, потому что я занимал важное место. Я сказал, что не сделаю политической ошибки,
	но пусть будет звонок, и моя подпись немедленно будет на
	этом документе.
	Может быть, я ошибался, но, во всяком случае, Люстерник не откажется подтвердить,
	что именно так было, т.е. что я согласен, но только после разрешения Наркоминдела, или
	оттуда идущего.

	..........
	
	ЛЮСТЕРНИК. Никаких санкций не могло быть, это не обращение Наркоминдела, а обращение общественности.
	{\bf Мы пошли в Наркоминдел с самого начала. Мы написали
	текст более резкий.} Они сказали, что до иностранных ученых это не дойдет, т.е. не дойдут
	некоторые места. Это была консультация, но не санкция, а консультация с одним из работников Наркоминдела.
	{\bf Этот документ не от Советского государства шел. Это был как
	бы частный, нет, не частный, а общественный документ — обращение к своим товарищам по профессии из другой страны.}
	
	ЛУЗИН. Я считаю, что Вы оговорку сделали по существу: именно — частный. Я оценивал, 
	что это был чересчур частный документ. Но если бы я знал, что он согласован с
	более или менее ответственным работником Наркоминдела, я бы абсолютно подписал.
\end{quotation}

\begin{quotation}
	ЛУЗИН. Ведь я, буквально, только что приехал. Я не мог еще обжиться в нашей стране.
	Ведь прошло только две недели, как я приехал.
\end{quotation}

Надо отметить, что Лузин, имевший большие связи во французских научных кругах,
в самом деле был особой фигурой, а Франция скоро (очень скоро) понадобилась в качестве союзника.
Общественные связи с Франции (и без того не обширные) рвать не стоило... С точки
зрения интересов советского государства прав был Лузин.

\sm

 %%%%%%%%%%%%%%%%%%%%%%%%%%%%%%%%%%%%%%%%%%%%%%%%%%%%%%%%%%%%%%%%%%%%%%%%%%%%%%%
 
 {\bf\punct Бегство Лузина из МГУ.%
 \label{ss:luzin-begstvo}}
 Так или иначе
 с 20.11.1930 Лузин  устраивается на работу в ЦАГИ, а 
  16.12.1930 - уходит из МГУ (по собственному желанию). Кажется, его легко можно понять.
  Есть два опубликованных письма, которые дают дополнительные штрихи к этой картине.
 
 \begin{quotation}
 Александров (С борта парохода Deutschland) -- Колмогорову, 1 июня 1931г.
 
  Николай Николаевич [Лузин], читающий теорию функций у Кагана — 
явление, в полном смысле слова, трагикомическое. Я вполне понимаю,
что Н.~Н. тяжело и невозможно было оставаться ни в 
Университете, ни в Институте ({\bf ведь и для меня мое дальнейшее пребывание
в Институте продолжает оставаться нерешенным вопросом}, 
независимо от каких бы то ни было окладов), 
но, все-таки, на какое 
лихо студентам Кагана нужны трансфинитные и строжайше, со 
всеми мыслимыми тонкостями изложенные иррациональные числа?
Я, конечно, представляю себе, что среди слушателей Кагана 
найдутся молодые люди, серьезно любящие математику и 
интересующиеся ею, но, как бы то ни было, все это очень печально.

С другой стороны, что же бедному Николаю Николаевичу делать!
Я отлично понимаю его психологию — <<lch bin ja auch ein armer
Musikant>> [Я и сам — всего лишь бедный музыкант (Р. Шуман)], а 
жизнерадостности, о которой ты писал по поводу Миши Лаврентьева,
предполагать в людях, менее молодых, довольно трудно.
 \end{quotation}

 \begin{quotation}
 Колмогоров--Александрову, 17 марта 1931г.
 
  Был у Николая Николаевича. Основное, во что никто не верит:
Н.~Н. действительно болен и выглядит очень плохо. Склероз аорты
и еще что-то — с трудом взбирается на свою лестницу. Бросил
университет именно поэтому и тяжко обижен, что ни 
<<академические круги>>, ни общественность этому не желают верить.
 \end{quotation}

 %%%%%%%%%%%%%%%%%%%%%%%%%%%%%%%%%%%%%%%%%%%%%%%%%%%%%%%%%%%%%%%%%%%%%%%%%%%%%%%
 
 {\bf\punct Учебник не может быть признан вредительским...%
 \label{ss:luzin-vygodsk}}
 Лузинские напасти  на этом не кончились. В том же 1931г. Выгодский \cite{Vygodsky-GL}  опубликовал
 рецензию на только что опубликованный учебник анализа Грэнвиль--Лузин.
 Начинается она как вполне деловая рецензия, с разными положительными словами 
 и с деловой критикой. Дальше довольно мирно: 
 \begin{quotation}
 И не об отдельных недостатках 
изложения мы хотим говорить, а о принципах положенных в основу
изложения теории.

Мы глубоко убеждены в том, что учебное руководство, вводящее
в новую дисциплину с ее новыми методами и проблемами, должно
прежде всего удовлетворять требованию, чтобы учащемуся с самого
начала было ясно, откуда возникает новая проблема, почему
создаются новые методы и как они развивают. Иначе говоря, 
процесс вырастания абстракций на почве конкретной действительности
должен быть осознан. Это осознание необходимо как для того, чтобы
студент научился плодотворно применять математический метод на
практике, так и для того, чтобы он смог получить живой интерес
к изучению теории и ключ к ее глубокому пониманию.

Академик Н. Н. Лузин по-видимому стоит на позиции, 
принципиально противоречащей этой точке зрения.
\end{quotation}

Лузина меньше всего можно было бы обвинить в том, что он не пытался <<вырастить>> понятия 
в голове студентов. Но выращивает он их неправильно....
 
 \begin{quotation}
  В основе такой трактовки вопроса [о производной], столь обычной для многих
авторов, лежат несомненно определенные идеологические мотивы.
Здесь перед нами явно выраженная идеалистическая концепция:
понятие рассматривается как нечто первичное; его практическая 
значимость есть лишь «приложение» логического к конкретному; 
содержание понятия формируется не практикой, а заложено где-то в 
смутных глубинах сознания, и задача науки — извлечь его из этих глубин
на ясную поверхность рассудочного мышления...

Мы в корне не согласны с тем, чтобы определение дифференциала
давать чисто формально, игнорируя историческое прошлое 
дифференциала как своеобразной, актуально бесконечно-малой величины; но если
это уже делать, то конечно так, как это делает Н. Н. Лузин. В \S 160
он вполне последовательно проводит тенденцию, красной нитью 
проходящую через всю книгу.
 \end{quotation}
 
 Окончательный вывод прокурора был положителен для подсудимого:
 \begin{quotation}
 Вообще, стоя на исторической точке зрения, мы не можем считать
преступным замалчивание факта, которому не поверил бы Ньютон или
Эйлер,—-по крайней мере, если бы об этом факте им было бы лишь
упомянуто. Но, как было показано, антиисторизм характерен для
книги Лузина.
\end{quotation}
 
Однако каковы категории мышления! И где грань между матаном  преступным%
\footnote{Можно, конечно, попытаться предположить, что слово <<преступный>> 
 сказано в шутку, но было это в атмосфере,
где сотоварищи Выгодского аппелировали к авторитету ОГПУ, а оная организация приходила за знакомыми Лузина...
Вряд ли это шутка.}
и матаном всего лишь  антиисторическим и идеалистическим?

Стоит иметь в виду, что Выгодский  был директором НИИММ (до мая 1932г.), а с 1932г. крупным издательским
деятелем (и тогда от оного деятеля издание лузинского учебника зависело)

В том же 1931г. Выгодский
выпускает свой собственный учебник <<Основы анализа бесконечно-малых>>,
о нем еще пойдет речь ниже. На это Лузин отвечает восхищенным письмом%
\footnote{Датировка письма не ясна.}:
\begin{quotation}
 Глубокоуважаемый Марк Яковлевич,
 
позвольте искренне поблагодарить Вас за Ваш чудесный и ценный
подарок: за присылку мне Вашего Курса Анализа. Я давно 
слышал о его появлении и слышал о страстных спорах, возбуждаемых
им. По-видимому, совсем нет, или есть очень мало, лиц, спокойно
к нему относящихся. Всякий, знакомящийся с ним, становится или
горячим его поклонником, или столь же страстным его 
противником. И это меня не удивляет ничуть, так как Вы мужественно 
коснулись самой болезненной точки Анализа вообще, современного -- в
особенности, метнув камень в «осиное гнездо».
\end{quotation}

Письма Лузина к Выгодскому были изданы А.~И.~Маркушевичем \cite{Markushevich} в 1976г. и потом неоднократно переиздавались
\cite{Demidov-pismo}, \cite{Laugwitz}%
\footnote{Автор данных записок должен признаться, что при всех удивительных вывертах
Выгодского, он испытывает к нему как к автору и даже как  деятелю большое уважение.
Он был искренний умный  фанатик (а вовсе  не шарлатан вроде Кольмана или Лейферта), верящий в разум 
и оному служивший. Он не был безграмотен (как тогда Яновская).
Однако письма эти не вызывают
у автора данных записок умиления, которое должен был испытывать читатель этих публикаций: как сам великий Лузин
пишет какому-то молодому Выгодскому...}. 

Лузин, надо сказать, и в этом
непростом жанре был нетривиальным
автором. Он, в частности, утверждал, что строгое введение актуальных бесконечно малых
в анализ возможно (оно и в самом деле возможно, как показал Робинсон (Abraham Robinson)
в 1960г., и, может, до того Хьюиттом (Edwin Hewitt)). 

Цель Лузина была достигнута, общий язык с Выгодским найти ему удалось.
Например, Выгодский весной 1932г. включил Лузина в комитет по изданию <<Математической энциклопедии>>
(проект не был осуществлен). Позже Лузин искал в Выгодском
союзника против Александрова, см. письмо от 14.09.1934, опубликованное в \cite{ErmTok}.

\sm

 %%%%%%%%%%%%%%%%%%%%%%%%%%%%%%%%%%%%%%%%%%%%%%%%%%%%%%%%%%%%%%%%%%%%%%%%%%%%%%%
 
  {\bf\punct Переворот в МатСборнике.%
  \label{ss:sbornik-perevorot}} Впрочем, жизнь в Москве 
  шла вперед и вперед. В день переворота в Обществе была сформирована новая
  редакция <<Мат. Сборника>>. В 1931гг вышли всего два выпуска Сборника, 1-2 и 3-4.
  Первый из выпусков открывался статьей \cite{Sbornik-revolution}
  <<От редакции.
Со следующего номера журнал «Математический сборник» будет
выходить под названием «Советский математический сборник».>>
Там сообщалось о перевороте (эта информация выше уже использовалась).
Завершалась статья словами: 
\begin{quotation}
Настоящий номер, сданный прежней
Редакцией, не успел еще отразить реорганизацию Об-ва.
Приступая к работе, редакция «Советского математического
сборника» призывает советских математиков сплотиться вокруг журнала
 и помочь обратить его в {\bf боевой орган} советской математики.
 \end{quotation}
 На обложке сообщалось \cite{Tok-white}, что 
\begin{quotation}
... материалы следует направлять ответственному редактору Л.~А.~Люстернику и
ответственному секретарю А.~О.~Гельфонду.
   \end{quotation}

   В следующем номере слово <<Советский>> в названии не появилось, зато появился список новой редколлегии:
  \begin{quotation} 
   Редакция:
О.~Ю.~Шмидт (отв. редактор), П.~С.~Александров, А.~Ф.~Бермант (отв. секретарь), М.~А.~Лаврентьев, 
Л.~А.~Люстерник, В.~И.~Смирнов, Ф.~И.~Франкль, А.~А.~Холщевников, Н.~Г.~Чеботарев.
   \end{quotation}

   Редколлегия изменилась  по сравнению с неизвестной редакцией от 12.11.1930 (ответственный редактор
   теперь не Люстерник, а Гельфонд вообще выпал,
   Александров с середины 1930г по середину 1931г. находился за границей и едва ли мог входить в первоначальный
   список обновленной редакции).
   Из незнакомых нам лиц в списке присутствует Бермант Анисим Федорович (1904--1959) (который в тот момент
   только что закончил аспирантуру у Лаврентьева), человек, судя по всему, тоже радикальных взглядов, и Холщевников%
   \footnote{Я нашел у А.~А.~Холщевникова публикацию о <<Математических рукописях>> Маркса (Фронт науки и техники, 1933, 2, 100-106), ультрафинитистская статья в <<Вестнике Коммунистической академии>> за 1928г (вып. 27)
   а также доклад на Втором Всесоюзном
   математическом съезде по секции истории и философии математики. Доклад содержал какую-то муть
   на тему континуум-гипотезы.}. Смирнов и Чеботарев, очевидно, представляли
   Ленинград и Казань (в следующем номере сообщалось,
   что сборник издается Московским, Ленинградским и Казанским
   математическими обществами). Стоит сопоставить  представителей иных городов в Шмидтовском Сборнике и  в кольмановском
   Мат.обществе 1931г. (см. ниже п.\ref{ss:obshchestvo-1931}), там
   иногородние деятели представляли прежде всего политику, а здесь математику. Московская часть списка,
   если оставить загадочного А.~А.~Холщевникова, тоже выглядит нормально (А.~Ф.~Бермант, скорее всего,
   занимался делопроизводством, а в А.~А.~Холщевникове естественно предполагать политического представителя,
   кстати его выпадение из редакции в 1935 совпало с появлением там Б.~И.~Сегала).
   
   В этом выпуске также содержалась новое заявление:
  \begin{quotation}
   Советские математики, поддерживайте свой
журнал!

{\bf Среди большинства советских математиков сохранилась традиция печатать
свои лучшие работы в иностранных журналах}. Больше того, существовала и
пользовалась распространением точка зрения, усматривавшая в факте печатания
большого количества наших работ за границей положительное явление — {\it преодоление
советской наукой культурной блокады}. Этот взгляд, конечно, неправилен:
рассыпанная по журналам Германии, Франции, Италии, Америки, Польши и других
буржуазных стран советская математика не выступает как таковая, не может
показать собственного лица.

Рост научных кадров внутри СССР, поворот советской математики лицом
к социалистическому строительству ставят перед нами задачу создания журнала,
отражающего эти сдвиги и организующего советскую математику в направлении
активного участия в соцстроительстве.

Советская математика может и должна иметь журнал международного значения.
Поэтому мы продолжаем обычай снабжать иностранными резюме статьи,
написанные на русском языке, и печатаем статьи на иностранных языках.
Опыт показал, что и математические статьи, написанные на русском языке, доходили
до иностранного читателя.

{\bf Группа московских математиков} обратилась в редакцию с письмом, в котором {\bf принимает 
на себя обязательство печатать свои статьи, в первую очередь, в
«Математическом сборнике»} и призывает к этому других математиков Советского
Союза.
  \end{quotation}
 Курсив из оригинала. Жирный текст мой (мы эти идеи где-то уже встречали? не так ли?).

   Что касается <<группы московских математиков>>, то ее состав  понятен, это
   Люстерник, Шнирельман и Гельфонд. Они свое обещание выполнили и за границей после 1931г.   не публиковались
   (у Гельфонда была одна заметка в Compt.Rendus  1934г.).
   
   \sm

 Значительно интересней, что новая редакция (она в начале 1933г. пополнилась Выгодским, Хинчиным 
 и Степановым%
 \footnote{Александров (Гёттинген)--Колмогрову, 9 ноября 1932 г.
 \newline
 {\it
 Редакцию я себе представляю приблизительно в таком составе:
Гл. редактор — О.~Ю.~Шмидт.
Члены редакции: Александров, Выгодский, Степанов.
Секретарь редакции — Бермант.
(Обязанности ясны: О.~Ю. — для больших оказий, М.~Я. — связь с 
издательством и партийный авторитет, я и В.~В. — фактическое
руководство работой, А.~Ф.~Бермант — исполнение.).....}
  \newline
  Это личная точка зрения Александрова на то, как надо.
 Понятно, что в трех случаях обязанности такими и были,  
 а кто фактически руководил работой журнала, история умалчивает.}) 
 выполнила другое свое обещание и превратила Мат.сборник
 в один из лучших мировых математических журналов (правда в редакции остался лишь один трех
 перечисленных молодых математиков -- Люстерник).
   
 %%%%%%%%%%%%%%%%%%%%%%%%%%%%%%%%%%%%%%%%%%%%%%%%%%%%%%%%%%%%%%%%%%%%%%%%%%%%%%%
 
\sm
 
 %%%%%%%%%%%%%%%%%%%%%%%%%%%%%%%%%%%%%%%%%%%%%%%%%%%%%%%%%%%%%%%%%%%%%%%%%%%%%%%

 {\bf\punct Падение Шмидта.%
 \label{ss:schmidt-padenie}} Мы видим, что Шмидт в 1931г. стал главредом Мат.Сборника,
 добавив еще одну строчку к списку своих должностей. В действительности, в самом начале 1931г.
 Шмидт потерпел тяжелое поражение.

 В  президиуме Комакадемии 23.12.1930--6.1.1931 состоялась дискуссия
 <<О положении на фронте естествознания>>.
 Материалы <<заключительного звена первого этапа>> дискуссии были опубликованы в книге \cite{Za-povorot}.
 Как явствует из введения, это не полная стенограмма.

 Докладчик -- Шмидт.
 Содокладчик  А.~А.~Максимов. Незадолго до этого в Комакадемии
 происходили бурные события, связанные с обвинениями в адрес группы Деборина. Из речей Шмидта и Максимова видно, что дискуссии предшествовала какая-то резолюция
 в <<правлении ассоциацией>> (естественных наук). По-видимому, к этой резолюции относились слова одного из
 выступавших (Новинский)
 \begin{quotation}
 И 
 нужно сказать, что вся эта самокритика, которую здесь О.~Ю.~Шмидт развернул, явилась только заключительным аккордом нашей 
 дискуссии и еще за 2-3 минуты до того, как О.~Ю.~Шмидт решил голосовать за то, что позиция 
 естественно-научного руководства, в частности его позиция, является антимарксистской, он решительно
 возражал против такой формулировки [и т.д.]
 \end{quotation}
 
 Максимов развивает суровую товарищескую
 критику Шмидта. В его длинной речи были, в частности, такие слова:
 \begin{quotation}
  Если к этому добавить, что по линии БСЭ, где подавляющее большинство статей 
  оказываются немарксистскими, или даже ан\-ти-марк\-си\-ст\-ски\-ми (в таких статьях как В и т.д. <<Вещество>>,
  <<Внушение>>, <<Аксиома>> и т.д. -- здесь даны не только неправильные, но прямо антимарксистские формулировки), все это показывает, что нам еще придется поработать в этом направлении,
  чтобы вскрывать эти ошибки, чтобы производя положительную работу, в то же время эту критику 
  продолжать на основе самокритики....
 \end{quotation}

Дальше идет несколько выступлений, кто-то обвиняет И.~И.~Агола в связи с БСЭ. 
А потом выступает товарищ Яновская:
 \begin{quotation}
  ... {\bf получалось из его Шмидта высказываний и на философской дискуссии такое впечатление, что
  он является только жертвой.} Но мне кажется, что дело обстоит все-таки не так. И если мы посмотрим
  на ту конкретную работу, которую вел т.~Шмидт, а мне пришлось это видеть в последнее время на
  его работе в БСЭ, где ему принадлежат некоторые статьи, то мы убедимся в том, что {\bf т.Шмидт в
  этом отношении вовсе не является жертвой, что у него неправильная идеалистическая позиция, скажем, в отношении
  математики}...
 \end{quotation}

Там же она выражает восторг учебником статистики Хотимского (см. выше п.\ref{ss:hotimski}).
 Дальше мы цитируем Яновскую по ее статье \cite{Yanov-BSE}, опубликованной в 1931г.:
 \begin{quotation}
 {\bf Казалось бы при таких условиях основной задачей наших 
марксистских авторов и изданий в отношении математики должна была быть 
беспощадная борьба с математическим идеализмом, где бы и в какой бы
форме он ни проявлялся, неуклонное разоблачение его подчас очень
тонких хитросплетений, неусыпная работа над установлением и в 
области математики подлинного единства теории и практики на основе
борьбы за осуществление генеральной линии партии на основе 
практики социалистического строительства. Казалось бы в таких 
условиях все наши работы и издания, где речь идет о математике, должны
выть боевыми, воинствующе партийными.} Казалось бы в первую
очередь это относится к такому, рассчитанному на огромную аудиторию
активно участвующих в строительстве социализма трудящихся СССР
изданию, как Большая советская энциклопедия.....

{\bf Между тем в БСЭ нет ни одной боевой
воинствующей, партийной статьи по математике.} Более того, с ее страниц
иногда прямо, иногда в несколько замаскированном и 
«подчищенном» виде пропагандируется махизм, конвенционализм, идеализм.

Уже самый состав сотрудников БСЭ свидетельствует о том, что
редакция поставила себе в отношении естествознания и техники более
чем скромные цели. К работе в этих разделах почти совершенно не
были привлечены молодые кадры советских ученых; она строилась
исключительно в расчете на старую «левую», либеральную 
профессуру, марксистски совершенно неподготовленную, «беспартийную»
в своей политической, а следовательно и философской установке,
воспитанную на Канте, Махе и Пуанкаре и в лучшем случае 
пытающуюся «соединить» махизм с материализмом. Развернутое 
наступление социализма в СССР и связанное с ним бешеное сопротивление
кассового врага раскрыло кавычки этой пресловутой 
«беспартийности»: значительная часть основных сотрудников и редакторов БСЭ
оказалась по ту сторону баррикады, в лагере контрреволюционного
вредительства, в связи с международной буржуазией, 
социал-фашизмом и кулацко-эсеровской контрреволюцией...

{\bf Неудивительно, что Энгельс и Ленин просто выпали из статей
этого отдела [математики], что отдел в целом не может быть квалифицирован иначе
как антимарксистский и идеалистический, несмотря на то, что именно
тут бой с идеализмом должен был быть особенно жесток}...

....конвенционализмом, идеализмом в области естествознания и математики.
Перед нами статья главного редактора БСЭ тов. Шмидта об алгебре...
[дальше много страниц сурового философского анализа, попутно достается проф. Кагану...]
 \end{quotation}
 
 Тема БСЭ принимает в дискуссии особое значение...
 В заключительном слове товарищ Максимов говорил, в частности:
 \begin{quotation}
Борьба, по-видимому, будет дальше продолжаться, то, что, мы проделали по БСЭ, о чем скажет, очевидно О.~Ю.~Шмидт, показало, что нам  много 
 	 еще придется бороться. Мы проделали, большую конкретную работу по БСЭ, но это показывает, что мы
 	только маленький кусочек работы сделали. Нам придется  здесь
 	критиковать и бороться дальше. 
 \end{quotation}
 
 В своей заключительной речи Шмидт, разворачивая (как это было положено) самокритику, все еще пытается отбиваться.
 Принимается длинная резолюция, осуждающая <<естественно-научное руководство>>,
 в составе Шмидта, а также Левина, Левита, Агола, Гессена и еще Серебровского.
 Там есть и пункт, посвященный БСЭ.
 
 Решение Комакадемии -- одного из мыслительных центров Партии
 подтверждается постановлением ЦК. 

\sm

 Кроме того, в журнале <<За марксистско-ленинское естествознание>>
 -- так был переименован  журнал <<Естествознание и марксизм>> (как утверждается в \cite{Bogolyubov-Rozhenko}
  основанный Шмидтом)
 --
  была напечатана 
 статья
 <<Объединенное заседание правления Ассоциации
 естествознания Коммунистической академии, руководящего состава естественно-научного Отделения
 ИКПФиЕ и главной редакции БСЭ (сокращенная стенограмма)>>
\cite{Zasedanie}.  Порядок следования этих двух дискуссий не ясен,
скорее всего вторая дискуссия была в перерыве первой 
(и к ней относится вышеприведенная цитата о проделанной работе в заключительном слове Максимова).
 
 Здесь после вступительного слова Шмидта начинается общая атака на БСЭ,
 сначала Максимов, потом Яновская, которая в качестве знатока математики 
 с особой решительностью обличает лично Шмидта (изрекая мутный поток 
 марксистствующего богословия). Далее  Баткис и Ландис по медицине,
 Климовицкий по биологии, Балезин по химии, Сурта на тему о вообще%
 \footnote{Участники дискуссии: Баткис Григорий Абрамович,
 	(1895-1960), Ландис Михаил Моисеевич, Балезин  Степан Афанасьевич(1904-1982),
 	Климовицкий Виктор Абрамович (1899 – 1965).
 Сурта Иван Захарович (1893-1937).}. Собрание приняло резолюцию. Приведем отрывок (отражающий тонкости и переливы комакадемической мысли):
\begin{quotation}
	\small
1. Имели место случаи использования БСЭ для пропаганды идеологии враждебных классов, идеализма, махизма, витализма и т.д., напр. слова <<Антропология>>, <<Биология>>, <<Волны>>, <<Воля>>, <<Внушение>>, <<Вероятность>>.

2. В ряде случаев пропагандируется антимарксизм и меньшевиствующий идеализм. Например, слова <<Аксиома>>, <<Влечение>>, <<Вейсман>>, <<Гаусс>>, <<Геометрия>> и др.

3. Слова <<Гельмгольц>>, <<Гарвей>>, <<Вещество>>, и др. представляют явную пропаганду механистического мировоззрения.

4. {\bf  В целом ряде статей отсутствует четкая классовая линия <<Алгебра>>}, <<Галилей>> и др.

5. Как правило отсутствует использование работ Маркса, Энгельса, Ленина по кардинальнейшим вопросам естествознания (все статьи по математическим вопросам, ст. <<Восприятие>> и др.)

6. Для писания статей практиковалось привлечение немарксистов и нематериалистов, а также буржуазных ученых и выступающих против марксизма и ленинизма (привлечение Франка Ф., наследовавшего кафедру Маха).

7. Руководителями крупнейших подотделов,
отдела техники, являлись такие фигуры, как Рамзин, Осадчий, Чарповский, оказавшиеся вредителями и контрреволюционерами. 

8. Приглашавшимся специалистам не давалось определенного задания, четкой установки при написании той или иной статьи. Каждый из них излагал свое собственное мировоззрение, и редакция не всегда вносила существенные идеологические поправки в статьи, а многие из вносимых поправок являлись далеко недостаточными.

9. Привлечение имеющихся партийных сил для написания статей было совершенно недостаточно.

10. Практиковалось исключительно индивидуальное выполнение заданий без привлечения коллективов, без предварительной коллективной критики помещаемых статей. 

{\bf Совещание} считает, что совокупность указанных ошибок является проявлением правого оппортунизма на практике и {\bf признает работу БСЭ в области естествознания в его теперешнем виде антимарксистским по существу.}
\end{quotation}

\sm

 Главным редактором БСЭ Шмидт, как это ни странно после такой головомойки, остался. 
 Вместо Кагана редактором естественно-научного раздела становится А.~А.~Максимов,
 и на другом месте того же Кагана -- редактора отдела математики -- оказывается Кольман.
 БСЭ была украшена обширной статьей Кольмана и Яновской <<Дифференциал>> т.22, 602-608,
 где 4/5 текста были посвящены Марксу, а также Энгельсу и Ленину.
 
 Список  должностей Шмидта (см. п. \ref{ss:schmidt})  в тот год укоротился втрое.
 Он больше не член Президиума Комакадемии, не заведует секцией естественных наук Комакадемии,
 не является редактором трех научно-марксистских журналов (но оказывается редактором 
 Мат. Сборника). Не является он и директором НИИММ... 
 В тот год он даже не возглавлял полярных экспедиций (но директором Арктического
 института все же остался, как и членом Президиума Госплана). 
 В следующем, 1932г,
  после перехода <<Сибирякова>>
 из Архангельска в Петропавловск за одну навигацию Шмидт
 начнет становиться  всесоюзной знаменитостью.
 В тот же год при Совнаркоме
 появится управление ГлавСевМорПуть%
 \footnote{В 1932г. <<Сибиряков>> и <<Русанов>> также забрали с Северной Земли группу  Г.~А.~Ушакова и Н.~Н.~Урванцева.
 Они исследовали
 в 1930-1932гг.
 архипелаг и составили его карту, что само по себе было крупнейшим географическим достижением
 (до того о Земле Тайвай, она же Земля Николая II, она же Северная земля, была известна теорема существования,
 установленная Б.~А.~Вилькицким-мл. и его экспедицией в 1914-15гг.).
 Обеспечение автономного существования четырех человек в течение двух лет в Арктике требовало определенных средств.
 В 1930г. выделенного финансирования не хватало, Ушаков  тогда взял кредит под шкуры 100 неубитых медведей.
 После 1932г. выделяемые на Арктический проект средства станут совсем иными. Скоро  под руководством Шмидта
 в Арктике будет выстроена инфраструктура с портами, авиабазами и метеостанциями.
 В 1941г. СевМорПуть окажется важной экономической и военной коммуникацией.
 \newline
 Стоит отметить, что в истории бывали и другие математики с не только академическими биографиями, например, Лазар Карно (Lazare  Carnot)
 и Поль Пенлеве (Paul Painlev\'e).}.

 А статья Кольмана и Яновской про дифференциал осталась чудо-юдом, так и не став <<античным образцом>> 
 для всеобщего подражания, о чем ниже.
 Но мы забежали вперед, пока же
 на фронте наступления Комакадемии начинался  новый, 1931год. Вместо Шмидта
 секцию естественных наук Комакадемии возглавляет
 Кольман, и вместо другой ипостаси того же Шмидта директором НИИММ оказывается Выгодский.
 
 \sm

 {\bf \punct <<За материалистическую диалектику в математике>>.%
 \label{ss:fon}}
 Это название сборника, вышедшего под редакцией Кольмана в 1931г., который отчасти состоял из
 ранее опубликованных статей. Он уже многократно цитировался выше и будет цитироваться.
 Это квинтэссенция свистопляски 1931года, желающие могут раскопать эту 
 дикую книгу и попытаться ее читать,
 мы лишь перечислим авторов с указанием числа статей (если таковая не одна), 
 
 {\small
 Боярский 
 
 Гливенко
 
 Выгодский (7)
 
 Кольман (8)
 
 Люстерник
 
 В.Молодший
 
 Орлов
 
 М. М. Пистрак
 
 Смит (2)
 
 Ставровский
 
 Хотимский (2)
 
 З.Цейтлин
 
 Цыпкин 
 
 Яновская (2)
 
 Ястремский 
 }
 
 %Кроме того была статья непонятного автора, <<Об основных
 %вопросах преподавания математики в средней школе>>, которая была удалена из доступной автору отсканированной
 %копии
 %(интересно бы ее почитать, она могла бы пролить свет на некоторые дальнейшие события).
 
 В этом списке интересно видеть Люстерника и Гливенко со статьями <<О топологических методах анализа>> и
 <<Интеграл в математике и математическом естествознании>>. Следует подчеркнуть, что это нормальные
 популярные статьи по математике. Автору последнее обстоятельство представляется весьма существенным,
 и в связи с этим стоит еще раз напомнить о
 <<Декларации инициативной группы>>.
 
 \sm
 
 %%%%%%%%%%%%%%%%%%%%%%%%%%%%%%%%%%%%%%%%%%%%%%%%%%%%%%%%%%%%%%%%%%%%%%%%%%%%%%%

\sm
 
 {\bf\punct Бригадно-лабораторный метод.%
 \label{ss:brigada}}
 Про (ди)ректора МГУ Касаткина Николая Васильевича (1891 – 1961), 
 пришедшего на смену Удальцову, известно мало, наиболее подробно о 
 его деятельности в МГУ
 пишет Ильченко в биографии Петровского, а биография самого Касаткина есть в  \cite{Korsakov}. 
  В целом, его история  ректорства (3.7.1930-22.7.1934) выглядит как цепь реорганизаций и потрясений,
 что, впрочем, под напором сверху происходило по всей системе образования СССР.
 Похоже на то, что человеком он был дельным,
 а Университет (точнее та часть, которая от него осталась после всех отделений) 
 уже полуразгромленный к середине 1930г.,
 вышел из потрясений, начавшихся в 1929г., в нормальном виде (быть может, и более сильном, чем в них входил),
 разумными выглядят и кадровые назначения 1932-33г.

 На 1930-1931гг. пришлась очередная напасть. Из приводимых в \cite{Ilchenko} тезисов доклада Касаткина
 (доклад был уже позже описываемых событий):
 \begin{quotation}
  Весной в массовом масштабе осуществлён переход на «активный бри\-гад\-но-ла\-бо\-ра\-тор\-ный метод занятий», введённый в 1930 г.
  Академические группы разделены на бригады по 3–5 студентов, которые без систематического руководства квалифицированных 
  преподавателей должны были заниматься «проработкой» учебных заданий. Этот метод почти полностью ликвидировал лекционные курсы,
  читаемые крупными учёными. Чтение общих и специальных курсов заменено вступительными лекциями
  к самостоятельной работой студенческих бригад, а экзамены и зачёты как индивидуальная форма оценки знаний студентов – 
  коллективным отчётом бригад о «проработке» ими учебных заданий. Творческий подъём, энтузиазм охвативший студентов, 
  аспирантов и научных работников в годы первой пятилетки, некоторое время скрывал несостоятельность 
  «бри\-гад\-но-ла\-бо\-ра\-тор\-но\-го метода» занятий. Но уже в первый год применения этого метода были отмечены 
  серьёзные недостатки в преподавании: слабое участие профессоров и преподавателей в учебной работе,
  низкая успеваемость студентов и т.д. Сторонники «метода» объясняли их его новизной%
  \footnote{Как знакомо... Это звучало в 1970е в связи с новой
  	школьной программой по математике,
  это звучит в последнее десятилетие в связи с ЕГЭ.
  Различие в том, что тогда команда <<отбой>> была дана сравнительно быстро.},
  неумением его применять,
  а не следствием порочности самой основы этой системы преподавания в целом....
 \end{quotation}

Вот описание этого метода в статье \cite{Lapko}:
\begin{quotation}
 Чтение лекций было отменено, вся тяжесть проработки
материала была перенесена на самостоятельные занятия студентов, но не
каждым в отдельности, а бригадой. Был введен так называемый
 бри\-гад\-но-ла\-бо\-ра\-тор\-ный метод. Этот метод был основан на механическом перенесении
норм организации труда из промышленности в высшие учебные заведения.
Копирование в учебной работе бригадной организации труда на производстве
преподносилось под видом «борьбы» за качество подготовки пролетарских
специалистов. Бри\-гад\-но-ла\-бо\-ра\-тор\-ный метод предполагал отмену лекционной
системы обучения — основного способа ознакомления студентов с теоретическими вопросами. Лекции допускались лишь в небольшом количестве.
Всемерно поощрялась групповая, бригадная проработка материала без преподавателя
путем кратких докладов членов бригады по отдельным частям
задания или путем беседы о прочитанном между отдельными студентами.
Весь объем материала каждым студентом не изучался. Индивидуальный
контроль преподавателя за работой студентов отсутствовал, так как существовала
форма групповой коллективной беседы в часы занятий. Кроме
того, этот метод проработки теоретического материала предполагал наличие
большого количества разнообразных учебников. Обеспечение же студентов
литературой (в частности, математической) было очень плохое как в качественном,
так и, особенно, в количественном отношении; часть студентов
пользовалась устаревшими пособиями, выпущенными еще в прошлом веке,
да и их было очень мало. 
\end{quotation}

Есть упоминания попытки Наркомпроса ограничить применение этого метода летом 1931г.
Тем не менее, он продолжал внедряться. Далее его вроде бы 
отменили осенью 1932г. Но он не везде отменился, и было
отдельное постановление ЦК уже 1936г. об отмене его там, где он еще оставался. 

Есть такое известие о жизни профессуры в это время:
\begin{quotation}
Колмогоров -- Александрову, 16.04.1931

 К сожалению, особенной разгрузки от 40—60 часов в декаду
наших математиков при этом ожидать не приходится! В 
университете в будущем году организуется 12 групп (!), изучающих 
дифференциальные уравнения.
\end{quotation}

Отметим, что 60 часов в декаду при пересчете на рабочую неделю дают 40 часов.
Нагрузка, превосходящая невыносимую. По-видимому, все же имеется в виду работа в двух
или большем числе мест работы. Еще я не понимаю, как можно было набрать эту нагрузку при бри\-гад\-но-ла\-бо\-ра\-тор\-ном методе.

\sm

 %%%%%%%%%%%%%%%%%%%%%%%%%%%%%%%%%%%%%%%%%%%%%%%%%%%%%%%%%%%%%%%%%%%%%%%%%%%%%%%

 {\bf\punct Москва помогает Ленинграду.%
 \label{ss:moskva-leningrad}} Между тем, в Ленинграде продолжалась борьба 
группы Лейферта против Гюнтера и прочих <<правых>>.
Об этом повествовалось в уже упомянутой книге <<На Ленинградском математическом фронте>> \cite{Lenmatfront}.
События эти изучалось в статье Н.~С.~Ермолаевой \cite{Erm3}.
Победный конец описывается в \cite{Lenmatfront} так:
\begin{quotation}
 Большую идеологическую и организационную помощь в этой работе
об-ву математиков-материалистов оказал приезд из Москвы представителей секции математики Коммунистической академии тт. В.~И.~Хотимского и
С.~А.~Яновской.

8 февраля [1931г.] состоялся доклад т.~Яновской, организованный Обществом математиков-материалистов,
«О задачах математики в реконструктивный период», 
и с большой речью в прениях выступил В.~И.~Хотимский. На этом собрании профессором Г.~М.~Фихтенгольцем было сделано
публичное заявление о его желании работать на новых путях в
интересах социалистического строительства, и это заявление, как и выступление проф. Б.~Н.~Делоне, были правильно учтены как выступления
не только личного характера, но как отражающие настроение целого
круга математиков Ленинграда.

10 марта инициативная группа математиков в составе: академика
И.~М.~Виноградова, профессоров А.~Р.~Кулишера, Г.~М.~Фихтенгольца,
Б.~Н.~Делоне и других приняла декларацию, послужившую основой
организации нового Общества математиков Ленинграда.

В то же время появилось письмо проф. Н.~М.~Гюнтера с признанием
части своих главнейших ошибок.
\end{quotation}

Сколько-либо заметных научных фигур среди ленинградских лево-радикалов не было. В 1932г.
Виноградов становится директором Физико-математического института (им. Стеклова).
Точная дата небезынтересна (тогда воды за единицу времени утекало много), но история ее умалчивает.

На должности директора НИИ математики  при Университете в следующие годы
до 1936-37
заметен%
\footnote{Назначен 6.10.1930 \cite{Mihlin}. Или, например:
\newline{\it
Александров--Колмогорову (17 февраля 1936 г.,из Москвы
в Тифлис)
\newline
Президиумом Академии Наук ты утвержден членом
математической группы Отделения без предварительного рассмотрения
твоей кандидатуры в группе, в отличие от Гельфонда и от
Лаврентьева, которые будут обсуждаться. Очевидно, разница
вызвана тем, что ты есть директор Института} [НИИ математики при МГУ]. {\it Однако твой
коллега А.~Р.~Кулишер  в состав группы не введен. Зато введен
академик [АН Белорусской ССР; тоже директор института] Бурстин [Burstin].
Впрочем, заседание группы, как говорят, 
произойдет в ближайшие дни, и там мы выберем нашего великого
Гельфонда} [в математическую группу, член-кором Гельфонд стал в 1939г.].
\newline
Забавно, что по такой мелочи приходится заниматься медиевистикой и извлекать информацию путем умозаключений
из  косвенных источников.}
Кулишер Александр Рувимович (1875 — ?), один из  ленинградских лево-радикалов.
Он часто значится среди репрессированных в 1937г., но в мемориальских
списках его нет, а в \cite{Nagibin}
он упоминается среди <<выбывших по разным причинам>> из Вятского пед.института
в 1943-44гг.%
\footnote{Существовал также Кулишер Роман Маркович (1903—1935), вроде бы живший 
в одном доме с упоминаемым Кулишером, и имевший несчастье быть женатым на
Ольге Петровне Драуле, сестре роковой Милды. Впрочем, наш
Кулишер и после известной истории продолжал оставаться директором института.}.
В 1937г. директором института становится В.~И.~Смирнов (один из недавних <<правых>>).
О том, кто заведовал ленинградской частью Института Стеклова в 1934-1940гг.
история почему-то (и едва ли без причин) умалчивает.

\sm

%%%%%%%%%%%%%%%%%%%%%%%%%%%%%%%%%%%%%%%%%%%%%%%%%%%%%%%%%%%%%%%%%%%%%%%%%%%%%%%
 
 {\bf\punct Директор Выгодский.%
 \label{ss:vygodsk-dir}}
Директорство Выгодского в Институте механики и математики замечательно
тем, что о нем не принято было упоминать, ни вообще, ни в биографиях Выгодского.
Эта традиция продолжалась и в 90е и в нулевые годы, и лишь совсем недавно
оно всплыло на свет Божий. Что именно такого интересного 
тогда творилось (а очевидно, что просто так предмет табуирован быть не может),
автор не знает. Вот три случайный отсвета:
 \begin{quotation}
 Колмогоров--Александрову (8 марта 1931)
 
 В январе М.Я. Выгодский 
назначен директором Института,
а А.Я. — заместителем. Но
Выгодский в тот же день скрылся
в неизвестном направлении.
Вместо него был оставлен 
Барсуков%
\footnote{Барсуков Александр Николаевич (1891-1958) -- крупный педагог.
Революционер, член РСДРП(б) с 1917г.,   был деканом Физмата МГУ, замректора,
главредом Учпедгиза. Основатель журнала <<Математика и физика в школе>>(1934),
позднее <<Математика в школе>> (1937).
Автор стабильного школьного учебника по алгебре, бывшего в ходу в 1956-1966.
}, которым и было дано
предписание мне вернуться. На
днях М.Я. приехал и скоро 
вступит в должность.
 \end{quotation}
 
 \begin{quotation}
  Оргдеятельность же Выгодского и его заместителя (зав. учебной
частью Берманта) производит впечатление хорошее: они отдают
Институту много времени и стараются все наладить самым 
деловым образом. Разгруппировывают вновь полученных аспирантов по
степени подготовки и не боятся заставить неграмотных изучать
дифференциальное исчисление. Отстаивают и достаточно 
теоретические специальности, хотя и намерены направлять на них 
небольшое число хорошо подготовленных аспирантов.
Восстанавливается индивидуальное руководство аспирантами....

В профессорской висит плакат:

Долой егоровщину в математике!

Математику — на службу социализму!
 \end{quotation}
 
 \begin{quotation}
 Колмогоров - Александрову 16 апреля 1931 г.
 
  Вообще же, многое в Институте тебе не понравится, но теперь
ясно, что он не погиб (это и В.~В.~[Степанов] теперь признает) и работать
можно, и с пользой. Бермант (а отчасти, и сам Выгодский) очень
увлечены тем, чтобы перешибить Европу (например, заявляют,
может быть, и не без оснований, что через год-два наша
Lesezimmer должна стать лучше Гёттингенской). Деньги им дают
в большом количестве. Я, в частности, с 1 марта получил по 300
руб., и обещают по 400. Буду спрашивать о тебе.
 \end{quotation}

 \begin{quotation}
 Александров (по-видимому, с судна, отправляющегося из США) -- Колмогорову, 1 июня 1931
 
  По поводу Института. Вопрос об окладах, конечно, пустяки.
Я согласен, что его не следует так освещать, как я это сделал в
первом письме. Но, вообще, я от Института жду мало приятного. 
{\bf Приказы же и выговоры Выгодского значительно превосходят даже мои,
казалось бы, достаточно подготовленные ожидания.} Уж не знаю,
как Ястремский, но я желал бы все-таки себя от выговоров
Выгодского гарантировать, даже в случае, если по тем или иным причинам
не выполню учебного плана. Ну, там видно будет!...

...{\bf ведь и для меня мое дальнейшее пребывание
в Институте продолжает оставаться нерешенным вопросом, 
независимо от каких бы то ни было окладов}...
 \end{quotation}

 %%%%%%%%%%%%%%%%%%%%%%%%%%%%%%%%%%%%%%%%%%%%%%%%%%%%%%%%%%%%%%%%%%%%%%%%%%%%%%%

{\bf\punct Конференция по планированию математики.%
\label{ss:planirovanie}}
Процитируем публикацию \cite{Sbornik-planirovanie-0} в Мат.сборнике, номер 3-4 за 1931г.
\begin{quotation}
С 5 до 9 июня 1931 г. заседала в Москве I Всероссийская конференция по
планированию математики. При организации конференции имел место ряд недочетов, который привел, например, к тому, что некоторые специальности (геометрия,
статистика) не были представлены и недостаточно широко была привлечена общественность. Но, несмотря на эти недостатки, она знаменует исторические сдвиги.
Конференция по планированию математики -- это только одна из целого ряда
научно-плановых конференций, которые созывались, начиная с конференции по
планированию науки, созванной НИС ВСНХ. Они выражают факт вступления в
период социализма, рост пролетарского коммунистического ядра среди научных
работников и, наконец, решительный поворот, который в последнее время происходит в рядах интеллигенции СССР, которая, выбросив из своей среды вредителей-контрреволюционеров, решительно становится на сторону пролетариата,
строящего социализм, и большевистской партии.
\end{quotation}

Тов. Кольман сделал 5 июня доклад  «Современный кризис математики и основные линии ее реконструкции»,
открывавший конференцию. Автор не видел текст этого доклада (скорее всего он был опубликован в каком-либо 
тогдашнем журнале).
По этому докладу 9 июня конференция приняла
исполненную  в кольмановском
стиле резолюцию «О кризисе буржуазной математики и о реконструкции математики в СССР»
(она была опубликована в Мат.сборнике \cite{Sbornik-planirovanie}.)
\begin{quotation}
	Даже некоторые представители буржуазной математики начинают уже сознавать теперь, как они страдают от стихийности и бесплановости развития, которые
	являются основной чертой буржуазной науки. Закон неравномерности развития
	капитализма имеет место и в науке%
	\footnote{Какие теоретические откровения!},
	 и это ведет к тому, что наряду с продолжающимся развитием новых отраслей и методов растет и математическая декадентщина, и все более широкие слои буржуазных учений обрекаются на гниение.
	Наряду с кризисом обоснования растет кризис математики, основанный на бесплановости, узости научного кругозора, отрыве теории от практики... 
	 Выход из создавшегося положения — подлинное
	внедрение методов плановости и коллективности, ликвидация отрыва теории от
	практики, разрешение методологических проблем математики, ускорение темпов
	развития математики и ликвидация моментов загнивания — возможен лишь в условиях строящегося социализма, лишь на основе метода марксизма-ленинизма.
\end{quotation}
В резолюции, среди прочего, бичуются Лузин, Гюнтер, Егоров.

 Конференция не ограничилась обсуждением доклада тов. Кольмана:
\begin{quotation}
	На конференции собрались представители целого ряда математических факультетов и институтов и, прежде всего, научно-технических институтов, которые
	представили конференции тезисы своих запросов в области математических исследований. Среди них были: Физико-технический институт — Ленинград, Военно-техническая академия — Ленинград, Центральный аэро-гидродинамический институт (ЦАГИ) — Москва, Геофизический институт — Москва, Институт гражданских
	сооружений — Москва и мн. др. 
\end{quotation}
Очевидно, представители технологических институтов сообщали об интересующих их 
задачах прикладной математики, группы математиков лоббировали свои собственные исследования для включения их в планы.

Конференция вынесла  резолюции 
<<по организационным вопросам>>,
<<по вопросам анализа>>,
<<по вопросам дискретной математики>>,
<<по вычислительной математике>>. 
 Например, по {\bf дискретной математике} были выдвинуты  следующие темы
 (Если на клетке слона написано «буйвол» – не верь глазам своим!):
\begin{quotation} \small
1. Алгебраическая геометрия функционального пространства
[под этим подразумевались уравнения Фредгольма как явствует
из пояснительной записки].

2. Топология многообразий и функционального пространства.

3. Предельные переходы в алгебраических процессах (качественные методы
решений дифференциальных и интегральных уравнений).

4. Вопросы абстрактной алгебры и теории групп, необходимые для решения
проблем топологии  и квантовой физики. 

5. Исчисление конечных разностей и интерполяция. Отмечалось захирение
теории конечных разностей в наших основных математических центрах и необходимость подготовки специалистов в этой области.

6. Прикладная теория чисел (отыскание и развитие тех областей теории
чисел, которые смыкаются либо прямо с техникой, либо через посредство
естествознания и анализа, например кристаллография с теорией тройничных положительных квадратических форм и вопросы пространственных форм с экстремальными задачами дискретной геометрии. 
\end{quotation}	

В резолюции   включались  кольмановские мотивы:
\begin{quotation}
	Теория функций
	действительного переменного играла в московской математике еще недавно доминирующую роль. С ней исторически связано возникновение ряда основных направлений математического исследования в Москве (топология, теория вероятностей,
	функциональный анализ). Эта теория выяснила ряд очень важных связей между
	свойствами функций, основанных на понятии предельного перехода, и создала
	аппарат для решения этих задач, но она находилась не только за границей, но
	и в СССР, в плену у буржуазной идеалистической философии и поэтому не могла
	решить исторической задачи обоснования анализа. В лице своих наиболее реакционных представителей она вступила на путь
	 солипсизма\footnote{В солипсизме Кольман неоднократно обвинял Лузина.}.
	
	{\bf Удовлетворительное решение проблемы непрерывности может быть дано при
	условии внедрения диалектического материализма в мировоззрение математиков,
	к чему в последнее время уже намечается сдвиг.} 
\end{quotation}

В резолюции по дискретной математике это было поставлено в план:
\begin{quotation}
	{\bf 
В качестве методологической проблемы было выдвинуто основательное выяснение соотношения между непрерывным и дискретным в математике.}
\end{quotation}	

Что касается резолюции по вычислительной математике, то она, по-ви\-ди\-мо\-му,
была целиком в руках передовых сил. Речь там шла вовсе не о том, что называется
<<численными методами>>:
\begin{quotation}
	\small
	...разработку «математики малых точностей», понимая под этим разработку
	методики быстрых ориентировочных вычислений с ограниченной точностью при
	помощи устного и письменного счета, а также применения простейших вычислительных пособий и приборов; .......
	
	развитие геометрографии с обращением особого внимания на подведение
	экспериментальной базы под меру оценки простоты и точности графических построений и с разработкой теории и методики приемов точного черчения;
	
	развитие экспериментальной математики, понимая под ней нахождение физи­ческих эквивалентов математическим операциям;....
	
	пересмотр программ преподавания математики, начиная со школы I ступени,
	в направлении включения в них теории и практики рациональных вычислений, главным образом, над приближенными величинами, и прекращение привития учащимся
	вредных навыков псевдоточных и нерациональных по приемам вычислений;....
	
	обеспечение возможности испытания, конструирования и изготовления пробных образцов вычислительных приборов%
	\footnote{Автору кажется очевидным из контекста, что речь идет не об арифмометрах
	(кстати, промышленность тогда производила, по меньшей мере, <<Железного Феликса>>) и
	не о подобных вычислительных
машинах. К тому же резолюция была заявкой на финансирование, очевидно, что специалисты по прикладыванию линейки не годились для разработки сложных механических и электрических устройств,
и деньги  просили не на это (опять-таки, обратное было бы видно из контекста).}
	 и специального чертежного инструмента,
	для чего конференция считает безусловно необходимым организацию при некоторых
		на\-уч\-но-ис\-сле\-до\-ва\-тель\-ских институтах специальных, надлежаще оборудованных мастерских, а равно и предоставление валюты для приобретения нужных
	образцов из-за границы. 

	организацию центрального органа, руководящего делом организации вычислительной работы с возложением на него наблюдательных и инспекторских
	функций;...

	организацию общества «Долой неграмотность счета», с тем чтобы оно организовало секцию при обществе «Техника массам»;
	
	организацию соцсоревнования и ударничества в области как 	на\-уч\-но-ис\-сле\-до\-ва\-тель\-ской, так и практической вычислительной работы; 
\end{quotation}

Понятно, что те, которые в дальнейшем продвигали приложения математики к технике и естественным
наукам -- это вовсе не те, кто задалбывали   прикладной математикой окружающих в 1930-1931гг.

В одной из цитат выше отмечалось, что геометрия была представлена <<недостаточно широко>>. Теория вероятностей вообще не была представлена. Очевидно, никаких
последствий по прошествию небольшого промежутка времени это не поимело, а сама эта
конференция была  старательно забыта.

\sm

 {\bf\punct Мат.Общество, год 1931.%
 \label{ss:obshchestvo-1931}}
 О том, что происходило с Московским математическим обществом после переворота  
 ноября  1930г.,
 история умалчивает. Из архива Мат.общества соответствующие страницы 
 были вырваны. Довольно странно, что не удается найти иных архивных данных на эту тему. Приходится заниматься
 медиевистикой и довольствоваться тем, что мы имеем.
 
 В сентябре 1931г. начал издаваться новый журнал
 <<Математическая наука - пролетарским кадрам>> под редакцией
 Кольмана,  Хотимского и Яновской. В общей сложности вышел один номер этого журнала
 (и автор, к сожалению, найти его не смог). В журнале была статья
 <<Реорганизация Московского математического общества>> \cite{Reorganizatsiya}, 
 объемом 1 стр., она 
 полностью приведена Т.~А.~Токаревой в \cite{Tok-white}. 
 После потока малоинформативной, но знакомой мути сообщается следующее
 \begin{quotation}
  Общество выделило ряд секций: по исследовательской работе, по педагогической, по массовой работе,
  по антирелигиозноной работе.

Общество выпускает два математических журнала: научный
(прежний «Математический сборник») и массовый журнал%
\footnote{По мнению Т.~А.~Токаревой {\it по всей
вероятности, имеется в виду журнал, публикующий эту информацию.}
Это похоже на истину.}.
Избран новый президиум Общества в следующем составе [дальнейшее деление на строчки - моё]:
\newline
Кольман, Выгодский, Хотимский, Гельфонд,
\newline
Орлов [Киев], Бурстин[Минск], 
\newline
Хинчин, Голубев, Райков, Франкль, Яновская,
Люстерник, Лаврентьев 
\newline
и др.
 \end{quotation}
 
 В тексте отсутствуют даты, более удивительно, что отсутствуют должности
 (в том числе, не указан президент), и, наконец присутствуют <<и др.>>, которых
 кольмановцы  сочли нужным не называть.
  То есть происходило что-то непростое
 с формированием этого руководства (что-то было переиграно...),
 и позднейшие вспоминающие не сочли это нужным вспоминать. 
 
 По-видимому, президентом Общества к моменту опубликования журнала был Кольман%
 \footnote{Как утверждается в \cite{Demidov-TT}:
 {\it Факт президентства Э.~Кольмана не имеет никаких документальных подтверждений — протоколы заседаний 
 Общества той поры уничтожены. О пребывании Э.~Кольмана на этом посту в указанный период было сообщено
 в 60-е годы Л.~А.~Люстерником на заседаниях семинара по истории математики
  и механики механико-математического факультета МГУ.}}.
 Список президиума делится на 3 части. Сначала, по-видимому, руководство общества,
 состоящее из трех диалектиков и Гельфонда%
 \footnote{Кстати, в цитируемом сборнике была статья Гельфонда <<Элементарный вывод формул Эйлера>>.
 	На всякий случай, Гельфонд был математик и ни в какой научно-диалектической мути не замечен.}.
 
 Далее представители важных центров. М.~Х.~Орлов - это харьковский диалектик,
 который сверг Бернштейна с должности директора института. Бурстин  Целестин Леонович (C.~Burstin),
 1888-1938, член компартии Австрии, из-за проблем с работой в 1929г. перебрался в СССР,
 с 1931г. - директор физико-технического института АН Белорусской ССР (видимо, тогда же и основанного).
 Автор многочисленных математических работ. Собственно именно две этих
 фамилии разделяют московскую часть списка на две части.
 
 Отсутствие Ленинграда не должно удивлять, ни Гюнтер, свергнутый в марте 1931г., ни кадет Виноградов
 в качестве представителей не подходили. А  представители Харькова и Минска были политически благонадёжны.
 
  Не упоминавшийся ранее в этих записках Франкль  Феликс Исидорович (Felix Frankl,
 1905-1961)  --  тоже австрийский математик, имигрировавший в 1929г. в СССР,
 работал в Комакадемии, однако в 1931г. перебрался в ЦАГИ.
 Почему-то в списке нет Шнирельмана (который в тот год работал в Москве в пединстите).
 
\sm

 Александров (см. выше) писал, что общество {\it фактически бездействовало}.
 Но доклады на нем все же были:
 \begin{quotation}
 Колмогоров -- Александрову, 1 апреля 1931г.
 
  Доклад Кольмана (<<Динамическая и статистическая закономерность>>)
  вначале содержал дельное изложение марксистской 
теории соотношения случайности и закономерности (которая, действительно, 
значительно превосходит соответствующие теории
других школ), {\bf но затем раскрылся в фантазерстве (о будущей математике),
не всегда грамотном и в устах столь авторитетного
человека не безвредном.}

4-го делаю доклад о своей командировке (тогда же доклад 
Глаголева) ;

29-го — <<Современное положение и проблемы теории вероятностей>>.

Намечено еще в Матем. обществе: 1) о работах Зарецкого;
2) теория вероятностей и изучение действительности 
(неизвестно, где); 3) факты из исследований по основаниям математики.
 \end{quotation}
 
 Про концовку письма ничего не понятно, но доклад 29.04.1931 все же состоялся.

 \begin{quotation}
 Колмогоров -- Александрову, 1 мая 1931г.
 
  Много разговаривал с Лаврентьевым — он все больше привлекает
меня своей деловитостью и, в то же время, жизнерадостной 
непосредственностью. Работает в ЦАГИ и, не в пример многим прочим, с
искренним интересом.

29-го делал в Обществе доклад <<Современные проблемы теории
вероятностей>>. О.~Ю.~Шмидт, вообще появляющийся очень редко,
ушел сейчас же после конца доклада. Гельфонд и Шнирельман 
говорили о приложениях теории вероятностей к теории чисел. А.~Я. [Хинчин] начал с
авторитетного заявления, что все такие приложения — чистый
вздор, но кончил тем, что, хотя в теории чисел нет случайности,
существует <<иррелевантность>> распределения признаков целых 
чисел (например, аддитивных и мультипликативных свойств). Вокруг
этих вопросов (должно быть, как более философских и 
диалектических) и развернулась дискуссия, но были и внимательные слушатели
самого моего доклада.
 \end{quotation}

%%%%%%%%%%%%%%%%%%%%%%%%%%%%%%%%%%%%%%%%%%%%%%%%%%%%%%%%%%%%%%%%%%%%%%%%%%%%%%% 
 
 {\bf \punct По городам и весям.%
 \label{ss:goroda}} Что в те годы творилось по разным городам - предмет для отдельного исследования, и автор
 ее не пытался это изучать.
 Упомянем два эпизода, связанных с героями настоящей истории.

 \sm

 {\sc Харьков.} Там атаке был подвергнут Бернштейн. Приведем отрывок
 из Орлова Михаила Хрисанфовича (1900-1936), статья
 \cite{Orlov} цитируется по \cite{Baranets}:
 \begin{quotation}
  «Акад.
Бернштейн ведёт активную борьбу против марксизма-ле\-ни\-низ\-ма,
прикрываясь лозунгами аполитичности и непартийности. Но, как всегда в таких случаях, эти лозунги припрятывают
враждебную нам политическую линию. И действительно,
обосновывая аполитичность, непартийность и надклассовость
математики, акад. Бернштейн становится на вполне определённые идеологические позиции,
характеризуемые как реакционная философия воинствующего эклектицизма. Ещё на Всероссийском
математическом съезде 1927 года акад. Бернштейн проявил свои методологические ярко антимарксистские взгляды.
 \end{quotation}

 О харьковской истории можно прочитать в
 \cite{Bogolyubov-Rozhenko}.
 Бернштейн покинул пост директора основанного им  в 1928 г. Украинского института математических наук,
 а в 1933г. перебрался из Харькова в Ленинград (в Стекловку), а потом в составе Стекловки -- в Москву.
 В  \cite{Bogolyubov-Rozhenko} сообщается, что Орлов впоследствии стал жертвой репрессий,
 (я не нашел подтверждений этому, быть может, он умер сам, так и  не дожив до 1937 года).

 \sm

{\sc Cаратов.}
Голубев выпал из нашей истории, когда он стал ректором Саратовского университета
(март 1921-январь 1923). Процитируем его дневниковые записи <<Тетрадь Омега>>,
\begin{quotation}
Неудачником, в сущности, был и я. Вместо ученого пришлось много сил и времени
отдать организационной работе, к которой у меня никогда не лежало сердце, но к которой
были способности, может быть, больше, чем к непосредственному научному творчеству,
которое меня всегда влекло. Но и в организации, в сущности, ничего не вышло: большая
работа в Саратове, куда была вложена вся энергия молодости, вся вера и весь её
энтузиазм, в конце концов, разрушена Каценбогенами%
\footnote{Каценбоген Соломон Захарович  (1889-1946) -- ректор Саратовского
 университета в 1928-1932.} и иже с ними, а в дальнейшем уже
не было той энергии и доверия.
Дальнейшая работа учебная и организационная тоже шла при постоянном
противодействии конкурирующих идей и дала очень мало.

Канун моих именин.

В изгнании, в Свердловске%
\footnote{Стоит иметь в виду дату, происходившие события не могли не влиять на настроение писавшего.}, 27.VIII, 1942г. 
 \end{quotation}

Где-то в 1924г. Голубев  за отсутствием в университете математической литературы
стал читать то, что было под рукой \cite{Golubev-Prot-Tyulina}:
\begin{quotation}
Я занялся изучением бывших у меня 
оттисков работ по теории крыла Н.~Е.~Жуковского и
С.~А.~Чаплыгина и совершенно неожиданно для себя
обнаружил, что в этих работах содержится то, о чем
я мечтал: приложение теории функций комплексного
переменного к изучению явлений природы и техники.
При этом приложение непосредственное,
дифференциальные уравнения

Теорию функций комплексного переменного 
можно было буквально видеть, осуществить на модели
течением жидкости. Это было претворением 
гениальных идей Римана в области техники и естествознания Эти идеи совершенно очаровали меня...
\end{quotation}
 Уже в 1927 году он уже публикует в трудах ЦАГИ книгу <<Теория крыла аэроплана в   плоскопараллельном потоке>>
а в 1931 - <<Теория крыла аэроплана конечного размаха>>. 
Цитируем \cite{LL}:
\begin{quotation}
 Первые годы работы в Саратове для В.~В.~Голубева сложились хорошо:
в университете, помимо основных курсов, можно было читать интересные
специальные курсы, появились способные молодые ученые из окончивших
студентов и т. д. Однако превращение физико-математического факультета
в отделение педагогического факультета привело к постепенному сокращению
объема программы по математике. В конце концов объем преподавания
настолько сократился, что возможность проводить в стенах факультета чтение
специальных курсов по теории аналитических функций и по гидромеха­нике пропала.
В 1929 г. курс по теории крыла конечного размаха В.~В.~Голубеву
пришлось читать уже не на факультете, а от организации Осоавиахима
в Саратове. 
\end{quotation}

В марте 1930 года ректорат назначает смотр кафедры Голубева%
\footnote{Была такая байка: Начальство задает вопрос: <<Верно ли, что Ваш отец был священником?>>
--<<Да, верно, как и отец Чернышевского, имя которого носит наш университет.>>}.

Летом 1930  года потерпевший поражение в Саратове Голубев приезжает в Москву и устраивается к Чаплыгину в ЦАГИ.
Через три года он окажется первым деканом Мехмата.

%%%%%%%%%%%%%%%%%%%%%%%%%%%%%%%%%%%%%%%%%%%%%%%%%%%%%%%%%%%%%%%%%%%%%%%%%%%%%%%
%%%%%%%%%%%%%%%%%%%%%%%%%%%%%%%%%%%%%%%%%%%%%%%%%%%%%%%%%%%%%%%%%%%%%%%%%%%%%%%
%%%%%%%%%%%%%%%%%%%%%%%%%%%%%%%%%%%%%%%%%%%%%%%%%%%%%%%%%%%%%%%%%%%%%%%%%%%%%%%
%%%%%%%%%%%%%%%%%%%%%%%%%%%%%%%%%%%%%%%%%%%%%%%%%%%%%%%%%%%%%%%%%%%%%%%%%%%%%%%
%%%%%%%%%%%%%%%%%%%%%%%%%%%%%%%%%%%%%%%%%%%%%%%%%%%%%%%%%%%%%%%%%%%%%%%%%%%%%%%
%%%%%%%%%%%%%%%%%%%%%%%%%%%%%%%%%%%%%%%%%%%%%%%%%%%%%%%%%%%%%%%%%%%%%%%%%%%%%%%

\section{От диалектической республики к мехмату%
\label{s:mechmath}}

\COUNTERS

\epigraph{Я не устаю повторять, что только позже, когда мы начали ездить по свету, мы поняли, что таких математических факультетов, как в ЛГУ, в мире было очень мало, – а такого, как мехмат в МГУ, в мире просто не было нигде – по концентрации и по охвату всей математики, существующей на то время; по научному молодежному потенциалу.}
{А. М. Вершик}

Изложение чудесных событий в предыдущем параграфе не претендует на полноту.
Продолжение  в том же духе без всякого сомнения неминуемо привели бы к концу математики
в отдельно взятой стране. Мы знаем, что этого не случилось, а случилось совсем наоборот.

О годах 1932-1933  история умалчивает еще более, чем о предшествующих.
Вот отрывки, которые автору удалось найти.

\sm

{\bf\punct Происходила дискуссия\dots.%
\label{ss:discussion}}
 Из Академической стенограммы.
\begin{quotation}
ЛЮСТЕРНИК. Насчет теоретических изысканий происходила довольно длительная дискуссия, но в 1930 году, я помню, уже было решение, чтобы не было никаких перегибов в
этом отношении, и что теоретические исследования должны проводиться в полной силе.
Например, топологический кружок ни на минуту в то время не прерывал своей деятельности.
 То же самое и все остальные теоретические семинары. Так что такой установки
не было. Если кто-нибудь из молодых наших товарищей иногда высказывал отдельные
мысли, то к тому моменту, когда посылалось это письмо, уже эти вопросы были ясны, и
ни один теоретический кружок, ни один семинар в то время не прерывал занятий.
\end{quotation}

Процитируем также Александрова \cite{Alex-Shmidt2}:
\begin{quotation}
	... речь 
	шла о сохранении нашего основного математического журнала <<Математический сборник>>.
	
	Раздавались голоса, что журнал совершенно устарел, что никакой чистой математики не должно быть, должна быть математика только прикладная, которая отвечает непосредственным потребностям сегодняшнего дня. Отто Юльевич решительно этому воспрепятствовал. Встав в этот момент во главе <<Математического сборника>>, О.~Ю.~Шмидт твёрдо и решительно вел борьбу с всякими возникавшими тогда вульгаризаторскими тенденциями и придал журналу тот характер серьёзной научности, который он сохраняет до сих пор.
\end{quotation}

Несколько отголосков тогдашнего противостояния есть в цитатах из Академической стенограммы в
п.\ref{za-granitsej}.
Так или иначе в 1930-1931гг при участии лево-радикальных математиков чистую математику удалось отстоять.
Н.~М.~Бескин \cite{Besk29} датирует уход странных прикладников с Мехмата 1934 годом.

\sm

{\bf\punct Ветер меняется.}
Из \cite{Ilchenko}:
\begin{quotation}
Январь, 13, 1932

На Коллегии НКП [Наркомпрос%
\footnote{Наркомом просвещения РСФСР 9.1929-10.1937 был А. С. Бубнов (1884-1938), член ЦК ВКП(б)
в	1917—1918 и 1924—1937гг., начальник ГлавПУРа в 1924-1929гг.,
то есть все эти потрясения и подвижки происходили при одном министре.}] РСФСР заслушан доклад об МГУ. В постановлении говорилось: «Учитывая огромное значение МГУ,
являющегося крупнейшим научно-учебным центром по подготовке научно-исследовательских кадров и старейшим в стране университетом...:

1. Превратить Московский университет в образцовый университет, отвечающий полностью предъявляемым к нему требованиям.
\end{quotation}

\sm

{\bf\punct Смена лиц.%
\label{ss:hinchin-1932}} 
\begin{quotation}
Александров (Москва) -- Колмогорову (Днепропетровск), 20 мая 1932г.

 Новости:
 
1.	Директором Института назначен А.~Я.~Хинчин.
Весь кабинет пока не сформирован, но одним из заместителей 
директора уже назначен Г.~К.~Хворостин.
\end{quotation}

Смена лиц 1932 года, очевидно, не носила <<революционного>> характера, как падение Егорова.
Так или иначе, на месте Егорова во главе НИИ механики и математики оказывается
крупнейший математик. В декабре следующего 1933г. директором становится Колмогоров.

\sm

Что касается героя, вокруг которого собрана наша история,
то Хинчин сразу пригласил Лузина вернуться на МГУ.
В \cite{ErmTok} приведено длинное письмо Лузина с отказом
(по содержанию видно, что предложение было отправлено не позднее первых чисел июня).
Хинчин на этом не остановился (см. п. \ref{ss:luzin-do-1936}). Вот еще документ:
\begin{quotation}
Колмогоров--Александрову, 30 сентября 1932г.
\newline
Теорию функций действительного переменного начинает 
читать А.~Я.~Хинчин с передачей Ник. Ник. по его приезде.
\end{quotation}

Лузин, как известно, в итоге уклонился.

\sm

{\bf\punct Матобщество.%
\label{ss:alexandrov-1932}}
Весной 1932г. Александров восстанавливает Московское математическое общество.
О подробностях этого ничего не известно.
В 1949г. он сам писал:
\begin{quotation}
Новый период работы Общества начинается, таким образом, лишь с весны 1932 г., 
когда не только восстанавливаются научные заседания Общества,
но и начинают постепенно развиваться различные новые виды его
деятельности....

Руководство Обществом последовательно осуществлялось следующими лицами: 
..........
С 1932 г. по настоящее время П.~С.~Александров (президент). Вице-президенты
(последовательно): М.~Я.~Выгодский; И.~И.~Привалов и Л.~Г.~Шнирельман, И.~И.~Привалов
и С.~Л.~Соболев; А.~Н.~Колмогоров и В.~В.~Степанов. Секретари: В.~Л.~Гончаров,
Н.~Е.~Кочин, В.~В.~Степанов, А.~И.~Мальцев, С.~А.~Гальперн.

Новый устав Общества, действующий с 1932 г., вполне обеспечивает охват
Обществом всех активно работающих московских математиков.

\noindent
Своими почётными членами Общество избрало: О.~Ю.~Шмидта, С.~Н.~Бернштейна,
А.~Н.~Крылова (скончался в 1945 г.), Н.~М.~Гюнтера (скончался в 1941 г.)
и С.~А.~Чаплыгина (скончался в 1942 г.). 
\end{quotation}

Итак, новое руководство 1932г. состояло из Александрова, Выгодского и Гончарова%
\footnote{Гончаров Василий Леонидович (1896-1955),
	ученик Бернштейна, работы по комплексному анализу и комбинаторике.
По-видимому,  статья \cite{Goncharov}, 1944 (аннонсировалась в <<Докладах>> в 1942-1944гг.), была первой работой по статистике 
перестановок (исследовалось распределение длины  максимального
цикла перестановки и числа циклов). По-видимому, также это была последняя 
научная работа Гончарова, далее он стал членом-корреспондентом только что организованной
АПН РСФСР и в дальнейшем занимался педагогикой.}. Выгодский вскоре
(видимо, 1933г.) сменился на Привалова. Отметим  список почетных членов, и в числе
их Гюнтера, недавнего руководителя разгромленных было реакционеров на Ленинградском математическом фронте.

Выше приводился список Александрова-Головина заседаний Мат.Общества за 1930г, оборвавшийся 10 октября 1930
года. Вот его продолжение до марта 1933 года включительно.

\sm

{\small
\hangindent=1cm\noindent
21 мая 1932 г.  А. Я. Хинчин «{\it Математическая теория стационарной очереди}».

\hangindent=1cm\noindent
28 мая 1932 г.  М. Я. Выгодский «{\it Архимед, его эпоха и методология}».

\hangindent=1cm\noindent
11 декабря 1932 г.  А. Н. Колмогоров «{\it О геометрических идеях Plucker'a}
и Klein'а».

\hangindent=1cm\noindent
23 декабря 1932 г.  П. С. Александров «{\it О дальнейших перспективах развития топологии}».

\hangindent=1cm\noindent
5 февраля 1933 г.  А. Я. Хинчин «О некоторых новых результатах в общей
динамике».

\hangindent=1cm\noindent
11 марта 1933 г.  Г. М. Шапиро «О соответствии между поверхностями по
параллельности касательных плоскостей».

\hangindent=1cm\noindent
23 марта 1933 г.  С. А. Яновская «Маркс и математика».
}

Мы видим, что заседания редки (или составители списка не захотели их все упоминать), и видим два доклада, посвященных
истории и диалектике (Выгодский и Яновская).
С апреля 1933г. общество входит в обычный режим.

\sm

Александров, по-видимому, вскоре решает взять на себя методологическую инициативу.
Согласно сообщению Куроша \cite{Kurosh} 1964 года,
\begin{quotation}
..В феврале [правлением ММО] 1934 г. обсуждался план
специального выпуска журнала «Фронт науки и техники» [5—6, 1934г.], содержавшего
математические статьи методологического характера; позже этот выпуск
был переиздан в качестве - отдельного сборника.
\end{quotation}

Сборник вышел  под редакцией Яновской
c названием <<О философии математики>> \cite{Filosofiya-matematiki},
его участники  Колмогоров, Александров, Курош,  В.~М.~Молодший, Гливенко, А.~М.~Фишер (Alexander Fischer)
 и сама Яновская (три статьи%
 \footnote{Судя по библиографиям, лишь одна из них входила в журнал. Статья Фишера была взята из журнала <<Под знаменем марксизма>>,
 1934,5, а статья Молодшего из то же журнала, 1935, 3}).
Само слово <<философия математики>> с выходом этого сборника приобретает мирный характер...

Общество расширяет формы своей деятельности, появляется секция средней школы,
премии для молодых ученых, московские математические олимпиады, оно занимается выдвижением в Академию.
Общество быстро становится влиятельной организацией, с ним  считался Наркомпрос
\cite{Kurosh}, в частности оно активно участвовало в улучшении школьных учебников,
о чем речь пойдет в п. \ref{ss:kiselev}.

\sm

{\bf \punct Об упрощенцах и упрощенчестве.%
\label{ss:stetski}} Это название статьи Стецкого%
\footnote{Стецкий Алексей Иванович, 1896-1938, член ЦК ВКП(б) с 1927г.},
крупного партийного деятеля в <<Правде>> от 4 июня 1932 года. Приведем ее полностью.
Жирный текст - это то, что было выделено в Правде. Курсив - мой.

\begin{quotation}
 ОБ УПРОЩЕНЧЕСТВЕ И УПРОЩЕНЦАХ

Никак нельзя забывать указание Ленина о том, — «... что без солидного философского обоснования
никакие естественные науки, никакой материализм не может выдержать
борьбы против натиска буржуазных идей и восстановления буржуазного миросозерцания. 
Чтобы выдержать эту борьбу и провести ее до конца с полным успехом, естественник 
должен быть современным материалистом,
сознательным сторонником того материализма, который представлен Марксом, 
т. е. должен быть диалектическим материалистом» (Ленин, т. XXVII, стр. 187).

Наши научные работники должны тщательно, систематически изучать диалектический материализм, 
они должны брать под огонь всякую пропаганду поповщины, идеализма во всех их видах и проявлениях.
Как бы ни старались сторонники механицизма снять задачу борьбы за диалектику,
эта задача — борьба за диалектический материализм против буржуазного мировоззрения — 
{\bf стоит и будет стоять} у нас в порядке дня.

Нельзя однако эту борьбу вести в порядке деклараций, заклятий, словесной шумихи. 
{\it От научного работника требуется прежде всего всерьез овладеть своей специальностью,
знать свою научную область, упорно работать в ней, — лишь при этом условии возможны
действительная борьба с буржуазным мировоззрением и победа диалектического метода.
Без этого условия декларация о партийности в науке {\bf ничего не стоит}.
Какой толк из декларации «за партийность в математике», если люди, выбрасывающие этот лозунг,
не знают математики. Можно десятки раз клясться в верности диалектическому материализму 
и требовать его осуществления в математике, но все эти клятвы и заявления — {\bf пустой звук},
если они не подкрепляются знанием самого предмета, его метода, его проблем,} пониманием того,
как и каким путем {\bf конкретно в данной дисциплине} проявляется буржуазное мировоззрение.

Пора и среди беспартийных специалистов поставить более квалифицированно пропаганду
марксизма-ленинизма и материалистической диалектики.

{\it Показать беспартийным, как вести доменный процесс на основе марксизма-ленинизма 
или как писать картины или строить дома на основе материалистической диалектики, мы не можем.
И тот, кто берется за это, тот {\bf шарлатан}}. Но дать им лучшие произведения основоположников марксизма,
поставить дельные лекции по вопросам диалектики, показать, как Маркс, Энгельс, Ленин и Сталин 
в своих работах проводили диалектический метод, — это мы не только можем, {\bf 
но и должны непременно делать}.
А те научные работники, которые искренно переходят к нам, учась у Маркса и Ленина,
они позаботятся и подумают над тем, как в своей области применять диалектический метод
на основе имеющихся у них знаний, фактического материала, действительного понимания проблем своей науки.

И, наконец, о научных обществах. Вокруг Комакадемии образовалась порядочная сеть обществ:
физиологов-марксистов, врачей-марк\-сис\-тов, математиков и пр.

Все эти общества имеют почти сплошь коммунистический состав. Они замкнулись в своей скорлупе
и отгородились от широких обществ советских специалистов. Не пора ли здесь эти перегородки снять,
коммунистам войти в широкие общества советских специалистов, чтобы рука об руку с ними, и споря,
когда нужно, работать над разрешением не на словах, а на деле, научных проблем? 
\end{quotation}

Автор - член ЦК, он за лекции по марксизму, за то чтобы научные работники их слушали.
Всецело за...
Но статья не про это, а про то, что наука должна заниматься наукой. В качестве примера приводится
математика (надо думать, какой-то математик поучаствовал в спуске этой телеги). И
мы видим
ласковое слово
<<шарлатан>>, обращенное лично к Кольману (а также Яновской, Хотимскому  и  их сотоварищам из других
наук).

Интересны и два последних абзаца. Комакадемическая программа (которой Кольман, видимо, и занимался) создания
профессиональных обществ на базе марксизма более не представляет интереса. Профессиональные общества
должны заниматься профессиональными проблемами...

Остается не ясным, как произошла смена ветра, и что за этим стояло.
Но диалектическая свистопляска 1931г. (по крайней мере, в математике, где у нее не было внутренних
профессиональных
источников, см. наше обсуждение <<Декларации инициативной группы>>) завершилась.
Понятно, что разговоры о диалектике, споры о диалектике, а также попытки самовыразиться
на базе диалектики, случались и позже,
но это уже было не о том.

Вскоре Кольман публикует в журнале <<Под знаменем марксизма>> статью <<На текущие
темы>> \cite{Kolman-temy}, в которой среди прочих тем он поддерживает статью Стецкого, выступает <<против ленинизма в кролиководстве>>
(вариациями на подобные темы был славен издававшийся Кольманом и в том же 1932 году почивший в бозе
журнал <<За марксистско-ленинское естествознание>>) и пишет о том, что надо поднимать уровень 
диалектического материализма. Насколько автор может судить по нескольким публикациям Кольмана
следующих  лет, тема диалектической кузькиной матери для математиков более-менее выветривается из  его сочинений
(но он еще вспомнит Лузина в своей книжке \cite{Kolm}).

А вот цитата из одной из апологий Яновской \cite{Levin-2012}:
{\it ... которая всегда защищала науку и ученых и выступала против упрощенчества в науке, свойственного марксистам-догматикам.}%
\newline
Она тоже вняла...

\sm

{\bf \punct  Смена ветра. Образование.%
\label{ss:wind-2}} Уходили в легенду и удивительные образовательные
изобретения. 

\begin{quotation}
Колмогоров (Днепропетровск)-- Александрову (Цюрих), 22.9.1932

\sm

Вчера появилось Постановление ЦК о высшей школе. Все
благоразумные предложения, о которых много говорилось
(80—85\% времени на специальные дисциплины, введение 
факультативных курсов по выбору, восстановление лекций, некоторое
ограничение практики в технических вузах, дипломные работы,
дифференцирование оплаты по квалификации преподавателей и
по роду занятий (лекция, семинар, упражнения), вступительные
экзамены даже для окончивших рабфаки, набор аспирантуры по
отзывам кафедр, индивидуальные планы работы аспирантов,
запрещение досрочных выпусков и т.д.) широко развиты. Говорится
еще о восстановлении ученых степеней, об учреждении 
<<инженерной академии>> и т.п.

Как обычно, правовое направление мысли законодателя, кроме
того, выразилось в

1.	предоставлении директорам университетов права налагать
взыскания на профессоров, вплоть до запрещения преподавать в
вузах на 1 — 5 лет (и именно за пропуск занятий без 
уважительных причин);

2.	установлении шестидневной недели с 48 днями летнего перерыва
и 12 днями зимнего (при этом профессорам и преподавателям
обеспечивается 2 месяца отпуска летом).

В силу того же постановления в Днепропетровске открывается
Университет.
\end{quotation}

Дальнейшая история образования 30х годов была сложной (см. кое-что
в статьях А.~Ф.~Лапко \cite{Lapko} и Люстерника  \cite{Lyu-2}, см. также книгу
\cite{Chanbarisov}, которую мне не удалось найти), 
но дальше уже была конструктивная деятельность со сложностями,
ошибками и с поисками решений, и в целом успешная%
\footnote{В тех же статьях 1967г. Люстерник  высказывал тревогу в связи с объявленной тогда
школьной реформой. Кажется, это был один из немногих таких голосов.}.

 При всех элементах фантасмагоричности 1929-31гг. итоги  были не столь уж однозначны.
В ходе катавасии 1929-31гг.
только в Москве
возникли МАИ, МИСиС,  МЭИ, Керосинка, МИСИ, которые вошли в группу ведущих технологических
вузов.
Кстати, такого рода вузы давали хорошие рабочие места для математиков. 

Московское механико-математическое отделение Физмата вышло из потрясений с лучшим составом профессоров,
чем в него входило.%
\footnote{Автор пишет об этом, не имея пофамильных списков, которых мне не удалось найти, а по мелькающим именам.}
Конечно, это был прежде всего итог развития математики предыдущих
лет, но важно ведь не только есть ли люди, но и устраиваются ли они на работу
(это же не всегда так...).
Например, в ходе загадочного директорства Выгодского профессором становится
Колмогоров%
\footnote{В 1932г. Выгодский выдвигал (нереализовавшийся) проект Математической энциклопедии
\newline
{\it 
Александров -Колмогорову, 20 мая 1932 г.
\newline
Состав редакции еще не намечен, но Выгодский говорил мне
(в частном порядке), что он будет стремиться во главе всего 
издания поставить ТЕБЯ, <<нашего советского Эйлера>> (его слова)!}}.
В 1930-1931гг. в качестве профессоров появляются Гельфонд (один из участников переворота в матобществе),
Голубев и Лаврентьев, в 1930г. Петровский становится ассистентом.
<<Исключенный из комсомола Курош>>, по-прежнему, остается наблюдаемым.
Так что издали кадровая политика 1930-31гг. выглядит как довольно успешная.

Надо, впрочем, иметь в виду увольнения 1933г. Вот что писалось в официальном издании 
\cite{Narkompros}
Наркомпроса в 1935г.:
\begin{quotation}
	Выше уже отмечалось, что в период перестройки университета и превалирования количественных показателей над качественными в ряды преподавательского состава МГУ просочилась некоторая, правда незначительная, прослойка людей с недостаточной научной квалификацией; в особенности это имело место по физическим и математическим дисциплинам. Уже в течение 1933 г. университет освободился от этих преподавателей, заменив их квалифицированными научными работниками. При этом чрезвычайно знаменательно, что проведенный пересмотр преподавательского персонала под углом зрения повышения требований к его теоретическому уровню не только не снизил, но, наоборот, повысил процент партийно-комсомольской прослойки кафедр (в особенности среди старшего преподавательского состава), что находит свое объяснение в большом росте научных кадров из аспирантского молодняка.
\end{quotation}

{\bf \punct Мехмат.%
\label{ss:mechmath}}
В мае 1933 г.  была восстановлена факультетская структура МГУ. Были организованы шесть факультетов:
механико-математический, физический, химический, биологический, почвенно-географический и рабфак. 
Через год был восстановлен исторический факультет.

Что касается названия мехмат, а не матмех, то в 30е на Мехмате было больше механиков, а не математиков.
Похоже на то, что существенное влияние на тогдашний Мехмат тогда оказывал ЦАГИ, где 
квалифицированные механики и математики были в самом деле нужны, а ЦАГИ в самом деле был нужен стране.
Приведем цитату из А.~И.~Некрасова (\cite{Golubev-Prot-Tyulina})
\begin{quotation}
Влияние ЦАГИ на подготовку кадров (на 
стороне) наиболее ярко проявилось в постановке 
преподавания на механическом отделении 
физико-математического факультета Московского университета, 
отмечал.  Это влияние сказалось 
даже на самой структуре механического отделения..
Причем в задачи отделения была поставлена 
подготовка кадров на широкой научной базе.
\end{quotation}
В свою очередь, из известных математиков в ЦАГИ в 30х работали, в частности, Лаврентьев, М.~Келдыш, Гантмахер
(как, впрочем, и  герой, вокруг которого собрано наше повествование). 

1 мая 1933г. деканом новообразованного механико-математического факультета
стал Голубев%
\footnote{Одна из причуд  истории этих бурных лет: началась она с падения Егорова, пострадавшего за христианские убеждения,
Голубев же стал христианином и вроде бы того не скрывал...}. 

Впрочем, это уже была другая история.

\section*{Добавление к \S \ref{s:mechmath}. Исторический фон}

{\small

{\bf\punct  Дело имяславцев и арест Егорова.%
\label{ss:losev}} Арестован он был по делу 
<<Истинно православной церкви>>. 

 Известно, что с начала 20х Православная церковь (находившаяся под атаками со стороны властей и революционных кругов)
была расколота по вопросу о политике по отношению к Советской власти.
С одной стороны была группа <<обновленцев>> во главе с А.~И.~Введенским, которая поддерживалась властью.
Были также группы, которая с середины 20х встали на курс политической 
лояльности к Советской власти. В 1927г. митрополит Сергий Старогородский,
 заместитель Патриаршего местоблюстителя опубликовал декларацию
 о лояльности Советской власти (позже, в период восстановления церкви
Сергий стал патриархом). Лояльные Советской власти верующие, по-видимому, 
должны были оказаться на стороне Сергия (и, в частности, это должно
было оказаться началом конца обновленчества).  С другой стороны немедленно возникли 
 группы, осуждавшие  политику Сергия и обвинявшие его в соглашательстве.
 Разногласия этих групп с Сергием имели политическую природу, при формулировке
 возражений непросто было воздержаться от политических обвинений,
 после формулировки каковых эти группы должны были 
 стать предметом особого интереса ОГПУ. Впрочем, эти замечания носят
 умозрительный характер.

  Некоторые современные православные авторы считают
<<Истинно православную церковь>> реальной организацией. Так это было или не так, в любом случае удар ОГПУ 1929-1931гг. был направлен против нелояльных
церковных групп.

Имяславцы (в число которых входили Егоров и А.~Ф.~Лосев, о котором сейчас и пойдет речь) были одной из таких групп%
\footnote{Известный историк церкви Д.~В.~Поспеловский относил имяславцев к <<крайним группам>>.}, 
едва ли связанной с большой организацией
(если даже таковая в самом была).
% Аресты имяславцев начались весной 1930г., тогда, в частности был арестован
%друг Егорова и знакомый Лузина А.Ф.Лосев. Согласно В.Н.Щелкачёву Егорова (как и его самого) арестовали в ночь
%с 9 на 10 октября 1930г. Дело имяславцев
% \footnote{В следственном деле фигурировала фамилия Бухгольца. Данных о том, что его арестовывали, %однако не известно.}
Документы по делу имяславцев 
видели, по крайней мере, А.~А.~Тахо-Годи (вторая жена Лосева) \cite{Taho-Godi}, В.~Н.~Щелкачёв
\cite{ShCh}, проходивший по этому делу (и, видимо, не бывший имяславцем), и историк С.~С.~Демидов \cite{Demidov}.

\sm

	 Арест Лосева был связан с публикацией
	его (ставшей впоследствии знаменитой) книги <<Диалектика мифа>> \cite{Losev-mif}.
	Процитируем
	биографию Лосева, составленную  Тахо-Годи:
	\begin{quotation}
		Книга, где Лосев раскрыл действенность мифов научных%
			\footnote{Вот образец: {\it Все
			эти бесконечные физики, химики, механики и астрономы имеют совершенно
			богословские представления о своих <<силах>>, <<законах>>, <<материи>>,
			<<электронах>>, <<газах>>, <<жидкостях>>, <<телах>>, <<теплоте>>, <<электричестве>> и
			т.д. Если бы они были чистыми физиками, химиками и т.д., они ограничились бы
			выводом только самих законов и больше ничего, да и всякие <<законы>>, даже
			самые основные и непоколебимые, толковались бы у них исключительно лишь как
			гипотезы. Это было бы чистой наукой \dots  Но, конечно, надо помнить, что тут речь идет исключительно о
			чистой науке и что реально никогда такой чистой науки не существует, что это
			есть анализ не реально-исторической науки, но лишь ее теоретически-смысловых
			основ и структур. С этой стороны видным делается как мифологическое засилие
			в современной науке у наивных ее <<практиков>>, у всяких экспериментаторов и
			философски не мыслящих  ее работников, так и полное несходство существа науки
			с существам мифологии ее работников, так и полное несходство существа науки
			с существам мифологии.}
		\newline
	При том, что книга Лосева является яркой и высокохудожественной
	(а автор настоящих записок относится к тем, кто считает, что
	книга заслуживает прочтения и продумывания),
стоит отметить, что логическая аргументация здесь (и не только здесь) сводится к тому, что
взгляды всех этих <<бесконечных физиков>> не соответствуют взглядам самого Лосева.}%
		, философских и литературных%
		\footnote{\it Я лично терпеть не
			могу женщин с непокрытыми головами. В этих последних есть некоторый тонкий
			блуд, - обычно мужчинам нравящийся. Также нельзя быть христианином и любить
			т.н. <<изящную литературу>>, которая на 99\% состоит из нудной жвачки на тему о
			том, как он очень любил, а она не любила, или как он изменил, а она осталась верной, или как он, подлец, бросил ее, а она повесилась или повесилась не
			она, а кто-то еще третий и т.д. и т.д. Не только <<изящная литература>>, но и
			все искусство, с Бетховенами и Вагнерами, есть ничто перед старознаменным
			догматиком <<Всемирную славу>> или Преображенским тропарем и кондаком; и
			никакая симфония не сравнится с красотой и значением колокольного звона.
			Христианская религия требует мифологии колокольного звона. Христианин,
			если он не умеет звонить на колокольне или не знает восьми церковных гласов
			или, по крайней мере, не умеет вовремя развести и подать кадило, еще не
			овладел всеми тонкостями диалектического метода. Колокольный звон, кроме
			того, есть часть богослужения; он очищает воздух от духов злобы поднебесной.
			Вот почему бес старается, чтобы не было звону.}%
		, а главное, социальных – в эпоху “великого перелома” и “построения социализма в одной стране”, – была запрещена цензурой, выбросившей все идеологически опасные места. А.~Ф. не убоялся запрета и вставил в печатавшийся текст ряд мест, которые были исключены цензурой. Предлог для ареста книги и ее автора был найден. А поскольку все издательские дела с чиновниками и типографиями вела супруга А.~Ф., то и она попала в тюрьму, а затем и в лагерь. Но иного выхода, кроме как высказать вслух заветные свои идеи, у философа не было. В одном из лагерных писем к жене он справедливо писал (22/III–1932): “В те годы я стихийно рос как философ, и трудно было (да и нужно ли?) держать себя в обручах советской цензуры”. “Я задыхался от невозможности выразиться и высказаться. Этим и объясняются {\bf контрабандные вставки в мои сочинения после цензуры, и  в том числе (и в особенности) в “Диалектику мифа”}. Я знал, что это опасно, но желание выразить себя, свою расцветавшую индивидуальность для философа и писателя превозмогает всякие соображения об опасности”...
	\end{quotation}
	
		Книга  известна, и представляется весьма маловероятным, что Лосев мог предлагать цензуре именно этот эпатажный текст%
		\footnote{ {\it  С точки зрения коммунистической
			мифологии, не только <<призрак ходит по Европе, призрак коммунизма>> (начало
			<<Коммун. манифеста>>), но при этом <<копошатся гады
			контр-революции>>, <<воют шакалы империализма>>, <<оскаливает зубы гидра
			буржуазии>>, <<зияют пастью финансовые акулы>>, и т.д. Тут же снуют такие
			фигуры, как <<бандиты во фраках>>, <<разбойники с моноклем>>, <<венценосные
			кровопускатели>>, <<людоеды в митрах>>, <<рясофорные скулодробители>>... Кроме
			того, везде тут <<темные силы>>, <<мрачная реакция>>, <<черная рать мракобесов>>;
			и в этой тьме - <<красная заря>> <<мирового пожара>>, <<красное знамя>>
			восстаний... Картинка! И после этого говорят, что тут нет никакой
			мифологии.}}.
		 Так или иначе власть жестко отреагировала на подмену согласованной для печати
	книги. Цитируем Тахо-Годи:
	\begin{quotation}
	В Деле Алексея Федоровича за No 100256, которое я 
	обозрела в июне 1995 года, находится справка о «Диалектике
	мифа», составленная начальником IV Отделения ИНФО
	ОГПУ Соловьевым (л. 193—196). Там указано, что в книгу
	размером в 153 страницы  (7 1/2 п. л.) автор «без 
	согласования с Главлитом%
		\footnote{Так называлось цензурное ведомство.} внес ряд принципиальных исправлений и
	дополнений (на с. 7, 8, 11, 17, 18, 19, 22, 70, 71, 74, 75, 77,
	78, 84, 87, 90, 92 и 95)». Кроме того, {\bf  он вставил страницы с
	98 по 134, с 241 по 263}%
\footnote{Отметим, что 153 машинописных страницы выросли до 268 книжных.}...
	
	Эта замечательная справка была подписана 18 апреля
	1930 года. Власти действовали мгновенно. Лосева 
	арестовали в ту же ночь.	
		\end{quotation}
		5 июня последовал арест жены Лосева
	
	Разумеется, было приказано уничтожить весь отпечатанный тираж <<Диалектики мифа>>. 
	Книга к тому времени уже начала продаваться, кроме того, известно о нескольких конкретных
	экземплярах, оказавшихся в библиотеках. Как будто один из этих экземпляров долгое
	время лежал в 
	Библиотеке имени Ленина, и как будто с него была даже снята копия.

	История с рукописью получила широкую огласку и упоминалась  в двух речах,
	произносившихся  на  XVII съезде ВКП(б), проходившем 26.6-13.7.1930г \cite{16}.
	О Лосеве говорил аж сам Каганович Лазарь Моисеевич.
	\begin{quotation}
	А вот пример из области философской литературы. В Правде была помещена рецензия о семи книгах философа-мракобеса Лосева.  Но последняя книга этого реакционера и черносотенца под названием <<Диалектика мифа>>, разрешенная к печатанию Главлитом, является самой откровенной пропагандой наглейшего классового врага (несмотря на то, что она разрешена Главлитом, она не увидела света). Приведу лишь несколько цитат из этого реакционного и мракобесовского произведения.
	
	{\it <<Дыромоляи, говорят, и сейчас не перевелись в глухой Сибири. А я, по грехам своим, никак не могу взять в толк, что земля движется, и неба никакого нет.>>%
		\footnote{Продолжение в оригинале:
			\newline
			{\it <<Учебники читал, когда-то хотел сам быть астрономом, даже женился на астрономке. Но вот до сих пор никак не могу себя убедить, что земля движется и что неба никакого нет. Какие-то там маятники да отклонения чего-то куда-то, какие-то параллаксы... Неубедительно. Просто жидковато как-то. Тут вопрос о целой земле идет, а вы какие-то маятники качаете. А главное, все это как-то неуютно, все это какое-то неродное, злое, жестокое. То я был на земле, под родным небом, слушал о вселенной, «яже не подвижется»...>>}}\dots
		
		 <<Католичество, которое хотело спасти живой и реальный мир, имело полное логическое право сжечь Джордано Бруно>>%
		 \footnote{{\it Если
		 	что-нибудь не имеет конца, - следовательно, оно не имеет границы и формы.
		 	Если что-нибудь не имеет границы и формы, это значит, что оно ничем не
		 	отличается от всего прочего. Но если оно ничем не отличается от всего
		 	прочего, то, следовательно, невозможно установить, существует ли оно вообще
		 	или нет. Итак, если мир бесконечен, то это значит, что ровно никакого мира
		 	не существует. Нигилизм Нового времени так, в сущности, и думает. Восхвалять
		 	бесконечность миров заставляло тут именно желание убить всякий мир; и
		 	католичество, которое хотело спасти живой и реальный мир, имело полное
		 	логическое право сжечь Дж. Бруно. А уничтожать мир необходимо было
		 	тем, кто выставлял на первый план свою личность.}}
		
		<<Диалектический материализм есть вопиющая нелепость.>>
	
	 <<Сжигать людей на кострах красивее, чем расстреливать, так же как готика красивей и конкретней новейших казарм, колокольный
	звон автомобильных воплей и платонизм -- материализма.>>%
	\footnote{\it Бить по воздуху плетью -
		пустое и отвлеченное занятие; но, например, порка крестьян и рабов есть
		проявление конкретной идеи, ибо тут - реальное общение личностей, понимающих
		друг друга.}
	
	<<Только монах один -- не мещанин. Только монах понимает правильно и достаточно глубоко половую жизнь, и только он знает глубину и красоту женской души.>>
	
	<<Говорили: идите к нам, у нас полный реализм, живая жизнь вместо ваших фантазий и мечтаний. Оказывается, полный обман и подлог. Нет, дяденька, не обманешь. Ты, дяденька, с меня шкуру хотел спустить, а не реалистом сделать. Ты, дяденька, вор и разбойник.>>
	
<<Коммунистам нельзя любить искусство. Раз искусство, значит гений. Раз гений, значит -- неравенство. Раз неравенство, значит -- эксплоатация.>>
	
	<<Иной раз вы с пафосом долбите: социализм возможен в одной стране, не чувствуете ли вы в это время, что кто-то или что-то на очень высокой ноте пищит у вас на душе: не-ет.>>}
	
	И это выпускается в Советской стране. О чем это говорит? Это говорит о том, что у нас все еще недостаточно бдительности. ({\sc Голоса:} Кто выпускает? Где выпущено?). Разрешено Главлитом.
 ({\sc Голоса:}
	Чье издание?) Это выпущено самим автором, но ведь
	вопрос заключается в том, что у нас, в Советской стране, в стране пролетарской диктатуры, на частном авторе должна лежать узда пролетарской диктатуры.
	А тут узды не оказалось. Очень жаль.  ({\sc Голоса:} Правильно.)
	\end{quotation}

%Pod_Skansenem

	А вот из речи товарища Киршона%
	\footnote{Киршон, Владимир Михайлович, 1902-1938,
	успешный советский драматург и деятель РАПП. Расстрелян по обвинению в троцкизме.}
	\begin{quotation}
	Или, вот, например, товарищи, место из книжки,
	которую цитировал уже здесь товарищ Каганович	
	<<Диалектика мифа>> Лосева -- она тоже вышла в 1930г. и пропущена Главлитом. Этот Лосев, помимо всяких
	других откровенно-черносотенных монархических высказываний, между прочим сообщает:
	
	{\it Пролетарские идеологи или ничем не отличаются от капиталистических гадов и шакалов, или отличаются,
	но еще им неизвестно, чем собственно они отличаются.}

Коммунист, работник Главлита, пропустивший эту книжку, в которой нас в лицо называют капиталистическими гадами и шакалами, мотивировал необходимость ее разрешения тем, что это <<оттенок философской мысли>>%
\footnote{Из \cite{Taho-Godi}:
	\newline
	{\it Кстати сказать, в заключении Главлита, сделанном цензором (он же баснописец) Басовым-Верхоянцевым, значилось: автор “трактата” – “совершенно чуждый марксизму (идеалист)”. Но разрешение было все-таки дано, “разве только в интересах собирания и сбережения оттенков философской мысли”.}}	
({\sc Смех. Голос.} Оппортунизм на практике.)
{\bf А я думаю, что не мешает нам за подобные оттенки ставить
к стенке} ({\sc Апплодисменты. Смех})
	\end{quotation}

	Продолжим цитирование Тахо-Годи:
	\begin{quotation}
		Новосёлов%
		\footnote{Новосёлов Михаил Александрович (1864-1938), известный 
			деятель православной церкви, последовательный противник Советской власти 
			и антисергианец. Сторонник создания катакомбной церкви. Арестован в марте 
			1928г. и более на свободу, видимо, не выходил.} активно распространял пропагандистские брошюры, писал письма своей пастве, вел переговоры с епископом Дмитрием, устанавливал связь с приходами Ленинграда, рассылал документ под названием “Большое имяславие”, в  котором соединились цели имяславия и антисергиевские.
		А.~Ф.~Лосева с Д.~Ф.~Егоровым этот документ “покоробил” некоторыми антисоветскими местами. Но тем не менее В.~М.~Лосева отправила “Большое имяславие” в Ленинград с одобрительным отзывом (3/Х–1930, л. 128). А это уже означало практическую деятельность, и круг Лосева – Новосёлова объективно начинал играть роль некоего руководящего центра...
		
		Если вдуматься поглубже, то уже только один документ “Большое имяславие” мог бы послужить поводом для ареста Лосевых. Тем более что этот “программный документ” для отложившихся от Сергия “в значительной мере” содержал выписки из сочинений А.~Ф. (8/Х-1930, л. 128). и его “заметок” (22/VII–1930, л. 117). Лосев не только просмотрел этот документ по просьбе Новоселова, но и был с этим документом “согласен” (там же). Кроме того, распространялись и другие брошюры, “Ташкентская” или “Воронежская”, тоже антисергиевские и тоже составленные имяславцами (причем автор последней якобы Новоселов), и она “небезызвестна” и Лосеву (22/VI–1930 от ПП ОГПУ по ЦЧО Алексеева, п/нач. ИНФО ОГПУ т. Герасимовой, л. 203). В Деле Лосева находился и документ, составленный священником Павлом в Дивеевской пустыни, начинавшийся словами: “Св. Нифонт предрекал...” По словам Лосева, этот документ содержит идеи, “соответствующие моим мыслям о советской власти, как сатанинской, о Церкви, как борце с ней, о капитализме, как грехе”. (15/ХП-1930,л. 142)...
		
		Документ [Большое имяславие], можно сказать, опасный и антисоветский, а вот А.~Ф.~Лосев, признавая свою вину, считает, что участвовал в этой антисоветской организации “без особой практической деятельности” (8/Х–1930, л. 129)... 
	\end{quotation}
	
	Егоров, по-видимому, был кем-то вроде религиозного главы московских имяславцев.
	Во всяком случае, Тахо-Годи упоминает и отчасти цитирует имяславческую
	декларацию 1922г,  первой там стояла подпись Егорова. Собрания имяславческого 
	кружка тоже происходили на квартире Егорова. Понятно, что Егоров 
	должен был стать объектом повышенного интереса ОГПУ. Арестован он был 
	в ночь с 9 на 10 октября 1930г. (одновременно с несколькими
	другими фигурантами дела).
	
	Так или иначе, в руках  ОГПУ был ряд  реальных эпизодов%
	\footnote{Вот отголоски из \cite{Taho-Godi}:
		\newline {\it Н.М.Соловьев -- враг компромиссов -- требовал <<политики
	патриарха Гермогена>>. Патриарх Тихон ... ответил, что
он <<Гермогеном быть не хочет>>.}
\newline
{\it 
Лосев был сторонником имяславия как чисто религиозной идеи,
и политика в этом деле была ему чужда, поэтому крайности
современных имяслацев на Кавказе, активно занимающихся политикой,
были для него неприемлимы (17/VII/1930) (Они, например, отказывались получать паспорта, распространяли антисоветские листовки, совершали другие акты гражданского неповиновения).}
\newline
Было бы весьма странным, если бы в Церкви не было людей, требующих <<политики патриарха Гермогена>>
(надеюсь, что читателю ясно, что это значит), их позиция вполне понятна, и странно таких людей
было бы сейчас осуждать. Но еще более странна
(если задуматься)
современная позиция интеллигенции и властей, о том, что никаких таких 
людей не было, и  все это придумали и  фальсифицировали в ВЧК/ОГПУ/НКВД.
И вообще в других странах революционные эпохи были полны заговорами,
и лишь у нас было не как у людей. 
 Это не отменяет наличия
массовой фальсификации судебных дел (что вполне очевидно в отношении 1937-38г),
так и фальсификаций в конкретных  делах и превентивных арестных кампаний
(что тоже более-менее понятно в отношении 1929-31гг.), так же как не отменяет 
жесткости законодательства и жестокости его применения.}, наказуемых с точки
	зрения тогдашнего  права. Так или иначе, 
	чекисты начали разработку разветвленной организации 
	“Политический и административный Центр всесоюзной контрреволюционной организации церковников <<Истинно-Православная Церковь>>. К делу были привлечены ряд деятелей Церкви из
	других городов и имяславцы на Кавказе.

	Дело, по-видимому, проходило без применения <<незаконных методов следствия>>
	(есть воспоминания Щелкачева, Лосев же впоследствии называл следовательницу,
	ведшую его лично, <<ласковой коброй>>). Так или иначе, дело о большом
	заговоре было скроено.
	
	В общей сложности по  делу привлекалось 48 человек,  часть их них отпала%
	\footnote{В деле фигурирует фамилия Бухгольца, но вроде бы никакие источники
		не сообщают нам о том, что он арестовывался. От имяславства он к тому времени,
		по словам Тахо-Годи, отошел.},
	часть была вскоре освобождена. Осуждено 33 человека, из них 17 отправились в места заключения
	(сроки от 5 до 10 лет), остальные были высланы. 
	
	Против Егорова реальных обвинений, скорее всего, не было (судя по
	 дальнейшим упоминаниям его имени, в нем активного контрреволюционера
	 не видели).
	 
	 Что касается Лосева, то
	 по приговору от 3 октября 1932 получил 10 лет лагерей, в конце сентября
	 был отправлен в Кемь, далее на Свирьстрой,
	  освобожден 8 сентября 1932г. постановлением
	 Коллегии ОГПУ от 7 сентября 1932 года (подробности, видимо, не известны). Судимость снята 4.08.1933.

\sm

{\bf\punct Пиррова победа Комакадемии.%
\label{ss:pirr}}
Комакадемия 
создавалась как центр социалистической мысли и как <<кузница>>
интеллектуальных марксистских кадров.
Ее руководитель - историк М.~Н.~Покровский(1868-1932) -- 
пользовался загадочно большим влиянием в советских верхах.

 Процитируем воспоминания Кольмана  \cite{Kolman-vospominaniya}.
\begin{quotation}
 Комакадемия представляла собой научно-исследовательское учреждение при ВЦИКе, 
 созданное для того, чтобы развивать на ди\-а\-лек\-ти\-ко-ма\-те\-ри\-али\-сти\-че\-ской основе общественные и
естественные науки. Этим она должна была противостоять Академии
Наук, от которой при тогдашнем ее составе нечего было ожидать.
\end{quotation}
На этот счет в верхах имелись разные мнения, была также идея оставить естественно-научное
отделение, а остатки гуманитарного присоединить к Комакадемии.

Наряду с этим существовал Институт красной профессуры, которым руководил тот же Покровский.

В принципе в этих структурах работали и  крупные ученые, из героев нашего повествования
это Шмидт, Гельфонд,
Люстерник. Из числа других героев рассказа, там училась Яновская, работали Кольман, Хотимский,
 а также один из лузитан Лихтенбаум. 

Так или иначе, кадры были взрощены и отобраны (взрощенных не хватало%
\footnote{Покровский \cite{Pokrovski}:
	{\it Из 96 красных профессоров, работающих сейчас в Москве, если не ошибаюсь, только 43 занимаются преподавательской деятельностью, а остальные 53 работают в разных наркоматах, работают в партийных учреждениях и т. д.}}, нужно было отбирать
на местах), и в 1929-1931гг. они
двинулись на овладение высшей школой и научными институтами.
В принципе важнейшая цель властей - усиление просоветских элементов и <<общественных организаций%
\footnote{В смысле ВКП(б), Комсомол, советские профсоюзы.}>>
в вузах была достигнута. 

Но кадры были неоднозначные,  далеко уехать на них было нельзя.
Кстати, забавно наблюдать как эти герои сцеплялись между собой.

\sm

Позволим себе лирическое отступление.
Не столь уж сложно обнаружить, что статья по алгебре лишена классового анализа.
Но уж если классовый анализ ввести, то статья все равно окажется несовершенной с точки зрения истинных
ценителей, а неизбежное (или почти неизбежное)  изменение качества статьи едва ли можно будет объяснить
чем-либо, кроме недостаточной глубины классового или диалектического анализа.

 Автор случайно наткнулся на несколько
сочинений такого рода, не имея намерений их разыскивать.
В.~Н.~Молодший \cite{Molodshy-Leifert} напал на учебник <<Методика математики для педагогических техникумов>>
под редакцией Лейферта (М.~А.~Белецкая, И.~Д.~Киселев, А.~Р.~Кулишер, Л.~А.~Лейферт, Е.~И.~Отто). 
Б.~Кедров%
\footnote{Кедров, Бонифатий Михайлович (1903-1985), учившийся в то время в Институте красной профессуры,
сын одного из наиболее печально известных деятелей  Гражданской войны.
Впоследствии -- видный советский марксистко-ленинский философ.}
обнаружил, что в учебнике статистики под редакцией Хотимского, <<борьба с вредительскими теориями ведется крайне слабо>>,
и написал по этому поводу 10-страничную рецензию (хотя, казалось бы в этом отношении авторы учебника были невинны).
Мария Смит-Фалкнер тоже пыталась нападать на этот учебник. Однако (см. Предисловие к \cite{Hotimsky-book})...
\begin{quotation}
 Однако дискуссия с полной очевидностью вскрыла, что эта критика велась  М.~Н.~Смит с меньшевистско-богдановских позиций
 и, что в работах самой М.~Н.~Смит содержится целая система антимарксистких взглядов. В принятой Комфракцией Общества 
 статистиков-марксистов при Комакадемии 2/III.1932г. резолюции по поводу статистических работ М.~Н.~Смит 
 зафиксировано, что <<вплоть до 1931г. в них под флагом стопроцентного марксизма проводились буржуазно-махисткие
 и просто идеалистические теории. В этих работах марксистско-ленинская теория подменена эклектической смесью
 из меньшевиствующего идеализма, богдановщины и рабски повторяемых идей буржуазной статистики,
 преимущественно английской (Боули и др.). Под флагом марксизма в этих теориях
 проступает идеалистическая теория <<истинных величин>>. Качественные различия низводятся М.~Н.~Смит 
 до количественных вариаций определенного размера. Место единого метода материалистической диалектики
  у М.~Н.~Смит  занимает ряд методов формальной логики, обычно представляющий для буржуазной теории истинно
  научную методологию. Диалектика в работах М.~Н.~Смит фигурирует лишь как один из многих методов и систематически извращается.
  В предисловии к <<Основам статистической методологии>> М.~Н.~Смит, ни словом не упомянув о Марксе, Энгельсе, Ленине,
  откровенно поставила перед собой задачу распространения среди <<русской>> учащейся молодежи идеи Боули....
  На протяжении ряда лет работы М.~Н.~Смит  засоряли умы советской молодежи элементами классово чуждой идеологии.>>%
  \footnote{Не стоит преувеличивать значения этих слов, это просто принципиальная товарищеская критика,
   приводных ремней к ОГПУ такие речи пока не имели... Это случится потом в 37-38гг.}
\end{quotation}
История умалчивает, за кем было последнее слово в этом научном споре. Мы знаем лишь, что 1933г. Хотимский оказался начальником политотдела
Завитинской моторно-тракторной станции на Дальнем Востоке.

В том же 1932г. товарищ Яновская подвергла суровой большевистской критике книгу товарища Выгодского
<<Основы исчисления бесконечно малых>> (впрочем, у Яновской уже тогда появились  признаки протрезвления,
и ее рецензию мы обсудим чуть ниже).

\sm

Часто пишут, что Комакадемия оказалась рассадником партийной оппозиции и вольномыслия. Возможно, что отчасти это было и правдой (в принципе это высказывание носит статистический характер
и может быть проверено). В любом случае принципиальные  научные споры комакадемических
бойцов и сопровождающие их <<оттенки философской мысли>>  производят непростое впечатление.

 Выше мы 
обсуждали события на одном отдельно взятом факультете одного (но главного) университета.
Зрелище отдает  фантастикой.  Примерно то же происходило по всей стране.
 Так продолжаться не могло (это было несовместимо с существованием
 образования и науки в стране), и перед властью вставала объективная задача
или снизить значение марксистской идеологии, или загнать ее в регламентированные рамки
(и уменьшить возможности марксизма как способа самовыражения).
Смягчение идеологии (по этой причине, и не только по ней) вскоре началось. Марксизм  был жестко регламентирован,
а Комакадемия как  интеллектуальный проект творческого марксизма%
\footnote{Стоит иметь в виду, что когда Шмидт заведовал естественно-научной частью
Комакадемии, в ее научных изданиях публиковались, например, разные молодые математики.
В виду экзотичности этих изданий я ни разу до них не добирался, но можно уверенно сказать,
что там бывали статьи высшего класса (скажем, Гельфонда или Шнирельмана,
они потом были перепечатаны). Скорее всего, это касалось не только математики...
Но в конце 1930го <<естественно-научное руководство>> подверглось суровому товарищескому осуждению.} провалилась.
C 1932гг. начинается ее закат...

\sm

{\bf\punct Математики и нематематики.%
\label{ss:biologi}} Возможно, что у современного читателя поведение математиков, описанное в предыдущих
параграфах, вызовет восклицание: <<Да как же такое могло быть?>> Автор,
понимая эти чувства и напомнив еще раз про зеркало, не меняющее голову с ногами,
все же должен заметить, что  те герои рассказа, которые были именно математиками, при всех своих качествах
оставались во время этих событий профессионалами. Напомню, что во время истории 1936
политизированные нападавшие  разделяли оценки  политические  и оценки научные. 
Люди, которые писали диалектическую муть в 1930-31гг., не были математиками, а Люстерник и Гливенко,
участвовавшие
 диалектическом сборнике, написали туда нормальные статьи. <<Декларация инициативной
группы>> при всей политизированности осталась неподверженной диалектической свистопляске
и содержит конструктивные предложения, которые потом сами авторы декларации начнут реализовывать.

{\bf Будучи расколотым политическим,
московское математическое сообщество сохранилось как профессиональное сообщество 
с профессиональной этикой.}

\sm

Часто математиков, включая в обсуждение общеисторический фон,
 сравнивают с биологами%
 \footnote{Реконструкции событий 30х годов в советской
 биологии предлагались  Н.~П.~Дубининым \cite{Dubinin}
 (очевидцем и деятельным участником событий) в 1973г.,
 генетиком В.~Н.~Сойфером в 1989-90гг, \cite{Soyfer}, и историком 
 науки Н.~Л.~Кременцовым \cite{Krementsov}
 в 1997. Взгляды этих авторов на движущие силы тогдашних событий заметно различаются. Отметим, что В.~Н.~Сойфер 
опубликовал целую книгу \cite{Soyfer-Dubinin}, 431с., обличающую 
Н.~П.~Дубинина. Публикация Н.~Л.~Кременцова (Nikolai Krementsov) в российском сборнике вызвала гневную отповедь \cite{Lebedev-Krementsov}
в том же сборнике. Автор этих записок меньше всего ставит своей
целью реконструкцию данной сложной истории (это и не его специальность),  цель этого пункта - лишь сравнение
стартовых условий математиков и биологов перед 30ми годами.}. Попробуем и мы провести это сравнение. 
Укажем четыре фактора, сложившихся к началу 30х и сыгравших в 
дальнейшей истории советской биологии резко отрицательную роль.

\sm

1)  {\sc Диалектизация биологии}. Ей занимались и настоящие биологи тоже.
Автор, работая над этими записками, наткнулся на подобные сочинения
И.~И.~Агола, 
см. \cite{Agol},  \cite{Agol-1}, С.~Г.~Левита \cite{Levit}, А.~Н.~Баха%
\footnote{Бах был также руководителем ВАРНИТСО и входил в 1929г. в группу захвата 
	Академии. Напомним однако, что на заседаниях Комиссии по делу Лузина 1936г. он же осаживал
зарвавшихся математиков.}
\cite{Bach}. Мы воздержимся от цитирования, всё было как в математике, только 
говорилось о биологии, и говорили это настоящие биологи%
\footnote{Еще о диалектике. Цитируем Тахо-Годи \cite{Taho-Godi}.
\newline
{\it <<Диалектическая разработка математики>> -- новая задача Лосева.
Он осуществит ее позже, вернувшись из лагеря [1933] сочинит книгу <<Диалектические
основы математики>>. Валентина Михайловна [жена Лосева] напишет
к ней интереснейшее предисловие...}
\newline
 Книга эта в итоге в 1997 была издана, желающие
могут ее изучить.}.

\sm
2) {\sc Чрезмерные обещания.} Товарищ М.~Б.~Митин \cite{Mitin}, главный советский философ 1930гг., на дискуссии по генетике 1939г., обращаясь
к А.~С.~Серебровскому, говорил:
\begin{quotation}
 Весь Ваш «генофонд» с пятилеткой в два с половиной года - реакционнейшее измышление,
  коренящееся в Ваших ошибочных теоретических построениях.
\end{quotation}
Не знаю, какие именно точно высказывания Серебровского Митин имел в виду.
Возможно, что вот это \cite{Serebrovski}:
 \begin{quotation}
 Если подсчитать, какое количество сил, времени, средств освободилось бы, если бы нам удалось очистить население нашего Союза от различного рода наследственных страданий, то наверное пятилетку можно было бы выполнить в 2 1/2 года. 
 \end{quotation}
Это не буквально обещание.
Но
процитируем также <<Резолюцию о задачах генетических исследований в области животноводства>>
в <<Трудах Всесоюзной конференции по планированию генетико-селекционных 
исследований>> \cite{Serebrovski-Vavilov}, проходившей под руководством А.~С.~Серебровского и Н.~И.~Вавилова в 1932г.
\begin{quotation}
	Конференция считает, что... на протяжении второй пятилетки...
	за счет улучшения генотипа сельскохозяйственных животных
	может быть достигнуто по СССР увеличение мясности не менее 10\% по крупному рогатому скоту,
	15\% по овцам, 20\% по свиньям, кроликам и курам, повышение молочности на 10\%
	процентов с сохранением \% жира, повышение яйценоскости на 20-25, шерстности
	на 50-70\%, улучшение тяговых свойств лошадей и пушных свойств.
\end{quotation}

Автор, не будучи биологом, не может сам оценить, насколько обоснованы были эти надежды.
Но ответ кажется понятным из слов   Серебровского (в трудах той же конференции):
\begin{quotation}
 Невозможность для генетиков капиталистических стран вести сколько-нибудь систематическую, широко поставленную работу с сельскохозяйственными животными поневоле заставляет их сосредоточиваться на работе с различными лабораторными зоологическими объектами -- кроликами (в лучшем случае), с морскими свинками, мышами, бабочками, жуками и, наконец, со знаменитой дрозофилой. Характер объектов накладывает определенный отпечаток и на характер тематики, и мы имеем чрезвычайно пышно разработанные главы генетики, тесно связанные, например, с дрозофилой, и полную неразработанность таких глав, которые бы имели особое значение для нашего народного хозяйства. Достаточно, например, сравнить степень разработанности такой типично дрозофильной главы генетики, как кроссинговер или инверзии и др. хромозомные аберрации, аллеломорфизм и пр. и, наоборот, такую главу, как генетика роста, генетика скороспелости, иммунности, не говоря уже о молочности, яйценоскости и т. п., -т. е. типично не дрозофильные темы, чтобы увидеть, как однобоко идет развитие зоологической генетики, оторванной от животноводства, запертой в академические лаборатории на Западе. И поэтому то -- несмотря на все блестящие достижения генетики, которые трудно преувеличить -- несмотря на это, мы должны признать состояние этой западной генетики совершенно неудовлетворительной....
\end{quotation}
 
 Основная эффект предполагалось получить с помощью метизации. Однако 
 тот же Серебровский говорил следующее:
 \begin{quotation}
 Коль скоро мы являемся почти полными младенцами в области  частной генетики сельскохозяйственных животных, не знаем чем различаются породы по генотипу, постольку мы, как генетики, не в состоянии оказать сколько-нибудь существенной помощи метизационной работе. 
\end{quotation}

Еще цитата из Серебровского:
\begin{quotation}
Сейчас же я остановлюсь на следующем вопросе: я считаю совершенно
 необходимым (хотя против этого раздаются голоса) нам, генетикам
	и селекционерам, взять на себя определенные обязательства в разрешении
	тех или иных народнохозяйственных проблем. Скажем, по линии проблемы
	мяса мы должны сказать, что мы берем на себя разработку таких планов
	и методов работы, чтобы обеспечить 12\% прироста мясной продукции
	в 1937 г. в смысле увеличения веса туш. Ответственность за остальные
	проценты прироста должны взять другие науки. 	
\end{quotation}
И заключение его выступления:
\begin{quotation}
	Спланируем же и будем вести нашу работу так, чтобы к назначенному сроку иметь честь услышать нетерпеливый звонок: «Алло, говорит Соцстроительство. Готова ли твоя работа?» — и иметь право ответить: «Готова и открыты новые многообещающие перспективы».
\end{quotation}

Стратегически идея подъёма сельского хозяйства на основе научной генетики
была правильной. Но  научной базы для улучшения животноводства
на ближайшие годы (на которые давались обещания) тогда  не было
(и, что забавно, отсутствие этой базы очевидно из самого выступления Серебровского%
\footnote{Думаю, что биолог, читая его речь, найдет там разные
	 фантастические предложения, я же  должен уклониться от подобных рассмотрений.}). С растениеводством
ситуация была получше, но и в этом случае не было оснований для чрезмерного оптимизма.
Из выступления Вавилова (там же)
 \begin{quotation} Этот отрыв генетики от селекции особенно свойственен западноевропейским селекционерам, а также работникам в Канаде и САСШ, где селекция главным образом сосредоточена в руках семенных фирм...
 	\end{quotation}
 Закончим тему заключительными словами Вавилова на этой конференции:
 \begin{quotation}
 	Я лично нисколько не сомневаюсь в том, что, не взирая ни на какие
 	трудности, условия нашей работы улучшаются с каждым днем, и тот могучий рост и тот энтузиазм, который безусловно имеется в нашей среде,
 	который все нас, вопреки всем трудностям, заставляет работать, -
 	доказывает, что нас ждут большие достижения, что, когда мы с вами
 	соберемся через несколько лет, мы с честью будем говорить о том, что
 	план, намеченный нами, выполнен.
 	
 	Заканчивая эту конференцию, мы дадим друг другу слово, что со
 	своей стороны сделаем все, чтобы этот план не только выполнить, но
 	и перевыполнить. 
 \end{quotation}
 
 \sm
 
 3) {\sc Евгеника.} Это учение об улучшении человеческой породы путем тем или иных
 мероприятий. Оно развилось в с конца XIX века, в ряде протестантских
 стран (Германия, англосаксы, Скандинавия) различные мероприятия по облагораживанию
 рода людского начали (особенно в 30е годы) проводиться в жизнь.
 В основном применялись принудительные стерилизации (в некоторых местах
 тайные, Северная Каролина), в Германии вскоре научились решать проблему более радикально.
  
  В России интерес к евгенике возник в начале 1920х под впечатлением
  массовой гибели людей в 1914-1922г. Возглавил евгеников Н.~К.~Кольцов,
  он также издавал <<Русский евгенический журнал>>%
  \footnote{Были также <<Известия Бюро по евгенике при Российской Академии наук>>
  	(Ю.А.Филипченко) и
  <<Медико-биологический журнал>> (С.~Г.~Левит).}.
  В то время евгеники имели поддержку в правительственных кругах
  (обычно называют наркомов Н.~А.~Семашко и А.~В.~Луначарского). Где-то в 1929-30г. власти сформулировали
  свое резко отрицательное отношение к евгенике, и Кольцов сразу свернул евгеническую
  деятельность.
  
  Сам Кольцов занимался теоретическими исследованиями и, возможно, относился 
  к этой науке с чрезмерным энтузиазмом. Но едва ли его можно было
  бы упрекнуть в какой-либо склонности к людоедству. Однако,
  были авторы, высказывавшие сочувствие стерилизационным проектам,
  а Серебровский [надо помнить, что он занимал весьма высокое положение
  в ученом мире] выступил и с такими научными идеями:
  \begin{quotation}
  	Социализм, разрушая частнокапиталистические отношения в хозяйстве, разрушит и современную семью, а в частности разрушит в мужчинах разницу в отношении к детям от своего или не своего сперматозоида. Точно так же, может быть несколько труднее, будет разрушено стыдливое отношение женщины к искусственному осеменению, и тогда все необходимые предпосылки к организации селекции человека будут даны\dots
  	
  	Совершенно очевидно, что при наличии такой обстановки евгеника – наука займет одно из самых почетных мест в системе человеческих наук, так как получаемые на ее базе успехи будут самыми драгоценными из всех мыслимых. В самом деле, при свойственной мужчинам громадной спермообразовательной деятельности и при современной отличной технике искусственного осеменения (находящего сейчас широкое применение лишь в коннозаводстве и овцеводстве), от одного выдающегося и ценного производителя можно будет получить до 1 000 и даже до 10 000 детей. При таких условиях селекция человека пойдет вперед гигантскими шагами. И отдельные женщины и целые коммуны будут тогда гордиться не «своими» детьми, а своими успехами и достижениями в этой несомненно самой удивительной области – в области создания новых форм человека%
  	\footnote{
  		Из воспоминаний Н.~П.~Дубинина \cite{Dubinin}:
  		\newline
  		{\it 
  		Острое столкновение произошло после выхода в свет статьи А.~С.~Серебровского <<Антропогенетика и евгеника в социалистическом обществе>>. Эта статья была крайней формой евгенических предложений о практических методах разведения людей с целью получения новой породы человека. Как-то в лаборатории в присутствии ряда сотрудников А.~С.~Серебровский сам поднял вопрос о своей статье. Я четко и без обиняков сказал, что статья эта реакционная, антимарксистская и ничего, кроме вреда, генетике принести не может, что заявление о том, будто бы путем генетического улучшения наших людей страна сможет выполнить пятилетку в два с половиной года, ничего не имеет общего с учением марксизма о личности и об обществе. Чувствовалось, что А.~С.~Серебровский возмущен, но он молчал. 
  		\newline
  		Тогда крайне рассерженным заступником выступил С.~Г.~Левит. Ведь он был главным редактором <<Трудов медико-генетического института>>,
  		 в которых увидела свет статья А.~С.~Серебровского. Левит обрушился на меня за то, что я посмел в такой форме высказать свое мнение, и заявил, что все попытки связать возражения против статьи с методологическими вопросами,- это, по его мнению, и есть кустарный, как он выразился, <<голоштанный>> марксизм. Этим и закончилась наша дискуссия при полном молчании остальных участников встречи.}}.
  \end{quotation}
  
 Так или иначе, до практического применения евгеники в СССР дело не дошло.
 В дальнейшем бывших евгеников в буквальном смысле слова не преследовали,
 но при случае тому или другому биологу коллеги поминали евгеническое
 прошлое (обострилось это в конце 30х, когда евгеника широко пошла в дело в
 фашистской Германии). Вспомнили его и Кольцову \cite{Bach-Koltsov},
  когда тот в 1939г. выставил свою
 кандидатуру в действительные члены АН СССР. Кольцов отказался каяться, и в итоге лишился поста директора
 института.
 
 \sm
 
 4) {\sc Особая роль дарвинизма.}
 Дарвинизм был одной из идеологических подпорок марксизма
 (так же, как и сейчас является подпоркой либерализма).
 Это способствовало повышенному интересу официальных философов к
 биологам (особенно к генетикам), с другой -- биологи получали
 возможность не только выяснять, кто из них лучший ученый
 (что неизбежно), но и кто лучший дарвинист. В последнем случае была возможность
 обращаться к высшим авторитетам и вовлекать их во внутринаучные дискуссии.
 
 \sm
 
 Так или иначе, генетики входили в 30е со славой Вавилова и   Кольцова,
 но общественные позиции вавиловцев и кольцовцев по ряду причин оказались крайне уязвимыми%
 \footnote{Понятно, что одним из факторов успешного развития генетики в СССР
 в 20х и первой половине 30х годов было как хорошее государственное финансирование,
так и более конкретно поддержка государственными деятелями конкретных сильных
ученых.  Но генетики перестарались с обещаниями,
и в следующем раунде это естественным образом повлекло недоверие к их заявлениям. 
\newline
Обещания, очевидно, не ограничивались 
Совещанием 1932г.
Процитируем Дубинина \cite{Dubinin}:
\newline
{\it
Устами Н.~И.~Вавилова, А.~С.~Серебровского и других генетики того времени торжественно обещали, что их наука приступила к выполнению своих задач по ускоренному строительству новых основ подъема урожайности сортов и продуктивности пород в условиях социалистического колхозного строя. Поскольку это обещание шло от людей, работавших по теории генетики, оно в какой-то мере задевало практиков селекции, имевших свои методы и свои достижения. Крупные селекционеры растений шли в основном самостоятельной дорогой. Совершенно самобытен был И.~В.~Мичурин. Знаменитые животноводы М.~Ф.~Иванов и Е.~А.~Богданов по ряду вопросов спорили с А.~С.~Серебровским и его учениками.	
\newline	
Однако Н.~И.~Вавилов и А.~С.~Серебровский в эти годы обладали наибольшими возможностями. Н.~И.~Вавилов был первым президентом Всесоюзной академии сельскохозяйственных наук имени В.~И.~Ленина. А.~С.~Серебровский встал во главе генетики и селекции животных во Всесоюзном институте животноводства. Н.~К.~Кольцов, А.~С.~Серебровский и М.~М.~Завадовский были академиками ВАСХНИЛ. После их избрания они выступили с докладами на коллегии Народного комиссариата по земледелию, в которых развили картину скорых блестящих успехов генетики в практике сельского хозяйства.}
\newline
Похоже, что генетики оказались также в неладах с селекционерами...
}.
 Это само по себе ничего не предрешало, но в последующее десятилетие
 это было постоянно действующим
 отрицательным фактором в потоке споров, конфликтов, случайностей,
 а также борьбы за денежные потоки и должности.
 
  \sm
 
 Пара слов о послевоенных событиях в биологии. Математическая свара 1936г. с участием
 партийных кругов кончилась в целом благополучно, но могла кончиться и очень плохо.
 Что касается биологов, то в 1945-1948гг. генетики активно атаковали Т.~Д.~Лысенко.
 Тот оставался непоколебленным,
 но, с другой стороны, окликов  <<Прекратить!>> по отношению к нападавшим генетикам
 не видно \cite{borba}.
 И тут генетикам крупно повезло, они обратили на свою сторону молодого аппаратчика,
  заведующего Сектором науки Отдела пропаганды и агитации ЦК ВКП(б) Ю.~А.~Жданова,
  каковой 10 апреля 1948г. прочитал на семинаре лекторов обкомов партии
   лекцию <<О положении в советской генетике>> \cite{borba}, по-видимому, выйдя 
   за пределы своих аппаратных полномочий.
  Молодой человек был сыном некого А.~А.~Жданова, а также кавалером некой
  С.~И.~Алиллуевой. Последнее, возможно, и не поимело значения,
  а член Политбюро и секретарь  ЦК КПСС А.~А.~Жданов
  имел в верхах много недоброжелателей.
  Идея расследовать историю сессии ВАСХНИЛ, не задумываясь, какое отношение
  она могла бы иметь  к А.~А.~Жданову (одному из самых высших лиц страны)
  и
  к его смерти, которая почему-то последовала через несколько дней после 
  завершения сессии ВАСХНИЛ, кажется автору несколько экстравагантной.%
    \footnote{Так или иначе, источники, подтверждающие, что вместе с лекцией Жданова-мл.
    	в гости к биологам из бутылки  вышел Большой Джин борьбы за власть в Политбюро ЦК КПСС, которому 
    	было не так уж много дела до науки, имеются. Прежде всего это воспоминания
    	\cite{Shepilov} Шепилова Дмитрия Трофимовича (1905-1995),
  	кандидата в члены Президиума ЦК КПСС (1956—1957), секретаря ЦК КПСС (1955—1956, 1957) \cite{Shepilov}, опубликованные лишь в 2001г.
  См. также удивительный документ, опубликованный в 1989г. \cite{borba}, часть 3, с.112-113,
   где секретарь ЦК пишет проект бумаги
  (видимо, постановления ЦК) против своего сына.
См. также запись от 31 мая 1948г. в дневниковых записях наркома В.~А.~Малышева \cite{Malyshev}.
\newline
К сожалению, при обилии обсуждений темы преследований генетики, количество опубликованных
документов невелико.}%
$^,$%
\footnote{Кажется странной и рассказываемая история \cite{Soyfer} о том, как
Сталин и Партия 15 лет боролись с генетикой и все не могли ее одолеть.}.
В любом случае это уже другое время и  другая история.}

\section{В следующем раунде%
\label{s:heroes}}

\COUNTERS

Здесь нас интересует судьба не научной части участников переворота 1929-1931 года.

\sm

{\bf\punct Хворостин.%
\label{ss:hvorostin-2}}
Мы видели его поступающим в аспирантуру в декабре 1929г., а в мае 1932г.  - зам. директора НИИ механики и математики.
В 1935ом году он становится (ди)ректором Саратовского университете. В Саратовских университетских хрониках
(например, \cite{Avrus}, \cite{Avrus-1})
 он - 
<<блестящий администратор>>, стремившийся превратить Саратов в <<Геттинген на Волге>>. Он пытается собрать в 
университете серьезных ученых, в том числе из числа сосланных, и подающую надежды молодежь 
(и он в самом деле собрал много
известных имен из разных наук). Он не в ладах 
с местным парткомом... В апологии Хворостина \cite{Avrus-1}
рассказывается

\begin{quotation}
 На партийном собрании университета, проходившем 5 дней
(28 марта – 1 апреля 1937 г.), и в докладе зав. отделом школ и науки
Саратовского обкома партии Голяндина, и в большинстве выступлений (их было 23) Гавриил Кириллович был подвергнут острой и
часто несправедливой критике. Его обвиняли в том, что он ставил
науку выше политики, допустил серьёзные ошибки при подборе
кадров. Так преподаватель физмата С.~А.~Суслов заявил, что у Хворостина не случайные ошибки, а их система, которую нужно вырвать с корнем. Выступившие на собрании Н.~И.~Усов и Д.~И.~Лучинин (будущие ректоры СГУ) обвинили Гавриилу Кирилловича 
в отсутствии большевизма, потере классовой бдительности, зажиме критики в свой адрес и в отношении близких ему людей.
\end{quotation}

Судя  высказываниям многих авторов, Хворостин имел поддержку кого-то в ЦК, и до поры до времени 
он был сильнее провинциального парткома.

В <<Математической хронике>> \cite{Uspehi-36-1} в первом томе только что появившегося сериала <<Успехи математических наук>>
сообщается:

\begin{quotation}
	\small
 В Саратовском государственном университете до последнего времени по математике
 работал только один профессор тов. Г.~П.~Боев. {\bf С начала текущего учебного года
университет развернул работу пяти математических кафедр, пригласив из Москвы для
руководства четырьмя новыми кафедрами профессоров А.~Я.~Хинчина, И.~Г.~Петровского,
А.~Г.~Куроша, В.~В.~Вагнера}. Кроме руководителей кафедр приглашен ряд лучших из
окончивших в 1935 г. аспирантуру математического института при Московском университете 
молодых математиков — тт. Симонов, Базилевич, Гордон и Барабанов. Таким
образом в Саратове с этого года сосредоточена первоклассная группа математиков...

Помимо привлечения на постоянную работу университет широко применяет практику
приглашения крупнейших специалистов для прочтения отдельных специальных курсов,.
Так, в зимнем и весеннем семестрах 1935/36 уч. года в университете читали лекции профессора
П.~С.~Александров, Л.~С.~Понтрягин и др., познакомившие математиков Саратова
со своими последними научными достижениями....

Одновременно ряд молодых научных работников университет...
направлен университетом в научную командировку в Москву... Один из командированных,
доцент Саратовского университета Н.~Г.~Чудаков, сделал блестящую работу по одному
из труднейших разделов теории чисел — вопросу о распределении нулей римановой
дзета-функции...

Не приходится сомневаться, что все эти мероприятия, с большой энергией проводимые директором 
Саратовского университета Г.~К.~Хворостиным, приведут к возникновению
у нас в Союзе еще одной хорошей математической школы, призванной играть
немалую роль в научной жизни страны. 
\end{quotation}

Если я не ошибаюсь, в биографиях Куроша его работа в Саратове не упоминается, в одной из биографий Хинчина
я видел слова <<Саратовский период>>, в математических биографиях Петровского
упоминаний о работе в Саратове тоже не было. В большой книге Е.~В.~Ильченко \cite{Petrovsky3} о Петровском
сообщается

\begin{quotation}
 И.Г.Петровский выезжает в другие города. В ноябре 1933 г. занимает
 должность профессора математики в Днепропетровском государственном университете, 
 в 1936/37 учебном году в первый раз читает новый курс по теории обыкновенных
 дифференциальных уравнений в Саратовском университете.
\end{quotation}

Работали ли ли Курош и Петровский в Саратове постоянно или <<вахтовым методом>>, сохраняя должности в Москве,
история умалчивает.

31.07.1937 (с марта-месяца 1937 уже много воды утекло) партком Саратовского университета объявил Хворостина виновным во всех
антибольшевистских грехах и постановил:
\begin{quotation}
 За засорение университета троцкистскими
и другими враждебными советской власти элементами, за зажим
самокритики, за связь с врагами народа, за активное пособничество и укрытие от партийной организации Хворостина~Г.~К. …
вывести из состава партийного комитета и исключить из рядов
ВКП(б).
Просить ЦК~ВКП(б) и Саратовский обком ВКП(б) снять Хворостина с поста директора университета.
\end{quotation}

Дальше все было как по писанному.
Второго августа последовал арест, 
следствие (разумеется) установило, что Хворостин был участником <<антисоветской, правотроцкистской, 
террористической, диверсионно-вредительской организации, действовавшей в Саратовской области>>.
Он был осужден выездной сессией Военной коллегии Верховного суда СССР и 21.01.1938 расстрелян.

\sm

В 1956 году  Петровский возбудил ходатайство о реабилитации Хворостина....

\sm

Как писал в А.~Н.~Боголюбов в <<Истории Отечественной математики>> \cite{Stokalo-4-2}, т.3,
\begin{quotation}
 Деятельность московских математиков в Саратовском университете
продолжалась около двух лет.
\end{quotation}

После ареста Хворостина математики поразбежались из Саратова (из упоминавшихся
зав.кафедрами остался В.~В.~Вагнер), поразбежались 
и люди других специальностей (а кто-то и был арестован). Но похоже, что несмотря на происшедший обвал,
усилия Хворостина не полностью ушли в песок, и дальше университет сохранил неплохой уровень...

Вот такая странная  у человека была судьба.

\sm

{\bf \punct Райков.%
\label{ss:raikov-2}}
Я несколько раз встречал высказывание, что в 1933-35 году он находился  в ссылке в Воронеже за троцкизм,
но  не знаю надежного подтверждающего
это источника.
Сам Райков, выступая на Втором Всесоюзном математическом съезде в 1934г., говорил
\begin{quotation}
	Товарищи, в прошлом году, закончив аспирантуру в Москве, я был послан на довольно ответственную работу заведующего кафедрой  в Воронеже...
\end{quotation}

 Из Академической стенограммы:

\begin{quotation}
ЛЮСТЕРНИК. ... и сам Райков сел за работу по теории чисел и читал лекции в Воронежском университете по курсу теории чисел [и функциям действительного] переменного. 
\end{quotation}

Первая опубликованная научная работа Райкова, которую мне удалось найти,
датируется 1936 годом \cite{Raikov-first}, Райкову было 30 лет, то есть в этом 
отношении он превзошел сам\'ого Лузина (но, работа, видимо, была сделана чуть раньше, а пребывание в Воронеже
принесло Райкову большую пользу).

Работа посвящена эргодическим свойствам отображения $z\to z^n$ из единичной окружности в себя,
она, очевидно, являлась продолжением работ Хинчина.
Далее Райков начинает активно публиковаться, выходит несколько работ по теории чисел, по-видимому, под влиянием Хинчина.

Самостоятельной научной фигурой он становится в 1937-38 годах, в связи с работами \cite{Raikov-tomsk},
\cite{Raikov-poisson} по <<арифметике>>
характеристических функций. Программа этой <<арифметики>> была посвящена изучению  разложение
 вероятностного распределения в сверточное произведение (см. книгу Е.~Лукача \cite{Lukach},
 главы 5-6,
 позже эта тема, впрочем, тоже развивалась).
 Естественно при этом перейти к характеристическим функциям, и ставить задачу о разложении
 положительно определенной функции на прямой в произведение положительно определенных функций.
 В принципе, к этому сюжету можно отнести
 и представления Колмогорова, Хинчина, Леви (Paul L\'evy), безгранично делимых распределений. 
 В более широкой постановке задача изучалась в конце тридцатых Хинчиным и Крамером
 (Harald Cram\'er), тогда были получены разные красивые результаты, Хинчин показал, что 
 у данной характеристической функции существует разложение в произведение безгранично делимой функции
 и конечного или счетного набора неделимых функций (так сказать, <<простые множители>>),
 но разложение это не единственно.
 Крамер показал, что если гауссовское распределение  разложить  в произведение,
 то множители обязательно будут гауссовыми. Райков начинает изучение <<арифметики>>
положительно определенных функций,  голоморфных  в полосе $-a<\Im t<b$, в частности, показывает что множители
 голоморфной функции обязательно голоморфны в той же полосе. С помощью этого он получает аналог
 теоремы Крамера для пуассоновских распределений (множители пуассоновского распределения - обязательно пуассоновские).
 
 Работа эта приобрела известность, но главное, после нее, Райков начинает заниматься положительно определенными функциями
 и анализом на группах. В частности доказывает аналог теоремы Бохнера (S.~Bochner) для локально компактных абелевых
 групп (примерно тогда же ее получили А.~Вейль (A.~Weil) и А.~Я.~Повзнер, см. \cite{Gorin-Raikov}). Здесь пути Райкова
 перекрещиваются с путями И.~М.~Гельфанда.
 
В 1937 году М.~Стоун (Marshall Stone) опубликовал работу
{\it Applications of the theory of Boolean rings to general topology,} \cite{Stone}. 
Работа соответствует своему названию, но больше известна тем, что там доказано
известное обобщение теоремы Вейерштрасса о равномерном приближении непрерывных функций. Там, однако, было утверждение
(Теорема 86) о том, что компактные топологические пространства определяются алгебрами непрерывных функций на них
\footnote{Two bicompact [тогдашний термин для слова <<компактный>>] $H$-spaces [хаусдорфовы пространства] 
are topologically equivalent if and
only if their function-rings are analytically isomorphic. One bicompact $H$-space
is a continuous image of another if and only if its function-ring is analytically
isomorphic to an analytical subring of the function-ring of the other.}.
Почему-то в это высказывание вцепились в Москве, лично Колмогоров, И.~М.~Гельфанд (ученик Колмогорова)
и Г.~Е.~Шилов (ученик Гельфанда). Из их интереса воспоследовало понятие нормированного кольца (1939-1940) и далее
$C^*$-алгебры (уже с участием М.~А.~Наймарка%
\footnote{Который, кстати, не происходил из Московской математической школы, 
а был учеником М.~Г.~Крейна, а тот учеником Чеботарёва.}).

В 1940 Гельфанд и Райков \cite{GR-1} вводят алгебру $L^1$ на локально компактной абелевой группе $\Gamma$
и показывают что все гомоморфизмы этой алгебры в алгебру $\C$ происходят
из характеров (одномерных представлений) группы  $\Gamma$ (тем самым групповые алгебры на локально компактных
группах были введены,
и могли использоваться как орудие). 

В 1943г. они же \cite{GR-2} показывают, что у локально компактной группы $G$ есть <<достачно много>> неприводимых
унитарных представлений, в том смысле, что для любого элемента $g\in G$ есть неприводимое унитарное представление
$\rho$ такое, что $\rho(g)\ne 1$. Доказательство: множество положительно определенных функций
на группе компактно, по теореме Крейна--Мильмана (кстати, опубликованной в 1940 году) оно есть выпуклая оболочка
своих крайних точек. Но положительно определенные функции являются матричными элементами 
представлений (отголосок теоремы Бохнера и характеристических функций, впрочем тривиальный),
а крайние точки соответствуют неприводимым представлениям.
Остается заметить, что положительно определенные функции на группе есть, потому что на локально компактной группе есть
мера Хаара и есть представление группы в $L^2$ на себе.

Вход в теорию бесконечномерных
представлений и некоммутативный гармонический анализ был проделан. Но мы уже вышли за временн\'ые рамки этих
записок, а потому вернемся в начало 30х.

\sm

{\bf \punct Выгодский.%
\label{ss:vygodski-2}}
Наивысший взлет его карьеры пришелся  на 1931-1934гг.
В 1931-32гг. мы видели его директором НИИ механики и математики, в
 1932-1934гг. он был директором Объединенного научно-технического издательства.

В 1934 году  опубликовал книгу <<Галилей и инквизиция>>.
В 1935 году   был исключен из Партии. В течении многих лет  утверждалось, что это исключение
было связано с тем, что он получил за эту книгу медаль Папской академии. В статье Демидова, Петровой и Токаревой
2015 года говорится:
\begin{quotation}
Книга получила положительную оценку папского
престола [и далее в сноске:] На заседании Московского математического общества, посвященного памяти Выгодского,
говорили даже о медали Папской академии, которая была присуждена ему за эту книгу.
К сожалению, никаких документальных подтверждений этому пока мы не имеем.%
\footnote{Автору настоящих записок также неизвестно, есть ли подтверждение положительной оценки папского престола.
	\newline
	Из интервью Юшкевича, человека осведомленного, \cite{Yushkevich-Interview}, 1988:
{\it Выгодский был по неизвестной
мне причине исключен из ВКП(б),}  то же самое говорится в явно неподцензурных воспоминания
Юшкевича \cite{Yush-school}.  В первом номере журнала <<Под знаменем марксизма>> 1935г. 
была поносная статья А.~A.~Максимова \cite{Maximov-Vygodsky} по поводу книги Выгодского.
\newline
{\it Он в своей книге совершенно обошел вопрос об об’ективной роли Галилея,
	о том, какому отряду подымавшейся буржуазии служил Галилей... 
	\newline
	Но ничего, абсолютно ничего не усвоил Выгодский из богатой сокровищницы
	 марксизма-ленинизма. Наоборот, он, претендующий быть марксистом, подменил марксизм-ленинизм таким подходом к историческим событиям эпохи Галилея,
	который ничего общего с марксизмом-ленинизмом не имеет и враждебен последнему.}
\newline
Стоит все же иметь в виду, что в этом жанре Выгодский был и сам не промах, а сентенции
уровня вроде цитированной мы слышим ежедневно с небольшой с заменой идеологических штампов.}
\end{quotation}

Так или иначе, исключен он был. Цитируем тех же авторов.
\begin{quotation}
В 1930–1931 гг. развязалась кампания борьбы с егоровщиной, в 1936 разразилось
«дело академика Лузина». Выгодский в начале этого периода — личность
в сообществе очень заметная. Такие люди часто становятся мишенью разного
рода интриг, в этот период интриг с политическим подтекстом....

 В 1935 г. он был исключен из партии «как оторвавшийся от партийной
жизни и в своих работах протаскивающий чуждые марксизму взгляды». Разумеется, он немедленно был удален со всех руководящих постов. Отныне он
был человеком, который в каждой анкете в соответствующей графе должен
был указывать: исключен из партии и обозначать причину исключения. В реалиях
советской жизни он становился человеком, пораженным в правах. И
эта «черная метка» преследовала его всю жизнь. Он был вынужден оставить
механико-математический факультет Московского университета и в 1952 г. переехать в Тулу, где до самого выхода на пенсию в 1960 г. работал профессором
Тульского педагогического института. С 1963 г. и до самой смерти он состоял
профессором-консультантом Тульского горного института
\end{quotation}

Так что и Выгодский оказался жертвой. На самом деле, биографические данные заставляют усомниться
в том, что в 1930х он был <<поражен в правах>>.

В 1938 году защитил докторскую <<Математика древних вавилонян>>. Это была первая 
защищенная к тому моменту докторская
по истории математики (вторая - Юшкевич). В тот момент, по-видимому,  не было и защищенных
кандидатских (хотя люди, получившие за это кандидатскую степень без защиты были).

Согласно <<Летописи Московского университета>> \cite{Ilchenko} в 1939г был
\begin{quotation}
организован научно-исследовательский семинар по истории математики для профессорско-преподавательского состава
и студентов (февр., рук. проф. М.~Я.~Выгодский)%
\footnote{Семинар этот под руководством Выгодского и Яновской существовал с 1932-33 года, то есть речь идет о
какой-то его реинкарнации. В данном контексте важно, что Выгодский в 1939 году был руководителем идеологического семинара.}. 
\end{quotation}

Стоит иметь в виду, что <<история математики>> всегда была идеологическим предметом, а уж в тридцатые годы (см. \S \ref{s:history}) ...

\sm

Согласно справочнику \cite{professora}, Выгодский работал в МГУ  профессором в 1933-1941 и 1945-1948.
С конца 30х он параллельно работал в МИТХТ.

Поиски по библиотечным каталогам показывают, что во второй половине 30х продолжали выходить переводы историко-математической литературы
с участием Выгодского.

В 1941 году выходит первое издание его <<Справочника по элементарной математике>> и  тогда же
книга «Арифметика и алгебра
в древнем мире».

\sm

Это отступление не  только о Выгодском, но и том, как пишется сейчас история математики. В любом случае в 30е годы
в правах он поражен не был (и слава богу), но после 1948г. с ним случилось  что-то еще.
Из воспоминаний Б.~А.~Розенфельда \cite{Rozenfeld}:
\begin{quotation}
  С большим трудом, уже после войны, он восстановился в партии, доказав,
 что работал в мусаватистских учреждениях, выполняя задания Красной Армии, 
 а затем тихо «механически выбыл» из партии. Я часто беседовал с Марком Яковлевичем и он разъяснил мне,
 в чем состоят его расхождения с КПСС и почему он не хочет в ней состоять.
\end{quotation}

\bigskip

А вообще-то, кроме своего не слишком славного директорства,
он сумел в своей жизни совершить много полезного и для математического образования и для математики. 
Его исторические работы не лишены искренних
элементов экстремизма, но это все же это научные исследования. 
Итогом его организационно-издательской деятельности 30х были многочисленные издания старинных  математических
сочинений, а также историко-математических книг (Цейтен, Нейгебаэур,  Клейн, Вилейтнер). 
Он был автором справочников и вузовских
учебников. <<Справочник по элементарной математике>>, изданный в 1941 году, был переиздан в 1948, а в 1986 вышло
его 27ое издание. Дальше номера изданий не считались, но были, например, издания 2012, 2014, 2015 годов...
Эта живучесть связана с исключительной удачностью книги (авторов книг по элементарной математике 
всегда хватало). Еще одна цитата
из Б.~А.~Розенфельда (предисловие к \cite{Rozenfeld-Vygod}): 
\begin{quotation}
 Это, конечно, не простые справочники, а учебные пособия нового типа, очень полезные при самостоятельном изучении
математики.
\end{quotation}

Кстати, автор данных записок придерживается того же мнения.

Да и с его директорством, как мы видели, все было не вполне однозначно.

\sm

{\bf \punct Яновская.\label{ss:yanov2}} У Яновской было много учеников, и о ней больше всего написано.
В связи с этим нам тоже придется остановиться на 
ней подробнее

\sm

a) {\sc <<Математические рукописи Маркса.>>} Несколько цитат из Кольмана. Сначала из Воспоминаний \cite{Kolman-vospominaniya}, опубликованных
в 1982г.:
\begin{quotation}
 В марте 1931 года был арестован Рязанов, директор ИМЭЛ [Институт Маркса--Энгельса--Ленина], и был назначен Адоратский, 
 старый большевик, ученый-марксист. Я был послан туда как член дирекции института.
 Адоратский поручил мне заведование кабинетом Маркса. С первого же дня работники института обратились в дирекцию
 с утверждением, что Рязанов публиковал далеко не все документы, которыми он располагал, что часть из них он, 
 якобы, прятал не в несгораемых шкафах, а в каком-то личном тайнике. И в самом деле, с помощью специалиста,
 вызванного из Угрозыска, этот тайник был обнаружен в стене кабинета Рязанова. В нем хранились фотокопии нескольких
 писем Маркса и Энгельса, а также математических тетрадей Маркса и его же записей по вопросам естествознания.
\end{quotation}

А вот что он же произносил с трибуны в апреле 1931г. \cite{Kolm-boevye-voprosy}:

\begin{quotation}
Специально по математике — по арифметике, геометрии, алгебре, но дифференциальному исчислению — имеется 50 работ Маркса.

Причем это не только заметки к прочитанным книгам, не только упражнения на отдельные темы,
но иногда почти что подготовленные к печати наброски по таким узловым, коренным проблемам математики, 
по которым Маркс, — это уже сейчас видно, — сказал много такого, что не потеряло решающего значения и до настоящего времени.
Труды, которые Маркс изучал по математике, являлись основными фундаментальными работами.
Он прорабатывал Ньютона, его «Метод Флексий», «Анализ», «Математические начала натуральной философии». 
Он прорабатывал Тейлора, Мак-Лорена, Лагранжа, Даламбера. Сейчас мы знаем обо всем этом.

{\bf Все это пряталось, что является преступлением не меньшим, 
чем преступление шахтинцев, которые прятали угольные пласты.} 
\end{quotation}

Снова воспоминания Кольмана:
\begin{quotation}
 Что же касается математических трудов Маркса, то, по-видимому, как я понимаю теперь,
 Рязанов потому тянул с работой над ними, с их расшифровкой и опубликованием, что не был уверен в их научной ценности,
 боялся повредить авторитету Маркса. Немецкий посредственный математик Гумбель%
 \footnote{Emil Julius Gumbel, он опубликовал статью на эту тему в сериале <<Летописи марксизма>> за 1927г.
 	\cite{Gumbel}.
Я, к сожалению, этой статьи  не видел.}, 
 которому Рязанов дал их на отзыв,
 не понял их методологического значения и не советовал издавать их%
 \footnote{По-видимому, были и другие отзывы, см.  \cite{Katolin}.}. Здесь произошло нечто подобное тому, 
 что и с <<Диалектикой природы>> Энгельса. Эйнштейн, которого Рязанов попросил ознакомиться с рукописью, 
 не понял громадного философского значения этих набросков (что объяснялось тем, 
 что Эйнштейн вообще не понимал диалектики), и указал - совершенно правильно - 
 что с чисто физической стороны они устарели, а поэтому рекомендовал воздержаться от их издания%
 \footnote{Лев Католин  в \cite{Katolin} утверждал в точности противоположное
 	(что Эйнштейн поддержал <<Диалектику природы>>). Оба источника не слишком надежны,
 	более правдоподобно известие Кольмана.}.

К чести Рязанова надо отметить, что тогда он не послушался совета гениального физика, 
и <<Диалектика природы>> была в 1929 году издана. Но в случае с математическими рукописями получилось иначе.
\end{quotation}

Рукописи эти представляли записи Маркса 1873–1883 при изучении учебников анализа. 
Набор этих учебников был выяснен К.~А.~Рыбниковым
уже после Войны. Цитируем  рецензию Н.~Х.~Розова 1968 года \cite{Rozov}. 

\begin{quotation}
 К сожалению, Марксу, как это видно из его рукописей, остались неизвестными
новейшие для второй половины XIX века математические теории (К.~Вейерштрасса,
Р.~Дедекинда, Г.~Кантора, и других). Все источники, которыми пользовался Маркс,
не выходят за пределы методологии математики XVIII века. К этому были существенные
причины. Маркс, живший в 70—80-х годах прошлого века (в период, к которому как
раз и относится большая часть его математических рукописей) в Англии, вынужден
был ориентироваться главным образом на те книги, которые были «в ходу» в английских
университетах. А господствовавший там традиционализм, беспрекословное преклонение
перед авторитетом Ньютона создали глухую стену, отгородившую Англию от всех новых
идей математического анализа.
\end{quotation}

Математик, прочитав это, может сделать  умозаключение (и оно будет правильным), что Маркс
не был знаком с  работами
Коши и Больцано по обоснованию анализа, как, впрочем, и к открывавшим математический XIX
век открытиям 
Фурье и Гаусса%
\footnote{Понятно, что архаичность источников должна была быть указана с самого начала людьми,
говорившими, что рукописи не надо публиковать. Это подтверждается словами Кольмана в книге 1936 года:
\newline
{\it Хотя с чисто математической стороны работы Маркса ограничены характером бывших в его распоряжении источников.}
\newline
и продолжение (не могу удержаться)
\newline
{\it 
\dots
и полученные им результаты не могут быть поэтому механически
перенесены на современную математику, Маркс дает такое логическое обоснование дифференциального исчисления, 
которое является вместе с тем и историей его развития. Он
вскрывает подлинный смысл дифференциальных символов и 
сущность того диалектического перехода, который приводит от 
алгебраического дифференцирования к дифференциальному 
исчислению в собственном смысле слова. Выясняя оперативный
смысл основных формул дифференциального исчисления, Маркс
предвосхитил точку зрения, к которой пришли такие видные
современные математики, как Адамар, лишь в наши дни.
\newline
Самая главная, величайшая ценность работ Маркса состоит
в том, что на конкретном примере, на одной из важнейших 
проблем обоснования математики, в них показано применение
метода материалистической диалектики в действии, что
в противовес всем математикам, вводящим новые понятия
чисто формально, здесь благодаря методу, 
устанавливающему единство логического и 
исторического, показана необходимость введения именно этих, а не
других понятий. Тем самым не посредством общих 
рассуждений, а путем показа дано доказательство возможности
перестройки математики на диалектико-материалистической
основе, дан прообраз такой перестройки, а вместе с тем и 
пример того, как должна выглядеть подлинно научная история
математических наук.}
\newline
Чувствуется, были возражения. И нечего было Эрнесту Яромировичу на этот счет сказать.}.
В этом не было ничего особенного, как не было ничего особенного и в отставании
учебников матана от развития науки (они всегда отставали, и это неизбежно).
Забавна история с публикацией и развеличиванием записей, которые сам  Маркс
не  намеревался публиковать. Снова слово Кольману.
\begin{quotation}
Я немедленно подсказал Адоратскому создать бригаду для их расшифровки - были приглашены Яновская и ее ученики,
Райков и Нахимовская, усердно взявшиеся за это каторжное дело%
\footnote{О чудесном  обретении рукописей из тайника  известно только из воспоминаний Кольмана. 
	О том, что Кольман извлек рукописи из небытия в их издании говорится
в первом издании книги \cite{Katolin} Льва Католина	(1973).
Но почему-то об этом не сказала
 Яновская в 1933г. в \cite{Yanov-PZM-2}.}.
\end{quotation}

В 1933г. рукописи были частично изданы в журнале <<Под знаменем марксизма>>
(номер 1) и  переизданы в сборнике
<<Марксизм и естествознание>>. Яновская \cite{Yanov-PZM-2} в связи с этим писала:
\begin{quotation}
	Лишь самому Марксу удалось вскрыть подлинно диалектический характер
	дифференциального исчисления. Ибо он сначала вскрывает {\it происхождение}
	дифференциальных символов, а затем, когда они, так сказать, возникли
	{\it оборачивает} точку зрения и делает исходным пунктом не {\it реальные}
	<<алгебраические>> соотношения, отражением которых служат
	дифференциальные символы, а самые эти символы, рассматриваемые как
	{\it уже данные и готовые}...
\end{quotation}

Определенная неприятность для математиков могла бы случиться, если бы преподавателей матана
на занятиях
обязали бы <<вскрывать подлинно диалектический характер>> дифференциала. Насколько я понимаю,
ничего подобного не произошло,  <<рукописи>> в 30х годах в дело не пошли.

В 1932г. Кольман выступил на Математическом конгрессе в Цюрихе с докладом 
<<Новое обоснование дифференциального исчисления Карлом Марксом>>.
По-видимому, отчет об этом докладе был победоносным, цитирую книгу Льва Католина
\cite{Katolin}, посвященную истории публикации <<Математических рукописей>> Маркса:
\begin{quotation}
 Хотя 
пятница девятого сентября была для конгресса днем, 
рекордным по насыщенности программы, само имя Маркса
в заголовке доклада привлекло множество слушателей.
В помещении, отведенном для занятий не слишком 
популярной на конгрессе и довольно слабо посещаемой
секции философии, стало тесно. Пришлось перенести
заседание в другую, большую аудиторию, которая,
впрочем, тоже оказалась переполненной.
\end{quotation}

Сам Кольман рассказывал так:
\begin{quotation}
 На конгрессе я сделал два секционных сообщения: «О новом обосновании
 дифференциального исчисления Карлом Марксом» и «О
функциях кватернионального переменного». Оба они не остались
без внимания. Первое, конечно, из-за имени Маркса, причем бывшие
на конгрессе русские эмигранты не преминули устроить во время
этого моего сообщения небольшую обструкцию, подняли шум в
аудитории, но их вежливо уняли, 
\end{quotation}

А вот рассказ де Рама (Georges  de Rham) 1978г., переданный  Дюгаком
\cite{Dug}:
\begin{quotation}
  На заседании председательствовал А.~Реймон [Arnold Reymond], профессор философии Лозаннского университета, 
  а в зале было много народу. В конце доклада Кольмана, П.~Бернайс [Paul I. Bernays], 
  профессор Геттингенского университета, 
с интересом следивший за изложением Кольмана, задал ему следующий
вопрос: «А как Карл Маркс брал производную от $\sin x$?». Тогда раздался оглушительный общий хохот,
и председатель заседания, который должен был делать доклад
после Кольмана, объявил: «Господа, так как уже поздно, я сделаю свой доклад на
следующем конгрессе математиков». Однако ввиду протеста из зала, он все-таки
сделал доклад%
\footnote{Серпинский (очевидец), \cite{Dug}: {\it Шеф советской
делегации Кольман говорил на Конгрессе глупости о взглядах Карла Маркса,
касающихся основ дифференциального исчисления.}
Ср. Юшкевич \cite{Yush2}, 1991 {\it Между прочим доклад Кольмана о математических рукописях
Маркса, незадолго до того опубликованных и прокомментированных С. А. Яновской,
был сделан очень неудачно и вызвал недоумение многих участников
конгресса.}
\newline
Колмогоров -- Александрову 30 сентября 1932 г.:
\newline
\it
Затем был второй доклад [Чеботарёва], о съезде, очень бледный (после 
вытягивали из докладчика более интересные ответы на вопросы). Впрочем, 
следующее место имело успех: <<На докладе тов. Кольмана я не был, так как
хотел послушать Мюнтца, докладывавшего одновременно, однако, сам
Мюнтц каким-то образом оказался на докладе тов. Кольмана>>.}.
\end{quotation}

Заданный вопрос о дифференцировании синуса может показаться странным, но на самом 
деле рассказ де Рама весьма правдоподобен.
Процитируем статью В.~Н.~Молодшего 1969г. \cite{Molodshy-Marx} о рукописях Маркса:
\begin{quotation}
 В рукописи о производной Маркс не показывает, как найти его методом производные функций 
 $\sin x$ и $\cos x$. Это, по-видимому, обусловлено тем, что в учебниках, которыми он располагал,
 <<алгебраические>> разложения функций  $\sin x$ и $\cos x$ в степенные ряды не рассматривались.
 Как указывалось выше, такие разложения были осуществлены Лагранжем
 и нашли впоследствии более строгое
 обоснование в некоторых курсах алгебры. Если бы этот факт был известен Марксу, он без труда
 нашел бы с помощью своего метода производные функций  $\sin x$ и $\cos x$.
\end{quotation}

Скорее всего, это отголосок того же, о чем говорил де Рам
(а иначе к чему все это?).

\sm

Так или иначе в 30х годах <<Математические рукописи>> в дело не пошли 
(такой статьи нет в БСЭ, что достаточно характерно, я не нашел упоминания о них и в длинной статье про Маркса,
хотя мог быть и был недостаточно внимателен). Но дело на этом не кончилось.

Снова цитируем Кольмана:
\begin{quotation}
Львиную долю работы над рукописями проделала сама
Яновская, давшая ценный комментарий к ним; они были фактически полностью подготовлены ею к печати. 
Но в ИМЭЛе, после смерти Адоратского[1878-1945], менялись директора,
и все они в вопросе об издании Марксовых математических рукописей занимали одинаково нерешительную,
колеблющуюся позицию страховщиков%
\footnote{Для контраста, цитата из И.~Р.~Шафаревича: <<{\it К несчастью для Маркса,
эти его «математические рукописи» популяризировались
одно время в СССР и стали широко известны. Они имеют все черты безнадежного дилетантизма
и ума, работающего вхолостую}>>.
\newline
Гливенко утверждал \cite{Glivenko}, что Маркс предвосхитил какое-то понимание дифференциала математиками.
Желающие могут изучить этот вопрос (в любом случае статья  Гливенко, прав он или нет, является  рациональным
текстом).}. Потребовались десятки лет, пока рукописи были опубликованы,
причем включившимся в это дело позже бездарным, но ловким карьеристом Рыбниковым,
одним из многочисленных учеников уже скончавшейся Яновской, присвоившим себе всю заслугу
и получившим за это ученое звание доктора.
\end{quotation}

Упомянутый  К.~А.~Рыбников в 1954 году
защитил по этим рукописям докторскую диссертацию (физико-математические науки).
Согласно Католину \cite{Katolin},
\begin{quotation}
Оппонентами
его диссертации выступили такие известные ученые,
как А.~О.~Гельфонд и А.~П.~Юшкевич. Академик
А.~Н.~Колмогоров прислал положительный отзыв.

Но
тем не менее от публикации самих рукописей Маркса
решено было пока воздержаться — в них оставались еще
неясные места...
\end{quotation}

К. А. Рыбников (который ранее работал начальником секретного  криптографического учреждения),
защитив диссертацию, увел у Яновской кафедру <<Истории математики>> Мехмата.
Что касается рукописей, то они в итоге были полностью изданы в 1968г. (через два года после смерти Яновской).
Согласно предисловию ИМЭЛ, рукописи были подготовлены к печати Яновской.

Как писал Н. Х. Розов в уже цитированной рецензии,
\begin{quotation}
 Математики, философы, историки, экономисты получили ценнейший подарок. 
\end{quotation}
Студенты-математики (по крайней мере, на Мехмате) в самом деле получили подарок, а именно 
<<Рукописи>> были включены в программу курса <<Марк\-си\-ст\-ско-ле\-нин\-ской философии>>.
Насколько я помню (примерно 77-78гг) тема эта спускалась на тормозах, то ли в связи определенным
благоразумием преподавателей с философского факультета, то ли в связи с тем, что 
они уже имели возможность ознакомиться с реакциями студентов на диалектическую природу дифференциала.

В этой истории забавно  адекватное отсутствие желания руководства
Института марксизма-ленинизма публиковать рукописи (при всей необходимой осторожности в отношении показаний
Кольмана, стоит отметь, что для публикации рукописей понадобилось 38 лет), и
проталкивание проекта около-математическими партийцами от   снизу. 

\sm

б) {\sc О позднейшей биографии.} Мы видели, что герои 1929-31гг.
Хворостин, Райков и Выгодский в дальнейшем оказались людьми, способными к положительной
деятельности. Похожее случилось и с Яновской. Например, известно, что в 1941-43гг. она, будучи в эвакуации,
работала в Пермском университете (в то время в <<Молотовском университете>>), и привезла
с собой оттуда студентов  Е.~Б.~Дынкина, О.~А.~Олейник, М.~М.~Постникова.
В 1947 году под ее редакцией были изданы «Основы теоретической логики» Гильберта и Аккермана
(Hilbert, Ackermann, Grundz\"uge der theoretischen Logik.),
в 1948  «Введение в логику и методологию дедуктивных наук» Тарского
(Tarski, Introduction to Logic and to the Methodology of Deductive Sciences), что было 
существенным вкладом в развитие логики в СССР... Биограф Яновской Б.~В.~Бирюков пишет
о ее <<духовном перерождении>>. Не берусь обсуждать вопросы духовности,
но в любом случае после Войны она выглядела совсем иначе.
Но мы уже  за  временн\'ыми
рамками этих записок, поэтому вернемся в 30е годы.

В юбилейной статье 1966 года в <<Успехах>> \cite{Yanov-UMN} приведен список сочинений Яновской. 
Первые ее статьи называются  <<Категория количества у Гегеля и сущность математики>> (1928),
 <<Закон единства противоположностей в математике>> (1929),  <<Идеализм в современной философии математики>> (1930).
 Дальше продолжается в том же духе, есть статья с на первый взгляд мирным названием <<Дифференциал>> 
 в соавторстве с Кольманом БСЭ (после кузькиной матери, показанной Шмидту),
 еще одна статья с Кольманом про Гегеля и т.д. 

 Названия резко меняются в 1937 году:
 
 \sm
 
 {\small
--- «Геометрия» Декарта (К 300-летию со времени выхода в свет), Фронт науки и техники, № 6, 25—35.

--- О рядах, Матем. в школе, № 4, 17—21. (Поправки см. в № 6, 127.)

--- Как древние вавилоняне четыре тысячи лет назад вычисляли квадратные корни
Матем. в школе, № 6, 71—73.

--- И. Л. Гейберг, Естествознание и математика в классической древности, Книга
 и пролетарская революция, № 3, 132—134.

--- Леонард Эйлер, Введение в анализ бесконечно малых, Книга и пролетарская революция, № 7, 124—127. }

\sm

До этого историко-математических статей у Яновской 
в списке работ не видно. При всей идеологизированности тогдашней истории математики,
эта наука тогда была более мирной, чем философия,
а с точки зрения 37 года, этот род занятий был и более благоразумен.

Кстати, степень {\bf доктора физико-математических наук}
Яновская получила  в 1935г без защиты диссертации. Работ по истории математики у нее еще не было.
То есть степень была получена за что-то иное: или за диалектику (но это философия), или, скорее, за
<<Математические рукописи>> Маркса, что при наличии желания, можно отнести к истории математики.
Еще раз напомним, что в  1934-35гг сразу после восстановления ученых степеней,
их присуждение получение без защиты было относительно обычным (интересно бы найти полный список таких докторов по математике, боюсь, что Яновская будет выглядеть в нем несколько необычно).

\sm

Впрочем, постепенный отход от образа
диалектической фурии у нее произошел до 37 года.
Я не занимал последовательным изучением ее тогдашних сочинений, но в целом видно,
что со временем ее философствования принимают все менее  воинственный характер. Как мы отмечали выше,
к 36 году мирное значение принимает 
и понятие  <<философия математики>>, а вышедший под редакцией Яновской сборник  \cite{Filosofiya-matematiki}
совсем не напоминает диалектический сборник 1931 г.\cite{Dialektika}).

Цитируем одну из апологий, \cite{Yanov-BDU}.
\begin{quotation}
Страх 30-х (особенно 37-38 гг.) сохранился у С.~А.~Яновской на всю
	жизнь, заставляя её быть осторожной даже в ситуациях, которые по меркам 60-х гг. выглядели абсолютно невинными.
\end{quotation}	
Про 37-38гг. понятно, Яновская как раз принадлежала к слою революционных фанатиков,
который оказался прорежен в 1937 году. Но похоже на истину и то, что бояться она начала раньше:
вряд ли уютно было оказаться одним из адресатов статьи в <<Правде>> за подписью члена ЦК. Да и ее
соратник Райков, быть может, не по своей воле в Воронеж съездил. В любом случае ее поведение
в середине тридцатых становится более мирным.

\sm

в) {\sc Принципиальный товарищеский спор.}
В сентябре 1932 году она опубликовала 12-страничную 
рецензию <<{\it Проблема учебника для втузов еще не решена}>> \cite{Yanov-KIPR}
на учебник анализа Выгодского (тот самый, который
хвалил Лузин, и в котором были актуальные бесконечно малые).
Рецензия выглядит довольно непривычно, это смесь  разумной критики с поисками дьяволов на 
острие иглы. Есть в рецензии  пассаж, где упоминаются и Лузин, и Фихтенгольц, которого еще недавно
громили на Ленинградском математическом фронте (и в придачу Каган, которого Яновская недавно громила лично):
\begin{quotation}
 ... у нас, особенно за последние годы, появилось довольно значительное количество учебников по математике,
 которые совершенно сознательно ставят себе задачей достижение этого органического единства теории и практики,
 сначала либо путем внесения (часто весьма большого как у Фихтенгольца) числа технических, физических
 или механических (вообще конкретных) задач в учебники в основном первого типа [типа Бибербаха,
 полностью изгнавшего из изложения дифференциал]
 или путем придания большей теоретичности учебникам второго (переработка, например, Лузиным учебника
 Грэнвиля, Каганом - учебника Филлипса).
\end{quotation}

Что касается учебника Выгодского, то он 
\begin{quotation}
... {\it наиболее оригинальный, быть может, наилучший из всех}%
\footnote{в этой цитате и ниже выделение  текста - авторское, (н а и б о л е е....)}...
\end{quotation}
Но это за здравие, рецензия  не благостна:
\begin{quotation}
 Первые же главы учебника т.~Выгодского -- не метод интегрального исчисления, а вычисление интегралов степенных
 функций (для целого, дробного и отрицательного показателя), основанное на ряде искусственных
 приемов и уловок.
\end{quotation}

Автор настоящих заметок согласен с тем, что это очень неудачно, но  рука нетривиального автора в книге Выгодского
 все же чувствуется.
Приведем еще несколько цитат.:
\begin{quotation}
 Если бы в следующих изданиях автору удалось рассеять тот туман, которым опутан у него дифференциал,
 и действительно вскрыть суть современных, более совершенных математических методов на основе подлинного разбора
 (а не пропаганды) исторически им предшествовавших, его учебник по праву заслужил бы название нового,
 по существу преодолевшего противоположность между строгим, но формальным, оторванным
 от практики буржуазным учебником одного сорта, и деляческим, эмпирическим, буржуазным же учебником другого
 сорта. {\it В настоящей же своей форме он содержит в себе элементы того и другого в непреодоленном виде}...
 В настоящей свое форме, несмотря на все свои достоинства и успех, он не может быть назван
 шагом вперед в деле марксистско-ленинского построения учебника по математике.
\end{quotation}

Впрочем,  это вполне джентельменский спор двух философов, учебник Выгодского открывается эпиграфом из 
<<Материализма и эмпириокритицизма>> (кстати, удачно подобранным),
а в  предисловии, среди прочих рассуждений на марксистские темы, говорится: 
\begin{quotation}
Я думаю, что изложенная выше установка, определяющая лицо этой
книги, является правильной с точки зрения методологии
{\it марксизма - ленинизма.}
\end{quotation}

Вот Яновская эту правильность диалектически  отрицает. И концовка:
\begin{quotation}
 .... раньше всего автору следует отказаться от своей манеры обучать не средствами серьезного научного разбора и критики,
 а путем пропаганды {\it положительно неправильных} с научной точки зрения методов, построенных на упрощенчестве и извращения 
 подлинного существа дела.
\end{quotation}

Кстати, {\it положительно неправилен} - это Яновская цитирует предисловие Выгодского, а Выгодский в этом 
словосочетании цитировал Энгельса.

У учебника были издания 1931, 1932, 1933. Больше он не издавался.
А то, что этот автор вскоре  научился писать хорошие учебники, мы знаем.

\sm

г) {\sc Яновская и дело Лузина.} В 1936 году Яновская, конечно не могла остаться в стороне.
Она было основным выступающим на митинге в МГУ 9 июля 1936 в связи с опубликованной в тот же день статьей
в <<Правде>>. Все, что
мы сейчас знаем -- текст  в анонимной статье в журнале
«Фронт науки и техники». 1936. № 9, \cite{Front1936-1}. Приведем сообщение полностью:
\begin{quotation}
	\small
 В своем докладе С.~А.~Яновская отметила, что присутствующие на данном собрании
могли бы очень многое прибавить к тому, что писалось в «Правде» о «деятельности» и
личности H.~H.~Лузина.

Печатая все свои оригинальные работы заграницей и помещая в советских изданиях
лишь малоценные статьи, издеваясь при этом над собственными похвальными отзывами
и работами, помещаемыми в советских журналах, лицемерно льстя в глаза советской
научной молодежи и сообщая «по секрету» друзьям, что время этой молодежи подходит
к концу, — Н.~Лузин думал, что ему удастся долго одурачивать нашу научную общественность.

Он действовал бесцеремонно, нечистоплотно, вредительски, рассчитывая на полную
свою безнаказанность. Нечистоплотное отношение Лузина к чужим работам выразилось, в частности, в плагиате у своих учеников (у т.~Новикова). {\bf Возмутительное вредительство его сказалось при переработке известного учебника Гренвиля по математике.}
Переработка свелась к тому, что текст подлинника в 450 стр. вырос до 750—800 стр. В
первой части был еще сохранен до некоторой степени систематический порядок, который имелся в учебнике Гренвиля, но во второй части изложение ведется таким образом,
что оно может дезориентировать читателя. Вся переработка книги пестрит дефектами и
ошибками.

Когда редактора ОНТИ почтительно напоминали Лузину о необходимости исправить
ошибки, вплотную заняться вторичной переработкой отдельных мест учебника, то получали
подобного рода ответы:
«Видите ли, я переезжаю на новую квартиру, у меня протекает потолок, и я не могу
заняться этим учебником».

В своем докладе тов. Яновская приводит между прочим такой характерный штрих:
В 1930 г. H.~H.~Лузин председательствовал на том собрании ученых, которое приняло
обращение к французским ученым — протест против интервентов — в связи с делом
{\bf Промпартии}. Но Лузин уклонился от того, чтобы собственноручно подписываться под
этим воззванием. Напрасно тогда {\bf молодой аспирант Рабинович} стучался в двери Лузина.
Узнав, что он пришел из института математики за подписью, Лузин заявил, что болен,
что ни принять, ни подписать обращение не может.

Вот эта манера — демонстрировать свое «советское лицо» перед советской общественностью
и сохранить свое подлинное лицо перед заграницей — весьма характерна
для Лузина. Это двурушничество так откровенно выпячивалось, что странно было его
не замечать.

Лузин тогда же под первым попавшимся предлогом уходит с математического факультета

И как констатирует в своем докладе т. Яновская, особенно после ухода Лузина из
университета, как из рога изобилия посыпались похвалы, расточаемые им всяческим
ничтожествам математической науки, похвалы, граничащие с издевательством и вредительством.

Без всякого урока для себя Н.~Н.~Лузин возвращается в 1935 г. в университет. Профессора
и студенты встречают его хорошо. Правда, некоторые профессора (например,
проф. Александров) нерешительно указывают на деканате, что Лузину поручаются курсы
не по специальности, что это будет только вредно для курсов. Но все осталось так, как
хотел Лузин.
\end{quotation}

Нельзя сказать (во всяком случае на основании этого текста), что Яновская сколько-либо выделялась бы
на фоне математиков. Ее позиция близка к позиции Шнирельмана, Люстерника,  Гельфонда (что не удивительно).
Выше выделены три фразы. На заседании Академической комиссии 
слово <<Промпартия>> появилась у Горбунова, зачитывавшего сообщение Шнирельмана,
а про <<аспиранта Рабинович>> в тот же день (видимо, с утра)  вспомнил Гельфонд.
% (Александров присутствовал
%на обоих мероприятиях, судя по репликам 7 июля Академическая комиссии заседала с утра). 
Про Гренвиля-Лузина вспоминал 11 июля Люстерник, в том же стиле, но без слов <<вредительство.>>
Но <<приносит вред>> математики тоже произносили, и без  слова <<вред>> Лузина не смогли бы исключить
из Академии.

\sm

д) {\sc Апологии и реальность.} Процитируем одну из апологий \cite{Levin-2011}--\cite{Levin-2012},
впрочем несколько радикальную:
\begin{quotation}
	\small
Уже в 1943 г. С.~А.~Яновская организовала первый в СССР семинар по математической логике при МГУ,
которым она руководила совместно с И.~И.~Жегалкиным 
и П.~С.~Новиковым (а позже – вместе с А.~А.~Марковым)...

 Важно отметить, что С.~А.~Яновской удалось то, что не удается большинству ученых 
 даже с высоким уровнем научных работ: она стала основателем двух крупных научных школ:
 советской школы истории математики и советской школы математической логики.
 
 Мы видим, что у С.~А.~Яновской имеются очень большие заслуги перед наукой
 и высшей школой страны. Почему же страна не воздала ей должное – не присвоила
 ей никаких почетных научных и педагогических званий, не наградила премиями,
 не избрала в Академию наук? Ответ очень прост: Софья Александровна
была исключительно скромным человеком и никогда не только не предпринимала
каких бы то ни было усилий, чтобы такие звания и премии получить, но даже не
заикалась об этом (в отличие от большинства современников своего круга, активно
занимавшихся такой деятельностью). В этих условиях для ее заслуженного награждения
требовалось лишь вмешательство одного-двух авторитетных ученых-
математиков, осознающих важность сделанного этим человеком и готовых совершить Поступок.
К сожалению, таких людей не нашлось, хотя многие крупные
ученые (П.~С.~Александров, А.~Н.~Колмогоров, А.~А.~Марков, П.~С.~Новиков, И.~Г.~Петровский,
С.~Л.~Соболев и др.) высоко ценили ее профессиональные достижения и
относились к ней с глубоким уважением. В связи с этим заметим, что, когда такие
люди находились, подобные вопросы решались достаточно просто даже в эпоху
Сталина. Например, когда в 1946 г. на общем собрании АН СССР, обсуждавшем
прием новых членов АН, президент АН академик С.~И.~Вавилов заявил с трибуны,
что ему стыдно, что он – академик, а Л.~Д.~Ландау – нет, вопрос о членстве Ландау
в АН СССР был решен в течение считанных минут. 
\end{quotation}

Яновская действительно занималась пропагандой логики и ее организационной поддержкой,
и нет оснований сомневаться в разных высказываниях о сыгранной ей в этом качестве положительной роли.
Но рассказы о ее значении представляются несколько преувеличенными. 
Логические исследования велись в России и до Революции, произносимые имена \cite{Gastev}  говорят нам не много
(во всяком случае,  автору этих записок), однако, кроме Жегалкина, в стране были и другие люди,
интересовавшиеся логикой. Решающей предпосылкой
к возникновению советской логической школы были работы Лузина и его учеников по дескриптивной теории множества.
Логические работы Колмогорова и
Гливенко появляются уже в 20е годы, в середине 30х появляются первые логические статьи А.~И.~Мальцева
(он был учеником Колмогорова). Новиков получил свои первые  знаменитые результаты по логике в начале сороковых.
Жаловаться на недостаточную развитость логики в Москве того времени как-то смешно, да и огромное научное
влияние Яновской на этих
деятелей вызывает определенные сомнения (более вероятно, что кто-то из левых математиков
сумел оказать на нее вразумляюшее воздействие, наиболее вероятен в этой роли Гливенко). 

Не все так просто было и семинаром, и с курсами логики. Во всяком случае, Ю.~А.~Гастев и Ю.~Т.~Медведев
\cite{Gastev} говорят следующее:
\begin{quotation}
Длительный процесс накопления логической культуры,... ведет свое начало от {\it основанного в начале 30-х годов
И.~И.~Жегалкиным семинара.}
\end{quotation}

Из статьи С.~И.~Адяна и В.~А.~Успенского о Новикове:
\begin{quotation}
	В 1943 г. на механико-математическом факультете МГУ был организован первый в нашей стране
	научно-исследовательский семинар по математической логике. С этого времени до конца жизни П.~С.~Новиков бессменно
	был одним из руководителей этого семинара...
\end{quotation}

Ершов писал \cite{Ershov}:
\begin{quotation}
В разное время это семинар возглавляли И.~И.~Жегалкин, П.~С.~Новиков, С.~А.~Яновская, А.~А.~Марков
\end{quotation}

Представляется, что роль Яновской была несколько преувеличена ее благодарными учениками и наследниками
(хотя сомневаться в том, что эта роль после Войны она была в целом положительна, нет оснований).

Что касается создания школы истории математики, то  комментаторы, наверно,
правы. Стоит лишь отметить, что первые работы Яновской по истории математики, как мы отметили, относятся 
к 1937г., а Выгодский занимался этим (как  основной научной профессией к 1930г.). 
О самой же <<истории математики>> и ее роли в истории речь пойдет в следующем параграфе.

\sm

{\bf \punct Лейферт.}
В 1979 г. С.~Г.~Михлин \cite{Mihlin} сделал выписки из документов архива Ленинградского
университета, относящиеся сначала к физико-математическому, а потом и
математико-механическому факультетам с 1918 г. по 60-ые годы. В \cite{Erm3}
сообщается, что Михлин не был взят в аспирантуру благодаря Лейферту 
(так что Лейферт должен быть предметом его отдельного интереса).

Лейферт (и, кстати, А.~Р.~Кулишер) впервые появляются в этих выписках  26.05.1924. О том, как безграмотный человек,
 оказался на факультете, история умалчивает, скорее всего это случилось под каким-то давлением извне.
Согласно выпискам Михлина, в весеннем семестре 1930г. Лейферт заведовал кафедрой аналитической геометрии в ЛГУ,
в осеннем - кафедрой приближенных вычислений. В октябре 1931г. список профессоров по математике такой:
{\it Фихтенгольц, Мюнтц, Франк~М.~Л., Смирнов~В.~И., Кошляков, Лейферт (!), Гюнтер,
	Акимов-Перетц~Д.~Я.(?)} (восклицательный знак, наверно, от Михлина, вопрос от публикатора выписок).

Дальше ограничимся несколькими цитатами. Статья Н.~С.~Ермолаевой 
\cite{Erm3}:
\begin{quotation}
	Прошло немного времени, и вот на заседании партийно-комсомольской фракции Общества математиков-материалистов
	21 февраля 1932 г. собравшиеся с изумлением узнали, что Л.~А.~Лейферт снят со всех руководящих постов за то,
	что только прикрывался диалектическим материализмом>>, проводя <<поверхностную критику установки буржуазных
	математиков>>, что им <<в самом толковании диамата допущены грубейшие ошибки>>, что по существу 
	<<никакой работы не было, за исключением брошюры <<На Ленинградском математическом фронте>>, 
	которая по своему содержанию представляет весьма скудную в идеологическом отношении почву>>,
	а также за допущение «левацких загибов>> и за то, что <<борьба с буржуазно-реакционной частью профессоров
	носила декларативный характер>>.
	
	Оказывается, Лейферт был снят директоратом и партийной ячейкой Института естествознания
	при Ленинградском отделении Комакадемии, где он работал по совместительству так же,
	как и в Педагогическом институте.
\end{quotation}

Из статьи Н.~Я.~Виленкина \cite{Vilenkin}:
\begin{quotation}
	Согласно спискам участников 
	Второго Всесоюзного съезда математиков, 
	Л.~А.~Лейферт в 1934 г. уже проживал в 
	Ростове-на-Дону, и вряд ли этот переезд был 
	добровольным. В дальнейшем он поменял 
	несколько городов и в конце 30-х годов 
	оказался в Воронежском университете, 
	где продолжал заниматься травлей ученых. 
%	По дошедшим до автора рассказам, Лейферт 
%	пытался преследовать известного советского 
%	геометра Н. В. Ефимова, работавшего тогда 
%	в Воронеже. Однако ночью накануне  
%	собрания, где должно было разбираться дело 
%	офицерского сына Н. В. Ефимова, Лейферт 
%	был арестован и, судя по всему, расстрелян.	
\end{quotation}

Воспоминания Е. Н. Ефимовой (дочери знаменитого геометра и декана Мехмата МГУ Н.~В.~Ефимова,
 \cite{Efimov} (версия описываемой истории есть у Н.~Я.~Виленкина \cite{Vilenkin}):
\begin{quotation}
	Рассказал мне дядя Боря%
	\footnote{Левин Борис Яковлевич (1906-1993,  позднее работал в Харькове,
		автор знаменитых работ по целым функциям), 
	 учился в Ростове вместе с Ефимовым. В 1931 или 32г. Ефимов переехал в аспирантуру в Москву,
	а Левин остался в аспирантуре в Ростове. Скорее всего, он имел там возможность лично наблюдать
	Лейферта.} и историю
	о папе, которую я не знала. Как только Борис Левин узнал, что некий Лейферт,
	после попытки разгромить ленинградских математиков, был направлен в Воронеж, он написал папе письмо
	и просил папу быть осторожным. Папа письмо получил и все понял.
	
	Лейферт решил опубликовать как свою научную работу статью с названием <<О проведении общего перпендикуляра
	к двум скрещивающимся прямым>>. Она попала на рецензию к Ефимову. 
	Николай Владимирович написал в отзыве <<Это была бы хорошая задача для контрольной 
	работы на первом курсе, но печатать такую работу невозможно, потому что это было известно еще Пифагору>>. 
	Лейферт настаивал и добился того, что возникло дело о травле советского математика Лейферта
	со стороны профессоров Ефимова и Фабриканта (Фабрикант, видимо, тоже дал отрицательный отзыв на эту работу).
	Было назначено собрание%
	\footnote{В принципе на рассмотрении на собрании ничего фатального не было, за исключением 
		известного временн\'ого промежутка, см. дату чуть ниже.}, на котором должны были рассматривать это дело. Папа и профессор Фабрикант пришли.
	Ждали долго, а собрание все не начинается. Мимо проходил кто-то из парткома,
	спросил их, что они тут делают, и только тут они узнали, что Лейферт арестован, 
	а собрание не состоится.
\end{quotation}

Лейферт был арестован 4 февраля 1938 г. по обвинению в участии в право-троцкистской террористической диверсионно-вредительской организации.  Расстрелян в г. Воронеж 22 апреля 1938 г.
Реабилитирован 22 мая 1958 г.
<<за отсутствием состава преступления>>.

%{\bf \punct Ленинградские левые.}
%Мрочек, Вацлав Ромуальдович (1879-1937).
%С сайта
% https://socialist-revolutionist.ru
 
 %\begin{quotation}
 %С 1930 по 1937 годы состоял в штате отделения прикладной астрономии Научного института им. П.Ф. Лесгафта, 
 %выполняя там ряд исторических исследований по заданию Н.А. Морозова.
 %С 1925 по 1930 годы Мрочек возглавлял «Кружок по истории и методологии точного знания»,
 %позднее влившийся в Общество математиков-материалистов при Комакадемии.
 %С 1930 по 1931 годы он входил в президиум этого Общества. 
 %В 1931 году организовал и возглавил в Доме ИТР им. В.М. Молотова секцию марксистской истории %техники (СМИТ). 
 %В 1934 году его назначили заместителем директора по учебной части организованного
 %в это время университета Истории науки и техники (при доме техпропаганды НКТП).
 %В 1933–1934 годах Мрочек был заместителем председателя комиссии по технической математике в Академии Наук.	
 %	\end{quotation}
 
 %Арестован 05.08.1937, расстрелян 27.08.1937	

%Двое были расстреляны в 37-38гг (???), Милинский умер в тюрьме в 1942г.
%Кулишер и Дрозд упоминаются в литературе как репрессированные,
%но в соответствующих базах их не видно (о Кулишере - ниже). 

   \section{Историки и математики%
   \label{s:history}}
   
   \COUNTERS
   
      \epigraph{Поэт в России должен жить долго, чтобы пережить всех мемуаристов.}
      {Федор Сологуб}
   
   \epigraph{...история, писавшаяся этими господами, ничего иного, кроме  политики, опрокинутой в прошлое, не представляет.}{М. Н. Покровский,
   председатель Коммунистической академии и ректор Института Красной профессуры.}

   Читатель, возможно, обратил внимание, что часть московских деятелей 1929-1931гг
   были участниками Гражданской войны (Хворостин, Яновская, 
   Выгодский, Хотимский, Кольман)%
   \footnote{Люстерник в \cite{Mints}   писал в воспоминаниях о Яновской:
   \newline{\it
   Она в частности, назвала руководителей Одесского подполья
 в годы гражданской войны, своего мужа  И.~И.~Яновского  и известного математика
 А.~Я.~Хинчина.}
\newline
   К известию стоит отнестись с осторожностью, потому что
   по опубликованным биографиям
   Хинчина во время Гражданской войны
   его видно в Иваново-Вознесенске и Нижнем Новгороде.}, и, скорее всего для них, как и для многих других людей,
  Год великого перелома был продолжением все той же Гражданской войны, лишь в новых условиях. Но отметим в числе героев
  этого года
   двух историков математики, Выгодского и Яновскую (кстати некоторые из ленинградских ультра-левых тоже ей были не чужды).
   
   \sm

   {\bf\punct История математики.%
   \label{ss:istmat}}
   До появления  в 30е годы диалектиков история математики как самостоятельная профессия в России была,
   по-видимому, представлена лишь В.~В.~Бобыниным  (1849—1919). Но это не значит, что математики не занимались историей
   своей науки в качестве части профессиональной  деятельности. Труды на эту тему в 20-30е годы
   оставили, в частности, А.~В.~Васильев, А.~Н.~Крылов, Д.~М.~Мордухай-Болтовской, Н.~Н.~Лузин, Н.~Г.~Чеботарёв. 
   В послевоенном Советском Союзе историко-научные  книги публиковали, в частности, Б.~В.~Гнеденко, В.~В.~Голубев,
   Б.~Н.~Делоне, Б.~А.~Розенфельд,  В.~Ф.~Каган, А.~И.~Маркушевич,
   в сравнительно недавнее время -- С.~Г.~Гиндикин, С.~М.~Тихомиров, В.~Г.~Мазья, А.~Н.~Ширяев, М~.И.~Монастырский, Д.~В.~Аносов
     (автор пишет эти перечни навскидку,
   и едва ли они полны). Многие математики (этого уже не счесть и не перечислить)
   публиковали историко-научные статьи. Сказать, что наши математики обделяли интересом историю своей науки,
   нельзя.

   Так или иначе, история математики как самостоятельная профессия появилась в СССР лишь в начале 30х в лице Выгодского.
   Профессия эта непростая, вроде как человек должен знать математику (а  это такая наука, 
   что ее и сами математики-то не знают), и одновременно быть историком,
   что тоже является сложной специфической формой профессиональной деятельности. Очень сложной зоной оказывается история математики
   XX века. Мне известен единственный советский/российский историк, который ей серьезно и интересно занимался 
   Ф.~А.~Медведев (1923-1993), его биографию см. в \cite{Zaitsev} (но и его исследования относятся
  к  эпохе Лебега).
   
   Так или иначе, школа профессиональных историков математики
   в СССР пошла от Выгодского и Яновской в 30е годы.
   Но даже по ее основателям видно, что это была не совсем историко-математическая школа.
   
   \sm

   {\bf \punct История и диалектика.%
   \label{ss:ist-dialektika}}
   Из статьи Яновской 1930 года \cite{Yanov-PZM-1}:
   \begin{quotation}
   А аcпиранты-математики, если не вслух, то потихоньку называли курс истории и философии естествознания, организованный для
   для них, «красным богословием», а старые профессора т. н. «Московской
   школы», авторитет которых в среде математиков был несокрушим, прилагали
   все усилия к тому, чтобы спасти «автономию» «чистой», «единой» для «всех
   времен и народов» математики от злостных покушений на нее со стороны
   материалистической философии, не стесняющейся открыто заявлять о своей
   партийности и классовом, пролетарском характере.
   \end{quotation}
   
   История математики (и вообще история естественных наук) была предметом  идеологическим%
   \footnote{Это признавали позднее в 2006г. С.~С.~Демидов и Т.~А.~Токарева \cite{DemidovTokareva},
   {\it<<идеологические установки нового общественного строя, возводившие историю науки, а следовательно
   и историю математики в ранг почти идеологической дисциплины.>>}
   }, кто и когда в точности
   это придумал - неизвестно.
   % (в числе лиц, начавших идеологизацию 
   %истории наук, был Шмидт).
   
   Но вообще это делалось на весьма высоком государственном
   уровне. В 20х годах при Академии работала <<Комиссия по истории знаний>> во главе с Вернадским>>. В 30ом году
   она была преобразована, руководителем стал аж сам академик Н.~И.~Бухарин, а <<в бюро по предложению Бухарина>>
   вошли
  А.~М.~Деборин, В.~И.~Вернадский и А.~Ф.~Иоффе (то есть два представителя Партии и два представителя науки).
  В начале 1932 года она была преобразована в  <<Институт истории науки и техники>>, список его последовательных
  директоров (1932-1937гг) говорит сам за себя: Бухарин%
  \footnote{Если, кто не помнит - один из лидеров Советского государства второй половине двадцатых
  	(вместе со Сталиным и А.~И.~Рыковым).
  В момент  своего избрания в Академию (на каковое академики согласились под нешуточными угрозами) еще был таковым.
  Время активной работы Бухарина в Академии почему-то совпало с наихудшим для академической науки
  промежутком за все время Советской власти и  одним из двух реально плохих моментов. 
  Возможно, что эти события напрямую не связаны, но отметить их временн\'ое
  совпадение необходимо. Задача, требующая отдельного исследования, -- понять была ли
здесь причинно-следственная связь.}, уже упоминавшийся нами А.~А.~Максимов
(своего рода Кольман-light), В.~В.~Осинский-Оболенский [кандидат в члены ЦК ВКП(б)].
  Впрочем,  институт издавал  разные книжки, и в том числе и вполне качественные.
  
%  Как видно из слов Яновской, идеологизация истории науки началась раньше 1930года.
%  Одним из ее инициаторов был Шмидт. Согласно
  
%\begin{quotation}  
%  Выступая на Общем собрании Коммунистической академии 29 ноября 1924 г.,
%Шмидт подчеркнул значение систематического изучения истории науки для 
%разработки узловых вопросов методологии и для материалистического понимания 
%сущности самой науки и указал на необходимость оформления истории науки в 
%самостоятельную научную дисциплину. По его предложению весной 1925 г. в 
%Коммунистической академии была создана Секция естественных и точных наук, с целью - 
%«содействовать построению чисто материалистической системы знаний на основе 
%диалектического материализма»
%\end{quotation}
  
 % Не ясно, был ли Шмидт изобретателем этой идеи.
 % Само выступление Шмидта было в струе тогдашнего политического момента:
 % 2 января 1925 года вышло постановление
 % ЦК ВКП(б) «Об общественном минимуме и пропаганде ленинизма в вузах» (я не нашел текста).
 % Согласно Летописи Московского университета \cite{Ilchenko}, в
  % том же  январе на медицинском и физико-математическом факультетах открылись кафедры диалектического материализма.
  %и по-видимому тогда же начали читаться циклы общественных наук.
  %Под 3 сентября 1928 сообщается
  %\begin{quotation}
  %Отменены курсы общественных дисциплин: политической экономии; истории революционного движения; истории партии.
  %В качестве обязательного оставлен курс исторического материализма. 
%\end{quotation}
  
На вторую половину 20х
 приходится начало внедрения идеологических курсов в вузы
(история эта автору не известна, и она, как и вся прочая тогдашняя
 история, была сложной и сопровождалась приливами
и отливами).
Иделогизированная история естествознания была частью этого процесса.
  
  Известно, что    Шмидт какое-то время читал  курс истории
  естествознания в Московском Университете, но  при  всей
  своей политизации, он со своим ес\-тест\-вен\-но-на\-уч\-ным фанатизмом вряд ли мог тянуть на богослова.
 
  С другой стороны, цитируем воспоминания Юшкевича \cite{Yush-school} о том же 1925 годе: 
  \begin{quotation}
   {\bf В 1925 г. в Московском университете был создан
семинар «Введение в историю и философию естествознания»; 
зачет и экзамен по этому предмету были обязательными.} Кадры 
специалистов в данной области еще только готовились, и {\bf руководила
семинаром С.~А.~Яновская}, тогда слушательница 
Института красной профессуры.
  \end{quotation}

  Яновская потом (начиная, с 1937г.) в самом деле занималась исследованием истории математики. То, что она могла вещать в году 1925
  (и еще долго после этого)
  могло быть лишь <<богословием>>, см. образцы в \cite{Yanov-BSE} (отрывки выше в п. \ref{ss:schmidt-padenie}), \cite{Zasedanie}, \cite{Kolman-differential}.
  Скорее всего,  предметом нежной любви аспирантов была она лично.

%  \begin{quotation}
%  Организация семинара отвечала стремлению молодежи
%  факультета осмыслить положение в тогдашней математике.
%  Необходимо было исследовать сложный комплекс
%  методологических проблем и оценить противоборствующие течения с позиций материалистической диалектики.
 % Для Софьи Александровны эти задачи тогда стали основными.
 % Не давая быстрых готовых решении, она привлекала
 % участников семинара к их поиску, направляла обсуждение
  %по методологически правильному пути. Софья
  %Александровна вела семинар с большим тактом, подкупающей
  %мягкостью и вместе с тем решительностью, когда
  %речь шла критике идеалистической трактовки тех или
  %иных вопросов. Семинар имел большой успех, занятия
  %проходили увлекательно. 
  %\end{quotation}	

  \sm
  
  {\bf\punct О первых ученых степенях по истории математики.%
  \label{ss:ist-stepeni}}
 Цитируем  А.~Н.~Боголюбова, Б.~А.~Розенфельда и Юшкевича
  \cite{Ist-shtokalo}:
  \begin{quotation}
  В это же время [довоенное] организуется подготовка специалистов по истории
математики. В МГУ защитили кандидатские диссертации
первые 
аспиранты: ученики Выгодского — Г.~Б.~Петросян, С.~Е.~Белозеров и 
ученики Яновской — К.~А.~Рыбников, Э.~Я.~Бахмутская; первые докторские
диссертации по этой специальности были защищены М.~Я.~Выгодским
и А.~П.~Юшкевичем.
\end{quotation}

Выгодский стал кандидатом в то время, когда степени давали без защиты диссертации в году за реальные,
хотя и странные,
работы по истории математики. Яновская в те же времена стала доктором физико-математических
наук, скорее всего, за <<Математические рукописи>> Маркса,
Но никак не за историю математики.

Поиск по современным базам и \cite{Diss-MGU} дает об остальных лицах следующие сведения:

\sm

{\small
Белозеров Семен Ефимович (1904-1987), защитил кандидатскую 26.09.1939 
по теме <<О некоторых вопросах из истории функций комплексного переменного>>,
{\bf ректор} Ростовского университета в 1938 - 1954 гг.

Петросян Гарегин Бахшиевич (1902-1997) , {\bf ректор} Ереванского государственного университета  1938—1941
Согласно Википедическим данным, он действительно историк математики, но аспирантуру в Москве
кончал по теоретической физике,  степень кандидата получил в Армении в 1935 году без защиты диссертации.

Рыбников, Константин Алексеевич (1913-2004), защитил диссертацию 
<<Ранняя история вариационного исчисления>> 26 июня 1941 года,
в 1945-1953 начальник какого закрытого криптографического учреждения

Бахмутская Эсфирь Яковлевна (1916-1972), диссертации в библиотечных каталогах РГБ, РНБ, МГУ, а также в \cite{Diss-MGU} не обнаружено
(что ничего не означает).
}

\sm

В первых двух случаях обратите внимание на скорость карьерного роста!

\sm

{\bf \punct Забытый историк математики.} В связи с историей математики в 1930гг. в последние десятилетия
упоминаются Яновская и Выгодский. Оказывается, что  в 1960х годах признавался еще один историк.
Цитируем статью Б.~М.~Кедрова (один из ведущих советских философов), Б.~А.~Розенфельда
(известный математик и историк математики) и А.~П.~Юшкевича (историк математики) \cite{Kedrov-Kolman},
посвящённоё 75-летия со дня рождения <<известного философа и историка математики>>:
\begin{quotation}
	Первая историко-математическая работа Э. Кольмана «Божественная эволюция геометрической мысли» была опубликована в журнале «Естествознание и марксизм» в 1929 г. В 1930 г. в том же журнале он публикует статью «Современные задачи математиков и фи\-зи\-ков-ма\-те\-ри\-а\-лис\-тов-диалектиков», а в 1931 г. в журнале «Физика и химия» — статью «Задачи математики в социалистическом строительстве в реконструктивный период».
	
	В 1931 г. Э. Кольман участвует в работе Международного конгресса по истории науки в Лондоне и выступает там с сообщениями о современном кризисе математики и основных линиях ее реконструкции и о неопубликованных работах Карла Маркса по математике, естествознанию, технике и их истории...
	
	В 1932 г. выходит советское издание книги одного из основоположников топологии И. Б. Листинга «Предварительные исследования по топологии» в переводе Э. Кольмана с его большой вводной статьей. Эта статья была первым исследованием по истории одной из важнейших ветвей современной математики.
	
	В книге «Предмет и метод современной математики» (1936) Э. Кольман широко использовал свою богатую эрудицию в области истории науки. Здесь, в частности, он предложил известное деление истории математики на три периода в соответствии с достигнутой степенью абстракции, считая в математике XVI—XVII вв. главным переход к переменным символам, а в новейшей математике — выдвижение на первый план переменных операций.	
\end{quotation}

Я воздерживаюсь от комментариев к упомянутым сочинениям и  цитирования (читатель, открывший их, узнает <<руку мастера>>).  Отмечу лишь, что Кольман после Войны был одним из штатных историков математики,  и вообще он <<автор
более 400 печатных работ, в том числе 50 книг и брошюр>>. Так или иначе, эмигрировав в 1976г. в  Швецию,
он не стал объектом дальнейших историко-математических апологий, а в 1991г. ему предстояло стать козлом историческо-математического отпущения.

  \sm
  
  {\bf\punct Социальная история отечественной математики.%
  \label{ss:soc-hist}} 
  Слабо-иде\-о\-ло\-ги\-че\-ский характер предмет <<История математики>> сохранял до 80х годов (впрочем, к тому времени
  уже давно идеологию
  было принято проповедовать не без фиги в кармане. Что случилось с этой наукой дальше?
  Процитируем статью С.~С.~Демидова и Т.~А~Токаревой
   \cite{DemidovTokareva}, посвященной 100-летию Юшкевича: 
  \begin{quotation}
  В конце XX - начале  XXI века социальная история математики становится одним
  из приоритетных направлений исследований ... школы.
  \end{quotation}

  <<Социальной истории математики>> в самом деле посвящены многочисленные статьи,
  социальной историей науки занимается Институт истории естествознания и техники им. С.~И.~Вавилова РАН,
  его публикации, попадающиеся на глаза, посвящены в основном этой теме.
 
 Встает вопрос о том, что входит в понятие <<социальная история>>? Математический социум, как и любой другой,
  обладает определенными внутренними нравами, нельзя сказать, что очень симпатичными. Понятно, что его
  история полна самых разнообразных склок и клановых разборок. Математический социум 
  определенным образом действовал  по отношению к внешним социумам, включая власти.
  Надо сказать, что  действия эти бывали достойными и успешными. Здесь
  можно вспомнить и участие математиков в развитии высоких технологий, начавшееся с Жуковского,
  Чаплыгина, Голубева еще в 1920е годы. Можно вспомнить и хорошо поставленное математическое
  образование, не только профессиональное, но также школьное и вузовское с конца 30х по начало 70х годов
  (над чем много потрудились  и педагоги, и профессиональные математики).
  Оно имело широкое культурное воздействие, последствия которого не вполне изжиты до сих пор.
 С другой стороны, были и малоблаговидные действия, процессы разложения, глупости и стратегические ошибки.
 
  Это все скорее <<история математиков>>, чем <<история математики>>, но все же история.
  
  Оказывается, что все перечисленное оказывается за рамками <<социальной истории>>. В оную историю входит главным
  образом преследование математиков властями. Нужно сказать, что таковые случались, но
  вообще математика процветала под гнетом Советской власти. А время, когда <<социальная история>>
  стала писаться, совпало с погромом прикладных наук и удушением наук фундаментальных
 (что заставляет задуматься о преемственности традиций <<истории  естествознания>>).

 Ну ладно, <<социальная история>> в этом смысле -- это тоже история, ее тоже надо изучать.  Интересно не это.
 
\sm

 {\bf\punct  К социальной истории <<социальной истории>>.%
 \label{ss:sochist-sochist}}
 В этом изучении неизбежно должно было возникнуть <<диалектическое противоречие>>. Самые яркие события 
 <<социальной истории>> в упомянутом смысле случились в
   1929-1931гг. (и о них речь шла выше). В них участвовали 3 культовых фигуры советской математики (Люстерник, Шнирельман
   и Гельфонд),
   а потому об этих событиях до конца 80х было принято молчать.
   А в числе вождей тогдашнего переворота стояли основатели Советской школы истории математики
  Яновская и Выгодский (а также Кольман). Как социальная история математики должна была глядеть на свою собственную
  социальную историю?
  
  Яновской и Выгодскому посвящены многочисленные тексты. Выше уже замечалась, что упоминание директорства Выгодского
  в самый неприятный момент 
  <<социальной истории>> из оной истории выпало.
  Потом об этом (а также не очень нежно о самом Выгодском) стали писать
  М.~Ю.~Колягин и О.~А.~Саввина \cite{Kol-savvina-1}, \cite{Kol-savvina-1}, о нем  упоминается в 
  \cite{professora} (ну и так далее, тайное становится явным).
  Первая сколько-либо правдоподобная биография Выгодского
  со стороны историков математики появилась лишь 2015 году \cite{Demidov-Vygodsky} (29 страниц). Но о директорстве
  Выгодского мы узнаем так:
  \begin{quotation}
В сентябре 1930 г. был арестован Егоров  — недавний директор
 Научно-исследовательского института математики и механики Московского университета,
действующий президент Московского математического общества и глава редакции 
основного математического журнала того времени — «Математического сборника». 
Московское математическое сообщество вступило в эпоху
кризиса, все названные московские математические институты оказались под
ударом. Выгодский оказался одним из наиболее видных и активных деятелей
тогдашнего математического сообщества: он входил в состав президиума
Всесоюзной математической ассоциации (1934–1956), руководил Институтом (1931–
1932), являлся вице-президентом Общества (1932–1933), ответственным редактором 
(вместе с С.~А.~Чаплыгиным) 1 и 2 выпусков 38 тома «Математического
Сборника» (1931). В 1933–1935 гг. он работал (по совместительству) в Институте
истории науки и техники АН СССР в Ленинграде, наконец, занимал важные позиции в Государственном технико-теоретическом издательстве — в 1932–
1933 гг. даже был главным его редактором
  \end{quotation}
  
 Это все, что сказано. Причем, о том каким институтом руководил Выгодский, читатель может установить
 лишь путем умозаключения (<<Институт>> с большой буквы), а что  значит  директорство в этот год,
 может понять лишь въедливый
 человек, специально интересовавшийся  историей Московского Физмата.

Автор этих заметок приводил выше
 две жесткие цитаты о Яновской, одна из статьи \cite{Yanov-BDU} И.~Г.~Башмаковой, 
 С.~С.~Демидова и В.~А.~Успенского
 (подчеркну, что последний автор - не историк, а математик), а вторая из статьи Юшкевича \cite{Yush1}, 
 которая никогда не цитируется, и которую можно найти,  лишь  занимаясь преднамеренными розысками. В целом Яновская  превратилась в иконописного персонажа (с отдельными фразами, вкрапляемыми в апологетическую литературу.
 Впрочем Б.~В.~Бирюков \cite{Biryukov-yanovsk} написал ее нормальную биографию, но он и не историк математики. 

   Есть обширная статья С.~С.~Демидова и В.~Д.~Есакова \cite{DeEs}, которая считается основным текстом по <<Делу Лузина>>.
   В ней много и с большим моральным напряжением рассказывается о свержении Егорова и 
   тогдашней математической смуте. Выгодский в статье вообще не упоминается. Яновская   упоминается, но не
   в связи
   с этими событиями.
Зато упоминаются действия неназванных <<красных профессоров>>:
   \begin{quotation}
   Отстранение Д.~Ф.~Егорова было достигнуто в результате длительной осады «пролетарского студенчества» 
   и решительных действий «красных профессоров». 
  \end{quotation}
 
\begin{quotation}
  Весной 1930 года его удаляют с поста директора Института математики и механики Московского университета.
  На его место приходит «красный профессор%
  \footnote{<<Красный профессор>> стоит в кавычках, что подразумевает ссылку на
  язык того времени. Кажется, 
  по тогдашней терминологии <<красный профессор>> -- это термин, обозначающий выпускника Института красной профессуры,
  каковым Шмидт не был. Хотя был профессором и был красным.}» О. Ю.
Шмидт, открывший свою деятельность на новом поприще призывом к сотрудникам Института перестроить работу
на марксистской основе и обвинением во вредительстве тех,
кто попытается этому препятствовать. Взявший после него слово уже бывший директор
Института Д.~Ф.~Егоров заявил, что истинным вредительством является навязывание
всем стандартного мировоззрения.
\end{quotation}

\begin{quotation}
  Так что власть в московском математическом сообществе сама упала в руки «молодежи» — 
  с 1932 года президентом Московского математического общества стал П.~С.~Александров (сохранивший этот пост до 1964 года!)%
\footnote{Смысл восклицательного знака непонятен, Александров, знаменитый математик и умнейший человек,
соответствовал занимаемой должности.}, а главным редактором «Математического сборника» —
«красный профессор» член партии
О.~Ю.~Шмидт.
\end{quotation}

Заметьте, что <<красные профессора>> упоминаются во множественном числе. А кто еще?
Читатель не слишком внимательный поймет прочитанное так, что свергал 
Егорова Шмидт. Но никаких конкретных доводов не приводится. Источник о полемике Егорова со Шмидтом не приводится,
но об этом разговоре см. выше п.\ref{ss:1930-spring}, и поведение Шмидта было вовсе  не однозначным.

 \sm
 
 Понятно, что выгораживание своих учителей и отцов-ос\-но\-ва\-те\-лей по-че\-ло\-ве\-чес\-ки понятно, но в данном случае 
 предмет <<социальной истории>>
 как науки улетучивается и начинается (как мы сейчас увидим) диффамация людей, которым
 не выпала честь служить историками математики. 
 
 \sm
 
{\bf\punct Юшкевич Адольф Павлович%
\label{ss:yushkevich-1}} (1906-1993). Заведующими кафедры/кабинета истории математики
на Мехмате последовательно были исследователи <<Математических рукописей>> Маркса Яновская и Рыбников%
\footnote{Согласно \cite{Gelfond-VIET}, на Мехмате некоторое время существовала 
Кафедра теории чисел и истории математики, заведующий Гельфонд.}.
 Но научным лидером советской истории математики
с конца с 1940х годов считался Юшкевич.
Он учился на Физмате МГУ в 1923-1929гг.
В 1936 году без защиты получил степень кандидата физ.-мат. наук,
в 1940 защитил докторскую <<Математика и ее преподавание в России в XVIII веке>>.
В 1930-1952гг.
 работал в Московском высшем техническом училище (<<Бауманка>>) -- сначала ассистентом, с 1934 г. --
доцентом, с 1940 г.-- профессором, в 1941-1952гг. зав. кафедрой высшей мате­матики. 
С 1945 работал в  Институте
истории естествознания и техники АН СССР.

Он автор многочисленных работ по истории математики,
также был организатором и соорганизатором разных важных и успешных издательских проектов
(например, сериал <<Историко-математические исследования>> вместе с Г.Ф.Рыбкиным,
один из заместителей ответственного редактора в цунном многотомном издании <<История отечественной математики>> \cite{Stokalo-4-2}, изданном в Киеве).
Он -- один из интересных персонажей уходящей в прошлое блестящей математической эпохи. 
Но здесь нас интересуют не  научные заслуги Юшкевича (заслуги бесспорные)
 в 
исследованиях математики разных стран от глубокой древности до XIX века, а 
 его биография в 1920-30г, когда эти достижения были еще впереди.

\sm

Но перед этим приведем еще одна цитату из \cite{DemidovTokareva}:
\begin{quotation}
Он [Юшкевич] требовал от каждого из своих учеников
хорошего знания литературы по изучаемой теме и серьезной ее проработки. В связи с этим вставал вопрос 
о знании иностранных  языков, чему Юшкевич, владевший основными живыми европейскими языками и латынью,
придавал особое значение... разумеется во главу угла ставилось хорошее владение математической 
стороной дела - задача должна была быть  рассмотренной c позиций современной науки.
Это было для него, так сказать, и правилом хорошего тона, и необходимым жизненным требованием --
диссертации, по правилам ВАК'а шли по разряду <<математика>>, а Адольф Павлович не 
хотел уронить своего достоинства в глазах коллег-математиков. Он был воспитанником
знаменитой математической школы и рассматривал, прежде всего, ее представителем.
\end{quotation}

\sm

{\bf\punct Юшкевич: эпизоды ранней биографии.%
\label{ss:yushkevich-2}}
Вернемся к нашей теме.
Его отец, Юшкевич Павел Соломонович (1873- 1945), был философом и меньшевиком.
Согласно воспоминаниям Юшкевича-мл. \cite{Yush-school}:
\begin{quotation}
 Его интересовала математика и в не меньшей
степени общественная жизнь. Кончая гимназию, он вступил, около
1890 г., в революционный кружок, руководимый студентом
Мих[аилом] Владимировичем] Морозовым [1868-1938]. Вскоре
их всех арестовали и, после полутора тюремных лет, отца 
отправили в ссылку в Кишинев. Там было немало подобных ему 
ссыльных, среди них Дав[ид] Бор[исович] Гольдендах (псевдоним 
Рязанов), будущий большевик и директор Института Маркса и 
Энгельса, а также академик....
\end{quotation}

Учился математике в Париже (в Сорбонне), видимо, в 1891-1892 и 1898-1803гг. 
Продолжаем цитирование:
\begin{quotation}
В революции 1905—1907 гг. отец принял активное участие.
Часто выступая на рабочих и матросских митингах, как 
представитель РСДРП [Российской социал-демократической рабочей 
партии], в которую вступил в начале века, а в 1903 г., после ее 
раскола [на II съезде], примкнул к фракции меньшевиков....
\end{quotation}

Книга Ленина <<Материализм и эмпириокритицизм>> (1909) была реакцией на
сборник
статей социал-демократов В.~Базарова, А.~Богданова, А.~В.~Луначарского, Я.~А.~Бермана,
О.~И.~Гельфонда,
 П.~С.~Юшкевича и С.~А.~Суворова под названием
«Очерки по философии марксизма»(1908). Само название книги Ленина является 
перефразировкой названия книги 
П.~Юшкевича  «Материализм и критический реализм>>. Юшкевич, в свою очередь, ответил  Ленину брошюрой
<<Творцы философской ортодоксии>> (1910)...

Согласно примечаниям к 18 тому Полного собрания сочинений Ленина \cite{Lenin}, издание 1968г.,
\begin{quotation}
 В 1917-1919 годах сотрудничал на Украине [видимо, в Одессе, но почему-то этого не сказано]
 в меньшевистском журнале «Объединение>>
и других антибольшевистских изданиях;
\end{quotation}

Из воспоминаний Юшкевича-мл.:
\begin{quotation}
Отношения между большевиками и меньшевиками обострялись,
хотя и постепенно — многие местные большевики были 
вчерашними меньшевиками. Отцу поручались курсы истории философии в
местной партийной школе. Он был в добрых отношениях с 
редактором официальной газеты И.~И.~Яновским...

Как-то в Одессу приехал Луначарский и, встретив отца, стал
уговаривать вернуться на север, где интеллектуальная жизнь 
оживает.  Отец с большими приключениями добрался в 1921 г. до 
Москвы, повидал Рязанова и [Л.Б.] Каменева
знал, и от последнего получил такой же совет.
\end{quotation}

В дальнейшем П.~С.~Юшкевич отошел от политической деятельности.
Он переводил на русский сочинения  Маркса--Энгельса, в частности он перевел <<Диалектику природы>>.
Позже переводил
 различных философов (Гельвеций, Гольбах, Дидро) и историко-математические книги.

А в только что цитированном отрывке небезынтересен набор упоминаемых лиц:

\sm

{\small  Яновский  Исаак Ильич -- это муж С.~А.~Яновской.

 Луначарский Анатолий Васильевич -- Нарком просвещения
 (октябрь 1917 - сентябрь 1929), впоследствии - культовая фигура советской истории
 (и, кстати, один из авторов сборника, который
 разносил Ленин).
 
 Каменев (Розенфельд) Лев Борисович  -- в 1924-25гг. один из лидеров Советского государства
 (вместе со Сталиным и Г.~Е.~Зиновьевым). 
 
Рязанов Давид Борисович (Давид-Симха Зельман-Берович Гольдендах). В числе меньшевиков-межрайонцев
(группа Троцкого)
вступил в 1917г. в РКП(б). Основатель
Института Маркса-Энгельса и его директор 1921-1931. Когда в начале 1929г. проводилась
операция по советизации АН СССР, Политбюро предполагало сделать Рязанова
вице-президентом (всего вице-президента было два). Рязанов  отказался по состоянию
здоровья. Тогда был выставлен Кржижановский.}

\sm

Среди авторов разнесенного Лениным сборника был и
Гельфонд Осип Исаакович. Согласно все тем же примечаниям в \cite{Lenin},
\begin{quotation}
	\small
Гельфонд~О.~И. (1863-1942) - один из авторов ревизионистского
сборника <<Очерки по философии марксизма>> (1908); по профессии
врач; с конца 80-х годов принимал участие в революционном движении,
в 1905 года содействовал распространению социал-демократической
литературы в Киеве. После Октябрьской социалистической революции
работал врачом (в 1922-1928 годах - в Институте К.~Маркса
и Ф.~Энгельса и Комакадемии) . Написал ряд работ по вопросам медицины
и несколько статей по философии: <<Философия Диогена и современный
позитивизм>> ( 1908 ), <<Об эмпириокритической теории познания>>
( 1908) и другие. Философские взгляды Гельфонда В.~И.~Ленин
назвал «окрошкой из материализма и агностицизма>>.
\end{quotation}
 
 А еще О. И. Гельфонд был отцом знаменитого математика, одного из героев настоящих записок.

 Продолжим цитировать Юшкевича-мл.:
\begin{quotation}
 Мы провели несколько недель в семье математика
Я.~С.~Дубнова [1887-1957], приходившегося нам далеким 
свойственником. Член Президиума ВЦИК П.~Г.~Смидович [1874-1935]
распорядился дать нам большую комнату (бывший актовый зал 
какой-то школы) на Старой Басманной, д. 19.... Отец занялся переводами, которые 
первоначально ему поручал знакомый по Кишиневу Д.~Б.~Рязанов, сестре
предложил занятия на медицинском. Я регулярно печатал на 
машинке переводы отца под его диктовку (один час в день). Стал 
обзаводиться знакомствами - прежде всего с детьми В.~А.~Базарова
(Руднева) [1874—1939], с сыном И.~И.~Гельфонда
 и младшей 
дочерью М.~В.~Морозова.
\end{quotation}

Отметим, что Базаров и Морозов упоминались чуть выше. Позже А.~О.~Гель\-фонд-мл.
	женился на дочери  Морозова.

Уже о более позднем времени:
\begin{quotation}
 В 1930-е гг. Гельфонд\footnote{Из введения С.~С.~Демидова и Т.~А.~Токаревой
 к воспоминаниям Юшкевича: {\it Особенное место в публикуемых воспоминаниях Адольфа Павловича
занимает личность его друга - выдающегося советского математика 
Александра Осиповича Гельфонда}.
\newline Это действительно так.} особенно подружился с Л.~Г.~Шнирельманом.
Это был чрезвычайно умный математик и интересный человек.
Нередко мы собирались тесной компанией в Доме ученых — оба
друга с женами и Шнирельман (редко кто-либо еще) — и он 
рассказывал нам фантастические, придуманные им истории с элементом
мистики и ужасов.
\end{quotation}

Вернемся к профессиональной биографии Адольфа Павловича.
По-ви\-ди\-мо\-му, от своего отца он унаследовал и философские, и историко-математические интересы. В 1923
году поступил на московский Физмат:
\begin{quotation}
 Собственно, я хотел на 
философский, но отец убедил меня в том, что философ неграмотный в 
математике и физике — невежда, и я решил поступить на математическое отделение...

В 1925-1926 гг. меня заинтересовали споры по вопросам 
обоснования математики и, прежде всего, анализа [отголоски дискуссий первых трех десятилетий XX
века приводились выше в \S \ref{s:filosofiya}].
\end{quotation}

и дальше

\begin{quotation}
 В 1925-1926 гг. мы с Гельфондом искали заработка...
 В эти же годы мы стали, как я сказал,
часто посещать Секцию естественных наук Комакадемии, и здесь
нам неожиданно повезло в том смысле, что благодаря Л.~М.~Лихтенбауму
и Л.~А.~Люстернику, мы получили первую постоянную 
работу: нас приняли младшими научными сотрудниками Отдела 
философии и истории естествознания Научно-исследовательского 
института имени К.~А.~Тимирязева [Институт входил в систему Комакадемии \cite{Dubinin}]. Институт этот, помещавшийся на
Пятницкой ул., д.48, был биологический... Математику там сперва
представляли Лихтенбаум и Люстерник и когда, по предложению
О.~Ю.~Шмидта, они перешли в Комакадемию, в Отделе 
образовались вакансии, на которые нас и зачислили по их рекомендации
(на полставки каждого)...
\end{quotation}

В 1926-1929г работал в 
Отделе конъюнктуры Центрального статистического управления
СССР, где должен был заниматься вопросами текущей 
экономической конъюнктуры. Что касается занимавшимся им должностей, есть случайное упоминание того, 
что он в 1938г.
работал старшим редактором издательства ОНТИ%
\footnote{В.~С.~Кирсанов. ВИЕТ. 2005. № 4. С. 105-124.}.

Продолжим цитирование:
\begin{quotation}
В конце 1920-х гг. я стал часто встречаться с Яновской, и у нас 
установились хорошие отношения. Когда в 1933 г, начал работать 
основанный ею и М.~Я.~Выгодским  университетский научный
семинар по истории математики%
\footnote{Колмогоров--Александрову, 19 октября 1932 г.
\newline
{\it 
Часы от 7 до 9 в четвертый день шестидневки, о которых ты
спрашиваешь, уже заняты философским семинаром, коего ты
состоишь руководителем (вкупе с С.~А.~Яновской, О.~Ю.~Шмидтом,
А.~Я.~Хинчиным, Л.~С.~Понтрягиным и мною). Это не означает, 
впрочем, что ты должен там часто бывать. Открывается семинар
моим докладом об основаниях теории вероятностей, но далее идет
обширная программа исторических докладов (Варьяш, Выгодский,
Яновская ...), в которой я не принимаю активного участия.}
}, я стал одним из его постоянных
участников. (Часто встречающееся утверждение, что семинар этот
возник в 1935 г. - ошибочно.) С.~А.~Яновская руководила или 
принимала участие в руководстве этим семинаром до конца жизни.
Особенностью ее частых выступлений, как и печатных работ, было
то, что к вопросам истории математики она подходила с позиций
марксистской методологии и вместе с тем рассматривала их, чем
дальше, тем чаще, под логико-математическим углом зрения.
\end{quotation}

Б. А. Розенфельд \cite{Rozenfeld} писал: 
\begin{quotation}
 Адольф Павлович хотел поступить в аспирантуру к С.~А.~Яновской, с которой был знаком еще в Одессе, 
 но она не решилась принять
его.  А.~П.~Юшкевич занимался историей математики самостоятельно
и вскоре стал одним из крупнейших историков математики.
\end{quotation}

О каких-либо общественных действиях Юшкевича в 20-30е годы автору не известно. От него остались публикации,
которые было бы небезынтересно разыскать%
\footnote{Из тех же воспоминаний: {\it Это была моя первая публикация по истории математики}
[Статья про Лазара Карно \cite{Yushkevich-Karno}. {\it Я зря 
пытался в ней усмотреть элементы гегелевской диалектики.}}. Например, остались предисловия редактора (Яновской) и
переводчика (Юшкевич) к книге Германа Вейля [Hermann Weyl] <<Философия математики>>, 1934 (это перевод нескольких
статей).
На неподготовленного читателя (автор помнит себя в качестве такого  в 1980г.)
предисловия эти могут произвести сильное впечатление. На фоне многого того, что автору пришлось прочесть,
работая над данными записками, они скорее выглядят даже  беззубо (все ж уже был 34ый год, а не 31ый).
Переводчик заключает свои речи такой сентенцией:
\begin{quotation}
Из настоящей работы читатель увидит все же, что интуиционизм
ставил ряд важнейших вопросов в своей критике формально-логического
направления в.математике и теории континуума. В этом нет, пожалуй,
ничего удивительного. <<Когда один идеалист ругает другого, на этом
выигрывает материализм>> (Ленин). И значение работ Вейля именно в этой
их критической стороне. 
\end{quotation}

При чтении этого текста видно серьезное отношение переводчика к философии, а также к основаниям
математики (за исследованиями в этой области переводчик, как видно из его текста, внимательно следил). 
Из предисловия ясно, что было какое-то общее решение об издании серии сборников 
по буржуазной философии, с такой целью:
\begin{quotation}
 Между тем 
выработка марксистских воззрений немыслима без знания и понимания этих —
в основном нематериалистических — теорий%
\footnote{Марксизм тогда еще не боялся дискуссий, это еще не ритуальный марксизм 
	1970х годов.}.
\end{quotation}

Ну и видно, что инициатором издания именно работ Вейля, и именно этих работ был сам Юшкевич:
\begin{quotation}
 При выборе материала для этого сборника я счел полезным остановиться на статьях Вейля,
 а не главы интуиционистской школы Броуера [Luitzen  Brouwer],
потому что работы последнего доступны лишь очень ограниченному
кругу читателей.
\end{quotation}

В списке работ Юшкевича 30х годов присутствуют предисловия к многим книгам, всего 10
(<<предисловия>>, <<предисловия редактора>>, <<предисловия к переводу>>, <<вступительная статья>>),
что показывает, что он был тогда не последним человеком.

Что касается его докторской диссертации, то С.~С.~Демидов и Т.~А.~Токарева говорят следующее:
\begin{quotation}
 ...он решил заняться подготовкой докторской, предполагая выбрать в качестве предмета изучения
 историю оснований анализа, но С.~А.~Яновская, учитывая идеологическую ситуацию, порекомендовала ему
 выбрать тему, связанную с историей математики в России. Защита состоялась 28 мая 1940 г. 
 на механико-математическом факультете МГУ; оппонентами выступали А.~О.~Гельфонд, М.~Я.~Выгодский,
 С.~А.~Яновская.
\end{quotation}

\sm

{\bf\punct Статьи Юшкевича 1989 и 1991 года.%
\label{ss:yushkevich-89-91}}
Во время Перестройки начали  публиковаться письма \cite{Kap} знаменитого физика
Капицы высшим государственным деятелям в защиту различных ученых,
 в частности, стало известно его письмо Молотову от 6 июля 1936г. в защиту Лузина (оно приведено ниже).
Юшкевич в связи с эти упоминает статью  в газете «Советская культура» от 21.05.1988.

В апрельском номере Вестника Академии наук за 1989г. выходит статья Юшкевича <<Дело Лузина>>
(отметим, что в тот момент Стенограмма уже была обретена и цитируется в статье). В 1991
году в Сборнике <<Репрессированная наука>> появляется статья того же Юшкевича с тем же названием (статьи имеют пересечения,
но тексты их различаются, что не удивительно, а интересно, что есть 
существенные  различия в их смыслах). В статьях впервые открыто пишется о событиях
на московском Физмате 1929-31гг., пишется  кратко и четко (особенно в первой), там также рассказывается о <<деле Лузина>> 1936г.
Юшкевич формулирует идею, что дело Лузина было хорошо продуманной атакой на науку со стороны Сталина,
и (во второй статье) что непосредственно курировал эту атаку Кольман. Мы обсудим это чуть ниже, а пока должны
отметить три обстоятельства.

\sm

Первое.
Юшкевич родился в 1906 году, во время потрясений московского Физмата ему было 23-25 лет, а во время дела Лузина
- 30. Он был очевидцем обеих историй (а также логически не исключено его участие). Он был членом
тогдашней лево-радикальной научной и около-научной тусовки, а также просто лево-радикальной тусовки. Более того, он входил в нее каким-то
семейным образом, в этой тусовке были дорогие и близкие ему люди. 

\sm

Второе. После Войны Юшкевич был довольно высокопоставленным историком и имел доступ ко многим
архивам. Скорее всего (особенно учитывая предыдущее, а также его знакомства с ведущими
математиками), он знал, что произошло на самом деле (наличие у него интереса к проблеме
видно, например из пп.\ref{ss:raspad-yushkevich}, \ref{ss:last}).

\sm

Третье. Юшкевич неоднократно видел, как история творится на глазах: 
не в том смысле, что при нем происходили исторические события, а 
в том смысле, что радикально менялось прошлое. Это наблюдалось, по крайней мере, в 1938-39гг. и в 1953-57гг. Он имел
возможность видеть,
что этот творчество имело какие угодно цели, кроме установления истины. В 1987-1991гг. он наблюдал 
этот процесс снова.

\sm

В связи со всем этим, у нас не должно быть априорной уверенности, в том,
что Юшкевич был беспристрастным летописцем и старался донести до нас именно истину.
Поэтому мы должны проанализировать его  статьи так, как положено анализировать
исторические источники.

\sm

Юшкевич в 1991г. сообщает о Стенограмме следующее:
\begin{quotation}
{\bf Никто из участников заседаний не предъявил Лузину политических обвинений},
все, близко знавшие Лузина, говорили, что он вообще никогда
не высказывался по политическим вопросам. 
\end{quotation}

Читатель может оценить степень истинности этого высказывания...

А так рассказывается о Хинчине:
\begin{quotation}
Заслуживает внимания выступление А.~Я.~Хинчина, возражавшего против ряда обвинений,
содержавшихся в правдинской статье «О врагах в советской маске».
«Суслина, — сказал Хинчин, — называют учеником Н.~Н.~Лузина,
загубленным Н.~Н. Ну, когда человек умирает от сыпного тифа, то это слишком
резкое выражение. Ведь он мог заболеть сыпным тифом и в Иванове. Но общее мнение таково, 
что из Иванова Н.~Н. выжил Суслина. Однако самый перевод из Москвы
в Иваново я считаю услугой, оказанной Суслину Н. Н., тогда
еще не враждебно настроенному по отношению к Лузину».
Далее А.~Я.~Хинчин возражал против того пункта той же правдинской
статьи, в котором Лузин объявлялся представителем «бесславной царской
московской математической школы, философией которой было черносотенство
и движущей силой — киты российской реакции». Дело в том, что, когда
Н.~Н.~Лузин был молод и складывались его академические и, вероятно, политические убеждения,
«...в это время в Московском университете происходила борьба между реакционной, действительно черносотенной группой, и другой группой,
которую возглавлял Д.~Ф.~Егоров и которая стремилась к европеизации
в буржуазном смысле этого слова Московского университета. И Лузин целиком принадлежал 
к этой второй группе. Он отнюдь не был связан с самодержавным
правлением».
\end{quotation}

Как мы видим, сказанное о Хинчине,
это правда, только правда, тщательно отобранная правда, но это не вся правда.

Подробно (несколько смягчено) обсуждаются научно-этические обвинения. 
 Выходит, что математики кое-как
приняли участие в атаке на Лузина сверху и предъявляли не вполне обоснованные обвинения.
Были, правда какие-то <<некоторые>>:
\begin{quotation}
На протяжении всех этих заседаний вырабатывалась и резолюция комиссии.
Некоторые участники заседаний, как Б.~И.~Сегал и другие, предлагали
сохранить в ней крайне резкие формулировки, содержавшиеся в статьях
«Правды».
\end{quotation}

Б. И. Сегала  к 1991 году мало кто не помнил (и, быть может, не все, кто помнил,
 поминали его добрым словом).

Перед обсуждением Стенограммы в статье 1991г. Юшкевич пишет:
\begin{quotation}
 К сожалению, эти стенограммы были обнаружены, когда данная статья была готова, и они учтены только
в той части, в которой уточняют отдельные положения статьи. 
\end{quotation}

Однако Стенограмма уже упоминается в статье за апрель 1989 года и там обсуждается. Кроме того,
сравнив реальную стенограмму и рассказ Юшкевича, мы можем понять, что работа по ретуши была проведена
большая и тщательная.

По-человечески, автора можно понять, он скрывает малосимпатичное поведение целого ряда уважаемых и
заслуженных людей, в том числе его
собственных давно умерших
друзей (и, может быть, всем нам было бы лучше этого не знать).
С другой точки зрения он действовал в рамках вечной идеологии российской интеллигенции:
интеллигенты прекрасны, а все, что когда-либо было и есть плохо, исходит от властей.
Точка зрения вполне удобная, никто ей возражать не стал, и, в частности не стали возражать ученики
лузинских учеников, которые знали больше, чем было в статьях о Лузине
 написано. А зря не стали...
При следующем обороте истории эти лузинские ученики превратятся в полных уродов...

\sm

{\bf Интересно сравнить и две версии статьи Юшкевича}. В первой, например, сообщается, что
<<на собрании [в МГУ 9 июля] с критикой Лузина выступили среди других Александров, Колмогоров, Люстерник, Понтрягин>
(а также кратко говорилось о выступлениях).
Это исчезло.  В первой статье была такая сентенция:
\begin{quotation}
 Как видим, молодое поколение математиков энергично выдвигалось на руководящие позиции. 
 Процесс этот, сам по себе естественный, в тогдашней обстановке приобрел весьма зловещий характер. 
 Часть выдвинувшихся на рубеже 20-х и 30-х годов молодых ученых считала, что старшее поколение —
 Д.~Ф.~Егоров и Н.~Н.~Лузин в особенности — отстают от новых направлений науки, 
 не должны далее возглавлять математическую жизнь ни в Москве, ни в Советском Союзе вообще%
  \footnote{Кстати, это свидетельство очевидца, подтверждающее нашу точку зрения на выступление
  <<Инициативной группы>>}. 
 Тут играли роль и чисто личные моменты,
 например, всем известные напряженные отношения между Н.~Н.~Лузиным и П.~С.~Александровым. 
\end{quotation}
Сентенция
была преобразована (увеличена), но упоминание об Александрове исчезло. 

В первой статье говорилось:
\begin{quotation}
 В шельмовании представителей физико-математических наук, о котором писала недавно «Советская культура»,
 с 1930 г. весьма активно участвовали Э.~Кольман, тогда сотрудник Коммунистической академии, и 
 С.~А.~Яновская, руководитель методологического семинара на физико-математическом факультете МГУ. 
 27 апреля 1931 г. в Ленинграде и 5 июня в Москве на I Всероссийской конференции
 по планированию математики Кольман сделал доклад «Современный кризис математики и основные линии ее реконструкции».
 По этому докладу 9 июня конференция приняла резолюцию «О кризисе буржуазной математики и о реконструкции математики в СССР». 
\end{quotation}
Во второй статье Яновская из этого контекста исчезла.

Мы видим, что переход от первой статьи ко второй (кстати, более обширной) не был движением
в сторону развития истины. Наоборот усиливалась те свойства изложения, о которых было написано выше.
Посмотрим однако, что нового (кроме более подробного обсуждения стенограммы) появилось в статье 1991г.

\sm

{\bf \punct О возрастании роли Кольмана.%
\label{ss:kolman-growth}}
В первой статье Кольман упоминается в трех местах: в цитированном выше абзаце, в биографической справке,
\begin{quotation}
	Кольман~Э.— уроженец Чехии, воспитанник Пражского университета. Попал в русский плен во время первой мировой войны. Вступил в ВКП(б), долгое время придерживался официальной линии партии. Закончил свою жизнь в политической эмиграции в Швеции. 
\end{quotation}

Третье место - сноска:
\begin{quotation}
  Заключение статьи очень близко к характеристике Н.~Н.~Лузина на с. 290 книги 
  Э.~Кольмана «Предмет и метод современной математики» (М.: Соцэкгиз, 1936), 
  подписанной к печати 14 июля того же года. Участие Кольмана в кампании против Лузина несомненно, 
  хотя под конец жизни он заверял всех, кто к нему обращался, что вообще ничего не помнит о «деле Лузина». 
  Однако {\bf  с очень большой долей вероятности можно утверждать,
  что  Кольман принимал активное участие} в подготовке цитированной выше статьи «Правды». 
\end{quotation}

Это весьма интересное место, и мы его детальным образом обсудим ниже.
Во второй статье сноска остается, но появляется такой абзац:
 \begin{quotation}
   Неизвестно, по чьей инициативе оно [дело Лузина] началось. Но {\bf можно с уверенностью сказать,
   что самое активное участие в организации «дела» принял Э.~Кольман}, в 1935—1938 гг. занимавший должность
заведующего отделом науки Московского комитета КПСС, когда его первым
секретарем был Н.~С.~Хрущев. {\bf Активное закулисное участие Кольмана не вызывает никаких сомнений:
без его санкции «дела» вообще
не было бы и контролировал ход «дела» с самого начала и даже до его начала, привлекая при том
некоторых математиков, конечно он. Все статьи, шельмовавшие Н.~Н.~Лузина
(и содержавшие часто ложную информацию),  все собрания, созываем для
осуждения Лузина, проходили под наблюдением Кольмана.} Позднее, в 60-е гг.,
он постарался забыть свое участие в «деле» Н.~Н.~Лузина. Я поверхностно
был знаком с Кольманом еще до Отечественной войны, а с 1953 г. он несколько
лет работал в Институте истории естествознания и техники АН СССР, и здесь
мы встречались регулярно. Знакомство наше продолжалось и после его ухода
из Института, и в разговорах со мной он был, вообще говоря, довольно откровенным. 
Но на мои неоднократные вопросы о том, что ему известно о «деле»
Н.~Н.~Лузина, он неизменно отвечал, что когда-то слышал о нем, но толком ничего не помнит. 
Такой же ответ он прислал уже будучи в эмиграции на вопрос
моего французского коллеги профессора П.~Дюгака. Память у Кольмана и на
старости лет была очень крепкая. Скорее всего, ему было неприятно вспоминать
о своем бессовестном поведении при организации антилузинской кампании.
Впрочем, стыдно было не ему одному, но и многим других лицам, принявшим
в этой кампании деятельное участие, как свидетельствуют их статьи о Н.~Н.~Лузине,
написанные по тому или иному поводу после его кончины.
 \end{quotation}
 Рассказывается о погромных статьях Кольмана%
 \footnote{{\it Так, в журнале «Большевик» (ныне «Коммунист») за 1931 г., № 2, он поместил статью «Вредительство в науке», где на нескольких страницах (73—81) «разоблачил» примеры
 	«вредительства» в самых разных сферах знания — от психологии до политической экономии и финансового дела, в частности в ряде методологически ошибочных, по его мнению, статей издания Большой советской энциклопедии, та-
 	ких как «Величина», «Вероятность», «Волны», «Гидромеханика» и даже в кратких
 	 биографиях Галилея, Гарвея и Гаусса. Тут досталось и главному редактору
 	О. Ю. Шмидту и руководителю математического отдела известному геометру
 	профессору В. Ф. Кагану.}.},
 о его сближении с Ленинградским обществом математиков-материалистов, про его доклад в Цюрихе. Все это так, и все это 
 проверяемо. Делается предположение, что именно он донес властям о поведении Бернштейна,
 Егорова и Гюнтера на Первом всесоюзном математическом съезде (может и сообщал,
 но вообще-то скандал был, на поверхности и тайной не был).
 
 \begin{quotation}
 	... не без его деятельного
 	участия ...
 	появилась «Декларация инициативной группы о реорганизации Математического общества», подписанная такими
 	выдающимися молодыми математиками, как будущие члены-корреспонденты
 	АН Л. А. Люстерник, Л. Г. Шнирельман, А. О. Гельфонд и будущий акад.
 	Л. С. Понтрягин....
 	\end{quotation}

 Мы анализировали
 выше этот документ, очевидно, что этот текст не мог быть продуктом творчества Кольмана,  Яновской, или Хотимский (хотя участие кого-либо из них  возможно%
 \footnote{А ленинградская декларация как раз была их сочинением, и по ней это видно.}). Кстати, свое утверждение
  Юшкевич не аргументирует. В качестве кого он выступает,
 очевидца или исследователя истории%
 \footnote{Кстати, тот же вопрос встает и про упомянутые выше <<все собрания>>...}? Его рассказ о декларации кончается так
 \begin{quotation}
 	 	Такую декларацию мог бы написать и Кольман, консультативное участие которого в ее составлении весьма вероятно.  
 \end{quotation}
 Интересно, что две сентенции не вполне согласованы. 
 
 И главное, дело Лузина, в котором Кольман из соучастника написания статьи в  Правде  в подстрочной сноске
 в первой публикации
 вырос в куратора этого дела во второй. Что мешало Юшкевичу сообщить об этом (и о всем прочем!) двумя годами ранее?
 
 \sm
 
 Тут автор должен напомнить про <<эталонных демонов>>, которых немало в изложениях российской
 (и особенно советской) истории.
 На этих демонов принято навешивать грехи всех и вся и делать их них источники всех бед.
 Те, кто сомневаются, подвергаются общественному остракизму: как? ты оправдываешь такого-то?
 В итоге история превращается в демонологию и в источник жестоких фантомных болей. О таких
 фантомных болях
 речь пойдет далее в этом параграфе.
 
\sm

   {\bf\punct Культ личности Кольмана.%
   \label{ss:kolman-kult1}}
Отметим, что {\it Юшкевич говорил о том, что Кольман участвовал в сочинении статьи в Правде, но не говорил,
что он был автором.}
Автором он был объявлен в статьях Левина \cite{Levin-1990} 1990
(на основе общефилософских соображений и цитированной выше сноски 
из статьи Юшкевича) и С.~С.~Демидова и Т.~А.~Есакова \cite{DeEs}.
 В дальнейшем  роль Кольмана в истории науки начала расти до небес, о нем 
 гуляют многочисленные высказывания типа  <<Виднейший идеолог сталинской эпохи>>,
 правда скорей по вторичным текстам, чем научным работам.

 Впрочем,  Н.~С.~Ермолаева \cite{Erm2} писала  <<{\it В Москве главным математическим идеологом становится Э.Я.Кольман}>>
(в отношении 1931г., и это не вызывает сомнений),
она же в другой статье
  \cite{Erm2} <<{\it Как известно Кольман был главным идеологом советской науки}>> (о 1936г.).
  
  В статье
С.~С.~Демидова и Т.~А.~Есакова он характеризуется как 
<<{\it одно из
влиятельнейших в сфере науки того времени лиц}>>, или косвенно
во фразе
\begin{quotation}
 Вождь, взявший
в это время Академию под свой непосредственный контроль и пытавшийся организовать
её работу в соответствии со своим видением того, как должен действовать «штаб советской науки»,
меньше всего был заинтересован в отдании Академии во власть партийных
функционеров среднего звена — мехлисов и кольманов. 
\end{quotation}
То есть Кольман изящно равняется с Мехлисом, а примерное иерархическое положение Мехлиса
российскому читателю было известно. Ну, хотя бы он был кандидатом в члены ЦК ВКП(б), то есть входил в
140
высших лиц СССР. При этом российский читатель скорее вспомнит о Мехлисе в более высоких ролях,
чем просто роль ЦКиста.

В статье С.~С.~Демидова \cite{Demidov} 1999г. Кольман назван <<{\it черным ангелом московской математики}>>.

Автор настоящих заметок, не испытывая к Кольману ни малейшей симпатии, позволит себе
высказать сомнения в величии сыгранной им роли.

\sm

{\bf\punct Снова о биографии Кольмана.%
\label{ss:kolman-biography}}
Мы начали обсуждать  биографию этого деятеля в \ref{ss:kolman-memoirs}
(и в частности усомнились в достоверности его воспоминаний).
 Во второй половине 20х он занимает разные, в общем
не очень высокие должности
(типа заведующего издательством, издательств тогда, как показывают справочники
<<Вся Москва>>, было много), не связанные с наукой.
Дальше  в 1929-1930 работал в отделе Агитпропа ЦК (возможно, меньшее время,
чем он говорит). В апреле 1929г. Кольман выступал на всесоюзной конференции марксистско-ленинских научных 
учреждений.

В январе-феврале 1931 года после изгнания Шмидта он оказывается руководителем Ассоциации
Институтов естествознания Комакадемии%
\footnote{В доступных автору публикациях
	\cite{Zasedanie}, \cite{Za-povorot} о разборке в Комакадемии
	Кольман присутствует, но какой-то особой его роли не видно. Яновская, напротив, выделяется.
В статьях 1931г. Кольман крыл Шмидта и Кагана (кажется, Кагана больше), следуя 
резолюции Комакадемии.},
это весьма высокая в отношении науки должность, причем именно в тот момент времени. 
Посмотрим на его достижения в этой должности.

Он становится тогда главредом  журнала
<<Естествознание и марксизм>>. Журнал
(как утверждается в \cite{Bogolyubov-Rozhenko}) был основан Шмидтом в 1929 году, с приходом Кольмана 
издание начинает называться
«За марксистско-ленинское естествознание»  (название достойно главреда). 
Судя по библиотечному каталогу МГУ, в конце 1932 года журнал завершает свое существование
(судя по разнообразным ссылкам, журнал кольмановского периода был фееричен, но автор не удосужился его посмотреть).

Был также  журнал <<Научное слово>>, где Шмидт был главным редактором. 
Согласно \cite{Bogolyubov-Rozhenko},
он переходит к Кольману, а (согласно библиотечным каталогам) через два выпуска 
он оказывается присоединенным к журналу СОРЕНА.

В 1931г. Кольман начинает издавать журнал <<Математическая наука -- пролетарским кадрам>>.
Вышел лишь один его выпуск (см. \cite{Tok-white}).

 В  1932 году прикрывается руководимая Кольманом <<Ассоциации естествознания Комакадемии>>. 
 
 Наш герой становится
 директором Института красной профессуры естествознания. Единственный след
 существования этого учреждения, который автору удалось найти, - информация о том, что Институт
 был распущен до июня 1933г (см. выше п.\ref{ss:kolman-memoirs}).
 
 Он, видимо, руководил  операцией по формированию
 марксистских ес\-тест\-вен\-но-на\-уч\-ных обществ.
 Операция проваливается, эти уродцы оказываются
 никому не нужными%
 \footnote{Вот образец продукции Ленинградского общества математиков-материалистов,
 	\cite{Tsigler}:
\newline 
{\it 
\dots основным условием в деле перестройки методики }[преподавания математики]{\it должно быть следующее: руководящий орган народного просвещения определяет лишь направляющую программу знаний, а социально производственное окружение (завод, шахта, колхоз) заполняет программу реальным содержанием.
\newline
Избираемые математические работы должны способствовать выработке целостного материалистического мировоззрения, развитию у учащихся творческих способностей, воспитанию в них четкости, исполнительности и аккуратности, а также необходимых организационно-плановых навыков
\dots
\newline
Прежде всего это }[учебник по математике]{\it будет не
универсальная <<хрестоматия>>, начиненная отвлеченными данными, а серия небольших математических брошюр (листовок), составленных по принципу отдельных социально-экономических, политических и производственных тем.
\newline
Примером такой математической разработки может служить тема: <<Прямоугольник, поверхность фигур и гипербола на изучении свойств индикатора>>. Социальный заказ  -- удовлетворить потребность нашего хозяйства в квалифицированных лаборантах теплотехнического производства.\dots 
\newline
Можно назвать и другие темы как-то:
<<Счет, измерение и построение в посевной кампании>>,
<<Число, отношения, и построение диаграммы на данных Октябрьской революции>>, <<Окружность, круг и цилиндр в жестяном производстве>>, <<Число, отношения и формулы в вопросах нормирования и рационализации труда>> и ряд других.
\newline
Совокупность таких тем, оформленных в виде брошюр в 2-4 печатных листа каждая, даст возможность учителю построить на их разработке в течение года
всю программу по математике для данной группы учащихся.%
}
\newline
Не стоит слишком строго судить этот гейзер бреда,  лучше задуматься еще раз о зеркале, которое не меняет голову с ногами, и о том не видели ли мы похожих социальных прожектов
(в мечтах и в жизни) на своей
памяти, и не поддерживало\slash поддерживает ли их  образованное общество \dots
}.

 Он лично проваливает даже ново-Московское математическое общество.
 При том, что в Москве можно было опереться на блестящих
 лево-радикальных математиков (очень активные Гельфонд, Люстерник, Шнирельман, впрочем, левыми были 
 также
 Гливенко и Хинчин).
 Они  имели вкус к общественно-организационной деятельности
 и с немалой пользой  прилагали на этом поприще свои силы. Был Выгодский,
 который при всех своих небесспорных качествах был человеком деятельным и очень талантливым... 
 
 \sm
 
 По-видимому, длинный список мест, где Кольман работал, был связан как с талантом пускать
 пыль в глаза, так и с особым талантом проваливать любое дело. 
 
 Как мы видели, в июне 1932 член ЦК ВКП(б) А.~И.~Стецкий опубликовал в <<Правде>> телегу
 на Кольмана.
 
 Уже после этого (см. п. \ref{ss:yanov2}.а) (сентябрь 1932г.) Кольман устроил из Карла Маркса всесветское математическое посмешище в Цюрихе (и вряд ли это могло остаться тайной
 для его начальства).

 \sm
 
 В 1933-1935гг. Кольман продолжал работать в Комакадемии.
 Но каких-либо заметных должностей, 
  связанных с наукой,
 в это время за ним  не видно.
 
 Он, однако, сохранял должности в редакциях. 
 В частности, он  входил в редакцию  <<Под знаменем марксизма>> (видимо в 1931- 1942г., до 1935г я смотрел), 
 до 1935г. журнал издавался Комакадемией, далее, по-видимому, каким-то ее наследником.
 Автор просмотрел статьи Кольмана в этом журнале 1932-35гг. После статьи <<На текущие темы>> \cite{Kolman-temy}
 он в целом уклоняется от обсуждения математики и  уходит от показа советским математикам
 (и вообще советским ученым) диалектической кузькиной 
 матери, см. его статьи \cite{Kolman-1933}--\cite{Kolman-ver}.
 Они%
 \footnote{Я не утверждаю, что они радуют глаз.} заметно отличаются от его сочинений 1930-1932гг.,
 см. \cite{Kolman-Losev}--\cite{Kolm-boevye-voprosy}, восемь статей в сборнике   \cite{Dialektika},
 резолюцию конференции по планированию математики в \cite{Sbornik-planirovanie}.

  И вообще  дискуссии 1933--35гг. (очень многостатеен был А.~А.~Максимов) вряд ли имели то значение,
 какое имела идеология в момент натиска 1930-31гг.
 
 Кроме того,  Кольман долго сохранял должность редактора отдела математики БСЭ.
 
 \sm
 
 {\small  До 1931г. (тома 1-21) редактором отдела математики был Каган, он же заведовал всем естественно-научным
 направлением (хороший был редактор).% Выход 22 тома задержался {\it ввиду сложности работы по статьям
% <<Диалектика>> и <<Диалектический материализм>>} до 1935г. Остальные тома тем временем продолжали выходить.
 В 23 томе впервые не были указаны редакторы отделов (это понятно), позже такое несколько раз повторялось.
 Укажем редакторство по годам сдачи томов: в 1931-32 Кольман, в 1934 Кольман, Хинчин,
 в 1935 Кольман, Колмогоров, в 1936 Колмогоров, Кольман, в 1937 Колмогоров. 
 Я, возможно, недостаточно внимательно и последовательно просматривал
 математические статьи Энциклопедии, но, по-моему, кольманский дух наложил
 на нее не такой уж сильный отпечаток. Выше приводилась
 поносная статья про Егорова, но ветер вскоре изменился.
 
 Еще раз цитируем Стенограмму:
  \begin{quotation}
 ЛУЗИН.	Когда я послал статью из Парижа и сам приехал, то узнал, что эта статья
 	[<<Дифференциальное исчисление>> для Большой советской энциклопедии] 
 	вызвала жестокие нарекания, и тогда Вениамин Федорович [ Каган] заявил, что эту статью печатать нельзя и т.д.
 	Так
 	продолжалось около двух-трех лет. После этого ко мне обратились с просьбой эту статью развить. Когда я спросил, зачем это,
 	раз она вызвала такое противодействие, то
 	мне сказали, что я должен развить ее в том же направлении. Это мною было сделано,
 	и я ее развил.
 \end{quotation}
 
 Выше упоминалась
 статья <<Дифференциал>> Кольмана и Яновской, посвященная, на 80\% Марксу, Энгельсу и Ленину. По-видимому,
  это сочинение должно было стать образцом для подражания. Но - увы! -
  нет в мире совершенства. Как сообщается во вкладке
  к 23 тому БСЭ, <<В виду сложности работы по статьям <<Диалектика>> и <<Диалектический материализма>>,
  очередной 22-ой
  том Б.~С.~Э не смог выйти в срок>>%
  \footnote{В 1930г. были заклеймены А. М. Деборин и И. К. Луппол, и теперь диалектику нужно было срочно менять} (т.е. в 1931г.). Вышел он лишь в 1935г., когда эта тема 
  уже не была столь актуальной. Вслед за образцовым сочинением шла втрое б\'ольшая статья Лузина
  <<Дифференциальное исчисление>>%
  \footnote{Объем - 18
  	страниц современного журнала большого формата (см. статьи Лузина в <<Математическое
  	образование>>,  2005:2,  2005:3). В 1934г. вышла еще б\'ольшая статья Лузина <<Функция (в математике)>>.
  Впрочем,  к этой статье была приписка некоего (не скрывавшегося) Э.К. с философическими, но не воинственными комментариями, приведем конец:
\newline
{\it Функциональная зависимость не может исчерпать и заменить собой причинную связь.
Она отражает объективные закономерности, присущие самим вещам, но в еще менее
полной форме, чем причины. Функциональная зависимость отражает лишь количественную сторону
связей и отношений.}},  нормальные 
  обширные статьи Степанова <<Дифференциальные уравнения>> и
  Бурстина, Сретенского <<Дифференциальная геометрия>>, а также статья
  Кагана <<Дифференциальные инварианты>>...}

\sm

Посмотрим на участие Кольмана в математических съездах.

Летом 1930г. Кольман присутствовал на Первом Всесоюзном съезде математиков. 
Там было четыре общесъездовских идеологических доклада: Шмидт, Выгодский, Яновская, Каган.
Кстати там было и три секционных доклада по философии и истории математики. Кольмана среди докладчиков
не было. 

Весной 1931г. Кольман руководил конференцией по планированию математики. Тогда 
он в самом деле занимал важную должность.

 В 1934г. наш герой присутствует на Втором всесоюзном математическом съезде.
  Он  огласил приветствие съезду от
 Комакадемии. Вошел в Президиум съезда (всего там было 60 человек). Был одним из 5 (правда первым по счету)
 из последовательных председателей на секции <<Истории и философии математики>> (другие: Жегалкин, Нахимовская, Выгодский, Струве).
 Несколько скромная роль для <<Идеолога математики Всея Руси>>, к тому же не отличавшегося личной скромностью. На съезде состоялось пленарное  заседание,
 <<посвященное обсуждению вопроса о развитии научной работы по математике
 в СССР и задачах научно-исследовательских институтов>> (председателем был Александров). Казалось бы
 самое место для приложения Кольмановских талантов. Но и там его не было среди выступающих%
 \footnote{Съезд избрал Совет Всесоюзной математической
 	ассоциации в составе 41 человека, среди них был и Кольман (скорее всего, в качестве представителя Комакадемии). Президиум ассоциации был избран в составе: О.~Ю.~Шмидт (председатель), М.~Я.~Выгодский, А.~Н.~Колмогоров, Б.~И.~Сегал (секретарь),
 	В.~И.~Смирнов. Итак, Шмидт снова на коне
 	(его только что торжественно встречали Москва и Ленинград, но Съезде он по болезни
 	отсутствовал) и снова главный
 политик среди математиков и представитель математиков среди политиков. Характерно, что в Президиуме --
  два комиссара,
Выгодский и Сегал. Но не Кольман.}.

\sm

%В цитиремой ниже в п.\ref{ss:oslo} бумаге
В 1936г.
предлагалось включить Кольмана в советскую делегацию
на Международный математический конгресс в Осло \cite{Politburo}, документ 
Политбюро от 3 июля 1936%
\footnote{Резолюция была: в отправлении советской делегации  отказать. Обратите внимание на дату,
и см. даты в следующем параграфе. Не исключено, что делегация на Конгресс (он был 12-18 июля 1936)
не была отправлена в связи начавшимся делом Лузина. Но этот вопрос требует дополнительного продумывания.}.
Но здесь он был идеальным кандидатом: соглядатай с владением иностранными
языками. 
 
 Кольман всплывает в 1936 году в качестве зав.отдела науки Московского
 городского комитета ВКП(б). Это
 довольно высокая должность, но не стоит преувеличивать ее значения. 
 То, что он действительно ее занимал подтверждается статьей \cite{Kedrov-Kolman},
 двое из трех авторов этой статьи не могли не быть знакомыми с Кольманом в 30х годах.
 С другой стороны было бы интересно проверить даты, они могут и не точными (1936-1937гг. у Кольмана
 и 1936-1938гг. в \cite{Kedrov-Kolman}, не обязательно хотя бы один из этих вариантов точен).
 
 С.~С.~Демидов и В.~Д.~Есаков характеризуют Кольмана как
  <<{\it одного из влиятельнейших в сфере науки того времени лиц — ... заведующего Отделом науки МК партии тов. Э.~Кольмана}.>> Для начала обсудит вопрос, является ли первая часть фразы следствием второй. 
  Отдел науки, научно-технических изобретений и открытий ЦК ВКП(б) был организован в мае
  1935г. \cite{knowbysight} и дальше продолжал свое существование под меняющимися вывесками. По-видимому, Кольман
  был не первым по счету человеком на  должности главного ученого МГК ВКП(б), и уж точно не последним. 
  Можем ли мы сказать про других таких деятелей, что они   -- <<влиятельнейшие в сфере науки
  лица>>? Я предпринял разнообразные поиски по доступным электронным источникам, книгам, справочникам и ... не смог найти
  никаких следов этих деятелей (и даже не смог найти их имен). Во всяком случае есть основания
  усомниться в подразумевающейся импликации.
  
  Кольман описывает свою деятельность так:
  \begin{quotation}
  	 Работа в МК была интересная, но крайне трудная, сложная. У меня, правда, было двое помощников: зам, инженер Каплун, хороший, умный, рассудительный товарищ, занимавшийся главным образом изобретательскими делами, и Крапивинцев, только что окончивший врач-дерматолог, на которого я старался перекладывать все биолого-медицинские вопросы. Непонятно почему, Крапивинцева%
  	 \footnote{Я пытался проверить эти данные. Крапивинцев Павел Николаевич  был содокладчиком на 
  	 	VI пленуме Центрального совета Всесоюзного общества изобретателей (26-31 января 1937 г.).
  	 	Далее он действительно был полпредом в Литве 30.12.1937 — 8.10.1938, потом
  	 	работал в в полпредстве  в Париже, дважды был на приеме у Сталина  (5.10.1938 и
  	 	21.04.1939).
  	 	Больше ничего найти не удалось. Даже при наблюдавшемся в конце 1937г. кадровом кризисе
  	 	в СССР
  	 	назначение послом человека, не занимавшего руководящей должности, выглядит
  	 	странно.}, честного, но порядком ограниченного, малокультурного парня, никогда не соприкасавшегося с международными проблемами, внезапно направили на дипломатическую работу в буржуазную Литву.
  	 
  	  Но отдел не имел ни одного инструктора или референта, и справляться с огромным объемом работы — с сотнями научно-исследовательских институтов, с научной работой вузов столицы, с научными обществами и издательствами - а ведь все это (за исключением центральных, которыми ведал отдел науки ЦК) подлежало нашему руководству...
  	  было просто физически невозможно... 
  \end{quotation}
То, как было устроено партийно-государственное управление в 1930гг., по-ви\-ди\-мо\-му, является
одной из самых мало-исследованных исторических тем. У автора настоящих записок нет идей
о функциях, которые исполнял Отдел науки МГК ВКП(б) (а также, как делились
тогда сферы влияния МГК и наркоматов). Но картины, которую описывает Кольман,
быть не могло.

Список отделов 
Московского горкома ВКП(б)  можно найти в справочнике <<Вся Москва. Адресно-справочная книга>>, 1936.
Они соответствуют отделам ЦК ВКП(б), каковые можно найти в базе \cite{knowbysight}.
Учебными учреждениями ведал Отдел школ. Слово школа  в тогдашнем понимании
включало в себя младшую, среднюю и высшую школу (которая в свою очередь делилась по типам).
Так вот над Университетом и вузами находился Отдел школ. Издательства по партийной линии подчинялись
Отделу печати. Военизированные научные-технические учреждения по партийной линии должны были подчиняться
ГлавПУРу (трудно себе представить,  что кому-либо еще). 

Менее очевидна ситуация с академическими институтами. Но Академия была важнейшей стратегической структурой,
а достижение  лояльности и управляемости Академии и научных институтов было задачей, решаемой на самом высоком уровне
(есть сборник <<{\it  Академия наук в решениях Политбюро}>>, там виден прямой интерес
Политбюро и структур ЦК, в том числе и в частных вопросах).
Академия имела единую партийную организацию%
\footnote{Таковая упоминается в \cite{Zankevich} под 1934-35 годом со ссылкой на книгу
	<<Академия наук СССР VII-му
	Всесоюзному съезду советов>>, 1935, по приводимым там данным выходит, что это была
	общая партийная организация Академии, и у этой организации был парторг.
Стоит иметь в виду, что в академических институтах тогда было мало партийных и тогда была общая борьба
за усиление партийного влияния в них.} Кажется более чем правдоподобным, что
академические институты подчинялись по партийной линии Отделу науки ЦК, чем
персонажу, у которого не было ни одного референта.

У автора нет идей, чем на самом деле занимался в то время отдел науки,  но он  не руководил
(и не имел физической  возможности руководить) 
сотнями организаций упомянутых типов. А иначе представьте, какие горы взаимных ученых кляуз (это без референтов-то прочитать...),
каков объем  отчетов, какие труды по призывам к рационализации...

\sm

Добавим, что
 в 1931г. во время Великого перелома Кольман и Комакадемия имели 
свободу рук. В 1936г. Кольман входил в сложную структуру, где эта свобода была ограничена,
а должность, которую он занимал, скорее предполагала влияние, чем руководство. Сами
научные учреждения  параллельно зависели от быстро усиливавшихся в 30е годы наркоматов. 

\sm

Осенью 1936г. в связи с биологическими дискуссиями дух Кольмана воспрял,
а его таланты снова оказалось ко двору.
В журнале <<Под знаменем марксизма>> выходит его погромная статья 
<<{\it Черносотенный бред фашизма и наша медикобиологическая наука}>> \cite{Kolman-black}
(я ее, к сожалению, не читал).
Тогда же он атаковал С.~Г.~Левита, который, видимо не без его стараний, был исключен из партии,
однако директором института он оставался еще 6 месяцев до июля 1937г.

В марте 1937г. Кольман сам был атакован в Правде статьей А. А. Максимова \cite{Maximov-Kolman},
и в том же 1937    (по его словам%
\footnote{Он утверждает, что не был переизбран в МК партии, перевыборы были
в мае 1937г. Это можно проверить по газетам, я этого не делал.})  становится безработным. Опять всплывает в 1938-ом, вскоре участвует в нападках на 
вавиловцев и кольцовцев.
В частности, он участвовал в дискуссии по генетике 1939г., см. \cite{Genetika}, \cite{Mitin},
 был не последним участником этой дискуссии, но и никак не первым (стоит отметить, что к тому моменту число 
философов-марксистов несколько поубавилось, и на оставшихся был повышенный спрос).
В 1940г. опубликовал какую-то чушь в ДАН СССР, представленную академиком Т.~Д.~Лысенко
(я не читал).

Судя по тому, что мы знаем о   Кольмане, он 
в той или иной степени умел приносить вред человечеству, находясь в любой должности
(такого рода научные таланты бывают в любые времена).
Лично он мог иметь существенное влияние
на математику лишь в 1931г. в Звездный час Комакадемии. Дальше он мог   портить кровь тем или иным
конкретным лицам, но общая политика в отношении науки определялась  другими людьми и руководствовалась 
другими целями.

О деятельности Кольмана в 30е годы нам кое-что известно. Забавно, что
об этом
он в своих воспоминаниях почти не рассказывает... Если в чем-то (факте работы в Комакадемии и 
в отделе науки МК ВКП(б)) его слова соответствуют истине, то уровень своей значимости он высоко поднял над реальным,
а степень его контактов с высочайшими особами сильно преувеличена. Какие-то подобные контакты 
были, но опять-таки, едва ли стоит серьезно воспринимать его описания, даже когда мы предполагаем 
реальность таких контактов%
\footnote{Биография Кольмана 1917-1925гг., описанная в его воспоминаниях и юбилейной статье
\cite{Kedrov-Kolman}, похожи (хотя детали  расходятся), и, может, примерно так оно и было, с преувеличениями того же типа. Но не исключено, что авторы статьи знали о биографии Кольмана этого периода от самого Кольмана.}.

Впрочем, сейчас мы услышим оценку значения Кольмана из авторитетного источника 1936г.

\sm

{\bf\punct  Слова Кржижановского, обращенные к Лузину.%
\label{ss:krzhizh}}
Из Стенограммы:
\begin{quotation}
	КРЖИЖАНОВСКИЙ.
  Я не математик, но имею некоторое отношение к кругу ученых. Я
не знаю такой области познания, где бы серьезный ученый отрицал значение высокой
теории. {\bf Тем более мне кажется странным, что в области математики, где такое бесконечное совершенство теорий, 
громадный размах теоретических познаний, которые, в
конце концов, на наших глазах создают новые области и новые дисциплины,  чтобы серьезный ученый на основании толкований 
с Кольманом и с Каганом пришел к выводу, что
эти теоретические изыскания могут быть вредны моей родине.} Не знаю, как у Вас, у
крупного ученого (а Вы являетесь крупным ученым) могла зародиться мысль, что это
может принести вред и поэтому Вы не будете печатать. Это смехотворное толкование.
Вы, крупнейший теоретик, Вы должны знать и знаете, какие практические интересы
могут быть связаны с высочайшей теорией. Разве мы не видим, что математика доходит
до таких практических знаний, до которых она раньше не доходила? А Вы толкуете, что
{\bf Вы побеседовали с Каганом и Кольманом и решили, что это вред может принести. 
 И Вы
хотите, чтобы мы поверили этому? Это несерьезно.} 
\end{quotation}

И по смыслу то же самое, в более жесткой форме:
 \begin{quotation}
 	КРЖИЖАНОВСКИЙ.
  {\bf Вы являетесь блюстителем высокого звания советского академика.} Вы понимаете, какое
   орудие стране дает наука молодой, рвущейся вперед к знанию, страны? {\bf Вы облечены
   высоким доверием и говорите так спокойно — «по трусости, но малодушию сдал эту
   позицию».}
   
  {\bf Положение обязывает Вас. Вот в чем дело.
   Как мы можем сказать нашему ученому миру, что вы испугались Райкова, Кольмана,
   Кагана, испугались настроений.}
   \end{quotation}

Слова Кржижановского о Кольмане (да еще в перечне с Райковым и Каганом)
как раз соответствует тем должностям, которые  занимал реальный Кольман.
Добавим, что эти слова - индикатор действительного положения Кольмана, независимо от 
наших прочих изысканий.

\sm
 
{\bf \punct Довод Юшкевича.%
\label{ss:yushkevich-argument}}
Упомянутое Юшкевичем место книги Кольмана звучит так:
   \begin{quotation}
 Известно, что так называемая {\bf «Московская математическая
школа»} — Цингер, Бугаев, Некрасов — проповедовала, будто
«арифмология» (теория чисел и непрерывных функций) 
обосновывает индивидуализм, анализ с его непрерывностью 
направлен против революционных идей, теория вероятностей
подтверждает беспричинность явлений и свободу воли, а вся
математика в целом находится в соответствии с принципами
философии Лопатина -- {\bf православием,  самодержавием} и 
народностью».  
   
 Этот черносотенный образ мыслей был полностью донесен
до наших дней одним из «столпов» этой школы Лузиным,
который придал ему некоторую более {\bf «современную» 
фашистскую} окраску. Вместе с тем Лузин «исправил» эту идеологию
в деталях, заменив открытую проповедь  православия более
тонким дурманом — субъективным идеализмом и солипсизмом.
Так, например, натуральный ряд чисел провозглашается
Лузиным «функцией головы того математика, который в 
данном случае говорит о натуральном ряде»
\end{quotation}

Параллельное место в Правде такое:
\begin{quotation}
   Мы знаем, откуда вырос академик
Лузин: мы знаем, что он один из стаи бесславной царской {\bf  «Московской математической школы»},
философией которой было черносотенство и движущей идеей — киты российской реакции:
{\bf православие и  самодержавие}. Мы знаем, что и сейчас он недалек от
подобных взглядов, может быть, чуть-чуть {\bf фашистски модернизированных}.
Но социальная почва, взрастившая Лузиных, исчезла, ушла и уходит из-под ног. Академик Лузин
мог бы стать честным советским ученым, каких из старого поколения много. Он не захотел этого; 
он, Лузин, остался врагом, рассчитывая на силу социальной мимикрии, на
непроницаемость маски, им на себя напяленной.
Не выйдет, господин Лузин!.....  
\end{quotation}

Далее интересны выходные данные книги Кольмана:
<<Сдано в набор 2/IV 1936 г. Подписано к печати 14/VII 1936 г.>>. Статья в <<Правде>>
датирована 3/VII 1936 г., то есть книга Кольмана была еще не сдана в набор.
Итак, сентенция, которая была в находившейся в тот момент в печати 
книге, оказывается также в газете <<Правда>>, что ставит вопрос о связи этих событий
(теоретически возможно и обратное воздействие).

По мнению автора данных записок правдинский пассаж по смыслу восходит к Кольману 
(и чуть ниже приведен неочевидный довод в пользу этого тезиса). Будем из 
этого исходить.

Влечет ли это авторство Кольмана, как утверждается в \cite{DeEs}?
Приходят в голову два сценария.

\sm

{\it Сценарий} 1. Правдист, собирающий материал, обращается к политически надежным  математическим
кадрам,
а также к партийцам, которые на Лузина ругались.
Обратиться к Кольману было достаточно естественно.
 
 \sm
 
{\it Сценарий} 2. Математик  пишет статью, сам формулирует разные околонаучные обвинения,
а вот профессионально извергать политические ругательства он не умеет. Тогда он обращается
к Кольману, чтобы тот поучил его уму-разуму.

\sm

В обоих случаях получается участие Кольмана (что, кстати, утверждал сам Юшкевич),
но не авторство.

Однако и с прямым участием все не совсем очевидно.

\sm

1) (тривиальное). Как сообщается в предисловии в книге Кольмана, рукопись, по просьбе автора,
смотрели Яновская, Колмогоров, Александров (последний - 7-ую главу). Ее могли смотреть и другие люди, например, кто-то
мог оказывать помощь в редактировании книги.
Мы видели, что Александров и Колмогоров не участвовали в политических нападках на Лузина,
а их отношение к <<Московской математической школе>> было положительным (см. их статьи \cite{Alexandrov-1937},
\cite{Kolmogorov-matematika} 1937-1938, которые частично цитировались в пп. \ref{ss:detail}.2, \ref{ss:egorov-arest}).
Поэтому передача через них выглядит маловероятной, а вот про Яновскую такого не скажешь.

\sm

2) (более интересное) Кольман напечатал много статей (и надо думать, он много публично выступал,
а книга написана на основе его лекций в Комакадемии).
Статьи его звонкие, но несколько однообразные (ну типа того, что он знал философскую истину и ее повторял).
Выше, в п. \ref{ss:kolman-memoirs} приведен отрывок  о Лузине из статьи Кольмана марта 1931г. Там говорилось 
тоже самое... Сентенции Кольмана вполне могли попасть в газету и независимо от его рукописи.

\sm

Довод Юшкевича   серьезен, однако, не будучи ничем подкрепленным,
не несет  доказательной силы.

\sm

{\bf \punct Обсуждение статьи в Правде.%
\label{ss:pravda}}
1) Кольман оставил обширное литературное наследие. В его сочинениях 1930г.
 очень высока плотность фи\-ло\-соф\-ско-ме\-то\-до\-ло\-ги\-че\-ских сентенций
и связанных с этим слов (они,  наряду с обвинениями во вредительстве, составляют  основное содержание сочинений Кольмана).
В тексте из <<Правды>> нет ни одной философской или методологической фразы,
и лишь одно философское слово: <<философия Московской школы>>. Поэтому авторство Кольмана выглядит крайне неправдоподобно.
Крайне маловероятно и то, что он мог быть буквальным автором политизированной части текста (последние 2/5).

Кроме того, в тексте не высказывается сомнений в ценности  лузинского направления в науке%
\footnote{На этот счет Кольман любил высказаться, упоминая Лузина, типа <<{\it В этом отношении характерны
некоторые работы Бореля и особенно Лузина, уходящие в пустые абстракции}>> \cite{Kolman-listing}.
В Стенограмме и в газетной кампании такого рода взгляды не наблюдались.}, нет классового анализа,
 лексикон ругательств как-то бедноват, да и автор статьи не имел художественного таланта Кольмана. 
 Читатель может открыть любую погромную статью Кольмана
 и восхититься изяществом ее слога.

 В принципе здесь можно провести лингвистический анализ
 (хотя бы провести анализ лексикона ругательств и слов, используемых как ругательства).
 Если кому интересно - карты в руки.
  
 \sm
 
2) Текст из Правды  выглядит как составленный из двух частей. Первая часть выглядит, как написанная 
математиком, или по <<рыбе>> от какого-то математика. 
 Взгляды этого математика близки к взглядам Люстерника, Шнирельмана и Гельфонда%
 \footnote{Не следует в связи с этим выдвигать гипотез об авторстве.
 	Стоит отметить, что подобные взгляды имели и  не столь известные
 математики, так же, как и участники около-математической передовой тусовки.
Добавлю, что гипотезы об авторстве, скажем, Шнирельмана или Люстерника
кажутся мне крайне маловероятными
 (насколько я могу это оценивать  по их текстам, попадавшимся мне на глаза).},
высказывавшимися ими в 1936г. и в Декларации из Мат. Сборника 1930, номер 3-4,
приведенной выше в п \ref{ss:sbornik-perevorot}.
 И взгляды эти отличны от взглядов Александрова с Колмогоровым, Хинчина и Соболева
 (высказывавшихся в том же 1936г.).
 
 Дальше идет поток политизированных ругательств.
 Последние два абзаца действительно выглядят как переработка кольмановских сентенций,
 но переработка другим человеком, вычистившим оттуда философию.

 \sm
 
 {\sc Вероятные сценарии.} Откуда-то с ЦК
 приходит указание (причем провокация с предложением Лузину написать
 статью в Известия могла  иметь и значение проверки <<сигнала>>%
 \footnote{Слово было такое.}, а то, что Лузин написал могло и дополнительно разозлить 
 правдистов), дальше или газетчики собирают материал,
 или просят  какого-то математика (из числа не вызывающих сомнений) написать статью. Собственно, это те два сценария,
 которые указаны в предыдущем пункте. Что касается непосредственного обращения к Кольману
 за умом-разумом, то оно могло быть, а могло и не быть.

\sm

{\sc Следующие  статьи в <<Правде>>.} Следующие анонимные статьи в <<Правде>>
вышли 9 июля, 14 июля и 6 августа. Стоит отметить, что никаких следов философии и методологии в этих статьях нет,
то есть авторство Кольмана, на мой взгляд, исключено. Что касается прочих второстепенных в данном случае
 вопросов,
 скажу свою точку зрения, не пытаясь ее аргументировать. Три статьи имеют внутриправдинское происхождение, частью от журналиста, присутствовавшего
на заседаниях комиссии и ученых митингах, частью от правдинского главреда, который 9 июля задал главную
идею, а в статьи 14 июля и 6 августа вставил усилительные сентенции (вмешательство главреда 
и направление этого вмешательства
очевидны из обсуждения следующего параграфа, надо думать  ему лично газетная кампания обязана такой яростью).
По-видимому, есть некоторая возможность проанализировать эти тексты на предмет единого авторства.
Я бы осторожно предположил, что первая статья (от 3 июля) была написана другим человеком.

Выше в пп.\ref{ss:detail}.2 мы видели, что философы от математики  с некоторым запозданием
решили понападать на Лузина.

\sm

{\bf \punct  Дальнейшие доводы в пользу авторства Кольмана.%
\label{ss:kolman-tsagi}}
В 1990г появилась
 статья А.~Ф.~Левина \cite{Levin-1990}, где подстрочное примечание Юшкевича было объяснено
в развернутой форме. 

В статье С.~С.~Демидова и В.~Д.~кова \cite{DeEs} было много раз повторено, что Кольман, судя по всему,
 был автором статьи в <<Правде>>,
каковое повторение само по себе доказательной силы не имеет. Также была опубликована следующая бумага,
писанная Кольманом в его  звездный час:
\begin{quotation}
\small
«СЕКРЕТНО

Член Академии наук математик Н.~Лузин, избранный в 1929 г. по кафедре философии,
отказался подписать обращение советских учёных к заграничным по поводу процесса
Промпартии и в знак протеста против реорганизации Московского Математического
Института и Московского Математического Общества, президент которого ЕГОРОВ
арестован, ЛУЗИН оставил работу в Московском Математическом Институте и ушёл в
ЦАГи. Так как ЛУЗИН является специалистом по абстрактнейшей части теории множеств,
не имеющей никаких практических приложений, и в качестве руководителя так
называемой Московской Математической школы, он хвастает, что «никогда не решил ни
одного конкретного уравнения», то вряд ли в ЦАГи он может принести большую пользу.
Нужно подчеркнуть, что ЛУЗИН близко связан с виднейшим французским математиком 
БОРЕЛЕМ, активным сотрудником Французского военного ведомства. В бытность
свою в 1929 г. в Париже ЛУЗИН гостил у БОРЕЛЯ.

О воинствующем идеализме ЛУЗИНА красноречиво говорит следующая выдержка
из отчёта на заседании Академии о его заграничной поездке: «повидимому, 
натуральный ряд чисел не представляет из себя абсолютно объективного образования. По-видимому,
он представляет собой функцию головы того математика, который в данном случае
говорит о натуральном ряде. {\bf По-видимому, среди задач арифметики есть задачи 
абсолютно неразрешимые}%
\footnote{Даже эта цидуля несет содержательную
информацию по истории математики как истории идей. Воистину удивительны бывают
исторические источники. Впрочем, теорема Гёделя (K.~F.~G\"odel) о неполноте
({\it если формальная арифметика непротиворечива, то в ней существует невыводимое и неопровержимое высказывание})
была опубликована в том же, 1931 году. Лузин вряд ли рассматривал формулу Гёделя как задачу арифметики,
и имел в  виду задачи осмысленные. Но вскоре после Гёделя  стало ясно, что так оно и есть.
Возможно, что упомянутый Кольманом отчет (он есть в описи лузинского архива) содержит больше информации на
эту тему, чем обсуждаемый документ.}».
На эту тему ЛУЗИНЫМ во время командировки во Францию
написана книга и там же издана.

Кроме ЛУЗИНА в самое последнее время стал работать в ЦАГи профессор КАСТЕРИН,
ушедший демонстративно из Института Физики при 1 МГУ, где он вел разлагающую
антиобщественную работу.

22/2.31. 
\end{quotation}

К этой бумаге в статье прилагается странная сентенция:
\begin{quotation}
Донос Кольмана имел и реальную практическую направленность в условиях начавшейся подготовки
к Международному конгрессу математиков в Цюрихе, который намечался
на 4—12 сентября 1932 г., и, устраняя Лузина от участия в нем, призван был
обеспечить присутствие в Цюрихе молодых руководителей Московского математического
общества П.~С.~Александрова, Э.~Кольмана и др. 
\end{quotation}

При чем тут Конгресс, до которого было еще полтора года? при чем тут Александров?
Какое отношение присутствие Александрова на Конгрессе
(он сам был знаменитый математик)  имело отношение к присутствию или неприсутствию Лузина? 
Что касается бумаги, то смысл ее предельно ясен: Кольман хочет
помешать работе в ЦАГИ двум профессорам, бежавшим туда от погрома в
МГУ\footnote{Кастерину Николаю Петровичу
(1869-1947), кстати, члену Матобщества c 1901г., он, скорее всего,  помешал.
Во всяком случае Кастерин тогда же ушел из ЦАГИ,
но через какое-то время (установить которое мне не удалось) оказался в МГУ. Лузина Чаплыгин, у которого дела у самого были
в тот момент
неважные, по-видимому, смог отстоять.}. Бумага добавляет штрих к прекрасному портрету Кольмана%
\footnote{Ниже в п. \ref{ss:rasputie} мы добавим еще один похожий штрих.}, но 
ни имеет никакого отношения к вопросу об авторстве правдинской статьи пять лет спустя..

\sm

{\bf\punct Новый  виток.%
\label{ss:moral}}
Чуть выше стоит отрывок из \cite{DeEs} про Лузина, Кольмана, Александрова и Конгресс в Цюрихе.
Смысл той сентенции  очевиден - косвенно прицепить Александрова к 
 Кольману. 
К вопросительным знакам, поставленным выше, следует добавить некоторые замечания:

\sm

а) Александров
 с мая 1930 по июнь 1931г. был за границей (за это время свершился переворот в Мат.Обществе
и пал Шмидт).

\sm

б) по дороге в Москву, Александров (как мы видели) не испытывал восторга от происходившего
и сомневался, оставаться ли ему в МГУ.

\sm

в) Александров (чье присутствие на Конгрессе должна была, по мнению
авторов \cite{DeEs} обеспечить упомянутая бумага) не упомянут в Президиуме Математического 
общества сентября 1931г. (когда первым в списке был Кольман).
А из кольмановского президиума Мат.общества на Конгресс не поехал никто%
\footnote{%
Согласно \cite{DeEs} или \cite{Politburo}, 23.07.1932 Политбюро утвердило делегацию СССР в составе 
 Бернштейна,  Чеботарева,  Александрова,  Хинчина и тов. Э. Кольмана.
  Согласно
 трудам конгресса, из Советского Союза присутствовали
 Александров, Кольман, Чеботарев, Мюнц (Ch. H. M\"untz), М.~Акимов, М.~Куренский, М.~Кравчук,
 Н.~М.~Крылов, М.~Г.~Пфейфер, М.~Орлов, Н.~Парфентьев.
 Кольмановский список руководителей Мат.общества приведен в п.\ref{ss:obshchestvo-1931}. Из этого списка
 на Конгрессе оказался только сам Кольман.
Напомним, что  Лузину 
в итоге выезд был  разрешен, а не поехал он сам, скорее всего, по состоянию здоровья,
см. п. \ref{ss:luzin-do-1936}. Почему не поехали Хинчин и Бернштейн, история умалчивает. 
Отметим также странность самого списка и наличие в нем лишь двух представителей Москвы, из коих один был Кольман 
(пять представителей от Украины, два от Ленинграда,
два от Казани).\label{fo:zurich}}, кроме самого Кольмана.

\sm

Продолжаем цитирование \cite{DeEs}, делая попутные реплики в сносках:
\begin{quotation}
 П.~С.~Александров
владычествовал в Московском математическом обществе — основной в Москве
 (а, следовательно, в СССР%
 \footnote{Логика восхищает. А Мехмат? А Стекловка? А математическая группа АН СССР?
 А Матмех?})
математической организации. Он находил понимание у одного из
влиятельнейших в сфере науки того времени лиц — у заведующего Отделом науки МК
партии тов. Э.~Кольмана%
\footnote{Напомним, Александров <<владычествовал>> с мая 1932, а Кольман был
	 на упомянутой должности с  1936
года.}....

Движущей силой интриги выступил главный редактор «Правды» Л.~З.~Мехлис. Ему
была подброшена идея — использовать собранные ({\bf Кольманом и его «молодыми» друзьями})
антилузинские материалы для развёртывания кампании по борьбе за воспитание
советского патриотизма у отечественных учёных.
\end{quotation}
\begin{quotation}
Отметим также, что тема раболепия и низкопоклонства перед Западом, составлявшая
суть социального заказа, утонула в различных обвинениях, представлявшихся, конечно,
очень важными {\bf «молодым советским математикам» и идеологическим соратникам
Э.~Кольмана}. 
\end{quotation}
\begin{quotation}
Как
только положение Кольмана на партийном Олимпе
пошатнулось — он потерял свой пост в
МК партии — математики (А.~О.~Гельфонд и Л.~Г.~Шнирельман) поторопились нанести
удар по его математической репутации, опубликовав в «Успехах математических наук»
разгромную рецензию
 на вышедшую ещё в 1936 г. (!)%
 \footnote{Интересно, что это такая  за репутация? и в чем смысл восклицательного знака? 
 почему Гельфонд со Шнирельманом должны были немедленно писать рецензию?
 Шнирельман поучаствовал в 1937 году в борьбе за качественные школьные учебники, что было важнее.
 И, кстати, в 1937г. был лишь один выпуск <<Успехов>>.} его книгу «Предмет
 и метод современной математики».

 А так как существовать без штатных идеологов в советское время не
полагалось никакой науке, они постарались максимально ограничить их круг в математическом 
сообществе небольшим числом лиц, которым они доверяли. Наиболее влиятельным
идеологом в московском математическом сообществе, особо поддерживаемым университетскими
математиками, стала С.~А.~Яновская, известный логик, философ и историк математики,
много сделавшая для развития в СССР исследований по математической логике%
\footnote{На этом основную теорему <<социальной истории математики>> можно считать
полностью доказанной.}.
\end{quotation}

{\bf Если кто не понял, то <<молодые друзья Кольмана>> и <<идеологические соратники Кольмана>>, это те, кто выступали против Лузина на Академической комиссии. То есть 
Александров, Гельфонд, Колмогоров, Люстерник,  Хинчин, Шнирельман.}

\bigskip

{\bf \punct  Диалектическая спираль.%
\label{ss:spiral}} Подытожим наше историографическое исследование.

Как известно, все люди не без греха,  в том числе люди выдающиеся, у которых к тому больше соблазнов и больше возможностей.
В Московской математической сваре 1936г. не в самых привлекательных ролях участвовал ряд знаменитых
мужей, имевших огромные научные заслуги и в своей прочей жизни (которая, впрочем,
 тоже была не без греха, как и у всех смертных),
большинство из них совершили много положительного и оставили по себе добрую память.
Некоторых их них помнят ныне живущие.
Да и сама великая эпоха советской математики осталась памятником их деятельности.

В число участников свары на рьяных ролях входили и друзья Юшкевича -- человека весьма осведомленного, и, почти наверняка знавшего,
что было на самом деле.
В первой своей статье он (использовав публикацию письма Капицы)
превратил эту  свару в <<очередной урок смирения научной интеллигенции>>,  преподанный Сталиным
(что вполне соответствует  перестроечной исторической стилистике перевешивания на Сталина грехов всех и вся).
По-видимому, это единственное место, где в своей статье он  серьезно отступил от истины.
Рассказ о Кольмане там вполне реален.

Во второй статье Юшкевич решил воспользоваться хлестаковскими воспоминаниями Кольмана,
дабы превратить его в героя, вертевшего судьбами советской математики, и навесить 
на него прегрешения всех прочих действующих лиц (кроме, разумеется, неприкосновенного в этом отношении Сталина). 
По-человечески, эта идея понятна. Кольмана уж точно не жалко, а с людей, по большому счету заслуживающих
огромное уважение, снимались их прегрешения.  Понятно, что те, кто знал больше, чем
писалось в <<Успехах>>, помалкивали, кто ж будет оправдывать
  Кольмана, при этом  обвиняя, скажем, 
Александрова, Соболева, Хинчина, Колмогорова или Люстерника?

Так или иначе, Кольман в новом изложении истории начал превращаться в эпического батыра, исшедшего из ада,
и
образ этот стал дорог пост-советским математикам.

Все было бы ничего, но теория Юшкевича  могла существовать лишь 
при неопубликованной стенограмме Академической комиссии. Дальше нужно было что-то менять.
Или попробовать отказаться от образа батыра Кольмана, заодно сняв со Сталина прегрешение, которого он не совершал
(каковое снятие  само по себе было бы грехом смертным).
Это было невозможно,  такой издательский проект не был бы поддержан (какая же это тогда <<социальная история>>?),
и такое не
приняла бы научная общественность 90х годов. Или сделать то, что было сделано, и
к широкому удовольствию деградировавших наследников великой эпохи
превратить Александрова, Гельфонда, Колмогорова, Люстерника, Хинчина, Шнирельмана
и примкнувших к ним Соболева и Шмидта в прислужников Батыра, участвовавших в преподании Сталиным урока советской интеллигенции.

     \section{Игра с огнем%
     \label{s:fire}}
     
     \epigraph{Поражен неожиданными совершенно незаслуженными газетными нападками на Вас.
     	 Ваш высокий всемирно признанный научный авторитет не может
     	быть поколеблен. Твердо надеюсь, Вы найдете в себе силы спокойно отнестись
     	к малоавторитетной критике Ваших трудов. О совершенно необоснованных
     	обвинениях другого порядка не говорю}{Телеграмма Чаплыгина Лузину, начало июля 1936г.}
     
     \COUNTERS
     
     Вернемся, наконец, в горячее лето 1936г. 
     Так или иначе, в игру должны были включиться партийные круги. Опубликованные
     на сегодняшний день дополнительные партийные материалы \cite{DeEs} состоят из двух писем Мехлиса Сталину.
     
     \sm
     
     {\bf\punct Мехлис.%
     \label{ss:mehlis}} Мехлис Лев Захарович (1889-1953) в это время был главным редактором  
     <<Правды>>. 
     
    {\small
     Основные вехи его биографии.
     
     $\bullet$ Член  социал-демократической сионистской партии «Поалей Цион» (1907-1910),
 
  $\bullet$   с 3.1918
     член РКП(б), в Гражданскую войну на различных комиссарских должностях,
  
   $\bullet$   в 1922-30гг. -- помощник Генерального секретаря ВКП(б).
   
    $\bullet$ В 
     1931 - 12.1937
     главный редактор газеты «Правда»,
     
      $\bullet$ с 10.2.1934 -- кандидат в члены ЦК ВКП(б),
     с 12.10.1937 -- член ЦК,
     
      $\bullet$ 30.12.1937 - 6.9.1940
     начальник Политического управления РККА (и зам. наркома обороны),
     
      $\bullet$ 19.1.1938 - 5.10.1952
     член Организационного бюро ЦК ВКП(б), 6.9.1940 -1941, 1946 - 27.10.1950
     
      $\bullet$ 6.9.1940 - 27.10.1950 народный комиссар советского контроля%
     \footnote{Наркомат советского контроля, бывший Рабкрин, список последовательных
     руководителей до Мехлиса (весьма внушительный):
     Сталин, Цюрупа, Куйбышев, Орджоникидзе, Андреев, Рудзутак, Куйбышев, Антипов, Косиор, Беленький, Землячка.
 Напомним, что до 1950г. это была главная антикоррупционная структура в СССР.},
    к этой должности прилагалась должность
    <<заместитель председателя Совнаркома СССР>>.  позднее министр государственного контроля СССР.

    $\bullet$ 10.7.1941 - 4.6.1942
    заместитель народного комиссара обороны СССР, армейский комиссар I-го ранга
    
    $\bullet$
    .1 - 19.5.1942
    представитель Ставки Верховного Главнокомандования на Крымском фронте, армейский комиссар I-го ранга
    
     $\bullet$ Далее член Военных советов ряда фронтов
    }

\sm
  
     Военные его очень не любили. Основное предъявлявшееся ему впоследствии обвинение
     -- разгром Крымского фронта в мае 1942г.
      Как утверждает современный военный историк А.~В.~Исаев, изучавший советские и немецкие
      документы того времени, степень ответственности Мехлиса
       за эту катастрофу была сильно преувеличена.
     
     Человек жестокий,  по-видимому, принципиальный, а также преданный лично Сталину,
       он много
     находился на контрольно-надзорных должностях, наделал там много дел
     и нажил много врагов. Вскоре после смерти стал одним из эталонных демонов советской
     истории. При вероятной преувеличенности его дурной славы, он был человеком,
     от которого стоило держаться как можно дальше.
     
     Однако, в обсуждаемый момент времени Мехлис был редактором газеты <<Правда>>, а  его
     надзорные функции и вызываемый ими страх у надзираемых был еще впереди,  математики
     не могли понимать, в чьих руках, отчасти, оказалось <<дело Лузина>>.
     
     \sm

     {\bf \punct <<Традиции раболепия>>.%
     \label{ss:rabolepie}}  Первая статья в <<Правде>>  вышла 3 июля,
     в тот же день состоялся митинг в Институте Стеклова. C.~C.~Демидов и В.~А.~Есаков
     \cite{DeEs}
     приводят следующий документ:
     \begin{quotation}
     	
     	3 июля 1936 года

     	Товарищам Сталину, Кагановичу, Андрееву, Жданову, Ежову
     
     	Тов. Молотову

    	Материалы, собранные редакцией «Правды» в связи с делом академика Н.~Лузина,
     	выявили, между прочим, один серьезного значения недостаток в работе научных организаций. Сводится этот недостаток к тому, что большинство ученых наиболее интересные
     	свои работы считают нужным публиковать главным образом и раньше всего не в СССР,
     	а в заграничной печати. Вызывается это двоякого рода причинами:
     	
     	во-первых, ненадежностью издания научных книг и журналов у нас, в СССР,
     	
     	во-вторых, тем ореолом уважения, которым до сих пор окружена в научной среде
     	(даже среди многих коммунистов — научных работников) любая, хотя бы и малозначительная работа, если она напечатана за границей.
     	
     	Считая такое положение совершенно ненормальным, прошу ЦК ВКП(б) санкционировать развернутое выступление по этому вопросу на страницах «Правды».
     	
     	Редактор «Правды» Л.~Мехлис

     	<резолюция: Молотову! Кажется, можно разрешить. Сталин
     	роспись В.~М.~Молотова и помета «Сообщ[ено] т.~Мехлису. П[оскребышев]»>
     	ЦК ВКП(б) —
     \end{quotation}

     Мы видим, что уже 3 июля 1936 года (день публикации в <<Правде>>) об истории услышали высшие лица страны. Стоит заметить, что инициатива идет от Мехлиса
     к Сталину, резолюция от последнего -- <<кажется, можно>>. То есть это не выглядит
     как исполнение генерального плана вождя.
     
     Последовала статья в <<Правде>> 9 июля:
     \begin{quotation}
     	До сих пор считается почему-то нормальным,
     	естественным в научной среде печатать труды советских ученых, прежде всего,
     	за границей или даже только за границей.\dots
     	
     	Дошло до того, что даже популярные работы (по топологии, теории вероятностей)
     	профессоров Александрова, Хинчина, Колмогорова впервые напечатаны были за границей на немецком языке, а затем только был «поднят вопрос» о переводе этих работ советских ученых на русский язык и переиздании в СССР%
     	\footnote{Были ли такие работы? Кстати книга Alexandroff-Hopf не была переведена на русский,
     		и быть может это связано с данными инвективами.}. Чувство глубокого возмущения
     	охватывает, когда знакомишься с такого рода фактами. Вскормленные советской страной и советской наукой подобные ученые продолжают раболепствовать перед любой чужой
     	страной, любым чужим языком. Или они не обладают ни каплей национальной гордости?
     	
     	Как
     	можно характеризовать, например, линию поведения физика Акулова, который до сих
     	пор публикует свои работы в германском фашистском журнале, где даже либеральные
     	буржуазные ученые в виде протеста против мракобесов отказались печататься!

     	Как же может потерпеть советская общественность, чтобы плоды работ наших,
     	советских, ученых были отдалены от десятков и сотен тысяч советской молодежи, жаждущей знаний! Как же можем мы терпеть, чтобы сами наши советские ученые распространяли и внутри страны и за границей нелепый взгляд о второразрядности советской
     	науки!
     	
     	Нередко слышатся ссылки на нашу полиграфическую бедность, на нехватку бумаги
     	для научных изданий и тому подобное. Эти ссылки не имеют под собой реальной почвы.
     	{\bf Советское государство хочет и может обеспечить целиком, на сто процентов, своевременное напечатание всех ценных работ всех советских научных работников.}
     [выделено в газете] Для этого
     	важнейшего дела есть и сейчас, будет все больше и бумаги и полиграфической техники.
     \end{quotation}
     
     Далее шла критика издательств.
     
     Стоит заметить, что идеи этого типа высказывались не только в партийных верхах. 
     Вот что говорил директор издательства ОНТИ Аршон на заседании Комиссии
     9 июля (кстати, он был одним из выступавших на митинге 3 июля в Стекловке, \cite{Front1936}):    
     \begin{quotation}
     	АРШОН. Мне кажется, что сводить вопрос к тому, что у нас мало бумаги и
     	медленно печатается, а поэтому требовалось печатать за границей, нельзя. Очень приятно слышать, когда Павел Сергеевич говорит, что, получивши щелчок от «Правды», он
     	перестроится, а раньше это было результатом недомыслия. {\bf Надо сказать, что, в частности, наше издательство в течение трех лет неоднократно ставило этот вопрос, и мы получали забавные доводы за то, что именно так и должно быть, т.е. что надо печатать за границей. Мы ставили вопрос так: все то, что имеется ценного в науке, созданного нашими
     	работниками, должно печататься на русском языке; это не должно помешать печатанию
     	на других языках.} И вот интересно, какие доводы приводились против этого: если мы будем печатать наши работы на русском языке, то наши молодые аспиранты и научные
     	работники не будут иметь стимула к изучению иностранных языков, поэтому надо печатать на иностранных языках. Второй довод: по специальным вопросам слишком узкий
     	круг ученых, эти ученые главным образом за границей и, конечно, они на русском языке
     	читать не будут. Как будто бы тот, кто интересуется своей областью, не прочтет то, что
     	напечатано, хотя бы на китайском языке, как это делают наши ученые: когда они чем-нибудь интересуются, они прочитают это на любом языке, и, тем более, что никогда не ставился вопрос так, что не нужно печатать одновременно и на иностранном языке, но на
     	русском — обязательно. Вопрос о том, что нужно создавать для огромной массы людей, научных кадров капитал на русском языке, чтобы это было в наших библиотеках,
     	отводился тем, что зачем это делать, когда он все равно прочтет на иностранном языке.
     	И последнее, что всегда так было и так останется.
     	
     	Мне кажется, что здесь в целом, в отношении большинства это неверно, 
     	но у меньшинства, к которому принадлежит Лузин, безусловно, было стремление рассматривать
     	СССР, как задний двор культуры, а культура там, на Западе.
     \end{quotation}
     
     \begin{quotation}
     	АРШОН. Мы с удовольствием печатаем большие труды на иностранном языке, но
     	с тем, чтобы они были напечатаны и на русском.
     \end{quotation}
     
     Выступление Аршона неоднозначно. Оно содержит очевидные элементы
     демагогии, однако вопрос о запрете печатанья за границей не стоит (не стоит он и у Мехлиса). А огромная научная литература на русском языке была в итоге создана,
     и русский (по-крайней мере, в математике) становится одним из важнейших научных языков
     (что, впрочем, в неменьшей степени связано с уровнем советской математики, кстати издание собственных
     журналов было неотъемлимым элементов самоорганизации математического сообщества).

     Но вернувшись к нашей теме, отметим, что  
      идеи, высказанные Мехлисом, были не только у Аршона (см. выше п. \ref{ss:sbornik-perevorot}),  
      если  <<Правда>> в самом деле собирала 
     <<материалы>>, она должна была такие идеи услышать от весьма авторитетных
     и наиболее близких к партийным кругам лиц.
     
     В общем, взгляд Мехлиса на этот вопрос оказался сходным взглядам
     леворадикальных математиков. Что ж, бывает.
     
     \sm

     {\bf \punct  Реакция академических кругов.%
     \label{ss:akademiki}}
     Если математики с энтузиазмом решали свои собственные задачи,
     то представители других наук должны были испытать страх 
     и должны были начать поиск средств к спасению.
        Реакция началась незамедлительно уже 3 июля.

  \begin{quotation}
   Вернадский -- Ферсману,	3 VII 1936, Узкое
  	
  	Дорогой Александр Евгеньевич,
  	
  	Мне кажется, Академия не может пройти мимо той <статьи>, которая напечатана в <<Правде>> о Лузине. Крупнейший математик и человек глубокой честности будет заниматься плагиатом у своих учеников! Неверно, что этот человек не работает. Я не знаю, в какой форме надо выступить, но первым делом надо заставить газету напечатать ответ Лузина. Напрасно он, конечно, тратит время на ревизию школы, и, очевидно, его отзыв
  	не отвечал оценке педагогов! Но прежде всего нельзя забывать, что это один крупнейших математиков не только нашей страны. Я пишу С.~А.~Чаплыгину, который сейчас в Ессентуках.
  	
  	Так радостно видеть Вас бодрым и полным,
  	настоящим человеком.

\phantom{.}\hfill  	Ваш В.~Вернадский
  \end{quotation} 
  
  \begin{quotation}
  	[Сопроводительное письмо]
  	Узкое
  	7. VII. 1936
  	
  	Дорогой Александр Евгеньевич,
  	
  	направляю к Вам письмо (к которому присоединился Н.~С.~Курнаков), к Вам как к председателю нашего Отделения. Дайте ему ход только в случае, если Вы найдете это нужным
  	и невредным. Мы не можем, конечно, издалека учесть точно положения. Если Вы в этом
  	откажете <неразборчиво>, верните письмо мне.
  	
  	Я думаю, что подобная история может оказаться, в конце концов, гибельной для Академии,
  	 если она приведет к удалению H.~H. [Лузина] из Академии или чему-нибудь подобному. Мы покатимся вниз по наклонной плоскости.
  	
  	Ваш
  	В.~И.~Вернадский
  	
  	P.S. Копию моего заявления пошлите Горб[унову], А[лександр] Евгеньевич], если Вы
  	не будете против этого возражать.
  	
  	Прилагаю письмо Насонова [Н.~В.] (он забыл поставить число), он также считает, что
  	если Вы считаете неубедительными наши письма — то верните мне.
  \end{quotation} 
  
    <<Отделение>> -- это отделение математических и естественных наук АН СССР
    (ОМЕН).  Заявление Вернадского и письмо Насонова приводятся в \cite{Delo}. 
  
  Из дневника Вернадского:
  \begin{quotation} 
  	9.VII.936 
  		Узкое 
  		
  		История с Лузиным не дает покоя. Был 6го Н~Д.~Зел[инский] и хотел подать общее заявление. Общее заявление при нын[ешних] обстоятельствах немыслимо. Я предложил отдельные (письма), он хотел Курнак[ова] [1]. Курнаков подписал мое председателю] Отделения] Ферсм[ану]. Послал свое и Насон[ов] 7го, вчера говорил с Зел[инским]. Он говорит, что А~Е. (Ферсман) против так[ого] заяв[ления] считает его неверным. Согласно указанию Курнакова {\bf я и раньше дал в этом отношении carte blanche А.Е.} [Ферсману]. Совершенно ясно, что в таких случаях, если презид[иум] АН и решит выступить члены АН ничего не смогут сделать. Ист[ория] с Луз[иным] не ясна. Очевидно, какие-то интриги. Раньше уже было. {\bf  Тогда К[ольман]  идеолог математики%
  			\footnote{Напомним, что он им в самом деле был в 1931г.},
  			 теперь тоже, очевидно, математики. Говорили об Александрове, который, хотя ученик Лузина, по-видимому, из карьерных соображений действует.}
  \end{quotation}
  
 Ответ Чаплыгина Вернадскому опубликован в \cite{Vernad-perepiska} и (с небольшим сокращением) Юшкевичем в \cite{Yush2}. 
Документы, связанные с Вернадским, усиленно публиковались. Естественно думать, что
должны были быть другие такие академические письма и другие академические усилия.

Мы видим, что в оппозиции к правдинской статье оказались весьма влиятельные члены академии.
Среди них был 
П. Л. Капица [Peter Kapitza]
Приводим полностью
его знаменитое (впоследствии) письмо : 
  \begin{quotation}
  	\small
 <резолюция: За ненадобностью вернуть гр-ну Капице. Молотов>
 
 Председателю СНК СССР В.~М.~Молотову
 
 6 июля 1936, дер. Жуковка дача 33
 
 Товарищ Молотов,
 
 Статья в «Правде» о Лузине меня озадачила, поразила и возмутила, и как советский
 ученый я чувствую, что должен сказать Вам, что я думаю по этому поводу.
 Лузин обвиняется во многом, я не знаю, справедливы ли эти обвинения, но я вполне
 допускаю, что они полностью обоснованы, но и в этом случае мое отрицательное отношение к статье не изменится.
 
 Сперва начну с некоторых обвинений Лузина мелкого характера. Он печатал свои
 лучшие работы не в Союзе. Это делают многие ученые у нас главным образом по двум
 причинам: 1) у нас скверно печатают — бумага, печать; 2) по международному обычаю,
 приоритет дается только (в том случае), если работы напечатаны по-французски, немецки или английски. Если же Лузин печатал в Союзе плохие работы, то в этом виноваты
 редакции журналов, которые их принимали.
 
 То, что он завидовал своим ученикам, и поэтому были случаи несправедливого отношения к ним,
 то, к сожалению, это явление встречается даже среди самых крупных ученых...
 Итак, остается одно обвинение против Лузина, он скрывал за лестью свои антисоветские настроения, хотя каких-либо больших преступлений не указывается. Тут стоит по
 существу очень важный и принципиальный вопрос: как относиться к ученому, если
 морально он не отвечает запросам эпохи.
 
 Ньютон, давший человечеству закон тяготения, был религиозный маньяк. Кардан,
 давший корни кубического уравнения и ряд важнейших открытий в механике, бьл кутила
 и развратник. Что бы Вы с ними сделали, если бы они жили у нас в Союзе?
 Положим, у Вас заболел близкий Вам человек. Позвали бы (Вы) гениального врача,
 если бы даже его моральные и политические убеждения Вам были противны?
 
 Возьмем пример ближе. Гениального Клода, изобретателя процессов ожижения,
 разделения и получения целого ряда газов, которыми пользуются теперь всюду в мире и
 у нас в Союзе. Он — французский фашист. Что бы Вы с ним стали делать, если бы он был
 советским гражданином и не захотел менять своих убеждений?
 
 Я не хочу защищать моральные качества Лузина. ...Но нет сомнений, что он наш
 крупнейший математик, один из четырех самых лучших наших математиков, его вклад в
 мировую науку признается всеми математиками, как у нас, так и за границей. К тому же,
 он сделал больше чем кто-либо другой из наших математиков, чтобы собрать и воспитать ту плеяду советских математиков, которую мы сейчас имеем в Союзе.
 Я считаю, что страна, имеющая (таких) крупных ученых, как Лузин, должна первым
 делом сделать все, чтобы его способности были наиболее полно использованы для человечества.
 
 Людей типа Лузина, идеологически нам не подходящих, во-первых, надо поставить в
 такие условия, чтобы они, продолжая работать в своей научной области, не имели широкого общественного влияния, во-вторых, нужно сделать все возможное, чтобы их перевоспитать в духе эпохи и сделать из них хороших советских граждан.
 Начнем с первого. Что Лузин не социалист, об этом, конечно, знали все в Академии,
 и таких там не мало, и, конечно, это не было неожиданно открыто директором 16-й школы после того, как Лузин разразился льстивыми комплиментами. Но, несмотря на это,
 его выбирают выполнять целый ряд общественных обязанностей, его просят рецензировать, ему поручают руководство Математической группой Академии наук ...
 
 Во-вторых, было ли сделано все возможное для того, чтобы перевоспитать Лузина и
 людей типа Лузина в Академии наук и можно ли это достигнуть такими методами, как
 статья в «Правде»?
 
 Я утверждаю, что нет, а как раз наоборот — этим затрудняется воспитание не только
 самого Лузина, но и целого ряда других ученых.
 Как вообще Вы взялись за перестройку Академии? Вы, первое, начали выбирать в
 академики партийных товарищей. Это был бы лучший метод, если бы у нас были крупные ученые среди партийцев. Оставляя в стороне общественные науки, наши партийные академики куда слабее старых, их авторитет поэтому мал.
 
 Вырастить новых ученых из молодежи тоже пока не удается. Я это объясняю совсем
 неправильным подходом с Вашей стороны к науке, чересчур узко утилитарным и недостаточно внимательным. Поэтому главный научный капитал у нас все же лежит в старом
 поколении людей, доставшихся (нам) по наследству. Поэтому следовало бы, казалось, все
 сделать, чтобы их перевоспитать, приучить и пр. Но то, что Вы делаете, совсем не достигает цели. Когда-то арестовали Лазарева, прогнали Сперанского, а теперь обрушились на
 Лузина. Немудрено, что от такого нежного обращения (такие) ученые как Успенский,
 Чичибабин, Ипатьев и другие, сбежали. Я по себе знаю, как бездушно вы можете обращаться с людьми.
 
 Возьмем, далее, тех партийных товарищей, которых Вы посылаете работать с учеными и которые, если хорошо подобраны, могли бы прекрасно перевоспитать нашу ушедшую от жизни научную среду. Ведь все время среди них обнаруживаются такие товарищи, за которых краснеть приходится. Я это по своему опыту знаю. Ведь того заместителя,
 которого мне в начале дали, я не могу иначе назвать, как совсем беспринципным человеком. Теперешнего же моего «зама», лучше которого я себе не желаю, я сам себе нашел. Правда, когда я его попросил, было сделано все, чтобы я его заполучил. И я уверен,
 если бы у всех директоров институтов Академии наук были бы такие же замы, как у
 меня, то дух Академии наук совсем бы изменился.
 
 Что Вы сделали, чтобы перевоспитать Лузина? Ничего. А чего достигнет эта статья в
 «Правде»? Либо он еще слащавее заговорит, либо у него произойдет нервное расстройство,
 и он прекратит научно работать. Только перепугаете, больше ничего. Пугать надо опасных
 врагов. А разве Лузины опасны Советскому Союзу? Новая конституция лучше, чем что-либо другое, показывает, что Союз достаточно мощен, чтобы не бояться Лузиных.
 Но вот, имея в руках все хозяйственные достижения, политические завоевания, которыми располагает наш Союз, я не понимаю, как можно не перевоспитать любого академика, каким бы он не был, стоит только внимательно взяться и подойти индивидуально.
 Пример — хотя бы Павлов. А крупных ученых у нас не так много, чтобы за это дело
 трудно было взяться.
 
 Из всех этих соображений, я не могу понять, какой тактический смысл статьи в
 «Правде», и вижу в ней только вредный шаг для нашей науки и для Академии, так как это
 не перевоспитывает наших ученых и не подымает их престиж в стране.
 А если к этому прибавить, что имя Лузина достаточно хорошо известно на западе,
 чтобы такая статья не прошла незамеченной. Благодаря своей слабости и неубедительности (она) может быть комментирована самым разнообразным и нелепым образом. Видя
 вред всего случившегося для науки в Союзе, я считаю, что должен об этом написать Вам.
 
 П. Капица
 \end{quotation}
 
 Капица за свою жизнь написал немало писем руководителям страны в защиту разных
 ученых. Известно, что в ряде случаев к нему прислушивались. Было ли так в 
 деле Лузина, мы не знаем.
 
 \sm

  {\bf\punct  Кржижановский в вышестоящих инстанциях.%
  \label{ss:krzhizh1}}
  Напомним, что 11 июля прозвучали слова о внесении в резолюцию формулировки 
  «\dots что
  H.~H.~Лузин своей деятельностью за последние годы приносил вред советской науке и
  Советскому Союзу». На следующем заседании, 13 июля Кржижановский
  сообщает: 
  \begin{quotation}
  	КРЖИЖАНОВСКИЙ.
   {\bf Та резолюция, которую мы писали, признана правильной, выдержанной, и нет никаких оснований ее переделывать. Но некоторые пожелания здесь есть, законные пожелания.} Вот что желают от
  нас: в самой резолюции нужно дать побольше фактического материала. По поводу каких
  пунктов? Например, по пункту третьему, там, где говорится о низкопоклонстве. Здесь
  надо непременно сделать несколько цитат, таких, которые это иллюстрируют. Я, конечно,
  это не в состоянии сделать, но вам это очень легко сделать. Затем в пункте относительно
  отзывов — привести с десяток фамилий, что можно сделать после того разбора, который
  мы сделали.\dots
  {\bf Здесь обращают внимание на то заключение, которое мы сделали: «...полностью
  подтверждает характеристику Лузина, данную в газете «Правда», как врага в советской
  маске». Есть совет это заключение переделать в таком духе}, чтобы здесь была совершенно самостоятельная мысль, чтобы «не плагиировать» из «Правды» и «переноса» не
  делать, а {\bf сказать таким образом: поступок Лузина является недостойным советского
  ученого, к тому же действительного члена Академии наук, а также несовместим с достоинством, которое должно быть у каждого советского гражданина.} Это уже не будет плагиатом.
 % 
 % Переходя же к остальной части резолюции, право, я затрудняюсь сказать, что здесь
 % можно изменить. Вот п[ункт] 5 надо, пожалуй, изменить. {\bf Здесь написано — «прямой
 % плагиат». Прямого плагиата здесь нет, здесь дело обстоит тоньше.} Вот мы с вами убедились, что Лузин обыкновенно очень мало выражается, а на самом деле принадлежит к
 % числу очень многословных ораторов. Почему? Он защищает здесь позицию реакционного крыла академиков, со страстью защищает его плоть и кровь. И вот нам нельзя дать в
 % резолюции ни одного момента, который дал бы ему возможность продолжать разговоры
  %в таком же духе. Например, в отношении прямого плагиата — он будет опять говорить: я
  %крупный ученый, зачем мне плагиат? Но мы должны подчеркнуть, что он обкрадывает
  %учеников. {\bf Надо сказать эту мысль его же словами «делает перенос»}, это очень тонко.
  
  %Мне вообще казалось бы, что эту сторону разрыва с учениками, которой мы касаемся, не мешало бы подразвить, потому что это уже констатировано, и самим им, в конце
  %концов, признано. Вспомните сегодняшнюю часть его выступления, где он говорит о
  %Лаврентьеве и Новикове. Это документ очень важный.
\end{quotation}

Решение основного вопроса становится благополучным для Лузина. Происходит некоторое смягчение и в других отношениях (непонятно, по инициативе ли партийной инстанции,
самого Кржижановского или академических верхов).

С.~С.~Демидов и В.~А.~Есаков \cite{DeEs} предлагают такое объяеснение:
  \begin{quotation}
  	Как явствует из реплик Г.~М.~Кржижановского, он имел разговор «наверху», судя по
  	всему с самим И.~В.~Сталиным, к которому единственным из Академии он имел прямой
  	допуск.
  	
  	Таким образом, и это главное (!), наверху (И.~В.~Сталиным) был поддержан мягкий
  	вариант резолюции.
  \end{quotation}
  
  {\small В книге Грэхема и Кантора \cite{GrK}, обсуждавшейся выше в \S \ref{s:filosofiya} эти сентенции трансформируются
  таким образом:
  \begin{quotation}
  	У главы комиссии Глеба Кржижановского возникли сомнения по поводу жуткого исхода, к которому двигались заседания. Он был явно смущен и даже напуган травлей, которую
  	предприняли коллеги Лузина, готовые ради своей карьеры избавиться от него.
   Кржижановский был старым революционером, он знал Сталина лично и решил лично поговорить с ним о деле Лузина. 
   
   {\bf Беседа Кржижановского и Сталина состоялась 11 или 12 июля 1936 года. Никакого свидетельства о ней не найдено, но ясно, что оба пришли к выводу}, что коллеги Лузина зашли
   слишком далеко. Сталин с самого начала не испытывал восторга от расследования дела Лузина, потому что не слишком высоко ставил Кольмана....   
\end{quotation}
     
  Последнее -- скорее отрывок из романа, чем из  исторической работы.  
  Никаких известий о беседе Сталина с Кржижановским нам не известно. Их отношение
  к Лузину тоже не известно (Кржижановский был вполне ж\"есток на заседаниях комиссии).
  Не известно, как высоко Сталин ставил Кольмана, и знал ли вообще Сталин о его существовании%
  \footnote{Кольман в \cite{Kolman-vospominaniya} 
  	рассказывает о своих встречах со Сталиным. Даже, если не брать
  в расчет сомнительную достоверность его повествования
  и того, что Кольман  сильно преувеличивает степень своих контактов с высочайшими персонами,
  последний вопрос возникает при чтении рассказа Кольмана  естественно.
  Автор не может удержаться от того, чтобы лишний раз процитировать. Среди бумаг, собранных Дюгаком
  и частично опубликованных в \cite{Dug}, есть и такое:
    \newline
  {\it К.~Каплан (K.~Kaplan)
по этому поводу пишет:
«Именно во время своего пребывания в Крыму Готвальд получил рапорт о резкой
критике Кольмана.
Сталин, который уже был в курсе этого по своим собственным
каналам, хотел все это выяснить и вызвал Готвальда для объяснений. Этот последний
прибыл к нему с большой опаской — в его окружении даже говорили, 
что он рисковал больше не вернуться обратно. Сталин хотел прежде всего узнать, кто такой
этот Кольман. Узнав, что речь идет не просто о полемисте, боровшемся в тридцатых
годах в СССР против буржуазных тенденций в естественных науках, но что он также
проходил как троцкист у Берии (присутствовавшим на встрече), Сталин решил:
"Поскольку это советский гражданин (что на самом деле было не так), вам надо 
только отправить его сюда, мы приведем его в чувство", — пообещал он».}
\newline
Автор не может проанализировать степень достоверности этого известия.}.}
  
  До нас дошел журнал посещений Сталина \cite{Stalin}%
  \footnote{Многие широко признаваемые за истину байки
  	рассыпаются при открывании этой книги. Например, есть  литература, о том, как
  	Н.~И.~Вавилов ходил в этот кабинет, и о том, как он там был в последний раз, например, \cite{Lebedev-Kolchinsky}. Между тем Вавилов ни разу
  	в этом кабинете не был.}.
  За  1924-1953гг Кржижановский был на приеме в кабинете Сталина лишь 6 раз,
  все в 1928-31гг. Конечно, Кржижановский со Сталиным мог бы, например,
  говорить по телефону, но речь шла о редактировании довольно длинной
  бумаги. Так что наличие разговора вызывает сомнение и с этой стороны.

  В той же статье \cite{DeEs} приведен следующий документ, который, по-моему, содержит 
  отрицательный ответ на заданный вопрос:
  \begin{quotation}
  	13 июля 1936
  	
  	<на бланке редакции газеты «Правда»>
  	
  	Товарищам СТАЛИНУ
  	
  	МОЛОТОВУ
  	
  	На состоявшемся 13 июля 1936 г. заседании комиссии Академии наук по делу Лузина
  	председатель комиссии тов. Г.~М.~Кржижановский нашел нужным сделать несколько замечаний по окончательной редакции резолюции.
  	
  	Ссылаясь весьма определенно, невзирая на присутствие беспартийных членов комиссии и беспартийных научных работников из математических кругов, тов. Г.~М.~Кржижановский предложил несколько отредактировать текст резолюции, в основном, по его словам, одобренный в соответствующих инстанциях.
  	
  	В частности, тов. Г.~М.~Кржижановский предложил изменить последний пункт резолюции в том смысле, чтобы не называть Лузина врагом в советской маске, как это сделано в «Правде». Тов. Кржижановский решил выразить свое отношение к линии поведения
  	Лузина «по-академически», сформулировав это поведение недостойным звания советского ученого.
  	
  	Прилагаю соответствующую выдержку из стенограммы выступления Г.~М.~Кржижановского, полученную в Академии наук.
  	
  	Редактор «Правды» Л.~Мехлис
  \end{quotation}
  
  Документ подтверждает, что ситуация была действительно опасной, и что Мехлис,
  в полном соответствии с его сложившейся впоследствии репутацией, жаждал крови.
  Но из текста ясно, что <<инстанции>> -- это не Сталин, вряд ли Мехлис стал бы
  обвинять Сталина перед Сталиным. Куда ходил Кржижановский с бумагами, Мехлис
  знал,  а мы не знаем. Скорее всего, <<инстанции>> -- это отдел науки ЦК
  (зав. К.~Я.~Бауман%
  \footnote{Бауман на приеме у Сталина в 1936 году тоже не был.  
  Мехлис  был 18.06, 1.07, 20.07, 11.12, 13.12. Можно лишь гадать,
  	задавался ли в течение 15-минутного пребывания Мехлиса в кабинете Сталина
  	 20.07 вопрос о Лузине.}%
  ), но, быть может, не только он. <<Инстанция>> (быть может, с участием Кржижановского) была не столь воинственна как Мехлис, а кроме того
  отменила формулировку из <<Правды>>, что уязвляло Мехлиса как главного редактора <<Правды>>.
  
  Так или иначе, Мехлис обвиняет перед Сталиным Кржижановского и <<инстанции>>.
  Реакция Сталина нам неизвестна, мы лишь видим, что благоприятное для Лузина решение
  осталось в силе. У математиков была легенда (версия есть у Ефремовича, см. п. \ref{ss:efremovich}), что Сталин решил остановить дело.
  Цитируем С.~П.~Новикова-младшего \cite{NovS}:
  \begin{quotation}
  	Это дошло до Сталина (кажется, через Капицу, который якобы писал Молотову, но читал Сталин). Он удивился: этого дела в академическом истеблишменте он пока не планировал. Говорят Сталин сказал: <<Газета Правда, товарищи не ошибается. Но здесь же сказано, что вредительство своеобразное. Значит и наказание должно быть своеобразное. Ограничиться обсуждением.>>
  \end{quotation}

    Вообще-то обвинения Лузина должны были выглядеть сверху
  как Батрахомиомахия.   
  
  Постановление Президиума АН СССР от 5 августа было довольно сдержанным,
  а завершалось так 
\begin{quotation}
	Президиум Академии Наук полагает, что поведение акад. Н.~Н.~Лузина несовместимо
	с достоинством действительного члена Академии Наук, и что наша научная 
	общественность имеет все основания ставить вопрос об исключении его из состава академиков.
	
	Однако, учитывая значение H.~H.~Лузина как крупного математика, взвешивая всю
	силу общественного воздействия, выявившегося в столь широком, единодушном и справедливом осуждении поведения H.~H.~Лузина, и исходя из желания предоставить Лузину
	возможность перестроить все его дальнейшее поведение и его работу, -
	Президиум считает возможным ограничиться предупреждением H.~H.~Лузина, что
	при отсутствии решительного перелома в его дальнейшем поведении Президиум вынужден будет неотложно поставить вопрос об исключении H.~H.~Лузина из академических рядов».
	\end{quotation}
    
    <<Правда>> 6 августа пересказала это своими словами:
    \begin{quotation}
    Это меньше всего были академические споры. Это вылилось в движение гневного
    протеста против затесавшегося в научную среду классового врага, пытавшегося привить
    советской научной общественности чуждые ей, заимствованные из буржуазного мира,
    грязные приемы и манеры. Дело Лузина дало возможность очень широким научным
    кругам почувствовать с особой остротой, в чем заключается достоинство советской науки. В страстном негодовании, которое объединило и старых академиков и советскую
    научную молодежь, было гордое сознание силы, свободы и независимости науки в советской стране. В своем постановлении президиум Академии наук признал, что научная советская общественность имеет основание требовать исключения Лузина из академических рядов. {\bf Презренный раб буржуазии}, подобострастно прикладывавшийся к ручке всякого иностранного «генерала» от науки и пренебрежительно относившийся к советской
    науке, упразднившей генеральство, — не достоин занимать место среди академиков.	
    
    Комиссия Академии наук полностью подтвердила ту характеристику, которая была
    дана Лузину в статьях «Правды». Это — враг в советской маске. Он принадлежал к правым кругам профессуры до революции. Советская власть простила ему прошлое. Она
    предоставила ему все возможности для научной работы. Лузину не приходилось испытывать материальной нужды, как его собратьям по науке за границей. Ему не грозило изгнание. Он пользовался полной свободой в области своей научной деятельности.
    Но он остался верен реакционной буржуазии. {\bf Свое высокое положение академика
    он использовал для вредительской подпольной работы. Он старался засорить ряды ученых, проталкивая заведомо невежественных и бездарных людей.} Циничным издевательством были его подхалимские высказывания о средней школе, напечатанные в «Известиях». Как торгаш по натуре, он обкрадывал своих учеников. На науку он смотрел как на
    базар, и преследованием молодых талантливых ученых он оберегал свое монопольное
    место.\dots
    
    Но пролетарская революция добралась и до этого места. Овладевая всеми высотами
    научными знаниями, советская научная молодежь взяла в свои руки и высшую математику. В борьбе с талантливыми учениками Лузин обнаружил свое внутреннее бессилие.
    Он боролся самыми грязными средствами, в буквальном смысле слова сживая со света
    опасных конкурентов. Но он был разоблачен своими же учениками. Беспартийная советская научная молодежь помогла сорвать маску с классового врага. Нет и не может быть
    буржуазных гнезд в советской науке!
    	\end{quotation}
    	
    	На этом шквал стих. Худшее было позади (хотя в тот момент это могло быть не вполне ясным).
    	
    	\sm

 {\bf \punct Публичная кампания.%
 \label{ss:public}}  Статьи в Правде 
часто имели инструктивный характер: <<Делай как я!>>. 
Те, кого такие публичные инструкции могли касаться,
должны были сами это понимать. Если статья имела продолжение
(а 3 июля была уже вторая статья), то сомнений уже ни у кого не должно было
остаться. По научным организациям страны прошли собрания и митинги, посвященные обличению
Лузина и лузинщины, а также поиску и искоренению лузинщины в своей среде
(что естественным образом становилось инструментом местных разборок).
Несколько таких статей из тогдашней печати приведено в \cite{Delo}. Тема исследовалась также в 
\cite{AlexandrovDA}. Как будто, никаких серьезных
последствий это не имело.
В Харькове в то время издавался на немецком языке физический журнал 
 «Physikalische Zeitschrift der Sowjet Union».
Д.~А.~Александров \cite{AlexandrovDA}  сообщает, что местный обком в порядке
борьбы с лузинщиной постановил выпускать этот журнал на украинском языке.
Москва это постановление отменила.

Параллельно с делом Лузина происходила свара в Пулковской обсерватории,
которая по мнению того же Д.~А.~Александрова \cite{AlexandrovDA}, который ее разбирал, была 
обусловлена внутренними причинами (со словами {\it достаточно очевидно}). Интересно, что директор обсерватории
Б.~П.~Герасимович, тоже обвинявшийся в низкопоклонстве перед Западом, был назначен 
членом академической комиссии по борьбе с лузинщиной в физике.

Через месяц после того, как баталии в Москве (и не только в Москве) уже отгремели, кампания против лузинщины
 началась
в Томске. Цитируем статью \cite{luzin-v-sibiri}, где она подробно исследовалась:
\begin{quotation}
Сигналом к началу кампании послужил приезд в Томск 5 сентября 1936 г. секретаря крайкома Р.~И.~Эйхе. Выступая 11 сентября на собрании партактива Томской организации, он указал на пассивность коммунистов и научной общественности при обсуждении «дела Лузина» как на серьезнейший «провал в работе». «В Томске, где насчитывается около 800 научных работников, почти никто на дело Лузина, на статьи в <<Правде>> никак не отозвался, словно дело Лузина никакого политического значения не имеет. Неужели вы полагаете, что в Томске нет лузинщины, нет отдельных проявлений раболепия перед буржуазной наукой?»
\end{quotation}
   
   Товарищ Эйхе был кандидатом в члены Политбюро ЦК ВКП(б) и сыграл какую-то
   особую роль в наступлении 1937 года.
   В Томске (где был хороший университет) 15-19 сентября 1936г. прошла жесткая кампания по борьбе с лузинщиной с горячими собраниями и с обороной
   обвиняемых. Как будто, прямые последствия этой истории были не велики
   (случайно или нет, Стефан Бергман вскоре оказался в Тбилиси, а П.~С.~Тартаковский в Ленинграде).
   Звучало обвинение, что <<Известия НИИ математики и механики Томского университета>>
   издавались по-немецки. Но они как издавались на русском, немецком 
   и французском языке до 1936г., так и продолжили издаваться в 1937-38г.
  
    По мнению авторов \cite{luzin-v-sibiri}  кампанию
   загасил сам товарищ Эйхе. Закончим цитирование:
   \begin{quotation}
Эти положения вкратце таковы: наряду с «правильными» выступлениями имелись «неправильные» (те, кто поначалу упорствовал и отрицал элементы «лузинщины» в своей работе – Тартаковский, Иваненко и др.); ряд ученых проделали эволюцию, и если на внутривузовских собраниях они пытались «отрицать ошибки», то на общегородском собрании большинство из них «правильно признавали свои ошибки» (хотя и с оговорками); в то же время двое – Иваненко и Фукс%
\footnote{У Б.~А.~Фукса была одна заметка в Compt. Rendus, 
то есть он вроде и не подпадал под обвинение газеты <<Правда>>.} – «не выступили с признанием свои ошибок».
\end{quotation}

   \section{После шквала%
   \label{s:after}}

      \COUNTERS
   
   Так или иначе, история 1936г. могла кончиться  очень 
   плохо, но кончилась она, в целом, благополучно. Лузин был свергнут
   с его должности, и это было положительным явлением, он остался в 
   Академии, и это тоже было хорошо. В 1936г. главный математик СССР
   не был избран, а во главе Математической группы встал Президиум в составе
   Бернштейна, Виноградова и Сегала (то есть два знаменитых математика-академика
   и один комиссар). В 1939г. должность председателя
   занял Колмогоров.
   
   Другие последствия, приписываемые этой истории, сильно преувеличены.

 \sm

 {\bf \punct Журналы, новые и старые.%
 \label{ss:journals}}  
  В течение 30х происходило постепенно улучшается издание советских математических
  журналов. До 1930г. был один центральный математический журнал, <<{\it Математический сборник}>>.
   C 1932г. начинается увеличение его объема. 
   $$
   \begin{matrix}
                           &1922&1924&1926&1927&1928&1929 \\
 \text{Количество выпусков}&   3   & 4  & 4  & 4  & 5  & 3  \\
 \text{Количество страниц} & 151 & 272& 412& 383&483 & 424\\
   \text{Количество статей}&  10  & 24 & 26 & 13 & 29 & 32
   \end{matrix}
   $$
   
   $$
   \begin{matrix}
   %Год
   %1900—1918
   %(в среднем
   %за год)
   1930&
   1931&
   1932&
   1933&
   1934&
   1935&
   1936&
   1937&
   1938&
   1939&
   1940\\
   %\text{Количество выпусков}&
   %3
   2&
   2&
   3&
   4&
   4&
   6&
   6&
   6&
   6&
   6&
   6\\
 %\text{Количество страниц}&
   %305
   212&
   240&
   450&
   528&
   668&
   745&
   1000&
   1253&
   %677+570
   1247&
  % 660+541
  1201&
 % 597+519
 1118
   \\
  % \text{Количество статей}&
   %13
   12&
   15&
   29&
   37&
   44&
   54&
   80&
   64&
 %  51+36
   87&
   %40+30
   70&
  % 30+34
  64
    \\
   \end{matrix}
   $$
  Еще Егоров  начал превращать <<Сборник>> в международный журнал, он 
  стал им в 30е годы и сохранял эту роль 
  до начала Второй мировой войны. В 1946г. Люстерник \cite{Lyusternik-1946}
  писал:
  \begin{quotation}
  	К началу Отечественной войны «Сборник» и по качеству материала и по
  	объёму был одним из лучших международных математических журналов. Возросло
  	сотрудничество крупных иностранных учёных в «Сборнике». 
  \end{quotation}
  В 1941 - 1947гг в журнале еще изредка печатались
  иностранные авторы, среди них P.~R.~Halmos, J.~E.~Littlewood, L.~J.~Mordell, W.~J.~Trjitzinsky.
  До 1945г. в журнале было  много статей советских авторов на иностранных языках,
  публикация их прекращается к 1948г.
  
  \sm
  
  В 1933г. оживляются <<{\it Доклады Академии наук}>>. Вот данные по числу статей, видимых Zentralblatt Math.
  $$
  \begin{matrix}
  	1933 & 1934 & 1935 &1936 &1937 &1938 & 1939&1940\\
  	 12  &  101 & 97   & 82 &  104 & 132 & 151 &143
  \end{matrix}
  $$
  Журнал публиковал короткие сообщения, до Войны значительная часть публикаций
  выходила на иностранных языках.

 В 1937г. начал выходить следующий новый журнал <<{\it Известия Академии наук СССР.  Серия математическая}>>.
 До этого было общее издание по всем естественным наукам: <<Известия Академии наук СССР. VII серия. Отделение математических и естественных наук>>. 
 
 В том же 1936г. появился журнал <<{\it Прикладная математика и механика}>>.
 
 Наряду с регулярными журналами появилось три центральных сериала.
 С 1931г. начали выходить (нерегулярно) <<{\it Труды института Стеклова}>>, всего
 до 1940г. вышло 14 выпусков (часть из них очень тонких, 10 стр.).
  
    С 1933г. началось издание сериала <<{\it Труды семинара по векторному и тензорному анализу}>>, в 1941г вышел
  5ый выпуск.
  
  В 1936г. начал выходить сериал <<{\it Успехи математических наук}>>. Всего до Войны
  вышло 9 выпусков. В регулярный журнал издание превратилось в 1946г. 
  
  Все эти издания возобновили выпуск после Войны. После Войны оборвался замечательный
  сериал <<{\it Известия НИИ математики и механики Томского университета}>>, основанный 
  в 1935г. Чуть ниже мы скажем о нем несколько слов.
  
 \sm

В 30е годы в провинции выходил ряд хороших математических изданий,
<<{\it Ученые записки Казанского государственного университета}>>, <<{\it Бюллетень
Средне-Азиатского государственного университета}>> (издававшийся
В.~И.~Романовским%
\footnote{Романовский Всеволод Иванович (1879—1954), яркий деятель и в научном и организационном плане, основные труды по статистике
и теории вероятностей (но не только, например, он обнаружил
три новых системы гипергеометрических ортогональных многочленов
\cite{Romanovsky}), один их основателей Ташкентского (Средне-Азиатского) университета.}, там печатались также иностранные авторы),
различные издания в Киеве.
По-моему, тогдашние провинциальные издания тянулись  вверх и 
и не были выражением своеобразного духа местных научных тусовок как в послевоенные десятилетия
Советского Союза
(когда исключением из этого правила был лишь <<Сибирский математический журнал>>).

\sm
 
 {\small
 	Теперь о Томских известиях.
 Из <<Математической хроники>> в первом томе сериала Успехи мат.наук (1936), \cite{Uspehi-36-1}:
 \begin{quotation}
  В 1935 г. вышли первые два выпуска%
  \footnote{Zentralblatt показывает второй выпуск в 1937 году.} <<Известий научно-исследовательского института 
  математики и механики при Томском гос. университете им. В. В. Куйбышева>>
(отв. редактор Ф. Молин, редколлегия — С. Бергман, Ф. Молин. Ф. Нётер). 
 \end{quotation}
 
 Молин уже упоминался выше, Стефан Бергман и Фриц Нётер - политэмигранты из Германии.
 В 1938 году вышел второй том <<Известий>> (два выпуска).
   Среди его авторов были, например,  Эйнштейн~А. [Albert Einstein], Колмогоров~А.~Н., Хинчин~А.~Я.,  фон Нейман~Дж. [von Neumann J.],
   Привалов~И., Куфарев~П.~П., Райков~Д.~А.,  Фукс~Б.~А., Нётер~Ф.[Noether F.],  Серпинский~В., 
   Бергман~С.[Bergman S.], Молин~Ф.~Э., Эрдеш~П. [Erd\"os P], Туран~П [Tur\'an P].
   
   Западная Сибирь  была своего рода эпицентром Тридцать Седьмого Года, и Томскому университету тогда досталось.
   В частности, Нётер был арестован в ноябре 37 года, находился в Орловском политизоляторе 
   и вместе с другими заключенными был расстрелян в сентябре 1941
   при приближении немцев к Орлу\footnote{Впрочем, есть саратовская легенда (источник -- Б.~Шайн), выглядящая 
   самоправдоподобно, 
   что
   саратовский механик С.~В.~Фалькович встретил Нётера в Московском метро в декабре 1941г, я не берусь анализировать
   эти известия. Я, кстати, тоже эту легенду слышал, только встреча случилась в ЦАГИ.}. Бергман в 1936 году перебрался в Тбилиси, и его судьба сложилось благополучно.
   
   После успокоения тамошней жизни Молин попытался продолжить издание сериала и в августе 1940 собрал следующий том.
   Он был издан лишь после Войны в 1946 году \cite{Uspehi-Tomsk}. На этом издание оборвалось.}

В тогдашнем математическом сообществе многократно озвучивалась идея издания популярного журнала.
В 1934-38гг. вышли 13 выпусков сериала <<{\it Математическое просвещение}>>%
\footnote{Издатель Бончковский Ростислав Николаевич, (1905–1942), погиб под Сталинградом.}.
В 1938 А.~Ф.~Бермант
  \cite{Bermant} писал
\begin{quotation}
	Наконец, с особенной настойчивостью следует пожелать, чтобы наступил
	конец мытарствам элементарного журнала <<Математическое просвещение>>, до сих
	пор существующего в виде непериодических сборников. 
\end{quotation}
 Однако сериал  закрылся в том же году, новое просветительское
  издание%
  \footnote{В 1957-1961 под тем же названием <<Математическое просвещение>> выходил издававшийся Я.~С.~Дубновым, А.~А.~Ляпуновым, и А.~И.~Маркушевичем,
  сериал, целью которого была подготовка реформы школьного математического образования. Под тем же названием с 1997 выходит популярный сериал.}$^,$
\footnote{C 1937 выходит специализированный педагогический журнал <<Математика в школе>>.} -- журнал <<Квант>> для школьников был основан Колмогоровым и И.~К.~Кикоиным лишь в  1970г.
  
  Так или иначе, не исключено (но, может, это все происходило само собой), что начало издания журнала <<Известия. Серия математическая>>
  и некоторое увеличение емкости <<Докладов>> после 1936г. были (или были отчасти)
  последствием дела Лузина.
  
  \sm
   
   {\bf \punct Заграничные публикации.%
   \label{ss:zagranitsa}}   С делом Лузина часто связывают 
   прекращение публикаций советских ученых за рубежом. Обсудим, насколько это
   соответствует истине.
   
    На академической комиссии
   таких призывов  не было. В буквальном мысле слова их не было и в печати. 
 \begin{quotation}
 	КРЖИЖАНОВСКИЙ. Может быть, там более изящные издания, но не в этом причина.
 {\bf 	Никто, конечно, не может сказать, что мы, говоря о том, чтобы печатали в нашей стране,
 	этим самым хотим оторвать научные работы от общей мировой продукции.} Но мы констатируем, что неслучайно [у Лузина] есть большая разница в работах, печатаемых за границей и
 	печатаемых здесь.
 \end{quotation}   
 
 Есть известная статья Д.~А.~Александрова <<Почему советские учёные перестали печататься за рубежом: становление самодостаточности и изолированности отечественной науки>>
 \cite{AlexandrovDA}:
 \begin{quotation}
 Казалось бы, такая кампания, прокатившаяся по стране [обличение лузинщины], и определила переломный момент, когда власти изолировали научное сообщество от Запада. Легко заключить, что ученые, под давлением сверху, должны были перестать печататься в зарубежных изданиях и общаться со своими западными коллегами. Однако это, казалось бы, самоочевидное утверждение не подтверждается ни количественными библиографическими данными, ни архивными свидетельствами...

 Спад доли зарубежных публикаций в 1936–1937 гг. оказался выражен разве лишь для сообщества математиков; впрочем, тенденция перехода на публикации в отечественных изданиях заметна у математиков и в более ранний период. При этом математики и после 1936 г. продолжали публиковаться за рубежом, хотя и менее интенсивно.
 \end{quotation}
 
 Автор данных записок провел поиск%
  \footnote{Первичный список публикаций лежит по адресу
  	\newline
  	http://www.mat.univie.ac.at/$\sim$neretin/misc/luzin/list.txt
  	\newline
  	обработка в
\newline 	
  	http://www.mat.univie.ac.at/$\sim$neretin/misc/luzin/schet.txt
  	\newline
  	http://www.mat.univie.ac.at/$\sim$neretin/misc/luzin/luzin-publications.html
  	\newline
  	При составлении списков есть довольно много проблем, с обилием транскрипций и с определением
  	гражданства (например, Дюбюк -- московский математик, а из двух Ивановых один оказался из Сан-Франциско, а другой предположительно был болгарин). Скорее всего, в списке
  	остались лакуны.} зарубежных публикаций советских математиков за 1930-1943 годы с помощью различных современных электронных баз
 (Zentralblatt, MathSciNet, MathNet) и библиографического справочника \cite{za30}.
 Политимигранты в статистику не включались.

   Во-первых, вот список математиков, которые публиковались за границей после дела Лузина. На всякий случай мы считаем лишь публикации, {\bf начиная с 38 года по 43 год, плюс Comptes Rendus 37 года}, для осторожности мы не учитываем 
    журнальные публикации
   37 года (статьи могли быть поданы в 36 году), кроме быстро публиковавшегося Comptes Rendus. 
   
   \sm
   
   П.~С.~Александров, Н.~А.~Артемьев, А.~Ф.~Бермант,  С.~В.~Бахвалов, С.~Н.~Бернштейн, Н.~Н.~Боголюбов, Б.~В.~Булгаков,   Н.~Б.~Веденисов, Б.~З.~Вулих,
    Я.~Л.~Геронимус, Н.~А.~Глаголев, В.~И.~Гливенко, Н.~М.~Гюнтер, Б.~В.~Гнеденко, А.~М.~Данилевский, А.~П.~Дицман, П.~Е.~Дюбюк, Л.~В.~Канторович, М.~В.~Келдыш, А.~Н.~Колмогоров, А.~А.~Конюс, М.~Г.~Крейн, Н.~М.~Крылов, А.~А.~Кулаков, В.~Д.~Купрадзе, М.~А.~Ларентьев, Б.~М.~Левитан, Е.~М.~Ливенсон, И.~М.~Максимов, А.~А.~Марков, И.~Е.~Огиевецкий, И.~Г.~Петровский, А.~Г.~Пинскер, Л.~С.~Понтрягин, И.~И.~Привалов, В.~И.~Романовский, Е.~Е.~Слуцкий, А.~Н.~Тулайков, В.~К.~Туркин, С.~П.~Фиников, С.~В.~Фомин, Г.~Ф.~Хильми, Н.~Г.~Чеботарев, В.~Л.~Шмульян 
    
   \sm
   
   На запрет это не очень похоже.
   
     \sm
   
   Несколько  активно публиковавшихся за границей математиков этого больше не делали: И.~М.~Виноградов, Н.~Н.~Лузин (разумеется), Д.~Е.~Меньшов (не очень активный), П.~К.~Рашевский, А.~Я.~Хинчин (две публикации 37 года), Г.~В.~Пфейфер, С.~Д.~Россинский, С.~П.~Слугинов, М.~К.~Куренский (одна публикация 37 года). Зато в вышеприведенном списке много людей, ранее не публиковавшихся за границей. 
   
   Вот даты заграничных публикаций персонально  обличенных <<Правдой>> товарищей 
   
   \sm
   
   Александров 1937, 1937, 1938, 1941, 1943 
   
   Бернштейн 1937, 1937, 1937,1937,1938, 1938, 1938, 1942 
   
   Колмогоров 1937, 1937, 1938, 1941 
   
   Хинчин 1937, 1937
   
   \sm

  Общая статистика заграничных публикаций советских математиков  по годам
  $$
  \begin{array}{cccccccc}
1930 &1931&1932&1933&1934&1935&1936
\\
   96 (81)&
   78 (75)&
   88 (73)&
   76&
   57&
   83&
   89
   \end{array}
  $$
  
  $$ 
   \begin{array}{ccccccc}
  1937& 1938&1939&1940&1941&1942&1943\\
      53&
   25&
   11&
   14&
   8&
   1&
   7
   \end{array}
   $$
   
   В скобках стоит число публикаций 1930-1932 без трудов Конгрессов в Болонье и Цюрихе.
   Мы видим, что в 1930--1936гг общее число заграничных публикаций советских математиков
   за границей держалось примерно на одном и том же уровне%
   \footnote{Емкость  советских журналов росла, поэтому относительная доля иностранных
   	публикаций должна была уменьшаться.}
    (яма 1934 года, возможно, объясняется кризисом в 1933 года в СССР). В 1937 году происходит спад и дальше количество публикаций убывает.
   
   Стоит иметь в виду, что для этого убывания были более веские причины, чем <<дело Лузина>>.
   
   Во-первых, в сентябре 1939 года началась Мировая война, а с апреля 39 года Европа к этой войне готовилась. СССР в 1938 году пережил две военных тревоги, в августе  на востоке и в сентябре  на западе.
   
   Во-вторых, с лета 37 по осень 38 в СССР шла весьма заметная охота на шпионов. Понятно, что многие граждане могли опасаться заграничных публикаций.
   
   В-третьих, с приходом к власти фашистов в Германии, советские авторы естественным образом
   стали там реже печататься. Резко отрицательное отношение к публикациям в фашистских странах
   советские власти недвусмысленно высказали в 1936 во время <<дела Лузина>>. В связи с гражданской 
   войной в Испании, начавшейся в  июле 1936 года, резко ухудшились советско-итальянские
   отношения, и еще одна важная с точки зрения математики страна, вслед за Германией, закрылась для публикаций.
   
   Однако был еще один фактор, сработавший в 37 году, который сразу не приходит в голову.
   Он становится виден, если составить статистику по разным странам.

 \begin{longtable}{l|lll| ll|lll}
 	\multicolumn{9}{c}{Публикации за рубежом 1930--1943 (по странам)}
 	\\
 	\hline
 	                 & Герма & Авст  & Ита  & Франция & Поль & США & Гол     & Другие \\
 	                 & ния   & рия   & лия  & C.R+пр. & ша   &     & ланд    &  \\ \hline
 	\endfirsthead    & Герм. & Авст. & Итал & Франция & Пол. & США & Голланд & Другие \\ \hline
 	\endhead   
    1930             & 8     & 2     & 13   & 33+8    & 5    & 2   &         & 9      \\
 	1931             & 20    &       & 12   & 24+5    & 4    & 3   &         & 4      \\
 	1932             & 13    &       & 18   & 9+5     & 10   & 4   &         & 5      \\
 	1933             & 17    & 1     & 11   & 24+7    & 7    & 1   &         & 6      \\
 	1934             & 10    & 1     & 6    & 24+1    & 6    & 3   & 1       & 5      \\
 	1935             & 11    &       & 8    & 42+2    & 6    & 6   & 3       & 4      \\
 	1936             & 4     &       & 10   & 30+8    & 11   & 10  & 4       & 10     \\ \hline
 	1937             & 3     &       & 2    & 16+3    & 6    & 9   & 7       & 6      \\
 	1938             &       &       &      & 6+10    & 1    & 1   & 4       & 3      \\
 	1939             &       &       &      & 2+2     &      & 3   & 2       & 3      \\
 	1940             & 1     &       &      & 0+2     &      & 5   & 1       & 4      \\ \hline
 	1941             &       &       &      &         &      & 3   &         & 4      \\
 	1942             &       &       &      &         &      & 1   &         &  \\
 	1943             &       &       & 1    &         &      & 6   &         &
 \end{longtable}  

По Франции отдельно записаны публикации в Comptes Rendus и в других журналах.
Столбец <<Другие>> -- это Швейцария, Швеция, Бельгия, Великобритания, Румыния, Чехословакия,
Латинская Америка, Япония, Индия,  Турция (вклад каждой из этих стран по-отдельности невелик)%
\footnote{Странная единичка в Италии за 1943 год -- это упоминавшийся в сноске \ref{fo:maximov}
И.~М.~Максимов, он же трижды публиковался в Tohoku (Япония) в 1940-41гг
 (вообще же   публикации
других советских математиков в
Японии кончились к 1937 году).}.

Мы видим, что в 1936 году на Францию+Польшу+Италию пришлось $38+11+10=59$ публикаций,
а на все остальные страны 28 публикаций. С фашистской Италией в 1936 году вконец испортились отношения, а во Франции и Польше, как мы видели, атака советских математиков на Лузина
вызвала возмущение   (напомним слова Серпинского%
\footnote{В 1948 году советская делегация должна была поехать на VI
	съезд польских математиков. Оказалось, что там должно было быть, в частности,
	чествование Серпинского. Советская делегация (Колмогоров, Александров
и секретарь парткома Стекловки	К.~К~ Марджанишвили) не поехала, ей запретили.
Запретило высокое партийное руководство. По просьбе К.~К.~Марджанишвили \cite{DeEs} или \cite{Politburo}.	
	И с чего бы это? И очень  ли хотели два других члена предполагаемой
	делегации чествовать Серпинского?}%
: <<{\it присутствие господ Александрова, Хинчина,
Колмогорова, Шнирельмана, которые самым нечестным образом
выступили против своего бывшего учителя и ложно обвинили
его, — нельзя терпеть ни в каком собрании честных людей}.>>). Соответственно,
по Франции и Польше в 1937 году произошло резкое сокращение числа публикаций
(возможно, по Comptes Rendus оно началось уже в 1936 году).

Интересно, что по Франции в 1938 году произошел откат, в обычных
журналах (не Comptes Rendus) было опубликовано 10 советских статей вместо 3 в 1937.
Однако экологическая ниша Comptes Rendus для советских математиков уже была занята
Докладами АН СССР.

Так или иначе, прекращение публикаций советских математиков за рубежом
произошло после Отечественной Войны и было связано с Холодной войной
и Железным занавесом, оно оставалось сильно затрудненным вплоть до Перестройки.
А сокращение числа публикаций в связи с <<делом Лузина>> в 1937
в самом деле было, но в основном, в связи с возмущением во Франции и Польше
(о чем многим советским математикам потом вспоминать было, по-видимому, не
вполне приятно), и, отчасти, с опасениями части математиков, которые прекратили
публикации сами.
   
   \sm
   
 \def\oslo{{\bf\punct Математический конгресс в Осло, 1936.%
   \label{ss:oslo}}
   В 20х годах советские математики, по-видимому, свободно ездили 
   за границу, если могли раздобыть средства на эти поездки
   (в связи с этим часто упоминаются <<рокфеллеровские стипендии>>, 
   также 
   упоминается год, когда их советским математикам не давали).
   В книге \cite{Politburo} <<Академия наук в решениях Политбюро>>
   приводится обширный документ от 28.02.1928, в шапке стоят
    тт. Бухарин, Рязанов и Покровский. Далее идет длинная записка
    от Покровского, что делать с Академией. Дальше, согласно
    публикации <<К протоколу фракции был приложен проект постановления Политбюро, в котором предусматривалось:.......>>
    \begin{quotation}
    	ДОПОЛНИТЕЛЬНЫЕ ПРЕДЛОЖЕНИЯ

    	2.	В состав фракции коммунистов-академиков включить т.т. Енукидзе, Криницкого, Литвинова, Луначарского и Горбунова. Считать фракцию комиссией Политбюро по руководству Академией Наук и по регулированию международных научных связей СССР (участие в международных научных съездах и конференциях).
    	
    	Предоставить комиссии право окончательного решения относящихся сюда вопросов о заграничных командировках советских ученых, а также всех командировок по линиии Академии Наук.	
    \end{quotation}

Первые четверо из перечисленных лиц находились весьма высоко в тогдашней
партийной иерархии, и не каждому смертному дано было предлагать
их куда-либо включать. Да и не каждому смертному дано готовить постановление Политбюро. Естественно предполагать, что в этой роли выступал член Политбюро Бухарин.

Что касается научных командировок, то в дальнейшем в книге \cite{Politburo} присутствуют 
многочисленные решения Политбюро на тему выпускать или не выпускать.
Все ли вопросы такого типа решались на уровне Политбюро, и в самом
ли деле Политбюро решало или штамповало чьи-то решения, я из публикации
понять не смог. Но так или иначе, с этого времени заграничные командировки
ученых требовали политического одобрения.

    14-18 июля 1936 в Осло
   проходил международный математический конгресс. В списке membres du Congr\`es
   довольно много советских участников, Гельфонд и Хинчин были приглашенными 
   докладчиками. Однако советские участники на конгрессе, судя по всему, в Осло не появились%
   \footnote{В трудах Конгресса в списке числа участников по странам около СССР стоит 
   цифра 11, но, по-видимому, это все же не так. В \cite{Dug} приводятся
   строчки из письма Г.~Фрейденталя [Hans Freudenthal] Дюгаку (1978) <<{\it Единственным математиком, приехавшим из СССР на Конгресс в Осло, был
   Ф.~Нетер [F.~Noether], немецкий эмигрант и, возможно, обладатель немецкого паспорта...}».}.

Была такая легенда, переданная Е.~Б.~Дынкиным \cite{Yaglom}:
\begin{quotation}
	 Ну, конечно, он [Б. Н. Делоне] любил с молодежью сплетничать. Он даже рассказывал, что
	советских математиков не пустили на конгресс в Осло из-за того, что там был Троцкий.
	Решался на это в 40-х годах...
\end{quotation}

Однако, есть  и такой опубликованный документ \cite{Politburo} Политбюро.
\begin{quotation}
	\footnotesize
	3 июля 1936 г.
	
	О Международном математическом конгрессе в Осло
	
	Отклонить предложение т. Баумана о командировании советских ученых на Международный математический конгресс в Осло.
	
	Протокол № 41, п. 55.	Д. 979. Л. 14.
	
	Принято в связи со следующей запиской К.Я. Баумана:
	
	«СЕКРЕТАРЮ ЦК ВКП(б) тов. И.В.СТАЛИНУ
	О КОМАНДИРОВАНИИ СОВЕТСКИХ УЧЕНЫХ НА МЕЖДУНАРОДНЫЙ
	МАТЕМАТИЧЕСКИЙ КОНГРЕСС В ОСЛО (НОРВЕГИЯ)
	
	С 13-го по 18-е июля с.г. в Осло (Норвегия) состоится очередной Международный математический конгресс. В работе конгресса принимают участие ученые всех крупнейших стран мира (Франция, Англия, Италия, Германия и др.).
	
	На предыдущем Конгрессе в 1932 году (Цюрих) все государства посылали многочисленные делегации (Франция — 69 чел., Германия — 118 чел., США — 66 чел. и т.д.). Советский Союз был представлен делегацией только в составе 3-х математиков%
	\footnote{Мы видели, см. сноску
		\ref{fo:zurich}, что советских участников в Цюрихе было 11. То ли были еще какие-то делегации,
		кроме общей советской делегации (например, делегация Украины), то ли некоторые
		участники поехали туда сами.}. Между тем работы советских математиков занимают одно из первых мест в мировой математической науке.
	
	Академия Наук СССР, Наркомпрос РСФСР, Академия Наук Белорусской ССР обратились с просьбой в ЦК ВКП(б) командировать делегацию на XI Международный математический конгресс в составе 29 чел.
	
	Отдел Науки, научно-технических изобретений и открытий предлагает командировать делегацию в 10 человек в следующем составе: Бернштейн~С.~Н., действительный член Академии Наук СССР, И.~М.~Виноградов — действительный член Академии Наук СССР, Бурстин~Ц.~Л. — действительный член Академии Наук Белорусской ССР, проф. Александров~П.~С., проф. Гельфонд~А.~О., проф. Хинчин~А.~Н., проф. Ковнер~С.~С., проф. Кольман~Э., проф. Колмогоров~А.~Н., проф. Мусхилишвили~Н.~И.
	
	Кроме того Академия Наук просит разрешить советской делегации внести предложение о созыве будущего XII Международного математического конгресса в СССР (Москва — Ленинград) в 1940 г.
	
	Отдел Науки ЦК ВКП(б) поддерживает предложение Академии Наук о созыве очередного математического Конгресса в СССР.
	
	Приложение: Проект постановления и список ученых.
	
	Зав. Отделом науки, научно-технических изобретений и открытий ЦК ВКП(б) Бауман».
	
	На записке резолюция (простым карандашом): «Против Ст[алин]», ниже роспись (синим) — В.Молотов и пометы, что «т. Микоян за предложение т. Сталина, т. Калинин — за».
	
	Выписки посланы: Молотову, Кржижановскому.
	К слову «список» (в приложении) дана сноска: «Список на 3 л. уничтожен. 2.IV.69 [подпись неразборчива]» (РГАСПИ. Ф. 17. On. 163. Д. 1113. Л. 157—159).
\end{quotation}

На бумаге стоит дата 3 июля.
До начала Конгресса оставалось 12 дней. Понятно, что Советская делегация уже должна была
быть утверждена. По-видимому, смысл записки Баумана состоял в том, чтобы эту делегацию увеличить. Однако,  дело Лузина уже разгорелось, и доходило ли оно
до высших эшелонов власти до того, или не доходило, но к этому дню уже дошло
(см. п.\ref{ss:rabolepie}). Похоже на то, что запрет на участие
в Конгрессе в Осло был связан с начавшимся
делом Лузина (а  советские математики в июле 1936г с точки зрения верхов,
видимо,  не должны были отвлекаться от  обсуждения статьи в газете <<Правда>> и не должны были сгоряча обсуждать это с заграничными коллегами).}
   
   \sm
   
 {\bf \punct Второе пришествие Киселёва.%
 \label{ss:kiselev}}  Здесь мы вкратце упомянем еще одно  общественное выступление советских математиков в конце 1936-1937 года (с историей Лузина оно не связано). До нас дошло несколько опубликованных резолюций математических организаций
  с жесткой критикой тогдашних школьных учебников по математике, главным образом, учебника <<Геометрии>> Гурвица%
  \footnote{Гурвиц Юлий Осипович (1882 — 1953) работал в дальнейшем доцентом
  	в ряде московских педагогических вузов, одновременно был учителем в школе,
  	входил в редакцию <<Математики в школе>>, см. некролог
  	в <<Математике
  	в школе>>, 1953, 4, за подписью  И.~К.~Андронова (где есть намеки на обсуждаемую
  	историю).
  	\newline
  	Гангнус Рудольф Вильгельмович (1883-1949).
  	На старом сайте <<Мемориала>> говорится, что он работал а железнодорожном
  	техникуме и в 1937г. получил 5 лет ссылки.
Далее работал учителем в  Каргополе, после отбывания ссылки  работал учителем в
Муроме, умер в Москве. На новом сайте этих данных нет.
	\newline
 Выступление математиков 1936-37гг., насколько
я могу судить по опубликованным текстам, выглядит чисто деловым, без политической
составляющей, хотя присутствуют многочисленные намеки на нечистоту
длительного удерживания этих учебников Наркомпросом.}
 и Гангнуса, см. \cite{Gangnus-1}--\cite{Gangnus-5}, \cite{Uspehi-38-5} и рецензии
  \cite{Gant-Gangnus}, \cite{Lyutin}, \cite{Shnirelman-pravda}, см. также \cite{Fihtengolts}. 
  В опубликованных заявлениях мы видим фамилии: Александров,  Гантмахер, В.~Л.~Гончаров, Люстерник,
  Привалов, Соболев,
  В.~А.~Тартаковский,
  Степанов, Фихтенгольц, С.~А.~Христианович,  Шнирельман, А.~Р.~Эйгес%
  \footnote{Фамилия  героя, вокруг которого собрано наше повествование, не встречается в резолюциях, которые мне пришлось читать. Однако напомню, см. \ref{ss:nachalo-oblavy}, что незадолго до того  он писал в Наркомпрос,
  	<<{\it особенно подчеркивая низкое качество
  	общепринятого тогда учебника геометрии Гангнуса}>>.}

 Высказывания были весьма жесткими (см.\cite{Gangnus-2}, \cite{Uspehi-38-4}):
  \begin{quotation}
  ... Точно так же и в отношении стабильных учебников постановление партии и правительства по существу не выполнено. Смысл введения стабильных учебников заключается, между прочим, в поднятии уровня школьных учебников и изгнании халтурных учебников. Вместо этого Наркомпрос в ряде случаев канонизировал в качестве стабильных безграмотные и халтурные учебники даже по тем предметам, по которым уже имелись на русском языке грамотные учебники (например по геометрии). 
  
  Наркомпросом ничего не сделано для изменения указанного положения вещей, несмотря на многочисленные сигналы как со стороны практических работников школ, так и со стороны отдельных ученых и научных учреждений. 
  
  {\bf Приходится признать, что часть вины все же несет и научная общественность, которая слишком поздно и недостаточно энергично реагировала на положение преподавания математики в средней школе}........ 
  
  9. Группа математики Академии наук СССР надеется, что Наркомпрос сделает все необходимые оргвыводы из сказанного выше и что лица, виновные в грубых ошибках в области руководства начальной и средней школами, будут заменены людьми, способными справиться с огромными задачами, стоящими перед советской школой. 
  
  В частности, необходимо произвести расследование всех обстоятельств, при которых в качестве стабильного учебника по геометрии были утверждены книги Гангнуса и Гурвица, и установить всех лиц, виновных в этом преступлении перед советской школой. 
  \end{quotation}
 
 В частности, 4 марта 1936 года в Правде была опубликована статья Шнирельмана <<Безграмотный учебник геометрии>>. Лев Генрихович с его политической активностью вызывает сложные чувства, но эта статья (к сожалению, у меня нет под рукой ее копии) была образцовым представителем жанра деловой разгромной рецензии, и, скорее всего, сыграла немалую роль
 в деле борьбы за качественные учебники. Приведем риторическую концовку
 \begin{quotation}
 	Быть может, и сейчас найдутся охотники продолжить <<исправление>>
 	негодного учебника Гурвица и Гангнуса... во-первых, нужно изменить
 	порядок следования глав, во-вторых, следует заменить почти все определения, в 
 	третьих заменить несколько сот неверных утверждений верными, в четвертых,
 	следовало бы написать учебник простым ясным языком. Само собой разумеется,
 	что напрашивается необходимость также заменить авторов учебника...
 	
 	Утвердив книги Гангнуса и Гурвица в качестве стабильных учебников,
 	Наркомпрос РСФСР по существу канонизировал безграмотную халтуру...
 	
 	Мало этого, по мнению руководства Управления средней школы Наркомпроса, в этих, с позволения сказать, учебниках все, якобы,
 	построено на движении, они диалектичны. Вряд ли можно придумать  
 	большее издевательство над диалектикой, от которого страдают 
 	и педагоги, и учащиеся...
 	\end{quotation}
 
 Через  дня газета <<Правда>> \cite{Pravda-Shnirelman} 
 сообщила%
 \footnote{Люстерник \cite{Lyusternik-1938} упоминает <<четырехкратные сигналы Правды>>,
 	  я нашел лишь два из них.}, что собрание профкома Мехмата и Институтов механики и математики с участием
 Александрова, Соболева,  Понтрягина, Тихонова, Глаголева, Люстерника, Кагана, Степанова,
 Бари, А.~Кулакова, Колмогорова, Петровского, Гливенко, Немыцкого, Ф.~И.~Франкля, Г.~Э.~Проектора --
 <<всего 65 человек>>
 \begin{quotation}
 	... присоединяется к выводам этой статьи и настаивает на исключении этой книги
 	из числа стабильных учебников.
 \end{quotation}

 Вопросу об учебниках, в частности, были посвящены два заседания Московского математического общества, 3 и 9 апреля 1937г.:
 \begin{quotation}
 РЕЗОЛЮЦИЯ, ПРИНЯТАЯ НА ЗАСЕДАНИИ МОСКОВСКОГО МАТЕМАТИЧЕСКОГО ОБЩЕСТВА 9 АПРЕЛЯ 1937 г.
 
 По вине невежественного руководства со стороны Управления средней школы Наркомпроса,
 в частности по вине А.~И.~Абиндера, учебная литература по математике находится в настоящее время
 на чрезвычайно низком уровне. Управление средней школы Наркомпроса, получая в течение ряда лет и со стороны научных организаций и со стороны учительства сигналы о безграмотности стабильного учебника геометрии Гурвица и Гангнуса, никакой подготовительной работы для замены этой безусловно вредной книги не вело. Книги Гурвица и Гангнуса должны быть изъяты и ни в каком случае не переиздаваемы. {\bf Временно стабильным учебником геометрии должен быть объявлен курс Киселева под редакцией Н.~А.~Глаголева.}
 
 Учитывая, что составление оригинального учебника геометрии потребует времени и что курс Киселева на ближайшее время может удовлетворить потребности преподавания (книга написана безусловно грамотно и усовершенствовалась на протяжении 40 изданий), собрание считает необходимым объявление конкурса с длительным сроком (не менее трех лет) на 1) стабильный учебник по геометрии, 2) курс элементарной геометрии для учителей. Собрание высказывает пожелание о переиздании Учпедгизом различных учебников по всем отделам математики, в частности курсов геометрии Давыдова и Извольского, а также о переводе ряда иностранных учебников. Вместе с тем собрание считает своевременным поставить вопрос о пересмотре учебной программы средней школы по курсу геометрии. 
 
 Председатель И.~И.~Привалов 
 
 Секретарь: В.~Л.~Гончаров 
 \end{quotation}
 
 Напомним, что в конце XIX-начале XX века в России основными школьными учебниками
 математики были легендарные <<Арифметика>>, <<Алгебра>> и <<Геометрия>> Андрея Петровича Киселёва
 (их первые издания вышли соответственно 1884, 1888, 1892, автор многократно их перерабатывал),
 а также <<Тригонометрия>> Николая Александровича Рыбкина. В 1920х годах многократно издавался 
 также учебник <<Алгебра и элементы анализа>> Киселева (<<для трудовой школы>>). 
 После потрясений начала 30х учебники Киселева выпали из школы%
 \footnote{Вопрос этот не изучен. В книге на <<Ленинградском математическом фронте>>\cite{Lenmatfront}
 в качестве
 ругательства упоминается <<киселёвщина>>. В 1934г. Киселё получил Орден Трудового красного знамени.} (в ходу оставалась лишь <<Алгебра>>).
 
 В итоге выступления математической общественности три учебника Киселева в переработанном виде%
 \footnote{ <<Геометрия>> была переработана  Н.~А.~Глаголевым (были добавлены симметрии и проектирования),
 <<Алгебра>> А.~Н.~Барсуковым, а <<Арифметика>>  самим Хинчиным. Сейчас встречаются интернет-филиппики с обличением этих переработок, но учебники в самом деле были уже несколько архаичны, особенно, <<Арифметика>> (см., например, статью Хинчина, \cite{Hinchin-Kiselev}).} были приняты Наркомпросом РСФСР в 1938-1940 как стабильные учебники.

Самому
Киселеву Андрею Петровичу (1852-1940) во время всей этой катавасии было 85 лет, и едва ли он мог в этом принимать участие. 
 
 Математики рассматривали эту меру как временную \footnote{Отчасти это, видимо, было связано с желанием
 	вернуть в школу элементы анализа.
 	Это начали делать  вскоре после Войны, а в 1965--76гг. в старших классах
 	школы основным учебником алгебры был  Е.~С.~Кочетков, Е.~С.~Кочеткова <<Алгебра и элементарные функции>>,
 	простой учебник, включавший, в частности,
 	пределы, производную, комплексные числа и математическую индукцию.
 	\newline
 	 В 20-30 годы
 наблюдается также идея, что
геометрию следует опереть на геометрические преобразования. В конце 1950х эта идея
возродится под влиянием Н.~Бурбаки, в начале 60 последуют неудачные попытки воплотить
ее в жизнь (массовый учебник В.~Г.~Болтянского и И.~М.~Яглома и экспериментальный учебник А.~И.~Фетисова).
Неудачи не будут должным образом оценены, и последует введение новых учебников
геометрии со сжиганием за спиной мостов...},
  однако фактически три учебника Киселева 
 и <<Тригонометрия>> Рыбкина остались стабильными учебниками вплоть до 1955 года. 
 <<Стереометрия>> Киселева пережила реформу 1955г.,  последний школьный выпуск по ней состоялся в 1976г. С выпуска 1977 началась деградация геометрии в советской средней школе, но это - совсем другая история, подробнее о школьной реформации 1960-1980, см. \cite{Ner-reform}.

 \section{Осень%
 \label{s:autumn}}
 
 \COUNTERS
  
 {\bf \punct Витязь на распутье.%
 \label{ss:rasputie}} После июля 1936г. Лузин покинул Московский университет и  Стекловку
 (по крайней мере, ее московское отделение). Он также перешел из математической группы Академии 
 в техническую (и, по-видимому, оставался там до конца жизни).
 
 Что касается ЦАГИ, то процитируем Бари:
 \begin{quotation}
 В эти же годы он пытался прочно связаться с работой в
ЦАГИ; по случайным причинам эта связь не наладилась.
(Тогдашний начальник ПУР Красной Армии, которому был подчинен
ЦАГИ, в письме начальнику ЦАГИ указал, что он считает, что
Н.Н. - ярко выраженный идеалист, так как где-то написал, что
«ряд натуральных чисел непознаваем». Это вызвало отказ Н.~Н. от
работы в ЦАГИ).
 \end{quotation}
 Если смотреть по окружающему  тексту, это должно было относиться к первой половине 30х.
 С другой стороны, в 1936г. Лузин вроде бы еще работал в ЦАГИ.
 Обвинение Лузина очень похоже на то, которое содержалось  бумаге Кольмана 1931 (см. п.\ref{ss:kolman-tsagi}).
Представляется маловероятным, что начальник Главного политического управления Красной Армии%
\footnote{Гамарник, Я. Б., 01.10.1929--31.05.1937, Смирнов П.А.  06.--12.1937, Мехлис Л.З. 30.12.1937--06.09.1940.}
мог сам по себе заинтересоваться проблемой неразрешимости арифметики. Вероятно, Кольман
написал еще одну телегу, или в связи с июльской историей 1936 года старая бумага была извлечена на свет божий.
Скорее всего, Лузин ушел из ЦАГИ тогда же в 1936г.  Впрочем, такая интерпретация известия Бари
не является единственно возможной (например, речь могла идти о ставке или полставке).

Мы знаем, что у Лузина была запрлата академика (очень высокая по общим меркам), и что он руководил семинаром по УрЧП
в технической группе академии. Были  ли у него  с осени 1936г. по конец 1937г. другие места работы, мы не знаем
(вроде бы он числился завотделом в Ленинградской Стекловке).

Процитируем Стенограмму:
 \begin{quotation}
 	ЛУЗИН. Затем у меня готовится большая работа по проективным множествам, целая книга, работа
 	была заказана два года тому назад. Работа эта уже осуществилась. Мне нужно вписать
 	только две главы, работа сделана на русском языке. Книга эта в рукописи находится у
 	Нины Карловны [Бари] и скоро появится...
 	
 	ЛУЗИН. Но, что касается того, почему я тотчас же этого не потребовал, то я скажу, что я имел
в виду написать пять нот, из которых я написал две и одновременно послал их Борелю, а
что касается остальных, то я предполагал осветить этот вопрос уже не тремя строчками,
а более подробно, использовав для этого три ноты, мною еще не опубликованные. Эти
три ноты я не опубликовал, так как я это время был занят другими заботами, т.к. мне
пришлось отложить все научные дела и прекратить публикацию.
 \end{quotation}
 
Нам ничего не известно ни о судьбе рукописи, ни об этих <<нотах>> (в Compt. Rend.,
они по понятным причинам не вышли). 
Но вряд ли случилась гибель важных текстов. Лузин имел возможность издавать свои <<ноты>> в ДАН СССР (и издавал).
Кроме того, похоже на то, что подготовка  русского перевода
парижской монографии Лузина была начата при его жизни, и здесь был сделан выбор: переводить монографию
или издавать оригинальный текст (выбран был перевод). 

Лузин после 1936г. еще издал по русски две статьи по дескриптивной теории множеств (очевидно, что не те, о которых
говорилось в 
Стенограмме), но по-крупному пути в дескриптивной теории множеств ему не было, и он сам понимал это лучше, чем кто-либо другой.
Из ТФДП он ушел уже за 20 лет до того, сразу после своей диссертации. Можно лишь фантазировать, могла ли иметь 
продолжение деятельность Лузина по периодограммам, но статья была зарублена Хинчиным, см.
\ref{ss:haltura}...

Так или иначе, перед Лузиным встал вопрос о том, что делать дальше.  В разных областях математики
царствовали его ученики или ученики его учеников, с которыми он сам давно рассорился...

Процитируем еще раз  Люстерника:
   \begin{quotation}
    Он ведь еще был не старым человеком и сильным математиком. В предвоенные
годы он «покончил одним ударом» с тематикой, которой занимались
московские математики в течение ряда десятилетий — «изгибаниями на
главном основании» (Н.~Н.~Лузин показал сравнительную узость этих преобразований).
А в послевоенные годы он «тряхнул стариной» и провел на высоком
уровне семинар по теории функций комплексного  переменного, в котором
поставил ряд задач теории функций двух переменных, и этот семинар послужил
отправной точкой серии последовавших друг за другом замечательных
исследований московских математиков по теории функций многих переменных
(начиная с работы участника этого семинара А.~С.~Кронрода).
   \end{quotation}
 
  \sm
  
{\bf \punct Изгибания на главном основании.%
\label{ss:izgibaniya}}
Так или иначе <<{\it по приглашению С.~П.~Финикова}>> \cite{Luzin-osnovaniya},
не позднее начала 1937г. (см. письмо Лузина Колмогорову, цитированное в п. \ref{ss:detail})
Лузин занялся задачей <<об изгибании на главном основании>>. 
Поясним о чем идет речь.

{\small
Рассмотрим двумерную поверхность в $\R^3$. На касательном пространстве в каждой точке
определена риманова метрика
(<<первая квадратичная форма>>),
а также <<вторая квадратичная форма>> (см. любой учебник дифференциальной геометрии).
Две линии на поверхности, пересекающиеся в данной точке называются {\it сопряженными},
если они ортогональны относительно второй квадратичной формы%
\footnote{Эта же старая терминология  встречается в учебниках
аналитической геометрии в словосочетании {\it <<сопряженные диаметры эллипса>>.}}.
Пусть дана сеть линий на поверхности (т.е. два семейства кривых, которые в окрестности точки
образуют картинку, диффеоморфную прямоугольнику с прямоугольной сеткой). Эта сеть
называется {\it сопряженной}, если кривые одного семейства сопряжены кривым 
другого семейства.

Пусть теперь две поверхности $S_1$, $S_2$ взаимно однозначно отображены
друг на друга с сохранением римановой метрики. Тогда, как обнаружил в 1866г.
Петерсон\footnote{Петерсон Карл Михайлович (1828—1881),
 автор красивых  и тонких  работ по теории поверхностей, которые были
оценены уже после его смерти. О его биографии известно мало. В 1853 по окончанию Дерптского
университета защитил кандидатскую диссертацию (в которой были, в частности, были впервые
получены условия Петерсона-Мейнарди--Кодацци). С 1865 до конца жизни 
работал учителем в Петропавловском мужском училище (Москва). 
Один из основателей Московского математического общества.} существует единственная сопряженная сеть на одной поверхности,
которой соответствует сопряженная же сеть на другой поверхности%
\footnote{Забыв про риманову метрику на $S_1$, мы можем считать, что в каждой точке $x$
поверхности $S_1$ определены две вторых квадратичных формы $L_x$, $M_x$,
одна ее собственная, а другая -- перенесенная с $S_2$. 
Мы пишем уравнение $\det(L-\lambda M)=0$. В точках общего положения это уравнение
имеет два различных корня $\lambda_1$, $\lambda_2$. Соответственно в каждой точке возникает два собственных направления 
$\ker (L-\lambda_1 M)$, $\ker (L-\lambda_1 M)$. Обе формы диагонализуются в 
получаемой системе координат. Остается провести интегральные кривые
к этим касательным направлениям. Отметим, что может случится так, что $\lambda_1$, $\lambda_2$
- это пара комплексно сопряженных чисел Тогда и сопряженная сеть будет комплексной. Что ж поделать.}
(при некоторых условиях невырожденности, которые ясны из сноски).
 
Эту сеть линий Петерсон назвал <<основанием изгибания>>. Если при сохранении этой сети возможно непрерывное семейство
изгибаний, то такая сеть была названа главным основанием (на самом деле как показал Петерсон, наличие двух
изгибаний на данном основании влечет наличие непрерывного семейства). В наши цели не входит
ни обсуждение дальнейших исследований Петерсона по изгибанию поверхностей, 
ни продолжение тематики главных оснований в Московской 
геометрической школе, где она была очень популярна, такие работы были и в Западной Европе (см. \cite{Finikov}).

Хотя и самим Петерсоном, и его последователями было построено много примеров изгибаний на главном основании,
оставалось неизвестным, для всех ли поверхностей возможны такие изгибания. Этой задачей, которая несколько
десятилетий не поддавалась усилиям Московской школы, в 1936-1937гг. занялся Лузин. Он быстро
сумел войти в новую для него тематику, в 1938г. опубликовал  несколько заметок в Докладах,
а 1939 - большую работу в Известиях АН СССР, однако, не по математическому, а по техническому
отделению.

Он пишет (весьма замысловатую) систему нелинейных УрЧП-ов на сети, которые могут быть главными основаниями,
в терминах коэффициентов римановой метрики. Далее  показывает, что система не всегда решается,
а условие ее совместности -- система УрЧП-ов на коэффициенты римановой метрики.
Уравнения получаются приравниванием к нулю определенных многочленов от частных
производных (каковые многочлены явно не выписаны). Таким образом получалось,
что на  поверхности с фиксированной римановой метрикой, вообще говоря, главных оснований нет.

Лузин пошел дальше. Рассмотрим  поверхность с фиксированной
римановой метрикой (пусть у нее есть главное основание) и фиксируем ее вложение в $\R^3$. Оказывается, что для общего вложения
любой поверхности
главных оснований тоже не существует.}

Лузин снова показал свой класс (да кто ж в этом классе мог бы усомниться!),
но результат был отрицательным и близким к закрытию тематики
\footnote{Из статьи Лузина \cite{Luzin-osnovaniya}:
\newline
\it
 Все это заставляет смотреть на свойство поверхности «обладать 
главным основанием» как на весьма частное свойство, имеющееся лишь
в исключительных случаях. Не может быть никакого сравнения, по
нашему мнению, сети главного основания с сетью линий кривизны, 
асимптотических линий, семейством геодезических линий и т. д. Последние
имеются на каждой аналитической поверхности, и их изучение, по всей
справедливости, образует главы классической дифференциальной 
геометрии, тогда как феномен главного основания, заслуживающий названия
«феномена Петерсона», может составлять не главу, но лишь параграф
дифференциальной геометрии...
\newline
Но, с другой стороны, не следует и недооценивать смысл этого 
феномена. Если бы главное основание составляло особенность какого-нибудь
одного, ранее известного, классического семейства поверхностей, интерес
феномена Петерсона был бы сведен на нет, так как главное основание
было бы просто свойством этого семейства поверхностей. И в таком
случае это обстоятельство, без сомнения, было бы давно усмотрено. Но
на самом деле этого нет, и феномен Петерсона дает нам какую-то, еще
мало понятную «косую» классификацию поверхностей, при которой 
наличие главного основания оказывается затрагивающим многие классические
семейства поверхностей без того, чтобы содержать их целиком. Надо
учесть еще, что погоня за общностью осуществляется ценою утраты
глубины и открытия неожиданных, хотя и частных, но глубоких соотношений. 
Так, на долю теории чисел выпало бы весьма мало теорем, если
бы она сосредоточивалась на понятии общего иррационального числа...}$^,$%
\footnote{Незадолго до смерти Лузин вернулся к задаче об условиях возможность изогнуть
вложенную поверхность на главном основании, от этой деятельности остались 350 страниц рукописей,
отрывок из которых был издан \cite{Luzin-add3}, \cite{L-Sret}.}

\sm

Выбор Лузина оказался неудачным.

\sm

{\bf \punct Другие проекты.%
\label{ss:projects}}
Как будто больших новых исследовательских проектов (по крайней мере, успешных) со стороны Лузина
больше не было. Однако в целом его деятельность  после 1936г. отнюдь не была бесплодной и 
могла бы составить славу какого-либо другого математика.

Одновременно с дифференциально-геометрическими изысканиями Лузин пишет учебник 
\cite{Luz-real}
теории функций действительной переменной для педвузов. Эта книга дважды издавалась
в 1940 и 1948гг. Автор не пытался анализировать ни ее содержание, ни ее влияние.
Лузин также продвигал другой образовательный проект, речь о котором пойдет в следующем параграфе.

В середине 40х Лузин снова вел семинар по теории функций действительной переменной в МГУ. Среди 
участников семинара был Кронрод Александр Семенович (1921-1986), защитивший 1949
под руководством Лузина кандидатскую диссертацию по теории функций
двух вещественных переменных, работа сразу была признана докторской.
В дальнейшем Кронрод был одним из основателей Computer science в СССР, а также 
инициатором образовательной технологии, использовавшейся в разное время в  7, 57, 91, 179 и 444
 школах Москвы, эту технологию далее развивал Н.~Н.~Константинов.
Повлиял ли тут лузинский семинар, или это было личным изобретением Кронрода,
автору не известно\footnote{Одним из нетривиальных
	элементов технологии был бег на месте. В дополнение
к школьным (хотя и несколько углубленным) <<Алгебре>> и <<Геометрии>>
вводился предмет с названием <<Математический анализ>>,
где все делалось не спешно: аккуратно вводились действительные числа, 
элементарная теория множеств, теория пределов и непрерывность
в довольно изощренной форме,   для продвинутых была мера Лебега.
Главным был стиль обучения: школьникам давали отпечатанные на пишущей машинке листочки с задачами,
и каждый школьник должен был лично сдать решение, по решению предполагался разговор.
Фактически урок состоял в цепочке индивидуальных разговоров 4-6 преподавательствующих
с 24-30 школьниками, выступлений у доски не было. Фактически получался бег на месте.
Пройденный материал был минимальным, но школьники обучались навыкам самостоятельного
мышления. Вспомним, что ТФДП 20х было школой мышления для московских математиков,
и теперь ТФДП (уже совсем элементарная) снова использовалась, но уже в старших классах.
Пожалуй отмечу, что при чтении рассказов про Лузина-лектора мне все время вспоминался Константинов
у доски, произносивший речи для участников олимпиад или кружков.
Возможно, что технология Кронрода--Константинова была гибридом лузинской школы и старых московских
математических кружков.}.

\sm

{\bf\punct Последние годы.} Процитируем Бари:
      \begin{quotation}
      	Работа его в Институте имени В.~А.~Стеклова оборвалась в
      	самом конце 1947 года вследствие резкого столкновения с одним
      	из академиков. После этого он уже не вернулся в Институт 
      	имени В.~А.~Стеклова. (На заседании математической секции АН 
      	обсуждались кандидаты в академики. Их было три: 
      	П.~С.~Александров, которого поддерживал А.~Н.~Колмогоров, И.~Г.~Петровский,
      	которого поддерживал И.~М.~Виноградов, и <Н.~Г.~Чеботарев>, 
      	которого поддерживал С.~Н.~Бернштейн. После острых споров, при 
      	выходе А.Н.Колмогоров ударил Н.~Н. по лицу. - Авт. рук.)
      \end{quotation}
   
   \begin{quotation}
    Родители Н.~Н. были здоровья слабого, особенно мать, которая
страдала застарелым пороком сердца и умерла рано. Эта слабость
здоровья перешла и к Н.~Н. Впервые случился тяжелый сердечный
припадок в 1941 г., в самом начале войны; положение было 
настолько тяжелое, что Н.~Н. пришлось провести в постели в то время
несколько месяцев. Повторялись в слабой степени сердечные 
припадки, несмотря на внимательное лечение, и в последующие годы.
28 февраля 1950 года в результате острого сердечного припадка
(инфаркт) Н.~Н. внезапно скончался на 67 году жизни. Похороны
состоялись 3-го марта на Введенском кладбище.
   \end{quotation}

 \section{Математика для  нематематиков%
 \label{ss:non-math}}
 
 \COUNTERS

 Еще одна сторона деятельности Лузина -- издание учебников для нематематиков.
 
 \sm
 
 {\bf  \punct  Гренвиль, Гренвиль--Лузин, Лузин.%
 \label{ss:graneville}}
 Напомним (подробнее см. \cite{Kol-savvina-2}), что учебник анализа Гренвиля был переведен на русский
 Маракуевым Николаем Николаевичем (1847–1911) и издан в 1912г.
 По жанру это был хороший элементарный курс анализа, рассчитанный на инженеров
 (и, быть может, представителей естественных наук и, быть может, на педагогов-математиков).
  Следующее издание вышло 1922г.
 с участием Лузина и Тарасова%
 \footnote{Тарасов Николай Петрович, учился на Московском Физмате,
 не окончил его, преподавал на рабфаке Иваново-Вознесенского политеха, где познакомился с Лузиным,
 в 1924г. все же закончил Московский университет, в дальнейшем  работал преподавателем,
 издавал учебники, в частности «Курс высшей математики для техникумов», выдержавший 17 изданий,  заведовал кафедрой в МАТИ.} 
 
 Лузин в 1907--1930гг. работал на Московском Физмате, место это особое,
 но  ему довелось 
  преподавать также в Иваново-Вознесенском политехе
  1918--1920 и в Лесотехническом институте 1920--1925. Он занялся редактированием
  учебника Гренвиля (Грэнвиля), внося в него разные дополнения, 
  список трудов Лузина показывает издания 1922, 1924, 1926, 1927, 1928,
  1928-1929, 1930. C 1930г. учебник начинает выходить как
  Гренвиль--Лузин, сначала как <<Элементы дифференциального и интегрального исчислений>>,
  потом как  <<Курс дифференциального и интегрального исчислений>>.
  Были издания 1930, 1930, 1930, 1931, 1933, 1933, 1934, 1934, 1935, 1935-1937,
  1937, 1938, 1942, общий тираж 270 000 экземпляров. Мы видели, что учебник 
  подвергался критике и нападкам,  выше цитировались Выгодский, который счел
  его антиисторическим, Люстерник, а также Яновская, которая 
  применила слова <<возмутительное вредительство>>. Но учебник,
  очевидно, был востребован высшей технической школой.
  
  В 1946г. вышел новый учебник Лузин <<Дифференциальное исчисление>> и
  <<Интегральное исчисление>>. Он четыре раза издавался с грифом <<Учебник>>,
  1946, 1949, 1952, 1953 (тиражи <<Дифференциального исчисления>> соответственно
  50 000, 100 000, 50 000, 50 000) и трижды как <<Учебное пособие>>,
  1955, 1958, 1961 (60 000, 150 000, 100 000). Возможно, что изменение
  грифа связано с тем, что эта книга была больше рассчитана
  на то, чтобы по ней учиться, чем на то, чтобы по ней сдавать и принимать экзамены.
  
    Это  900-страничный элементарный учебник анализа,
  издалека он выглядит достаточно необычно. Лузин избегает разных
  ловушек, в которые попадали многие более поздние математические курсы. В нем отсутствуют элементы развития педагогики
  на своей собственной основе (когда плодятся искусственные трудности,
  которые преподавателям приятнее и легче преодолевать, чем естественные) и минимизируется возможный материал для заучивания. Одновременно автор учебника пытается всеми правдами и неправдами добиться понимания предмета
  и не заминать места, где возможно возникновение непонимания. Что касается формального уровня строгости,
  то Лузин временами ведет игру на грани фола, преднамеренно предлагая неаккуратные формулировки,
  рассуждает исходя из них, но далее говорит формально четко. Курс, по-моему, преднамеренно просчитан 
  на возможности разного уровня прочтения.
   Стоит отметить и то, что Лузин делает ставку
  на геометрическую восприятие.

  В 60е-начале 70х годов XX века была предпринята попытка модернизации 
  курсов математики для нематематиков в высшей школе. История была сложной 
  и драматичной, с успехами и неудачами, и вряд ли когда-нибудь эта история будет написана. 
  Но в
   среднем картина 
  к 80му году была  плачевной (см., например,
  лекцию В.~А.~Рохлина, 1981, \cite{Rokhlin}), дальше уже
  происходило дезорганизованное отступление
   без какого-либо видения стратегической перспективы  со стороны математического социума.
 
   Учебник Лузина пришелся на благополучную эпоху положительного отношения
  общества к математике и по-видимому относился к числу многих факторов, формировавшим это отношение (как и учебники Фихтенгольца и В.~И.~Смирнова, справочники Выгодского или 
  школьные учебники Киселёва).
   Учебники анализа почему-то сами собой перестают быть годными со временем,
   быстрее, чем становится архаичным их содержимое.
   (об этом писал сам Лузин, см. чуть ниже), вряд ли учебник Лузина
   может быть успешным сегодня. Но  для преподавателей, самостоятельно
   просчитывающих курсы или для авторов новых учебников книга Лузина  
   может быть серьезным и полезным источником для раздумий.

 Небезынтересное обсуждение педагогических взглядов Лузина есть в статье В.~Л.~Минковского 
 <<О методико-математических воззрениях, 1930 Лузина>>,
 \cite{Mink}. Мы лишь приведем цитаты из двух введений Лузина к учебникам
 анализа, где он, отчасти, сам высказывает свои взгляды.
 
 \sm
 
 {\bf \punct  Из введения Лузина к учебнику Грэнвиля, 1924.%
 \label{ss:graneville-preface}}
 \begin{quotation}
 	Задача составления «Курса Анализа» для высших технических школ 
 	есть одна из самых трудных и, вместе с тем, одна из наиболее привлекательных тем для лиц, занимающихся распространением математических знаний
  Трудность этой проблемы вырастает из тех заданий, в которых
 	она должна решиться.
 	
 	{\it Во-первых}, намечаемый круг читателей обуславливает помещение в книге лишь самых необходимых сведений из анализа, без которых
 	не может обойтись ни один инженер.
 	
 	{\it C другой стороны}, самое изложение этих сведений не может рвать
 	с современным состоянием математических знаний, отставая от него на половину или на
 	три четверти века.
 	
Составление всякого курса анализа должно следовать по ра\-вно\-дей\-ст\-ву\-ю­щей этих двух почти взаимно исключающих друг друга требований. Отсюда - трудность проблемы.

Каждый составитель руководства анализа знает на опыте ту силу искушения, которую ему приходится, преодолевать,
чтобы не подпасть под преобладающее влияние одного из перечисленных требований. И, в самом деле, 
здесь для составителя нужна большая осторожность. 	

Во-первых, негармоничное подчинение составителя первому требованию неизбежно влечет чрезмерное упрощение текста, 
переходящее порою в прямое его опрощение. Следуя этому пути, встречают, прежде всего, тот класс учебников анализа,
в которых говорится о непрерывных функциях таким языком, как если бы они всегда имели производную, и где обращаются 
с бесконечными пределами определенного интеграла так, как если бы пределы эти были числами, хотя это влечет 
уже к фактическим ошибкам счисления. На этом же пути, спускаясь несколькими ступеньками ниже, встречают учебники,
в которых пишут знаки многозначных функций, не упоминая о том, какое именно значение принимается за основное, и где разлагают функции в ряды, не указывая на то, где именно пригодны эти разложения. Дальнейшее следование по этому пути приводит к учебникам  с явственно выраженным распадом теоретического содержания. В книжках этого рода, обычно весьма тонких по объему, авторы стараются тем не менее удержать всю совокупность фактов, ценою удаления соединительной
 логической ткани между ними, не замечая, что именно эта операция и делает их книги непреодолимо трудными для учащихся, предоставленных своей одной лишь механической памяти, так как является удаленным именно то самое, что скрепляло эти факты и могущественно помогало памяти....
 
 {\it Во-вторых}, также негармоничное подчинение составителя второму 
 требованию обычно приводит к тому, что написанный составителем учебник имитирует  университетский курс анализа. 
 Часто бывает, что такая книга не удовлетворяет ни инженеров, ни учащихся университетов, будучи слишком теоретичной 
 для первых и слишком недостаточной для вторых%
 \footnote{Для сравнения, В.~А.~Рохлин, 1981, \cite{Rokhlin}:
 \newline
 \it
 ... Как быть с людьми, которые проявляют интерес к математике, но обучаются по программам, по которым научиться нельзя?....
  Однако часто и программы, и соответствующие учебники не являются самостоятельными, а представляют собой просто испорченные курсы,
  по которым готовят математиков....}. Чтобы понять, почему это происходит, достаточно обратить внимание
 на самый процесс составления таких руководств. Обычно при этом исходят из какого-нибудь университетского курса анализа,
 гарантирующего научность и  логичность  изложения, и затем подвергают его осторожному процессу сокращения,
 удаляя из него факты, ненужные для инженера, и оставляя с легкой модификацией теоретические  пояснения и взаимные связи нужных фактов. 
 Но при этом оставлении части теоретического материала представляется почти невозможным сохранить чувство меры: 
 сила логической связи фактов  в действительно научном курсе является настолько непреодолимой, что составитель, 
 принужденный следовать ходу логической мысли, продолжает двигаться в этом направлении, так сказать,
 но инерции и неизбежно вводит вещи, уже в самом деле чуждые для инженера. Именно вследствие этой причины,
 почти во всех хороших учебниках для высших технических школ мы находим доказательство иррациональности числа $e = 2, 71828 ... $,
 искусившее составителей своею краткостью и легкостью, но совершенно излишнее для инженера, которому нужны только три его десятичных знака.
 Другою
 причиною этого слишком близкого следования университетскому курсу является,
 по-видимому, недостаточно длительный опыт составителей в их личной работе над основами анализа,
 вследствие чего в учебник проскальзывают  многие весьма затруднительные для учащихся рассуждения только потому,
 что составителю они кажутся совершенными в научном отношении, хотя ближайшее рассмотрение их часто обнаруживает,
 что в смысле их строгости они недорого стоят и всегда могут быть заменены другими, более интуитивными и столь же научными.
 Этим может быть объяснено сохранение в учебниках многих рудиментарных понятий и теорий.
 \end{quotation}
 
 Если вспомнить  недавние времена, то в последние 30 лет курсы математики на нематематических специальностях институтов и университетов пребывают в состоянии постепенной деградации и распада.
 Комментарии Лузина выглядят весьма актуально, причем часто реализуются
 оба упомянутых сценария одновременно...
 
 \sm 
 
 {\bf\punct Из введения Лузина к учебнику анализа И. И. Жегалкина и М. И. Слудской.%
 \label{ss:zhegalkin-preface}}
 
 \begin{quotation}
Решающая реформа преподавания пришла, если говорить о Франции,
вместе с преподавательской деятельностью директора Высшей нормальной
школы в Париже Jules Tannery...
 	
 Основной идеей реформы преподавания Jules Tannery было систематическое расчленение на простые элементы тех умственных актов, которые 
 необходимы для овладения началами математического анализа, и доведение
 этих элементов до возможно более ясного понимания читателями и слушателями.
 
 Jules Tannery никогда не прибегал в своем преподавании к формальным
 фокусам, всегда считая их очень вредными и постоянно предупреждая о том,
 что всякое искусственное сокращение в процессах мысли оплачивается слишком дорогой ценой непонимания самого существа дела. Все его математические рассуждения всегда были прямыми и иногда длинными, но столь
 тщательно разработанными, расчлененными и до такой степени ясными, что
 овладение ими уже не представляло ни малейшего труда. Jules Tannery все
 время внимательнейшим образом следит за состоянием ума начинающего,
 предвидит возникновение в нем сомнений, иллюзий и заблуждений, предупреждает о них и каждое мгновение готов прийти на помощь читателю. И,
 действительно, читая его, всегда испытываешь странное чувство отождествления ума автора и читателя. «Только то и прочно, что понятно» было любимой
 фразой Jules Tannery, в противоположность известному изречению Даламбера, сказанному им в хаосе идей XVIII в.: «Идите дальше, потом когда-нибудь
 поймете»....

 Первая и основная идея И. И. Жегалкина — отрицательная, это есть
 совершенно ясно осознанная им невозможность исходить при составлении
 учебника от обычного представления об идеальном читателе. А между тем,
 большинство учебников именно и отправляется от этого представления, наделяя
 этого абстрактного читателя беспредельными внимательностью, понятливостью,
 догадливостью и сообразительностью. Для этого читателя нет ни в
 чем никаких затруднений и препятствий: достаточно автору хотя бы один
 раз указать на какое-нибудь обстоятельство, как этот идеальный читатель
 уже понимает его с полуслова и запоминает если не на всю жизнь, то во всяком
 случае на все время изучения книги. Поэтому-то авторы, следующие по
 этому пути, имеют необыкновенно экономные размеры учебника. Поэтому-то авторы, следующие по
 этому пути, имеют необыкновенно экономные размеры учебника. Книги этого
 рода обычно кажутся весьма привлекательными для издательств, довольных
 возможностью преподнести дифференциальное и интегральное исчисления уложенными в небольшое число страниц, и для учащихся, еще не знакомых
 с той истиной, что чем толще учебник математического анализа, тем он скорее будет прочтен и усвоен. Но уже очень скоро у учащихся наступают разочарование и охлаждение к книге или, что бесконечно хуже, к самому предмету.

 Когда вдумываются в причины возникновения иллюзии «идеального читателя»,
  то немедленно замечают, {\it что под таким читателем автор просто
 разумеет себя самого и именно то самое состояние своего ума, которое он имеет
 в момент создания учебника, но отнюдь не то состояние ума, которое было
 у автора, когда он сам впервые знакомился с излагаемыми им идеями}[выделение текста здесь и ниже - авторское]....

Таким образом, самой характерной чертой предлагаемого курса анализа, чертой, отличающей его от всех остальных курсов, является {\it исключительная ориентировка его на понимание учащимся всех процессов рассуждения}. При этой ориентировке на понимание нисколько не страшны дефекты
памяти, так как самый ход однажды понятого материала не позволяет утратиться существенному, деталь же легко восстановить и по любому справочнику. Ориентировка на понимание кажется громоздкой лишь в начале, но
на деле затраченное на понимание время с лихвой окупается в дальнейшем,
так как при правильно понятых основаниях дисциплины дальнейший материал принимает характер лишь упражнений в давно известном, чем создается уже экономия во времени.

Следует еще указать на то, что ориентировка курса на понимание имеет
в виду неизмеримо более важное соображение, чем простое облегчение учащемуся трудностей изучения: это -- {\it инициативу учащихся}. Эта драгоценная инициатива, т. е. {\it умение ориентироваться в новой обстановке, выходящей
из привычного шаблона}, может прийти и всегда приходит лишь от совершенного понимания материала, а вовсе не от твердости его усвоения, т. е. вбирания памятью.

Следует иметь в виду, что выпадение из памяти какой-нибудь важной
части материала приводит в расстройство и весь остальной, вообще
казавшийся хорошо усвоенным материал только в том случае, если он раньше
не был ориентирован на понимание. 	
 	\end{quotation}

\section{Эпилог%
\label{s:epilog}}

\COUNTERS

  \epigraph{ Возьмите историю развития советской математики.
   Вся она - несколько поколений! - началась из школы академика Николая Лузина,
   которая возникла во время гражданской войны.}
   {И. Р. Шафаревич}

   Стоит отметить, что только что процитированный Шафаревич был учеником Б.~Н.~Делоне, который не относился к школе Лузина.
   Впрочем учился он в Москве, где атмосферу определяли выходцы из Лузитании.
   	
   	\sm

{\bf \punct Последние из лузитан.%
\label{ss:last}} Один из героев нашего повествования,
Меньшов Дмитрий Евгеньевич, с 1941г. работал в должности заведующего
кафедрой Теории функций и функционального анализа МГУ. Если автору не изменяет память, в последний раз его видели на Мехмате весной 1977г. В какой-то из дней его увезли в академический санаторий <<Узкое>>. В марте 1980г.
А.~П.~Юшкевич и С.~С.~Демидов взяли у него интервью <<Воспоминания о молодых годах и о возникновении Московской школы теории функций>>, которое появилось в 1983г.
в 27 выпуске <<Историко-математематических исследований>> (предыдущий выпуск был в 1982г., т.е., воспоминания Меньшова - человека
вообще не последнего и не последнего из авторов, публиковавшихся в историко-математических исследованиях,
 - были задержаны с публикацией). Однако опасных вопросов интервью никак не коснулось.
Это и не удивительно. Из редакционного введения к интервью:
\begin{quotation}
В основу текста «Воспоминаний» чл.- корр. АН СССР Д.~Е.~Меньшова легла магнитофонная запись его беседы с А.~П.~Юшкевичем в марте 1980 г. в санатории АН СССР «Узкое». Текст записи, переписанный с ленты Е.~И.~Славутиным и отредактированный А.~П.~Юшкевичем, был в ноябре того же года дополнен и проверен Д.~Е.~Меньшовым. В беседе участвовал также С.~С.~Демидов. «Воспоминания» публикуются в форме вопросов и ответов, какую они имели вначале; вопросы и ответы помечены буквами В. и О. {\bf Редакция благодарна академикам П.~С.~Александрову и А.~Н.~Колмогорову, которые частично познакомились с согласия Д.~Е.~Меньшова, с его текстом и внесли в него несколько уточнений.}	
\end{quotation}

В день столетия  Лузина 8.12.1983 Колмогоров отправил в редакцию сериала такое письмо (опубликовано в 1985г, \cite{Kolmogor-letter}):
\begin{quotation}
В «Воспоминаниях о молодых годах и о возникновении Московской школы теории функций» Д.~Е.~Меньшова, напечатанных в XXVII выпуске, на с. 327—328 говорится о моих отношениях с Н.~Н.~Лузиным, среди прочего, следующее: «По существу А.~Н.~Колмогоров не нуждался в руководстве. Колмогоров не был в числе прямых учеников Н.~Н.~Лузина, а занимался с В.~В.~Степановым. В 1925—1929 гг. (в аспирантуре) официальным руководителем А.~Н.~Колмогорова был Н.~Н.~Лузин, которого во время его заграничной командировки заменял я. Вообще же А.~Н.~Колмогоров шел им самим выбранным путем» и т. д. Я хотел внести некоторые уточнения. На самом деле я являюсь прямым учеником Н.~Н.~Лузина. В 1921 г. он узнал о моем первом результате в теории тригонометрических рядов (существование рядов Фурье — Лебега со сколь угодно медленно убывающими коэффициентами). Тогда же Н.~Н.~Лузин пригласил меня регулярно бывать у него дома, на Арбате, 25, вместе с небольшой группой моих сверстников (Л.~В.~Келдыш, П.~С.~Новиков, И.~Н.~Хлодовский). Беседы эти были для меня чрезвычайно полезными и продолжались до моего поступления в аспирантуру в 1925 г., когда Н.~Н.~Лузин стал моим официальным научным руководителем. Продолжавшиеся во время аспирантуры (за исключением поездок Н.~Н.~Лузина за границу), беседы с ним имели для меня большое стимулирующее значение.

Разумеется, беседы и занятия с В.~В.~Степановым и Д.~Е.~Меньшовым также сыграли свою роль в формировании моих интересов.

Москва, 8 декабря 1983 г, А.~Н.~Колмогоров	
\end{quotation}

Там же был опубликован ответ редакции, истинный смысл которого остался ясен лишь участникам переписки:
\begin{quotation}
Помещая письмо академика А.~Н.~Колмогорова, редколлегия «Историко-математических исследований» должна со своей стороны сказать, что в том месте «Воспоминаний» чл.-корр. АН СССР Д.~Е.~Меньшова, которое упомянуто в этом письме, была допущена по недосмотру одного из членов редколлегии ошибка при публикации. На самом деле в оригинальном тексте Д.~Е.~Меньшова было написано: «По существу А.~Н.~Колмогоров не нуждался в чьем-либо руководстве. Вначале Колмогоров был приглашен в число учеников Н.~Н.~Лузина, а затем занимался с В.~В.~Степановым. В 1925—1929 гг. (в аспирантуре) официальным руководителем Колмогорова был Н.~Н.~Лузин, которого во время его заграничных командировок заменял я. Вообще же А.~Н.~Колмогоров шел им самим избранными путями» и т. д.

Письмо А.~Н.~Колмогорова существенно уточняет вопрос о его отношениях с Н. Н. Лузиным в 20-е годы нашего века. Что касается заграничных поездок Лузина, то они имели место в первой половине 1926 и 1927 гг. и затем с лета 1928 по лето 1930 г.

Редколлегия приносит извинения Д.~Е.~Меньшову и читателям за неточность, о которой идет речь, и благодарит А.~Н.~Колмогорова за его разъяснение%
\footnote{Стоит отметить, что в коротком тексте <<О первых шагах школы Лузина>> \cite{Menshov-L-2}, опубликованном в Стекловских трудах в 1983г., Меньшов просто упоминает Колмогорова через запятую в числе учеников Лузина.}.	
	\end{quotation}

Александров Павел Сергеевич - основной%
\footnote{или один из двух основных.} противник Лузина - умер в ноябре 1982 года. 

С 9 по 18 сентября 1983 г. в г.~Кемерово проходила Всесоюзная школа по теории
функций, посвященная 100-летию со дня рождения Лузина. За два дня до ее начала
кемеровская университетская многотиражка опубликовала интервью с Колмогоровым, взятое В.~А.~Успенским, 
`<<{\it Ученик об учителе}>> \cite{Kolmogor}. Это было громкое широковещательное заявление,
сделанное для заинтересованных узких кругов. 

Андрей Николаевич успел при жизни перепечатать это интервью трижды, в Вестнике АН СССР (1984), 
в Успехах мат. наук (1985) и в своей популярной книге <<Математика. Наука и профессия>> (1988).

Так завершилось сорокопятилетнее противостояние, связанное с именем Лузина. С осени 1983 года и по 1985г. 
советская математика широко отмечала 100-летие со дня рождения своего основателя. 
А вскоре начался закат пришедшей когда-то вместе с диссертацией Лузина великой эпохи... 

\sm

Мы закончим несколькими цитатами.

\sm

{\bf\punct Шафаревич.}  Из интервью, данного В.~Н.~Тростникову, \cite{Shaf-cult}.

\begin{quotation}
Т. Ну, скажем, Лузин был дореволюционный математик...

Ш. Главный успех его был как раз в эпоху после революции. Создание его школы произошло в эпоху Гражданской войны.
Один из учеников Лузина любил писать полукомические стихи.
И он писал о их школе, которую они называли Лузитания: «Катком замерзли коридоры, горячие здесь только споры».
Эпоха разрухи, голода.

Т. Вы знаете, я читал у Колмогорова какие-то записки, он сам ученик Лузина, но младший. 
Так вот, Колмогоров писал, что он сделал какую-то курсовую работу в университете на первом курсе и ему,
по-моему, мешок крупы дали. И он говорит, я чувствовал себя просто богатым, потому что я, мол, кашу ел эту.
Вот, смотрите, значит, советская школа математики зародилась в очень трудное время, безусловно. 
Причем интересен еще такой момент – эти математики свободно тогда ездили за границу, да?

Ш. Да. Были так называемые рокфеллеровские стипендии, которые получали многие на несколько месяцев, на полгода.
С семьей уехать можно было за границу и там поработать. Причем из лузинской школы там ни одного не осталось%
\footnote{Из письма Александрова Колмогорову от 30 марта 1931:
{\it Теперь о <<моих>> делах. Постарайся иметь разговор с Выгодским,
Хинчиным и др. о том, чтобы, по возможности, хотя бы за месяц до
моего приезда меня перевели в разряд персонально оплачиваемых 
работников. Это тоже не в очень серьезном порядке, но отказываясь
здесь от предложений профессуры в 625 долларов в месяц, я надеюсь,
что я в своем отечестве могу рассчитывать на профессуру
в 625 рублей в месяц.}}.

Т. То есть невозвращенцев не было?

Ш. Несколько невозвращенцев было из Ленинграда. Но мне кажется, что это второстепенные математики.
А вот из известных имен я ни одного не знаю. Во всяком случае, из лузинской школы.

Т. То есть люди ехали в относительно сытую Европу, благополучную, и возвращались в голодную Россию, да?
Вот как вы думаете, что же заставляло их так поступить?

Ш. Я думаю, прежде всего, что интересно было. Лузин сумел создать вокруг себя атмосферу духовного горения.
И не могли ничего эквивалентного найти на Западе, во что бы они могли вписаться. 
Ну, там были такие школы, но надо было в них войти. А тут другая была школа,
в которой они себя чувствовали своими людьми, говорили математически на своем языке, 
свои проблемы у них ценились и так далее. Ну, вообще мне кажется, что вопрос ваш, 
как вы его сформулировали, он более широкий. Почему в такую суровую эпоху у нас была яркая культура?
\end{quotation}

{\bf \punct Егоров.}
Из письма Лузину, 19 мая(1 июня) 1906:
\begin{quotation}
Что касается до способов научной работы, то мне думается, что здесь надо больше свободы
и непосредственности: работайте, как работается, и не задавайтесь предвзятыми идеями.
В конце концов, научная работа имеет много общего с поэтическим творчеством, и всякое принуждение 
(со стороны или от себя лично во имя какого-нибудь принципа) здесь ни к чему не приводит.	
	\end{quotation}
   
{\bf \punct Лузин.} Доклад в <<Узком>>, процитированный Бари в \cite{TyulinaA}: 
   \begin{quotation}
    Математики изумляются гармонии чисел и геометрических
форм. Они приходят в трепет, когда новое открытие открывает им
неожиданные перспективы. И та радость, которую они 
переживают, разве это не есть радость эстетического порядка, хотя обычные
чувства зрения и слуха здесь не участвуют. Мало, очень мало 
избранников способны полностью вкушать эту радость, это — правда,
но разве не то же самое происходит и с искусствами, самыми 
бесплодными? Математик есть подлинный художник, потому что он
созерцает ту гармонию идей, которая одна только и прекрасна, и,
следовательно, достойна наших усилий.

Математик изучает свою науку вовсе не потому, что она 
полезна. Он изучает ее потому, что она прекрасна. Если бы она не была
прекрасной, не стоило бы затрачивать труда знакомиться с ней, не
стоило бы на нее тратить своей жизни. Еще раз повторяю, что я
говорю не о той красоте, которая поражает наши чувства, т.е. о
красоте качества и видимого. И не потому, что мы пренебрегаем
этой красотой, но просто потому, что с нею нечего было бы делать
пауке. Я говорю о красоте более глубокой, проистекающей из 
гармонии и согласованности во едино всех частей, которую один лишь
чистый интеллект и сможет оценить. Именно эта гармония и дает
твердую основу тем красочным видимостям, в которых купаются
наши чувства. Без этой твердой поддержки красота этих бегущих
образов была бы несовершенной, потому что была бы 
колеблющейся и всегда изменяющейся. Наоборот, интеллектуальная красота
пребывает вечной и неизменной, и единственно ради нее математик
и осуждает себя на долгие и тягостные работы.
Нужно ли еще прибавлять, что в развитии этого чувства 
интеллектуальной красоты лежит залог всякого прогресса?
   \end{quotation}

{\bf\punct Лаврентьев.}
\begin{quotation}
   Заканчивая, мне бы хотелось еще раз подчеркнуть, что Лузин не только
многих из нас научил одержимости в преследовании намеченной цели, но он
также показал, как надо увлекать молодежь на научный подвиг.
\end{quotation}

%Лузин H. H. Об арифметических методах математиков XVII века / Публикация С. С. Деми-
%дова, вступительная статья и примечания И. Г. Башмаковой // Вопросы истории естествознания
%и техники. 1933. № 3. С. 25—37.

%5. Шафаревич И. Р. Пьер Ферма и развитие теории чисел (к выходу русского издания тео-
%ретико-числовых трудов П. Ферма) // Вопросы истории естествознания и техники. 1933. № 3. С.
%37—40.

%Гантмахер, редакторство?

%Э. Борель и
%А. Лебег обратились в фонд Рокфеллера, который согласился субсидировать Лузина до 31 мая
%1930 г.

%Третье письмо H. H. Лузина M. Я. Выгодскому и связанные с ним архивные документы /
%Публикация и примечания Н. С. Ермолаевой и Т. А. Токаревой // Историко-математические иссле-
%дования. Новая серия. М., 1999. Вып. 3(38). С. 100—118.

%47. Письмо Н. Н. Лузина к О. Ю. Шмидту / Публикация

%Бутягин ZB

%Я хочу сказать следующее:
%очень неприятно, что нет товарища Шмидта, но с ним я буду иметь особый разговор в
%порядке партийной дисциплины,

 %Советские ученые: очерки и воспоминания
 
 %Синкиевич
 
% http://www.spbgasu.ru/upload-files/vuz_v_licah/publish/sinkevich_gi/70.pdf

 % https://work-way.com/blog/2015/05/23/e-kolman-boevye-voprosy-estestvoznaniya-i-tex-niki-v-rekonstruktivnyj-period/`
 
 \def\Ru{(Russian)}
 \def\TR{(Translated from the Russian)}
 \def\Fr{(French)}
 \def\Ge{(German)}
 \def\En{(English)}
 
\def\RMS{Russ. math. surv} 
\def\UMN{Usp. mat. nauk}
\def\umn{Усп. мат. наук}

\def\msb{Мат. сборник}
\def\Msb{Mat. sb.}
\def\MSB{Sb. math}

\def\izv{Изв АН СССР, сер. мат.}
\def\izvv{Изв АН СССР}
\def\IZV{Math. USSR Izv.}
\def\Izv{Izv. Akad. Nauk SSSR, Ser. Mat.}
\def\Izvv{Izv. Akad. Nauk SSSR}

\def\vestnik{Вестник АН СССР}
\def\Vestnik{Vestnik Akad. nauk SSSR}

\def\CR{Comptes Rendus Acad. Sci.}

\def\shkol{Мат. в шк.}
\def\Shkol{Mat. v shkole}

\def\istmat{Ист.-мат. иссл.}
\def\Istmat{Ist.-mat.issled.}

\def\VIET{Vopr. Istor. Estestvozn. Tekh.} 
\def\viet{Вопр. ист. естествозн. и техн.}

\def\steklov{Тр. МИАН}
\def\steklovv{Trudy mat. inst. Steklov}
\def\Steklov{Proc. Steklov Inst. Math.}

 {\footnotesize
 
 \section*{Remarks to bibliography}
 
 The biblography consists of two parts. The first part [1]--\cite{Moscow} 
 is a collection of references to this work. The second part 
 contains bibliographical sources of non-wikipedic type about
 heroes of our story.
 
 \sm
 
 {\bf Some  Russian journals, periodicals, and newspapers of XX centure.}
 
 %\stepcounter{nomer}

 \nomer Автоматика и телемеханика//Avtomatika i telemekhanika 
 \newline
 (Automation and Remote Control)
 \newline
 A journal since 1936.
 
 \nomer Большевик//Bolshevik,
 \newline
  1924-1952, the main thick journal of the Communist party.
   \newline
 Continued as Коммунист//Kommunist.
 
 \nomer Бюллетень ВАРНИТСО%
 \footnote{Всесоюзная ассоциация работников науки и техники для содействия социалистическому строительству СССР (All-Union assiciation for a help
 	to the socialist construction in the USSR)}//Bulleten VARNITSO
 \newline
 Polytical and scientific-organizational
 journal,  1929-1931,
 continued as \ref{j:front}).

\nomer Вестник Академии наук СССР//Vestnik akademii nauk SSSR 
\newline
(The bulletin of the Academy
 of sciences of the USSR).
 \newline 
 Scientific-organiztional journal, since 1932. After 1992 continued as
  \newline
 Вестник Российской академии наук//Vestnik Rossi{\u \i}sko{\u \i} akademii nauk
   \newline
 (Bulletin of the Russian Academy of Sciences)
 
\nomer Вестник Московского Университета// Vestnik Moskovskogo universiteta
  \newline
 (Bulletin of the Moscow University). A journal since 1946 with numerous series.

\nomer 
Вопросы истории естествознания и техники//Voprosy istorii estestvoznaniya i techniki
\newline 
(Problems of history of natural sciences and technics)
\newline
A periodical since 1956, a journal since 1980.

\nomer Вопросы философии//Voprosy filosofii
  \newline
(Problems of philosophy). 
  \newline
  A journal since 1947

\nomer Доклады АН СССР//Doklady Akademii Nauk SSSR,
\newline
An old journal with numerous moments of branching and
 variations of titles.
\newline
Translated as
 Doklady: mathematics.

\nomer Естествознание и марксизм//Estestvoznanie i marksizm
\newline
(Natural sciences 
and Marxism)
\newline
 A scientific-philosophical journal, 1929-1930.
 Continued in 1931-1932 as (anti-scientific journal)
\newline 
 'За марксистско-ленинское естествознание//Za marksitsko-leninskoe estestvoznanie 
  \newline
 (For Marxist-Leninist natural sciences)', 

\nomer Историко-математические исследования//Istoriko-matematicheskie issledovaniya
 \newline
 (Historical-mathematical investigations).
  \newline A periodical since 1948.

\nomer Изв\ест\ия Императорской Академ\ии Наукъ//Izvestiya Imperatorsko{\u \i} akademii
nauk
 \newline
(Bulletin d'Acad\'emie Imp\`eriale des Sciences)
 \newline
 A journal 1894-1906. It was splitted
after 1906, see the next line.
 \newline
Изв\ест\ия Императорской Академ\ии Наукъ,  VI сер\ия//Izvestiya Imperatorsko{\u \i} akademii
nauk, VI series
 \newline
 (Bulletin d'Acad\'emie Imp\`eriale des Sciences, VI series)
  \newline A journal 1907-1916. Continued as 
   \newline
Изв\ест\ия Росс\ийской Академ\ии Наукъ. VI сер\ия
//Izvestiya Rossi{\u \i}sko{\u \i} Akademii Nauk, VI seriya
 \newline
 (Proceedings of the Russian
Academy of Sciences, VI series), 1917-1925;
 \newline
 Continued as
Известия Академии наук СССР. VI серия//Izvestiya Akademii nauk SSSR, VI series
 \newline
(Proceedings of the Academy of sciences of the USSR), 1926-1927. It was splitted
after 1927
 \newline
Известия Академии наук СССР. VII серия. Отделение физико-математических наук//
Izvestiya Akademii nauk SSSR, VII series. Otdelenie fiziko-matematicheskih nauk
 \newline
(Proceedings of the Academy of sciences of the USSR, VII series. Department of physical and mathematical sciences), 1928-1930.
 Continued as
  \newline
Известия Академии наук СССР. VII серия. Отделение математических и естественных наук//Izvestiya Akademii nauk SSSR, VII series. Otdelenie matematicheskih i estestvennyh nauk
 \newline
(Bulletin of the Academy of sciences of the USSR, VII series. Department of  mathematical and natural sciences), 1931-1935.
 It was splitted
after 1936
 \newline
Известия Академии наук СССР. Серия математическая//Izvestiya Akademii nauk SSSR, seriya matemeticheskaya, 1936-1991.
 \newline
 English translation: Mathematics of the USSR-Izvestiya since 1967.
  \newline
Continued as 
 \newline
Известия Российской академии наук. Серия математическая
 //Izvestiya Rossi{\u \i}sko{\u \i} akademii nauk, seriya matemeticheskaya].
  \newline
 English translation: Izvestiya: Mathematics.

\nomer Известия Научно-исследовательского инсти­тута математики и механики при Томском государственном университете//Izvestiya Nauchno-issledovatelskogo instituta matematiki
 i mehaniki pri Tomskom gosudarstvennom universitete
 \newline
  (Proceedings of the Research institute
 of Mathematics and Mechanics of the Tomsk State University)%
 \footnote{In Zentralblatt the periodical was reviewed as {\it Mitt. Forsch.-Inst. Math. Mech. Univ. Tomsk}.}, a periodical, 1935-1946.
 
\nomer История и методология естественных наук//Istoriya i metodologiya estestvennyh nauk 
 \newline
(History and methodology of natural sciences),
 \newline
a periodical 1960-?
 
\nomer Квант//Kvant
 \newline
 (Quant). A poplular journal for high-schoolers, since 1970

\nomer КИПР -- Книга и пролетарская революция //KIPR -- Kniga i proletarskaya revolutsiya
\newline
(Book and proletary revolution).
A bibliographical journal, 1932-1940

\nomer Математика в школе //Matematika v shkole
 \newline
 (Mathematics in school),
a  journal for teachers, since 1937.

\nomer Математика. Кибернетика.
 \newline
(Mathematics. Сybernetic). A scientific-popular periodical,
1967-1991.

\nomer Математическое просвещение//Matematicheskoe prosveshchenie
 \newline
 (Mathematical	enlightenment). A periodical 1934-1938. Continued in
1957-1961 and since 1997.

\nomer Математический сборник//Matematicheski{\u \i} sbornik, a journal since 1866.
 \newline
 English translation {\it Mathematics of the USSR. Sbornik}, 1967-1993;
{\it Sbornik. Mathematics} since 1993

\nomer Наука и жизнь//Nauka i zhizn
 \newline
 (Science and life), scientific-popular journal, since 1934.

\nomer Под знаменем марксизма//Pod znamenem marksizma
 \newline
 (Under banner of Marxism), a philosophical journal, 1922-1944. . 

\nomer Природа//Priroda 
\newline
(Nature). A popular scientific journal since 1912.

\nomer СОРЕНА%
\footnote{Социалистическая реконструкция и наука//Socialistic reconstruction and science}//SORENA
\newline
 Popular scientific and social journal, 1931-1936.

\nomer Теория вероятностей и ее применения//Teoriya veroyatnoste{\u \i} i ee
primeneniya]
\newline
 A journal since 1956, 
 \newline
 translated as
Theory of Probability and its Applications.

\nomer Труды Математического института имени В. А. Стеклова//
Trudy Matematicheskogo Instituta imeni V. A. Steklova
\newline 
A periodical since 1931,
\newline
 translated as 	Proceedings of the Steklov Institute of Mathematics since 1966.

\nomer Труды Московского математического общества// Trudy Moskovsogo matematicheskogo obshchestva.
\newline
A periodical since 1951
\newline
 Translated as
 'Transactions of the Moscow Mathematical Society'
 
 %\umn, \UMN, \RMS
\nomer Успехи математических наук//Uspekhi matematicheskih nauk,
\newline
 a periodical 1936-1944,
a regular journal since 1946;
\newline English translation {\it Russian mathematical surveys}
began 1960 

\nomer Учёные записки Московского университета//Uchenye zapiski Moskovskogo
universiteta
\newline
 (Transactions of the Moscow University).
 \newline
  A journal, 1833-1961
(with numerous variations of the title and series).

\nomer \label{j:front}
 Фронт науки и техники//Front nauki i techniki
 \newline
  (Front of Science and Technics).
 A polytical-scientific journal, 1931-1938.

}

%KBF3EXRP

 \tt
 
 Mathematical Department, University of Vienna;
 
 Institute for Theoretical and Experimental Physics, Moscow;
 
 MechMath.Department, Moscow State University;
 
 Institute for information transmission problems, Moscow.
 
 \sm
 
 e-mail: yurii.neretin(at)univie.ac.at 
 
 URL: http://www.mat.univie.ac.at/$\sim$neretin/

\end{document}